\documentclass[11pt,reqno,a4paper,oneside]{amsart}

\usepackage[utf8]{inputenc}
\usepackage[left=2.5cm,
            right=3cm,
            top=2.8cm,
            bottom=3cm]{geometry}

\RequirePackage{silence}
\usepackage{amsmath,amsthm,amssymb}
\usepackage{mathtools,esint,upref}
\usepackage{xfrac,xcolor}
\usepackage{bbm}
\usepackage{enumitem,comment}
\usepackage[hidelinks]{hyperref}
\usepackage{lmodern,cancel}
\usepackage[final]{pdfpages}
\usepackage{stmaryrd}

\newtheorem{theorem}{Theorem}[section]
\newtheorem*{theorem*}{Theorem}
\newtheorem{proposition}[theorem]{Proposition}
\newtheorem{lemma}[theorem]{Lemma}
\newtheorem{corollary}[theorem]{Corollary}
\newtheorem{assumption}[theorem]{Assumption}

\theoremstyle{definition}
\newtheorem{definition}[theorem]{Definition}
\theoremstyle{remark}
\newtheorem{remark}[theorem]{Remark}

\newtheorem{step}{Step}[subsection]
\newtheorem{case}{Case}[subsection]
\newtheorem{claim}{Claim}[subsection]

\numberwithin{equation}{subsection}

\DeclarePairedDelimiter{\norm}{\lVert}{\rVert}
\DeclarePairedDelimiter{\snorm}{\llbracket}{\rrbracket}

\newcommand{\R}{\mathbb{R}}
\newcommand{\N}{\mathbb{N}}
\newcommand{\Z}{\mathbb{Z}}
\newcommand{\T}{\mathbb{T}}
\newcommand{\bP}{\mathbb{P}}
\newcommand{\bS}{\mathbb{S}}
\newcommand{\fS}{\mathfrak{S}}
\newcommand{\balpha}{\boldsymbol{\alpha}}
\newcommand{\bbeta}{\boldsymbol{\beta}}
\newcommand{\bgamma}{\boldsymbol{\gamma}}
\newcommand{\bomega}{\boldsymbol{\omega}}

\newcommand{\one}{\mathbbm{1}}

\newcommand{\mr}{\mathring}

\newcommand*{\ini}{\ensuremath{\mathrm{in}}}
\newcommand*{\Id}{\ensuremath{\mathrm{Id}}}

\DeclareMathOperator{\divr}{div}
\DeclareMathOperator{\adj}{adj}
\DeclareMathOperator{\supp}{supp}

\newcommand{\resF}[1]{{\langle\kern-0.35ex\langle #1 
    \rangle\kern-0.35ex\rangle}}

\newcommand{\resN}[1]{{|\kern-0.2ex[ #1 
    ]\kern-0.2ex|}}

\newcommand{\vertiii}[1]{{\vert\kern-0.25ex\vert\kern-0.25ex\vert #1 
    \vert\kern-0.25ex\vert\kern-0.25ex\vert}}

\title[Optimal regularity in scalar anomalous dissipation]
      {Scalar anomalous dissipation and optimal regularity via iterated homogenization}

\author[J. Burczak]{Jan Burczak}
\address{Mathematisches Institut, Leipzig University, Augustusplatz 10, 04109 Leipzig, Germany and Max Planck Institute for Mathematics in the Sciences, Inselstrasse 22-26, 04103 Leipzig, Germany}
\email{burczak@math.uni-leipzig.de}

\author[L. Sz{\'{e}}kelyhidi]{ L\'{a}szl\'{o} Sz{\'{e}}kelyhidi, Jr.}
\thanks{LSz gratefully acknowledges the support of the  Deutsche Forschungsgemeinschaft (DFG, German Research Foundation) through GZ SZ 325/2-1.}
\address{Max Planck Institute for Mathematics in the Sciences, Inselstrasse 22-26, 04103 Leipzig, Germany}
\email{szekelyhidi@mis.mpg.de}

\author[B. Wu]{Bian Wu}
\address{Max Planck Institute for Mathematics in the Sciences, Inselstrasse 22-26, 04103 Leipzig, Germany}
\email{bian.wu@mis.mpg.de}

\begin{document}

\begin{abstract}
For any $\beta_0<1/3$ we construct divergence free vector fields in $ C_{x,t}^{\beta_0}$ and a sequence of diffusivities $\kappa_q \searrow 0$ such that, for an arbitrary initial datum from a low regularity class, the classical solution $\rho_q$ to the advection-diffusion equation exhibits anomalous dissipation along the sequence $\kappa_q$. At the same time $\rho_q$ remains uniformly bounded in  $C_t^{0} C_x^{\alpha_0}$, where $\beta_0 + 2\alpha_0<1$. Our result confirms a conjecture of Armstrong and Vicol \cite{ArmstrongVicol} and shows sharpness of the Obukhov-Corrsin threshold within the context of iterated homogenization. Our construction confirms time-homogeneity of the dissipation anomaly, as required in turbulence theory, and as a consequence we also obtain better time regularity for the scalar $\rho_q$ than the classical prediction of Yaglom. 
\end{abstract}

\maketitle

\setcounter{tocdepth}{1}

\tableofcontents


\section{Introduction}
We consider the advection-diffusion equation
\begin{equation}  \label{e:advection_diffusion}
\begin{aligned}
      \partial_t \rho_\kappa  +  u \cdot \nabla  \rho_\kappa  =&\, \kappa \Delta \rho_\kappa      &  \text{ in }\T^2 \times (0,1]\,, \\ 
      \rho_\kappa (\cdot, 0) =&\, \rho_{\ini} &  \text{ on } \T^2,
\end{aligned}
\end{equation}
where  $u: \T^2 \times [0,1] \rightarrow \R^2$ is a given divergence-free vector field and $\rho_{\ini}$ is an initial datum with zero mean. Throughout this paper we only consider data $u, \rho_{\ini}$ which is regular enough for global well-posedness of \eqref{e:advection_diffusion} in the energy class to hold (e.g.~$\rho_{\ini} \in L^2$ and $u \in L^\infty L^d$, $d=2$, suffices). In this case the strong energy identity  
\begin{equation}\label{e:LHidentity}
\frac12\| \rho_\kappa (T) \|^2_{L^2} + \kappa \int_{0}^{T}\| \nabla \rho_\kappa \|^2_{L^2} \,dt 
            = \frac12\| \rho_\ini \|^2_{L^2}\quad \textrm{ for all }T\in [0,1] 
\end{equation}
holds. Furthermore, the solution to \eqref{e:advection_diffusion} is unique and as smooth as the datum allows. In what follows, we refer to such solutions as 'classical solution' and note in passing that weaker notions of solution also exist (see e.g.\cite{Modena2018}).

\subsection{Main result}
In this paper we prove the following theorem:
\begin{theorem}   \label{t:mainHld}
Let $0<\beta_0<\sfrac{1}{3}$ and $\alpha_0>0$ such that $2\alpha_0 + \beta_0 < 1$. There exists a divergence-free vector field 
\begin{equation}\label{e:1:t:mainHld}
 u \in C^{\beta_0}(\T^2\times [0,1])
\end{equation}
and a sequence of diffusivities $\{\kappa_q\}_{q \in \N}$, $\kappa_q\searrow 0$, such that for any zero-mean $\rho_{\ini}\in H^s(\T^2)$ with $s>\frac{1+\beta_0}{2}$ the classical solution $\rho_{q}:=\rho_{\kappa_q}$ of \eqref{e:advection_diffusion} satisfies for any $T\in (0,1]$
\begin{equation}     \label{e:2:t:mainHld}
      \lim_{q \rightarrow \infty} \kappa_q \int_{0}^T\|\nabla \rho_q \|_{L^2}^2 \,dt 
            \geq c_0 \| \rho_{\ini} \|^2_{L^2}\,.
\end{equation}
Moreover, if in addition $\rho_{\ini}\in H^{1+\alpha'}(\T^2)$ for some $\alpha'>\alpha_0$, then 
\begin{equation}     \label{e:4:t:mainHld}
      \sup_{q \in \N} \sup_{t\in [0,1]}\| \rho_q \|_{C^{\alpha_0}} 
            \le c_1 \| \rho_{\ini} \|_{H^{1+\alpha'}}\,.
\end{equation}
The constant $c_0 > 0$ depends only on $u$, $T$ and the ratio $\frac{\|\rho_{\ini}\|_2} {\|\nabla^s \rho_{\ini}\|_2}$, whereas the constant $c_1>0$ depends only on $u$.
\end{theorem}

Theorem \ref{e:1:t:mainHld} is a substantial extension of the main result proved in \cite{ArmstrongVicol} via iterated homogenization: the first statement above \eqref{e:2:t:mainHld} is the subject of \cite[Theorem 1.1]{ArmstrongVicol}, but the second statement above \eqref{e:4:t:mainHld} has merely been conjectured to hold in \cite[Section 1.3]{ArmstrongVicol}, and is in fact the main novelty of our work. We also remark that the first statement above in the 3D case has been substantially improved in our prior work \cite{burczak2023anomalous}, where (i) $u$ is in addition a weak solution of the 3D Euler equations satisfying an $h$-principle, and (ii) using the explicit dependence of $c_0$ on the data and based on the fluctuation-dissipation relation introduced in \cite{Drivas_Eyink_2017}, estimate \eqref{e:2:t:mainHld} is complemented with a lower bound on the average rate of separation between two realizations of the associated stochastic differential equation of \eqref{e:advection_diffusion}, in agreement with the Richardson law of dispersion.

\subsection{Context}\label{s:context}
In hydrodynamic turbulence theory, the advection-diffusion equation
\eqref{e:advection_diffusion} serves as an important toy model for 'turbulent transport'. For example, the classical K41 theory of turbulence based on the 3D Navier-Stokes equations finds its analogy in the Corrsin-Obukhov-Yaglom theory for the advection-diffusion equation. There are two foundational postulates common to both of these theories \cite[Section 6]{Frisch_1995}\footnote{Here we focus on the possible mathematical implications of these postulates, notwithstanding well-founded doubts as to their validity in a physical context, such as scaling corrections due to intermittency or the recent work \cite{iyer2025zerothlawturbulence} challenging the zeroth law.}:
\begin{itemize}
    \item[(I)] finite non-vanishing mean rate of dissipation in the inviscid limit;
    \item[(II)] homogeneity and isotropy on small scales.
\end{itemize}
In the context of advection-diffusion \eqref{e:advection_diffusion} postulate (I) translates, in its most basic form, into the statement that 
\begin{equation}\label{f:ead}
\inf_{\kappa>0}\kappa \int_{0}^{1}\|\nabla \rho_\kappa \|_{L^2}^2\,dt 
            \geq c_0    
\end{equation}
for some constant $c_0>0$, and is often referred to as the 'zeroth law of scalar turbulence' or 'anomalous dissipation'. Examples of divergence-free vector fields $u$ for which one can verify \eqref{f:ead} have been constructed in \cite{DriElgIyeJeo22} and shortly after in \cite{CCS,ElgLis23} with optimal regularity properties (see Section \ref{ss:regularity} below for a more in-depth discussion). The drawback of these constructions is that they do not satisfy postulate (II). More precisely, the time-translation invariance of the underlying model in connection with postulate (II) requires, at the least, the strengthening of \eqref{f:ead} to 
\begin{equation}\label{f:ads}
\inf_{\kappa>0}\kappa \int_{0}^{T}\|\nabla \rho_\kappa \|_{L^2}^2\,dt 
            \geq c_0(T)\quad\forall T>0,    
\end{equation}
where $c_0(T)>0$. The works \cite{ArmstrongVicol,burczak2023anomalous} as well as present work almost achieve \eqref{f:ads}, the only downside being that the infimum $\inf_{\kappa>0}$ has to be replaced by $\limsup_{\kappa\to 0}$ (indeed, it seems very difficult if not impossible to upgrade the $\limsup$ to $\inf$ within the iterated homogenization framework used here). 

\subsubsection{Optimal regularity}\label{ss:regularity}
The central prediction of K41 theory for Navier-Stokes turbulence, based on the hypothesis of anomalous dissipation, is a precise scaling law for the energy spectrum in the inertial range: $E(k)\sim c_0^{2/3}k^{-5/3}$. A mathematical formulation of this prediction has famously been put forward by L. Onsager in 1949, which became known as Onsager's conjecture and is concerned with 'ideal turbulence': the threshold regularity for conservation of energy in the Euler equations is precisely $C^{1/3}$ on the H\"older scale (corresponding to monofractal scaling). We refer to \cite{Eyink_survey} for an excellent recent survey on Onsager's ideal turbulence and the current state of the art. The corresponding prediction in the context of passive scalar advection-diffusion is the Obukhov-Corrsin law and can be stated in the most basic form as follows \cite{Obu49,Cor51}, see also \cite{DriElgIyeJeo22}: if the advecting velocity field $u\in C^{\beta}$ is H\"older continuous with exponent $\alpha>0$, then the solution $\rho_\kappa$ of \eqref{e:advection_diffusion} should be $\rho_\kappa\in C^{\alpha}$ uniformly in $\kappa>0$, with  $\beta+2\alpha=1$, and not better. A more precise form of this prediction was obtained by Yaglom \cite{Yaglom1949}, see also \cite{CCS}, stating that the criticality threshold for anomalous dissipation is 
\begin{equation}\label{e:yaglom}
    u \in L_t^p C_x^\beta, \; \rho_\kappa  \in L_t^r C_x^\alpha, \quad \text{uniformly in } \kappa,  \qquad \frac{1}{p} + \frac{2}{r} = 1, \;\beta + 2 \alpha =1.
\end{equation}
The regular side of the Yaglom relation, namely the statement that no anomalous dissipation is possible provided
\begin{equation}\label{e:yaglom_regular}
    u \in L_t^p C_x^\beta, \rho_\kappa  \in L_t^r C_x^\alpha, \quad \text{uniformly in } \kappa,  \qquad \frac{1}{p} + \frac{2}{r} \le 1, \;\beta + 2 \alpha >1
\end{equation}
is known to hold, see \cite{Eyi96,ConPro93,ConPro94,CCS,DriElgIyeJeo22}, with interesting connections to non-classical transport theory. Conversely, sharpness of the Yaglom regime \eqref{e:yaglom_regular} has been obtained in \cite{CCS} by Colombo, Crippa and Sorella. More precisely, in \cite{CCS} the authors prove that whenever $\sfrac{1}{p} + \sfrac{2}{r} = 1$ and $\beta + 2 \alpha <1$ with $2\leq r\leq 4$, there exists initial data $\rho_{\ini}\in C^\infty$ 
and divergence-free $u\in L_t^p C_x^\beta$ such that the classical solutions of \eqref{e:advection_diffusion} satisfy \eqref{f:ead} and 
\begin{align}     
      \limsup_{\kappa \to 0} \| \rho_{\kappa} \|_{L^r_t C^{\alpha}_x (\T^2 \times [0,1])} 
            \le C_0. 
\end{align}
This result of \cite{CCS} has been extended by Elgindi and Liss in \cite{ElgLis23} (based on ideas in \cite{DriElgIyeJeo22}) at the endpoint time case $r=2$, $p=\infty$ (in fact with velocity even continuous in time) by allowing for arbitrary smooth initial data. 

Compared to these works, the novelty of our contribution is twofold: firstly, our estimates \eqref{e:4:t:mainHld}, \eqref{e:1:t:mainHld} not only saturate the Obukhov-Corrsin relation from below, but in fact yield better time regularity. 
More importantly, we are able to achieve this in combination with the stronger, time-homogeneous form of anomalous dissipation \eqref{e:2:t:mainHld}. In this respect the work of Armstrong and Vicol was a major step forward in the field as it not only achieved \eqref{e:2:t:mainHld} in agreement with postulate (II), but also identified iterated homogenization as the correct mechanism behind the phenomenon of anomalous dissipation (in agreement with the classical idea of eddy diffusivity). However, the authors of \cite{ArmstrongVicol} left open the question of uniform regularity of the density, and conjectured that with their approach, combined with “large-scale regularity” techniques developed in quantitative homogenization theory, one can reach $C^0_tC^\frac{1-\beta}{2}_x$ regularity of $\rho_q$ uniformly in $q$ (in agreement with the  Obukhov-Corrsin prediction). 

In this paper we confirm this conjecture. Our methods complement the ideas introduced in \cite{ArmstrongVicol}, but we do not follow the 'large scale regularity' approach suggested there. Instead, we develop a novel framework for iterated quantitative homogenization for high-order correctors. 

\subsubsection{Outlook}\label{ss:outlook}

We believe that this paper does not only contain a new result in the fine analysis of the advection-diffusion equation and in anomalous dissipation of scalar variance, but provides a set of useful and sharp technical tools. In particular, we believe these tools can allow to gain insights into the following themes, which we plan to address in a subsequent work: 
\begin{enumerate}[wide,labelwidth=!,labelindent=0pt,label=\textsc{(P.\arabic*)}]

\item \label{pro:a1} \textit{Spatial homogeneity of anomalous dissipation.} 
 We expect the anomalous dissipation constructed in Theorem \ref{t:mainHld} to be homogeneous both in space and in time, which would amount to complete justification of postulate (II). In technical terms, we expect that for any nonempty open parabolic cylinder $E \subset \T^2 \times [0,1]$ and generic initial datum,
\begin{align}     
      \liminf_{q \rightarrow \infty} \kappa_q \iint_{E} | \nabla \rho_q |^2 dxdt 
            \geq c(E,\rho_{\ini},u).
\end{align}
A further interesting question is then to characterize the dependence of $c(E,\rho_{\ini},u)$ on $E$, $\rho_{\ini}$ and $u$. We expect that the techniques in this paper will help such characterization.

\item \label{pro:a2} \textit{H-principle and Euler equations.} Based on the work \cite{burczak2023anomalous} bridging iterated homogenization in \cite{ArmstrongVicol} and convex integration in \cite{BDSV}, we expect the existence of a large class of turbulent vector fields $u$, exhibiting both homogeneous anomalous dissipation (in the sense of \ref{pro:a1}) and optimal H\"older regularity for scalar fields, which are dense in an appropriate sense. This would yield an h-principle for these vector fields. Moreover, we expect they can be constructed as H\"older continuous weak solutions to the incompressible Euler equations, given the insights from convex integration \cite{DSz13,Isett2018,BDSV,GiriRadu}.
\end{enumerate}

\subsubsection{Related works}
In the recent past there has been a surge of activity concerning the mathematical foundation of the phenomenon of anomalous dissipation and the closely related issue of mixing. Although a rigorous approach for the nonlinear Navier-Stokes case is still widely open, the case of the advection-diffusion equation \eqref{e:advection_diffusion} has seen several very interesting advances. Beyond the works already mentioned we wish to point out \cite{CZDelElg20,ElgLisMat23,Yao2017,HCRR,HuyTit,JohSor,DRDII,TaoZwo,Ziz}. Finally, we mention the excellent surveys \cite{Drivas_survey,Eyink_survey} for a wider overview and context of the subject.

\subsection{Challenges and main ideas}    \label{ss:challenges}

Our construction of the vector field $u$ in Theorem \ref{t:mainHld} essentially follows the approach of \cite{ArmstrongVicol,burczak2023anomalous}, arising iteratively from a sequence $\{u_q\}$, where $u_{q+1} = u_q + w_{q+1}$ and $w_{q+1}$ is a highly oscillatory perturbation with spatial frequency $\lambda_{q+1}$ and temporal frequency $\mu_{q+1}$. The perturbation $w_{q+1}$ is transported by $u_q$ locally in time intervals of length $\mu_{q+1}^{-1}$ - see for instance Section \ref{ss:vectorfield}. In a first approximation, we can think of such a perturbation as taking the form
\begin{align}     \label{e:2:genIntro}
      w_{q+1}(x,t) =&\, a_q(x,t) \eta(\mu_{q+1}t) W\big( \lambda_{q+1} \Phi_{q}(x,t) \big).
\end{align}
 In this way we obtain a sequence of advection-diffusion equations, with advecting vector field $u_q$ and diffusion coefficient $\kappa_q\sim\lambda_q^{-\theta}$ with $\theta\sim 1+\beta$, such that at each level $q$ the deformation due to transport and the diffusion are balanced. In oversimplified terms the corresponding solutions $\rho_q$ arise then from iteratively homogenizing equation $q+1$ to equation $q$, see \cite[Section 2]{burczak2023anomalous}. Looking at the explicit form of the first order corrector and basic energy estimates in \cite{ArmstrongVicol,burczak2023anomalous} suggests then the following naive approach to regularity: obtain inductive information on higher regularity of $\rho_{q+1}-\rho_q$ in the form of Campanato bounds, and then to apply the tools developed for the "large scale regularity" theory in stochastic homogenization \cite{ArmstrongBook}. However, this strategy seems to run into two fundamental difficulties: 
\begin{itemize}
\item on small scales $\sim\lambda_q^{-1}$ we have to deal with a lack of sufficiently large scale separation: The scale of separation is at most $\sfrac{\lambda_q}{\lambda_{q+1}}$, which is very mild compared to the degeneracy in ellipticity $\kappa_{q+1}$. In order to achieve the optimal H\"older regularity, we thus need to compute correctors in the homogenization procedure to a certain very high order $Q$. 
\item in contrast, on larger scales we see that the background transport $u_q$ remains at order $1$, while the ellipticity coefficient $\kappa_{q}$ degenerates to zero. As well shall see later, this "high contrast environment" means that energy estimates and Gr\"onwall or Schauder estimates are far from being sufficient for pursuing any useful uniform H\"older regularity of $\rho_q$.
\end{itemize}
Although the first issue is possibly more technical than fundamental, since high order expansions are by now well known in periodic homogenization and also in quasi-periodic or stochastic settings \cite{Gu2017high}, the second issue seems to be a fundamental feature of our setup, prohibiting the use of "large-scale regularity" techniques \cite{AvLin1987,AvLin1989}. Instead, we are forced to develop a version of iterated homogenization which includes high order expansions. In turn this leads to new challenges, which are of a different nature from those encountered in quantitative stochastic homogenization \cite{GloriaOtto,FischerOtto,ArmstrongBook,Armstrong2021homogenization,Armstrong2023Large}.

\subsection{Main ingredients of our approach}

To prove Theorem \ref{t:mainHld}, we establish a regularity theory for these turbulent vector fields via higher order correctors. Below, we summarize the main new ingredients in this work.

\subsubsection{Expanding the flow maps to reveal higher order periodicity in time}       \label{sss:genIntro}

Even though, as explained above, our velocity fluctuations $w_{q+1}$ are heuristically expected to be periodic in space and time as in \eqref{e:2:genIntro}, this is only true upto first order. In fact, the nature of temporal cutoffs and the use of local-in-time flow maps means that the actual fluctuations are not periodic, neither in space nor in time. For the purposes of \cite{ArmstrongVicol,burczak2023anomalous} this is sufficient, as effects from lower order terms have a negligible contribution in the energy norm $L^\infty_t L^2_x$. However, for the higher order corrector expansion these terms do become significant and they may lead to non-negligible contributions. 

To recover periodicity in space and time, we define a globally defined approximate flow map $\Phi_q$ arising from a temporal Taylor expansion of the local flow maps 
\begin{align}     \label{e:6:genIntro}
      \Phi_{q,\iota} = x + \sum_{m=1}^{\infty} \frac{1}{m!} ( t_{q,\iota} - t )^m D_{q,t}^{m-1} u_q,
\end{align}
where $D_{q,t} := \partial_t + u_q \cdot \nabla$ is the transport derivative associated to $u_q$ and $t_{q,\iota}$ is the center of the support of $\eta_\iota$. Truncating this expansion at order $m_*$, we can retrieve
\begin{align}     \label{e:10:genIntro}
      w_{q+1}(x,t) =&\, a_q(x,t) \sum_{m=0}^{m_*} \eta^{(m)} ( \mu_{q+1} t ) w \big( \lambda_{q+1} \Phi_{q}(x,t) \big) 
            + O ( \lambda_{q+1}^{-m_*\gamma_0} )
\end{align}
for some small positive exponent $\gamma_0 > 0$. The flow map $\Phi_q:\T^2 \times [0,1] \rightarrow \T $ is globally defined in time and $\eta^{(m)}: \T \rightarrow \R $ are globally periodic. 

In \eqref{e:10:genIntro}, $a_q$ contains large spatial scale $(\geq \lambda_q^{-1})$ and large temporal scale $(\geq \mu_q^{-1})$, $\eta^{(m)}(\mu_{q+1} t)$ contains fast temporal scale $\mu_{q+1}^{-1}$ and $w( \lambda_{q+1} \Phi_q )$ contains fast spatial scale $\lambda_{q+1}^{-1}$. The expansion above reduces our problem, at each iteration step, to a cleaner homogenization problem up to arbitrarily small perturbations, at least formally, which is very desirable.

\subsubsection{Single space-time hybrid homogenization problem}     \label{sss:introHom}

With above expansion on flow maps, at each iteration step, we can focus on the following over-simplified advection diffusion equation in $\T^2 \times [0,1]$,
\begin{equation}  \label{e:2:introHom}  
\begin{split}
      D_{q,t} \rho_{q+1} - \divr \bigg(  
            \begin{bmatrix}
                  \kappa_{q+1}, &\,\psi \big( \lambda_{q+1} \Phi_{q}(x,t), \mu_{q+1} t \big) \\
                  -\psi \big( \lambda_{q+1} \Phi_{q}(x,t), \mu_{q+1} t \big), &\, \kappa_{q+1}
            \end{bmatrix} \nabla \rho_{q+1} \bigg) =&\, 0,       \\ 
      \rho_{q+1} (\cdot, 0) =&\, \rho_{\ini}.
\end{split}
\end{equation}

From classical homogenization theory \cite{BenLioPap78,JiKoOl,ShenBook}, we know there are three types of correctors that we have to consider. The first type is the spatial corrector, which captures the effect of fast spatial oscillation $\lambda_{q+1}^{-1}$. The second type is the temporal corrector, which captures the effect of temporal oscillation $\mu_{q+1}^{-1}$. The third type is the boundary corrector, which is needed to enforce initial condition on parabolic boundary $\{t=0\}$. 

Obtaining a high order expansion for such equations seems highly delicate, because now we have 2 distinct small scales $\lambda_{q+1}^{-1}$ and $\mu_{q+1}^{-1}$ in addition to the large scale, with no apriori well-defined ordering among different powers. In other words, there is not a unqiue way to even compute a formal asymptotic expansion. We present an algorithmic way \textit{to simultaneously compute all three types of correctors up to arbitrarily high order.} The details are presented in Section \ref{ss:stateExpInertial} and Section \ref{ss:proof_expansion}.

\subsubsection{Chain of homogenization problems, renormalized differential operators and polynomial-derivative space}     \label{sss:intro:renormDiffOp}

As mentioned above, at each step $q$ we require an expansion of $\rho_{q+1}$ up to very high order $Q$, forcing us to solve $Q$ different homogenization problems at step $q$. But then this means we need to solve $Q^2$ homogenization problems at step $q-1$. When $q$ goes to infinity, this would require us to solve infinitely many homogenization problems at step $1$, which is not feasible.

To overcome this issue, we solve a chain of $Q^2$ problems at step $q+1$. We solve the first $Q$ problems by expansion techniques in Section \ref{sss:introHom} and the remaining $Q^2-Q$ problems by the energy method. This is possible because of a special structure of the chains, consisting of two parts: the divergence form of the $Q$ homogenized problems and a special algebraic structure of function coefficients, introduced by the deformation of the flow maps $\nabla \Phi_q$. In particular, the function coefficients in the homogenized problems are polynomials of $\nabla \Phi_q$ and its derivatives. To efficiently implement estimates and capture their algebraic structure, we use renormalized differential operators and introduce a special polynomial-derivative space. The purpose is to be able to efficiently propagate both their algebraic structure and the estimates on chains of homogenization problems when passing from step $q+1$ to step $q$. These tools are developed in Section \ref{ss:poly_space}, Section \ref{ss:r_estimates} and Section \ref{s:decomposition}. The homogenization schemes are then presented in Section \ref{s:homInertial} and \ref{s:homDissip}.

\subsubsection{Inertial range, dissipative range and scale-by-scale decomposition of initial datum}

In Section \ref{sss:genIntro}, Section \ref{sss:introHom} and Section \ref{sss:intro:renormDiffOp}, we have described the main new ingredients in a single homogenization step. Next, we briefly describe how to implement these iteratively in $q$ and how to deal with arbitrary initial datum in $H^{1+\alpha'}$ containing small scales. There are two key differences with respect to \cite{ArmstrongVicol,burczak2023anomalous}.

First, instead of viewing the initial datum containing large scales as a single entity, we decompose the initial datum into infinitely many energy shells a'la Littlewood-Paley adapted to the scales $\{\lambda_q\}_q$. For a single shell component, we inductively apply the homogenization expansion to resolve its energy spectrum. 

Secondly, we do homogenization separately in the inertial range and the dissipative range. The inertial range is where the diffusivity is significantly enhanced, amounting to a large eddy effect. The dissipative range contains scales below the Kolmogorov scale. This approach is required largely due to the challenges mentioned in Section \ref{ss:challenges}.

\subsection{Notations} 

We use usual Einstein summation convention.

\subsubsection{Multi-index}   \label{sss:ntt:multi_ind}
We use $\mathcal{I}$ to denote the set of multi-indices consisting of indices $x_1$, $x_2$ and $t$, i.e.
\begin{align}
      \mathcal{I} := \bigcup_{n \in \N} \{ 1, 2, t \}^{n}.
\end{align}
If $\balpha \in \mathcal{I}$, then there is a unique $n_0 \in \N$ with $\balpha \in \{ 1, 2, t \}^{n_0}$. We use $\balpha(i)$ to denote $i$-th element in $\balpha$, then
\begin{align}
      \balpha = ( \balpha(1), \balpha(2), \ldots, \balpha( n_0 ) ).
\end{align}
We use $| \balpha |_x$ to denote the cardinality of $x_1$ and $x_2$, and $| \balpha |_t$ to denote the cardinality of $t$. Moreover, we define
\begin{align*}
      |\balpha| :=&\, |\balpha|_x + |\balpha|_t,      \quad
      [\balpha] := |\balpha|_x + 2 |\balpha|_t,
\end{align*}
and
\begin{align*}
      \mathcal{I}_x =&\, \{ \balpha \in \mathcal{I} \mid |\balpha|_t = 0 \},  \\
      \mathcal{I}(i,j) =&\, \{ \balpha \in \mathcal{I} \mid 
            |\balpha|_x \geq i, |\balpha|_t \geq j \},
\end{align*}
and a special collection of differentiation indices via
\begin{equation}     \label{e:indexDNabla}
\begin{split}
      \mathcal{I}_* = \big\{ \balpha \in \mathcal{I} \mid \balpha( |\balpha| ) \in \{1,2\} \big\}.
\end{split}
\end{equation}
Here, $\mathcal{I}(i,j)$ consists of indices with at least $i$ $x$-elements and $j$ $t$-elements.

\subsubsection{Differentiation and transport derivative}    \label{sss:ntt:deri}
We use $\partial^{\balpha}$ or $\nabla^{\balpha}$ to denote the usual partial derivative with index $\balpha$ which only contains spatial derivatives. For a smooth vector $v: \T^2 \times [0,1]$, we define the transport derivative associated to $v$ to be $D_t := \partial_t + v \cdot \nabla$. Using $D_{1} := \partial_{1}$ and $D_{2} := \partial_{2}$, we define, for any $\balpha$ with $|\balpha| = n \in \N$,
\begin{align}
      D^{\balpha} = D^{\balpha(1)} D^{\balpha(2)} \ldots D^{\balpha(n)}.
\end{align}
We call $D^{\balpha}$ the generalized derivative associated to vector field $v$ with index $\balpha$.

\begin{remark}    \label{r:IStar}
For differentiation indices in $\mathcal{I}_*$ defined in \eqref{e:indexDNabla}, the last element is $x$-element, hence the first derivative being applied is always space derivative. 
\end{remark}

In dimension $2$, we define $\nabla^\perp := [\partial_2, -\partial_1]^T$. A related and important set of notations \textit{renormalized differential operator} is introduced in Section \ref{ss:poly_space}.
For a matrix-valued function $H : \T^2 \rightarrow \R^{2 \times 2}$, we define $\divr H := \partial_k H_{jk}$. Then if $H$ is antisymmetric, we have $\divr H = \nabla^\perp H_{12}$.
For a function $\Phi : \T^2 \rightarrow \R^2$, we define $\nabla \Phi: \T^2 \rightarrow \R^{2 \times 2}$ by $[\nabla \Phi]_{ij} := \partial_j \Phi_i$. 
We also define
\begin{align}
      \adj \begin{bmatrix} a & b \\ c & d  \end{bmatrix}
      = \begin{bmatrix} d & -b \\ -c & a  \end{bmatrix}.
\end{align}
As a consequence, $\adj \circ \adj$ is an identity operator and 
\begin{align}     \label{e:adjointPerp}
      (A v)^\perp = ( \adj^T A ) v^\perp,     \quad \text{for any } A \in \R^{2 \times 2}, v \in \R^2.
\end{align}

\subsubsection{Oscillatory variables in homogenization}     \label{sss:ntt:hom}
$( \cdot )$ and $[ \cdot ]$ are normal brackets. $\xi$ and $\tau$ are periodic variables $\xi \in \T^2$ and $\tau \in \T$. $\langle \cdot \rangle_\xi$ and $\langle \cdot \rangle_\tau$ are the integrals in $\xi$ and $\tau$ respectively. $\langle \cdot \rangle_{\xi,\tau}$ is the integral in both $\xi$ and $\tau$. 

For periodic variable $\tau \in \T$, we define $\partial_\tau^{-1}$ to be the inverse operation of differentiation, mapping a function with zero $\tau$-average to another function with zero $\tau$-average. The analogous definition applies to $\partial_{\xi_1}^{-1}$ and $\partial_{\xi_2}^{-1}$. We also define $\Delta_\xi^{-1}$ to be the inverse Laplacian, mapping a function with zero $\xi$-average to another function with zero $\xi$-average.

\subsubsection{Constants}     \label{sss:ntt:cst}
We use $\lesssim$ to denote $\leq$ up to an absolute constant. For any integer $p$, we use $\lesssim_p$ to denote $\leq$ up to an absolute constant to power $p$, i.e.
\begin{align}
      \cdot \lesssim_p \cdot
      \qquad \text{means} \qquad
      \cdot \leq c^p \cdot
\end{align}
for some absolute constant $c$.

\subsubsection{Norms}         \label{sss:ntt:norm}
For $1 \leq p \leq +\infty$, we use $\|\cdot\|_p$ to denote the respective $L^p$ norm of the function domain. We also use the following simplified norm notation for function $\rho: \T^2 \times [0,1] \rightarrow \R $,
\begin{align}
      \| \rho \| := \| \rho \|_{L_t^\infty L_x^2}
\end{align}

\newpage

\section{The vector field}    \label{s:vectorField}
In this section we construct our drift and provide certain useful preliminaries.

\subsection{Absolute parameters}       \label{ss:parameters}

By \textit{absolute parameters}, we mean parameters that merely depend on $\beta, \alpha_0, \alpha', s$ in the statement of Theorem \ref{t:mainHld}. They are fixed at the beginning of the paper and remain unchanged throughout the paper.

For $ \beta, \alpha_0, \alpha', s \in (0,1) $ with
\begin{align*}
      \beta < \frac{1}{3}, \quad
      2\alpha_0 + \beta < 1,   \quad
      \alpha' > \alpha_0, \quad 
      s > \frac{1+\beta}{2},
\end{align*}
we fix $b$ satisfying
\begin{align}     \label{e:bConsParameter}
      b \in \Big( 1, \frac{11}{10} \Big), \quad
      \alpha' \geq b(\alpha_0+1)-1, \quad
      s \geq \frac{b(1+\beta)}{(b+1)(2-b)},
\end{align}
and
\begin{align}     \label{e:2:bConsParameter}
      2\alpha_0 + \beta + \frac{1+\beta}{b+1} (b-1) \alpha_0 \leq 1.
\end{align}
We also fix a large enough constant $\lambda_0$ and define the following parameters:
\begin{align}     \label{e:Parameter}
      \lambda_q := \big\lceil \lambda_0^{(b^q)} \big\rceil, \quad 
      \mu_q := \lambda_q^{e}, \quad 
      \delta_q:= \lambda_q^{-2\beta}, \quad 
      \kappa_q := \lambda_q^{-\theta},
\end{align}
\begin{align}     \label{e:Exponent}
      \theta = \frac{2b}{b+1}(1+\beta),
\end{align}
then
\begin{align}     \label{e:0:mScaleRela}
      \frac{\delta_q^{\sfrac12}} {\kappa_q\lambda_q}
      = \lambda_q^{ \frac{1+\beta}{b+1} (b-1) }.
\end{align}

We also fix several relatively small exponents (absolute constants): $\gamma_I, \gamma_R > 0$ for renormalized differentiation and interpolation, $\gamma_S$ for auxiliary scale control in homogenization, $\gamma > 0$ for the finest scale gap in homogenization. They satisfy the following relation
\begin{align}     \label{e:smallGammaR}
      4(\gamma_I + \gamma_R) \leq \gamma_S \leq \frac{b-1}{b+1} ,  \quad 
      3 b \gamma \leq \min\{ \gamma_I, \gamma_R \} .
\end{align}
Then we define
\begin{align}     \label{e:mrParameter}
      \mr \lambda_q := \lambda_q^{1+b\gamma_I} \Bigg( \frac{\delta_q^{\sfrac12}} {\kappa_q\lambda_q} \Bigg)^{\sfrac12} , \quad 
      \mr \mu_q := \delta_q^{\sfrac12} \lambda_q^{1+2b\gamma_I}, \quad 
      \varepsilon_q := \frac{ \delta_q^{\sfrac12} \lambda_q^{1+2b\gamma_R} } {\mu_{q+1}}.
\end{align}
This leads to the following relation      
\begin{align}
      \mr \mu_q \lesssim \kappa_q \mr \lambda_q^2 \lesssim \mr \mu_q.   \label{e:1:mScaleRela}
\end{align}

We require the following relations:
\begin{align}     
      \frac{ \delta_{q+1}^{\sfrac12} \mr\lambda_q }{ \kappa_{q+1} \lambda_{q+1}^2 } \cdot \lambda_{q+1}^{2\gamma_S}
            \leq&\, 1,     \label{e:2:mScaleRela} \\ 
      \max \Bigg\{ \frac{\mu_{q+1}} {\delta_{q+1}^{\sfrac12}\mr{\lambda}_q}, 
            \frac{ \varepsilon_q \kappa_{q+1} \lambda_{q+1}^2 } { \delta^{\sfrac12}_{q+1} \mr{\lambda}_q }, 
            \frac{ \mr{\mu}_q } { \varepsilon_q \delta_{q+1}^{\sfrac12} \mr{\lambda}_q }, 
            \frac{ \kappa_{q+1}\lambda_{q+1}^2 } { \varepsilon_q \delta_{q+1}^{\sfrac12} \lambda_{q+1} } \Bigg\} \leq&\, \lambda^{\gamma_S},   \label{e:4:mScaleRela}
\end{align}
A sufficient condition for \eqref{e:2:mScaleRela} is that
\begin{align}     \label{e:12:mScaleRela}
      3b\gamma_S \leq \frac{1-(2b+1)\beta} {2(b+1)}.
\end{align}
A sufficient condition for \eqref{e:4:mScaleRela} is \eqref{e:smallGammaR}, \eqref{e:12:mScaleRela} and 
\begin{align}     \label{e:14:mScaleRela}
      1-b\beta + \frac{(1+\beta)(b-1)}{2(1+b)} - \frac{b\gamma_S}{2}
      \leq be \leq  1-b\beta + \frac{(1+\beta)(b-1)}{2(1+b)} + \frac{b\gamma_S}{2}.
\end{align}

We define $\alpha \in (0,1)$ such that
\begin{align}     \label{e:alphaDef}
      2\alpha + \beta = 1 + 2b\gamma_I - 2b^2 \gamma.
\end{align}

We also fix absolute integers $N_*, N, Q \in \N$ such that
\begin{align}     \label{e:NQ}
      N_* \gamma \geq 3,      \quad 
      Q \geq N \geq N_*^4.
\end{align}

\begin{remark}\label{r:auxiliaryScalesRela}
\eqref{e:12:mScaleRela} leads to the constraint $\beta < \frac{1}{3}$. The constraints \eqref{e:4:mScaleRela} and \eqref{e:14:mScaleRela} are more at the technical level. The essential constraint for parameters is
\begin{align}
      \mr \mu_{q+1} > \delta_{q+1}^{\sfrac12} \lambda_{q+1}
            > \kappa_{q+1}\lambda_{q+1}^2 > \mu_{q+1} > \mr\mu_q > \delta_{q}^{\sfrac12} \lambda_{q},
\end{align}
which also leads to the constraint $\beta < \frac{1}{3}$. From direct computations, one can verify above relation from our constraints on parameters.
\end{remark}

\begin{remark}    \label{r:auxiliaryParaRela}
We have the following relations on parameters. From \eqref{e:smallGammaR}, \eqref{e:mrParameter}, \eqref{e:12:mScaleRela} and \eqref{e:14:mScaleRela}, we have 
\begin{align*}
      \varepsilon_q < 1.
\end{align*}
From \eqref{e:0:mScaleRela}, \eqref{e:smallGammaR} and $\beta < \frac{1}{3}$, we have 
\begin{align*}
      \mr\lambda_{q+1} > \lambda_{q+1} > \mr \lambda_q > \lambda_q.
\end{align*}
From \eqref{e:smallGammaR}, \eqref{e:mrParameter} and \eqref{e:14:mScaleRela}, we have
\begin{align*}
      \mr\mu_{q+1} > \mu_{q+1} > \mr \mu_q > \mu_q.
\end{align*}
\end{remark}

\subsection{Cutoffs} \label{ss:cutoffs}
First, in the following lemma, we define some cutoff functions.
\begin{lemma}     \label{l:cutoff}
There exist smooth periodic functions $\eta_1, \eta_2: \T \rightarrow \R $ and non-periodic $\tilde\eta: \R \rightarrow \R$ satisfying the following
\begin{enumerate}[leftmargin=*,label=\textsc{(B.\arabic*)},align=left]
\item
\begin{align}     \label{e:2:cutoff}
      \int_\T \eta_1^2(\tau) d\tau \int_\T \sin^2(2\pi s) ds 
      = \int_\T \eta_2^2(\tau) d\tau \int_\T \sin^2(2\pi s) ds = 1.
\end{align}
\item $\supp \eta_1 \subset \bigcup_{\iota \in \Z} (\iota+\frac{1}{8}, \iota+\frac{3}{8})$ and $\supp \eta_2 \subset \bigcup_{\iota \in \Z} (\iota+\frac{5}{8},\iota+\frac{7}{8})$ and for any $p \geq 0$
\begin{align}     \label{e:4:cutoff}
      | \partial_\tau^p \eta_1| + | \partial_\tau^p \eta_2 | \lesssim p^{3p} 
\end{align}
\item $\tilde \eta(\tau) = 1$ for $\tau \in ( -\frac{3}{8}, \frac{3}{8} )$ and $\supp \tilde \eta \subset [-\frac{5}{8}, \frac{5}{8}]$ and for any $p \geq 0$
\begin{align}     \label{e:6:cutoff}
      | \partial_\tau^p \tilde\eta | \lesssim p^{3p} 
\end{align}
\item Define
\begin{align}     \label{e:8:cutoff}
      \tilde \eta_\iota(\tau) := \tilde \eta (\tau-\iota),
\end{align}
then $\{\tilde \eta_\iota\}_{\iota \in \Z}$ forms a partition of unity on $\R$.
\end{enumerate}
The constant in all estimates is absolute.
\end{lemma}

\begin{proof}[Proof of Lemma \ref{l:cutoff}]
For the standard (smooth, compactly supported) mollifier 
\[
\phi (\tau)= C \exp \left(\frac{1}{\tau^2-1} \right) \one_{(-1,1)}
\]
we have the following estimate
\begin{equation}\label{e:std_mol_est}
| \partial_\tau^p \phi (\tau)| \le C\left(\frac{8}{e} p\right)^{2p} \le c_0 p^{3p},
\end{equation}
where the former is a direct computation. 
Shifting, scaling and periodising $\phi$, we obtain $\eta_1$, $\eta_2$ with their estimates. 

The partition of unity $\tilde \eta(\tau)$ can be produced as follows. Take function $r(\tau)$ whose graph is a isosceles trapezoid, equal to $1$ on  $( -\sfrac{3}{8}-\varepsilon_0, \sfrac{3}{8}+\varepsilon_0)$ and vanishing outside $[-\sfrac{5}{8}+\varepsilon_0, \sfrac{5}{8}-\varepsilon_0]$. Copies of this function, appropriately shifted, yield (non-smooth) partition of unity. Mollification of $r(\tau)$, by convolution with $\varepsilon^{-1}_0 \phi (\frac{\cdot}{\varepsilon_0})$, produces function $\tilde \eta(\tau)$ with bounds on derivatives inherited from the bound \eqref{e:std_mol_est}. Observe that, thanks to symmetry, shifted copies of $\tilde \eta(\tau)$ remain a partition of unity (now smooth).
\end{proof}

\subsection{Construction of the vector field} \label{ss:vectorfield}

In this section, we construct two divergence-free vector fields $\{u_q\}_{q}$ and $\{v_q\}_{q}$ by induction. The $\{\psi_q\}_{q}$ are the stream functions of $\{u_q\}_{q}$. Then we define the limit vector field $u$ and its stream function $\psi$ via
\begin{align}
      u := \lim_{q \rightarrow \infty} u_q, \quad 
      \psi := \lim_{q \rightarrow \infty} \psi_q.
\end{align}

Define
\begin{align}
      u_0 = v_0 = 0.
\end{align}
Assume we have $u_q$ and $v_q$. Below we construct $u_{q+1}$ and $v_{q+1}$.
At the induction step, define $\{t_{q,\iota}\}_{\iota \in \Z}$ via $t_{q,\iota} := \frac{\iota}{\mu_{q+1}}$ and 
\begin{align}     \label{e:6:constructVF}
      \tilde \eta_{q,\iota}( t ) := \tilde \eta_\iota( \mu_{q+1} t ).
\end{align}
Here, $\tilde \eta_\iota$ is given in Lemma \ref{l:cutoff}. Then $\tilde \eta_{q,\iota} = 1$ on $( t_{q,\iota} - \frac{3}{8\mu_{q+1}}, t_{q,\iota} + \frac{3}{8\mu_{q+1}} )$ and $\supp \tilde \eta_{q,\iota} \subset ( t_{q,\iota} - \frac{5}{8\mu_{q+1}}, t_{q,\iota} + \frac{5}{8\mu_{q+1}} )$.

Next, we define two local flow maps $\bar \Phi_{q,\iota}, \Phi_{q,\iota}: \T^2 \times \R \rightarrow \T^2$. Here, $\bar \Phi_{q,\iota}$ solves the transport equation
\begin{equation}     \label{e:8:constructVF}
\begin{split}
      \partial_t \bar \Phi_{q,\iota} + ( u_q \cdot \nabla ) \bar \Phi_{q,\iota} =&\, 0, \\ 
      \bar \Phi_{q,\iota} (x, t_{q,\iota}) =&\, x. 
\end{split}
\end{equation}
And $\Phi_{q,\iota}$ is given by the following explicit formula
\begin{align}     \label{e:10:constructVF}
      \Phi_{q,\iota} = x + \sum_{m=1}^{N_*} \frac{1}{m!} ( t_{q,\iota} - t )^m D_{q,t}^{m-1} u_q.
\end{align}
We define the global flow maps $\bar\Phi_q, \Phi_q: \T^2 \times \R \rightarrow \T^2$ via
\begin{align}     \label{e:16:constructVF}
      \bar \Phi_q = \sum_{\iota} \bar \Phi_{q,\iota} \tilde \eta_{q,\iota},       \quad 
      \Phi_q = \sum_{\iota} \Phi_{q,\iota} \tilde \eta_{q,\iota}.
\end{align}

We define \textit{antisymmetric matrix}-valued function $H: \T^2 \times \T \rightarrow \R^{2 \times 2}$ via
\begin{align}
      H_{12}(\xi, \tau) :=&\, \eta_1(\tau) \Pi_1(\xi) + \eta_2(\tau) \Pi_2(\xi) 
            = - H_{21}(\xi, \tau),  \label{e:H} \\ 
            \Pi_1(\xi) =&\, \sin(2\pi \xi_1), \quad \Pi_2(\xi) = \sin(2\pi \xi_2),  \label{e:hom:cos}
\end{align}
with $\eta_1$ and $\eta_2$ given by Lemma \ref{l:cutoff}. Then we define the stream functions
\begin{align}
      \bar \psi_{q+1}(x,t) :=&\, \frac{\delta_{q+1}^{\sfrac12}}{\lambda_{q+1}} 
            H_{12} \big( \lambda_{q+1} \bar \Phi_q(x,t), \mu_{q+1} t \big),     \label{e:barpsiQ+1} \\ 
      \psi_{q+1}(x,t) :=&\, \frac{\delta_{q+1}^{\sfrac12}}{\lambda_{q+1}} 
            \det \nabla \Phi_q H_{12} \big( \lambda_{q+1} \Phi_q(x,t), \mu_{q+1} t \big),     \label{e:psiQ+1} 
\end{align}
We construct $\{u_q\}_q$ and $\{v_q\}_q$ inductively via 
\begin{align}
      u_{q+1} =&\, u_q + \bar w_{q+1},    \quad
            \bar w_{q+1} = \nabla^\perp \bar \psi_{q+1}     \label{e:24:constructVF} \\ 
      v_{q+1} =&\, u_q + w_{q+1},  \quad       
            w_{q+1} = \nabla^\perp \psi_{q+1}.          \label{e:32:constructVF} 
\end{align}
Here, $v_{q+1}$ is not iteratively defined by $v_q$. Instead, $v_{q+1}$ is defined to be a perturbation of $u_{q+1}$.

Define
\begin{align}
      \bar \Psi_{q+1} :=&\, \sum_{j=0}^{q+1} \bar \psi_j,     \label{e:barCumuPsiQ+1} \\ 
      \Psi_{q+1} :=&\, \bar \Psi_q + \psi_{q+1},     \label{e:CumuPsiQ+1}
\end{align}

then we have
\begin{align}
      u_{q} &\,= \nabla^\perp \bar \Psi_q,      \quad 
            v_{q} = \nabla^\perp \Psi_q,       \label{e:26:constructVF} \\ 
      \bar w_{q+1} =&\, \frac{\delta_{q+1}^{\sfrac12}}{\lambda_{q+1}} \divr \Big( H ( \lambda_{q+1} \bar \Phi_q, \mu_{q+1} t ) \Big),     \label{e:28:constructVF} \\ 
      w_{q+1} =&\, \frac{\delta_{q+1}^{\sfrac12}}{\lambda_{q+1}} 
      \divr \Big( \det \nabla \Phi_q H ( \lambda_{q+1} \Phi_q, \mu_{q+1} t ) \Big).       \label{e:36:constructVF} 
\end{align}

From the support information of the cutoff functions $\eta_1, \eta_2$ in Lemma \ref{l:cutoff} and the construction of the velocity fields, we have the following property
\begin{align}     \label{e:u_zero_initialT}
      \nabla^p u_q( x, 0 ) = \nabla^p v_q( x, 0 ) = 0, \quad \forall \, p,q \in \N, x \in \T^2.
\end{align}

We also define a support condition.
\begin{definition}[Support condition]      \label{d:dis_supp_TPR}
A function $\eta: \T \rightarrow \R$ has support disjoint to $\eta_1, \eta_2$ if 
\begin{align}
      \supp \eta \cap \supp \eta_1 = \supp \eta \cap \supp \eta_2 = \varnothing.
\end{align}
\end{definition}

\begin{definition}[Shear condition]      \label{d:shear_structure}
A pair of functions $(\chi, \eta)$ with $\chi: \T^2 \rightarrow \R$, $\eta: \T \rightarrow \R$ satisfy the shear condition if either 
\begin{align}
      \supp \eta \subset \supp \eta_1,    \quad 
      \partial_{\xi_2} \chi = 0,
      \qquad \text{or} \qquad
      \supp \eta \subset \supp \eta_2,    \quad 
      \partial_{\xi_1} \chi = 0.
\end{align}
A pair of functions $(\chi, \eta)$ satisfy the generalized shear condition if either it satisfies shear condition or
\begin{align}
      \nabla_\xi \chi = 0.
\end{align}
\end{definition}

\subsection{The polynomial-derivative space}\label{ss:poly_space}

We introduce now quantities useful for keeping track of forthcoming estimates of stream function, while high-order homogenisation is applied.

Recall the notations in Section \ref{sss:ntt:multi_ind} and Section \ref{sss:ntt:deri}. 

\begin{definition}[Renormalized derivatives]    \label{d:renormalizedDiff}
For any differentiation index $\bomega \in \mathcal{I}$, let $D_q^{\bomega}$ be the \textit{generalized derivative associated to} the vector field $u_q$ with index $\bomega$. Define a \textit{renormalized differential operator associated to} $u_q$ via
\begin{align}
      \mr D_q^{\bomega} = \mr \lambda_q^{-|\bomega|_x} \mr \mu_q^{-|\bomega|_t} D_q^{\bomega}.
\end{align}
Similarly, we define $O_q^{\bomega}$ to be the \textit{generalized derivative associated to} the vector field $v_q$ with $\bomega$. We also define the \textit{renormalized differential operator associated to} $v_q$ via
\begin{align}
      \mr O_q^{\bomega} = \mr \lambda_q^{-|\bomega|_x} \mr \mu_q^{-|\bomega|_t} O_q^{\bomega}.
\end{align}
The same definition also extends to the space derivatives $\mr \nabla_q$, $\mr \partial_q^{\bomega}$ and $\mr\divr_q$, i.e.
\begin{align}
      \mr \nabla_q := \mr \lambda_q^{-1} \nabla,      \quad 
      \mr \partial_q^{\bomega} := \mr \lambda_q^{-|\bomega|} \partial^{\bomega},    \quad 
      \mr \divr_q := \mr \lambda_q^{-1} \divr
\end{align} 
for $\bomega$ only containing space derivatives.
\end{definition}

Note that, for fixed $q$, the renormalized \textit{space} derivatives are the same for different vector fields $u_q$ and $v_q$.

From the construction, $u_q$ and $v_q$ are smooth for any $q \in \N$. Now we define the following algebra generated by derivatives of their stream functions.

\begin{definition}[Polynomial-derivative spaces $\mathcal{P}_q$ and $\mathcal{O}_q$]
For any $n_0 \in \N$, recalling \eqref{e:barCumuPsiQ+1} and \eqref{e:CumuPsiQ+1}, we define
\begin{equation}  \label{e:0:poly_Uq}
\begin{split}
      \mathcal{P}_q(n_0) := \Bigg\{ &\,  a_0 + \sum_{k=1}^{k_0} a_k \prod_{i=1}^{I_k} \mr D_q^{\balpha_{k,i}} \Bigg( \frac{\mr \lambda_{q}^2} {\delta_{q}^{\sfrac12}\lambda_{q}^{1+2\gamma_R}}  \bar \Psi_q \Bigg) \, \Bigg| \,
            k_0 \in \N. \, \{I_k\}_{1 \leq k \leq k_0} \subset \N^+.  \\
            &\, \{a_k\}_{0 \leq k \leq k_0} \subset \R.\, \{ \balpha_{k,i} \}_{1 \leq k \leq k_0, 1 \leq i \leq I_k} \subset \mathcal{I}(2,0) \cup \mathcal{I}(1,N_*), 
                  \sup_k \sum_i | \balpha_{k,i} | \leq n_0
      \Bigg\} 
\end{split}
\end{equation}
and
\begin{equation}  \label{e:2:poly_Oq}
\begin{split}
      \mathcal{O}_q(n_0) := \Bigg\{ &\,  a_0 + \sum_{k=1}^{k_0} a_k \prod_{i=1}^{I_k} \mr O_q^{\balpha_{k,i}} \Bigg( \frac{\mr \lambda_{q}^2} {\delta_{q}^{\sfrac12}\lambda_{q}^{1+2\gamma_R}} \Psi_q \Bigg) \, \Bigg| \, \\
            &\, k_0, I_k, a_k, \balpha_{k,i} \text{ as above}, \, 
                  \sup_k \sum_i | \balpha_{k,i} | \leq n_0
      \Bigg\}.
\end{split}
\end{equation}
\end{definition}

\begin{remark}    \label{r:renormal_factor}
The prefactor $\frac{ \mr \lambda_{q}^2 } { \delta_{q}^{\sfrac12}\lambda_{q}^{1+2\gamma_R} }$ is chosen such that
\begin{align}
      \frac{1}{\varepsilon_q \mu_q} \nabla u_q 
            = \frac{1} {\delta_{q}^{\sfrac12}\lambda_{q}^{1+2\gamma_R}} \nabla u_q 
            = \frac{\mr\lambda_{q}^2} { \delta_{q}^{\sfrac12} \lambda_{q}^{1+2\gamma_R} } 
                  \mr \nabla_q \mr \nabla^\perp_q \bar \Psi_q   \quad 
      \rightarrow \quad
      \vertiii {\frac{1}{\varepsilon_q \mu_q} \nabla u_q}_q \sim 1.
\end{align}
\end{remark}

We also define a norm-like quantity for these polynomial-derivative spaces. Note that these are nowhere rigorously close to the notion of norm in Banach spaces.
\begin{definition}[Norm-like quantity of polynomial-derivative space]
For any $f \in \mathcal{P}_q(n_0)$, we define a number $\vertiii{f}_q$ via
\begin{equation}  \label{e:0:polyNorm}
\begin{split}
      \vertiii{f}_q := 
      \inf \Bigg\{ \sum_{k=0}^{k_0} |a_k| \, \Bigg| \, f \text{ can be written } 
            f = a_0 + \sum_{k=1}^{k_0} a_k \prod_{i=1}^{I_k} \mr D_q^{\balpha_{k,i}} \Bigg( \frac{\mr \lambda_{q}^2} {\delta_{q}^{\sfrac12}\lambda_{q}^{1+2\gamma_R}}  \bar \Psi_q \Bigg) &\, \\
            \text{with } k_0,I_k,a_k, \balpha_{k,i}, n_0 \text{ satisfying the same conditions in } \eqref{e:0:poly_Uq} &\,
      \Bigg\}
\end{split}
\end{equation}
Similarly, we also define the number $\vertiii{f}_q$ for each element $f \in \mathcal{O}_q(n_0)$.
\end{definition}

The infimum in \eqref{e:0:polyNorm} is actually a minimum. This is explained in the following remark.

\begin{remark}    \label{r:PolyDeriP2}
An element in $f \in \mathcal{P}_q(n_0)$ has multiple representations of form
\begin{align}     \label{e:2:PolyDeriP2}
      f = a_0 + \sum_{k=1}^{k_0} a_k \prod_{i=1}^{I_k} \mr D_q^{\balpha_{k,i}} 
            \Bigg( \frac{\mr \lambda_{q}^2} {\delta_{q}^{\sfrac12}\lambda_{q}^{1+2\gamma_R}}  \bar \Psi_q \Bigg),
\end{align}
since one can switch the order of differentiation in indices $\balpha_{k,i}$. However, $\vertiii{f}_q$ is uniquely defined for fixed $q$, because the infimum is taken among all representations. For fixed $n_0$, such representations live in a space with finite dimension, therefore the infimum in \eqref{e:0:polyNorm} is actually a minimum.

\begin{definition}[Lower bound of differentiation order]    \label{d:polyOrderLBd}
We say the differentiation order of $f \in \mathcal{P}_q(n_0)$ is bounded from below by $n \in \N$ if there exists a representation \eqref{e:2:PolyDeriP2} satisfying 
\begin{align}     \label{e:2:polyOrderLBd}
      \sup_{1 \leq k \leq k_0} \sum_i | \balpha_{k,i} | \leq n_0,       \quad 
      \sum_{k=0}^{k_0} |a_k| = \vertiii{f}_q 
\end{align}
such that $a_0 = 0$ and
\begin{align}     \label{e:4:polyOrderLBd}
      \inf_{1 \leq k \leq k_0} \sum_i | \balpha_{k,i} | \geq n.
\end{align}
We also extend the same definition to $\mathcal{O}_q(n_0)$.
\end{definition}

Given above definitions, we have the following lemma on a truncation projection regarding the differentiation order.

\begin{lemma}[Truncation on differentiation order]    \label{l:polyOrderTrun}
For any $n_0 \in \N$ and any element $f \in \mathcal{P}_q(n_0)$ and any integer $n \in \N$, there exists a projection of $f$ denoted by
\begin{align}
      \bP_{n} f \in \mathcal{P}_q(n)
\end{align}
such that the differentiation order of $f - \bP_{n} f \in \mathcal{P}_q(n_0)$ is bounded from below by $n$. Moreover, we have
\begin{align}
      \vertiii{ \bP_{n} f }_q + \vertiii{ f - \bP_{n} f }_q
            = \vertiii{ f }_q.
\end{align}
The same result also holds for $\mathcal{O}_q(n_0)$.
\end{lemma}

\begin{proof}
We choose a representation of $f$ like \eqref{e:2:PolyDeriP2} satisfying \eqref{e:2:polyOrderLBd}, then we define
\begin{align*}
      \bP_{n} f :=&\, a_0 + \sum_{k \in \blacksquare_{n}} 
            a_k \prod_{i=1}^{I_k} \mr D_q^{\balpha_{k,i}} 
            \Bigg( \frac{\mr \lambda_{q}^2} {\delta_{q}^{\sfrac12}\lambda_{q}^{1+2\gamma_R}}  \bar \Psi_q \Bigg),       \\
      \blacksquare_{n} :=&\, \bigg\{ 1 \leq k \leq k_0 \,\bigg|\, \sum_i | \balpha_{k,i} | \leq n \bigg\}.
\end{align*}
\end{proof}

\begin{remark}
Here, we do not claim the truncation projection in Lemma \ref{l:polyOrderTrun} is unique.
\end{remark}

The elements in $ \mathcal{O}_q $ are approximations of elements in $ \mathcal{P}_q $. Now we introduce the notion of congruence to formalize this correspondence.

\begin{definition}[Congruence]      \label{d:polyCongr}
Fix $q \geq 1$. For some $g \in \mathcal{P}_q(n_0)$ and $h \in \mathcal{O}_q(n_1)$, we say $g$ and $h$ are \textit{congruent} if there exist $k_0, \{I_k\}_k, \{a_k\}_k, \{\balpha_{k,i}\}_{k,i}$ as in \eqref{e:0:poly_Uq}, satisfying that 
\begin{align}
      g =&\, a_0 + \sum_{k=1}^{k_0} a_k \prod_{i=1}^{I_k} \mr D_q^{\balpha_{k,i}} \Bigg( \frac{\mr \lambda_{q}^2} {\delta_{q}^{\sfrac12}\lambda_{q}^{1+2\gamma_R}}  \bar \Psi_q \Bigg),  \quad  
      h = a_0 + \sum_{k=1}^{k_0} a_k \prod_{i=1}^{I_k} \mr O_q^{\balpha_{k,i}} \Bigg( \frac{\mr \lambda_{q}^2} {\delta_{q}^{\sfrac12}\lambda_{q}^{1+2\gamma_R}} \Psi_q \Bigg). 
\end{align}
\end{definition}

\begin{remark}    \label{r:PolyDeriP1}
For $a \in \mathcal{P}_q(0) = \R$, 
\begin{align}
      \vertiii{a}_q := |a|.
\end{align}
\end{remark}

The same applies to $\mathcal{O}_q(n_0)$.
\end{remark}

\begin{remark}    \label{r:PolyDeriP3}
From the support of the cutoff function in Lemma \ref{l:cutoff} and the definition of $\bar \Psi_q$, any element $h$ in $\mathcal{P}_q(n_0)$ satisfies
\begin{align*}
      h(x,t) \text{ is constant } \forall \, (x,t) \in \T^2 \times \bigg[ 0 , \frac{1}{8\mu_q} \bigg],
      \quad \text{ and } \quad
      |h(x,t)| \leq \vertiii{h}_q, \, \forall \, (x,t) \in \T^2 \times \bigg[ 0 , \frac{1}{8\mu_q} \bigg].
\end{align*}
The same applies to $\mathcal{O}_q(n_0)$.
\end{remark}

\begin{remark}    \label{r:PolyDeriP4}
For any $n_0 \in \N$, $\mathcal{P}_q(n_0)$ and $\mathcal{O}_q(n_0)$ are linear spaces.
\end{remark}

\begin{remark}    \label{r:opPolynomial}
If $f \in \mathcal{P}_q(n_1)$ and $g \in \mathcal{P}_q(n_2)$, we have
\begin{align}
      f g \in&\, \mathcal{P}_q(n_1+n_2), \\ 
      f + g \in&\, \mathcal{P}_q( \max \{n_1,n_2\} ), \\ 
      \vertiii{f+g}_{q} \leq&\, \vertiii{f}_{q} + \vertiii{g}_{q}, \\ 
      \vertiii{fg}_{q} \leq&\, \vertiii{f}_{q} \vertiii{g}_{q}, \\ 
      \nabla f, D_{q,t} f \in&\, \mathcal{P}_q ( n_1+1).
\end{align}
The same applies to $\mathcal{O}_q$.
\end{remark}

\newpage

\section{Analysis of stream functions}     \label{ss:r_estimates}

In this section, we prove some estimates related to the stream functions $\{\bar\Psi_q\}_q$ and $\{\Psi_q\}_q$. 

\begin{lemma} \label{l:uEstimate} Recall $ N_*, Q$ and $N (N_*)$ fixed in \eqref{e:NQ}.
Recall $H$ defined in \eqref{e:H}, the flow map $\Phi_q$ defined in \eqref{e:10:constructVF} (which gives number $N_*$), the flow map $\bar \Phi_q$ defined in \eqref{e:10:constructVF} \eqref{e:16:constructVF} and the vector field $u_q$ defined \eqref{e:24:constructVF}. 

There exist smooth functions 
\[
\begin{aligned}
&B_{q,m}, E_{q,m}, S_{q,m}, \Omega_{q,m}: \T^2 \times [0,1] \rightarrow \R^{2 \times 2},\\
&\omega_{q,m}: \T^2 \times [0,1] \rightarrow \R^2,\\
&\varpi_{q,m}: \T^2 \times [0,1] \rightarrow \R, \\
&\phi_{m}, \varphi_{m}, z_{m}, \vartheta_{m}, \sigma_m, \zeta_m: \T \rightarrow \R
\end{aligned}
\]
indexed by $q \in \N$ and $1 \leq m \leq N (N_*)$, such that: \\
(I) We have the following identities for any $(x,t) \in \T^2 \times [0,1]$, $\xi \in \T^2$ and $\tau \in \T$ (below, if the argument on the left hand side is $(x,t)$, we omit it for brevity). 
\begin{align}
      \nabla \Phi_q =&\, \Id + \varepsilon_q \sum_{m=1}^N B_{q,m} (x,t) {\phi_{m}} (\mu_{q+1} t),     \label{e:6:uEstimate}  \\ 
      {\det \nabla \Phi_q} =&\, {1+ \varepsilon_q}{\sum_{m=1}^N \varpi_{q,m}(x,t) \zeta_{m} (\mu_{q+1} t),}    \label{e:7:uEstimate} \\ 
         \nabla \Phi_q \nabla \Phi_q^T - (\det \nabla \Phi_q)^2 \Id
            =&\, \varepsilon_q \sum_{m=1}^N S_{q,m} (x,t) \vartheta_{m} (\mu_{q+1} t),       \label{e:10:uEstimate}  \\
      \adj \nabla \Phi_q^{T} \det \nabla \Phi_q =&\, \Id + \varepsilon_q \sum_{m=1}^N E^T_{q,m} (x,t) \varphi_{m} (\mu_{q+1} t),   \label{e:8:uEstimate}\\
      {\big( (\det \nabla \Phi_q) H(\xi, \tau)}\big)_{ij}
            =&\, H_{ij}(\xi, \tau) + \varepsilon_q \sum_{m=1}^N \sigma_m(\tau) \Omega_{q,m,ij} (x,t) H_{12}(\xi, \tau),    \label{e:12:uEstimate} \\
      \sum_{\iota \in \Z} \tilde \eta_{q,\iota} D_{q,t} \Phi_{q,\iota} = &\, {\omega_{q,N_*}}(x,t) {\phi_{N_*}} (\mu_{q+1} t),       \label{e:2:uEstimate} \\ 
      \sum_{\iota \in \Z} \partial_t \tilde \eta_{q,\iota} \Phi_{q,\iota} =&\, \sum_{m=1}^{{N_*}} {\omega_{q,m}} (x,t) z_{m} (\mu_{q+1} t),      \label{e:4:uEstimate}
\end{align}
with
\begin{align}
      \omega_{q,N} \in (\mathcal{P}_q(N))^2, \quad \varpi_{q,m} \in \mathcal{P}_q(N),  
            \label{e:algbVarpiOmega} \\ 
      B_{q,m}, E_{q,m}, S_{q,m}, \Omega_{q,m} \in (\mathcal{P}_q(N))^{2 \times 2}.  
            \label{e:algbBESOmega}
\end{align}
Furthermore, for any $1 \leq m \leq N$, $z_m$ has support disjoint to $\eta_1$ and $\eta_2$. \\
(II) We have the following estimates.
For $p \lesssim 8 Q^3$ the functions defined in (I) satisfy
\begin{equation}
      \| \partial_\tau^p \phi_m \|_\infty + \| \partial_\tau^p \varphi_m \|_\infty + \| \partial_\tau^p z_m \|_\infty + \| \partial_\tau^p \vartheta_m \|_\infty + \| \partial_\tau^p \sigma_m \|_\infty + \| \partial_\tau^p \zeta_m \|_\infty 
            \lesssim\, 1.       \label{e:phiZetaThetaEst}  
\end{equation}
and 
\begin{align}
      \sum_{m=1}^N \vertiii{ B_{q,m} }_q + \vertiii{ E_{q,m} }_q + \vertiii{ S_{q,m} }_q
            + \vertiii{ \Omega_{q,m} }_q + \vertiii{ \omega_{q,m} }_q + \vertiii{ \varpi_{q,m} }_q 
            \lesssim&1. \label{e:mtrPolyEst} \\
                  {\vertiii{ \omega_{q,{N_*}} }_q} \lesssim&\, \lambda_{q+1}^{-3},       \label{e:trsptRsdPolyEst}
\end{align}
The constant in $\lesssim$ depends in fact only on $Q$ (since by choice \eqref{e:NQ}  $N\le Q$ and  $N_*\le Q$, and $Q$ is an absolute constant).
\end{lemma}
Lemma \ref{l:uEstimate} is proven in Section \ref{sec:l:uEstimate_pf}.

\begin{lemma}     \label{l:streamEst}
For the cumulative stream function $\bar \Psi_{q+1}$ defined in \eqref{e:barCumuPsiQ+1}, and for the stream function $\bar \psi_{q+1}$ defined in \eqref{e:barpsiQ+1}, we have
\begin{equation}
      \| \nabla^n \bar \Psi_{q+1} \|_{\infty} 
      \le\,C(n)  \frac{\delta_{{q+1}}^{\sfrac12}}{\lambda_{{q+1}}} \lambda_{{q+1}}^{n}       \qquad \text{for any } n \ge 2 \label{e:streamEs_rig1}
          \end{equation}
          and 
\begin{equation}\label{e:streamEs_rig1corr}
\norm{\nabla^n \bar \psi_{{q+1}}}_\infty \le C(n) \frac{\delta_{{q+1}}^{\sfrac12}}{\lambda_{{q+1}}} \lambda_{{q+1}}^{n} \qquad \text{for any } n \ge 0.
\end{equation}
For any fixed number $Y$we have also 
          \begin{equation}
      \| D^{\balpha}_{q+1} \bar \Psi_{q+1} \|_{\infty} 
      \lesssim\, 1 \vee \frac{\delta_{{q+1}}^{\sfrac12}}{\lambda_{{q+1}}} (\delta_{{q+1}}^{\sfrac12} \lambda_{{q+1}})^{|\balpha|_t} \lambda_{{q+1}}^{|\balpha|_x}      \quad \text{for any } |\balpha| \leq Y \label{e:streamEs_rig2} 
          \end{equation}
          and
\begin{equation}\label{e:streamEs_rig2_psi}
      \| D^{\balpha}_{{q+1}} \bar \psi_{q+1} \|_{\infty} 
      \lesssim \frac{\delta_{{q+1}}^{\sfrac12}}{\lambda_{{q+1}}} (\delta_{{q+1}}^{\sfrac12} \lambda_{{q+1}})^{|\balpha|_t} \lambda_{{q+1}}^{|\balpha|_x} \qquad \text{for any }  |\balpha| \leq Y. 
\end{equation}
The constant $C(n)= \frac{(n-1)!}{n^2} C^{n-1} \tilde C$, for certain absolute constants, and $\lesssim$ depends only on the upper bound $Y$. The constants can be taken $C=50 \pi$ and $\tilde C =2$. Importantly, none of the constants depend on $q$.
\end{lemma}
Lemma \ref{l:streamEst} is proven in Section  \ref{sec:streamEst_rig1_pf}.

\begin{lemma}     \label{l:streamEst_rig_nobar}
Take $Y$ fixed in Lemma  \ref{l:streamEst}. For the cumulative stream function $\Psi_{q+1}$  defined in \eqref{e:CumuPsiQ+1}, and the stream functions  $\bar \psi_{q+1}$, $\psi_{q+1}$ defined in  \eqref{e:barpsiQ+1}, \eqref{e:psiQ+1},   we have for any $|\balpha| \leq Y-N_*-1$
          \begin{equation}
      \| D^{\balpha}_{q+1} \Psi_{q+1} \|_{\infty} +     \| O^{\balpha}_{q+1} \Psi_{q+1} \|_{\infty} 
      \lesssim\, 1 \vee \frac{\delta_{q+1}^{\sfrac12}}{\lambda_{q+1}} (\delta_{q+1}^{\sfrac12} \lambda_{q+1})^{|\balpha|_t} \lambda_{q+1}^{|\balpha|_x} \label{e:streamEs_rig2_nobar} 
          \end{equation}
 \begin{equation}\label{e:streamEs_rig2_nobar_psi}
      \| D^{\balpha}_{q+1} \psi_{q+1} \|_{\infty} +   \| O^{\balpha}_{q+1} \psi_{q+1} \|_{\infty} 
      \lesssim \frac{\delta_{q+1}^{\sfrac12}}{\lambda_{q+1}} (\delta_{q+1}^{\sfrac12} \lambda_{q+1})^{|\balpha|_t} \lambda_{q+1}^{|\balpha|_x}.
\end{equation}
For their difference $\bar \psi_{q+1} -  \psi_{q+1}= \bar \Psi_{q+1} -  \Psi_{q+1}$ it holds, recalling $\gamma$ of section \ref{ss:parameters}, for any $|\balpha| \leq Y-N_*-1$
\begin{equation}\label{e:streamEs_rig2_nobar_psi_diff}
      \| D^{\balpha}_{q+1} (\bar \psi_{q+1} -  \psi_{q+1})\|_{\infty} +   \| O^{\balpha}_{q+1} (\bar \psi_{q+1} -  \psi_{q+1}) \|_{\infty} 
      \lesssim \lambda^{-\gamma N_*+1}_{q+1} \frac{\delta_{q+1}^{\sfrac12}}{\lambda_{q+1}} (\delta_{q+1}^{\sfrac12} \lambda_{q+1})^{|\balpha|_t} \lambda_{q+1}^{|\balpha|_x}.
\end{equation}
For the difference of derivative operators it holds for any $ |\balpha| \leq Y-N_*-1$
\begin{equation}\label{e:streamEs_psi_DOdiff}
      \| (D^{\balpha}_{q+1}- O^{\balpha}_{q+1}) \bar \psi_{q+1} \|_{\infty} +    \| (D^{\balpha}_{q+1}- O^{\balpha}_{q+1}) \psi_{q+1} \|_{\infty}
      \lesssim \lambda^{-\gamma N_*+1}_{q+1} \delta_{q+1} (\delta_{q+1}^{\sfrac12} \lambda_{q+1})^{|\balpha|_t} \lambda_{q+1}^{|\balpha|_x}.
\end{equation}
as well as
\begin{equation}\label{e:streamEs_Psi_DOdiff}
      \| (D^{\balpha}_{q+1}- O^{\balpha}_{q+1}) \bar \Psi_{q+1} \|_{\infty} +    \| (D^{\balpha}_{q+1}- O^{\balpha}_{q+1}) \Psi_{q+1} \|_{\infty}
      \lesssim \delta^{\sfrac12} _{q+1} \lambda^{-\gamma N_*+2}_{q+1}  (\delta_{q+1}^{\sfrac12} \lambda_{q+1})^{|\balpha|_t} \lambda_{q+1}^{|\balpha|_x}.
\end{equation}
The constant in $\lesssim$ depends only on the upper bound $Y$.
\end{lemma}

\begin{corollary}\label{l:strDiffEst2}
In fact, denoting by 
$\Theta^{\balpha}_{q+1}$ either   $D^{\balpha}_{q+1}$ or  $O^{\balpha}_{q+1}$, any sequence consisting of symbols $\Theta^{\balpha_i}_{q+1}$ and  $(D^{\balpha_j}_{q+1} - O^{\balpha_j}_{q+1})$ acting on either $\bar \psi_{q+1}$ or on $\psi_{q+1}$ can be estimated as \eqref{e:streamEs_psi_DOdiff}, and acting on either $\bar \psi_{q+1}$ or on $\psi_{q+1}$ can be estimated as \eqref{e:streamEs_Psi_DOdiff},  provided there is at least one term $(D^{\balpha_k}_{q+1} - O^{\balpha_k}_{q+1})$. More precisely, for $f$ being either $\bar \psi_{q+1}$ or $\psi_{q+1}$ and for $F$ being either $\bar \Psi_{q+1}$ or $\Psi_{q+1}$ we have, for any $|\balpha| \leq Y-N_*-1$
\begin{equation}\label{e:streamEs_psi_DOdiff_gen}
       \|  (D^{\balpha_1}_{q+1} - O^{\balpha_1}_{q+1}) \Theta^{\balpha_2}_{q+1}  (D^{\balpha_3}_{q+1} - O^{\balpha_3}_{q+1}) \dots \Theta^{\balpha_k}_{q+1} f \|_{\infty} 
      \lesssim \lambda^{-\gamma N_*+1}_{q+1} \delta_{q+1} (\delta_{q+1}^{\sfrac12} \lambda_{q+1})^{|\balpha|_t} \lambda_{q+1}^{|\balpha|_x} \text{ etc},
\end{equation}
\begin{equation}\label{e:streamEs_Psi_DOdiff_gen}
       \|  (D^{\balpha_1}_{q+1} - O^{\balpha_1}_{q+1}) \Theta^{\balpha_2}_{q+1}  (D^{\balpha_3}_{q+1} - O^{\balpha_3}_{q+1}) \dots \Theta^{\balpha_k}_{q+1} F \|_{\infty} 
      \lesssim \delta^{\sfrac12} _{q+1} \lambda^{-\gamma N_*+2}_{q+1}  (\delta_{q+1}^{\sfrac12} \lambda_{q+1})^{|\balpha|_t} \lambda_{q+1}^{|\balpha|_x}  \text{ etc},
\end{equation}
where 'etc' means that the innermost and outermost operator actiong on $f$ may be either $\Theta^{\balpha'}_{q+1}$ or  $(D^{\balpha'}_{q+1} - O^{\balpha'}_{q+1})$  (so there are three more cases for l.h.s. of \eqref{e:streamEs_psi_DOdiff_gen}).
\end{corollary}

\begin{corollary}\label{l:strDiffEst3}
Fix positive $M$. For any smooth function $f:\T^2 \times [0,1] \to \R$, satisfying 
\begin{equation}\label{e:f_alastream_assD}
   \emph{either}  \quad  \| D^{\balpha}_{q+1} f \|_{\infty} 
      \lesssim M (\delta_{q+1}^{\sfrac12} \lambda_{q+1})^{|\balpha|_t} \lambda_{q+1}^{|\balpha|_x} \qquad \text{for any }  |\balpha| \leq Y-N_*-n_0.
\end{equation}
\begin{equation}\label{e:f_alastream_assO}
   \emph{or} \quad    \| O^{\balpha}_{q+1} f \|_{\infty} 
      \lesssim M (\delta_{q+1}^{\sfrac12} \lambda_{q+1})^{|\balpha|_t} \lambda_{q+1}^{|\balpha|_x} \qquad \text{for any }  |\balpha| \leq Y-N_*-n_0.
\end{equation}
with any $n_0 \ge 1$, $n_0 \le Y-N_*-1$, the analogues of estimates \eqref{e:streamEs_psi_DOdiff}, \eqref{e:streamEs_psi_DOdiff_gen} hold for any $|\balpha| \leq Y-N_*-n_0$. More precisely, we have for any $|\balpha| \leq Y-N_*-n_0$
\begin{equation}\label{e:streamEs_psi_DOdiff_genf}
      \| (D^{\balpha}_{q+1}- O^{\balpha}_{q+1}) f \|_{\infty} 
      \lesssim M \delta^{\sfrac12} _{q+1} \lambda^{-\gamma N_*+2}_{q+1}  (\delta_{q+1}^{\sfrac12} \lambda_{q+1})^{|\balpha|_t} \lambda_{q+1}^{|\balpha|_x}.
\end{equation}
\begin{equation}\label{e:streamEs_psi_DOdiff_gen_genf}
       \|  (D^{\balpha_1}_{q+1} - O^{\balpha_1}_{q+1}) \Theta^{\balpha_2}_{q+1}  (D^{\balpha_3}_{q+1} - O^{\balpha_3}_{q+1}) \dots \Theta^{\balpha_k}_{q+1} f \|_{\infty} 
      \lesssim M \delta^{\sfrac12} _{q+1} \lambda^{-\gamma N_*+2}_{q+1} (\delta_{q+1}^{\sfrac12} \lambda_{q+1})^{|\balpha|_t} \lambda_{q+1}^{|\balpha|_x} \text{ etc}.
\end{equation}
\end{corollary}

Lemma \ref{l:streamEst_rig_nobar} and its Corollaries are proven in Section \ref{sec:streamEst_rig_nobar_pf}.

\begin{corollary}     \label{c:CalPEst}
Recall the space $\mathcal{P}_q$ defined in \eqref{e:0:poly_Uq}. For any $n_0 \in \N$, and any $f \in \mathcal{P}_q (n_0)$, we have the following estimates for $|\balpha| + n_0 \leq Y$,
\begin{align}     \label{e:2:CalPEst}
      \| \mr D^{\balpha}_q f \|_{\infty} 
            \lesssim \vertiii{f}_q \lambda_q^{-|\balpha|\gamma_R}.
\end{align}
If the differentiation order of $f$ is bounded from below by some $n_* \in \N$. Then we have
\begin{align}     \label{e:4:CalPEst}
      \| \mr D^{\balpha}_q f \|_{\infty} 
            \lesssim \vertiii{f}_q \lambda_q^{-(|\balpha|+n_*)\gamma_R}.
\end{align}
\end{corollary}

\begin{proof}
Without loss of generality, we consider $f \in \mathcal{P}_{q} (n_0)$ that is not constant, we consider the representation of $f$ that attains $\vertiii{f}_{q+1}$, i.e. 
\begin{align}     \label{e:10:CalPEst}
      f = \sum_{k=1}^{k_0} a_k \prod_{i=1}^{I_k} \mr D_{q}^{\balpha_{k,i}} \Bigg( \frac{\mr \lambda_{q}^2} {\delta_{q}^{\sfrac12}\lambda_{q}^{1+2\gamma_R}} \bar \Psi_{q} \Bigg), \quad 
       \sum_{k=1}^{k_0} |a_k| = \vertiii{g}_{q} .
\end{align}
Applying Lemma \ref{l:streamEst}, we have
\begin{align}     
      \Bigg\| \mr D_{q}^{\balpha_{k,i}} \Bigg( \frac{\mr \lambda_{q}^2} {\delta_{q}^{\sfrac12}\lambda_{q}^{1+2\gamma_R}} \bar \Psi_{q} \Bigg) \Bigg\|_\infty 
      \lesssim&\, \lambda_q^{-2\gamma_R}
            \bigg( \frac{\lambda_q}{\mr \lambda_q} \bigg)^{|\balpha_{k,i}|_x-2}
            \bigg( \frac{\delta_q^{\sfrac12}\lambda_q} {\mr \mu_q} \bigg)^{|\balpha_{k,i}|_t} 
      \lesssim \lambda_q^{-|\balpha_{k,i}|\gamma_R} .      
                  \label{e:14:CalPEst}
\end{align}
Here we use \eqref{e:0:mScaleRela}, \eqref{e:smallGammaR} and \eqref{e:mrParameter}.

Now \eqref{e:2:CalPEst} follows from \eqref{e:14:CalPEst}. Recalling Definition \ref{d:polyOrderLBd}, \eqref{e:4:CalPEst} follows similarly.
\end{proof}

\begin{corollary}     \label{c:CalOEst}
Recall the spaces $\mathcal{O}_q$ defined in \eqref{e:2:poly_Oq}. For any $n_0 \in \N$, and any $f \in \mathcal{O}_q (n_0)$, we have the following estimates for $|\balpha| + n_0 \leq Y$,
\begin{align}     \label{e:2:CalOEst}
      \| \mr O^{\balpha}_q f \|_{\infty} 
            \lesssim \vertiii{f}_q \lambda_q^{-|\balpha|\gamma_R}.
\end{align}
If the differentiation order of $f$ is bounded from below by some $n_* \in \N$. Then we have
\begin{align}     \label{e:4:CalOEst}
      \| \mr O^{\balpha}_q f \|_{\infty} 
            \lesssim \vertiii{f}_q \lambda_q^{-(|\balpha|+n_*)\gamma_R}.
\end{align}
\end{corollary}

\begin{corollary}       \label{c:CongruEst}
For any $f \in \mathcal{O}_q (n_0)$ and $g \in \mathcal{P}_q (n_0)$, suppose $f$ and $g$ are congruent. Then we have the following estimates for any $\balpha \in \mathcal{I}$ with $|\balpha| + n_0 \leq Y - N_*$,
\begin{align}
      \| O^{\balpha}_q f - D^{\balpha}_q g \|_{\infty} 
      \lesssim
            \vertiii{f}_q \lambda_q^{-|\balpha|\gamma_R - 3}.
\end{align}
\end{corollary}

\begin{proof}
This is a direct consequence of Lemma \ref{l:streamEst_rig_nobar} and Corollary \ref{l:strDiffEst2}.
\end{proof}

\subsection{Proof of Lemma  \ref{l:uEstimate}}\label{sec:l:uEstimate_pf}
Below $D^{m-1}_{q}$ is $m-1$-th advective derivative with flow $u_q$.
\subsubsection{Proof of identity \eqref{e:6:uEstimate}}
Formula \eqref{e:10:constructVF} gives
\[
     \nabla \Phi_{q,\iota} = \Id + \sum_{m=1}^{N_*} \frac{1}{m!} ( t_{q,\iota} - t )^m  \nabla D_q^{m-1} u_q.
\]
Using that $t_{q,\iota} = \frac{\iota}{\mu_{q+1}}$ and that $\tilde \eta_{q,\iota} =\tilde \eta( \mu_{q+1} t - \iota)$ (cf \eqref{e:8:cutoff},
\eqref{e:6:constructVF}) we obtain 
\begin{equation}\label{e:6:uEstimate_pre}
     \nabla  \Phi_q = \sum_{\iota}  \nabla \Phi_{q,\iota} \tilde \eta_{q,\iota}
      =\Id + \varepsilon_q \sum_{m=1}^{N_*} \underbrace{\frac{1}{\varepsilon_q \mu^m_{q+1}} \nabla D_q^{m-1} u_q}_{=:B_{q,m} (x,t)} \sum_{\iota}  \underbrace{ \frac{1}{m!}( \iota - \mu_{q+1} t)^m  \tilde \eta( \mu_{q+1} t - \iota)}_{=:\phi_{m, \iota} (\mu_{q+1} t)}.
\end{equation}
Defining $\phi_{m}:= \sum_{\iota} \phi_{m, \iota}$ we have \eqref{e:6:uEstimate}.

\begin{remark}[$N_*$ vs $N$]
    Observe that \eqref{e:6:uEstimate_pre} is a sum over $i=1,\dots, N_*$, whereas \eqref{e:6:uEstimate}  is a sum over $i=1,\dots, N$, with $N\ge  N_*$. This is realised by extending the sum in \eqref{e:6:uEstimate_pre} with dummy entries (zeros). The reason is notation: the forthcoming representation will consist of sums up to $N^2_*$ etc, and in order not to keep track of the precise value of the largest index in final formulas, we choose the upper bound for all of them, i.e.\ $N$, and extend shorter sums with dummy entries.
\end{remark}

\subsubsection{Proof of identity \eqref{e:7:uEstimate}} Since 
\[
\begin{aligned}
\det (\Id + \underbrace{\sum^{N_*}_{m=1}A_m}_{=:A} ) &= 1 + \det A + tr A\\
&= 1 + \sum^{N_*}_{m=1} \sum_{m'=1}^{N_*} (A_{m,11} A_{m',22} - A_{m,21} A_{m',12}) + \sum^{N_*}_{m=1} A_{m,11} A_{m,22}
\end{aligned}
\]
We have via \eqref{e:6:uEstimate}
\begin{equation}\label{e:6p:uE}
\begin{aligned}
\det  \nabla  \Phi_q 
=& 1 +  \varepsilon^2_q \left(\sum^{N_*}_{m=1} \sum_{m'=1}^{N_*} ({B}_{q,m,11} {B}_{q,m',22 }- {B}_{q,m,21}    {B}_{q,m',12})  \phi_{m}  \phi_{m'} \right)\\
&\quad +  \varepsilon_q \sum^{N_*}_{m=1} ({B}_{q,m,11} + {B}_{q,m,22})  \phi_{m} 
\end{aligned}
\end{equation}

\begin{remark}[Multiple indices into long sums]
We will order terms with two (and more) indices into terms with single index according to the following version of enumerating product of two (and more) countable sets: having element $a_{m,m'}$ we define $b_j$ as follows:
\[
b_j = \begin{cases}
 a_{l,l} &\text{ if } j=l^2 \\
  a_{i,l} &\text{ if } j = (l-1)^2+i, \; i=1,\dots, l-1 \\
   a_{l,i} &\text{ if } j = (l-1)^2+i, \; i=l,\dots, 2l-1
\end{cases}
\]

Importantly for precise control of constants in estimates
\begin{equation}\label{e:enum_remark}
\begin{aligned}
(i)& \text{
 this enumeration does not involve $N_*$ }\\
(ii)& m,m' \le j \text{ i.e. any of the double indices is not larger than the single index.} 
\end{aligned}
\end{equation}
The full extent of the constants control is not necessary in the current paper.
\end{remark}

In particular, define
$\tilde B_{q,k,l} = {B}_{q,k,11} {B}_{q,l,22 }- {B}_{q,k,21}    {B}_{q,l,12} $
and with this
\begin{equation}\label{e:de_sigmaOmega}
\begin{aligned}
&\Omega^{quad}_{q,j} = 
\begin{cases}
\tilde B_{q,m,m} &\text{ if } j=m^2 \\
\tilde B_{q,m,i} &\text{ if } j = (m-1)^2+i, \; i=1,\dots, m-1 \\
\tilde B_{q,i,m}  &\text{ if } j = (m-1)^2+i, \; i=m,\dots, 2m-1
\end{cases} \\
&\sigma_j =
\begin{cases}
\phi^2_{m}  &\text{ if } j=m^2 \\
\phi_{m}  \phi_{i} &\text{ if } j = (m-1)^2+i, \; i=1,\dots, m-1 \\
\phi_{i}  \phi_{m} &\text{ if } j = (m-1)^2+i, \; i=m,\dots, 2m-1
\end{cases} \\
&\Omega^{lin}_{q,j} =  ({B}_{q,j,11} + {B}_{q,j,22}) \\
\end{aligned}
\end{equation}
Formulas \eqref{e:de_sigmaOmega} used in \eqref{e:6p:uE} give 
\begin{equation}\label{e:7:uEstimate_0}
\det \nabla \Phi_q =1 +  \varepsilon^2_q \sum^{{N^2_*}}_{m=1} \Omega^{quad}_{q,m} \sigma_m +  \varepsilon_q \sum^{{N_*}}_{m=1} \Omega^{lin}_{q,m} \phi_m.
\end{equation}
Therefore, defining 
\[
    \zeta_m = \begin{cases}
        \phi_m &\text{ for } m=1, \dots, N^* \\
        \sigma_{m-N^*} &\text{ for } m=N^*+1, \dots, N^2_* + N^*+1\\
    \end{cases}
\]
and 
\[
    \varpi_{q,m} = \begin{cases}
        \Omega^{lin}_{q,m}  &\text{ for } m=1, \dots, N^* \\
        \varepsilon_q\Omega^{quad}_{q,m-N^*} &\text{ for } m=N^*+1, \dots, N^2_* + N^*+1\\
    \end{cases}
\]
we have \eqref{e:7:uEstimate}.

\subsubsection{Proof of identity \eqref{e:10:uEstimate}}
 
 The formula \eqref{e:6:uEstimate_pre} for $\nabla \Phi_q$ gives:
\begin{align*}
     \nabla \Phi_q \nabla \Phi^T_q
     &\, = \Id + \varepsilon_q \sum_{m=1}^{N_*} (B_{q,m} + B^T_{q,m})\phi_{m}  + \varepsilon^2_q \sum_{m, m'=1}^{N_*} B_{q,m} B^T_{q,m'}  \phi_{m}  \phi_{m'}  \\
     &=: \Id + \varepsilon^2_q \sum_{m=1}^{N^2_*} S^{quadr}_{q,m}\sigma_m +  \varepsilon_q \sum^{{N_*}}_{m=1} S^{lin}_{q,m} \phi_m
\end{align*}  where
\begin{equation}\label{e:de_S}
\begin{aligned}
&S^{quad}_{q,j} := 
\begin{cases}
{B}_{q,m} {B^T}_{q,m} &\text{ if } j=m^2 \\
{B}_{q,m} {B^T}_{q,i}  &\text{ if } j = (m-1)^2+i, \; i=1,\dots, m-1 \\
{B}_{q,i} {B^T}_{q,m}  &\text{ if } j = (m-1)^2+i, \; i=m,\dots, 2m-1
\end{cases} \\
&S^{lin}_{q,j} :=  {B}_{q,j} +  {B^T}_{q,j}
\end{aligned}
\end{equation}

The preceding two computations together 
give
\begin{equation}\label{e:10uEstimate_pre}
\begin{aligned}
     \nabla \Phi_q \nabla \Phi^T_q -(\det  \nabla  \Phi_q)^2 \Id 
    =  &\varepsilon^4_q \sum^{{N^4_*}}_{m=1} S^{1,quadric}_{q,m} \vartheta^1_m + \varepsilon^3_q \sum^{{N^3_*}}_{m=1} S^{1,cubic}_{q,m} \varphi_m  \\
    &+ \varepsilon^2_q \sum^{{N^2_*}}_{m=1} S^{1,quad}_{q,m} \sigma_m + \varepsilon_q \sum^{{N_*}}_{m=1} S^{1,lin}_{q,m} \phi_m 
    \end{aligned}
\end{equation}
where
\begin{equation}\label{e:de_S1t1}
\begin{aligned}
S^{1,quadric}_{q,j} &:= -
\begin{cases}
\Omega^{quad}_{q,m} \Omega^{quad}_{q,m} \Id &\text{ if } j=m^2 \\
\Omega^{quad}_{q,m} \Omega^{quad}_{q,i} \Id&\text{ if } j = (m-1)^2+i, \; i=1,\dots, m-1 \\
\Omega^{quad}_{q,i} \Omega^{quad}_{q,m} \Id &\text{ if } j = (m-1)^2+i, \; i=l,\dots, 2m-1
\end{cases} \\
\vartheta^1_j &:=
\begin{cases}
\sigma^2_{m}  &\text{ if } j=m^2 \\
\sigma_{m}  \sigma_{i} &\text{ if } j = (m-1)^2+i, \; i=1,\dots, m-1 \\
\sigma_{i}  \sigma_{m} &\text{ if } j = (m-1)^2+i, \; i=m,\dots, 2m-1
\end{cases} \\
S^{1,quad}_{q,j} &:= S^{quadr}_{q,j} - 4 \Omega^{quadr}_{q,j}\\ 
S^{1,lin}_{q,j} &:=  S^{lin}_{q,j} - 2 \Omega^{lin}_{q,j}
\end{aligned}
\end{equation}
(where we used in $S^{1,quad}_{q,j}$ that product of two $\Omega^{lin}_{q,j}$ is $\Omega^{quad}_{q,j}$). The remaining cubic term consists of products of quadratic and linear terms, so recalling \eqref{e:de_sigmaOmega},
we write
\begin{equation}\label{e:de_S1cubic}
\begin{aligned}
&S^{1,cubic}_{q,j}  := 
\begin{cases}
\tilde B_{q,m,m} \Omega^{lin}_{q,m} &\text{ if } j=m^3 \\
\tilde B_{q,k,i} \Omega^{lin}_{q,m} &\text{ accordingly } \\
\text{ etc, along each side of the cube}
\end{cases} \\
&\varphi_j :=
\begin{cases}
\phi^3_{m}  &\text{ if } j=m^3 \\
\phi_{k} \phi_{i} \phi_{m}  &\text{ accordingly } \\
\text{ etc, along each side of the cube.}
\end{cases}
\end{aligned}
\end{equation}
Formula \eqref{e:10uEstimate_pre} with the preceding explicitly defined ingredients yields \eqref{e:10:uEstimate}, by appropriate shifts of indices of respective ingredients (analogously to how \eqref{e:7:uEstimate_0} gives \eqref{e:7:uEstimate}).

\subsubsection{Proof of identity \eqref{e:8:uEstimate}}

Since $\adj A^{T} = (\adj A)^T$ and in 2d $\adj (\Id +A)^{T}= \Id + \adj A^{T}$ and $\adj$ is linear, we have

\begin{equation}
\label{e:8:uEstimate:pre}
\begin{aligned}
{\adj \nabla \Phi_q^{T} \det \nabla \Phi_q} 
= \Id + \varepsilon^3_q \sum^{{N^3_*}}_{m=1} (E^{cubic}_{q,m})^T  \varphi_m  + \varepsilon^2_q \sum^{{N^2_*}}_{m=1} (E^{quadr}_{q,m})^T  \sigma_m + \varepsilon_q \sum^{{N_*}}_{m=1} (E^{lin}_{q,m})^T  \phi_m 
\end{aligned}
\end{equation}
with (recall $ \Omega_{q,m}$ and  $\tilde B_{q,k,i}$ are numbers, cf \eqref{e:de_sigmaOmega})
\begin{equation}\label{e:de_E}
\begin{aligned}
&E^{cubic}_{q,j}  = 
\begin{cases}
\tilde B_{q,m,m} \adj B_{q,m}  &\text{ if } j=m^3 \\
\tilde B_{q,k,i} \adj B_{q,m}  &\text{ accordingly } \\
\text{ etc, along each side of the cube}
\end{cases} \\
&E^{quadr}_{q,j}  = 
\begin{cases}
\tilde \Omega^{quad}_{q,m} + \Omega^{lin}_{q,m} \adj B_{q,m}  &\text{ if } j=m^2 \\
\tilde B_{q,k,i} \adj B_{q,m}  &\text{ accordingly } \\
\text{ symmetrically}
\end{cases} \\
&E^{lin}_{q,j}  = 
\Omega^{lin}_{q,j}+ \adj B^T_{q,j}  
\end{aligned}
\end{equation}
Formula \eqref{e:8:uEstimate:pre} gives \eqref{e:8:uEstimate} (as \eqref{e:10uEstimate_pre} has given \eqref{e:10:uEstimate}).

\subsubsection{Proof of identity \eqref{e:12:uEstimate}}
Thanks to skew-symmetry of $H$ we have $AHA^T = \det(A) H$ and thus for a scalar function $c$ it holds $\divr (cH) = \nabla^\perp (c H_{12})$. In particular:
\begin{equation}\label{e:adj_form_simple}
\begin{aligned}
\adj \nabla \Phi_q H  \adj \nabla \Phi^{T}_q &= \nabla \Phi_q^{-1} \det \nabla \Phi_q H \nabla \Phi_q^{-T} \det \nabla \Phi_q = (\det \nabla \Phi_q)^2 \nabla \Phi_q^{-1}  H \nabla \Phi_q^{-T} \\
&= \det \nabla \Phi_q  H.  
\end{aligned}
\end{equation}
This and the identity \eqref{e:7:uEstimate} gives \eqref{e:12:uEstimate}, with 
\[
\sigma_m:= \zeta_{m} \quad \text{and} \quad \Omega_{q,m} :=\begin{bmatrix} 0 & 1 \\ -1 & 0  \end{bmatrix}\varpi_{q,m}.
\]

\subsubsection{Proof of identity \eqref{e:2:uEstimate}}
 The definition \eqref{e:10:constructVF} yields cancellations in the following computation:
\begin{align*}
      D_q\Phi_{q,\iota} =&\sum_{m=1}^{N_*} \frac{1}{m!} D_q (( t_{q,\iota} - t )^m D_q^{m-1} u_q)\\
      =&\sum_{m=0}^{N_*-1} \frac{-1}{m!} (( t_{q,\iota} - t )^{m} D_q^{m} u_q) + \sum_{m=1}^{N_*} \frac{1}{m!} (( t_{q,\iota} - t )^{m} D_q^{m} u_q) \\
       =& - D_q^{0}  u_q + \frac{1}{N_*!} (( t_{q,\iota} - t )^{N_*} D_q^{N_*} u_q)
\end{align*}
recalling that by convention $D_t^{0} =0$, we thus have
\begin{equation}\label{e:app_flow_loc}
\begin{aligned}
          \sum_\iota \tilde \eta_{q,\iota}   D_q\Phi_{q,\iota} =&     \sum_\iota \tilde \eta_{q,\iota}  \frac{1}{N_*!} (( t_{q,\iota} - t )^{N_*} D_q^{N_*} u_q) \\
       =& \underbrace{ \frac{1}{ \mu^{N_*}_{q+1}} D_q^{N_*} u_q}_{=:\omega_{q,{N_*}} (x,t)} \sum_{\iota}  \underbrace{\frac{1}{{N_*}!} ( \iota - \mu_{q+1} t)^{N_*}  \tilde \eta( \mu_{q+1} t - \iota)}_{=\phi_{{N_*}, \iota} (\mu_{q+1} t)}   
\end{aligned}
\end{equation}
which is indeed \eqref{e:2:uEstimate}. 

\subsubsection{Proof of identity  \eqref{e:4:uEstimate}}
Now compute, using again  $t_{q,\iota} = \frac{\iota}{\mu_{q+1}}$ and that $\tilde \eta_{q,\iota} =\tilde \eta( \mu_{q+1} t - \iota)$ 
\begin{equation}
\begin{aligned}
     &\sum_\iota \partial_t \tilde \eta_{q,\iota} \Phi_{q,\iota} \\
     &=  \sum_\iota \mu_{q+1} (\partial_\tau \tilde \eta)( \mu_{q+1} t - \iota) (x + \sum_{m=1}^{N_*} \frac{1}{m!} ( t_{q,\iota} - t )^m D_q^{m-1} u_q)\\ 
     &= \sum_\iota \mu_{q+1} (\partial_\tau \tilde \eta)( \mu_{q+1} t - \iota) x +   \sum_{m=1}^{N_*}  \sum_\iota \frac{1}{m!} ( \iota - \mu_{q+1} t)^m(\partial_\tau \tilde \eta)( \mu_{q+1} t - \iota) (\frac{1}{\mu_{q+1}^{m-1}}  D_q^{m-1} u_q) \\ 
     &= \mu_{q+1} \sum_\iota  \underbrace{(\partial_\tau \tilde \eta)( \mu_{q+1} t - \iota)}_{=:z_{0,\iota }} x +   \sum_{m=1}^{N_*}  \underbrace{(\frac{1}{\mu_{q+1}^{m-1}} D_q^{m-1} u_q)}_{=:\omega_{q,{m-1}} (x,t)} \sum_\iota   \underbrace{ \frac{1}{m!} ( \iota - \mu_{q+1} t)^m(\partial_\tau \tilde \eta)( \mu_{q+1} t - \iota)}_{=:z_{m,\iota }} \label{e:where_ommdef}.
\end{aligned}
\end{equation}
Defining for $m \ge 0$
\begin{equation}\label{e:def_z}
z_{m}:= \sum_{\iota} z_{m, \iota},
\end{equation}
with $z_{m, \iota}$ defined in the preceding formula,
we obtain \eqref{e:4:uEstimate}.

\subsubsection{ 
Statements \eqref{e:algbVarpiOmega}, \eqref{e:algbBESOmega}} They follow immediately from the respective definitions.

\subsubsection{Proof of estimates \eqref{e:phiZetaThetaEst}} 

We will use the following very lossy, but sufficient to prove \eqref{e:phiZetaThetaEst}, inequality 
\begin{equation}\label{e:vs_inter}
|\partial_\tau^p (fg)| \le 2^p|\sup_{i \le p} \partial_\tau^i f| |\sup_{j \le p}  \partial_\tau^j g|  .
\end{equation}
Using \eqref{e:vs_inter} for
\[
\phi_{m, \iota} (\tau) = \frac{1}{m!}  ( \iota - \tau)^m  \tilde \eta( \tau - \iota),
\]
yields via estimate \eqref{e:6:cutoff} for $ \tilde \eta$, and the fact that support of $  \tilde \eta$ yields $|\iota - \tau|\le 1$ 
\[
 | \partial_\tau^p \phi_{m, \iota} | \le 2^p c_0 p^{3p}
\]
and observing that in definition of $\phi_{m} (\tau)$ at most two $\phi_{m, \iota}$'s overlap, we have 
\begin{equation}
          \| \partial_\tau^p \phi_m \|_\infty \leq c_0    2^{p+1} p^{3p} \lesssim 1, \label{e:phiEst} 
\end{equation}
where $\lesssim 1$ follows by the assumed absolute bound on total amount of derivatives, $p\le 8Q^3$.
The estimate \eqref{e:phiEst} is the part of \eqref{e:phiZetaThetaEst} where $\| \partial_\tau^p \phi_m \|_\infty$ is estimated. 
Analogously, via the definition \eqref{e:def_z} of $z_m$ we have the part of \eqref{e:phiZetaThetaEst} where $\| \partial_\tau^p z_m \|_\infty$ is estimated. 
The product $\phi_{m} \phi_{m'}$ satisfies via \eqref{e:phiEst}
\[
\begin{aligned}
 | \partial_\tau^p (\phi_{m} \phi_{m'}) | &\le \sum_{i=0}^{p} {p \choose i}  | \partial_\tau^i \phi_{m} | |\partial_\tau^{p-i} \phi_{m'}| \le \sum_{i=0}^{p} {p \choose i}  2^{i+1}  c_0 i^{3i} 2^{p+1-i}  c_0 {(p-i)}^{3(p-i)} \\
 &\le  c_0^2 2^{p+2}  \sum_{i=0}^{p} {p \choose i}  ({p^3})^{i} ({p^3})^{(p-i)} = c_0^2 2^{2p+2}  p^{3p}
 \end{aligned}
\]
In particular, the quadratic cutoffs $\sigma_m$, in view of their definition \eqref{e:de_sigmaOmega}, compare remark \eqref{e:enum_remark}, satisfy 
the estimate $ \| \partial_\tau^p \sigma_m\|_\infty \lesssim 1$. The remainder of the required estimate \eqref{e:phiZetaThetaEst} follows similarly.

\subsubsection{Coefficients estimates \eqref{e:mtrPolyEst}, \eqref{e:trsptRsdPolyEst}}

Consider $B_{q,m}$ defined in \eqref{e:6:uEstimate_pre}. This and $u_{q} = \nabla^\perp \bar \Psi_q$ given by \eqref{e:26:constructVF} gives
\[
\begin{aligned}
B_{q,m} &= \frac{1}{\varepsilon_q \mu^m_{q+1}} \nabla D_q^{m-1} \nabla^\perp \bar \Psi_q \\
&\stackrel{\text{def }\ref{d:renormalizedDiff}}{=} \frac{1}{\varepsilon_q \mu^m_{q+1}} \mr \lambda_q^{2} \mr \mu_q^{m-1}  \mr D_q^{\bomega}  \bar \Psi_q \qquad \big(\text{ with }\bomega = ( \partial, \underbrace{D_q}_{m-1 \text{ times}}, \partial ) \big)\\
&= \frac{\mu_q}{ \mu^m_{q+1}} \mr \mu_q^{m-1}  \mr D_q^{\bomega} \Bigg( \frac{\mr \lambda_{q}^2} {\delta_{q}^{\sfrac12}\lambda_{q}^{1+2\gamma_R}} \bar \Psi_{q} \Bigg)
\end{aligned} 
\]
where the last identity holds just because    $\varepsilon_q \mu_q 
            = \delta_{q}^{\sfrac12}\lambda_{q}^{1+2\gamma_R}$, recall Remark \ref{r:renormal_factor}. The relation $\mu_{q+1} > \mr \mu_q > \mu_q$, recall Remark \ref{r:auxiliaryParaRela}, implies \eqref{e:mtrPolyEst} for $B_{q,m}$, via definition of the quantity $\vertiii{\cdot }_q$, see \eqref{e:0:polyNorm}. 

            Inspecting the formulas for the rest of r.h.s. of \eqref{e:mtrPolyEst} given in part (i) of the proof, we see that they are products related to $B_{q,m}$. Since the quantity $\vertiii{\cdot }_q$ behaves good under multiplication, see Remark \ref{r:opPolynomial}, we have the entire \eqref{e:mtrPolyEst}, except for estimate of $\omega_{q,{m}}$.

The proof of \eqref{e:mtrPolyEst} for $\omega_{q,{m}}$ and of \eqref{e:trsptRsdPolyEst} is similar. Indeed, recalling that by \eqref{e:where_ommdef}
\[
\omega_{q,{m}} =\frac{1}{ \mu^{m}_{q+1}} D_q^{m} u_q 
\]
we have 
\[
\omega_{q,{m}} = \frac{\varepsilon_q \mu_q}{ \mr \lambda_{q} \mu^m_{q+1}} \mr \mu_q^{m}  \mr D_q^{\bomega} \Bigg( \frac{\mr \lambda_{q}^2} {\delta_{q}^{\sfrac12}\lambda_{q}^{1+2\gamma_R}} \bar \Psi_{q} \Bigg)
\]
As for $B_{q,m}$, the prefactor before $D_q^{\bomega}$ is small, so we have \eqref{e:mtrPolyEst} for $\omega_{q,{m}}$. In particular the prefactor $\frac{\varepsilon_q \mu_q}{ \mr \lambda_{q} \mu^m_{q+1}} \mr \mu_q^{m} \le \frac{ \mr \mu_q^{m}}{\mu^{m}_{q+1}} $ is as small as we want for
 $m= N_*$, with $N_*$ sufficiently large.  Thus via the (safe) choice of $N_*$ in \eqref{e:NQ}, we obtain \eqref{e:trsptRsdPolyEst}.

 \subsection{Proof of Lemma \ref{l:streamEst}}\label{sec:streamEst_rig1_pf}

The wanted estimates will be proven via induction over $q \to q+1$.  Recall that \[
\bar  \Psi_{q+1} = \bar \Psi_q + \bar \psi_{q+1},  
\]
where by definition \eqref{e:psiQ+1}
\begin{equation}\label{e:recall_psiQ+1}
\bar \psi_{q+1} = \frac{\delta_{q+1}^{\sfrac12}}{\lambda_{q+1}} 
            H_{12} \big( \lambda_{q+1} \bar \Phi_q(x,t), \mu_{q+1} t \big).     
\end{equation}

\subsubsection{Initial remarks}\label{sssec:ir_ests}

Let us start with explaining why we propagate all space derivatives and a finite number of mixed derivatives. Very briefly, it is caused by (standard) loss of space derivative. More precisely, there is

\paragraph{(i. Loss of space derivative)}
Indeed, recall from its definition  \eqref{e:8:constructVF}, that $\bar \Phi_{q, \iota}$ (the local, backwards flow) satisfies $\partial_t \bar \Phi_{q,\iota} + ( u_q \cdot \nabla ) \bar \Phi_{q,\iota} =\, 0$, which we denote briefly as $D_q \bar \Phi_{q, \iota}=0$. Taking space derivative, we have
\begin{equation}\label{e:tr_hof}
\begin{aligned}
D_q \nabla \bar \Phi_{q, \iota} &= -\nabla  u_q \cdot  \nabla \bar \Phi_{q, \iota}\\
&=-\nabla  \nabla^T \bar \Psi_q \cdot  \nabla \bar \Phi_{q, \iota},
\end{aligned}
\end{equation}
the latter via the definition \eqref{e:26:constructVF} of $u_q$. Since $\nabla \bar \Phi_{q, \iota} = 1$ at its initial time $t_{q, \iota}$, the standard Gronwall-type estimate for the forced transport equation \eqref{e:tr_hof} yields:
\begin{equation}\label{e:loss_sd}
\| \nabla \bar \Phi_{q, \iota} (t)\|_\infty \le \exp \int_{t_0}^t \| \nabla  \nabla^T \bar \Psi_q (s) \|_\infty ds 
\end{equation}
In particular, estimates on $\nabla \bar \Phi_{q, \iota}$ require knowledge of estimates of $\nabla^2  \bar \Psi_q$. Taking space derivative in \eqref{e:recall_psiQ+1} we have
\[
\nabla \bar \psi_{q+1} = \delta_{q+1}^{\sfrac12} \nabla \bar \Phi_q(x,t) 
            (\partial_\xi H_{12}) \big( \lambda_{q+1} \bar \Phi_q(x,t), \mu_{q+1} t \big). 
\]
which together with \eqref{e:loss_sd} means that 
estimates on $\nabla \bar \psi_{q+1}$ and thus on $\nabla \bar \Psi_{q+1}$ require knowledge of estimates of $\nabla^2  \bar \Psi_q$. This is the mentioned before standard loss of derivative phenomenon. It persists for any higher space derivative, since the above reasoning via the computation
\begin{equation}\label{e:tr_hofHi}
\begin{aligned}
D_q \partial_x^n \bar \Phi_{q, \iota} &\sim_n \partial_x^n  u_q \cdot  \nabla \bar \Phi_{q, \iota} + \partial_x^{n-1}  u_q \cdot  \nabla \partial_x \bar \Phi_{q, \iota} + \dots + \partial_x  u_q \cdot  \nabla \partial^{n-1}_x \bar \Phi_{q, \iota} \\
 &= \partial_x^n  \nabla^T \bar \Psi_q \cdot  \nabla \bar \Phi_{q, \iota} + \partial_x^{n-1}  \nabla^T \bar \Psi_q \cdot  \nabla \partial_x \bar \Phi_{q, \iota} + \dots + \partial_x  \nabla^T \bar \Psi_q \cdot  \nabla \partial^{n-1}_x \bar \Phi_{q, \iota}
\end{aligned}
\end{equation}
and chain rule applied to \eqref{e:recall_psiQ+1}, shows that estimates on $\nabla^n \bar \Psi_{q+1}$ require knowledge of estimates of $\nabla^{n+1}  \bar \Psi_q$. 

\paragraph{(ii.  No loss in mixed derivatives)} However, if we want to estimate $D_{q+1}^{\bomega} \bar  \Psi_{q+1}$ with ${\bomega}$ containing at least one $t$ (i.e.\ with at least one advective derivative), the situation is different. Assuming for a moment that the commutator $[D_{q+1}^{\bomega}, D_{q}^{\bomega}]$ causes no problems, we want to estimate $D_{q}^{\bomega} \bar  \Psi_{q+1}$. In the simplest case $D_{q}^{\bomega} = D_q \nabla$ we see from \eqref{e:tr_hof} that when estimating directly its l.h.s. with r.h.s., there is no loss of (order of) derivatives. This holds for any $D_{q}^{\bomega}$. Indeed, any $D_{q}^{\bomega}$ is a sequence of blocks of space and transport derivatives. The case of the innermost derivative being the transport derivative yields $D_{q}^{\bomega} \bar  \Phi_{q, \iota}=0$, so we are always in the case $D_{q}^{\bomega}  = D_{q}^{\balpha}  D_{q} \nabla^n $, $n+1+ |\balpha| = |\bomega|$. The expression $ D_{q} \nabla^n  \bar  \Phi_{q, \iota}$ is as in \eqref{e:tr_hofHi}. In particular, we see that this (mixed) derivative  of $n+1$-order (l.h.s.) is equal to a sum of products of derivatives, where $\bar \Psi_q$ appears in at most $n+1$-order of derivative (r.h.s.). Outer $D_{q}^{\balpha}$ derivatives induce at most further $|\balpha|$ derivatives of $\bar \Psi_q$. All in all, we expect no loss of derivative in the mixed derivative case.\footnote{Of course, the estimate of $D_{q+1}^{\bomega} \bar  \Psi_{q+1} = D_{q+1}^{\bomega} (\bar \Psi_q + \bar \psi_{q+1})$ requires dealing  with $D_{q+1}^{\bomega}$ and not with $D_{q}^{\bomega}$, full flow $\bar\Phi_{q}$ and not $\bar\Phi_{q, \iota}$, etc. This initial discussion indicates only that in the case of mixed derivatives, there is no derivative loss mechanism, which appears when dealing with pure space derivatives.}

\paragraph{(iii. Propagation of derivatives)} To sum up, we have to deal with loss of derivative in case of pure space derivatives, and not in the case of mixed derivatives. 

There are two standard remedies to the derivative loss. Either one mollifies the quantities involved, which allows appropriately high derivative estimates to be given not by (lossy) inductive assumptions, but also by mollified quantities; or one chooses to propagate all derivatives. In the former case, one needs to quantify the difference between the original and mollified quantities, in the latter one needs to quantify and control the arising combinatorial coefficients  (appearing for instance in product rule, see e.g.\ \eqref{e:tr_hofHi} with its $\sim_n$, or in the chain rule). 
We will follow the latter approach, which has been used in \cite{ArmstrongVicol}.

Therefore we propagate all space derivatives, but we don't have to propagate all mixed derivatives, which allows us also to be less careful about constants in the mixed derivative case.

\begin{remark}[Amplitude and magnitude]
    We will use interchangeably 'magnitude of $f$' or 'amplitude of $f$' for $\| f \|_\infty$
\end{remark}

\subsubsection{Propagating pure space derivatives}
This section is analogous to section 2.3 of \cite{ArmstrongVicol}. We provide our proof for the sake of reader's convenience, in particular, since our language and that of \cite{ArmstrongVicol} differ, and it seems our proof is straightforward.
\begin{remark}(Comparison of name of scales with \cite{ArmstrongVicol})
    The name of the step in our case is $q$, in the case of \cite{ArmstrongVicol} is $m$. Keeping this in mind, the translation of scales is
    \[
    \varepsilon_m \sim \lambda^{-1}_q, \quad a_m =\delta^\frac{1}{2}_q \lambda_q,
    \]
    in particular $a_m  \varepsilon^2_m = \frac{\delta^{\sfrac{1}{2}}_q }{ \lambda_q}$ is the amplitude of the potential $\bar \psi_{q+1}$.
\end{remark}

We want to reproduce \eqref{e:streamEs_rig1}, i.e.\ \[
\| \nabla^n \bar \Psi_q \|_{\infty} 
      \le\, \frac{(n-1)!}{n^2} \tilde C (C{\lambda}_q)^{n-1}\delta_q^{\sfrac12}     \qquad \text{for any } n \ge 2
\] 
with  $q\to q+1$. Since
\[
\bar  \Psi_{q+1} = \bar \Psi_q + \bar \psi_{q+1},  \]
it is crucial to prove an appropriate estimate on space derivatives of \[
\bar \psi_{q+1} (x,t)=\frac{\delta_{q+1}^{\sfrac12}}{\lambda_{q+1}} 
            H_{12} \big( \lambda_{q+1} \bar \Phi_q(x,t), \mu_{q+1} t \big),\] recall \eqref{e:recall_psiQ+1}, where
the main point is to control carefully the combinatorial coefficients that arise when taking derivatives of an arbitrary order. To this end let us recall and tailor the setting and two results of \cite{ArmstrongVicol}.

\paragraph{(i. Tools for dealing with an arbitrary space derivative)}
Introduce analytic-norm-type quantity for $n\in\N$
\begin{align}
\label{e:antq}
\snorm{f}_{n,R} 
=  \frac{(n+1)^2}{ n!R^{n}} \sup_{|i| = n} \  \norm{\partial^i f}_{L^\infty(\R^2)} 
\qquad \mbox{and} \qquad 
\snorm{f}_{R} = \sup_{n\in\N} \ \snorm{f}_{n,R}
\,.
\end{align}
where $\partial^i$ denotes any space derivative of order $|i|$. 

In order to control space derivatives of $\bar \Phi_{q}$, let us use Lemma B.7 of \cite{ArmstrongVicol} with its parameters given as follows:
\[
Y := \bar \Phi_{q, \iota} - \Id, \quad \mathbf{f} := u_q, \quad \mathbf{g} := -u_q, \quad d:=2,
\]
initial time shifted to $t_0:=t_{q, \iota}$, and $N$ arbitrary. We get
\begin{lemma} 
\label{l:transport_regularity_analytic}
Assume there exist~$C, R$ and such that for any $n\ge 1$ 
\begin{equation}
\label{e:tra_ass}
\sup_{t\in \R} \
\snorm{ u_q(t,\cdot)}_{n,R} \leq C
    \end{equation}
Then, $\bar\Phi_{q, \iota}$ (solving \eqref{e:8:constructVF}) satisfies for any $n\ge 1$ 
\begin{align}
\sup_{t  \in [-T,T]} \ 
\frac{1}{|t|} \snorm{ \bar\Phi_{q, \iota} (t_{q, \iota} +t, \cdot)- Id}_{n,R_{\bar\Phi}(t)}    \leq 16 C
\end{align}
where 
\begin{equation}
R_{\bar\Phi}(t) := R(1+   \frac{|t|}{T}) 
\qquad \mbox{and} \qquad 
T := \frac{1}{8 C R  }.
\label{e:tra_consts}
\end{equation}
\end{lemma}
In order to deal with the composition appearing in the definition of $\bar \psi_{q+1}$, we recall Proposition B.6 of \cite{ArmstrongVicol} (using there  arbitrary $m$)

\begin{proposition}[Composition estimate]
\label{prop:compose:analytic}
Assume that there exist positive constants $C_h$, $C_g$, $R_h$, $R_g \in (0,\infty)$ such that 
\begin{equation*}
\snorm{h}_{n,R_h} \leq C_h
\qquad \mbox{and} \qquad 
\snorm{g}_{n_1,R_g} \leq C_g
\end{equation*}
for all $0 \leq n
$ and $1\leq n_1 
$. Then, for every  $0 \leq n$, we have 
\begin{align}
\snorm{h \circ g }_{n,R} \leq C_h \, 
\qquad \mbox{where} \qquad 
R
= R_g (1 + d C_g R_h)\,.
\label{e.barf.3}
\end{align}
\end{proposition}
\paragraph{(ii. Estimate on the flow map)}
Since $u_q = \nabla^\perp \bar \Psi_q$ our assumption 
\eqref{e:streamEs_rig1} is \[
\| \nabla^n \bar u_q  \|_{\infty} 
      \le\, \frac{n!}{(n+1)^2}({C\lambda}_q)^{n} \delta_q^{\sfrac12}  \tilde C    \qquad \text{for any } n \ge 1
\] 
i.e.
\begin{equation}\label{e:vel_ass_spaceders}
    \snorm{u_q}_{n,C \lambda_q} \leq \tilde C \delta_q^{\frac{1}{2}}  \qquad \text{for any } n \ge 1.
\end{equation}
This used in Lemma \ref{l:transport_regularity_analytic} gives that for $n\ge 1$
\begin{equation}
\snorm{ \bar\Phi_{q, \iota} (t_{q, \iota} +t, \cdot) - \Id}_{n,2 C \lambda_q}    \leq 16 \tilde C \delta_q^{\frac{1}{2}} |t|
\end{equation}
for any $|t| \le \frac{1}{8 \tilde C C \delta_q^{\frac{1}{2}} \lambda_q }$; where above we replaced $R_{\bar\Phi}(t)$ given by Lemma \ref{l:transport_regularity_analytic} with its upper bound:
\[R_{\bar\Phi}(t) =  C \lambda_q (1+   \sfrac{|t|}{T}) \le 2 C \lambda_q, 
\]
which is allowed because the norm-type-quantity $\snorm{ \cdot}_{n,R}$ is monotone in $R$ by definition. 
In particular any
\begin{equation}\label{e:grwl_ass_simple}
    \mu_{q+1} \ge 8 \tilde C C \delta_q^{\frac{1}{2}} \lambda_q
\end{equation}
is admissible as the upper bound for $|t|$ and gives
\begin{equation}
\snorm{ \bar\Phi_{q, \iota} (t_{q, \iota} +t, \cdot) - \Id}_{n,2 C \lambda_q}    \leq  16 \tilde C \delta_q^{\frac{1}{2}} \frac{1}{\mu_{q+1}}
\end{equation}
which translated into the standard norm is 
\begin{equation}
\norm{ \nabla^n( \bar\Phi_{q, \iota} (t_{q, \iota} +t, \cdot) - \Id)}_\infty   \leq 16 \tilde C
\frac{ n!}{(n+1)^2} (2 C \lambda_q)^{n} \delta_q^{\frac{1}{2}} \frac{1}{\mu_{q+1}} \leq \frac{ n!}{(n+1)^2} (2 C \lambda_q)^{n-1},
\end{equation}
where for the latter inequality suffices
\begin{equation}\label{e:grwl_ass_sep}
    \mu_{q+1} \ge 16 \tilde C  \delta_q^{\frac{1}{2}} \lambda_q.
\end{equation}
Both assumptions \eqref{e:grwl_ass_sep} and \eqref{e:grwl_ass_simple} on relation between $\mu_{q+1}$ and $\delta_{q}^{\sfrac12} \lambda_{q}$ are valid in our scheme, since they are our CFL (Gronwall) condition $\mu_{q+1} \gg \delta_{q}^{\sfrac12} \lambda_{q}$.

Recalling that the global flow map $\bar \Phi_q = \sum_{\iota} \bar \Phi_{q,\iota} \tilde \eta_{q,\iota}$ by its definition \eqref{e:16:constructVF}, and using the properties of $\tilde \eta_{q,\iota}$ (which here, unlike in certain convex integration schemes, does not depend on space variable at all), we therefore have for $n \ge 1$
\[
\begin{cases}
   \norm{ \nabla \bar\Phi_{q}}   &\leq 1+ \frac{1}{4} \le \frac{1}{4} 5 (2 C \lambda_q)^{1-1} , \\
   \norm{ \nabla^n\bar\Phi_{q}}   &\leq \frac{ n!}{(n+1)^2} (2 C \lambda_q)^{n-1}, \quad \text{for } n \ge 2,
\end{cases}
\]
i.e.\  for $n \ge 1$
\begin{equation}\label{e:space_prop_flow}
\snorm{\bar\Phi_{q}}_{n,2C \lambda_q}   \le 5 (2 C \lambda_q)^{-1} 
\end{equation}

\paragraph{(iii. Estimate on the potential)}
Having the estimate \eqref{e:space_prop_flow}, we use it now to estimate space derivatives of 
\[
    \bar \psi_{q+1} = \frac{\delta_{q+1}^{\sfrac12}}{\lambda_{q+1}} 
            H_{12} \big( \lambda_{q+1} \bar \Phi_q(x,t), \mu_{q+1} t \big).
\]
To this end choose in the composition Proposition \ref{prop:compose:analytic}
\[
h (\cdot, t):=  \frac{\delta_{q+1}^{\sfrac12}}{\lambda_{q+1}} 
            H_{12} (\cdot, t), \quad g (x,t) = \lambda_{q+1} \bar \Phi_q(x,t)
\]
From the definition of $H_{12}$ (observe that its space derivatives cost $2 \pi$ each and, importantly, they are much more tame than its time derivatives) and from \eqref{e:space_prop_flow} we have

\begin{equation*}
\snorm{h}_{n,2\pi} \leq  \frac{\delta_{q+1}^{\sfrac12}}{\lambda_{q+1}} 
\qquad \mbox{and} \qquad 
\snorm{g}_{n_1,2C \lambda_q}   \le 5 (2 C \lambda_q)^{-1} \lambda_{q+1} 
\end{equation*}
for all $0 \leq n\leq m$ and $1\leq n_1 \leq m$. Then,   for every  $0 \leq n\leq m$, Proposition \ref{prop:compose:analytic} gives for every $n \ge 0$
\begin{align}
\snorm{    \bar \psi_{q+1}}_{n,R} \leq \frac{\delta_{q+1}^{\sfrac12}}{\lambda_{q+1}}
\end{align}
with 
\[
R
= 2 C \lambda_q (1 + d 5 (2 C \lambda_q)^{-1} \lambda_{q+1}  2 \pi)= 2 C \lambda_q + \lambda_{q+1} 20 \pi \le C \frac{1}{2} \lambda_{q+1} 
\]
where the last inequality holds thanks to (initial) scale separation and choice of $C> 40 \pi$. In other words for $n \ge 0$
\[
\begin{aligned}
\norm{\nabla^n \bar \psi_{q+1}}_\infty &\leq \frac{ n!}{(n+1)^2} \delta_{q+1}^{\sfrac12} (C \frac{1}{2} \lambda_{q+1})^n  \lambda^{-1}_{q+1}\\
&\leq \underbrace{\frac{ n!}{(n+1)^2} 2^{-n}}_{\le \frac{ (n-1)!}{n^2}} \delta_{q+1}^{\sfrac12} (C \lambda_{q+1})^{n-1} C
\end{aligned}
\]
which is the wanted estimate \eqref{e:streamEs_rig1corr}. The estimate on $\bar \Psi_{q+1}$ follows by adding the inductive assumption (valid for $n\ge 2$) to \eqref{e:streamEs_rig1corr}, to obtain for $n\ge 2$
\[
\begin{aligned}
\norm{\nabla^n \bar \Psi_{q+1}}_\infty \le \norm{\nabla^n \bar \Psi_{q}}_\infty + \norm{\nabla^n \bar \psi_{q+1}}_\infty &\le \frac{(n-1)!}{n^2} \tilde C (C{\lambda}_q)^{n-1}\delta_q^{\sfrac12}  + \frac{ (n-1)!}{n^2} \delta_{q+1}^{\sfrac12} (C \lambda_{q+1})^{n-1} C \\
&\le \frac{(n-1)!}{n^2}  \delta_{q+1}^{\sfrac12} \lambda_{q+1}^{n-1} \left(\frac{1}{2}\tilde C C^{n-1}  +C \right) 
\end{aligned}
\]
where the $\frac{1}{2}$ in the second line above is given by the scale separation. Consequently, choosing $\tilde C$ so that $\left(\frac{1}{2}\tilde C C^{n-1}  +C \right) \le \tilde C C^{n-1}$, for which suffices $\tilde C=2$, we have the desired inductive step for $\bar \Psi_{q+1}$, i.e.\ \eqref{e:streamEs_rig1}.
\begin{remark}(lower bound on number of derivatives)\label{rem:psiPsi}
   The estimate on $\nabla^n \bar \Psi_{q+1}$ is restricted for $n \ge 2$, whereas the estimate on $\nabla^n \bar \psi_{q+1}$ holds for  $n \ge 0$. The lower bound $n \ge 2$ follows from the fact that estimate for sum of terms cannot be small ($\bar \Psi_q = \sum_{i \le q} \bar \psi_i$), whereas estimates for a single summand ($\bar \psi_{q+1}$) can be small. (We could also state the inequality without the restriction $n \ge 2$ as $ \| \nabla^n \bar \Psi_q \|_{\infty} 
      \lesssim_n 1 \vee {\lambda}_q^{n-2} (\delta_q^{\sfrac12} \lambda_q)$.)
\end{remark}

\subsubsection{Interlude. Non-commutative product and chain rule}
Recall our multi-index notation ${\balpha} = (\alpha(1), \dots,  \alpha(n_0))$, where each $\alpha(i)$ denotes one derivative (of certain type, in our case: partial spatial derivative or transport derivative). Denoting by ${\bbeta}\prec{\balpha}$ any subsequence of ${\balpha}$ without permutation (i.e.\ the order of entries of ${\balpha}$ is conserved) we have 
\begin{equation}\label{e:product_non_comm}
 D^{\balpha} (fg) = \sum_{{\bbeta}\prec{\balpha}}  (D^{\bbeta} f) (D^{\balpha \setminus \bbeta} g),
 \end{equation}
which is a 'primitive' of the binomial formula, without commutativity properties between derivatives.  

In case we have a smooth function $H : \R \to \R$, the following non-commutative chain rule (or 'Fa\'a di Bruno formula') holds
\begin{equation}\label{e:faa_non_comm_1d}
 D^{\balpha} (H (f)) = \sum_{i=1}^{|\balpha|}  H^{(i)} \prod_{\substack{j=1..i \\ |\balpha_j| \le |\balpha|+1-i \\ \cup_j \balpha_j = \balpha \\ \balpha_j \prec \balpha}} D^{\balpha_j} f,
 \end{equation}
where $H^{(i)}$ is the $i$-th derivative of $H$,  $\cup_j \balpha_j = \balpha$ denotes that the elements of the choice of sequences $\balpha_j$ constitute together, after possible rearrangement, the sequence $\balpha$, and $\balpha_j \prec \balpha$ denotes as before that the sequences  $\balpha_j$ are chosen from $\balpha$ respecting the order of elements of $\balpha$.

In case $H : \R^n \to \R$, the multivariate version of \eqref{e:faa_non_comm_1d} reads
\begin{equation}\label{e:faa_non_comm_nd}
 D^{\balpha} (H (f_1, \dots, f_n)) = \sum_{\zeta=1}^{|\balpha|}  
  \sum_{\substack{k_1=1..n \\ k_2=1..n\\ \dots \\ k_\zeta=1..n}}
 \partial_{k_\zeta} \dots \partial_{k_1} H \prod_{\substack{j=1..\zeta \\ |\balpha_j| \le |\balpha|+1-i \\ \cup_j \balpha_j = \balpha \\ \balpha_j \prec \balpha}} D^{\balpha_j} f_{k_j},
 \end{equation}
We will usually specify these formulas to simpler cases, when applied.

\subsubsection{Propagating mixed derivatives. I. Flowmap derivatives}\label{ssec:mixed_flowmap}
We begin proving estimate \eqref{e:streamEs_rig2}. Now, the total number of derivatives is bounded from above by $Y$, consequently the combinatorial constants resulting from composition and product estimates are unimportant: as usual in similar iterative schemes, they will be uniformly stabilised at the step $q \to q+1$, when applied.

Crucially, there is no loss of total number of derivatives in the step $q \to q+1$, recall Section \ref{sssec:ir_ests} (ii).

We already know, thanks to \eqref{e:space_prop_flow}, that for each local flow $\Phi_{q, \iota} $ and for the global flow $ \bar \Phi_{q}$
\begin{equation}\label{e:flow_map_est_ini}
\begin{aligned}
 |\nabla^i  \bar \Phi_{q, \iota}|_{supp\; \tilde \eta_{q, \iota} }\quad+ \quad |\nabla^i  \bar \Phi_{q}|
&\lesssim 1, \qquad i=0,1 \\
 |\nabla^i  \bar \Phi_{q, \iota}|_{supp\; \tilde \eta_{q, \iota} }\quad+ \quad |\nabla^i  \bar \Phi_{q}|
&\lesssim  {\lambda}^{i-1}_q , \qquad  i = 1, \dots, L
\end{aligned}
\end{equation}

Let us analyse $D_q^{\balpha}  \bar \Phi_{q, \iota}$ with respect to how many advective derivatives there are in $D_q^{\balpha}$. Recall \eqref{e:tr_hofHi} and consider the case of one advective derivative present in $D_q^{\balpha}$, i.e.\ $D_q^{\balpha} = \nabla^m D_q \nabla^n \bar \Phi_{q, \iota}$. We have
\begin{equation}
\begin{aligned}
D_q \nabla^n \bar \Phi_{q, \iota} &\sim_n \sum_{i=1}^{n} \underbrace{\nabla^{n+1-i}  \nabla^T \bar \Psi_q}_{\text{by assumption, costs } {\lambda}_q^{n+1-i}\delta_q^{\sfrac12}; }  \underbrace{\nabla^i \bar \Phi_{q, \iota}}_{\quad\text{costs } {\lambda}^{i-1}_q \text{ on }supp\; \tilde \eta_{q, \iota} } 
\end{aligned}
\end{equation}
Since for both terms the cost of one space derivative is ${\lambda}_q$, $\nabla^m D_q \nabla^n \bar \Phi_{q, \iota}$ costs \[{\lambda}^m_q {\lambda}_q^{n+1-i}\delta_q^{\sfrac12} {\lambda}^{i-1}_q = \delta_q^{\sfrac12} {\lambda}^{m+n}_q.\] In other words, we have
\begin{equation}\label{e:flowmap_single_adv_est}
 |D_q^{\balpha}  \bar \Phi_{q, \iota}|_{supp\; \tilde \eta_{q, \iota} } 
\lesssim {\lambda}_q^{|\balpha|_x-1} (\delta_q^{\sfrac12} \lambda_q)^{|\balpha|_t}, \qquad |\balpha|_t = 1, \quad |\balpha|_x = |\balpha|-1 \;(=m+n)
\end{equation}
 Remembering that the transport derivatives of the global flow involve also $\tilde \eta_{q, \iota}$, whose (time) derivatives cost $\mu_{q+1} > \delta_q^{\sfrac12} \lambda_q$ each, we obtain
\begin{equation}
 |D_q^{\balpha}   \bar \Phi_q| 
\lesssim {\lambda}_q^{|\balpha|_x-1} \mu_{q+1}^{|\balpha|_t} \qquad |\balpha|_t = 1, \quad |\balpha|_x = |\balpha|-1 
\end{equation}

Further
\begin{equation}
\begin{aligned}
\nabla^m D_q \nabla^n \bar \Phi_{q, \iota} &\sim \sum_{j=0}^{m} \sum_{i=1}^{n}\underbrace{\nabla^{n+1-i+m-j}  \nabla^T \bar \Psi_q}_{\text{costs } {\lambda}_q^{n+1-i+m-j}\delta_q^{\sfrac12}     }  \underbrace{\nabla^{i+j} \bar \Phi_{q, \iota}}_{\text{costs } {\lambda}^{i-1+j}_q \text{ on }supp\; \tilde \eta_{q, \iota} }.
\end{aligned}
\end{equation}
Applying above $D_q$, in view of the already established case of one advective derivative acting on the flowmap, i.e.\ \eqref{e:flowmap_single_adv_est}, we see that $\nabla^k D_q \nabla^m D_q \nabla^n \bar \Phi_{q, \iota}$ costs ${\lambda}^{m+n+k-1}_q (\delta_q^{\sfrac12} {\lambda}_q)^2$ on the supports of $\tilde \eta_{q, \iota}$, i.e.
\[
|D_q^{\balpha}  \bar \Phi_{q, \iota}|_{supp\; \tilde \eta_{q, \iota} } 
\lesssim {\lambda}_q^{|\balpha|_x-1} (\delta_q^{\sfrac12} \lambda_q)^{|\balpha|_t}, \qquad |\balpha|_t = 2, \quad |\balpha|_x = |\balpha|-2 \;(=m+n+k)
\]
Taking into account as before the derivatives of the cutoffs appearing for the global flowmap, we have
\begin{equation} 
 |D_q^{\balpha}   \bar \Phi_q| 
\lesssim {\lambda}_q^{|\balpha|_x-1} \mu_{q+1}^{|\balpha|_t} \qquad |\balpha|_t = 2, \quad |\balpha|_x = |\balpha|-2. 
\end{equation}
Repeating this reasoning for three, etc, up to $|\balpha|-1$ advective derivatives, we reach 
\begin{equation}\label{e:flow_est_interm}
 |D_q^{\balpha}   \bar \Phi_q| 
\lesssim {\lambda}_q^{|\balpha|_x-1} \mu_{q+1}^{|\balpha|_t} \quad |\balpha|_x \ge 1.
\end{equation}
(The case of $|\balpha|$-many advective derivatives, i.e.\ solely advective derivatives present in $D_q^{\balpha}   \bar \Phi_q$ yields nonzero quantities only from cutoffs, since $ \Phi_{q, \iota}=0$ is a flowmap. Consequently $\mu_{q+1}^{|\balpha|}$ appears as r.h.s.\ in \eqref{e:flow_est_interm}. We leave this case out, because it is not needed in what follows, and not consistent with \eqref{e:flow_map_est_ini}).

Observe that the estimate \eqref{e:flow_est_interm} is consistent with \eqref{e:flow_map_est_ini}, so we can write 
\begin{equation}\label{e:flow_est_interm2}
 |D_q^{\balpha}   \bar \Phi_q| 
\lesssim {\lambda}_q^{|\balpha|_x-1} \mu_{q+1}^{|\balpha|_t} \quad \text{ for any } \quad |\balpha| \le Y, \; |\balpha|_x \ge 1.
\end{equation}

\begin{remark}
    One can argue equivalently using commutator estimates. Most of them, available 'from the shelf', seem to induce however some loss of derivative.
\end{remark}

\subsubsection{Propagating mixed derivatives. II. Estimates for $ D_q^{\balpha} \psi_{q+1}$.}\label{ssec:pmdbarpsi1}
Recall
\[
\bar \psi_{q+1}  =    \frac{\delta_{q+1}^{\sfrac12}}{\lambda_{q+1}}  H_{12} \big( \lambda_{q+1} \bar \Phi_q(x,t), \mu_{q+1} t \big)
\]
First, we make the following structural observation
\[
\begin{aligned}
&D_{q} \bar \psi_{q+1} \\
&=  \frac{\delta_{q+1}^{\sfrac12}}{\lambda_{q+1}} 
          D_{{q}}   H_{12} \big( \lambda_{q+1} \bar \Phi_q(x,t), \mu_{q+1} t \big) \\
&= \frac{\delta_{q+1}^{\sfrac12}}{\lambda_{q+1}} \mu_{q+1} (\partial_\tau    H_{12}) \big( \lambda_{q+1} \bar \Phi_q(x,t), \mu_{q+1} t \big) + \delta_{q+1}^{\sfrac12} ( D_{q} \bar \Phi_q) \cdot (\nabla_\xi H_{12}) \big( \lambda_{q+1} \bar \Phi_q(x,t), \mu_{q+1} t \big).
\end{aligned} 
\]
Importantly, the last term above vanishes, because
\begin{equation}\label{e:trans:cancel}
(D_{{q}} \bar \Phi_q) \cdot (\nabla_\xi H_{12}) ()=  \left(   \sum_\iota \tilde \eta_{q,\iota} D_q  \bar \Phi_{q,\iota} +
      \sum_\iota \partial_t \tilde \eta_{q,\iota}  \bar \Phi_{q,\iota} \right) \cdot (\nabla_\xi H_{12}) () = 0,      
\end{equation}
indeed $D_q  \bar \Phi_{q,\iota} =0$ and the (time) supports of derivatives of partitions of unity $\tilde \eta$ are disjoint from the (time) supports of $H_{12}$, since $H_{12}$ involves in its definition appropriate  cutoffs $\eta_1, \eta_2$. Therefore we have

\begin{equation}\label{e:trans_Dcorr}
D_{{q}} \bar \psi_{q+1} =\frac{\delta_{q+1}^{\sfrac12}}{\lambda_{q+1}} \mu_{q+1} (\partial_\tau    H_{12}) \big( \lambda_{q+1} \bar \Phi_q(x,t), \mu_{q+1} t \big).
\end{equation}
Observe that the higher advective derivatives of $\bar \psi_{q+1}$ behave identically, thus 
\begin{equation}\label{e:trans_Dcorr2}
D^n_{{q}} \bar \psi_{q+1} =\frac{\delta_{q+1}^{\sfrac12}}{\lambda_{q+1}} \mu^n_{q+1} (\partial^n_\tau    H_{12}) \big( \lambda_{q+1} \bar \Phi_q(x,t), \mu_{q+1} t \big).
\end{equation}

Without the cancellation \eqref{e:trans:cancel}, and the resulting \eqref{e:trans_Dcorr}, \eqref{e:trans_Dcorr2}, the obtained estimates would prove insufficient.

For the general case of $D_q^{\balpha}$ we may use directly the multivariate chain rule \eqref{e:faa_non_comm_nd} specified to three dimensions, with choices:
\[
H (\cdot,\cdot,\cdot):= H_{12} (\lambda_{q+1} \cdot,\lambda_{q+1} \cdot, \cdot), \quad (f_1, f_2) :=\bar \Phi_q(x,t), \quad  f_3 := \mu_{q+1} t
\]
It is more tractable, however, to make several simplifying observations first. Recall we are now interested only in derivatives up to order $Y$. Consequently the size of any $\partial_{k_\zeta} \dots \partial_{k_1} H_{12}$ is bounded above by a uniform constant, and changes of size between different derivatives of $H_{12}$ are thus negligible. Moreover   \\
(i) While taking a (space or advective) derivative of 
\begin{equation}\label{e:simplH}
H_{12}( \lambda_{q+1} \bar \Phi_q(x,t), \mu_{q+1} t), \text{\quad  or of any \quad} (\partial_{k_\zeta} \dots \partial_{k_1} H)( \lambda_{q+1} \bar \Phi_q(x,t), \mu_{q+1} t )
\end{equation}
the space derivative acts only on the inner functions $\bar \Phi_q(x,t)$ and the advective derivative acts only on $\mu_{q+1} t$. This is the already observed cancellation \eqref{e:trans:cancel}.
\\
(ii) Further, the cost of an advective  derivative acting on \eqref{e:simplH} is $\mu_{q+1}$ (in view of (i), it is now pure time derivative), which is consistent with the cost of the advective derivative of $\bar \Phi_q(x,t)$, cf \eqref{e:flow_est_interm2}.\\
(iii) In view of the chain rule \eqref{e:faa_non_comm_nd}, 
\begin{equation}\label{e:simplH2}
D_q^{\balpha'} \left(H ( \lambda_{q+1} \bar \Phi_q(x,t), \mu_{q+1} t ))\right)
\end{equation}
is a sum of products of terms of type \eqref{e:simplH} and advective derivatives of $\bar \Phi_q(x,t)$. Consequently, via (ii), an advective derivative of \eqref{e:simplH2} costs $\mu_{q+1}$.\\
(iv) Each space derivative applied to \eqref{e:simplH2} results in the chain rule present only via space derivative of $\lambda_{q+1} \bar \Phi_q(x,t)$ (and costs up to $\lambda_{q+1}$, which will be clear in the following computation). \\

The above simplifications show that  $D_q^{\balpha} \bar \psi_{q+1}$ can be modelled with the following combination of (simpler) one-dimensional chain rules (compare \eqref{e:faa_non_comm_1d})
\begin{equation}\label{e:simple_chain_applied}
\begin{aligned}
| D^{\balpha} \bar \psi_{q+1} (x,t)|  &= 
\frac{\delta_{q+1}^{\sfrac12}}{\lambda_{q+1}} | D^{\balpha}   H_{12} (\lambda_{q+1} \bar \Phi_q(x,t), \mu_{q+1} t)|\\
&\lesssim_Y \frac{\delta_{q+1}^{\sfrac12}}{\lambda_{q+1}}  \sum_{\substack{i=0..|\balpha|_x,\\ k=0..|\balpha|_t,\\ i+j\ge 1}} \lambda^{i}_{q+1} \prod_{\substack{j=1..i+k \\ |\balpha_j| \le  |\balpha|+1-(i+k) \\ \cup_j \balpha_j = \balpha \\ \balpha_j \prec \balpha}} \|D^{\balpha_j} \bar \Phi_q \|_\infty\\
\text{use \eqref{e:flow_est_interm2}}\quad&\lesssim_Y \frac{\delta_{q+1}^{\sfrac12}}{\lambda_{q+1}}  \sum_{\substack{i=0..|\balpha|_x,\\ k=0..|\balpha|_t,\\ i+j\ge 1}} \lambda^{i}_{q+1} \prod_{\substack{j=1..i+k \\ |\balpha_j| \le  |\balpha|+1-(i+k) \\ \cup_j \balpha_j = \balpha \\ \balpha_j \prec \balpha}}   {\lambda}_q^{|\balpha_j|_x-1} \mu_{q+1}^{|\balpha_j|_t}\\
&\lesssim_Y \frac{\delta_{q+1}^{\sfrac12}}{\lambda_{q+1}}  \mu^{|\balpha|_t}_{q+1} \sum_{i=0}^{|\balpha|_x}  
\lambda^{i}_{q+1} \prod_{\substack{j=1..i \\ |\balpha_j| \le  |\balpha_x|+1-i \\ \cup_j \balpha_j = \balpha_x \\ \balpha_j \prec \balpha}} {\lambda}_q^{|\balpha_j|_x-1}\\
&\lesssim_Y \frac{\delta_{q+1}^{\sfrac12}}{\lambda_{q+1}}  \mu^{|\balpha|_t}_{q+1}   
\lambda^{|\balpha|_x}_{q+1},
\end{aligned}
\end{equation}
where the last inequality follows from the fact that in the sum (in the penultimate row) the last summand is the largest (since $\lambda_{q+1} \gg \lambda_{q}$, up to combinatorial coefficients, taken over by $\lesssim_Y$).
We obtained
\begin{equation}\label{e:slow_corr_est}
| D_q^{\balpha} \bar \psi_{q+1} (x,t)|
\lesssim_Y \delta_{q+1}^{\sfrac12}  \mu^{|\balpha|_t}_{q+1}   
\lambda^{|\balpha|_x-1}_{q+1}
    \end{equation}
for any $|\balpha| \le Y$. Observe that even though the used estimate \eqref{e:flow_est_interm2} has restriction $|\balpha|_x \ge 1$, it does not appear in \eqref{e:slow_corr_est}. Indeed, in \eqref{e:simple_chain_applied} application of chain rule does not need \eqref{e:flow_est_interm2} with $|\balpha|_x =0$, and there is eventually $\lambda^{-1}_{q+1}$ appearing thanks to $\sfrac{\delta_{q+1}^{\sfrac12}}{\lambda_{q+1}}$.
\subsubsection{Propagating mixed derivatives. III. Estimates for $ D_{q+1}^{\balpha} \psi_{q+1}$.}\label{ssec:pmdbarpsi2}

In order to deal with transport derivatives with respect to $u_{q+1} = u_q + \nabla^\perp \bar \psi_{q+1}$, let us shift the flow as follows
\begin{equation}\label{e:flowshift_0}
 D_{q+1} f =  D_q f + \nabla^\perp \bar \psi_{q+1} \cdot \nabla  f.
\end{equation}
We will use the estimate 
\begin{equation}\label{e:slow_corr_est_appl}
| D_q^{\balpha} \bar \psi_{q+1} (x,t)|
\lesssim_Y \delta_{q+1}^{\sfrac12}  ( \delta_{q+1}^{\sfrac12} \lambda_{q+1})^{|\balpha|_t}   
\lambda^{|\balpha|_x-1}_{q+1} \quad \text{ for any } |\balpha| \le Y
    \end{equation}
which follows from 
\eqref{e:slow_corr_est} and \eqref{e:streamEs_rig1corr}, with $\mu_{q+1} \le \delta_{q+1}^{\sfrac12} \lambda_{q+1}$, see Remark \ref{r:auxiliaryScalesRela}.
In particular, choosing $f=\bar \psi_{q+1}$ in \eqref{e:flowshift_0}, we have
\begin{equation}\label{e:trans_02hspace}
| D_{q+1} \bar \psi_{q+1}| \lesssim \delta_{q+1} = \frac{\delta_{q+1}^{\sfrac12}}{\lambda_{q+1}} (\delta_{q+1}^{\sfrac12} \lambda_{q+1})
\end{equation}
which is the wanted \eqref{e:streamEs_rig2} for $|\balpha|_x =0$ and $|\balpha|_t =1$. One more advective derivative on $\bar \psi_{q+1}$ yields:
\[
\begin{aligned}
&D^2_{{q+1}} \bar \psi_{q+1} 
\\
&=  ( D_q  + \nabla^\perp \bar \psi_{q+1} \cdot \nabla ) ( D_q \psi_{q+1} + \nabla^\perp \bar \psi_{q+1} \cdot \nabla \bar \psi_{q+1})\\
&=    D^2_{q} \psi_{q+1} +  (D_q  \nabla^\perp \bar \psi_{q+1}) \cdot \nabla \bar \psi_{q+1} + (\nabla^\perp \bar \psi_{q+1}) \cdot (D_q \nabla \bar \psi_{q+1}) +   \\& \; + \nabla^\perp \bar \psi_{q+1} \cdot \nabla   D_q \psi_{q+1} + \nabla^\perp \bar \psi_{q+1} \cdot \nabla \nabla^\perp \bar \psi_{q+1} \cdot \nabla \bar \psi_{q+1} + \nabla^\perp \bar \psi_{q+1} \cdot (\nabla^\perp \bar \psi_{q+1} \cdot \nabla^2 \bar \psi_{q+1})
\end{aligned}
\]
(where $\cdot$ needs to be interpreted appropriately). Using \eqref{e:slow_corr_est_appl} we consequently have
\[
\begin{aligned}
|D^2_{q+1} \bar \psi_{q+1}| \lesssim &   \frac{\delta_{q+1}^{\sfrac12}}{\lambda_{q+1}} \mu^2_{q+1}+ \delta_{q+1}^{\sfrac12} \mu_{q+1} \delta_{q+1}^{\sfrac12} + \delta_{q+1}^{\sfrac12}  \delta_{q+1}^{\sfrac12}  \lambda_{q+1} \delta_{q+1}^{\sfrac12} 
\end{aligned}
\]
which in view of the recently invoked  $\mu_{q+1} \le \delta_{q+1}^{\sfrac12} \lambda_{q+1}$ gives 
\[
|D^2_{q+1} \bar \psi_{q+1}| \lesssim  \frac{\delta_{q+1}^{\sfrac12}}{\lambda_{q+1}} (\delta_{q+1}^{\sfrac12} \lambda_{q+1})^2.
\]
We computed $D^2_{q+1}$ to clearly make the point the the expected cost of higher-order advective derivatives is not $\mu_{q+1}$ as one could expect from \eqref{e:slow_corr_est}, but $\delta_{q+1}^{\sfrac12} \lambda_{q+1}$, due to $\nabla^\perp \bar \psi_{q+1} \cdot \nabla$.

The general estimate for $D_{q+1}^{\balpha} \bar \psi_{q+1}$ follows the same lines. Indeed, let us observe first that any block of advective derivatives present in $D_{q+1}^{\balpha}$ is a block of $(D_q + \nabla^\perp \bar \psi_{q+1} \cdot \nabla)$. Application of $(D_q + \nabla^\perp \bar \psi_{q+1} \cdot \nabla)$ to any function $f$ with magnitude $M$, whose advective derivative cost is (no bigger than) $\delta_{q+1}^{\sfrac12}\lambda_{q+1}$ and space derivative cost is $\lambda_{q+1}$ results in a function whose  magnitude is 
\[
 M (\delta_{q+1}^{\sfrac12}\lambda_{q+1} + \|\nabla \bar \psi_{q+1}\|\lambda_{q+1} ) = M (\delta_{q+1}^{\sfrac12}\lambda_{q+1}  + \delta_{q+1}^{\sfrac12}\lambda_{q+1} ) \le 2 M \delta_{q+1}^{\sfrac12}\lambda_{q+1} 
\]
and whose derivative cost is the same as that of $f$. The constant $2$ above is consumed by $\lesssim_Y$. This observation can be used to any $D_{q+1}^{\balpha'} \bar \psi_{q+1}$ with $\balpha' < \balpha$, yielding  
\[
      \| D^{\balpha}_{q+1} \bar \psi_{q+1} \|_{\infty} 
      \lesssim_Y \frac{\delta_{q+1}^{\sfrac12}}{\lambda_{q+1}} (\delta_{q+1}^{\sfrac12} \lambda_{q+1})^{|\balpha|_t} \lambda_{q+1}^{|\balpha|_x} \]
which is the wanted 
      \eqref{e:streamEs_rig2_psi}. Taking into account that $\bar \Psi_{q+1}= \bar \Psi_{q}+  \bar \psi_{q+1} $, compare also Remark \ref{rem:psiPsi}, we have the remaining wanted \eqref{e:streamEs_rig2}.
This ends the proof of Lemma \ref{l:streamEst}.

 \subsection{Proof of Lemma \ref{l:streamEst_rig_nobar} and its Corollaries}\label{sec:streamEst_rig_nobar_pf} 
\subsubsection{Initial remarks}
Proof of Lemma \ref{l:streamEst_rig_nobar} uses in many places reasonings very similar to proof of Lemma \ref{l:streamEst}, which we invoke without repeating fully. In particular, we will use
\begin{remark}[Remark on commutativity]\label{rem:comm}
As observed in the proof of Lemma \ref{l:streamEst}, even though the space vs advective derivatives did not commute, since the estimates for each type of derivative were consistent regardless of the order of derivation, the resulting inequalities were as if the derivatives commuted. 
\end{remark}

There is, however, one important difference between the currently wanted estimates for $\Psi_{q+1}$ and the already proven in Lemma \ref{l:streamEst} estimates for $\bar \Psi_{q+1}$. Namely, the estimates on $\bar \Psi_{q+1}$ follow from an inductive procedure, so in essence from estimates on $\bar \Psi_{q}$. Therefore we were careful concerning derivative loss, since otherwise induction over $q\in \N$ fails (looses eventually all derivatives). However, the estimates for $\Psi_{q+1}$ involve essentially the estimate for $\bar \Psi_{q+1}$ and no induction is used. (So 'loss of derivative' $Y \to Y-N_*-1$ is not related to inductive step $\Psi_q \to \Psi_{q+1}$, but to a 'same step comparison' between $\Psi_{q+1}$ and  $\bar \Psi_{q+1}$).

\subsubsection{Estimates on flowmaps}

 Before proceeding to the estimates for $ \psi_{q+1}$, let us obtain estimates for the explicit flow $\Phi_{q}$ and for the difference between the exact and the explicit flows, i.e.\ $\bar \Phi_{q} - \Phi_{q} = \sum_{\iota} (\bar\Phi_{q, \iota} - \Phi_{q, \iota}) \tilde \eta_{q, \iota}$. 
In view of the formula directly preceding \eqref{e:app_flow_loc}
\[
  D_q\Phi_{q,\iota} = \frac{1}{N_*!} (( t_{q,\iota} - t )^{N_*} D_q^{N_*} u_q) = 
\frac{1}{ \mu^{N_*}_{q+1}} D_q^{N_*} \nabla^\perp \bar \Psi_q \frac{1}{{N_*}!} ( \iota - \mu_{q+1} t)^{N_*}=:G.
\]
Importantly, the velocity used above is $u_q =\nabla^\perp \bar \Psi_q$, whose estimates we know. In particular, in view of \eqref{e:streamEs_rig2}, for any $|\balpha|\le Y - N_*$
\begin{equation}\label{e:flows_comparison_rhs}
\begin{aligned}
|D^{\balpha}_q G|_{\supp \tilde \eta_{q, \iota}} \lesssim \frac{1}{ \mu^{N_*}_{q+1}} \frac{\delta_{q}^{\sfrac12}}{\lambda_{q}} (\delta_{q}^{\sfrac12} \lambda_{q})^{N_*+|\balpha|_t} \lambda_{q}^{|\balpha|_x+1} &=   (\frac{\delta_{q}^{\sfrac12} \lambda_{q}}{ \mu_{q+1}})^{N_*} \delta_{q}^{\sfrac12}(\delta_{q}^{\sfrac12} \lambda_{q})^{|\balpha|_t} \lambda_{q}^{|\balpha|_x} \\
&\le \delta_{q}^{\sfrac12} \lambda^{-\gamma N_*}_{q+1} (\delta_{q}^{\sfrac12} \lambda_{q})^{|\balpha|_t} \lambda_{q}^{|\balpha|_x},
\end{aligned}
\end{equation}
with the last uses the $\gamma$-gap in $\mu_{q+1} > \delta_{q}^{\sfrac12} \lambda_{q}$. 
So $\bar\Phi_{q, \iota} - \Phi_{q, \iota}$ solves 
$D_q(\bar\Phi_{q, \iota} - \Phi_{q, \iota}) =g$
on the support of $\tilde \eta_{q, \iota}$ with zero initial datum.
The solution of $\partial_t f + u \cdot \nabla f = g$ with zero initial datum enjoys
\begin{equation*}
    \|f(\cdot, t)\|_0 \leq \int_{t_0}^t \|g(\cdot, \tau)\|_0 d \tau, 
\end{equation*}
and, up to its Gronwall time, \begin{equation*}
|\nabla^n f(t)| \lesssim \int_{t_0}^t \big( |g(\cdot, \tau)|_{n} + (t-\tau) |u|_{n} |g(\cdot, \tau)|_1\big) d\tau \quad \text{ for } n \ge 1
\end{equation*} see for instance \cite{GiriRadu}, Proposition B.1. Consequently, using  \eqref{e:flows_comparison_rhs} and  \eqref{e:streamEs_rig2} we have for any $n = 0, \dots, Y - N_*$
\begin{equation*}
\begin{aligned}
|\nabla^n (\bar\Phi_{q, \iota} - \Phi_{q, \iota})| &\lesssim \int_{t_0}^t \big( \delta_{q}^{\sfrac12} \lambda^{-\gamma N_*}_{q+1} \lambda_{q}^{n}  + (t-\tau) \delta_{q}^{\sfrac12} \lambda_{q}^{n} \delta_{q}^{\sfrac12} \lambda^{-\gamma N_*}_{q+1} \lambda_{q} \big) d\tau\\
 &\le (\mu^{-1}_{q} \delta_{q}^{\sfrac12} + \mu^{-2}_{q} \delta_{q}\lambda_{q}) \lambda_{q}^{n}  \lambda^{-\gamma N_*}_{q+1} \\
 &\le \lambda^{-\gamma N_*}_{q+1} \lambda_{q}^{n},
\end{aligned}
\end{equation*}
where for the last inequality we used lossy inequality $(\mu^{-1}_{q} \delta_{q}^{\sfrac12} + \mu^{-2}_{q} \delta_{q}\lambda_{q}) \le 1$ for the scales. This means that we have the following counterpart of the estimate
\eqref{e:flow_map_est_ini}
\begin{equation}\label{e:flow_map_est_ini_diff}
 |\nabla^i  (\bar\Phi_{q, \iota} - \Phi_{q, \iota})|_{supp\; \tilde \eta_{q, \iota} }, \quad |\nabla^i  (\bar\Phi_{q, \iota} - \Phi_{q, \iota})|
\lesssim  \lambda^{-\gamma N_*}_{q+1} \lambda_{q}^{i} \le \lambda^{-\gamma N_*+1}_{q+1} \lambda_{q}^{i-1} \qquad  i = 0, \dots, L-N_*
\end{equation}
Reasoning now analogously to section \ref{ssec:mixed_flowmap} and observing that the cost of derivatives of the forcing term $G$ is commensurate, we have for $|\balpha| = 0, \dots, Y - N_*$
\begin{equation}\label{e:flow_map_est_diff}
 |D_q^{\balpha}   (\bar \Phi_q - \Phi_q)| 
\lesssim {\lambda}_q^{|\balpha|_x-1} \mu_{q+1}^{|\balpha|_t} {\lambda}^{-\gamma N_*+1}_{q+1}.
\end{equation}
Triangle inequality and \eqref{e:flow_est_interm2}, \eqref{e:flow_map_est_diff} yield for $|\balpha| = 1, \dots, Y - N_*$
\begin{equation}\label{e:app_flow_map_est}
 |D_q^{\balpha}  \Phi_q| 
\lesssim \frac{1}{{\lambda}_q} {\lambda}_q^{|\balpha|_x} \mu_{q+1}^{|\balpha|_t}. 
\end{equation}

\subsubsection{Estimates for $ D^{\balpha}_{q+1} \psi_{q+1}$}

Recall that by definition \eqref{e:psiQ+1}
\[ \psi_{q+1}(x,t) = \frac{\delta_{q+1}^{\sfrac12}}{\lambda_{q+1}} 
            (\det \nabla \Phi_q) H_{12} \big( \lambda_{q+1} \Phi_q(x,t), \mu_{q+1} t \big).     \] 
The term $H_{12} ( \lambda_{q+1} \Phi_q(x,t), \mu_{q+1} t )$ above behaves as the term $H_{12} ( \lambda_{q+1} \bar \Phi_q(x,t), \mu_{q+1} t )$ appearing in $\bar\psi_{q+1}(x,t)$, already analysed. Indeed, $ \bar \Phi_q(x,t)$ and  $\Phi_q(x,t)$ have identical estimates for any $|\balpha| \le Y - N_*$ (see, respectively, \eqref{e:flow_est_interm2} and \eqref{e:app_flow_map_est}) and the involved cutoff functions are the same. Consequently, along the proof of estimates for $\bar\psi_{q+1}(x,t)$, see sections \ref{ssec:pmdbarpsi1}, \ref{ssec:pmdbarpsi2}, we have 
\begin{equation}\label{e:to_psi_first}
      \left\| D^{\balpha}_{q+1}  \frac{\delta_{q+1}^{\sfrac12}}{\lambda_{q+1}} 
        H_{12} \big( \lambda_{q+1} \Phi_q(x,t), \mu_{q+1} t \big)    \right\|_{\infty} 
      \lesssim_Y \frac{\delta_{q+1}^{\sfrac12}}{\lambda_{q+1}} (\delta_{q+1}^{\sfrac12} \lambda_{q+1})^{|\balpha|_t} \lambda_{q+1}^{|\balpha|_x} 
 \quad  \text{ for any }  |\balpha| \le Y - N_*.
    \end{equation}
The estimate \eqref{e:app_flow_map_est}, taking into account Remark \ref{rem:comm} and quadraticity of $\det$, gives 
\begin{equation}\label{e:slowdetflow}
|D_q^{\bomega}  (\det \nabla\Phi_q)| 
\lesssim {\lambda}_q^{|\bomega|_x} \mu_{q+1}^{|\bomega|_t} 
\quad \text{ for any } |\bomega| \le Y - N_*-1 .
\end{equation}
In order to provide the wanted estimate for $D^{\balpha}_{q+1} \psi_{q+1}$, we need estimates for $D_{q+1}^{\bomega}  (\det \nabla\Phi_q)$, i.e.\ 'shift the flow' in the estimate \eqref{e:slowdetflow}.  Recall that  $D_{q+1} = (D_q + \nabla^\perp \bar \psi_{q+1} \cdot \nabla)$. We have full knowledge of estimates of $\bar \psi_{q+1}$, so repeating the reasoning of section \ref{ssec:pmdbarpsi2}, in particular its second part, now   for function  $\det \nabla\Phi_q$, we obtain for any $|\bomega| \le Y - N_*-1$.  
\begin{equation}\label{e:app_flow_map_est_fast}
|D_{q+1}^{\bomega}  (\det \nabla\Phi_q)| 
\lesssim \lambda_{q+1}^{|\bomega|_x} (\delta_{q+1}^{\sfrac12} \lambda_{q+1})^{|\bomega|_t}. 
\end{equation}
Putting together estimates \eqref{e:app_flow_map_est_fast} and \eqref{e:to_psi_first} for the two ingredients of  $\psi_{q+1}$ we reach via the product rule/Remark \ref{rem:comm} the wanted \eqref{e:streamEs_rig2_nobar_psi} for $D_{q+1}^{\balpha}$ derivatives, and consequently \eqref{e:streamEs_rig2_nobar} for $D_{q+1}^{\balpha}$ derivatives, with the upper bound on number of derivatives dictated by the bound in \eqref{e:slowdetflow}.

Let us note in passing, for further reference, that analogously to the way we obtained \eqref{e:app_flow_map_est_fast} from \eqref{e:slowdetflow}, we see that \eqref{e:slowdetflow} implies for any $|\bomega| \le Y - N_*-1$ (and for each component $\nabla\Phi^1_q$, $\nabla\Phi^2_q$, which we do not distinguish below)
\begin{equation}\label{e:app_flow_map_est_fast2}
|D_{q+1}^{\bomega}  ( \nabla\Phi_q)| 
\lesssim \lambda_{q+1}^{|\bomega|_x} (\delta_{q+1}^{\sfrac12} \lambda_{q+1})^{|\bomega|_t}
\end{equation}
and that \[
 |D_q^{\balpha}  \nabla (\bar \Phi_q - \Phi_q)| 
\lesssim {\lambda}_q^{|\balpha|_x} \mu_{q+1}^{|\balpha|_t} {\lambda}^{-\gamma N_*+1}_{q+1}
\]
given by \eqref{e:flow_map_est_diff} implies
\begin{equation}\label{e:app_flow_map_est_fast21}
|D_{q+1}^{\bomega}  ( \bar \nabla\Phi_q -\nabla\Phi_q)| 
\lesssim  \lambda^{-\gamma N_*+1}_{q+1} \lambda_{q+1}^{|\bomega|_x} (\delta_{q+1}^{\sfrac12} \lambda_{q+1})^{|\bomega|_t}. 
\end{equation}
\subsubsection{Estimates for $ D^{\balpha}_{q+1} (\bar \psi_{q+1}- \psi_{q+1})$} In view of definitions \eqref{e:barCumuPsiQ+1}, \eqref{e:CumuPsiQ+1}, and \eqref{e:barpsiQ+1}, \eqref{e:psiQ+1}, we have
$\bar \Psi_{q+1} - \Psi_{q+1} = \bar \psi_{q+1}  -  \psi_{q+1}$ and
\begin{equation}\label{e:flow_diff_start}
\begin{aligned}
& \bar\psi_{q+1}  -  \psi_{q+1}\\  &= (1-  \det \nabla \Phi_q) \bar \psi_{q+1} &:=I\\
&+ \frac{\delta_{q+1}^{\sfrac12}}{\lambda_{q+1}}
            (\det \nabla \Phi_q) (H_{12} \left(  \lambda_{q+1} \Phi_q(x,t), \mu_{q+1} t ) - H_{12} (\lambda_{q+1} \bar \Phi_q(x,t), \mu_{q+1} t )
            \right) 
            &=:II
\end{aligned}
\end{equation}
The term I in \eqref{e:flow_diff_start} involves $1-  \det \nabla \Phi_q = \det \nabla \bar \Phi_q-  \det \nabla \Phi_q$
valid in view of the incompressibility of the flow; therefore quadracity of $\det$ gives $1-  \det \nabla \Phi_q \sim (\nabla \bar \Phi_q - \nabla\Phi_q) (\nabla \bar \Phi_q + \nabla\Phi_q)$. Or, performing the computations, we have
\[
\begin{aligned}
1-  \det \nabla \Phi_q = & \det \nabla \bar \Phi_q-  \det \nabla \Phi_q \\
=  & \partial_{x_1}\bar\Phi^1_q  \partial_{x_2}\bar\Phi^2_q -  \partial_{x_1}\bar\Phi^2_q \partial_{x_2}\bar\Phi^1_q - ( \partial_{x_1}\Phi^1_q \partial_{x_2}\Phi^2_q - \partial_{x_1}\Phi^2_q \partial_{x_2}\Phi^1_q ) \\
=& (\partial_{x_1}\bar\Phi^1_q -\partial_{x_1}\Phi^1_q)  \partial_{x_2}\bar\Phi^2_q +  \partial_{x_1} \Phi^1_q (\partial_{x_2}\bar\Phi^2_q -\partial_{x_2}\Phi^2_q) \\
&+(\partial_{x_1}\Phi^2_q -\partial_{x_1}\bar\Phi^2_q)  \partial_{x_2}\Phi^1_q +  \partial_{x_1} \bar\Phi^1_q (\partial_{x_2}\Phi^1_q -\partial_{x_2}\bar\Phi^1_q)
\end{aligned}
\]
Therefore, invoking \eqref{e:app_flow_map_est_fast2} and product rule/Remark \ref{rem:comm} we have for any $|\bomega| \le Y - N_*-1$
\begin{equation}
\label{e:der_diff_1}
\|D_{q+1}^{\bomega} (\det \nabla \bar \Phi_q-  \det \nabla \Phi_q) \|_{\infty} 
      \lesssim  \lambda^{-\gamma N_*+1}_{q+1} \lambda_{q+1}^{|\bomega|_x} (\delta_{q+1}^{\sfrac12} \lambda_{q+1})^{|\bomega|_t}.   
\end{equation}
This, with the estimate
\[ \| D^{\balpha}_{q+1} \bar \psi_{q+1} \|_{\infty} 
      \lesssim \frac{\delta_{q+1}^{\sfrac12}}{\lambda_{q+1}} (\delta_{q+1}^{\sfrac12} \lambda_{q+1})^{|\balpha|_t} \lambda_{q+1}^{|\balpha|_x}\]
  given by    \eqref{e:streamEs_rig2_psi}
yields for $I$ of \eqref{e:flow_diff_start}, again via product rule
\begin{equation}\label{e:flow_diff_estI}
\| D^{\balpha}_{q+1} I \|_{\infty}
      \lesssim \lambda^{-\gamma N_*+1}_{q+1} \frac{\delta_{q+1}^{\sfrac12}}{\lambda_{q+1}} (\delta_{q+1}^{\sfrac12} \lambda_{q+1})^{|\balpha|_t} \lambda_{q+1}^{|\balpha|_x}
\end{equation}
for any $|\balpha| \le Y - N_*-1$.

The term $II$ of r.h.s. of \eqref{e:flow_diff_start} involves difference $H_{12} (\xi, \tau) - H_{12}(\xi', \tau)$. In view of the definition of $H_{12}$, see
\eqref{e:H},  \eqref{e:hom:cos}, we have
\[
H_{12} (\xi, \tau) - H_{12}(\xi', \tau) = \eta_1 (\tau) \left( \sin ( 2 \pi \xi_1)- \sin ( 2 \pi \xi'_1) \right)+ \eta_2 (\tau) \left( \sin ( 2 \pi \xi_2)- \sin ( 2 \pi \xi'_2) \right)
\]
since the two summands above have identical structure, we  suppress $i=1,2$ and write, with $f( \cdot) = \sin (2\pi \cdot)$,
\[
H_{12} (\xi, \tau) - H_{12}(\xi', \tau) \sim  \eta (\tau) \left( f ( \xi_1)- f (\xi'_1) \right)
\]
so 
\[
\begin{aligned}
d (x,t) &:= H_{12} (  \lambda_{q+1} \Phi_q(x,t), \mu_{q+1} t ) - H_{12} (\lambda_{q+1} \bar \Phi_q(x,t), \mu_{q+1} t )
           ) \\
           &\sim \eta (\mu_{q+1} t )  f (\lambda_{q+1} \Phi_q(x,t)) - f (\lambda_{q+1} \bar \Phi_q(x,t))               \end{aligned}
\]
The estimates of $ D^{\balpha}_{q+1}$ of the above quantity follow now the already known path. It may be sufficient to realise that the derivatives of outer functions $\sin (\pi \cdot)$, $\cos (\pi \cdot)$ play no role in estimates with restricted number of derivatives ($\le Y$). More precisely, all the simplifying observations (i)-(iv) of section \ref{ssec:pmdbarpsi1} hold also currently. Following the lines of section \ref{ssec:pmdbarpsi1} we use \eqref{e:simple_chain_applied} with $ \bar \psi_{q+1} (x,t)$ replaced with our $d (x,t)$, and $\|D^{\balpha_j} \bar \Phi_q \|_\infty$ replaced with $\|D^{\balpha_j} (\Phi_q -\bar \Phi_q)\|_\infty$, whose estimates we already know thanks to \eqref{e:flow_map_est_diff}. This means that as we obtained \eqref{e:slow_corr_est} we have now 
\[
| D_q^{\balpha} \left( H_{12} (  \lambda_{q+1} \Phi_q(x,t), \mu_{q+1} t ) - H_{12} (\lambda_{q+1} \bar \Phi_q(x,t), \mu_{q+1} t ) \right)|
\lesssim_Y \lambda^{-\gamma N_*+1}_{q+1} \mu^{|\balpha|_t}_{q+1}   
\lambda^{|\balpha|_x}_{q+1}
\]
    so
    \begin{equation}\label{e:slow_corr_est_diff}
\frac{\delta_{q+1}^{\sfrac12}}{\lambda_{q+1}} | D_q^{\balpha} \left( H_{12} (  \lambda_{q+1} \Phi_q(x,t), \mu_{q+1} t ) - H_{12} (\lambda_{q+1} \bar \Phi_q(x,t), \mu_{q+1} t ) \right)|
\lesssim_Y \frac{\delta_{q+1}^{\sfrac12}}{\lambda_{q+1}} \lambda^{-\gamma N_*+1}_{q+1} \mu^{|\balpha|_t}_{q+1}   
\lambda^{|\balpha|_x}_{q+1}
    \end{equation}
    Observe that \eqref{e:slow_corr_est_diff} is consistent with \eqref{e:slow_corr_est}, i.e. it is the smallness $\lambda^{-\gamma N_*+1}_{q+1}$ (following from the proximity between exact and approximate flows $\bar \Phi_q$, $\Phi_q$) times r.h.s.\ of \eqref{e:slow_corr_est}.
    
    The shift to estimates for fast derivatives $D_{q+1}^{\balpha}$ is now identical to section \ref{ssec:pmdbarpsi2}. Now, instead of starting with \eqref{e:slow_corr_est} we start with \eqref{e:slow_corr_est_diff} and obtain
\[
\begin{aligned}
\frac{\delta_{q+1}^{\sfrac12}}{\lambda_{q+1}} &| D_{q+1}^{\balpha} \left( H_{12} (  \lambda_{q+1} \Phi_q(x,t), \mu_{q+1} t ) - H_{12} (\lambda_{q+1} \bar \Phi_q(x,t), \mu_{q+1} t ) \right)|\\
&\lesssim_Y \frac{\delta_{q+1}^{\sfrac12}}{\lambda_{q+1}} \lambda^{-\gamma N_*+1}_{q+1} (\delta_{q+1}^{\sfrac12} \lambda_{q+1})^{|\balpha|_t} \lambda_{q+1}^{|\balpha|_x}.
    \end{aligned}
      \]
      This estimate, together with \eqref{e:app_flow_map_est_fast} yields 
      \begin{equation}\label{e:flow_diff_estII}
\| D^{\balpha}_{q+1} II \|_{\infty}
      \lesssim \lambda^{-\gamma N_*+1}_{q+1} \frac{\delta_{q+1}^{\sfrac12}}{\lambda_{q+1}} (\delta_{q+1}^{\sfrac12} \lambda_{q+1})^{|\balpha|_t} \lambda_{q+1}^{|\balpha|_x}.
\end{equation}
The two estimates \eqref{e:flow_diff_estI}, \eqref{e:flow_diff_estII} in \eqref{e:flow_diff_start} give the wanted \eqref{e:streamEs_rig2_nobar_psi_diff} for $ D^{\balpha}_{q+1}$ derivatives. (As before, for any $|\balpha| \le Y - N_*-1$, since the known derivatives of approximate flow $\Phi_q$ are restricted to $|\balpha| \le Y - N_*$, and zeroth-order formula \eqref{e:flow_diff_start} involves first derivative $\Phi_q$.)

\subsubsection{Estimates for $ O^{\balpha}_{q+1}$ derivatives}
The only difference between 
$ O^{\balpha}$ derivatives and $D^{\balpha}$ derivatives lies in the flow involved, since in view of their definition \eqref{e:26:constructVF}   
\begin{equation}\label{e:OvsD_dervs}
D_{q+1} = (D_q + \nabla^\perp \bar \psi_{q+1} \cdot \nabla) \quad \text{ whereas } \quad O_{q+1} = (D_q + \nabla^\perp \psi_{q+1} \cdot \nabla).
\end{equation}
At this stage of the proof we have full information on estimates for $D^{\balpha}_{q+1} \psi_{q+1}$, namely
\[      \| D^{\balpha}_{q+1} \psi_{q+1} \|_{\infty}
      \lesssim \frac{\delta_{q+1}^{\sfrac12}}{\lambda_{q+1}} (\delta_{q+1}^{\sfrac12} \lambda_{q+1})^{|\balpha|_t} \lambda_{q+1}^{|\balpha|_x} \qquad \text{for any }  |\balpha| \leq Y-N_*-1.
      \]
      This estimate is consistent with the estimate \eqref{e:slow_corr_est_appl} for $\psi_{q+1}$ used in section \ref{ssec:pmdbarpsi2} to deal with $D^{\balpha}_{q+1}$. Therefore, following the lines of the proof from section \ref{ssec:pmdbarpsi2} onwards, we have the inequalities \eqref{e:streamEs_rig2_nobar} \eqref{e:streamEs_rig2_nobar_psi} \eqref{e:streamEs_rig2_nobar_psi_diff} for $ O^{\balpha}_{q+1}$. Concerning the admissible highest order of derivatives, observe that the highest possible purely space derivative is of order $Y-N_*-1$, because for pure space derivatives $ O^{\balpha}_{q+1} = D^{\balpha}_{q+1}$, and the highest possible mixed derivative (i.e.\ with at least one advective derivative) is also of order $Y-N_*-1$, because as before there is no loss of total number of derivatives in this case, as one sees from \eqref{e:OvsD_dervs}.

\subsubsection{Estimates for $ D^{\balpha}_{q+1} - O^{\balpha}_{q+1}$ derivatives}
In view of \eqref{e:OvsD_dervs}
\begin{equation}\label{e:DOdiff1}
(D_{q+1} - O_{q+1})f =   (\nabla^\perp \bar \psi_{q+1} - \nabla^\perp \psi_{q+1}) \cdot \nabla f.
    \end{equation}
Let us make the following observations:\\
(i) The cost of derivative of $(\nabla^\perp \bar \psi_{q+1} - \nabla^\perp \psi_{q+1})$ is
$\delta_{q+1}^{\sfrac12} \lambda_{q+1}$ for an advective derivative (no matter if it is $O_{q+1}$ or $D_{q+1}$ derivative) and $\lambda_{q+1}$ for space derivative, and the amplitude of  $(\nabla^\perp \bar \psi_{q+1} - \nabla^\perp \psi_{q+1})$ is $\delta_{q+1}^{\sfrac12} \lambda^{-\gamma N_*+1}_{q+1}$, all thanks to \eqref{e:streamEs_rig2_nobar_psi_diff}. This means that $D_{q+1} - O_{q+1}$ applied to any $f$ with amplitude $M$ and
$\delta_{q+1}^{\sfrac12} \lambda_{q+1}$ cost for an advective derivative  and $\lambda_{q+1}$ cost for space derivative, results in a function with amplitude $ \delta_{q+1}^{\sfrac12} \lambda^{-\gamma N_*+1}_{q+1} M \lambda_{q+1}$ and
$\delta_{q+1}^{\sfrac12} \lambda_{q+1}$ cost for an advective derivative  and $\lambda_{q+1}$ cost for space derivative. \\
(ii) Functions $\bar \psi_{q+1}$, $\psi_{q+1}$ (or any of their derivatives) satisfy the conditions for $f$ of (i), thanks to \eqref{e:streamEs_rig2_nobar_psi}, with amplitude $M=\frac{\delta_{q+1}^{\sfrac12} }{\lambda_{q+1}}$ (and appropriately for derivatives).\\
(iii) It holds (for $D=D_{q+1}$, $O=O_{q+1}$)
\begin{equation}\label{e:DOdiffn}
(D^m - O^m) = (D - O) D^{m-1} + O (D - O) D^{m-2} + O^2 (D - O) D^{m-3} + \dots + O^{m-1} (D - O)
    \end{equation}
(iv) The observations (i)-(iii) imply in particular that $(D_{q+1}^m - O_{q+1}^m) \bar \psi_{q+1}$ or $(D_{q+1}^m - O_{q+1}^m) \psi_{q+1}$ is a function with amplitude  $m \delta_{q+1}^{\sfrac12} \lambda^{-\gamma N_*+1}_{q+1} \delta_{q+1}^{\sfrac12} \sim  \lambda^{-\gamma N_*+1}_{q+1} \delta_{q+1}$,
 with $\delta_{q+1}^{\sfrac12} \lambda_{q+1}$ cost of an advective derivative  and $\lambda_{q+1}$ cost for space derivative. (The combinatorial factor $m$ is negligible, since in the estimates we have $\lesssim_Y$, with $Y$ being the highest possible derivative order.)\\
 (v) Analogously to (iv), using (ii) for derivatives, we have that
 \[
 \begin{aligned}
 (D_{q+1}^m - O_{q+1}^m) D^{\bomega}_{q+1} \bar \psi_{q+1}, \\ (D_{q+1}^m - O_{q+1}^m) D^{\bomega}_{q+1}  \psi_{q+1}, \\ (D_{q+1}^m - O_{q+1}^m) O^{\bomega}_{q+1} \bar \psi_{q+1}, \\ (D_{q+1}^m - O_{q+1}^m) O^{\bomega}_{q+1} \psi_{q+1}
\end{aligned}
 \]
are functions with amplitude  $\sim \lambda^{-\gamma N_*+1}_{q+1} \delta_{q+1} (\delta_{q+1}^{\sfrac12} \lambda_{q+1})^{|\bomega|_t} \lambda_{q+1}^{|\bomega|_x}$ and the costs of derivatives as before.

(vi) Furthermore, any $D^{\balpha}_{q+1} - O^{\balpha}_{q+1}$ consists of blocks of advective derivatives and of space derivatives, for instance 
\begin{equation}\label{e:DOblocks}
\begin{aligned}
D^{\balpha}_{q+1} - O^{\balpha}_{q+1} &= D^{n_1}_{q+1} \nabla^{m_1} .. D^{n_k}_{q+1} - O^{n_1}_{q+1} \nabla^{m_1} .. O^{n_k}_{q+1} \\
&= D^{n_1}_{q+1} \nabla^{m_1} .. (D^{n_k}_{q+1} - O^{n_k}_{q+1}) +  \dots + (D^{n_1}_{q+1}  - O^{n_1}_{q+1}) \nabla^{m_1} .. O^{n_k}_{q+1}.
\end{aligned}
\end{equation}
Using observation (v) in \eqref{e:DOblocks} gives that 
\[
 (D^{\balpha}_{q+1} - O^{\balpha}_{q+1})  \bar \psi_{q+1}, \quad  (D^{\balpha}_{q+1} - O^{\balpha}_{q+1}) \psi_{q+1}
 \]
are functions with amplitude  $\sim \lambda^{-\gamma N_*+1}_{q+1} \delta_{q+1} (\delta_{q+1}^{\sfrac12} \lambda_{q+1})^{|\balpha|_t} \lambda_{q+1}^{|\balpha|_x}$ and their costs of derivatives as before. This justifies the wanted estimate \eqref{e:streamEs_psi_DOdiff}.

Using again (v) and reasoning along (vi) shows that, denoting by  $\Theta^{\balpha}_{q+1}$ either   $D^{\balpha}_{q+1}$ or  $O^{\balpha}_{q+1}$, the symbol
\begin{equation}\label{e:DOblocks_sym}
 (D^{\balpha_1}_{q+1} - O^{\balpha_1}_{q+1}) \Theta^{\balpha_2}_{q+1}  (D^{\balpha_3}_{q+1} - O^{\balpha_3}_{q+1}) \dots \Theta^{\balpha_k}_{q+1} 
 \end{equation}
acting on either $\bar \psi_{q+1}$ or on $\psi_{q+1}$ yield functions with \[
\text{amplitude  } \quad \sim \lambda^{-\gamma N_*+1}_{q+1} \delta_{q+1} (\delta_{q+1}^{\sfrac12} \lambda_{q+1})^{|\balpha|_t} \lambda_{q+1}^{|\balpha|_x}
\] and the costs of derivatives as before, i.e.\ $\delta_{q+1}^{\sfrac12} \lambda_{q+1}$ is cost of an advective derivative  and $\lambda_{q+1}$ is cost for a space derivative. Of course, the only requirement for the symbol of type \eqref{e:DOblocks_sym} is that it has at least one block of  $(D^{\balpha'}_{q+1} - O^{\balpha'}_{q+1})$, its outermost and the innermost part may be any of $(D^{\bomega}_{q+1} - O^{\bomega}_{q+1})$, $\Theta^{\balpha}_{q+1}$. Observe also that we do not claim further gain from more than one block of type $(D^{\bomega}_{q+1} - O^{\bomega}_{q+1})$ appearing, so a brutal estimate replacing all but one blocks $(D^{\bomega}_{q+1} - O^{\bomega}_{q+1})$ with the respective sum suffices. 

The upper bound of number of admissible derivatives is dictated by what is available from \eqref{e:streamEs_rig2_nobar} and \eqref{e:streamEs_rig2_nobar_psi}, which is $Y-N_*-1$ (and there is no further loss of derivatives, as before in the case of mixed derivatives; indeed otherwise, i.e.\ in case $\balpha$ are only pure space derivatives,  $(D^{\balpha}_{q+1} - O^{\balpha}_{q+1}) =0$.) We have reached \eqref{e:streamEs_psi_DOdiff_gen}.

Replacing observation (ii) with:\\
(ii') Any derivative of $\bar \Psi_{q+1}$ or of $\Psi_{q+1}$ satisfies the conditions for $f$ of (i), thanks to \eqref{e:streamEs_rig2_nobar} with (lossy) amplitude $M=1$ for (zeroth-order) $\bar \Psi_{q+1}$ or $\Psi_{q+1}$ and appropriately for derivatives.\\
We obtain \eqref{e:streamEs_Psi_DOdiff}, \eqref{e:streamEs_Psi_DOdiff_gen} along the preceding lines.

In fact, as already seen from the original observation (ii), for any smooth function $f$ such that 
\begin{equation}\label{e:f_alastream_ass}
      \| D^{\balpha}_{q+1} f \|_{\infty} +   \| O^{\balpha}_{q+1} f \|_{\infty} 
      \lesssim M (\delta_{q+1}^{\sfrac12} \lambda_{q+1})^{|\balpha|_t} \lambda_{q+1}^{|\balpha|_x} \qquad \text{for any }  |\balpha| \leq Y-N_*-n_0.
\end{equation}
we have the same conclusions, thus Corollary \ref{l:strDiffEst3} holds, under the weaker assumption that we have control on both $D$ and $O$ derivatives, i.e.\ taking \eqref{e:f_alastream_ass} in place of the 'either-or' of \eqref{e:f_alastream_assD}, \eqref{e:f_alastream_assO}.

In order to upgrade from assuming \eqref{e:f_alastream_ass} to using merely 'either-or' of Corollary \ref{l:strDiffEst3}, we proceed as follows. Assume we know only the behavior of $D$-derivatives, which means that we assume only \eqref{e:f_alastream_assD} (the other case is symmetric). But this assumption already contains knowledge of pure space derivatives, so we can bound 
$(D_{q+1} - O_{q+1})f$ in view of the identity \eqref{e:DOdiff1}. Analogously we can bound $(D_{q+1} - O_{q+1})^m f$ up to maximal order allowed by the assumption. But (with $D=D_{q+1}$, $O=O_{q+1}$)
\[
D^2 - O^2 = (D-O) D + \underbrace{O (D-O)}_{=D(D-O)-(D-O)^2}
\]
so we have also appropriate bound on $D^2 - O^2$, consequently on $O^2$. Repeating this reasoning inductively for higher-order terms we see that one recovers the lacking estimate \eqref{e:f_alastream_assO}
on $O^{\balpha}_{q+1}$. From this point onwards we have all the assumptions of Corollary \ref{l:strDiffEst3}.

\newpage

\section{Decomposition lemmas}      \label{s:decomposition}

As we see in Section \ref{ss:poly_space}, the vector fields at step $q$ generate oscillatory functions documented by polynomial-derivative space $\mathcal{P}_q$ and $\mathcal{O}_q$ at step $q$. In this section, we quantitatively decompose the elements in $\mathcal{O}_{q+1}$ into products of elements in $\mathcal{P}_q$ and known oscillatory profiles.

\subsection{Statements of the decomposition lemmas}

From \eqref{e:2:uEstimate} and \eqref{e:4:uEstimate}, we have 
\begin{align} 
      D_{q,t} \Phi_q (x,t) =&\, \omega_{q,N}(x,t) \phi_{N} (\mu_{q+1} t) + \sum_{m=1}^N \omega_{q,m} (x,t) z_{m} (\mu_{q+1} t).       \label{e:DqtPhiq} 
\end{align}
For any $\eta \in \T \rightarrow \R$ with
\begin{align}     \label{e:asmpSupp}
      \supp \eta \subset \supp \eta_1 \text{ or } \supp \eta \subset \supp \eta_2,
\end{align}
we have
\begin{align}     \label{e:cons:asmpSupp}
      \eta( \mu_{q+1} t) D_{q,t} \Phi_q (x,t) =&\, \omega_{q,N}(x,t) \phi_{N} (\mu_{q+1} t) \eta( \mu_{q+1} t).
\end{align}

Now we can state the decomposition lemmas.

\begin{lemma}     \label{l:decompPrep}
For any $\alpha \in \mathcal{I}$ with $|\balpha| \leq 2NQ$, we have 
\begin{align}     \label{e:l0:decompPrep}
      \mr O_{q+1}^{\balpha} \Bigg( \frac{\mr \lambda_{q+1}^2} {\delta_{q+1}^{\sfrac12}\lambda_{q+1}^{1+2\gamma_R}} \Psi_{q+1} \Bigg) 
            = g_0 + \sum_{k=1}^{k_*} \ddot\chi_{k}( \lambda_{q+1} \Phi_q ) \ddot\eta_{k} (\mu_{q+1} t) g_{k}
\end{align}
for some $k_* \in \N$, $\ddot{\chi}_k: \T^2 \rightarrow \R$, $\ddot{\eta}_k: \T \rightarrow \R$, $g_k \in C^\infty(\T^2 \times[0,1])$ with 
\begin{align}
      k_* \leq&\, 2N ( 3N^2 )^{|\balpha|},       \label{e:l2:decompPrep} \\ 
      \max_{k \geq 1} \| \nabla_\xi^p \ddot\chi_{k} \|_{\infty} 
            \lesssim&\, 1,        \quad p \leq 8Q^3,        \label{e:l3:decompPrep} \\ 
      \max_{k \geq 1} \| \partial_\tau^p \ddot\eta_{k} \|_{\infty} 
            \lesssim&\, 1.        \quad p \leq 8Q^3.       \label{e:l4:decompPrep} 
\end{align}
For any $k$, $(\ddot\chi_{k}, \ddot\eta_{k})$ satisfies the generalized shear condition. Moreover, we have 
\begin{align}
      g_{k} &\,\in \mathcal{P}_q ( 3N(|\balpha|+1) ),       
            \label{e:l5:decompPrep}  \\ 
      \vertiii{ g_0 }_q &\,\lesssim r_0 ( |\balpha|_x, |\balpha|_t ),       
            \label{e:l6:decompPrep}  \\ 
      \sum_{k \geq 1} \vertiii{ g_k }_q &\,\lesssim r_1 ( |\balpha|_x, |\balpha|_t ), 
            \label{e:l7:decompPrep}
\end{align}
with some $\{r_0 (i,j)\}_{i,j \in \N} \subset \R$, $\{r_1 (i,j)\}_{i,j \in \N} \subset \R$ satisfying
\begin{align}
      r_1 (0,0) =&\, \frac{\mr \lambda_{q+1}^2} {\lambda_{q+1}^{2(1+\gamma_R)}},   \quad 
      r_0 (0,0) = \frac{ \mr\lambda^2_{q+1} }{ \mr\lambda^2_{q} } \frac{\delta_{q}^{\sfrac12}\lambda_q^{1+2\gamma_R}} {\delta_{q+1}^{\sfrac12}\lambda_{q+1}^{1+2\gamma_R}},       \label{e:l9:decompPrep} \\ 
      r_1 (i+1,j) \leq&\, \frac{\lambda_{q+1}} {\mr \lambda_{q+1}} r_1 (i,j), \quad 
      r_0 (i+1,j) \leq \frac{\mr \lambda_{q}} {\mr \lambda_{q+1}} r_0 (i,j),        \label{e:l10:decompPrep}  \\ 
      r_1 (i,j+1) \leq&\, \max\bigg\{ \frac{\mu_{q+1}}{\mr \mu_{q+1}}, \frac{\delta_{q+1}^{\sfrac12}\lambda_{q+1}}{\mr \mu_{q+1}} \bigg\} r_1 (i,j)
            + \max \bigg\{ \frac{\mr \mu_{q}}{\mr \mu_{q+1}}, \frac{ \delta_{q+1}^{\sfrac12} \mr \lambda_{q} } {\mr \mu_{q+1}} \bigg\} r_0 (i,j),     \label{e:l11:decompPrep}  \\ 
      r_0(i,j+1) \leq&\, \frac{\mr \mu_{q}}{\mr \mu_{q+1}} r_0 (i,j).   \label{e:l12:decompPrep} 
\end{align}
\end{lemma}

\begin{lemma}     \label{l:decompRho}
For any $\chi: \T^2 \rightarrow \R$, any $\eta: \T \rightarrow \R$ and any $\rho: \T^2 \times [0,1] \rightarrow \R$ satisfying
\begin{align*}
      \| \nabla_\xi^p \chi \|_{\infty} \lesssim&\, 1,       \quad p \leq p_0,      \\ 
      \| \partial_\tau^p \eta \|_{\infty} \lesssim&\, 1,    \quad p \leq p_0
\end{align*}
for some $p_0 \in \N$. Define
\begin{align}     \label{e:0:decompRho}
      g(x,t) =  h(x,t) \chi \big( \lambda_{q+1} \Phi_q(x,t) \big) \eta (\mu_{q+1} t) \rho(x,t).
\end{align}
Suppose $(\chi, \eta)$ satisfies generalized shear condition and $h \in \mathcal{P}_q (Q)$. Then we have the following decomposition for any $\alpha \in \mathcal{I}$ with $|\balpha| \leq 3Q$
\begin{align*}
      \mr O_{q+1}^{\balpha} g = \sum_{|\bomega| \leq |\balpha|} \mr D_{q}^{\bomega} \rho \sum_{k=1}^{\Gamma^{(\bomega)}} h_{k}^{(\bomega)} \chi_{k}^{(\bomega)} ( \lambda_{q+1} \Phi_q ) \eta_{k}^{(\bomega)} (\mu_{q+1} t)
\end{align*}
with
\begin{align}
      h_{k}^{(\bomega)}       &\,\in \mathcal{P}_q \big( Q + 2(N+2) |\balpha| \big),       
            \label{e:3:decompRho} \\ 
      \sum_{k,\bomega} \vertiii{ h_{k}^{(\bomega)} }_q 
            &\,\lesssim \bigg( \frac{\lambda_{q+1}}{\mr \lambda_{q+1}} \bigg)^{|\balpha|_x} \bigg( \frac{\delta_{q+1}^{\sfrac12} \lambda_{q+1}}{\mr \mu_{q+1}} \bigg)^{|\balpha|_t} \vertiii{h}_q    \label{e:4:decompRho}
\end{align}
and
\begin{align}
      \Gamma^{(\bomega)}      &\,\leq ( 3N^2 )^{|\balpha|},       \label{e:2:decompRho} \\ 
      \| \nabla_\xi^p \chi_{k}^{(\bomega)} \|_{\infty}      
            &\,\lesssim 1,    \quad p \leq p_0 - |\bomega|,      \label{e:5:decompRho} \\ 
      \| \partial_\tau^p \eta_{k}^{(\bomega)} \|_{\infty}
            &\,\lesssim 1.    \quad p \leq p_0 - |\bomega|_t.       \label{e:6:decompRho}
\end{align}
Moreover, $(\chi_{k}^{(\bomega)}, \eta_{k}^{(\bomega)})$ satisfies the generalized shear condition for any $k$ and $\bomega$. Furthermore, if we assume $\chi$ is constant, then we have
\begin{align}
      \sum_{k,\bomega} \vertiii{ h_{k}^{(\bomega)} }_q 
            &\,\lesssim \bigg( \frac{ \mr \lambda_q } {\mr \lambda_{q+1}} \bigg)^{|\balpha|_x} \bigg( \frac{\delta_{q+1}^{\sfrac12} \lambda_{q+1}}{\mr \mu_{q+1}} \bigg)^{|\balpha|_t} \vertiii{h}_q.       \label{e:8:decompRho}
\end{align}
\end{lemma}

The proofs of the above two lemmas are postponed to the end of this section.

\begin{corollary}       \label{c:decomp}
For any $g \in \mathcal{O}_{q+1} (2NQ)$, we have the following decomposition
\begin{align}     \label{e:0:decomp}
      g = g_0 + \sum_{k=1}^{k_*} \ddot{\chi}_{k}( \lambda_{q+1} \Phi_q ) \ddot{\eta}_{k} (\mu_{q+1} t) g_{k}
\end{align}
with
\begin{align}
      g_{k} \in \mathcal{P}_q ( 9N^2Q ),        \quad
      \sum_{k=0}^{k_*} \vertiii{ g_{k} }_q 
            \leq \vertiii{ g }_{q+1}         \label{e:6:decomp} 
\end{align}
and
\begin{align}
      k_*         &\,\leq (10N)^{8NQ},    \label{e:2:decomp} \\ 
      \max_{k} \| \nabla_\xi^p \ddot{\chi}_{k} \|_{\infty} 
            &\,\leq 1,        \quad p \leq 8Q^3,      \label{e:8:decomp} \\ 
      \max_{k} \| \partial_\tau^p \ddot{\eta}_{k} \|_{\infty} 
            &\,\leq 1,        \quad p \leq 8Q^3.       \label{e:10:decomp}
\end{align}
Moreover, $(\ddot{\chi}_{k}, \ddot{\eta}_{k})$ satisfies the generalized shear condition for any $k$.
\end{corollary}

\begin{proof}
If $g$ is a constant, all statements are trivial. Assume $d_* \leq 2NQ$. For any $g \in \mathcal{O}_{q+1} (d_*)$ that are not constant, we consider the representation of $g$ that attains $\vertiii{g}_{q+1}$, i.e. 
\begin{align}
       g = \sum_{k=1}^{k_0} a_k \prod_{i=1}^{I_k} \mr O_{q+1}^{\balpha_{k,i}} \Bigg( \frac{\mr \lambda_{q+1}^2} {\delta_{q+1}^{\sfrac12}\lambda_{q+1}^{1+2\gamma_R}}  \Psi_{q+1} \Bigg), \quad 
       \sum_{k=1}^{k_0} |a_k| = \vertiii{g}_{q+1}.
\end{align}
Without loss of generality, we also assume $k_0 \geq 1$, $a_k \neq 0$ and $I_k \geq 1$ for any $k$.

For $1 \leq d \leq d_*$, we define
\begin{align}
      \Lambda(d) :=&\, \bigg\{ 1 \leq k \leq k_0 \,\bigg|\, \sum_{i} |\balpha_{k,i}| = d \bigg\}.
\end{align}
From elementary combinatorics, we have the following estimate on the cardinality of $\Lambda(d)$
\begin{align}
      |\Lambda(d)| \leq 6^{d}.      \label{e:18:decomp}
\end{align}
For a fixed $d$, we apply Lemma \ref{l:decompPrep} to those $k \in \Lambda(d)$.
\begin{align}
      \prod_{i=1}^{I_k} \mr D_{q+1}^{\balpha_{k,i}} \Bigg( \frac{\mr \lambda_{q+1}^2} {\delta_{q+1}^{\sfrac12}\lambda_{q+1}^{1+2\gamma_R}}  \Psi_{q+1} \Bigg)
      =&\, \prod_{i=1}^{I_k} \Bigg( g_{k,i,0} + \sum_{j=1}^{J_{k,i}} \ddot\chi_{k,i,j}( \lambda_{q+1} \Phi_q ) \ddot\eta_{k,i,j} (\mu_{q+1} t)
            g_{k,i,j} \Bigg)        \label{e:20:decomp} \\ 
      =&\, g_{k,0} + \sum_{j=1}^{J_{k}} \ddot\chi_{k,j}( \lambda_{q+1} \Phi_q ) \ddot\eta_{k,j} (\mu_{q+1} t) g_{k,j}       \label{e:22:decomp}
\end{align}
From \eqref{e:20:decomp} to \eqref{e:22:decomp}, we just expand the product in $i$, then we rewrite and renumber the product in the form of \eqref{e:22:decomp}. From \eqref{e:l2:decompPrep}, we have the following estimate for $k \in \Lambda(d)$
\begin{align}     \label{e:24:decomp}
      1 + J_k \leq \prod_{i=1}^{I_k} (1+J_{k,i}) \leq (3N^2)^{ \sum_{i=1}^{I_k}(|\balpha_{k,i}|+1) }
      \leq (3N^2)^{2d}.
\end{align}

From Lemma \ref{l:decompPrep}, we also have
\begin{align}     
      \max_j \| \nabla_\xi^p \ddot\chi_{k,i,j} \|_{\infty} 
            \lesssim&\, 1,    \quad p \leq 8Q^3,      \label{e:26:decomp} \\ 
      \max_j \| \partial_\tau^p \ddot\eta_{k,i,j} \|_{\infty} 
            \lesssim&\, 1,    \quad p \leq 8Q^3       \label{e:28:decomp}
\end{align}
and 
\begin{align}
      g_{k,i,j} \in&\, \mathcal{P}_q \big( 3N (|\balpha_{k,i}|+1) \big),      \label{e:30:decomp} \\ 
      \sum_{j=0}^{J_{k,i}} \vertiii{ g_{k,i,j} }_q \lesssim&\,
            \lambda_{q+1}^{-2\gamma_R}
            \bigg( \frac{\lambda_{q+1}}{\mr \lambda_{q+1}} \bigg)^{|\balpha_{k,i}|_x-2} 
            \bigg( \frac{\delta_{q+1}^{\sfrac12} \lambda_{q+1}}{\mr \mu_{q+1}} \bigg)^{|\balpha_{k,i}|_t}.        \label{e:32:decomp}
\end{align}
Note that $\balpha_{k,i} \in \mathcal{I}(2,0) \cup \mathcal{I}(1,N_*)$ from \eqref{e:0:poly_Uq}. Using $N\gamma_I > 2$ from \eqref{e:smallGammaR} and \eqref{e:NQ}, we absorb the constants and deduce
\begin{align}     \label{e:36:decomp}
      \sum_{j=0}^{J_{k,i}} \vertiii{ g_{k,i,j} }_q 
            &\,\lesssim \lambda_{q+1}^{-\frac{\gamma_R}{2}},   \quad 
      \prod_{i=1}^{I_k} \sum_{j=0}^{J_{k,i}} \vertiii{ g_{k,i,j} }_q 
            \lesssim \lambda_{q+1}^{-\frac{I_k\gamma_R}{2}}.
\end{align}

Then from the product rules, \eqref{e:26:decomp} and \eqref{e:28:decomp}, we have
\begin{align*}
      \max_j \| \nabla_\xi^p \ddot\chi_{k,j} \|_{\infty} \lesssim&\, 1,     \quad p \leq 8Q^3,      \\ 
      \max_j \| \partial_\tau^p \ddot\eta_{k,j} \|_{\infty} \lesssim&\, 1,  \quad p \leq 8Q^3.
\end{align*}

From Remark \ref{r:opPolynomial} and \eqref{e:36:decomp}, we have for $k \in \Lambda(d)$
\begin{align}     \label{e:38:decomp}
      \sum_{j=0}^{J_k} \vertiii{ g_{k,j} }_q
      \lesssim \prod_{i=1}^{I_k} \sum_{j=0}^{J_{k,i}} \vertiii{ g_{k,i,j} }_q 
      \lesssim \lambda_{q+1}^{-\frac{I_k\gamma_R}{2}}.
\end{align}

Now we sum up \eqref{e:22:decomp} in $k$. With slight abuse of notations, we renumber all terms in $k$ and $j$ in a single index, again denoted by $ 1 \leq k \leq k_*$, to arrive at the form \eqref{e:0:decomp}. From \eqref{e:18:decomp} and \eqref{e:24:decomp}, the cardinality satisfies the estimate which verifies \eqref{e:2:decomp},
\begin{align*}
      k_*   \leq \sum_{d=1}^{d_*} \sum_{k \in \Lambda(d)} (3N^2)^{2d} 
            \leq \sum_{d=1}^{d_*} (9N)^{4d} \leq (10N)^{4d_*}.
\end{align*}

Also, \eqref{e:38:decomp} yields the following estimate verifying \eqref{e:6:decomp},
\begin{align}     \label{e:40:decomp}
      \sum_{k} \vertiii{ g_{k} }_q 
      \lesssim |a_0| + \sum_{d=1}^{d_*} \sum_{k \in \Lambda(d)} |a_k| \sum_{j=0}^{J_k} \vertiii{ g_{k,j} }_q 
      \lesssim \vertiii{ g }_{q+1}.
\end{align}
From \eqref{e:30:decomp} and Remark \ref{r:opPolynomial}, we have
\begin{align}     \label{e:42:decomp}
      g_k \in \mathcal{P}_q \bigg( 3N \bigg( I_k + \sup_k \sum_{i=1}^{I_k} |\balpha_{k,i}| \bigg) \bigg).
\end{align}
Since $g \in \mathcal{O}_{q+1} (2NQ) $, we have $I_k \leq NQ$. Then \eqref{e:6:decomp} follows from \eqref{e:40:decomp} and \eqref{e:42:decomp}.
\end{proof}

\subsection{Proof of Lemma \ref{l:decompPrep} and Lemma \ref{l:decompRho}}

In the rest of this section, we prove Lemma \ref{l:decompPrep} and Lemma \ref{l:decompRho}. The proofs are similar inductive arguments.

\begin{proof}[Proof of Lemma \ref{l:decompPrep}]
All the indices $m,m_1,m_2$ are summed from $1$ to $N$. When $|\balpha|=0$, we have
\begin{align*}
      \frac{\mr \lambda_{q+1}^2} {\delta_{q+1}^{\sfrac12}\lambda_{q+1}^{1+2\gamma_R}} \Psi_{q+1} 
      = \frac{ \mr\lambda^2_{q+1} }{ \mr\lambda^2_{q} } 
            \frac{\delta_{q}^{\sfrac12}\lambda_q^{1+2\gamma_R}} {\delta_{q+1}^{\sfrac12}\lambda_{q+1}^{1+2\gamma_R}} 
            \Bigg( \frac{\mr \lambda_{q}^2} {\delta_{q}^{\sfrac12}\lambda_{q}^{1+2\gamma_R}}  \bar \Psi_q \Bigg) 
      + \frac{\mr \lambda_{q+1}^2} {\delta_{q+1}^{\sfrac12}\lambda_{q+1}^{1+2\gamma_R}} \psi_{q+1}
\end{align*}
and 
\begin{align}
      \psi_{q+1} =&\,   
            \frac{\delta_{q+1}^{\sfrac12}} {\lambda_{q+1}} \det \nabla \Phi_q \Big( \eta_1(\mu_{q+1} t) \Pi_1( \lambda_{q+1} \Phi_q ) + \eta_2(\mu_{q+1} t) \Pi_2( \lambda_{q+1} \Phi_q ) \Big)  \nonumber \\ 
            =&\, \frac{\delta_{q+1}^{\sfrac12}} {\lambda_{q+1}} \sum_{m} \varpi_{q,m} \zeta_{m} (\mu_{q+1} t) \Big( \eta_1(\mu_{q+1} t) \Pi_1( \lambda_{q+1} \Phi_q ) + \eta_2(\mu_{q+1} t) \Pi_2( \lambda_{q+1} \Phi_q ) \Big)      \nonumber \\ 
            =&\, \frac{\delta_{q+1}^{\sfrac12}} {\lambda_{q+1}} \sum_{m} \varpi_{q,m} \Big( [ \eta_1\zeta_m ] (\mu_{q+1} t) \Pi_1( \lambda_{q+1} \Phi_q ) + [ \eta_2\zeta_m ] (\mu_{q+1} t) \Pi_2( \lambda_{q+1} \Phi_q ) \Big).        \label{e:5:0:decompPrep}
\end{align}
Then we renumber and write the terms in \eqref{e:5:0:decompPrep} in the form of \eqref{e:l0:decompPrep}, giving
\begin{align}
      \frac{\mr \lambda_{q+1}^2} {\delta_{q+1}^{\sfrac12}\lambda_{q+1}^{1+2\gamma_R}} 
            \Psi_{q+1} =&\, g_{0,0} + \sum_{k=1}^{k_0} g_{0,k} \chi_{0,k}( \lambda_{q+1} \Phi_q ) \eta_{0,k} (\mu_{q+1} t),    \nonumber \\ 
      g_{0,0} =&\, \frac{ \mr\lambda^2_{q+1} }{ \mr\lambda^2_{q} } 
            \frac{ \delta_{q}^{\sfrac12}\lambda_q^{1+2\gamma_R} } { \delta_{q+1}^{\sfrac12}\lambda_{q+1}^{1+2\gamma_R} } 
            \bigg( \frac{ \mr \lambda_{q}^2 } { \delta_{q}^{\sfrac12}\lambda_{q}^{1+2\gamma_R} }  \bar \Psi_q \bigg),        \label{e:5:2:decompPrep} \\ 
      g_{0,m}= g_{0,N+m}= &\, \frac{ \mr \lambda_{q+1}^2 } { \lambda_{q+1}^{2(1+\gamma_R)} } \varpi_{q,m},
            \quad \text{for any } 1 \leq m \leq N,    \label{e:5:3:decompPrep}
\end{align}
with $k_0 = 2N$ and $\chi_{0,k}, \eta_{0,k}$ given by \eqref{e:5:0:decompPrep}. From Lemma \ref{l:uEstimate}, we have
\begin{align*}
      \max_{k \geq 1} \vertiii{ g_{0,k} }_q 
            \leq \Bigg( \frac{\mr \lambda_{q+1}} {\lambda_{q+1}^{1+\gamma_R}} \Bigg)^2.
\end{align*}
From \eqref{e:5:2:decompPrep}, we also have
\begin{align*}
      \vertiii{ g_{0,0} }_q \leq \frac{ \mr\lambda^2_{q+1} }{ \mr\lambda^2_{q} } 
            \frac{\delta_{q}^{\sfrac12}\lambda_q^{1+2\gamma_R}} {\delta_{q+1}^{\sfrac12}\lambda_{q+1}^{1+2\gamma_R}}.
\end{align*}
This verifies the statements when $|\balpha| = 0$. In particular, \eqref{e:l3:decompPrep} follows from \eqref{e:hom:cos}. \eqref{e:l4:decompPrep} follows from Lemma \ref{l:cutoff} and \eqref{e:phiZetaThetaEst}. (\ref{e:l5:decompPrep}-\ref{e:l9:decompPrep}) follow from \eqref{e:algbVarpiOmega}, \eqref{e:mtrPolyEst} and Remark \ref{r:opPolynomial}.

Let $n:= |\bbeta|$. Assume that for some $k_n \in \N^+$, we have
\begin{align}     \label{e:6:decompPrep}
      \mr O_{q+1}^{\bbeta} \Bigg( \frac{\mr \lambda_{q+1}^2} {\delta_{q+1}^{\sfrac12}\lambda_{q+1}^{1+2\gamma_R}} \Psi_{q+1} \Bigg) 
      = g_{\bbeta,0} + \sum_{k=1}^{k_n} g_{\bbeta,k} \chi_{\bbeta,k}( \lambda_{q+1} \Phi_q ) \eta_{\bbeta,k} (\mu_{q+1} t), 
\end{align}
And we assume (\ref{e:l2:decompPrep}-\ref{e:l7:decompPrep}) for $\bbeta$. From such inductive assumptions, our goal is to derive the information on
\begin{align}     \label{e:9:decompPrep}
      \mr O_{q+1}^{\balpha} \Bigg( \frac{\mr \lambda_{q+1}^2} {\delta_{q+1}^{\sfrac12}\lambda_{q+1}^{1+2\gamma_R}} \Psi_{q+1} \Bigg) 
      = g_{\balpha,0} + \sum_{k=1}^{k_{n+1}} g_{\balpha,k} \chi_{\balpha,k}( \lambda_{q+1} \Phi_q ) \eta_{\balpha,k} (\mu_{q+1} t)
\end{align}
for some $k_{n+1} \in \N$.

We consider the two cases: $\balpha = 1 \bbeta$ or $\balpha = t \bbeta$. The case of $\balpha = 2 \bbeta$ is completely symmetric to $\balpha = 1 \bbeta$.

\begin{case}

We consider the case $\balpha = 1 \bbeta$. Then
\begin{align}     \label{e:12:decompPrep}
      \mr O_{q+1}^{\balpha} =&\, \mr \partial_{q+1,1} \mr O_{q+1}^{\bbeta}
            = \frac{1}{\mr \lambda_{q+1}} \partial_{1} \mr O_{q+1}^{\bbeta}.
\end{align}
Recalling the notion of renormalized derivatives in Definition \ref{d:renormalizedDiff}, we have
\begin{align}
      \mr O_{q+1}^{\balpha} \Bigg( \frac{\mr \lambda_{q+1}^2} {\delta_{q+1}^{\sfrac12}\lambda_{q+1}^{1+2\gamma_R}} \Psi_{q+1} \Bigg) 
      =&\, \frac{1}{\mr \lambda_{q+1}} \partial_{1} \sum_{k=1}^{k_n} g_{\bbeta,k} \chi_{\bbeta,k}( \lambda_{q+1} \Phi_q ) \eta_{\bbeta,k} (\mu_{q+1} t) 
            + \frac{1}{\mr \lambda_{q+1}} \partial_{1} g_{\bbeta,0}    \nonumber \\ 
            =&\, \frac{\mr \lambda_q}{\mr \lambda_{q+1}} \sum_{k=1}^{k_n} \chi_{\bbeta,k}( \lambda_{q+1} \Phi_q ) \eta_{\bbeta,k} (\mu_{q+1} t) \mr \partial_{q,1} g_{\bbeta,k} 
                  + \frac{\mr \lambda_q}{\mr \lambda_{q+1}} \mr \partial_{q,1} g_{\bbeta,0}    \nonumber \\ 
            +&\, \frac{\lambda_{q+1}}{\mr \lambda_{q+1}} \sum_{k=1}^{k_n} \Big[ \Id + \varepsilon_q \sum_{m} B_{q,m} \phi_{m} (\mu_{q+1} t) \Big]_{i1} \partial_{\xi_i} \chi_{\bbeta,k}( \lambda_{q+1} \Phi_q ) \eta_{\bbeta,k} (\mu_{q+1} t) g_{\bbeta,k}      \nonumber \\ 
            =&\, \frac{\mr \lambda_q}{\mr \lambda_{q+1}} \sum_{k=1}^{k_n} \chi_{\bbeta,k}( \lambda_{q+1} \Phi_q ) \eta_{\bbeta,k} (\mu_{q+1} t) \mr \partial_{q,1} g_{\bbeta,k}
                  + \frac{\mr \lambda_q}{\mr \lambda_{q+1}} \mr \partial_{q,1} g_{\bbeta,0}    \label{e:14_1:decompPrep} \\ 
            +&\, \frac{\varepsilon_q \lambda_{q+1}}{\mr \lambda_{q+1}} \sum_{k=1}^{k_n} \sum_{m} \partial_{\xi_i} \chi_{\bbeta,k}( \lambda_{q+1} \Phi_q ) [ \phi_{m} \eta_{\bbeta,k} ] (\mu_{q+1} t) [B_{q,m}]_{i1} g_{\bbeta,k}     \label{e:14_2:decompPrep} \\ 
            +&\, \frac{\lambda_{q+1}}{\mr \lambda_{q+1}} \sum_{k=1}^{k_n} \partial_{\xi_1} \chi_{\bbeta,k}( \lambda_{q+1} \Phi_q ) \eta_{\bbeta,k} (\mu_{q+1} t) g_{\bbeta,k}.    \label{e:14_3:decompPrep}
\end{align}
Here, we use \eqref{e:6:uEstimate} for computing $\nabla \Phi_q$.

Now we renumber the contributions from (\ref{e:14_1:decompPrep}-\ref{e:14_3:decompPrep}), in the form of \eqref{e:9:decompPrep}. Then
\begin{align}
      k_{n+1} \leq (N+2) k_n,
\end{align}
which verifies \eqref{e:l2:decompPrep} for $\balpha$. We derive the following information for any $p \in \N$
\begin{align}
      \max_{k \geq 1} \| \partial_\tau^p \eta_{\balpha,k} \|_{\infty} \leq&\, \max_{k \geq 1,m} \big\{ \| \partial_\tau^p \eta_{\bbeta,k} \|_{\infty} , \| \partial_\tau^p(\phi_{m} \eta_{\bbeta,k}) \|_{\infty} \big\},   \label{e:14_4:decompPrep} \\ 
      \max_{k \geq 1} \| \nabla_\xi^p \chi_{\balpha,k} \|_{\infty} \leq&\, \max_{k \geq 1} \big\{
            \| \nabla_\xi^p \chi_{\bbeta,k} \|_{\infty}, \| \nabla_\xi^{p+1} \chi_{\bbeta,k} \|_{\infty} \big\},     \label{e:14_5:decompPrep}
\end{align}
which verifies \eqref{e:l3:decompPrep}, \eqref{e:l4:decompPrep} for $\balpha$.

We also have 
\begin{align}
      \vertiii{ g_{\balpha,0} }_q \leq &\,
            \frac{\mr \lambda_q}{\mr \lambda_{q+1}} \vertiii{ \mr \partial_{q,1} g_{\bbeta,0} }_{q},    \label{e:14_6:decompPrep} \\ 
      \max_{k \geq 1} \vertiii{ g_{\balpha,k} }_q \leq &\, \max_{k \geq 1,m} \bigg\{
            \frac{\mr \lambda_q}{\mr \lambda_{q+1}} \vertiii{ \mr \partial_{q,1} g_{\bbeta,k} }_q, 
            \frac{\varepsilon_q \lambda_{q+1}}{\mr \lambda_{q+1}} \vertiii{ B_{q,m} g_{\bbeta,k} }_q, 
            \frac{\lambda_{q+1}}{\mr \lambda_{q+1}} \vertiii{ g_{\bbeta,k} }_q \bigg\},     \label{e:14_7:decompPrep}
\end{align}
verifying \eqref{e:l6:decompPrep}, \eqref{e:l7:decompPrep} and \eqref{e:l10:decompPrep} for $\balpha$. Similarly, \eqref{e:l5:decompPrep} holds for $\balpha$. Here, we use Remark \ref{r:opPolynomial}, \eqref{e:algbBESOmega} and \eqref{e:mtrPolyEst}.

\end{case}

\begin{case}
We consider the case $\balpha = t \bbeta$, then
\begin{equation}        \label{e:15:decompPrep}
\begin{split}
      \mr O_{q+1}^{\balpha} =&\, \mr O_{q+1,t} \mr O_{q+1}^{\bbeta}
            = \mr \mu_{q+1}^{-1} O_{q+1,t} \mr O_{q+1}^{\bbeta} \\ 
            =&\, \mr \mu_{q+1}^{-1} D_{q,t} \mr O_{q+1}^{\bbeta}
                  + \mr \mu_{q+1}^{-1} \nabla^\perp \psi_{q+1} \cdot \nabla \mr O_{q+1}^{\bbeta}.
\end{split}
\end{equation}
Then we can compute
\begin{align}
      \mr O_{q+1}^{\balpha} \Bigg( \frac{\mr \lambda_{q+1}^2} {\delta_{q+1}^{\sfrac12}\lambda_{q+1}^{1+2\gamma_R}} \Psi_{q+1} \Bigg)
      =&\, \mr \mu_{q+1}^{-1} D_{q,t} \sum_{k=1}^{k_n} g_{\bbeta,k} \chi_{\bbeta,k}( \lambda_{q+1} \Phi_q ) \eta_{\bbeta,k} (\mu_{q+1} t)     \label{e:16:decompPrep} \\ 
            +&\, \mr \mu_{q+1}^{-1} \nabla^\perp \psi_{q+1} \cdot \nabla \sum_{k=1}^{k_n} g_{\bbeta,k} \chi_{\bbeta,k}( \lambda_{q+1} \Phi_q ) \eta_{\bbeta,k} (\mu_{q+1} t)         \label{e:18:decompPrep} \\ 
            +&\, \mr \mu_{q+1}^{-1} D_{q,t} g_{\bbeta,0} + \mr \mu_{q+1}^{-1} \nabla^\perp \psi_{q+1} \cdot \nabla g_{\bbeta,0}        \label{e:19:decompPrep} \\ 
            =&\, \mr \mu_{q+1}^{-1} \sum_{k=1}^{k_n} D_{q,t} g_{\bbeta,k} \chi_{\bbeta,k}( \lambda_{q+1} \Phi_q ) \eta_{\bbeta,k} (\mu_{q+1} t)        \label{e:20:decompPrep}  \\ 
            +&\, \lambda_{q+1} \mr \mu_{q+1}^{-1} \sum_{k=1}^{k_n} g_{\bbeta,k} D_{q,t} \Phi_q \cdot \nabla_\xi \chi_{\bbeta,k}( \lambda_{q+1} \Phi_q ) \eta_{\bbeta,k} (\mu_{q+1} t)      \label{e:22:decompPrep}  \\
            +&\, \mu_q \mr \mu_{q+1}^{-1} \sum_{k=1}^{k_n} g_{\bbeta,k} \chi_{\bbeta,k}( \lambda_{q+1} \Phi_q ) \partial_\tau \eta_{\bbeta,k} (\mu_{q+1} t)        \label{e:24:decompPrep}  \\
            +&\, \mr \mu_{q+1}^{-1} \sum_{k=1}^{k_n} \nabla^\perp \psi_{q+1} \cdot \nabla g_{\bbeta,k} \chi_{\bbeta,k}( \lambda_{q+1} \Phi_q ) \eta_{\bbeta,k} (\mu_{q+1} t)   \label{e:26:decompPrep} \\ 
            +&\, \lambda_{q+1} \mr \mu_{q+1}^{-1} \sum_{k=1}^{k_n} g_{\bbeta,k} \nabla^\perp \psi_{q+1} \cdot \nabla \Phi_{q,i} \partial_{\xi_i} \chi_{\bbeta,k}( \lambda_{q+1} \Phi_q ) \eta_{\bbeta,k} (\mu_{q+1} t)         \label{e:28:decompPrep} \\ 
            +&\, \mr \mu_{q+1}^{-1} D_{q,t} g_{\bbeta,0} + \mr \mu_{q+1}^{-1} \nabla^\perp \psi_{q+1} \cdot \nabla g_{\bbeta,0}         \label{e:29:decompPrep} \\ 
            :=&\, J_1 + J_2 + J_3 + J_4 + J_5 + J_6 + J_7.        \nonumber 
\end{align}
We also compute 
\begin{equation}  \label{e:32:decompPrep}
\begin{split}
      \nabla \psi_{q+1} 
            =&\, \delta_{q+1}^{\sfrac12} \lambda_{q+1}^{-1} \sum_{m} \nabla \varpi_{q,m} \zeta_{m} (\mu_{q+1} t) \Big( \eta_1(\mu_{q+1} t) \Pi_1( \lambda_{q+1} \Phi_q ) + \eta_2(\mu_{q+1} t) \Pi_2( \lambda_{q+1} \Phi_q ) \Big) \\ 
            +&\, \delta_{q+1}^{\sfrac12}   
            \Big( \Id + \varepsilon_q \sum_{m} \adj E_{q,m}^T \varphi_{m} (\mu_{q+1} t) \Big) \eta_1(\mu_{q+1} t) \nabla_{\xi} \Pi_1( \lambda_{q+1} \Phi_q ) \\ 
            +&\, \delta_{q+1}^{\sfrac12}  
            \Big( \Id + \varepsilon_q \sum_{m} \adj E_{q,m}^T \varphi_{m} (\mu_{q+1} t) \Big)  \eta_2(\mu_{q+1} t) \nabla_{\xi} \Pi_2( \lambda_{q+1} \Phi_q ) \\ 
            =&\, \delta_{q+1}^{\sfrac12} \lambda_{q+1}^{-1} \sum_{m} \nabla \varpi_{q,m} [\eta_1\zeta_{m}] (\mu_{q+1} t) \Pi_1( \lambda_{q+1} \Phi_q ) \\ 
            +&\, \varepsilon_q \delta_{q+1}^{\sfrac12} \sum_{m} \adj E_{q,m}^T [ \eta_1 \varphi_{m} ] (\mu_{q+1} t) \nabla_{\xi} \Pi_1( \lambda_{q+1} \Phi_q ) \\ 
            +&\, \delta_{q+1}^{\sfrac12} \eta_1 (\mu_{q+1} t) \nabla_{\xi} \Pi_1( \lambda_{q+1} \Phi_q ) \\ 
            +&\, \text{symmetric terms concering } \eta_2, \Pi_2
\end{split}
\end{equation}

The contributions from \eqref{e:20:decompPrep}, \eqref{e:24:decompPrep} and $J_6$ in \eqref{e:29:decompPrep} are respectively
\begin{align}
      J_1 =&\, \frac{\mr \mu_{q}}{\mr \mu_{q+1}} \sum_{k=1}^{k_n} \chi_{\bbeta,k}( \lambda_{q+1} \Phi_q ) \eta_{\bbeta,k} (\mu_{q+1} t) \mr D_{q,t} g_{\bbeta,k},       \label{e:42:decompPrep} \\ 
      J_3 =&\, \frac{\mu_{q+1}}{\mr \mu_{q+1}} \sum_{k=1}^{k_n} \chi_{\bbeta,k}( \lambda_{q+1} \Phi_q ) \partial_\tau \eta_{\bbeta,k} (\mu_{q+1} t) g_{\bbeta,k},         \label{e:43:decompPrep} \\ 
      J_6 =&\, \frac{\mr \mu_{q}}{\mr \mu_{q+1}} \mr D_{q,t} g_{\bbeta,0}.         \label{e:41:decompPrep}
 \end{align} 
Using \eqref{e:DqtPhiq}, the contribution from \eqref{e:22:decompPrep} is
\begin{align}     \label{e:44:decompPrep}
      J_2 = \frac{ \lambda_{q+1} }{ \mr \mu_{q+1} } \sum_{k=1}^{k_n} \omega_{q,N} \cdot \nabla_\xi \chi_{\bbeta,k}( \lambda_{q+1} \Phi_q ) \eta_{\bbeta,k} (\mu_{q+1} t) \phi_{N} (\mu_{q+1} t) g_{\bbeta,k}. 
\end{align}
Using \eqref{e:32:decompPrep} and \eqref{e:adjointPerp}, the contribution from $J_7$ in \eqref{e:29:decompPrep} is
\begin{align}
      J_7 =&\, \frac{ \delta_{q+1}^{\sfrac12} \mr\lambda_q^2 } { \mr \mu_{q+1} \lambda_{q+1} } 
            \sum_{m} [ \eta_1 \zeta_{m} ] (\mu_{q+1} t) \Pi_1( \lambda_{q+1} \Phi_q ) \mr\nabla_q^\perp \varpi_{q,m} \cdot \mr\nabla_q g_{\bbeta,0}        \label{e:45_1:decompPrep} \\ 
            +&\, \frac{ \varepsilon_q \delta_{q+1}^{\sfrac12} \mr \lambda_q } { \mr \mu_{q+1} } 
            \sum_{m}  [ \eta_1 \varphi_{m} ] (\mu_{q+1} t) \nabla_{\xi}^\perp \Pi_1( \lambda_{q+1} \Phi_q ) \cdot \big( \adj E_{q,m}^T \mr \nabla_{q} g_{\bbeta,0} \big)     \label{e:45_2:decompPrep} \\ 
            +&\, \frac{ \delta_{q+1}^{\sfrac12} \mr\lambda_q } {\mr \mu_{q+1}} 
            \eta_1 (\mu_{q+1} t) \nabla_{\xi}^\perp \Pi_1( \lambda_{q+1} \Phi_q ) \cdot \mr \nabla_{q} g_{\bbeta,0}     \label{e:45_3:decompPrep} \\ 
            +&\, \text{symmetric terms concering } \eta_2, \Pi_2         \label{e:45_4:decompPrep}
\end{align}
Using \eqref{e:32:decompPrep}, the contribution from \eqref{e:26:decompPrep} is 
\begin{align}
      J_4 =&\, \mr \mu_{q+1}^{-1} \sum_{k=1}^{k_n} \nabla^\perp \psi_{q+1} \cdot \nabla g_{\bbeta,k} \chi_{\bbeta,k}( \lambda_{q+1} \Phi_q ) \eta_{\bbeta,k} (\mu_{q+1} t)     \nonumber \\ 
      =&\, \frac{ \delta_{q+1}^{\sfrac12} \mr \lambda_q^2 } { \mr \mu_{q+1} \lambda_{q+1} }  
            \sum_{k=1}^{k_n} \sum_{m} [ \eta_1 \zeta_{m} \eta_{\bbeta,k} ] (\mu_{q+1} t) 
                  [ \Pi_1\chi_{\bbeta,k} ] ( \lambda_{q+1} \Phi_q )   
                  \mr \nabla_q^\perp \varpi_{q,m} \cdot \mr \nabla_q g_{\bbeta,k}      \label{e:46:decompPrep} \\ 
      +&\, \frac{\delta_{q+1}^{\sfrac12}\mr \lambda_q}{\mr \mu_{q+1}} 
            \sum_{k=1}^{k_n} [ \eta_1 \eta_{\bbeta,k} ] (\mu_{q+1} t)
            [ \partial_{\xi_1} \Pi_1 \chi_{\bbeta,k} ] ( \lambda_{q+1} \Phi_q )  \mr \partial_{q,1} g_{\bbeta,k}      \label{e:48:decompPrep} \\ 
      +&\, \frac{\varepsilon_q \delta_{q+1}^{\sfrac12}\mr \lambda_q} {\mr \mu_{q+1}} 
            \sum_{k=1}^{k_n} \sum_m [ \eta_1 \varphi_{m} \eta_{\bbeta,k} ] (\mu_{q+1} t) 
            [ \nabla_{\xi}^\perp \Pi_1 \chi_{\bbeta,k} ] ( \lambda_{q+1} \Phi_q ) \cdot \big( \adj E_{q,m}^T \mr \nabla_{q} g_{\bbeta,k} \big)      \label{e:49:decompPrep} \\ 
      +&\, \text{symmetric terms concering } \eta_2, \Pi_2.        \label{e:50:decompPrep}
\end{align}
Using \eqref{e:32:decompPrep} and \eqref{e:6:uEstimate}, the contribution from \eqref{e:28:decompPrep} is
\begin{align}
      J_5 =&\, \frac{\lambda_{q+1}}{\mr \mu_{q+1}} \sum_{k=1}^{k_n} \partial_{\xi_i} \chi_{\bbeta,k}( \lambda_{q+1} \Phi_q ) \eta_{\bbeta,k} (\mu_{q+1} t) g_{\bbeta,k} \nabla^\perp \psi_{q+1} \cdot \nabla \Phi_{q,i}        \nonumber \\ 
      =&\, \frac{ \delta_{q+1}^{\sfrac12} \mr\lambda_{q} }{\mr \mu_{q+1}}  
            \sum_{k=1}^{k_n} \sum_{m} \big[ \Pi_1 \partial_{\xi_i} \chi_{\bbeta,k} \big] ( \lambda_{q+1} \Phi_q ) \big[ \eta_1 \zeta_{m} \eta_{\bbeta,k} \big] (\mu_{q+1} t) g_{\bbeta,k} \mr\partial^\perp_{q,i} \varpi_{q,m}        \label{e:52:decompPrep} \\ 
      +&\, \frac{ \varepsilon_q \delta_{q+1}^{\sfrac12} \mr\lambda_q } {\mr \mu_{q+1}} 
            \sum_{k=1}^{k_n} \sum_{m_1,m_2} \big[ \Pi_1 \partial_{\xi_i} \chi_{\bbeta,k} \big] ( \lambda_{q+1} \Phi_q ) \big[ \eta_1 \zeta_{m_1} \phi_{m_2} \eta_{\bbeta,k} \big] (\mu_{q+1} t) g_{\bbeta,k} \mr\partial^\perp_{q,j} \varpi_{q,m_1} [ B_{q,m_2} ]_{ij}         \label{e:53:decompPrep} \\ 
      +&\, \frac{\delta_{q+1}^{\sfrac12}\lambda_{q+1}}{\mr \mu_{q+1}} 
            \sum_{k=1}^{k_n} 
            [ \partial^\perp_{\xi_i} \Pi_1 \partial_{\xi_i} \chi_{\bbeta,k} ] ( \lambda_{q+1} \Phi_q ) 
            [ \eta_1\eta_{\bbeta,k} ] (\mu_{q+1} t) g_{\bbeta,k}     \label{e:54:decompPrep} \\ 
      +&\, \frac{ \varepsilon_q \delta_{q+1}^{\sfrac12}\lambda_{q+1} } {\mr \mu_{q+1}} 
            \sum_{k=1}^{k_n} \sum_{m_2} 
            [ \partial^\perp_{\xi_j} \Pi_1 \partial_{\xi_i} \chi_{\bbeta,k} ] ( \lambda_{q+1} \Phi_q ) 
            [ \eta_1 \eta_{\bbeta,k} \phi_{m_2} ] (\mu_{q+1} t) 
            g_{\bbeta,k} [B_{q,m_2}]_{ij}       \label{e:55:decompPrep} \\ 
      +&\, \frac{ \varepsilon_q \delta_{q+1}^{\sfrac12} \lambda_{q+1} }{\mr \mu_{q+1}} 
            \sum_{k=1}^{k_n} \sum_{m_1} 
            [ \partial_{\xi_j}^\perp \Pi_1 \partial_{\xi_i} \chi_{\bbeta,k} ] ( \lambda_{q+1} \Phi_q ) 
            [ \eta_1\varphi_{m_1}\eta_{\bbeta,k} ] (\mu_{q+1} t) 
            g_{\bbeta,k} [E_{q,m_1}]_{ij}     \label{e:56:decompPrep} \\ 
      +&\, \frac{ \varepsilon_q^2 \delta_{q+1}^{\sfrac12} \lambda_{q+1} }{\mr \mu_{q+1}} 
            \sum_{k=1}^{k_n} \sum_{m_1,m_2}
            [ \partial_{\xi_l}^\perp \Pi_1 \partial_{\xi_i} \chi_{\bbeta,k} ] ( \lambda_{q+1} \Phi_q ) 
            [ \eta_1\varphi_{m_1} \phi_{m_2}\eta_{\bbeta,k} ] (\mu_{q+1} t) 
            g_{\bbeta,k} [E_{q,m_1}]_{jl} [B_{q,m_2}]_{ij}        \label{e:57:decompPrep} \\ 
      +&\, \text{symmetric terms concering } \eta_2, \Pi_2        \label{e:58:decompPrep}
\end{align}
Here, $\Phi_{q,i}$ refers to the $i$-th component of the flowmap $\Phi_q$.

Now we renumber all contributions from (\ref{e:20:decompPrep}-\ref{e:28:decompPrep}), computed in (\ref{e:42:decompPrep}-\ref{e:58:decompPrep}), in the form of \eqref{e:9:decompPrep}. Simply by counting all terms, we have
\begin{align*}
      k_{n+1} \leq 3N^2 k_n,
\end{align*}
verifying \eqref{e:l2:decompPrep}.

For the functions $\eta_{\balpha,k}$, we have 
\begin{align}     
      \max_{k \geq 1} \| \partial_\tau^p \eta_{\balpha,k} \|_{\infty} \leq&\, \max_{k \geq 1,m,m_1,m_2} 
      \big\{ \| \partial_\tau^p \eta_{\bbeta,k} \|_{\infty} , \| \partial_\tau^p (\eta_{\bbeta,k}\phi_N) \|_{\infty} , \| \partial_\tau^{p+1} \eta_{\bbeta,k} \|_{\infty} ,   \label{e:64:decompPrep} \\
            &\,\| \partial_\tau^p \eta_1 \|_{\infty}, \| \partial_\tau^p (\eta_1\zeta_m) \|_{\infty} , \| \partial_\tau^p (\eta_1\varphi_m) \|_{\infty} ,         \label{e:65:decompPrep}  \\ 
            &\, \| \partial_\tau^p (\eta_1 \eta_{\bbeta,k}) \|_{\infty} ,
            \| \partial_\tau^p (\eta_1 \zeta_m \eta_{\bbeta,k}) \|_{\infty}, 
            \|\partial_\tau^p (\eta_1 \varphi_m \eta_{\bbeta,k}) \|_{\infty} , 
            \|\partial_\tau^p (\eta_1 \phi_m \eta_{\bbeta,k}) \|_{\infty} ,   \label{e:66:decompPrep} \\ 
            &\, \| \partial_\tau^p (\eta_1 \zeta_{m_1} \phi_{m_2} \eta_{\bbeta,k}) \|_{\infty} , \| \partial_\tau^p (\eta_1 \varphi_{m_1} \phi_{m_2} \eta_{\bbeta,k}) \|_{\infty},       \label{e:67:decompPrep} \\ 
            &\, \text{symmetric terms concering } \eta_2  \big\}.       \label{e:68:decompPrep}
\end{align}
Here, \eqref{e:64:decompPrep} and \eqref{e:65:decompPrep} contain the contributions from (\ref{e:42:decompPrep},\ref{e:43:decompPrep},\ref{e:44:decompPrep}) and (\ref{e:45_1:decompPrep}-\ref{e:45_4:decompPrep}), respectively. Similarly, \eqref{e:66:decompPrep} and \eqref{e:67:decompPrep} contain the contributions from (\ref{e:46:decompPrep}-\ref{e:58:decompPrep}). From induction assumptions on $\bbeta$, (\ref{e:64:decompPrep}-\ref{e:68:decompPrep}) and \eqref{e:phiZetaThetaEst}, we can verify \eqref{e:l4:decompPrep} for $\balpha$.

Similar to the analysis on $\eta_{\balpha,k}$, we have, for the functions $\chi_{\balpha,k}$,
\begin{equation}        \label{e:81:decompPrep}
\begin{split}
      \max_{k \geq 1} \| \nabla_\xi^p \chi_{\balpha,k} \|_{\infty} \leq \max_{k \geq 1} 
      \big\{&\, \| \nabla_\xi^p \chi_{\bbeta,k} \|_{\infty} , \| \nabla_\xi^{p+1} \chi_{\bbeta,k} \|_{\infty} , \|\nabla_\xi^p (\Pi_1\chi_{\bbeta,k}) \|_{\infty} , \| \nabla_\xi^p (\nabla_\xi \Pi_1\chi_{\bbeta,k}) \|_{\infty}, \\
            &\,\|\nabla_\xi^p ( \Pi_1 \nabla_\xi \chi_{\bbeta,k}) \|_{\infty} , \|\nabla_\xi^p ( \nabla_\xi \Pi_1 \nabla_\xi \chi_{\bbeta,k}) \|_{\infty} , \\ 
            &\, \text{symmetric terms concering } \Pi_2 \big\}.
\end{split}
\end{equation}
Similarly to $\eta_{\balpha,k}$, we can verify \eqref{e:l3:decompPrep}.

For the functions $g_{\balpha,0}$, we have the contribution from $J_6$ giving the following
\begin{align}     \label{e:84:decompPrep}
      \vertiii{ g_{\balpha,0} }_q \leq 
      \frac{\mr \mu_{q}}{\mr \mu_{q+1}} \vertiii{ \mr D_{q,t} g_{\bbeta,0} }_q,
\end{align}
verifying \eqref{e:l6:decompPrep} and \eqref{e:l12:decompPrep} for $\balpha$.

For the functions $g_{\balpha,k}$, we have for any 
\begin{equation}        \label{e:90:decompPrep}
\begin{split}
      \max_{k \geq 1} \vertiii{ g_{\balpha,k} }_q 
            \leq&\, \max_{k \geq 1,m,m_1,m_2} \bigg\{ 
      \frac{\mr \mu_{q}}{\mr \mu_{q+1}} \vertiii{ \mr D_{q,t} g_{\bbeta,k} }_q,
      \frac{\mu_{q+1}}{\mr \mu_{q+1}} \vertiii{ g_{\bbeta,k} }_q,
      \frac{\lambda_q}{\mr \mu_{q+1}} \vertiii{ \omega_{q,N} g_{\bbeta,k} }_q,
            \\ 
      &\, \frac{\delta_{q+1}^{\sfrac12}}{\lambda_{q+1}} \frac{\mr \lambda_q^2}{\mr \mu_{q+1}} \vertiii{ \mr \nabla_q^\perp \varpi_{q,m} \cdot \mr \nabla_q g_{\bbeta,k}  }_q,
      \frac{ \varepsilon_q \delta_{q+1}^{\sfrac12} \mr \lambda_q } {\mr \mu_{q+1}} \vertiii{ \adj E^T_{q,m} \mr \nabla_{q} g_{\bbeta,k}  }_q,
      \frac{\delta_{q+1}^{\sfrac12}\mr \lambda_q}{\mr \mu_{q+1}} \vertiii{ \mr \nabla_{q} g_{\bbeta,k}  }_q, \\ 
      &\, \frac{\delta_{q+1}^{\sfrac12}\mr\lambda_q} {\mr \mu_{q+1}} \vertiii{ g_{\bbeta,k} \mr \nabla^\perp_q \varpi_{q,m} }_q, 
      \frac{\varepsilon_q \delta_{q+1}^{\sfrac12} \mr \lambda_q}{\mr \mu_{q+1}} \vertiii{ g_{\bbeta,k} \mr\nabla^\perp_q \varpi_{q,m_1} B_{q,m_2} }_q ,
       \\ 
      &\, \frac{\delta_{q+1}^{\sfrac12}\lambda_{q+1}}{\mr \mu_{q+1}} \vertiii{ g_{\bbeta,k} }_q,
      \frac{ \varepsilon_q \delta_{q+1}^{\sfrac12}\lambda_{q+1} } { \mr \mu_{q+1} } \vertiii{ g_{\bbeta,k} B_{q,m_2} }_q,
      \frac{ \varepsilon_q \delta_{q+1}^{\sfrac12}\lambda_{q+1} } { \mr \mu_{q+1} } \vertiii{ g_{\bbeta,k} E_{q,m_1} }_q, \\ 
      &\,\frac{ \varepsilon_q^2 \delta_{q+1}^{\sfrac12}\lambda_{q+1} } { \mr \mu_{q+1} } \vertiii{ g_{\bbeta,k} B_{q,m_2} E_{q,m_1} }_q  \bigg\},  \\ 
      +&\, \max_{m} \bigg\{ 
      \frac{\mr \mu_{q}}{\mr \mu_{q+1}} \vertiii{ \mr D_{q,t} g_{\bbeta,0} }_q, 
      \frac{\delta_{q+1}^{\sfrac12}}{\lambda_{q+1}} \frac{\mr \lambda_q^2}{\mr \mu_{q+1}} \vertiii{ \mr \nabla_q^\perp \varpi_{q,m} \cdot \mr \nabla_q g_{\bbeta,0} }_q, \\
      &\, \frac{ \varepsilon_q \delta_{q+1}^{\sfrac12} \mr \lambda_q } {\mr \mu_{q+1}} \vertiii{ \adj E^T_{q,m} \mr \nabla_{q} g_{\bbeta,0} }_q,
      \frac{\delta_{q+1}^{\sfrac12}\mr \lambda_q}{\mr \mu_{q+1}} \vertiii{ \mr \nabla_{q} g_{\bbeta,0} }_q
      \bigg\}
\end{split}
\end{equation}
The terms in \eqref{e:90:decompPrep} for $k \geq 1$ are respectively from \eqref{e:42:decompPrep}, \eqref{e:43:decompPrep}, \eqref{e:44:decompPrep}, (\ref{e:46:decompPrep}), (\ref{e:49:decompPrep}), (\ref{e:48:decompPrep}), (\ref{e:52:decompPrep}), (\ref{e:53:decompPrep}), \eqref{e:54:decompPrep}, \eqref{e:55:decompPrep}, \eqref{e:56:decompPrep}. The terms in \eqref{e:90:decompPrep} for $k = 0$ are respectively from \eqref{e:41:decompPrep}, \eqref{e:45_1:decompPrep}, \eqref{e:45_2:decompPrep} and \eqref{e:45_3:decompPrep}.

Similar to the computation in \eqref{e:90:decompPrep}, we can deduce \eqref{e:l5:decompPrep}. Here we apply Remark \ref{r:opPolynomial} and use (\ref{e:algbVarpiOmega}-\ref{e:algbBESOmega}) and induction assumptions on $g_{\bbeta,k}$.

Now we apply Remark \ref{r:opPolynomial}. From \eqref{e:90:decompPrep}, the $\mathcal{P}_q$ bounds on $B_{q,m}, E_{q,m}, \varpi_{q,m}, \omega_{q,N}$ in (\ref{e:mtrPolyEst}-\ref{e:trsptRsdPolyEst}), and the relations in Remark \eqref{r:auxiliaryParaRela}
we have
\begin{align}
      \max_{k \geq 1} \vertiii{ g_{\balpha,k}  }_q \leq&\,
            \max\bigg\{ \frac{\mu_{q+1}}{\mr \mu_{q+1}},
            \frac{\delta_{q+1}^{\sfrac12}\lambda_{q+1}}{\mr \mu_{q+1}} \bigg\}
                  \max_{k \geq 1} \vertiii{ g_{\bbeta,k} }_q, \\ 
      \max_{k \geq 1} \vertiii{ g_{\balpha,k}  }_q \leq&\,
            \max \bigg\{ \frac{\mr \mu_{q}}{\mr \mu_{q+1}},
            \frac{ \delta_{q+1}^{\sfrac12} \mr \lambda_{q} } {\mr \mu_{q+1}} \bigg\}
                  \vertiii{ g_{\bbeta,0} }_q.
\end{align}
This proves \eqref{e:l7:decompPrep} and \eqref{e:l11:decompPrep}.

\end{case}
\end{proof}

\begin{proof}[Proof of Lemma \ref{l:decompRho}]
We prove this by induction. For $\balpha = 0$, the decomposition is \eqref{e:0:decompRho} itself.

Let $n := |\bbeta|$. Assume we have, for some $\{ \Gamma(n,l) \}_{1 \leq l \leq n} \subset \N$, that
\begin{align}     \label{e:12:decompRho}
      \mr O_{q+1}^{\bbeta} \Big( h \chi ( \lambda_{q+1} \Phi_q ) \eta (\mu_{q+1} t) \rho \Big)
            = \sum_{|\bomega| \leq n} \mr D_{q}^{\bomega} \rho \sum_{k=1}^{\Gamma(n,|\bomega|)} h_{\bbeta,k}^{(\bomega)} \chi_{\bbeta,k}^{(\bomega)} ( \lambda_{q+1} \Phi_q ) \eta_{\bbeta,k}^{(\bomega)} (\mu_{q+1} t).
\end{align}
We also assume the estimates (\ref{e:2:decompRho}-\ref{e:6:decompRho}) for $\mr D_{q+1}^{\bbeta}$.

For $\balpha = t\beta$, $\balpha = 1\beta$ or $\balpha = 2\beta$, the goal is to derive estimates (\ref{e:2:decompRho}-\ref{e:6:decompRho}) for $\mr O_{q+1}^{\balpha}$ in the following decomposition
\begin{align}     \label{e:14:decompRho}
      \mr O_{q+1}^{\balpha} \Big( h \chi ( \lambda_{q+1} \Phi_q ) \eta (\mu_{q+1} t) \rho \Big)
            = \sum_{|\bomega| \leq n+1} \mr D_{q}^{\bomega} \rho \sum_{k=1}^{\Gamma(n+1,|\bomega|)} h_{\balpha,k}^{(\bomega)} \chi_{\balpha,k}^{(\bomega)} ( \lambda_{q+1} \Phi_q ) \eta_{\balpha,k}^{(\bomega)} (\mu_{q+1} t)
\end{align}
for some $\{ \Gamma(n+1,l) \}_{1 \leq l \leq n+1} \subset \N$.

We denote 
\begin{align}     \label{e:15:decompRho}
      r(n,i) := \sum_{|\bbeta|=n} \sum_{|\bomega|=i} \sum_{k=1} \vertiii{ h_{\bbeta,k}^{(\bomega)} }_q. 
\end{align}

\begin{case}

Consider $\balpha = 1 \bbeta$. From \eqref{e:12:decompPrep}, we have
\begin{align}
      \mr O_{q+1}^{\balpha} g =&\, \mr \lambda_{q+1}^{-1} \partial_1 \bigg( \sum_{|\bomega| \leq n} \mr D_{q}^{\bomega} \rho \sum_k h_{\bbeta,k}^{(\bomega)} \chi_{\bbeta,k}^{(\bomega)} ( \lambda_{q+1} \Phi_q ) \eta_{\bbeta,k}^{(\bomega)} (\mu_{q+1} t) \bigg) \\ 
      =&\, \frac{\mr \lambda_q}{\mr \lambda_{q+1}} \sum_{|\bomega| \leq n} \mr\partial_{q,1} \mr D_{q}^{\bomega} \rho \sum_k h_{\bbeta,k}^{(\bomega)} \chi_{\bbeta,k}^{(\bomega)} ( \lambda_{q+1} \Phi_q ) \eta_{\bbeta,k}^{(\bomega)} (\mu_{q+1} t)   \label{e:16:decompRho} \\ 
      +&\, \frac{\mr \lambda_q}{\mr \lambda_{q+1}} \sum_{|\bomega| \leq n} \mr D_{q}^{\bomega} \rho \sum_k \mr\partial_{q,1} h_{\bbeta,k}^{(\bomega)} \chi_{\bbeta,k}^{(\bomega)} ( \lambda_{q+1} \Phi_q ) \eta_{\bbeta,k}^{(\bomega)} (\mu_{q+1} t)   \label{e:18:decompRho} \\ 
      +&\, \frac{\lambda_{q+1}}{\mr \lambda_{q+1}} \sum_{|\bomega| \leq n} \mr D_{q}^{\bomega} \rho \sum_k h_{\bbeta,k}^{(\bomega)} \partial_{\xi_1} \chi_{\bbeta,k}^{(\bomega)} ( \lambda_{q+1} \Phi_q ) \eta_{\bbeta,k}^{(\bomega)} (\mu_{q+1} t)    \label{e:20:decompRho} \\ 
      +&\, \frac{\varepsilon_q \lambda_{q+1}}{\mr \lambda_{q+1}} \sum_{|\bomega| \leq n} \mr D_{q}^{\bomega} \rho \sum_{m,k} [B_{q,m}]_{i1} h_{\bbeta,k}^{(\bomega)} \partial_{\xi_i} \chi_{\bbeta,k}^{(\bomega)} ( \lambda_{q+1} \Phi_q ) \big[ \phi_m \eta_{\bbeta,k}^{(\bomega)} \big] (\mu_{q+1} t)   \label{e:22:decompRho}
\end{align}

Now we need to renumber the terms above and write them in the form of \eqref{e:14:decompRho}.

For the contributions from (\ref{e:16:decompRho}-\ref{e:22:decompRho}), we have
\begin{align}     \label{e:26:decompRho}
      r ( n+1, |\bomega| ) \lesssim 
            \frac{\lambda_{q+1}}{\mr \lambda_{q+1}} r (n, |\bomega|)
            + \frac{\mr \lambda_q}{\mr \lambda_{q+1}} r ( n, |\bomega|-1) .
\end{align}

To estimate $\Gamma(n+1,|\bomega|)$, we combine above computations and deduce
\begin{align}     \label{e:28:decompRho}
      \Gamma(n+1,|\bomega|) \leq 3N \Gamma(n,|\bomega|) + \Gamma(n,|\bomega|-1).
\end{align}

\end{case}

\begin{case}
Consider $\balpha = t \bbeta$. From \eqref{e:15:decompPrep}, we have (below $k$ is summed from $1$ to $\Gamma(n,|\bomega|)$ and $m,m_1,m_2$ are summed from $1$ to $N$)
\begin{align}
      \mr O_{q+1}^{\balpha} g =&\, \frac{1}{\mr \mu_{q+1}} D_{q,t} \bigg( \sum_{|\bomega| \leq n} \mr D_{q}^{\bomega} \rho \sum_k h_{\bbeta,k}^{(\bomega)} \chi_{\bbeta,k}^{(\bomega)} ( \lambda_{q+1} \Phi_q ) \eta_{\bbeta,k}^{(\bomega)} (\mu_{q+1} t) \bigg)      \nonumber \\ 
      +&\, \frac{1}{\mr \mu_{q+1}} \nabla^\perp \psi_{q+1} \cdot \nabla \bigg( \sum_{|\bomega| \leq n} \mr D_{q}^{\bomega} \rho \sum_k h_{\bbeta,k}^{(\bomega)} \chi_{\bbeta,k}^{(\bomega)} ( \lambda_{q+1} \Phi_q ) \eta_{\bbeta,k}^{(\bomega)} (\mu_{q+1} t) \bigg)       \label{e:31:decompRho} \\ 
      =&\, \frac{\mr \mu_q}{\mr \mu_{q+1}} \sum_{|\bomega| \leq n} \mr D_{q,t} \mr D_{q}^{\bomega} \rho \sum_k h_{\bbeta,k}^{(\bomega)} \chi_{\bbeta,k}^{(\bomega)} ( \lambda_{q+1} \Phi_q ) \eta_{\bbeta,k}^{(\bomega)} (\mu_{q+1} t)        \label{e:32:decompRho} \\ 
      +&\, \frac{\mr \mu_q}{\mr \mu_{q+1}} \sum_{|\bomega| \leq n} \mr D_{q}^{\bomega} \rho \sum_k \mr D_{q,t} h_{\bbeta,k}^{(\bomega)} \chi_{\bbeta,k}^{(\bomega)} ( \lambda_{q+1} \Phi_q ) \eta_{\bbeta,k}^{(\bomega)} (\mu_{q+1} t)        \label{e:33:decompRho} \\ 
      +&\, \frac{\lambda_{q+1}}{\mr \mu_{q+1}} \sum_{|\bomega| \leq n} \mr D_{q}^{\bomega} \rho \sum_k h_{\bbeta,k}^{(\bomega)} \omega_{q,N} \cdot \nabla_{\xi} \chi_{\bbeta,k}^{(\bomega)} ( \lambda_{q+1} \Phi_q ) \big[ \phi_N \eta_{\bbeta,k}^{(\bomega)} \big] (\mu_{q+1} t)       \label{e:34:decompRho} \\ 
      +&\, \frac{\mu_q}{\mr \mu_{q+1}} \sum_{|\bomega| \leq n} \mr D_{q}^{\bomega} \rho \sum_k h_{\bbeta,k}^{(\bomega)} \chi_{\bbeta,k}^{(\bomega)} ( \lambda_{q+1} \Phi_q ) \partial_\tau \eta_{\bbeta,k}^{(\bomega)} (\mu_{q+1} t)          \label{e:35:decompRho} \\ 
      +&\, J_1 + J_2.         \label{e:36:decompRho}
\end{align}
Here, $J_1$ and $J_2$ are defined to contain the contributions from \eqref{e:31:decompRho}. Specifically, $J_1$ contains the terms when $\nabla$ hits $\mr D_{q}^{\bomega} \rho$ and $h_{\bbeta,k}^{(\bomega)}$, and $J_2$ contains the terms when $\nabla$ hits $\chi_{\bbeta,k}^{(\bomega)} ( \lambda_{q+1} \Phi_q )$.

Recalling \eqref{e:32:decompPrep}, we have
\begin{align}
      J_1 =&\, \frac{\mr \lambda_q}{\mr \mu_{q+1}} \sum_{|\bomega| \leq n} \sum_k \nabla^\perp \psi_{q+1} \cdot \mr\nabla_q \big( h_{\bbeta,k}^{(\bomega)} \mr D_{q}^{\bomega} \rho \big) \chi_{\bbeta,k}^{(\bomega)} ( \lambda_{q+1} \Phi_q ) \eta_{\bbeta,k}^{(\bomega)} (\mu_{q+1} t)      \nonumber \\ 
      =&\, \frac{\delta_{q+1}^{\sfrac12} \mr\lambda_q^2} {\mr \mu_{q+1}\lambda_{q+1}} \sum_{|\bomega| \leq n} \mr\nabla_q \mr D_{q}^{\bomega} \rho \cdot \sum_{m,k} \mr \nabla_q^\perp \varpi_{q,m} h_{\bbeta,k}^{(\bomega)} \big[ \chi_{\bbeta,k}^{(\bomega)} \Pi_1 \big] ( \lambda_{q+1} \Phi_q ) \big[ \zeta_m \eta_1 \eta_{\bbeta,k}^{(\bomega)} \big] (\mu_{q+1} t)    \label{e:38:decompRho} \\ 
      +&\, \frac{\delta_{q+1}^{\sfrac12} \mr\lambda_q^2} {\mr \mu_{q+1}\lambda_{q+1}} \sum_{|\bomega| \leq n} \mr D_{q}^{\bomega} \rho \sum_{m,k} \mr \nabla_q^\perp \varpi_{q,m} \cdot \mr\nabla_q h_{\bbeta,k}^{(\bomega)} \big[ \chi_{\bbeta,k}^{(\bomega)} \Pi_1 \big] ( \lambda_{q+1} \Phi_q ) \big[ \zeta_m \eta_1 \eta_{\bbeta,k}^{(\bomega)} \big] (\mu_{q+1} t)    \label{e:39:decompRho} \\ 
      +&\, \frac{\delta_{q+1}^{\sfrac12} \mr\lambda_q}{\mr \mu_{q+1}} \sum_{|\bomega| \leq n} \mr\nabla_q \mr D_{q}^{\bomega} \rho \cdot \sum_k h_{\bbeta,k}^{(\bomega)} \big[ \chi_{\bbeta,k}^{(\bomega)} \nabla_{\xi}^\perp \Pi_1 \big] ( \lambda_{q+1} \Phi_q ) \big[ \eta_1 \eta_{\bbeta,k}^{(\bomega)} \big] (\mu_{q+1} t)         \label{e:40:decompRho} \\ 
      +&\, \frac{\delta_{q+1}^{\sfrac12} \mr\lambda_q}{\mr \mu_{q+1}} \sum_{|\bomega| \leq n} \mr D_{q}^{\bomega} \rho \sum_k \mr\nabla_q h_{\bbeta,k}^{(\bomega)} \cdot \big[ \chi_{\bbeta,k}^{(\bomega)} \nabla_{\xi}^\perp \Pi_1 \big] ( \lambda_{q+1} \Phi_q ) \big[ \eta_1 \eta_{\bbeta,k}^{(\bomega)} \big] (\mu_{q+1} t)         \label{e:41:decompRho} \\ 
      +&\, \frac{\varepsilon_q \delta_{q+1}^{\sfrac12} \mr\lambda_q}{\mr \mu_{q+1}} \sum_{|\bomega| \leq n} \mr\partial_{q,j} \mr D_{q}^{\bomega} \rho \sum_{m,k} [ E_{q,m} ]_{ji} h_{\bbeta,k}^{(\bomega)} \big[ \chi_{\bbeta,k}^{(\bomega)} \partial_{\xi_i}^\perp \Pi_1 \big] ( \lambda_{q+1} \Phi_q ) \big[ \eta_1 \eta_{\bbeta,k}^{(\bomega)} \big] (\mu_{q+1} t)       \label{e:42:decompRho} \\ 
      +&\, \frac{\varepsilon_q \delta_{q+1}^{\sfrac12} \mr\lambda_q}{\mr \mu_{q+1}} \sum_{|\bomega| \leq n} \mr D_{q}^{\bomega} \rho \sum_{m,k} [ E_{q,m} ]_{ji} \mr\partial_{q,j} h_{\bbeta,k}^{(\bomega)} \big[ \chi_{\bbeta,k}^{(\bomega)} \partial_{\xi_i}^\perp \Pi_1 \big] ( \lambda_{q+1} \Phi_q ) \big[ \eta_1 \eta_{\bbeta,k}^{(\bomega)} \big] (\mu_{q+1} t)       \label{e:43:decompRho} \\ 
      +&\, \text{symmetric terms concering } \eta_2, \Pi_2     \label{e:44:decompRho}
\end{align}
and
\begin{align}
      J_2 =&\, \frac{\delta_{q+1}^{\sfrac12} \mr\lambda_q}{\mr \mu_{q+1}} 
                  \sum_{|\bomega| \leq n} \mr D_{q}^{\bomega} \rho 
                  \sum_{m,k} h_{\bbeta,k}^{(\bomega)} \mr \nabla_q^\perp \varpi_{q,m} \cdot \big[ \nabla_{\xi} \chi_{\bbeta,k}^{(\bomega)} \Pi_1 \big] ( \lambda_{q+1} \Phi_q ) \big[ \zeta_m \eta_1 \eta_{\bbeta,k}^{(\bomega)} \big] (\mu_{q+1} t)     \label{e:46:decompRho} \\ 
      +&\, \frac{\varepsilon_q\delta_{q+1}^{\sfrac12} \mr\lambda_q}{\mr \mu_{q+1}} 
            \sum_{|\bomega| \leq n} \mr D_{q}^{\bomega} \rho \sum_{m_1,m_2,k} \mr \partial_{q,j}^\perp \varpi_{q,m_1} [B_{q,m_2}]_{ij} h_{\bbeta,k}^{(\bomega)} \big[ \Pi_1 \partial_{\xi_i} \chi_{\bbeta,k}^{(\bomega)} \big] ( \lambda_{q+1} \Phi_q ) \big[ \zeta_{m_1} \phi_{m_2} \eta_1 \eta_{\bbeta,k}^{(\bomega)} \big] (\mu_{q+1} t)         \label{e:47:decompRho} \\ 
      +&\, \frac{\delta_{q+1}^{\sfrac12} \lambda_{q+1}}{\mr \mu_{q+1}} 
            \sum_{|\bomega| \leq n} \mr D_{q}^{\bomega} \rho \sum_k h_{\bbeta,k}^{(\bomega)} \big[ \nabla_{\xi} \chi_{\bbeta,k}^{(\bomega)} \cdot \nabla_{\xi}^\perp \Pi_1 \big] ( \lambda_{q+1} \Phi_q ) \big[ \eta_1 \eta_{\bbeta,k}^{(\bomega)} \big] (\mu_{q+1} t)      \label{e:48:decompRho} \\ 
      +&\, \frac{\varepsilon_q\delta_{q+1}^{\sfrac12} \lambda_{q+1}}{\mr \mu_{q+1}} 
            \sum_{|\bomega| \leq n} \mr D_{q}^{\bomega} \rho \sum_{m,k} [B_{q,m}]_{ij} h_{\bbeta,k}^{(\bomega)} \big[ \partial_{\xi_i} \chi_{\bbeta,k}^{(\bomega)} \partial_{\xi_j}^\perp \Pi_1 \big] ( \lambda_{q+1} \Phi_q ) \big[ \phi_m \eta_1 \eta_{\bbeta,k}^{(\bomega)} \big] (\mu_{q+1} t)        \label{e:49:decompRho} \\ 
      +&\, \frac{\varepsilon_q\delta_{q+1}^{\sfrac12} \lambda_{q+1}}{\mr \mu_{q+1}} 
            \sum_{|\bomega| \leq n} \mr D_{q}^{\bomega} \rho \sum_{m,k} [E_{q,m}]_{ij} h_{\bbeta,k}^{(\bomega)} \big[ \partial_{\xi_i} \chi_{\bbeta,k}^{(\bomega)} \partial_{\xi_j}^\perp \Pi_1 \big] ( \lambda_{q+1} \Phi_q ) \big[ \varphi_m \eta_1 \eta_{\bbeta,k}^{(\bomega)} \big] (\mu_{q+1} t)      \label{e:50:decompRho} \\ 
      +&\, \frac{\varepsilon_q^2\delta_{q+1}^{\sfrac12} \lambda_{q+1}}{\mr \mu_{q+1}} 
            \sum_{|\bomega| \leq n} \mr D_{q}^{\bomega} \rho \sum_{m_1,m_2,k} [B_{q,m_2} E_{q,m_1}]_{ij} h_{\bbeta,k}^{(\bomega)} \big[ \partial_{\xi_i} \chi_{\bbeta,k}^{(\bomega)} \partial_{\xi_j}^\perp \Pi_1 \big] ( \lambda_{q+1} \Phi_q ) \big[ \varphi_{m_1} \phi_{m_2} \eta_1 \eta_{\bbeta,k}^{(\bomega)} \big] (\mu_{q+1} t)       \label{e:51:decompRho} \\ 
      +&\, \text{symmetric terms concerning } \eta_2, \Pi_2.     \label{e:52:decompRho}
\end{align}

For the contributions from (\ref{e:32:decompRho}-\ref{e:52:decompRho}), we have
\begin{align}     \label{e:54:decompRho}
      r ( n+1, |\bomega| ) \lesssim
            \frac{\delta_{q+1}^{\sfrac12} \lambda_{q+1}}{\mr \mu_{q+1}}  r (n, |\bomega|), 
            + \frac{ \delta_{q+1}^{\sfrac12} \mr\lambda_q + \mr \mu_{q} }{ \mr \mu_{q+1} } 
                  r ( n, |\bomega|-1).
\end{align}
Also, we have
\begin{align}     \label{e:56:decompRho}
      \Gamma(n+1,|\bomega|) \leq 3N^2 \Gamma(n,|\bomega|) + 4N \Gamma(n,|\bomega|-1)
\end{align}

\end{case}

From \eqref{e:26:decompRho}, \eqref{e:28:decompRho}, \eqref{e:54:decompRho}, \eqref{e:56:decompRho}, we deduce \eqref{e:2:decompRho} and \eqref{e:4:decompRho}. \eqref{e:3:decompRho}, \eqref{e:5:decompRho} and \eqref{e:6:decompRho} are proved in the same way as Lemma \ref{l:decompPrep}.

When $\chi$ is a constant function, we do not get fast spatial scale $\lambda_{q+1}$. Therefore, \eqref{e:8:decompRho} is a direct consequence of this observation.
\end{proof}

\newpage

\section{Energy-type estimates}
As before, we will write $D_u = \partial_t + u \cdot \nabla $. By
$ D^{\bbeta}$ ($D^{\balpha}$ etc) we denote generalised derivative, it has $|\bbeta|_x$-order of space derivatives and it has $|\bbeta|_t$-many advective derivatives $D_u$, all arranged according to $\bbeta$.

Consider 
\begin{equation}     \label{e:parabolic_const_gen}
\begin{split}
      D_u \rho - \kappa \Delta \rho =&  \divr f, \\ 
      \rho(\cdot,0) =&\, \rho_{\ini}.
\end{split}
\end{equation}
In this section we provide energy-type estimates for $ D^{\bbeta} \rho$. 

\subsection{Scales assumption and how to apply them}
We make the following 'scales assumption' on the drift $u$:
\begin{assumption}\label{ass:scales_foor_ener}
There are three real numbers $\mathcal{A}, \mathcal{S}, \mathcal{T}$ larger than $1$, such that for any  $|\balpha| \le 2 U$, where $U$ is a fixed natural number, it holds 
\begin{equation}\label{e:scales_u1}
\left( \frac{\|D^{\balpha} u\|_{L^\infty}}{ \kappa}  \right)^2 \le  \mathcal{A} \mathcal{S} \mathcal{S}^{|\balpha|_x} \mathcal{T}^{|\balpha|_t}
\end{equation}
\end{assumption}
One may think of $\mathcal{A} \mathcal{S}$ above as 'amplitude', $\mathcal{S}$ as 'cost of space derivative', and  $\mathcal{T}$ as 'cost of advective derivative'.

The general results of this section, i.e.\ Propositions \ref{prop:energy_est}, \ref{prop:energy_est_tech} will specified in section \ref{sec:energ_spec} to the setting of the rest of this paper. More precisely, in the rest of the paper (other chapters), we invoke the energy-type estimates for $u:=u_q$, $\kappa:=\kappa_q$ and with:
\begin{equation}\label{e:scales_transl_ener}
\begin{aligned}
    \mathcal{A} &:=C_U  \\
    \mathcal{S} &:=
     \delta_q^{\sfrac{1}{2}} \lambda_q \kappa^{-1}_q \\
         \mathcal{T} &:=
    (\delta_q^{\sfrac{1}{2}} \lambda_q)^2 
\end{aligned}
\end{equation}
\subsubsection{Consistency of our scales assumptions}
The introduced scales assumption are consistent with the rest of the paper, in view of
\begin{lemma}\label{rem:scales_consistent}
      The velocities $u_q$, $v_q$ of Section \ref{ss:r_estimates} enjoy Assumption \ref{ass:scales_foor_ener} with the choice \eqref{e:scales_transl_ener}, with $U= \lfloor \frac{Y-N_*-1}{2} \rfloor$.
\end{lemma}
\begin{proof}
In view of \eqref{e:scales_transl_ener},  \eqref{e:scales_u1} is
\begin{equation}\label{e:scales_connect}
\left( \frac{\|D^{\balpha} u\|_{L^\infty}}{ \kappa}  \right)^2 \le C_U   \left(\delta_q^{\sfrac{1}{2}} \lambda_q \kappa^{-1}_q \right)^{|\balpha|_x+1}  (\delta_q^{\sfrac{1}{2}} \lambda_q)^{2|\balpha|_t} 
\end{equation}
Estimates of Section \ref{ss:r_estimates} yield the upper bound 
\begin{equation}
\left( \frac{\|D^{\balpha} u\|_{L^\infty}}{ \kappa}  \right)^2 \lesssim   (\delta_q^{\sfrac{1}{2}}   \lambda^{|\balpha|_x}_q \kappa^{-1}_q)^2  (\delta_q^{\sfrac{1}{2}} \lambda_q)^{2|\balpha|_t} 
\end{equation}
which says that with choices \eqref{e:scales_transl_ener} in order to have \eqref{e:scales_u1} or equivalently \eqref{e:scales_connect}, it suffices to have 
\begin{equation}\label{e:scales_true}
  (\delta_q^{\sfrac{1}{2}}   \lambda^{|\balpha|_x}_q \kappa^{-1}_q)^2  (\delta_q^{\sfrac{1}{2}} \lambda_q)^{2|\balpha|_t} \lesssim \left(\delta_q^{\sfrac{1}{2}} \lambda_q \kappa^{-1}_q \right)^{|\balpha|_x+1} (\delta_q^{\sfrac{1}{2}} \lambda_q)^{2|\balpha|_t}  
\end{equation}
which is, taking into account definitions of $\delta_q$ and $\kappa_q$
\[
  \lambda^{-\beta}_q \lambda_q \lambda^\theta_q \ge (  \lambda^{-\beta}_q  \lambda^n_q \lambda^\theta_q)^\frac{2}{n+1}
\]
this last inequality holds iff $\theta \ge \beta +1$.
This, in view of \eqref{e:Exponent} which defines $\theta = \frac{2b}{b+1}(1+\beta)$, is always satisfied.     
\end{proof}
\begin{remark}[Informal version of Assumption \ref{ass:scales_foor_ener} with the choice \eqref{e:scales_transl_ener}]\label{rem:informal_scales}
Observe that the assumption \eqref{e:scales_u1} with the choice \eqref{e:scales_transl_ener} can be informally stated
as
\begin{equation}\label{e:scales_nonrig}
\left( \frac{\|D^{\balpha} u\|_{L^\infty}}{ \kappa}  \right)^2 \le  \left(\frac{\|\nabla u\|_{L^\infty}}{\kappa}\right)^{|\balpha|_x+1} \|\nabla u\|^{2|\balpha|_t}_{L^\infty},
\end{equation}
because the equivalent form \eqref{e:scales_connect} of \eqref{e:scales_u1} has the right hand side of \eqref{e:scales_nonrig}, if one 
imagines that the estimate for $\|\nabla u_q\|_{L^\infty}$ is saturated, i.e. that  $\|\nabla u_q\|_{L^\infty} \sim \delta_q^{\sfrac{1}{2}} \lambda_q$.
\end{remark}
\subsection{Notations of this section}\label{sec:ener_not}
To state our estimates and to proceed efficiently with the proof, introduce the following shorthand notations, for:\\
(i)dissipation
\begin{equation}     \label{n:sh_est1}
\mathcal{D}_\rho := \kappa \int_0^T \|\nabla \rho (s)\|_{L^2}^2 ds 
\end{equation}
(ii) 'initial datum':
 \begin{equation} \label{n:sh_est1.1}
ID_n  := \|\nabla^n \rho_{in}\|^2_{L^2}
\end{equation}
(iii) 'initial datum and forcing':
\begin{equation}     \label{n:sh_est2}
\begin{aligned}
IDF_n  := \|\nabla^n \rho_{in}\|^2_{L^2}  + \kappa^{-1} \int_0^T \|\nabla^n f (s)\|_{L^2}^2 ds
    \end{aligned}
\end{equation}
(iv) 'Leray-Hopf quantity of function $h$'
\begin{equation}     \label{n:sh_estLH}
LH_h (\balpha)  := \sup_{t\le T} \| D^{\balpha}  h (t) \|^2_2 
+ \kappa \int_0^T \| D^{\balpha}  \nabla h (t) \|^2_2 dt 
\end{equation}
(v) 'maximum of Leray-Hopf quantities with a certain order of derivatives'
\begin{equation}     \label{n:sh_estLHmax}
LH_h (l_x,l_t)  := \sup_{\bbeta: \substack{|\bbeta|_x = l_x, \\  |\bbeta|_t = l_t}} LH_h(\bbeta)
\end{equation}
(vi) renormalised version of Leray-Hopf quantities
\begin{equation} \label{n:sh_estLHmax_r}
{\mr{LH}_h}(l_x,l_t)  := \sup_{\bbeta: \substack{|\bbeta|_x = l_x, \\  |\bbeta|_t = l_t}} {\mr{LH}_h}(\bbeta) = \sup_{\bbeta: \substack{|\bbeta|_x = l_x, \\  |\bbeta|_t = l_t}} \left(\sup_{t\le T} \| \mr D^{\bbeta}  h (t) \|^2_2 
+ \kappa \mr\lambda_q^2 \int_0^T \| \mr D_q^{\bbeta} \mr \nabla_q h (t) \|^2_2 dt \right).
\end{equation}

\subsection{General energy-type estimates}
The main result of this section, which requires of scales merely Assumption \ref{ass:scales_foor_ener}, reads
\begin{proposition}\label{prop:energy_est}
Take smooth, divergence free $u$, smooth forcing $f$ and datum $\rho_{\ini}$. If the scales Assumption \ref{ass:scales_foor_ener} holds, then the smooth solution of \eqref{e:parabolic_const_gen} satisfies:\\
(I) For $n \ge 1$, $n \le 2U$ 
\begin{equation}\label{e:ener_gen_space}
\sup_{t \le T} \|\nabla^n \rho(t)\|^2_{L^2}+ \kappa\int_0^T \|\nabla^{n+1} \rho (s)\|_{L^2}^2 ds \lesssim_n \mathcal{A}^n \mathcal{S}^n  \mathcal{D}_\rho + IDF_n=:RHS(n)
 \end{equation}
(II) Using $RHS(n)$ defined in \eqref{e:ener_gen_space}, we have for the case $|\bbeta|_t>0$: For any $|\bbeta| \le U$
\begin{equation}\label{e:ener_gen_mixed}
\begin{aligned}
LH_\rho(\bbeta)  \lesssim_{|\bbeta|} & \kappa^{2|\bbeta|_t} RHS(|\bbeta|_x+2|\bbeta|_t) \\
&+  \sum_{j=1}^{|\bbeta|_t} \sum_{i=0}^{j-1}  \kappa^{2i} (\kappa^2 \mathcal{A}\mathcal{S}^{2})^{j-i} \mathcal{T}^{|\bbeta|_t-j} \sum_{\substack{k_x + l_x = |\bbeta|_x+2i,  l_x>0}} \mathcal{S}^{k_x} RHS(l_x) \\
&+
\sum_{j=0}^{|\bbeta|_t-1} \kappa^{2j} LH_{\divr f} (|\bbeta|_x+2j,|\bbeta|_t-j-1)\\
+&\sum_{n=1}^{|\bbeta|_t} \sum_{j=0}^{n-2} \kappa^{2j} (\kappa^2 \mathcal{A}\mathcal{S}^{2})^{n-j-1}
\sum_{\substack{k_x + l_x = |\bbeta|_x+2j,\\k_t + l_t = |\bbeta|_t-n,}} \mathcal{S}^{k_x} \mathcal{T}^{k_t} LH_{\divr f} (l_x,l_t),
\end{aligned}
\end{equation}
where we employ convention that sums over empty sets (of indices) are vanishing.
\end{proposition}
In fact, the above Proposition will be proven as a special case of even more comprehensive, technical Proposition \ref{prop:energy_est_tech}, where the quantities with which we estimate $LH_\rho(\bbeta)$ are more flexible, e.g.\ they are allowed to have less space derivatives, at the cost of more advective derivatives. These proofs are provided in Section \ref{ss:energy_est_gp}.

The generality of Propostions \ref{prop:energy_est},\ref{prop:energy_est_tech} makes them useful for diverse choices of scales. For the purpose of particular scales used in this paper, we specify Propostion \ref{prop:energy_est} in the following two sections.

\subsection{First specification of the general energy-type estimates}\label{sec:energ_spec}
In this paper, scales do not only satisfy Assumption \ref{ass:scales_foor_ener}, but they are also chosen according to \eqref{e:scales_transl_ener}. This means in particular that $\mathcal{A}$ can be consumed by $\lesssim$ (since it is a constant) and that $\kappa^2  = \mathcal{T}\mathcal{S}^{-2}$. It turns out that these two properties suffice to extract the estimates that we need, even without resorting to the specific choice \eqref{e:scales_transl_ener} of scales. More precisely, in addition to Assumption \ref{ass:scales_foor_ener}, we take
\begin{assumption}\label{ass:scales_foor_ener2}
The scales satisfy $\mathcal{A} = C$, $\kappa^2  = \mathcal{T}\mathcal{S}^{-2}$
\end{assumption}
Using Assumptions \ref{ass:scales_foor_ener}, \ref{ass:scales_foor_ener2}
and definition of $RHS(n)$, we see that for the second line of r.h.s.\ of \eqref{e:ener_gen_mixed} one has 
\begin{equation}\label{e:dens_simple}
\begin{aligned}
&\sum_{j=1}^{|\bbeta|_t} \sum_{i=0}^{j-1}  \kappa^{2i} (\kappa^2 \mathcal{A}\mathcal{S}^{2})^{j-i} \mathcal{T}^{|\bbeta|_t-j} \sum_{\substack{k_x + l_x = |\bbeta|_x+2i,  l_x>0}} \mathcal{S}^{k_x} RHS(l_x) \\
&\lesssim \sum_{j=1}^{|\bbeta|_t} \sum_{i=0}^{j-1} \mathcal{S}^{-2i} \mathcal{T}^{|\bbeta|_t} \sum_{\substack{k_x + l_x = |\bbeta|_x+2i,  l_x>0}} \mathcal{S}^{k_x} \left( \mathcal{S}^{l_x}  \mathcal{D}_\rho + IDF_{l_x} \right) \\
&\lesssim \mathcal{T}^{|\bbeta|_t} \mathcal{S}^{|\bbeta|_x}  \mathcal{D}_\rho + \mathcal{T}^{|\bbeta|_t} \mathcal{S}^{|\bbeta|_x} \sum_{i=0}^{|\bbeta|_t-1} \sum_{\substack{k_x + l_x = |\bbeta|_x+2i,  l_x>0}} \mathcal{S}^{k_x-2i-|\bbeta|_x}  IDF_{l_x}.
\end{aligned}
\end{equation}
Putting this estimate together with an analogous estimate for the first line of r.h.s.\ of \eqref{e:ener_gen_mixed}, we see that the both the first and the second line of r.h.s.\ of \eqref{e:ener_gen_mixed} are estimated as follows
\begin{equation}\label{e:p:ee:fr1}
\text{two first r.h.s. lines of \eqref{e:ener_gen_mixed}} \lesssim
\mathcal{T}^{|\bbeta|_t} \mathcal{S}^{|\bbeta|_x}  \mathcal{D}_\rho + \mathcal{T}^{|\bbeta|_t} \mathcal{S}^{|\bbeta|_x} \sum_{l=1}^{|\bbeta|_x+2|\bbeta|_t} \mathcal{S}^{-l}  IDF_{l}.
\end{equation}
Concerning the 'forcing terms' of r.h.s.\ of \eqref{e:ener_gen_mixed}, we first observe that for $|\bbeta|_t=1$ the last line of r.h.s.\ of \eqref{e:ener_gen_mixed} vanishes, so the only forcing term is $|\bbeta|_t=1$-case of the middle line of r.h.s.\ of \eqref{e:ener_gen_mixed}, being here simply  $LH_{\divr f} (|\bbeta|_x, 0)$ (consistently with the initialization of the induction in section \ref{ssec:pf_ener_tech}). For $|\bbeta|_t>1$, via $\kappa^2  = \mathcal{T}\mathcal{S}^{-2}$, we have
\begin{equation}\label{e:force_mid_simple}
\begin{aligned}
\text{third r.h.s\ line of \eqref{e:ener_gen_mixed}} = \sum_{j=0}^{|\bbeta|_t-1} &(\mathcal{T}\mathcal{S}^{-2})^{j} LH_{\divr f} (|\bbeta|_x+2j,|\bbeta|_t-j-1)
    \end{aligned}   
\end{equation}
and we have 
\begin{equation}\label{e:force_last_simple}
\begin{aligned}
 \text{fourth r.h.s\ line of \eqref{e:ener_gen_mixed}} &=   \sum_{n=1}^{|\bbeta|_t} \sum_{j=0}^{n-2} \mathcal{S}^{|\bbeta|_x} \mathcal{T}^{|\bbeta|_t-1}
\sum_{\substack{l_x =0, \dots, |\bbeta|_x+2j,\\l_t =0, \dots,  |\bbeta|_t-n}} \mathcal{S}^{-l_x} \mathcal{T}^{-l_t} LH_{\divr f} (l_x,l_t)\\
&\lesssim   \sum_{n=1}^{|\bbeta|_t} \sum_{j=0}^{n-1} \mathcal{S}^{|\bbeta|_x} \mathcal{T}^{|\bbeta|_t-1}
\sum_{\substack{l_x =0, \dots, |\bbeta|_x+2j,\\l_t =0, \dots,  |\bbeta|_t-(j+1)}} \mathcal{S}^{-l_x} \mathcal{T}^{-l_t} LH_{\divr f} (l_x,l_t) \\
&\lesssim \sum_{j=0}^{|\bbeta|_t-1} \mathcal{S}^{|\bbeta|_x} \mathcal{T}^{|\bbeta|_t-1}
\sum_{\substack{l_x =0, \dots, |\bbeta|_x+2j,\\l_t =0, \dots,  |\bbeta|_t-(j+1)}} \mathcal{S}^{-l_x} \mathcal{T}^{-l_t} LH_{\divr f} (l_x,l_t),
\end{aligned}
\end{equation}
where the inequalities are valid, respectively, by lossy extension of the range of sums, and then the fact the the summands do not depend on $n$ any more. Observe that now the  last line of \eqref{e:force_last_simple} contains, as its $l_x= |\bbeta|_x+2j$, $l_t=|\bbeta|_t-(j+1)$-term, the estimate \eqref{e:force_mid_simple}. Thus we have
\begin{equation}\label{e:force_last_simple2}
\begin{aligned}
 &\text{third and fourth r.h.s\ line of \eqref{e:ener_gen_mixed}} \\
& \lesssim \sum_{j=0}^{|\bbeta|_t-1} \mathcal{S}^{|\bbeta|_x} \mathcal{T}^{|\bbeta|_t-1}
\sum_{\substack{l_x =0, \dots, |\bbeta|_x+2j,\\l_t =0, \dots,  |\bbeta|_t-(j+1)}} \mathcal{S}^{-l_x} \mathcal{T}^{-l_t} LH_{\divr f} (l_x,l_t),
\end{aligned}
\end{equation}
Putting together \eqref{e:p:ee:fr1}, \eqref{e:force_last_simple2} for \eqref{e:ener_gen_mixed}, using the definition \eqref{n:sh_est2} of the quantity $IDF_n$ as well as providing basic energy estimate, we have
\begin{corollary}
\label{cor:energy_est_specifed}
Take smooth, divergence free $u$, smooth forcing $f$ and datum $\rho_{\ini}$. Then the smooth solution of \eqref{e:parabolic_const_gen} satisfies the base energy identity
\begin{equation}\label{e:ener_base_specifed}
\mathcal{D}_\rho := \sup_{t \le T} \| \rho(t)\|^2_{L^2}+ \kappa\int_0^T \|\rho (s)\|_{L^2}^2 ds \le  \|\rho_{in}\|^2_{L^2}  + \kappa^{-1} \int_0^T \|f (s)\|_{L^2}^2 ds.
\end{equation}
Furthermore, if the scales Assumptions \ref{ass:scales_foor_ener} and \ref{ass:scales_foor_ener2} hold, then we have:\\
(I. Case of pure space derivatives) For $n \ge 1$, $n \le 2U$ 
\begin{equation}\label{e:ener_gen_space_specifed}
\sup_{t \le T} \|\nabla^n \rho(t)\|^2_{L^2}+ \kappa\int_0^T \|\nabla^{n+1} \rho (s)\|_{L^2}^2 ds \lesssim \mathcal{S}^n  \mathcal{D}_\rho + \left(\|\nabla^n \rho_{in}\|^2_{L^2}  + \kappa^{-1} \int_0^T \|\nabla^n f (s)\|_{L^2}^2 ds \right)
 \end{equation}
(II. Case of mixed derivatives) For the case $|\bbeta|_t>0$: For any $|\bbeta| \le U$
\begin{equation}\label{e:ener_gen_mixed_redux}
\begin{aligned}
\mathcal{T}^{-|\bbeta|_t} \mathcal{S}^{-|\bbeta|_x} LH_\rho(\bbeta)  &\lesssim    \mathcal{D}_\rho + \sum_{l=1}^{|\bbeta|_x+2|\bbeta|_t} \mathcal{S}^{-l}   \left( \|\nabla^l \rho_{in}\|^2_{L^2}  + \kappa^{-1} \int_0^T \|\nabla^l f (s)\|_{L^2}^2 ds\right)\\
&+  \sum_{j=0}^{|\bbeta|_t-1} \mathcal{T}^{-1}
\sum_{\substack{l_x =0, \dots, |\bbeta|_x+2j,\\l_t =0, \dots,  |\bbeta|_t-(j+1)}} \mathcal{S}^{-l_x} \mathcal{T}^{-l_t} LH_{\divr f} (l_x,l_t),
\end{aligned}
\end{equation}
\end{corollary}
Several remarks related to Corollary \ref{cor:energy_est_specifed} are in order
\begin{remark}[Consistency with $u_q$, $v_q$ and generality of advective derivatives]\label{rem:uqvq_consist}
 Lemma \ref{rem:scales_consistent} shows that
the assumptions of Corollary \ref{cor:energy_est_specifed} hold for solutions of \eqref{e:parabolic_const_gen} in case of $D_u$ there being either $D$-derivative (i.e.\ $u:= u_q$) or $O$-derivative (i.e.\ $u:= v_q$), with  $U= \lfloor \frac{Y-N_*-1}{2} \rfloor$. Analogously, taking into account results of Section \ref{ss:r_estimates}  Corollary \ref{cor:energy_est_specifed} holds for any combination of $O^{\bbeta_1}$ and $D^{\bbeta_2}$ present in the quantity $LH_\rho (\bbeta)$, provided the same combination appears in $LH_{\divr f} (\bbeta)$. 
\end{remark}
\begin{remark}[Non-divergence form forcing]\label{rem:nondiv_force}
In case \eqref{e:parabolic_const_gen} is replaced with \begin{equation}    \label{e:parabolic_const_gen_ndiv}
\begin{split}
      D_u \rho - \kappa \Delta \rho =&  F, \\ 
      \rho(\cdot,0) =&\, \rho_{\ini}.
\end{split}
\end{equation}
we will content ourselves with possibly the estimates obtained as follows. The pure space estimate with $n=0$ uses  Poincare inequality to  estimate $\rho$ with its gradient, and similarly the  pure space estimate with $n>0$ uses one integration by parts 
\[
\int \nabla^n F \nabla^n \rho dx= -\int \nabla^{n-1} F \nabla^{n+1} \rho dx
\]
Consequently, for non-divergence forcing, Corollary \ref{cor:energy_est_specifed} holds with the following changes\\
 (i) Instead of $LH_{\divr f} (\bbeta)$ we have $LH_{F} (\bbeta)$, \\
 (ii) instead of \[\left(\|\nabla^n \rho_{in}\|^2_{L^2}  + \kappa^{-1} \int_0^T \|\nabla^n f (s)\|_{L^2}^2 ds \right)\] in right-hand sides, we have \[\left(\|\nabla^n \rho_{in}\|^2_{L^2}  + \kappa^{-1} \int_0^T \|\nabla^{0 \wedge (n-1)} F (s)\|_{L^2}^2 ds \right)\]
\end{remark}

\subsection{Further specification of the general energy-type estimates}

In order to make our estimates even more convenient and tailored to our needs, we split them into two cases: of nonzero initial datum and zero forcing, and of zero initial datum and nonzero forcing. Furthermore, we renormalise certain quantities. To this end, recall definition \eqref{n:sh_estLHmax} of quantity $LH_h (l_x,l_t)$ and definition \eqref{n:sh_estLHmax_r} of its renormalisation, i.e.\ $\mr LH_h (l_x,l_t)$.

\subsubsection{Setup}
Consider
\begin{equation}     \label{e:parabolic_const}
\begin{split}
      D_{q} \rho - \kappa \Delta \rho =&\, F = \mr \lambda_q^{-1} \divr f, \\ 
      \rho(\cdot,0) =&\, \rho_{\ini}.
\end{split}
\end{equation}
\begin{remark}
In view of Remark \ref{rem:uqvq_consist}, if we replace $D_{q}$ by $O_{q}$ in \eqref{e:parabolic_const}, or in definition of the Leray-Hopf quantities used below, the results of this section hold true.
\end{remark}

\subsubsection{Scales and renormalised scales}
We will use the scales
\begin{equation}\label{e:scales_transl_ener2}
\begin{aligned}
    \mathcal{S} &:=
     \delta_q^{\sfrac{1}{2}} \lambda_q \kappa^{-1}_q \\
         \mathcal{T} &:=
    (\delta_q^{\sfrac{1}{2}} \lambda_q)^2 
\end{aligned}
\end{equation}
as in \eqref{e:scales_transl_ener} and additionally their renormalised counterparts:
\begin{equation}\label{e:scales_transl_ener2ren}
\begin{aligned}
    \mathcal{\mr S} &:=  \frac{\mathcal{S}}{\mr \lambda^2_q} \\
         \mathcal{\mr T} &:=
\frac{\mathcal{T}}{\mr \mu^2_q} 
\end{aligned}
\end{equation}
One can informally identify these scales as follows 
\begin{equation}\label{e:scales_transl_ener3}
\begin{aligned}
    \mathcal{S} &\sim \frac{\| \nabla u_q \|_{\infty} } {\kappa_q}, \qquad \mathcal{\mr S} &\sim \frac{\| \nabla u_q \|_{\infty} } {\kappa_q \mr \lambda^2_q},  \\
         \mathcal{T} &\sim \| \nabla u_q \|^2_{\infty}, \qquad   \mathcal{\mr T} &\sim \frac{\| \nabla u_q \|^2_{\infty}}{\mr \mu^2_q},
\end{aligned}
\end{equation}
see also Remark \ref{rem:informal_scales}.

\subsubsection{Result for initial datum and no forcing}
Using Corollary \ref{cor:energy_est_specifed}, taking into account Remark \ref{rem:uqvq_consist}, as well as notation as in Section \ref{sec:ener_not},
we obtain directly the case with no forcing
\begin{lemma}     \label{l:L2_init_energy}
Consider \eqref{e:parabolic_const} with $f = 0$. For any $|\bbeta| \le \lfloor \frac{Y-N_*-1}{2} \rfloor$, $|\bbeta|>0$
we have
\begin{equation}  \label{e:2:L2_init_energy}
 \mathcal{S}^{-|\bbeta|_x}
           \mathcal{T}^{-|\bbeta|_t} {{LH}_\rho}(|\bbeta|_x,|\bbeta|_t)  \lesssim_{|\bbeta|}    \kappa \iint |\nabla \rho|^2 dxdt  + \sum_{p=1}^{|\bbeta|_x+2|\bbeta|_t}  \mathcal{S}^{-p} \| \nabla^p \rho_{\ini} \|_2^2.
\end{equation}
Reformulated in renormalized quantities, we have
\begin{equation}  \label{e:6:L2_init_energy}
\mathcal{\mr S}^{-|\bbeta|_x} \mathcal{\mr T}^{-|\bbeta|_t}  {\mr{LH}_\rho}(|\bbeta|_x,|\bbeta|_t) 
\lesssim_{|\bbeta|}   \kappa \mr \lambda_q^2 \iint |\mr\nabla \rho|^2 dxdt  + \sum_{p=1}^{|\bbeta|_x+2|\bbeta|_t}  \mathcal{\mr S}^{-p}  \| \mr \nabla_q^p \rho_{\ini} \|_2^2.
\end{equation}
\end{lemma}

\begin{remark}[Case of pure space derivatives]
In both estimates, in the case $|\bbeta|_t=0$, the sum $\sum_{p=1}^{|\bbeta|_x+2|\bbeta|_t}$ may be replaced with its last summand (so that after multiplying with the space scale present as the first r.h.s.\ term, the initial datum is bare $\| \nabla^{|\bbeta|_x} \rho_{\ini} \|_2^2$, or respectively $\| \mr \nabla_q^{|\bbeta|_x} \rho_{\ini} \|_2^2$). For instance \eqref{e:2:L2_init_energy} for $|\bbeta|_t=0$ reads
\begin{equation}  \label{e:2:L2_init_energy_ps}
    \mathcal{S}^{-|\bbeta|_x}
    {{LH}_\rho}(|\bbeta|_x,|\bbeta|_t)  \lesssim_{|\bbeta|} 
       \kappa \iint |\nabla \rho|^2 dxdt  +  \|\nabla^{|\bbeta|_x} \rho_{\ini} \|_2^2.
\end{equation}
\end{remark}

\subsubsection{Results with forcing and zero initial datum (non-renormalised)}

Let us now specify Corollary \ref{cor:energy_est_specifed} to the case with no initial datum, but forced. Our forcing that fits Corollary \ref{cor:energy_est_specifed} equals $\mr\lambda_q^{-1}f$, see \eqref{e:parabolic_const}, or in non-divergence case simply $F$ (then we take into account Remark \ref{rem:nondiv_force}). In the former case, observe that the quantities appearing in the estimates are quadratic, thus ${\mr\lambda_q}^{-2}$ appears accordingly on right-hand sides. Taking also into account Remark \ref{rem:nondiv_force}, Corollary \ref{cor:energy_est_specifed} yields 

\begin{lemma}     \label{l:L2_F_energy_pre}
Consider \eqref{e:parabolic_const} with $\rho_{\ini} = 0$. For any $|\bbeta| \le \lfloor \frac{Y-N_*-1}{2} \rfloor$, $|\bbeta|>0$ we have
\begin{equation}  \label{e:2:L2_F_energy_pre}
\begin{aligned}
       \mathcal{T}^{-|\bbeta|_t} \mathcal{S}^{-|\bbeta|_x} {{LH}_\rho}(|\bbeta|_x,|\bbeta|_t)  &\lesssim \mr\lambda_q^{-2}   \sum_{j=0}^{|\bbeta|_t-1} \mathcal{T}^{-1}
\sum_{\substack{l_x =0, \dots, |\bbeta|_x+2j,\\l_t =0, \dots,  |\bbeta|_t-(j+1)}} \mathcal{S}^{-l_x} \mathcal{T}^{-l_t} LH_{\divr f} (l_x,l_t) \\
       &+ \mr\lambda_q^{-2}  \kappa^{-1}\sum_{l=1}^{|\bbeta|_x+2|\bbeta|_t}  \mathcal{S}^{-l}   \int_0^T \|\nabla^l f (s)\|_{L^2}^2 ds
\end{aligned}
\end{equation}
and
\begin{equation}  \label{e:6:L2_F_energy_pre}
\begin{aligned}
      \mathcal{T}^{-|\bbeta|_t} \mathcal{S}^{-|\bbeta|_x} {{LH}_\rho}(|\bbeta|_x,|\bbeta|_t)   &\lesssim   \sum_{j=0}^{|\bbeta|_t-1} \mathcal{T}^{-1}
\sum_{\substack{l_x =0, \dots, |\bbeta|_x+2j,\\l_t =0, \dots,  |\bbeta|_t-(j+1)}} \mathcal{S}^{-l_x} \mathcal{T}^{-l_t} LH_{F} (l_x,l_t),  \\
      &+ \kappa^{-1}\sum_{l=1}^{|\bbeta|_x+2|\bbeta|_t}  \mathcal{S}^{-l}   \int_0^T \|\nabla^{l-1} F (s)\|_{L^2}^2 ds.
\end{aligned}
\end{equation}
As in the previous lemma, in case of pure space derivatives the first lines of r.h.s.'s vanish, and the remaining sum (over $l$'s) may be replaced by its last summand.
\end{lemma}

\subsubsection{Results with forcing and zero initial datum (renormalised and simplified)}

Finally, let us renormalise the quantities in Lemma \ref{l:L2_F_energy_pre} and simplify the formulas.

Consider first the l.h.s's in Lemma \ref{l:L2_F_energy_pre}.
Since the order of $\bbeta$ agrees with the powers of scales appearing there, the l.h.s's in Lemma \ref{l:L2_F_energy_pre} are equal to their renormalised counterparts, i.e.\
\begin{equation}\label{e:lhs_ener_nrnrm}
    \mathcal{T}^{-|\bbeta|_t}  \mathcal{S}^{-|\bbeta|_x} {LH}_{\rho}(\bbeta) =     \mathcal{\mr T}^{-|\bbeta|_t}  \mathcal{\mr S}^{-|\bbeta|_x} \mr{LH}_{\rho}(\bbeta).
\end{equation}
(i. Case of forcing in divergence form.) Let us now consider the r.h.s.\ of \eqref{e:2:L2_F_energy_pre}. 
Concerning the first term of the r.h.s.\ of \eqref{e:2:L2_F_energy_pre}, we have
\begin{equation}\label{e:2:L2_F_energy_pre_n}
\begin{aligned}
&\mr\lambda_q^{-2}   \sum_{j=0}^{|\bbeta|_t-1} \mathcal{T}^{-1}
\sum_{\substack{l_x =0, \dots, |\bbeta|_x+2j,\\l_t =0, \dots,  |\bbeta|_t-(j+1)}} \mathcal{S}^{-l_x} \mathcal{T}^{-l_t} LH_{\divr f} (l_x,l_t) \\
&=  \mr\lambda_q^{-2} \underbrace{\mathcal{S}^2 \mathcal{T}^{-1}}_{:=\kappa^{-2}} \sum_{j=0}^{|\bbeta|_t-1} \sum_{\substack{l_x =0, \dots, |\bbeta|_x+2j,\\l_t =0, \dots,  |\bbeta|_t-(j+1)}} \mathcal{S}^{-l_x-2} \mathcal{T}^{-l_t} LH_{\divr f} (l_x,l_t) \\
&\le \mr\lambda_q^{-2}\kappa^{-2}\sum_{j=0}^{|\bbeta|_t-1} \sum_{\substack{\tilde l_x =1, \dots, |\bbeta|_x+2j+1,\\l_t =0, \dots,  |\bbeta|_t-(j+1)}} \mathcal{S}^{-\tilde l_x-1} \mathcal{T}^{-l_t} LH_{f} (\tilde l_x,l_t) \\
&= \mr\lambda_q^{-4}\kappa^{-2}\sum_{j=0}^{|\bbeta|_t-1}  \sum_{\substack{\tilde l_x =1, \dots, |\bbeta|_x+2j+1,\\l_t =0, \dots,  |\bbeta|_t-(j+1)}} \mathcal{\mr S}^{-\tilde l_x-1} \mathcal{\mr T}^{-l_t} \mr{LH}_{f} (\tilde l_x,l_t).
\end{aligned}
\end{equation}
Concerning the second term of the r.h.s.\ of \eqref{e:2:L2_F_energy_pre}, we have
\begin{equation}\label{e:2:L2_F_energy_pre_n2}
\begin{aligned}
\mr\lambda_q^{-2}  \kappa^{-1}\sum_{l=1}^{|\bbeta|_x+2|\bbeta|_t}  \mathcal{S}^{-l}   \int_0^T \|\nabla^l f (s)\|_{L^2}^2 ds &= \mr\lambda_q^{-2}  \kappa^{-1}\sum_{l=1}^{|\bbeta|_x+2|\bbeta|_t}  \mathcal{\mr S}^{-l}   \int_0^T \|\mr\nabla_q^l f (s)\|_{L^2}^2 ds \\
&= \kappa^{-2} \mr\lambda_q^{-4} \sum_{ l = 1 }^{|\bbeta|_x+2|\bbeta|_t}  \mathcal{\mr S}^{-l} 
             \underbrace{\kappa \mr\lambda_q^{2} \int_0^T \|\mr\nabla_q^l f (s)\|_{L^2}^2 ds}_{\le  \mr{LH}_{f} ( l-1,0)} \\
&\le \kappa^{-2} \mr\lambda_q^{-4} \sum_{ \tilde l = 0}^{|\bbeta|_x+2|\bbeta|_t-1}  \mathcal{\mr S}^{- \tilde  l-1} 
  \mr{LH}_{f} ( \tilde  l,0)
\end{aligned}
\end{equation}
where the under-braced inequality  holds since $\mr\nabla_q^l = \mr\nabla_q^{l-1} \mr\nabla_q$.

Observe that if in the last line of \eqref{e:2:L2_F_energy_pre_n} we extend the lower bound of the sum from $\tilde l_x=1$ down to $\tilde l_x=0$, we accommodate as the $l_t=0$, $j=|\bbeta|_t-1$ term of the sum the entire \eqref{e:2:L2_F_energy_pre_n2}. Thus, putting together \eqref{e:lhs_ener_nrnrm}, \eqref{e:2:L2_F_energy_pre_n}, \eqref{e:2:L2_F_energy_pre_n2} into \eqref{e:2:L2_F_energy_pre}, we have 
\begin{equation}     \label{e:6:L2_F_energy}
      \mathcal{\mr T}^{-|\bbeta|_t}  \mathcal{\mr S}^{-|\bbeta|_x} \mr{LH}_{\rho}(\bbeta) 
            \lesssim_{|\bbeta|} \mr\lambda_q^{-4}\kappa^{-2}\sum_{j=0}^{|\bbeta|_t-1}  \sum_{\substack{\tilde l_x =0, \dots, |\bbeta|_x+2j+1,\\l_t =0, \dots,  |\bbeta|_t-(j+1)}} \mathcal{\mr S}^{-\tilde l_x-1} \mathcal{\mr T}^{-l_t} \mr{LH}_{f} (\tilde l_x,l_t).
\end{equation}
(ii. Case of forcing in non-divergence form.)
Analogously we deal with the \eqref{e:6:L2_F_energy_pre}. More precisely, we have for the first r.h.s. term of \eqref{e:6:L2_F_energy_pre}
\begin{equation}\label{e:2:L2_F_energy_pre_nF}
\begin{aligned}
&\sum_{j=0}^{|\bbeta|_t-1} \mathcal{T}^{-1}
\sum_{\substack{l_x =0, \dots, |\bbeta|_x+2j,\\l_t =0, \dots,  |\bbeta|_t-(j+1)}} \mathcal{S}^{-l_x} \mathcal{T}^{-l_t} LH_{F} (l_x,l_t) \\
&= \underbrace{\mathcal{S}^2 \mathcal{T}^{-1}}_{:=\kappa^{-2}} \sum_{j=0}^{|\bbeta|_t-1} \sum_{\substack{l_x =0, \dots, |\bbeta|_x+2j,\\l_t =0, \dots,  |\bbeta|_t-(j+1)}} \mathcal{S}^{-l_x-2} \mathcal{T}^{-l_t} LH_{F} (l_x,l_t) \\
&= \mr\lambda_q^{-4}\kappa^{-2}\sum_{j=0}^{|\bbeta|_t-1}  \sum_{\substack{l_x =0, \dots, |\bbeta|_x+2j,\\l_t =0, \dots,  |\bbeta|_t-(j+1)}} \mathcal{\mr S}^{-l_x-2} \mathcal{\mr T}^{-l_t} LH_{F} (l_x,l_t) 
\end{aligned}
\end{equation}
and for the second r.h.s. term of \eqref{e:6:L2_F_energy_pre}
\[
\begin{aligned}
  \kappa^{-1}\sum_{l=1}^{|\bbeta|_x+2|\bbeta|_t}  \mathcal{S}^{-l}   \int_0^T \|\nabla^{l-1} F (s)\|_{L^2}^2 ds &= \kappa^{-2} \mr\lambda_q^{-4} \sum_{ l = 1 }^{|\bbeta|_x+2|\bbeta|_t}  \mathcal{\mr S}^{-l} 
             \kappa \mr\lambda_q^{2} \int_0^T \|\mr\nabla_q^{l-1} F (s)\|_{L^2}^2 ds \\
        &= \kappa^{-2} \mr\lambda_q^{-4} \sum_{ \tilde l = 0 }^{|\bbeta|_x+2|\bbeta|_t-1}  \mathcal{\mr S}^{-\tilde l-1} 
             \kappa \mr\lambda_q^{2} \int_0^T \|\mr\nabla_q^{\tilde l} F (s)\|_{L^2}^2 ds
\end{aligned}
\]
Splitting in the preceding identity the sum over $\tilde l$ into the summand for $\tilde l=0$ and the rest, we have   
\begin{equation}\label{e:2:L2_F_energy_pre_n2F}
\begin{aligned}
  \kappa^{-1}&\sum_{l=1}^{|\bbeta|_x+2|\bbeta|_t}  \mathcal{S}^{-l}   \int_0^T \|\nabla^{l-1} F (s)\|_{L^2}^2 ds \\
  &= \kappa^{-2} \mr\lambda_q^{-4} \sum_{ \tilde l = 1 }^{|\bbeta|_x+2|\bbeta|_t-1}  \mathcal{\mr S}^{-\tilde l-1} 
             \kappa \mr\lambda_q^{2} \int_0^T \|\mr\nabla_q^{\tilde l} F (s)\|_{L^2}^2 ds + \kappa^{-2} \mr\lambda_q^{-4} \mathcal{\mr S}^{-1} 
             \kappa \mr\lambda_q^{2} \int_0^T \|\mr\nabla_q^{0} F (s)\|_{L^2}^2 ds \\
             &\le \kappa^{-2} \mr\lambda_q^{-4} \sum_{ \tilde l = 1 }^{|\bbeta|_x+2|\bbeta|_t-1}  \mathcal{\mr S}^{-\tilde l-1} \mr{LH}_{F} (\tilde l-1,0) + \kappa^{-2} \mr\lambda_q^{-4} \mathcal{\mr S}^{-1} 
             \kappa \mr\lambda_q^{2} \int_0^T \|\mr\nabla_q^{0} F (s)\|_{L^2}^2 ds \\
             &= \kappa^{-2} \mr\lambda_q^{-4} \sum_{  l = 0}^{|\bbeta|_x+2|\bbeta|_t-2}  \mathcal{\mr S}^{-l-2} \mr{LH}_{F} (l,0) + \kappa^{-1} \mr\lambda_q^{-2} \mathcal{\mr S}^{-1} 
             \int_0^T \|F (s)\|_{L^2}^2 ds.
\end{aligned}
\end{equation}
Observe that the first term of the last line of \eqref{e:2:L2_F_energy_pre_n2F} is the  $l_t=0$, $j=|\bbeta|_t-1$ term of the sum in the last line of \eqref{e:2:L2_F_energy_pre_nF}. Taking this into account and 
putting together \eqref{e:lhs_ener_nrnrm}, \eqref{e:2:L2_F_energy_pre_nF}, \eqref{e:2:L2_F_energy_pre_n2F} into \eqref{e:6:L2_F_energy_pre}, one has 
\begin{equation}     \label{e:6:L2_F_energy_nd}
\begin{aligned}
      \mathcal{\mr T}^{-|\bbeta|_t}  \mathcal{\mr S}^{-|\bbeta|_x} \mr{LH}_{\rho}(\bbeta) 
            &\lesssim_{|\bbeta|}  \mr\lambda_q^{-4}\kappa^{-2}\sum_{j=0}^{|\bbeta|_t-1}  \sum_{\substack{l_x =0, \dots, |\bbeta|_x+2j,\\l_t =0, \dots,  |\bbeta|_t-(j+1)}} \mathcal{\mr S}^{-l_x-2} \mathcal{\mr T}^{-l_t} LH_{F} (l_x,l_t) \\ 
            &+\kappa^{-1} \mr\lambda_q^{-2} \mathcal{\mr S}^{-1} 
             \int_0^T \|F (s)\|_{L^2}^2 ds.
\end{aligned}
\end{equation}
The formulas \eqref{e:6:L2_F_energy}, \eqref{e:6:L2_F_energy_nd} are our key estimates. In order to shorten the sums appearing on their r.h.s.'s, let us use  the notation $[\balpha] := |\balpha|_x + 2 |\balpha|_t$, see section \ref{sss:ntt:multi_ind}, and rewrite the sum in \eqref{e:6:L2_F_energy} according to
\[
\sum_{j=0}^{|\bbeta|_t-1}  \sum_{\substack{\tilde l_x =0, \dots, |\bbeta|_x+2j+1,\\l_t =0, \dots,  |\bbeta|_t-(j+1)}} (\cdot)\le  \sum_{\substack{\tilde l_x+2 l_t \le |\bbeta|_x+2|\bbeta|_t-1,\\l_t =0, \dots,  |\bbeta|_t-1}}  (\cdot)\le \sum_{ \substack{|\balpha|_x + 2|\balpha|_t \le  |\bbeta|_x+2|\bbeta|_t-1, \\ |\balpha|_t \le |\bbeta|_t-1} }(\cdot) = \sum_{ \substack{[\balpha]  \le  [\bbeta] -1, \\ |\balpha|_t \le |\bbeta|_t-1} }(\cdot),
\]
where the penultimate sum means that we sum over all $\balpha \in \mathcal{I}$ such that  $0 \leq |\balpha|_x + 2|\balpha|_t \leq |\bbeta|_x+2|\bbeta|_t-1$ and $0 \le |\balpha|_t \le |\bbeta|_t-1$,  and analogously in the last sum. Similar rewriting of the sum of \eqref{e:6:L2_F_energy_nd} differs only in the range being now $[\balpha]  \le  [\bbeta] -2$.
Observe in both cases the ranges are admissible as long as $|\bbeta|>0$, where in the case of pure space derivatives, as long as $|\bbeta|>0$, we have only the range $[\balpha]  \le  [\bbeta] -1$ or, respectively, $[\balpha]  \le  [\bbeta] -2$, present. With these remarks,  formulas \eqref{e:6:L2_F_energy}, \eqref{e:6:L2_F_energy_nd} yield

\begin{lemma}\label{l:L2_F_energy}
Consider \eqref{e:parabolic_const} with $\rho_{\ini} = 0$. For any $\bbeta \in \mathcal{I} $ such that $|\bbeta| \le \lfloor \frac{Y-N_*-1}{2} \rfloor$, $|\bbeta|>0$ we have:\\
In case of the forcing term of \eqref{e:parabolic_const} being $\lambda_q^{-1} \divr f$:
\begin{equation}\label{e:2:L2_F_energy}
    \mathcal{\mr T}^{-|\bbeta|_t}  \mathcal{\mr S}^{-|\bbeta|_x} \mr{LH}_{\rho}(\bbeta) 
            \lesssim_{|\bbeta|} \mr\lambda_q^{-4}\kappa^{-2} \sum_{ \substack{[\balpha]  \le  [\bbeta] -1, \\ |\balpha|_t \le |\bbeta|_t-1} } \mathcal{\mr T}^{-|\balpha|_t}  \mathcal{\mr S}^{-(|\balpha|_x+1)} 
            \cdot \mr{LH}_{f}(\balpha).
    \end{equation}
In case of the forcing term of \eqref{e:parabolic_const} being $F$:  
\begin{equation}\label{e:4:L2_F_energy}
\begin{split}
      \mathcal{\mr T}^{-|\bbeta|_t}  \mathcal{\mr S}^{-|\bbeta|_x} \mr{LH}_{\rho}(\bbeta) 
            \lesssim_{|\bbeta|} &\, \kappa^{-1} \mr\lambda_q^{-2}  \mathcal{\mr S}^{-1}   \int_0^T \| F (s)\|_{L^2}^2 ds \\
            &+\mr\lambda_q^{-4}\kappa^{-2} \sum_{ \substack{[\balpha]  \le  [\bbeta] -2, \\ |\balpha|_t \le |\bbeta|_t-1} }\mathcal{\mr T}^{-|\balpha|_t}  \mathcal{\mr S}^{-(|\balpha|_x+2)} 
            \cdot \mr{LH}_{F}(\balpha).
\end{split}
\end{equation}
In case of pure space derivatives, the ranges of r.h.s. sums are, respectively, $[\balpha]  \le  [\bbeta] -1$, $[\balpha]  \le  [\bbeta] -2$.
\end{lemma}

\begin{corollary}       \label{c:L2_FHom_energy}
Suppose for some $n_0 \leq 12Q^3$ and some $\balpha \in \mathcal{I}_*$ with $[\balpha] \leq 3Q$,  satisfies the following estimates
\begin{align}
      \sup_t \int | \mr D_q^{\bomega} g (t) |^2 dx
      + \kappa \mr \lambda_q^2 \iint | \mr D_q^{\bomega} \mr \nabla_q g |^2 dxdt
      \lesssim&\, \lambda_q^{- 4 [\bomega] b \gamma},     
            \quad \forall \, [\bomega] \leq n_0,      \label{e:2:L2_FHom_energy} \\ 
      h \in \mathcal{P}_q(n_1),    \quad 
      \vertiii{h}_q \lesssim&\, \kappa \mr\lambda_q^2.    \label{e:4:L2_FHom_energy}
\end{align}
The force $f$ is given by
\begin{align}     \label{e:para_cstF}
      f := h \mr D_q^{\balpha} g,
\end{align}
Then we have
\begin{equation}  \label{e:6:L2_FHom_energy}
\begin{split}
      \| \mr D_q^{\bbeta} \rho \|_{L^\infty_t L^2_x} + 
      \kappa^{\sfrac12} \mr\lambda_q \| \mr D_q^{\bbeta} \mr \nabla_q \rho \|_{L^2}
      \lesssim \lambda_q^{ -2 ([\bbeta]+[\balpha]) b\gamma },
            \quad \forall \, [\bbeta] \leq n_0 - 3Q.
\end{split}
\end{equation}
\end{corollary}

\begin{proof}
Note that we have
\begin{align*}
      | \mr D_q^{\bbeta} f |^2
            \sim | \mr D_q^{\bbeta} h \mr D_q^{\balpha} g |^2 
                  + | h \mr D_q^{\bbeta} \mr D_q^{\balpha} g |^2 
\end{align*}
and from Corollary \ref{c:CalOEst}
\begin{align*}
      \| \mr D_q^{\bbeta} h \|_\infty \lesssim \kappa\mr \lambda_q^2 \lambda_q^{-[\bbeta]\gamma_R}, 
            \quad |\bbeta| \leq Y - n_1.
\end{align*}
Then we have
\begin{align*}
      \mr{LH}_f(\bbeta) \lesssim \kappa^2 \mr \lambda_q^4 \cdot \lambda_q^{-4([\bbeta]+[\balpha])b\gamma},      \quad 
            [\balpha] + [\bbeta] \leq n_0
\end{align*}
Here, we require $\gamma_R \geq 2 b\gamma$. Then from Lemma \ref{l:L2_F_energy}, we have
\begin{align*}
      \mr{LH}_{\rho}(\bbeta) 
            \lesssim_{|\bbeta|} \kappa^{-2} \mr\lambda_q^{-4} \sup_{[\balpha] \leq [\bbeta]-1} \mr{LH}_{f}(\balpha) 
            \lesssim_{|\bbeta|} \lambda_q^{-4([\bbeta]+[\balpha])b\gamma}.
\end{align*}
Here, we use the relation 
\begin{align*}
      \mathcal{\mr S}^2 \lesssim \mathcal{\mr T} \lesssim \mathcal{\mr S}^2,
\end{align*}
which follows from the definition of $\mathcal{\mr S}$ and $\mathcal{\mr T}$, and the parameter relations in Section \ref{ss:parameters}.
\end{proof}

\subsection{Proof of Proposition \ref{prop:energy_est}}\label{ss:energy_est_gp}

\subsubsection{Notational remarks}\label{ssec:note_energ}

We will use at times slightly informal language, to avoid too cumbersome notation. In particular\\
(i. Concerning derivatives) Instead of an arbitrary partial space derivative of order $k$, we simply write $\nabla^k$; consequently '$\cdot$' between high-order derivatives needs to be interpreted appropriately. \\
(ii. Concerning combinatorial coefficients)
We will write $a \sim_{\bbeta} \sum_{i} b_i$ if $a = \sum_{i} c_{\bbeta}(i) b_i$, where $c_{\bbeta}(i)$ are coefficients uniformly bounded in terms of $|\bbeta|$ (or $a \lesssim_{|\bbeta|} \sum_{i} b_i$ in case of estimates).
Since the coefficient $c_{\bbeta}(i)$ will originate from combinatorial coefficients of a product rule, we will refer to  $c_{\bbeta}(i)$'s as to 'combinatorial-type coefficients'.

\subsubsection{Case of pure space derivatives}

    The estimate \eqref{e:ener_gen_space} for pure space derivatives can be readily obtained along the lines of  \cite{burczak2023anomalous}. Indeed, inspecting the proof of Lemma 4.1  \cite{burczak2023anomalous} (see the last inequality in the proof) and taking into account Corollary 4.2 \cite{burczak2023anomalous} (which includes forcing), together with Young inequality to estimate the forcing term, we have
\begin{equation}\label{e:ener_un_f_2}
\begin{aligned}
 &\sup_{t \le T} \|\nabla^n \rho(t)\|^2_{L^2}+ \kappa\int_0^T \|\nabla^{n+1} \rho (s)\|_{L^2}^2 ds \le \\ &\|\nabla^n \rho_{in}\|^2_{L^2} + C_n \left(  \left(\frac{\|\nabla u\|_{L^\infty}}{\kappa}\right)^n + \left( \frac{\|\nabla^{n} u\|_{L^\infty}}{ \kappa}  \right)^\frac{2n}{n+1}   \right)  \kappa \int_0^T \|\nabla \rho (s)\|_{L^2}^2 ds  \\
 &+ \kappa^{-1} \int_0^T \|\nabla^n f (s)\|_{L^2}^2 ds.
\end{aligned}
 \end{equation}
 (The forcing term was left un-estimated in \cite{burczak2023anomalous}, since it was subsequently treated in 'time-averaging' section by further integration by parts. Now we don't need it, since our homogenisation Ansatz is higher order.)

Using in \eqref{e:ener_un_f_2} the scales relationship \eqref{e:scales_u1} we arrive at 
\[
\begin{aligned}
 &\sup_{t \le T} \|\nabla^n \rho(t)\|^2_{L^2}+ \kappa\int_0^T \|\nabla^{n+1} \rho (s)\|_{L^2}^2 ds \le \\ &\|\nabla^n \rho_{in}\|^2_{L^2} + 2 C_n  \mathcal{A}^n  \mathcal{S}^n \kappa \int_0^T \|\nabla \rho (s)\|_{L^2}^2 ds + \kappa^{-1} \int_0^T \|\nabla^n f (s)\|_{L^2}^2 ds
\end{aligned}
\]
which is \eqref{e:ener_gen_space}, taking into account the shorthands \eqref{n:sh_est1}, \eqref{n:sh_est2}.

\subsubsection{Formula for replication in case of advective derivatives}

Consider the case of a single advective derivative, i.e.\ $|\bbeta|_x = |\bbeta| -1$, $|\bbeta|_t = 1$. If the advective derivative is the innermost one, we have from the equation \eqref{e:parabolic_const_gen}
\begin{equation}\label{e:energy_mixed_1}
    \nabla^k  D \rho = \kappa \nabla^k \Delta \rho +  \nabla^k \divr f.   
\end{equation}
If the advective derivative is the outermost one, then
\begin{equation}\label{e:energy_mixed_2}
\begin{aligned}
    D \nabla^k \rho &=  \nabla^k  D  \rho + [D,  \nabla^k] \rho =  \nabla^k  D  \rho + [u \cdot \nabla ,  \nabla^k] \rho \\
    &=  \nabla^k  D  \rho + \sum_{i=1}^k c_{k,i} \nabla^{i} u \cdot \nabla^{k-i} \nabla \rho,
\end{aligned}
\end{equation}
where $c_{k,i}$ are combinatorial coefficients. 
Acting with $\nabla^l$ on \eqref{e:energy_mixed_2} we have
\begin{equation}\label{e:energy_mixed_3}
\begin{aligned}
  \nabla^l  D \nabla^k \rho &=  \nabla^{k+l}  D  \rho + \sum_{j=0}^l  \tilde c_{l,j} \sum_{i=1}^k c_{k,i} \nabla^{i+j} u \cdot \nabla^{k-i+l-j} \nabla \rho \\
  &=  \kappa \nabla^{k+l} \Delta \rho +  \nabla^{k+l} \divr f + \sum_{i=1}^{k+l} c_{k,i,l} \nabla^{i} u \cdot \nabla^{k+l-i} \nabla \rho,
  \end{aligned}
\end{equation}
where for the last identity we used 
\eqref{e:energy_mixed_1}.

Together, we see from \eqref{e:energy_mixed_1} and \eqref{e:energy_mixed_3} with $k+l=\bbeta$ that in the case $|\bbeta|_x = |\bbeta| -1$, $|\bbeta|_t = 1$ it holds
\begin{equation}\label{e:energy_mixed_4}
\begin{aligned}
 D^{\bbeta}  \rho &\sim_{\bbeta} \kappa \nabla^{|\bbeta|_x} \Delta \rho +  \nabla^{|\bbeta|_x} \divr f + \sum_{i=1}^{|\bbeta|_x} \nabla^{i} u \cdot \nabla^{|\bbeta|_x-i} \nabla \rho\\
  &= \kappa \nabla^{|\bbeta|_x} \Delta \rho +  \nabla^{|\bbeta|_x} \divr f + \sum_{j=0}^{|\bbeta|_x-1} \nabla^{j} \nabla u \cdot \nabla^{|\bbeta|_x-j} \rho
 \end{aligned}
\end{equation}
Importantly, we see that we do not have advective derivatives on r.h.s. of \eqref{e:energy_mixed_4}, so we can use the case of pure space derivatives to estimate  $D^{\bbeta}  \rho$ for any $|\bbeta|_x = |\bbeta| -1$, $|\bbeta|_t = 1$. Observe also that we 'paid' for this single advective derivative with two space derivatives: the highest space derivative on r.h.s. is of order $\bbeta|_x+2$ (the dissipative term).

In \eqref{e:energy_mixed_4} we see also our 'informal notation' of Section \ref{ssec:note_energ} (ii) is useful: comparing \eqref{e:energy_mixed_1} and \eqref{e:energy_mixed_2} we see that certain orders of taking advective and space derivatives produce less terms than other ones, but they are all captured by \eqref {e:energy_mixed_4}.

Next, consider the case of two advective derivatives and the rest being space derivatives, i.e.\ $|\bbeta|_x = |\bbeta| -2$, $|\bbeta|_t = 2$. In other words, we consider now $ \nabla^m D \nabla^l  D \nabla^k \rho$. Using \eqref{e:energy_mixed_3} we see that 
\begin{equation}\label{e:2adv_ener}
\begin{aligned}
 D \nabla^l  D \nabla^k \rho &= \kappa   D  \nabla^{k+l} \Delta \rho +   D \nabla^{k+l} \divr f +  D \sum_{i=1}^{k+l} c_{k,i,l} \nabla^{i} u \cdot \nabla^{k+l-i} \nabla \rho \\
 &\sim_{\bbeta} \kappa   D  \nabla^{k+l+2}  \rho +   D \nabla^{k+l} \divr f +  \sum_{i=1}^{k+l}  D \nabla^{i} u \cdot \nabla^{k+l-i} \nabla \rho + \sum_{i=1}^{k+l}  \nabla^{i} u \cdot D \nabla^{k+l+1-i} \rho.
  \end{aligned}
\end{equation}
Acting on \eqref{e:2adv_ener} with $\nabla^m$ we have
\begin{equation}\label{e:energy_mixed_3_sec}
\begin{aligned}
 \nabla^m D \nabla^l  D \nabla^k \rho &\sim_{\bbeta} \kappa  \nabla^m   D  \nabla^{k+l+2}  \rho + \nabla^m  D \nabla^{k+l} \divr f +  \\
 &+\nabla^m \left(\sum_{i=1}^{k+l}  D \nabla^{i} u \cdot \nabla^{k+l-i} \nabla \rho + \sum_{i=1}^{k+l}  \nabla^{i} u \cdot D \nabla^{k+l-i} \nabla \rho \right)
  \end{aligned}
\end{equation}
so, analogously to \eqref{e:energy_mixed_4}, we have now (via $k+l+m = |\bbeta|_x$)
\begin{equation}\label{e:energy_mixed_5}
 D^{\bbeta}  \rho \sim_{\bbeta} \kappa   \sum_{\balpha:   \substack{|\balpha|_x = |\bbeta|_x +2, \\  |\balpha|_t = 1 }} D^{\balpha}  \rho +   \sum_{\balpha:   \substack{|\balpha|_x = |\bbeta|_x, \\  |\balpha|_t = 1 }} D^{\balpha}  \divr f + \sum_{\balpha, \bomega:   \substack{|\balpha|_x + |\bomega|_x = |\bbeta|_x, |\bomega|_x >0\\  |\balpha|_t + |\bomega|_t = 1 }}  D^{\balpha} \nabla u \cdot D^{\bomega} \rho.
\end{equation}

We can obtain now recursively the general case. Indeed, analyzing structure of \eqref{e:energy_mixed_4}, \eqref{e:energy_mixed_5} (where $1 = |\bbeta|_t-1$) we see that in the general case it holds
\begin{equation}\label{e:energy_mixed_gen}
 D^{\bbeta}  \rho \sim_{\bbeta} \kappa   \sum_{\balpha:   \substack{|\balpha|_x = |\bbeta|_x +2, \\  |\balpha|_t = |\bbeta|_t-1 }} D^{\balpha}  \rho +   \sum_{\balpha:   \substack{|\balpha|_x = |\bbeta|_x, \\  |\balpha|_t = |\bbeta|_t-1 }} D^{\balpha} \divr f + \sum_{\balpha, \bomega:   \substack{|\balpha|_x + |\bomega|_x = |\bbeta|_x, |\bomega|_x >0\\  |\balpha|_t + |\bomega|_t = |\bbeta|_t-1 }}  D^{\balpha} \nabla u \cdot D^{\bomega}\rho.
\end{equation}
We want to estimate the Leray-Hopf quantity
\[
LH_\rho(\bbeta)  := \sup_{t\le T} \| D^{\bbeta}  \rho (t) \|^2_2 
+ \kappa \int_0^T \| D^{\bbeta}  \nabla \rho (t) \|^2_2 dt ,
\]
compare its definition \eqref{n:sh_estLH}. So via \eqref{e:energy_mixed_gen} for $\rho$ and for $\nabla \rho$ we have

\begin{equation}\label{e:ener_est_gen_forrec}
\begin{aligned}
LH_\rho(\bbeta) \lesssim_{|\bbeta|}& \kappa^2   \sum_{\sup\, \balpha:   \substack{|\balpha|_x = |\bbeta|_x +2, \\  |\balpha|_t = |\bbeta|_t-1 }} LH_\rho(\balpha) +   \sum_{\sup\, \balpha:   \substack{|\balpha|_x = |\bbeta|_x, \\  |\balpha|_t = |\bbeta|_t-1 }} LH_{\divr f} (\balpha) \\
&+ \sum_{\sup\, \balpha, \bgamma:   \substack{|\balpha|_x + |\bgamma|_x = |\bbeta|_x, |\bgamma|_x >0\\  |\balpha|_t + |\bgamma|_t = |\bbeta|_t-1 }} \| D^{\balpha} \nabla u \|^2_\infty  LH_\rho(\bgamma),   
\end{aligned}
\end{equation}
importantly, as already observed in the special case of two advective derivatives, also in the general formula \eqref{e:energy_mixed_gen} there is one advective derivative less on r.h.s.\ compared to l.h.s. (and two space derivatives more). The formula \eqref{e:ener_est_gen_forrec} applied recursively, allows to decrease number of advective derivatives from $|\balpha|_t$ down to $0$, at the cost of  $2|\balpha|_t$ more space derivatives. Instead of going all the way back, we will use \eqref{e:ener_est_gen_forrec} within induction step.
\subsubsection{$h$-depth  formula}\label{ssec:gf_repl}
Since $|\bbeta| \le U$, we can replace $\lesssim_{|\bbeta|}$ of \eqref{e:ener_est_gen_forrec} with $\lesssim_U$, for which we will write simply $\lesssim$ ($U$ is pre-fixed). This and our considerations thus far
allow to rewrite \eqref{e:ener_est_gen_forrec}, using the quantity \eqref{n:sh_estLHmax},
as
\begin{equation}\label{e:ener_est_gen_mn}
\begin{aligned}
LH_\rho(m,n) \lesssim& \kappa^2 LH_\rho (m+2,n-1)+  LH_{\divr f} (m,n-1) \\
&+ \mathcal{\tilde A} \sum_{   \substack{k_x + l_x = m, l_x>0\\  k_t + l_t = n-1 }} \mathcal{S}^{k_x} \mathcal{T}^{k_t} LH_\rho(l_x, l_t),    
\end{aligned}
\end{equation}
where to estimate $\| D^{\balpha} \nabla u \|^2_\infty$ we used the assumption \eqref{e:scales_u1} that gives 
\[
\| D^{\balpha} \nabla u \|^2_\infty \le  \kappa^2 \mathcal{A}\mathcal{S}^{k_x+2} \mathcal{T}^{k_t} =:\mathcal{\tilde A}\mathcal{S}^{k_x} \mathcal{T}^{k_t}
\]
where by definition $\mathcal{\tilde A} = \kappa^2 \mathcal{A}\mathcal{S}^{2}$.
The inequality \eqref{e:ener_est_gen_mn} allows to conjecture
\begin{proposition}\label{prop:energy_est_tech}
    Within assumptions of Proposition \ref{prop:energy_est} one has for an arbitrary $h=1,\dots, n$
\begin{equation}\label{e:ener_est_gen_conj}
LH_\rho(m,n)  \le RLH_\rho(m,n,h) + RLH_{\divr f}(m,n,h),
\end{equation}
where 
\begin{equation}\label{e:ener_est_gen_conj_rhsr}
\begin{aligned}
RLH_\rho(m,n,h) &:= \kappa^{2h} LH_\rho (m+2h,n-h)+ \sum_{i=0}^{h-1} \kappa^{2i} \mathcal{\tilde A}^{h-i} \sum_{   \substack{\\k_x + l_x = m+2i, l_x>0 \\  k_t + l_t = n-h,  l_t >0 }} \mathcal{S}^{k_x} \mathcal{T}^{k_t} LH_\rho(l_x, l_t) \\
&+  \sum_{j=1}^{h} \sum_{i=0}^{j-1} \kappa^{2i} \mathcal{\tilde A}^{j-i} \mathcal{T}^{n-j} \sum_{\substack{k_x + l_x = m+2i, l_x>0}} \mathcal{S}^{k_x} \left( RHS(l_x)\right),
\end{aligned}
\end{equation}
with $RHS(l_x)=\mathcal{A}^{l_x} \mathcal{S}^{l_x} \mathcal{D}_\rho + IDF_{l_x}$ as defined in \eqref{e:ener_gen_space},
and 
\begin{equation}\label{e:ener_est_gen_conj_rhsforce}
\begin{aligned}
RLH_{\divr f}(m,n,h) :=&  \sum_{j=1}^{h} \kappa^{2(j-1)} LH_{\divr f} (m+2(j-1),n-j)   \\
&+ \sum_{j=1}^{h} \sum_{i=0}^{j-2} \kappa^{2i} \mathcal{\tilde A}^{j-i-1} \sum_{\substack{k_x + l_x = m+2i,\\k_t + l_t = n-j,}} \mathcal{S}^{k_x} \mathcal{T}^{k_t} LH_{\divr f} (l_x,l_t),
\end{aligned}
\end{equation}
where we use the convention that sum over empty set is vanishing, so for $j=1$ the second line of \eqref{e:ener_est_gen_conj_rhsforce} equals $0$.
\end{proposition}
Even though the formulas in Proposition \ref{prop:energy_est_tech} seem cumbersome, their structure should be clear: using the replication formula \eqref{e:ener_est_gen_mn} recursively, one lowers the required-on-$\rho$ advective derivatives at the cost of twice as many space derivatives, paying attention to have $\kappa$, $\mathcal{\tilde A}$, $\mathcal{S}^{k_x}$, $\mathcal{T}$ in appropriate powers. In particular, the first line of \eqref{e:ener_est_gen_conj_rhsr} can be expanded further using the replication formula \eqref{e:ener_est_gen_mn}, while the second line of \eqref{e:ener_est_gen_conj_rhsr} is the 'initial term'; it contains: (i)  the case of $l_t=0$ instead of  $l_t>0$ in the first line of \eqref{e:ener_est_gen_conj_rhsr}, which is the second line of \eqref{e:ener_est_gen_conj_rhsr} for $j=h$, (ii) the previous scales of initial data, i.e.\ the second line of \eqref{e:ener_est_gen_conj_rhsr} for $j<h$. (Of course we can write both terms: involving $l_t=0$ and  involving $l_t>0$ together, but it is useful for application of Proposition \ref{prop:energy_est_tech} to separate them.)
Similarly, the forcing-related term \eqref{e:ener_est_gen_conj_rhsforce} appears, with the difference compared \eqref{e:ener_est_gen_conj_rhsr} being that forcing acts via 'each scale': we can't lower the needed advective derivatives from $n$ down to $n-h$, but we have the full sum $\sum_{j=1}^h$.

\subsubsection{Proof of Proposition \ref{prop:energy_est_tech}}\label{ssec:pf_ener_tech}
We prove now \eqref{e:ener_est_gen_conj} inductively over $h$. For $h=1$ we have it directly from \eqref{e:ener_est_gen_mn}, splitting there 
\[
\sum_{   \substack{k_x + l_x = m, l_x>0  \\  k_t + l_t = n-1 }} \quad \text{ into } \quad \sum_{   \substack{k_x + l_x = m, l_x>0 \\  k_t + l_t = n-1, \\ l_t>0 }} \quad \text{ and } \quad \sum_{   \substack{k_x + l_x = m, l_x>0 \\  k_t + l_t = n-1, \\ l_t=0 }} 
\]
which holds, because it is exactly \eqref{e:ener_est_gen_mn}.

Assume now that \eqref{e:ener_est_gen_conj} holds for any $h\le h_0-1$; this is our inductive assumption and we want to prove that \eqref{e:ener_est_gen_conj} holds for $h_0$. This is a straightforward computation which involves merely shifts in summing indices, but we provide it for completeness.
The inductive assumption is
\begin{equation}\label{e:ener_est_gen_iass}
\begin{aligned}
LH_\rho(m,n) &\lesssim \kappa^{2(h_0-1)} LH_\rho (m+2(h_0-1),n-(h_0-1)) &=:I\\
&+ \sum_{i=0}^{h_0-2} \sum_{   \substack{k_x + l_x = m+2i, l_x>0 \\  k_t + l_t = n-(h_0-1), \\ l_t>0 }} \kappa^{2i} \mathcal{\tilde A}^{(h_0-1)-i} \mathcal{S}^{k_x} \mathcal{T}^{k_t} LH_\rho(l_x, l_t) &=:II\\
 &+\sum_{j=1}^{h_0-1} \sum_{i=0}^{j-1} \kappa^{2i} \mathcal{\tilde A}^{j-i} \mathcal{T}^{n-j} \sum_{\substack{k_x + l_x = m+2i, l_x>0 }} \mathcal{S}^{k_x} RHS(l_x)&=:III\\
&+\sum_{j=1}^{h_0-1}  \kappa^{2(j-1)} LH_{\divr f} (m+2(j-1),n-j)  &=:IV \\
&+ \sum_{j=1}^{h_0-1}  \sum_{i=0}^{j-2} \kappa^{2i} \mathcal{\tilde A}^{j-i-1} \sum_{\substack{k_x + l_x = m+2i,\\k_t + l_t = n-j,}} \mathcal{S}^{k_x} \mathcal{T}^{k_t} LH_{\divr f} (l_x,l_t)  &=:V,
\end{aligned}
\end{equation}
Using \eqref{e:ener_est_gen_mn} with $m:=m+2(h_0-1), n:=n-h_0+1$ we have
\begin{equation}\label{e:ener_est_gen_defI}
\begin{aligned}
I \lesssim& \kappa^{2h_0}LH_\rho (m+2h_0,n-h_0) &=:I^1_\rho\\
&+ \kappa^{2(h_0-1)} LH_{\divr f} (m+2h_0-2,n-h_0) &=:I_{\divr f} \\
&+ \kappa^{2(h_0-1)} \mathcal{\tilde A} \sum_{   \substack{k_x + l_x = m+2h_0-2, l_x>0 \\  k_t + l_t = n-h_0}} \mathcal{S}^{k_x} \mathcal{T}^{k_t} LH_\rho(l_x, l_t)&=:I^2_\rho, 
\end{aligned}
\end{equation}
 Using in $II$ \eqref{e:ener_est_gen_mn} with $m:=l_x, n:=l_t$ we have
 \[
\begin{aligned}
II \lesssim  
\sum_{   \substack{i=0, \dots, (h_0-1)-1, \\k_x + l_x = m+2i, l_x>0 \\  k_t + l_t = n-(h_0-1), l_t>0 }}& \kappa^{2i} \mathcal{\tilde A}^{(h_0-1)-i} \mathcal{S}^{k_x} \mathcal{T}^{k_t}  \kappa^2 LH_\rho (l_x+2,l_t-1) &=:II^1_{\rho} \\
+  \sum_{   \substack{i=0, \dots, (h_0-1)-1, \\k_x + l_x = m+2i, \\  k_t + l_t = n-(h_0-1), l_t>0 }}& \kappa^{2i} \mathcal{\tilde A}^{(h_0-1)-i} \mathcal{S}^{k_x} \mathcal{T}^{k_t} LH_{\divr f} (l_x,l_t-1) &=:II_{\divr f}\\
+ \sum_{   \substack{i=0, \dots, (h_0-1)-1, \\k_x + l_x = m+2i, l_x>0 \\  k_t + l_t = n-(h_0-1), l_t>0 }}& \kappa^{2i} \mathcal{\tilde A}^{(h_0-1)-i} \mathcal{S}^{k_x} \mathcal{T}^{k_t}  \mathcal{\tilde A} \sum_{   \substack{\tilde k_x + \tilde l_x = l_x, \tilde l_x>0 \\  \tilde k_t + \tilde l_t = l_t-1 }} \mathcal{S}^{\tilde k_x} \mathcal{T}^{\tilde k_t} LH_\rho(\tilde l_x, \tilde l_t) &=:II^2_{\rho}
\end{aligned}
\]
By shifting indices of $II^1_{\rho}$ as follows: $i+1$ to $i$, $l_x+2$ to $l_x$ and $l_t-1$ to $l_t$ (the latter allowed since $II^1_{\rho}$ involves sums only over $l_t >0$) we have
\begin{equation}\label{e:ener_est_gen_iass3_r1}
\begin{aligned}
 II^1_\rho \le& \sum_{   \substack{i=1, \dots, h_0-1, \\k_x + l_x = m+2i,  l_x>0\\  k_t + l_t = n-h_0 }}   \kappa^{2i} \mathcal{\tilde A}^{h_0-i} \mathcal{S}^{k_x} \mathcal{T}^{k_t} LH_\rho(l_x, l_t) \\
 &= \sum_{   \substack{i=1, \dots, h_0-1, \\k_x + l_x = m+2i,  l_x>0\\  k_t + l_t = n-h_0, l_t>0}}   \kappa^{2i} \mathcal{\tilde A}^{h_0-i} \mathcal{S}^{k_x} \mathcal{T}^{k_t} LH_\rho(l_x, l_t) \\
 &+ \sum_{   \substack{i=1, \dots, h_0-1, \\k_x + l_x = m+2i,  l_x>0}}   \kappa^{2i} \mathcal{\tilde A}^{h_0-i} \mathcal{S}^{k_x} \mathcal{T}^{n-h_0} \underbrace{LH_\rho(l_x,0)}_{\lesssim RHS(l_x)}
    \end{aligned}
\end{equation}
In the double sum of $II^2_{\rho}$ we compute $\tilde k_x + \tilde l_x = l_x = m+2i-k_x$, so $(\tilde k_x +k_x)+ \tilde l_x = m+2i$ and analogously $(\tilde k_t +k_t)+ \tilde l_t = n-h_0$. Therefore
\begin{equation}\label{e:ener_est_gen_iass3_r2}
 II^2_\rho \lesssim \sum_{   \substack{i=0, \dots, h_0-2, \\ (\tilde k_x + k_x) + l_x = m+2i,  l_x>0 \\  (\tilde k_t + k_t) + l_t = n-h_0 }}   \kappa^{2i} \mathcal{\tilde A}^{h_0-i} \mathcal{S}^{k_x+ \tilde k_x} \mathcal{T}^{k_t +\tilde k_t} LH_\rho(\tilde l_x, \tilde l_t),
\end{equation}
which is already contained in the r.h.s.\ of \eqref{e:ener_est_gen_iass3_r1}.

Together, definitions in \eqref{e:ener_est_gen_iass}, \eqref{e:ener_est_gen_defI} (splitting $I_\rho^2$ into $l_t=0$ and $l_t>0$ part) and inequalities \eqref{e:ener_est_gen_iass3_r1}, \eqref{e:ener_est_gen_iass3_r2} give 
\begin{equation}\label{e:ener_est_gen_iass3_rr1}
\begin{aligned}
I^1_\rho + I^2_\rho + II^1_\rho + II^2_\rho + III &\lesssim \kappa^{2h_0}LH_\rho (m+2h_0,n-h_0) \\
&+ \kappa^{2(h_0-1)} \mathcal{\tilde A} \sum_{   \substack{k_x + l_x = m+2(h_0-1),  l_x>0 \\  k_t + l_t = n-h_0, l_t>0}} \mathcal{S}^{k_x} \mathcal{T}^{k_t} LH_\rho(l_x, l_t) \\
&+ \kappa^{2(h_0-1)} \mathcal{\tilde A} \sum_{   \substack{k_x + l_x = m+2(h_0-1),  l_x>0}} \mathcal{S}^{k_x} \mathcal{T}^{n-h_0} RHS(l_x)\\
&+ \sum_{   \substack{i=1, \dots, h_0-1, \\k_x + l_x = m+2i,  l_x>0\\  k_t + l_t = n-h_0, l_t>0}}   \kappa^{2i} \mathcal{\tilde A}^{h_0-i} \mathcal{S}^{k_x} \mathcal{T}^{k_t} LH_\rho(l_x, l_t) \\
&+ \sum_{   \substack{i=1, \dots, h_0-1, \\k_x + l_x = m+2i,  l_x>0}}   \kappa^{2i} \mathcal{\tilde A}^{h_0-i} \mathcal{S}^{k_x} \mathcal{T}^{n-h_0} RHS(l_x) \\
&+\sum_{j=1}^{h_0-1} \sum_{i=0}^{j-1} \kappa^{2i} \mathcal{\tilde A}^{j-i} \mathcal{T}^{n-j} \sum_{\substack{k_x + l_x = m+2i,  l_x>0}} \mathcal{S}^{k_x} RHS(l_x).
\end{aligned}
\end{equation}
Compare r.h.s.\ of \eqref{e:ener_est_gen_iass3_rr1} with $RLH_\rho(m,n,h_0)$ as defined by \eqref{e:ener_est_gen_conj_rhsr}. We see that the first r.h.s.\ term of \eqref{e:ener_est_gen_iass3_rr1} equals the first term of $RLH_\rho(m,n,h_0)$. The two terms of r.h.s.\ of \eqref{e:ener_est_gen_iass3_rr1} that are sums involving $l_t>0$ are estimated by the former of these two terms (the latter is $i=h_0-1$ summand in the former), which equals the second term of $RLH_\rho(m,n,h_0)$. The 'initial terms', i.e.\ the remaining three terms of r.h.s.\ of \eqref{e:ener_est_gen_iass3_rr1} that have $RHS(l_x)$ sum up to (a constant times) the last term of $RLH_\rho(m,n,h_0)$. So
altogether we see that the terms not involving forcing, which estimate $LH_\rho (m,n)$ are controlled by $RLH_\rho(m,n,h_0)$. 

We are left with estimating the 'forcing terms', namely
\[
\begin{aligned}
I_{\divr f} + II_{\divr f} + IV +V &= \kappa^{2(h_0-1)} LH_{\divr f} (m+2(h_0-1),n-h_0) \\
&+ \sum_{   \substack{i=0, \dots, h_0-2, \\k_x + l_x = m+2i, \\  k_t + l_t = n-(h_0-1), l_t>0 }} \kappa^{2i} \mathcal{\tilde A}^{(h_0-1)-i} \mathcal{S}^{k_x} \mathcal{T}^{k_t} LH_{\divr f} (l_x,l_t-1) \\
&+\sum_{j=1}^{h_0-1}  \kappa^{2(j-1)} LH_{\divr f} (m+2(j-1),n-j)\\
&+ \sum_{j=1}^{h_0-1}  \sum_{i=0}^{j-2} \kappa^{2i} \mathcal{\tilde A}^{j-i-1} \sum_{\substack{k_x + l_x = m+2i,\\k_t + l_t = n-j,}} \mathcal{S}^{k_x} \mathcal{T}^{k_t} LH_{\divr f} (l_x,l_t), 
\end{aligned}
\]
whose first and third r.h.s.\ term summed are equal to the first term of $RLH_{\divr f}(m,n,h_0)$ as defined in \eqref{e:ener_est_gen_conj_rhsforce}; shifting indices as before in the second term we can estimate it with 
\[
\sum_{   \substack{i=0, \dots, h_0-2, \\k_x + l_x = m+2i, \\  k_t + l_t = n -h_0, l_t>0 }} \kappa^{2i} \mathcal{\tilde A}^{(h_0-1)-i} \mathcal{S}^{k_x} \mathcal{T}^{k_t} LH_{\divr f} (l_x,l_t),
\]
which is $j=h_0$ summand of the second term of $RLH_{\divr f}(m,n,h_0)$, whereas the remaining to be estimated term is the $\sum_{j<h_0}$ part of the second term of $RLH_{\divr f}(m,n,h_0)$.

Altogether, we have verified 
\[
LH_\rho(m,n)  \le RLH_\rho(m,n,h_0) + RLH_{\divr f}(m,n,h_0)
\]
which is what we wanted for our inductive step. Consequently Proposition \ref{prop:energy_est_tech} is proven.
\subsubsection{Concluding remarks for proof of Proposition \ref{prop:energy_est}}
Take Proposition \ref{prop:energy_est_tech} with $h=n$, which gives ($l_t>0$ is now empty condition), since $LH_\rho (m+2n,0) \lesssim RHS(m+2n)$
\[
\begin{aligned}
LH_\rho(m,n)  &\lesssim \kappa^{2n} RHS(m+2n) +  \sum_{j=1}^{n} \sum_{i=0}^{j-1} \kappa^{2i} \mathcal{\tilde A}^{j-i} \mathcal{T}^{n-j} \sum_{\substack{k_x + l_x = m+2i,  l_x>0}} \mathcal{S}^{k_x} RHS(l_x) \\
&+
\sum_{j=1}^{n} \kappa^{2(j-1)} LH_{\divr f} (m+2(j-1),n-j)\\
&+\sum_{j=1}^{n} \sum_{i=0}^{j-2} \kappa^{2i} \mathcal{\tilde A}^{j-i-1} \sum_{\substack{k_x + l_x = m+2i,\\k_t + l_t = n-j,}} \mathcal{S}^{k_x} \mathcal{T}^{k_t} LH_{\divr f} (l_x,l_t).
\end{aligned}
\]
Via definitions of $\mathcal{\tilde A}$ and  of $LH_\rho(m,n)$ as in Section \ref{ssec:gf_repl}, it gives  \eqref{e:ener_gen_mixed}. 

Let us remark that the upper bound $|\bbeta| \le U$ is a careless estimate of the amount of derivatives needed: in the worst case, when estimating $U$-many advective derivatives, we need $2U$-many space derivatives in formula \eqref{e:ener_gen_mixed}, available via Assumption \ref{ass:scales_foor_ener}. Of course, all our estimates have a huge constant in $\lesssim$, which is however merely $U$-dependent.

\newpage

\section{Homogenization in inertial range}      \label{s:homInertial}

In this section, we prove the homogenization estimates from step $q$ to step $q+1$.

\subsection{Characterization of initial datum}  \label{ss:char_datum}

We introduce a few notions to characterize initial datum. We need two types of notions, to measure the initial datum at $L^2$ level, and to control the derivatives of initial datum up to order $14Q^3$.

At $L^2$-level, we define energy dissipation as follows.
\begin{definition}[Energy dissipation]
For a divergence-free vector field $u \in C^\infty ( \T^2 \times [0,1] )$ and $\kappa > 0$, we use $| \cdot |_{\mathfrak{D}(u,\kappa)}$ to denote energy dissipation for initial datum $\rho_{\ini} : \T^2 \rightarrow \R$
\begin{align}     \label{e:dissipDef}
      | \rho_{\ini} |_{\mathfrak{D}(u,\kappa)} := \bigg( \kappa \iint | \nabla \rho |^2 dxdt \bigg)^{\sfrac12}
      \quad \text{with} \quad 
      \left\{
      \begin{aligned}
            \partial_t \rho + u \cdot \nabla \rho - \kappa \Delta \rho =&\, 0, \\
            \quad \rho(\cdot,0) =&\, \rho_{\ini}.
      \end{aligned}
      \right.
\end{align}
\end{definition}
\noindent
Then we have the following trivial estimate
\begin{align}     \label{e:trivialDissipEst}
      | \rho_{\ini} |_{\mathfrak{D}(u,\kappa)} \leq \| \rho_{\ini} \|_{L^2}.
\end{align}

To control the (macroscopic) derivatives efficiently, we introduce \textit{composed initial datum}.

\begin{definition}[Composed initial datum]
Given some $n_* \in \N$ and $\rho_{\ini} : \T^2 \rightarrow \R$, define
\begin{equation}     \label{e:dIniDatNorm}
\begin{split}
      | \rho_{\ini} |_{\mathfrak{S}(q,n_*)} :=
            \inf \bigg\{&\, R > 0 \,\bigg|\,  \rho = \sum_{0 \leq k \leq n_*} \rho_{\ini}^{(k)};  \\
            &\, \| \mr \nabla_q^p \rho_{\ini}^{(k)} \|_{2} 
                  < \mr\lambda_q^{ - 2(k+p) b\gamma } R , \, \forall \, 1 \leq p \leq 14Q^3 - k, \, k \leq n_*
      \bigg\}.
\end{split}
\end{equation}
Naturally, $\mathfrak{S}(q,n_*)$ is the space of functions for which $| \cdot |_{\mathfrak{S}(q,n_*)}$ is finite.
Here, we call $\{\rho_{\ini}^{(k)}\}_k$ \textit{atoms} of the composed initial datum $\rho_{\ini}$. For brevity, we also introduce shorthand notations
\begin{equation}     \label{e:d3IniDatNorm}
\begin{split}
      \mathfrak{S}(q)  := \mathfrak{S} \big( q, Q^3 \big),     \quad
      \mathfrak{S}'(q) := \mathfrak{S} \big( q, 2Q^3 \big). 
\end{split}
\end{equation}
\end{definition}

\begin{remark}    \label{r:datumProj}
For any $\rho_{\ini} \in C^\infty(\T^2)$ and $\balpha$ only containing space derivatives, we have
\begin{align}
      | \mr \partial_q^{\balpha} \rho_{\ini} |_{\mathfrak{S}(q,n+[\balpha])}
            \leq&\, | \rho_{\ini} |_{\mathfrak{S}(q,n)}.       \label{e:4:datumProj}
\end{align}
\end{remark}

In composed initial datum, the bounds are weighted by certain scale in $k$. Below we also introduce an unweighted version.

\begin{definition}[Composed initial datum without weight]
For $\rho_{\ini} \in \mathfrak{S}'(q)$, we define
\begin{equation}     \label{e:unw_IniDatNorm}
\begin{split}
      | \rho_{\ini} |_{\mathfrak{P}(q)} :=
            \inf \Big\{ R > 0 \,\big|\, 
            \| \mr \nabla_q^p \rho_{\ini} \|_{2} 
                  < \mr\lambda_q^{ - 2 b p \gamma } R , \, \forall \, p \leq 12Q^3 \Big\}.
\end{split}
\end{equation}
\end{definition}

Comparing $\mathfrak{P}$ to $\mathfrak{S}'$, we have the following immediate consequence
\begin{align}    \label{e:2:unw_IniDatNorm} 
      | \rho_{\ini} |_{\mathfrak{P}(q_1)} \lesssim | \rho_{\ini} |_{\mathfrak{P}(q_2)} \lesssim | \rho_{\ini} |_{\mathfrak{S}'(q_2)},
            \quad \forall \, q_2 \leq q_1, \rho_{\ini} \in C^\infty(\T^2).
\end{align}

Every element of $\mathfrak{S}'(q)$ is a sum of atoms for which we have derivative estimates of different order. In the following definition and remark, we introduce an operator $\bS$ to classify these atoms into two classes: one with derivative estimates of higher order, denoted by $\mathfrak{S}(q)$, and the other with derivative estimates of possibly lower order, denoted by $\mathfrak{S}'(q)$.
\begin{definition}      \label{d:SProj}
We define a projection $\bS: \mathfrak{S}'(q) \rightarrow \mathfrak{S}(q)$ with
\begin{align}     
      | \bS \rho_{\ini} |_{\mathfrak{S}(q)} + | \rho_{\ini} - \bS \rho_{\ini} |_{\mathfrak{S}'(q)} 
            \lesssim&\, | \rho_{\ini} |_{\mathfrak{S}'(q)},       \label{e:d4IniDatNorm} \\
      | \rho_{\ini} - \bS \rho_{\ini} |_{\mathfrak{P}(q)} 
            \lesssim&\, \mr\lambda_q^{ -2 Q^3 b \gamma } | \rho_{\ini} |_{\mathfrak{S}'(q)}.    \label{e:d6IniDatNorm}
\end{align}
\end{definition}

\begin{remark}          \label{r:2:SProj}
This projection is designed to do some bookkeeping and one has multiple choices for defining it. In all instances where we shall use $\bS$, it acts on elements of form
\begin{align}     \label{e:4:r2_SProj}
      f = \mr \partial_q^{\balpha} \rho_{\ini}^{(k_0)} 
            \in \mathfrak{S}'(q)
\end{align}
for some function $\rho_{\ini} \in \mathfrak{S}(q)$ and its atoms $\{ \rho_{\ini}^{(k)} \}_k$. We choose $\bS$ such that
\begin{align*}
      \bS f := \left\{
      \begin{aligned}
            \mr \partial_q^{\balpha} \rho_{\ini}^{(k_0)},     \quad \text{if } [\balpha] + k_0 \leq Q^3, \\
            0,    \quad \text{if } [\balpha] + k_0 > Q^3.
      \end{aligned}
      \right.
            \in \mathfrak{S}(q).
\end{align*}
Then we have
\begin{align}     \label{e:2:datumProj}
      \| \mr \nabla_q^p (f - \bS f) \|_{L^2}
            \lesssim&\, \mr\lambda_q^{ -2 (Q^3+p) b \gamma } | \rho_{\ini} |_{\mathfrak{S}'(q)}.
\end{align}
When $\bS$ acts on linear combinations of terms with form \eqref{e:4:r2_SProj}, we extend the definition of $\bS$ linearly. This verifies \eqref{e:d6IniDatNorm}.

\end{remark}

For functions defined in $\T^2 \times [0,1]$, we introduce a technical condition for its behavior at time $t=0$, which essentially says the transport derivatives associated to $u_q$ can be conveniently expressed by space derivatives.

\begin{definition}[Transmission condition]
For functions $\rho_n: \T^2 \times [0,1] \rightarrow \R $, we say the collection $\{\rho_n\}_{0 \leq n \leq n_*}$ satisfies transmission condition at step $q$ if for any $p \leq 10Q$
\begin{align}
      \mr D_{q,t}^p \rho_n (x,0) = \sum_{a=0}^{n_*} \sum_{|\balpha|_x=1, \balpha \in \mathcal{I}_x }^{Q^3+Q^2+2p} c_{n,a}^{(p,\balpha)} \mr \partial_q^{\balpha} \rho_a (x,0),
\end{align}
for a set of constants $c_{n,a}^{(p,\balpha)}$ satisfies
\begin{align}
      | c_{n,a}^{(p,\balpha)} | \lesssim 1.
\end{align}
In $\lesssim$ is an absolute constant.
\end{definition}

\subsection{Main results of this section}       \label{ss:mainR_inerH}

In our inductive scheme from $u_q$ to $u_{q+1}$, we need to solve the operator $D_{q,t} - \kappa \Delta$ multiple times ($2Q^2$ times) with rather general forces. We also need to propagate a large collection of estimates. For this purpose, we introduce the following notions for estimates of chain type.

\begin{definition}[Chain estimate I]   \label{d:resolve}
Fix $\kappa>0$ and initial datum $\rho_{\ini} \in C^\infty (\T^2)$. Consider 
\begin{align}
      \{ s_n \}_{0 \leq n \leq 2Q^2} \subset&\, \R_{\geq 2},    \label{e:2:resolve} \\
      \{ g_n \}_{0 \leq n \leq 2Q^2} \subset&\, \mathcal{P}_q(2NQ),    \label{e:3:resolve} \\ 
      \{ \bomega_n \}_{0 \leq n \leq 2Q^2} \subset&\, \mathcal{I}_*,   
            \quad [\bomega_n] \leq 3Q     \label{e:4:resolve}
\end{align}
with
\begin{align}     \label{e:5:resolve}
      \vertiii{ g_n }_q \lesssim \kappa \mr \lambda_q^2 \lambda_{q}^{ - s_n b \gamma }, 
            \quad \forall \, n.
\end{align}
Define the functions $\{\rho_{n}\}_{0 \leq n \leq 2Q^2}$ to be the solutions of
\begin{equation}  \label{e:resolveRho}
\begin{split}
      D_{q,t} \rho_0 - \kappa \Delta \rho_0 =&\, 0, \\ 
      \rho_0 (\cdot, 0) =&\, \rho_{\ini}, \\ 
      \ldots&\, \\
      D_{q,t} \rho_{n+1} - \kappa \Delta \rho_{n+1} =&\, \mr\lambda_q^{-1} \divr( g_{n} \mr D_q^{\bomega_{n}} \rho_n ), \\ 
      \rho_{n+1} (\cdot, 0) =&\, 0, \\ 
      \ldots&\, 
\end{split}
\end{equation}
Then we define the notion $\resN{ \cdot }_{u_q,\kappa}$ via
\begin{equation*}
\begin{split}
      \resN{ \rho_{\ini} }_{u_q,\kappa} := \inf \bigg\{ R>0 \, \bigg| &\,
            \| \mr D_q^{\bbeta} \rho_n \|_\infty \leq
                  \lambda_q^{ - \big( [\bbeta] + \sum_{i=0}^{n-1} (s_i-1) \big) b \gamma } R, 
                  \quad \forall\, \bbeta \in \mathcal{I}_*, [\bbeta] \leq 11Q^3 - 3Qn, \\ 
            &\, \text{for any } n, \{s_n\}_n, \{g_n\}_n, \{\bomega_n\}_n, \{\rho_n\}_n \text{ satisfying }
                  \eqref{e:2:resolve}-\eqref{e:resolveRho}
      \bigg\}
\end{split}
\end{equation*}
We remind that, according to Section \ref{sss:ntt:norm}, $\| \cdot \|$ denotes $L_t^\infty L_x^2$-norm and $\| \cdot \|_\infty$ denotes $L_{x,t}^\infty$ norms.

To propagate $\alpha_0$-H\"older estimates, we also define the notion $\resN{ \cdot }^{(H,\alpha_0)}_{u_q,\kappa}$ via
\begin{equation*}
\begin{split}
      \resN{ \rho_{\ini} }^{(H,\alpha_0)}_{u_q,\kappa} := \inf \bigg\{ &\, R>0 \, \bigg| \,
            \| \rho_n \|_{L^\infty_t C^{\alpha_0}_x} \leq
                  \lambda_q^{ - \big( \sum_{i=0}^{n-1} (s_i-1) \big) b \gamma } R, \\ 
            &\, \text{for any } n, \{s_n\}_n, \{g_n\}_n, \{\bomega_n\}_n, \{\rho_n\}_n \text{ satisfying }
                  \eqref{e:2:resolve}-\eqref{e:resolveRho}
      \bigg\}
\end{split}
\end{equation*}
We extend above definitions to $v_q$, denoted by $\resN{ \cdot }_{v_q,\kappa}$ and $\resN{ \cdot }^{(H,\alpha_0)}_{v_q,\kappa}$.
\end{definition}

Now we can state the main results of this section.

\begin{theorem}   \label{t:homDissip}
Fix some $q_* \leq q$. Suppose that, for any $\rho_{\ini} \in C^\infty (\T^2)$, we have
\begin{align}     
      \resN{ \rho_{\ini} }_{u_q,\kappa_q} 
            \leq&\, \lambda_q^{-\alpha} \bigg( | \rho_{\ini} |_{\mathfrak{D}(u_q,\kappa_q)} 
                  + \sum_{j=q_*}^{q} \lambda_j^{-\sfrac{\gamma}{2}} | \rho_{\ini} |_{\fS(j)} \bigg),   \label{e:1:homDissip} 
\end{align}
then
\begin{align}     \label{e:2:homDissip}
      \left| | \rho_{\ini} |_{\mathfrak{D}(u_q,\kappa_q)} 
            - | \rho_{\ini} |_{\mathfrak{D}(u_{q+1},\kappa_{q+1})} \right|
            \leq&\, \lambda_{q+1}^{-\gamma} \bigg( | \rho_{\ini} |_{\mathfrak{D}(u_q,\kappa_q)} 
                  + \sum_{j=q_*}^{q+1} \lambda_j^{-\sfrac{\gamma}{2}} | \rho_{\ini} |_{\fS(j)} \bigg).
\end{align}
\end{theorem}

\begin{theorem}   \label{t:homRslv}
Fix some $q_* \leq q$. Suppose that, for any $\rho_{\ini} \in C^\infty (\T^2)$, we have
\begin{align}     
      \resN{ \rho_{\ini} }_{u_q,\kappa_q} 
            \leq&\, \lambda_q^{-\alpha} \bigg( | \rho_{\ini} |_{\mathfrak{D}(u_q,\kappa_q)} 
                  + \sum_{j=q_*}^{q} \lambda_j^{-\sfrac{\gamma}{2}} | \rho_{\ini} |_{\fS(j)} \bigg),   \label{e:2:homRslv} \\ 
      \resN{ \rho_{\ini} }^{(H,\alpha_0)}_{u_q,\kappa_q} 
            \leq&\, \left( 1 - \frac{4}{q} \right) \bigg( | \rho_{\ini} |_{\mathfrak{D}(u_q,\kappa_q)} 
                  + \sum_{j=q_*}^{q} \lambda_j^{-\sfrac{\gamma}{2}} | \rho_{\ini} |_{\fS(j)} \bigg), \label{e:4:homRslv}
\end{align}
then for any $\rho_{\ini} \in C^\infty (\T^2)$, we have
\begin{align}
      \resN{ \rho_{\ini} }_{u_{q+1},\kappa_{q+1}} 
            \leq&\, \lambda_{q+1}^{-\alpha} \bigg( | \rho_{\ini} |_{\mathfrak{D}(u_{q+1},\kappa_{q+1})} 
                  + \sum_{j=q_*}^{q+1} \lambda_j^{-\sfrac{\gamma}{2}} | \rho_{\ini} |_{\fS(j)} \bigg),      \label{e:6:homRslv} \\
      \resN{ \rho_{\ini} }^{(H,\alpha_0)}_{u_{q+1},\kappa_{q+1}} 
            \leq&\, \left( 1 - \frac{4}{q+1} \right) \bigg( | \rho_{\ini} |_{\mathfrak{D}(u_{q+1},\kappa_{q+1})}
                  + \sum_{j=q_*}^{q+1} \lambda_j^{-\sfrac{\gamma}{2}} | \rho_{\ini} |_{\fS(j)} \bigg).     \label{e:8:homRslv}
\end{align}
\end{theorem}

\begin{remark}    \label{r:linearityHomRslv}
From the definitions of $\resN{ \cdot }_{\cdot,\cdot}$ and $\resN{ \cdot }^{(H,\alpha_0)}_{\cdot,\cdot}$, we see that the estimates in Theorem \ref{t:homRslv} and Theorem \ref{t:homDissip} are linear in $\rho_{\ini}$.
\end{remark}

The proofs of above theorems are postponed to Section \ref{ss2:estQ+1}. To validate the estimates (\ref{e:6:homRslv}-\ref{e:8:homRslv}) for $D_{q+1,t} - \kappa_{q+1} \Delta$, we need two ingredients. One ingredient is unconditional energy estimates in Section \ref{ss:unconditional_inertial}. The other ingredient is homogenization which deduces estimates on $D_{q+1,t} - \kappa_{q+1} \Delta$ from \eqref{e:2:homRslv} on $D_{q,t} - \kappa_{q} \Delta$. This homogenization scheme is implemented in Section \ref{ss:stateExpInertial}, Section \ref{ss:proof_expansion}, Section \ref{ss:reEXP1}, Section \ref{ss:generalExpF} and Section \ref{ss2:estQ+1}.

For technical convenience, we also introduce the following notion for chain estimates.

\begin{definition}[Chain estimate II]   \label{d:FResolve}
Fix $\kappa>0$, $0 \leq n_* \leq 2Q^2$ and $f \in C^\infty(\T^2 \times [0,1])$. Consider $\{ s_n \}$, $\{ g_n \}$ and $\{ \bomega_n \}$ satisfying \eqref{e:2:resolve}-\eqref{e:5:resolve}. Define the functions $\{\rho_{n}\}_{1 \leq n \leq 2Q^2-n_*}$ to be the solutions of
\begin{equation}  \label{e:FResolveRho}
\begin{split}
      D_{q,t} \rho_1 - \kappa \Delta \rho_1 =&\, \mr\lambda_q^{-1} \divr( g_0 \mr D_q^{\bomega_0} f ), \\ 
      \rho_0 (\cdot, 0) =&\, 0, \\ 
      \ldots&\, \\
      D_{q,t} \rho_{n+1} - \kappa \Delta \rho_{n+1} =&\, \mr\lambda_q^{-1} \divr( g_{n} \mr D_q^{\bomega_{n}} \rho_n ), \\ 
      \rho_{n+1} (\cdot, 0) =&\, 0, \\ 
      \ldots&\, 
\end{split}
\end{equation}
Then we define the notion $\resF{ \cdot }_{u_q,\kappa,n_*}$ via
\begin{equation*}
\begin{split}
      \resF{ f }_{u_q,\kappa,n_*} := \inf \bigg\{ &\, R>0 \, \bigg| \,
            \| \mr D_q^{\bbeta} f \|_\infty \leq 
            \lambda_q^{ - [\bbeta] b \gamma } R,
            \quad \forall\, \bbeta \in \mathcal{I}_*, [\bbeta] \leq 11Q^3 - 3Qn_*,      \\
            &\,\| \mr D_q^{\bbeta} \rho_n \|_\infty \leq
            \lambda_q^{ - \big( [\bbeta] + \sum_{i=0}^{n-1} (s_i-1) \big) b \gamma } R, 
            \quad \forall\, \bbeta \in \mathcal{I}_*, [\bbeta] \leq 11Q^3 - 3Q(n+n_*),    \\ 
            &\, \text{for any } n, \{s_n\}_n, \{g_n\}_n, \{\bomega_n\}_n, \{\rho_n\}_n \text{ satisfying }
                  \eqref{e:2:resolve}-\eqref{e:5:resolve}, \eqref{e:FResolveRho}
      \bigg\}
\end{split}
\end{equation*}
Similarly, we also define the alternative for propagating H\"older estimates
\begin{equation*}
\begin{split}
      \resF{ f }^{(H,\alpha_0)}_{u_q,\kappa,n_*} := \inf \bigg\{ &\, R>0 \, \bigg| \,
            \| f \|_{L^\infty_t C^{\alpha_0}_x} \leq 
            \lambda_q^{ - [\bbeta] b \gamma } R,      
            \| \rho_n \|_{L^\infty_t C^{\alpha_0}_x} \leq
            \lambda_q^{ - \big( \sum_{i=0}^{n-1} (s_i-1) \big) b \gamma } R,    \\ 
            &\, \text{for any } n, \{s_n\}_n, \{g_n\}_n, \{\bomega_n\}_n, \{\rho_n\}_n \text{ satisfying }
                  \eqref{e:2:resolve}-\eqref{e:5:resolve}, \eqref{e:FResolveRho}
      \bigg\}
\end{split}
\end{equation*}
We also extend above definitions to $v_q$, denoted by $\resF{ \cdot }_{v_q,\kappa,n_*}$ and $\resF{ \cdot }^{(H,\alpha_0)}_{v_q,\kappa,n_*}$.
\end{definition}

We have following immediate consequences for the notions above.

\begin{lemma}    \label{l:P0resolveToF}
For any $n_1 \leq n_2$, we have 
\begin{align*}
      \resN{ \rho_{\ini} }_{u_q,\kappa} \leq&\, \resN{ \rho_{\ini} }_{u_q,\kappa},  \\ 
      \resF{ f }_{u_q,\kappa,n_2} \leq&\, \resF{ f }_{u_q,\kappa,n_1}.
\end{align*}
\end{lemma}
This lemma is a direct consequence of the definitions.

\begin{lemma}    \label{l:P1resolveToF}
For $ \rho_{\ini} \in C^\infty (\T^2) $, we define $\{ \rho_n \}_n$ as in \eqref{e:resolveRho}. Then, for any $0 \leq n \leq 2Q^2$, we have
\begin{align*}
      \resF{ \rho_n }_{u_q,\kappa,n} 
            \leq&\, \lambda_q^{ - \big( \sum_{i=0}^{n-1} (s_i-1) \big) b \gamma } \resN{ \rho_{\ini} }_{u_q,\kappa}.
\end{align*}
\end{lemma}

\begin{proof}
Comparing Definition \ref{d:FResolve} and Definition \ref{d:resolve} concludes the proof.
\end{proof}

\begin{lemma}    \label{l:P2resolveToF}
Suppose $\rho_{\ini} \in C^\infty (\T^2)$ and $g \in \mathcal{P}_q(2NQ)$ satisfies
\begin{align*}
      \vertiii{ g }_q \leq \kappa \mr \lambda_q^2 \lambda_{q}^{-s_*b\gamma}
\end{align*}
for some $s_* \in \R_{\geq 2}$. Define
\begin{equation}  \label{e:2:P2resolveToF}
\begin{split}
      D_{q,t} \rho_* - \kappa \Delta \rho_* 
            =&\, \mr\lambda_q^{-1} \divr \big( g \mr D_q^{\bomega} \rho_F \big), \\ 
      \rho_F (\cdot, 0) =&\, \rho_{\ini}
\end{split}
\end{equation}
for some $\bomega \in \mathcal{I}_*$ with $[\bomega] \leq 3Q$, then we have
\begin{align}     \label{e:6:P2resolveToF}
      \resF{ \rho_* }_{u_q,\kappa,n_0+1} 
            \leq&\, 
            \resF{ \rho_F }_{u_q,\kappa,n_0} \lambda_{q}^{-(s_*-1)b\gamma}
            +  \resN{ \rho_{\ini} }_{u_q,\kappa}.
\end{align}
\end{lemma}
The proof of this lemma is immediate and analogous to Lemma \ref{l:P1resolveToF}.

\begin{remark}
Analogous results to Lemma \ref{l:P0resolveToF}, Lemma \ref{l:P1resolveToF} and Lemma \ref{l:P2resolveToF} hold for H\"older estimates $\resN{ \rho_{\ini} }^{(H,\alpha_0)}_{u_q,\kappa}$ and $\resF{ f }^{(H,\alpha_0)}_{u_q,\kappa,n_*}$ and the $v_q$ alternatives.
\end{remark}

\subsection{Notations of this section}    \label{ss:nota_inerH}

In this section, we introduce the following simplified notations.

\begin{definition}[Parameters]   \label{d:homParaI} 
Define
\begin{align}     \label{d:mainPara}
      \lambda := \lambda_{q+1},            \quad
      \delta := \delta_{q+1},            \quad
      \mu := \mu_{q+1},            \quad 
      \bar\lambda := \lambda_{q},      \quad
      \bar\delta := \delta_{q},            \quad
      \bar\mu := \mu_{q},  
\end{align}
\begin{align}     \label{d:circleLambda}
      \mr {\lambda} := \mr {\lambda}_{q+1},      \quad
      \mr {\mu} := \mr {\mu}_{q+1},       \quad
      \mr {\bar\lambda} := \mr {\lambda}_{q},      \quad
      \mr {\bar\mu} := \mr {\mu}_{q},  
\end{align} 
\begin{align}     \label{d:eddyDiff}  
      \bar \kappa := \kappa \bigg( 1 + \frac{\delta}{\kappa^2\lambda^2} \bigg),     \quad 
      \varepsilon := \frac{\bar\delta^{\sfrac12} \bar\lambda^{1+2b\gamma_R}} {\mu}.
\end{align}
\end{definition}

\begin{definition}[Coefficient fields]   \label{d:homFlowmap}
Recall the velocity fields $u_q$ and $v_q$, the flow map $\bar\Phi_q$ and $\Phi_q$ constructed in Section \ref{ss:vectorfield} and Lemma \ref{l:uEstimate}. Define the following functions,
\begin{align*}
      &\,\bar u := u_q,          \quad
      u := u_{q+1},           \quad
      v := v_{q+1},           \quad
      \Phi := \Phi_q,      \quad
      \phi_* := \phi_N,       \quad
      \omega_m := \omega_{q,m},     \\
      &\,B_m := B_{q,m},               \quad 
      E_m := E_{q,m},               \quad
      S_m := S_{q,m},            \quad
      \Omega_m := \Omega_{q,m},     \quad
      \omega_* := \lambda^3 \omega_{q,N}.
\end{align*}
We introduce the following shorthand notations for operators, quantities and spaces,
\begin{align*}
      O := O_{q+1},                 \quad 
      \mr O := \mr O_{q+1},         \quad 
      \bar D :=&\, D_{q},           \quad 
      \mr {\bar D} := \mr D_{q},    \\
      \bar{\mathcal{P}} := \mathcal{P}_q,       \quad
      \mathcal{P} := \mathcal{P}_{q+1},         \quad 
      \vertiii{\cdot} :=&\, \vertiii{\cdot}_{q},   \quad 
      \vertiii{\cdot}_{+1} := \vertiii{\cdot}_{q+1}.
\end{align*}
We also define $Z, \Omega, \Xi: \T^2 \times [0,1] \times \T$ via
\begin{align}
      \Omega (x, t, \tau) :=&\, \omega_*(x,t) \phi_*(\tau),       \label{e:6:homFlowmap}  \\ 
      Z (x, t, \tau) :=&\, \sum_{m} \omega_m(x,t) z_m(\tau),      \label{e:8:homFlowmap}  \\ 
      \Xi (x, t, \tau) :=&\, \sum_{m} S_m (x,t) \vartheta_m (\tau).       \label{e:10:homFlowmap} 
\end{align}
In the whole Section \ref{s:homInertial}, $\sum_m$ always means $\sum_{m=1}^N$ for the integer $N$ defined in Section \ref{ss:parameters}.
\end{definition}

Consider the following equation for $\varrho: \T^2 \times [0,1] \rightarrow \R$,
\begin{equation}  \label{e:q+1:meq}
\begin{split}
      O_{q+1,t} \varrho - \kappa_{q+1} \Delta \varrho =&\, 0, \\ 
      \varrho (\cdot, 0) =&\, \rho_{\ini}.
\end{split}
\end{equation}
From above definitions, we can rewrite \eqref{e:q+1:meq} in the following form
\begin{equation} \label{e:hom:meq}
\begin{split}
    \bar D_t \varrho - \divr ( A \nabla \varrho ) =& 0, \\ 
    \varrho (\cdot, 0) =& \rho_{\text{in}}, 
\end{split}
\end{equation}
with $A: \T^2 \times [0,1] \rightarrow \R^{2 \times 2}$ given by
\begin{align}   \label{e:hom:mat} 
      A(x,t) :=&\, \kappa \Id + \frac{\delta^{\sfrac12}}{\lambda} \adj \nabla \Phi (x,t) H ( \lambda \Phi(x,t), \mu t ) \adj \nabla \Phi^{T} (x,t)
\end{align}
with $H$ an anti-symmetric matrix
\begin{align}
      H_{12}(\xi, \tau) =&\, \eta_1(\tau) \sin(2\pi \xi_1) + \eta_2(\tau) \sin(2\pi \xi_2) 
            = - H_{21}(\xi, \tau),  \label{e:hom:MatEnt} \\ 
      \Pi_1(\xi) =&\, \sin(2\pi \xi_1), \quad \Pi_2(\xi) = \sin(2\pi \xi_2).  \label{e:hom:sin}
\end{align}

\begin{remark}
From above definitions, we have
\begin{align}     \label{e:6:homParaI}
      \gamma > 0, \quad
      \varepsilon < 1, \quad 
      \lambda_q < \lambda, \quad
      \delta_q > \delta.
\end{align}
From \eqref{e:4:mScaleRela}, we have
\begin{align}     \label{d:lambda_gamma}
      \lambda^{\gamma_S} \geq 
      \max \Bigg\{ 
            \frac{\mu} {\delta^{\sfrac12}\mr{\bar\lambda}}, 
            \frac{\varepsilon \kappa\lambda^2} {\delta^{\sfrac12}\mr{\bar\lambda}}, 
            \frac{\mr{\bar\mu}} {\varepsilon\delta^{\sfrac12}\mr{\bar\lambda}},
            \frac{\kappa\lambda^2} {\varepsilon\delta^{\sfrac12}\lambda}
      \Bigg\}.
\end{align}
From \eqref{e:smallGammaR} and \eqref{e:2:mScaleRela}, we have
\begin{align}     \label{e:homGammaRela}
      \frac{ \delta^{\sfrac12}\mr{\bar\lambda} } {\kappa\lambda^2} 
            \cdot \lambda^{2b\gamma+\gamma_S} \leq 1.
\end{align}
From \eqref{e:mrParameter} and \eqref{e:14:mScaleRela}, we also have
\begin{align}     \label{d:bigMu}
      \frac{\delta^{\sfrac12}\mr{\bar\lambda}}{\mu} \leq \lambda^{-\gamma}. 
\end{align}
\end{remark}

\begin{remark}    \label{r:2:homFlowmap}
Recalling Lemma \ref{l:uEstimate}, above definitions lead to
\begin{align}     
      \bar D_t \Phi =&\, Z(x, t, \mu t) + \lambda^{-3} \Omega (x, t, \mu t),  \label{e:12:homFlowmap}  \\ 
      \Xi (x, t, \mu t) =&\, \frac{\kappa}{\varepsilon} \big( \nabla \Phi \nabla \Phi^T - (\det \nabla \Phi)^2 \Id \big) = \sum_{m=1}^N S_m (x,t) \vartheta_m (\mu t),       \label{e:14:homFlowmap} \\ 
      \adj \nabla \Phi \det \nabla \Phi =&\, \Id + \varepsilon \sum_{m=1}^N E_m \varphi_m (\mu t), \quad 
      \nabla \Phi = \Id + \varepsilon \sum_{m=1}^N B_m \phi_m (\mu t),    \label{e:16:homFlowmap} \\ 
      ( \adj \nabla \Phi H(\xi, \tau) &\, \adj \nabla \Phi^{T} )_{ij}
            = H_{ij}(\xi, \tau) + \varepsilon \sum_{m=1}^N \sigma_m(\tau) \Omega_{m,ij} H_{12}(\xi, \tau).     \label{e:18:homFlowmap}
\end{align}
From Lemma \ref{l:uEstimate}, we have, for any $(x,t) \in \T^2 \times [0,1]$,
\begin{align}
      \supp Z(x,t,\cdot) \cap \supp \eta_1 = \supp Z(x,t,\cdot) \cap \supp \eta_2 = \varnothing.      \label{e:20:homFlowmap}
\end{align}
Recalling \eqref{e:10:uEstimate}, we have, from \eqref{e:hom:mat} and \eqref{e:14:homFlowmap},
\begin{align}   \label{e:hom:mat_alt} 
      A(x,t) = \adj \nabla \Phi \bigg( \kappa \Id + \frac{\delta^{\sfrac12}}{\lambda} H ( \lambda \Phi, \mu t )  \bigg) \adj \nabla \Phi^{T} 
            + \kappa \big( \Id - \adj \nabla \Phi \adj \nabla \Phi^T \big) .
\end{align}
We also have
\begin{align}     \label{e:24:homFlowmap}
      \Phi(x,0) = x.
\end{align}
\end{remark}

\begin{remark}[Estimates on velocity]   \label{r:homParaII}
Under Definition \ref{d:homParaI} and Definition \ref{d:homFlowmap}, we deduce from Lemma \ref{l:uEstimate}
\begin{align}
      \omega_N, \omega_* \in ( \mathcal{\bar P}(N) )^2, \quad
            B_m, E_m, S_m, \Omega_m \in&\, ( \mathcal{\bar P}(N) )^{2 \times 2}.  \label{e:mtrStrHom}
\end{align}
Moreover, for any $\balpha$ with $[\balpha] \leq 8Q^3$ and $p$ with $p \leq 8Q^3$,
\begin{align}
      \vertiii{ \mr{\bar D}^{\balpha} B_m }, \vertiii{ \mr{\bar D}^{\balpha} E_m }, \vertiii{ \mr{\bar D}^{\balpha} S_m }, \vertiii{ \mr{\bar D}^{\balpha} \Omega_m }, \vertiii{ \mr{\bar D}^{\balpha} \omega_m }, \vertiii{ \mr{\bar D}^{\balpha} \omega_* } \lesssim&\, 1,       
            \label{e:mtrEstHom} \\ 
      \| \partial_\tau^p \phi_m \|_{\infty}, \| \partial_\tau^p \phi_* \|_{\infty}, \| \partial_\tau^p \varphi_m \|_{\infty}, \| \partial_\tau^p z_m \|_{\infty}, \| \partial_\tau^p \sigma_m \|_{\infty}, \| \partial_\tau^p \vartheta_m \|_{\infty} \lesssim&\, 1.
            \label{e:prdEstHom}
\end{align}
\end{remark}

\subsection{Energy and stability estimates for chains}     \label{ss:unconditional_inertial}

In this section, we prove unconditional energy estimates for equations like \eqref{e:resolveRho}. We also prove stability estimates between $D_{q+1,t} - \kappa_{q+1} \Delta$ and $O_{q+1,t} - \kappa_{q+1} \Delta$.

\begin{assumption}[Step $q+1$ with $\mathcal{O}$ coefficient fields]   \label{a:OEq}
For initial datum $\rho_{\ini} \in C^\infty (\T^2)$, suppose that $\{ g_n \}_{0 \leq n \leq 2Q^2} \subset \mathcal{O}_{q+1}(2NQ)$ and $\{ \bomega_n \}_{0 \leq n \leq 2Q^2} \subset \mathcal{I}_*$ with $[\bomega_n] \leq 3Q$ satisfying
\begin{align}     \label{e:2:OEq_Idct}
      \vertiii{ g_n }_{q+1} \lesssim&\, \kappa \mr \lambda^2 \lambda^{-s_nb\gamma}
\end{align}
for some $\{s_n\}_{0 \leq n \leq 2Q^2} \subset \R_{\geq 2}$. Define the functions $\{\varrho_{n}\}_{0 \leq n \leq 2Q^2}$ to be the solutions to
\begin{align}
      O_t \varrho_0 - \kappa \Delta \varrho_0 =&\, 0,       \label{e:5:OEq_Idct} \\ 
      \varrho_0 (\cdot, 0) =&\, \rho_{\ini},    \label{e:6:OEq_Idct} \\ 
      \ldots,     \nonumber \\ 
      O_t \varrho_{n+1} - \kappa \Delta \varrho_{n+1} =&\,\mr \divr ( g_n \mr O^{\bomega_n} \varrho_n ),       \label{e:7:OEq_Idct1} \\ 
      \varrho_{n+1} (\cdot, 0) =&\, 0.    \label{e:8:OEq_Idct}
\end{align}
\end{assumption}

For $\{\varrho_n\}$ defined in Assumption \ref{a:OEq}, we have the following energy estimates.

\begin{proposition}     \label{p:unconEnergyHom}
Taking Assumption \ref{a:OEq}, we have that, for $\bbeta \in \mathcal{I}_*$ with $[\bbeta] \leq 12Q^3 - 3Qn$,
\begin{equation}     \label{e:2:unconEnergyHo} 
\begin{split}
      \| \mr O^{\bbeta} \varrho_n \|
            + \kappa^{\sfrac12} \mr\lambda \| \mr O^{\bbeta} \mr \nabla \varrho_n \|_2
            \lesssim \lambda^{ - \big( 2[\bbeta] + \sum_{i=0}^{n-1} s_i \big) b\gamma } 
            \left( | \rho_{\ini} |_{\mathfrak{D}(v,\kappa)} + | \rho_{\ini} |_{\mathfrak{P}(q+1)} \right).
\end{split}
\end{equation}
From interpolation, we also have
\begin{align}     
      \resN{ \rho_{\ini} }_{v,\kappa} 
            \lesssim&\, \mr \lambda^{1-2\gamma} 
                  \left( | \rho_{\ini} |_{\mathfrak{D}(v,\kappa)} + | \rho_{\ini} |_{\mathfrak{P}(q+1)} \right),    \label{e:i1:unconEnergyHo} \\ 
      \resN{ \rho_{\ini} }_{v,\kappa}^{(H,\alpha_0)} 
            \lesssim&\, \mr \lambda^{(1-2\gamma)(1+\alpha_0)} 
                  \left( | \rho_{\ini} |_{\mathfrak{D}(v,\kappa)} + | \rho_{\ini} |_{\mathfrak{P}(q+1)} \right).    \label{e:i2:unconEnergyHo}
\end{align}
Moreover, for $\bbeta \in \mathcal{I}_*$ with $[\bbeta] \leq 11Q^3 - 3Qn$,
\begin{align}
      \| \mr O^{\bbeta} \varrho_n \|_\infty
            \lesssim \mr\lambda^{1-2\gamma} \lambda^{ - \big( 2[\bbeta] + \sum_{i=0}^{n-1} s_i \big) b\gamma }
                  \left( | \rho_{\ini} |_{\mathfrak{D}(v,\kappa)} + | \rho_{\ini} |_{\mathfrak{P}(q+1)} \right).    
                  \label{e:4:unconEnergyHo}
\end{align}
\end{proposition}

\begin{remark}    \label{r:energySobolevLoss}
The above result is presented for $\T^2$. For higher space dimensions, the additional power $\mr \lambda^{1-2\gamma}$ in \eqref{e:i2:unconEnergyHo} and \eqref{e:4:unconEnergyHo} changes.
\end{remark}

\begin{proof}
From the smallness of $\gamma$ given by \eqref{e:smallGammaR}, we apply energy estimates in Lemma \ref{l:L2_init_energy} to deduce that, for $\bbeta \in \mathcal{I}_*$ with $[\bbeta] \leq 12Q^3$,
\begin{align*}
      \| \mr O^{\bbeta} \varrho_0 \|
            + \kappa^{\sfrac12} \mr\lambda \| \mr O^{\bbeta} \mr \nabla \varrho_0 \|_2
            \lesssim \lambda^{-2[\bbeta]b\gamma} \left( | \rho_{\ini} |_{\mathfrak{D}(v,\kappa)} + | \rho_{\ini} |_{\mathfrak{P}(q+1)} \right).
\end{align*}
From \eqref{e:2:OEq_Idct} and Corollary \ref{c:L2_FHom_energy}, we have $\bbeta \in \mathcal{I}_*$ with $[\bbeta] \leq 12Q^3 - 3Q$,
\begin{equation*}
\begin{split}
      \| \mr O^{\bbeta} \varrho_1 \|
            + \kappa^{\sfrac12} \mr\lambda \| \mr O^{\bbeta} \mr \nabla \varrho_1 \|_2
            \lesssim \lambda^{ - ( 2[\bbeta]+2[\bomega_0]-2+s_0 ) b\gamma } \left( | \rho_{\ini} |_{\mathfrak{D}(v,\kappa)} + | \rho_{\ini} |_{\mathfrak{P}(q+1)} \right).
\end{split}
\end{equation*}
Inductively, we have for $\bbeta \in \mathcal{I}_*$ with $[\bbeta] \leq 12Q^3 - 3Qn$,
\begin{equation*}
\begin{split}
      \| \mr O^{\bbeta} \varrho_n \|
            + \kappa^{\sfrac12} \mr\lambda \| \mr O^{\bbeta} \mr \nabla \varrho_n \|_2
            \lesssim \lambda^{ - \big( 2[\bbeta] + \sum_{i=0}^{n-1} (s_i + 2[\bomega_i] - 2) \big) b\gamma } \left( | \rho_{\ini} |_{\mathfrak{D}(v,\kappa)} + | \rho_{\ini} |_{\mathfrak{P}(q+1)} \right).
\end{split}
\end{equation*}

To deduce \eqref{e:i1:unconEnergyHo} and \eqref{e:4:unconEnergyHo}, we apply the interpolation inequality in Lemma \ref{l:a:interpolate} to obtain
\begin{align}     \label{e:12:unconEnergyHo}
      \| \nabla^p \varrho_n \|_{L^\infty_{x,t}}
            \lesssim&\, \| \nabla^{k+p} \varrho_n \|^s_{ L^\infty_t L^2_x }
                  \| \nabla^p \varrho_n \|^{1-s}_{ L^\infty_t L^2_x },
\end{align}
with $s = Q^{-3}$, $k = Q^3$ and $p \leq 11Q^3 - 3Qn$. From \eqref{e:2:unconEnergyHo} and \eqref{e:12:unconEnergyHo}, we have
\begin{align}     \label{e:14:unconEnergyHo}
      \| \mr\nabla^p \varrho_n \|_{L^\infty_{x,t}}
            \lesssim&\, \mr \lambda^{1-2\gamma} \lambda^{ - \big( 2p + \sum_{i=0}^{n-1} s_i \big) b\gamma } 
            \left( | \rho_{\ini} |_{\mathfrak{D}(v,\kappa)} + | \rho_{\ini} |_{\mathfrak{P}(q+1)} \right).
\end{align}
When $\bbeta$ contains time derivatives, the following $L_{x,t}^\infty$ estimate follows from the equation \eqref{e:5:OEq_Idct} or \eqref{e:7:OEq_Idct1}, commutator equalities and \eqref{e:14:unconEnergyHo},
\begin{align}     \label{e:16:unconEnergyHo}
      \| \mr O^{\bbeta} \varrho_n \|_{L^\infty_{x,t}}
            \lesssim&\, \mr \lambda^{1-2\gamma} \lambda^{ - \big( 2[\bbeta] + \sum_{i=0}^{n-1} s_i \big) b\gamma } 
            \left( | \rho_{\ini} |_{\mathfrak{D}(v,\kappa)} + | \rho_{\ini} |_{\mathfrak{P}(q+1)} \right).
\end{align}

To deduce \eqref{e:i2:unconEnergyHo}, we apply the interpolation inequality \eqref{e:l:8:interpolate} in Lemma \ref{l:a:interpolate}.

\end{proof}

\begin{assumption}[Step $q+1$ with $\mathcal{P}$ coefficient fields]   \label{a:PEq}
Fix some initial datum $\rho_{\ini} \in C^\infty (\T^2)$. Assume $\{ \bomega_n \}_{0 \leq n \leq 2Q^2} \subset \mathcal{I}_*$ with $[\bomega_n] \leq 3Q$ and $\{ h_n \}_{0 \leq n \leq 2Q^2} \subset \mathcal{P}_{q+1}(2NQ)$ satisfying
\begin{align}     \label{e:2:PEq_Idct}
      \vertiii{ h_n }_{q+1} \lesssim&\, \kappa \mr \lambda^2 \lambda^{-s_nb\gamma}
\end{align}
for some $\{s_n\}_{0 \leq n \leq 2Q^2} \subset \R_{\geq 2}$ Define the functions $\{\ddot\varrho_{n}\}_{0 \leq n \leq 2Q^2}$ to be the solutions to
\begin{align*}
      D_t \ddot{\varrho}_0 - \kappa \Delta \ddot{\varrho}_0 =&\, 0,     \quad 
      \ddot{\varrho}_0 (\cdot, 0) = \rho_{\ini}, \\ 
      \ldots,     \nonumber \\ 
      D_t \ddot{\varrho}_{n+1} - \kappa \Delta \ddot{\varrho}_{n+1} =&\,\mr \divr ( h_n \mr D^{\bomega_n} \ddot{\varrho}_n ),       \quad 
      \ddot{\varrho}_{n+1} (\cdot, 0) = 0. 
\end{align*}
\end{assumption}

We also have the following energy estimates for $\{\ddot\varrho_n\}$ defined in Assumption \ref{a:PEq}. The proof is the same as the proof of Proposition \ref{p:unconEnergyHom}. We omit it here.

\begin{proposition}     \label{p:unconEnergyP}
Taking Assumption \ref{a:PEq}, we have that, for $\bbeta \in \mathcal{I}_*$, $[\bbeta] \leq 12Q^3 - 3Q n$,
\begin{align}     
      \| \mr D^{\bbeta} \ddot\varrho_n \|
            + \kappa^{\sfrac12} \mr\lambda \| \mr D^{\bbeta} \mr \nabla \ddot\varrho_n \|_2  
      \lesssim&\, \lambda^{ - \big( 2[\bbeta] + \sum_{i=0}^{n-1} s_i \big) b\gamma }
            \left( | \rho_{\ini} |_{\mathfrak{D}(u,\kappa)} + | \rho_{\ini} |_{\mathfrak{P}(q+1)} \right).
                  \label{e:2:unconEnergyP}
\end{align}
From interpolation, we also have
\begin{align}     
      \resN{ \rho_{\ini} }_{u,\kappa} \lesssim&\, \mr \lambda^{1-2\gamma} \left( | \rho_{\ini} |_{\mathfrak{D}(u,\kappa)} + | \rho_{\ini} |_{\mathfrak{P}(q+1)} \right),    \label{e:4:unconEnergyP} \\ 
      \resN{ \rho_{\ini} }_{u,\kappa}^{(H,\alpha_0)} 
            \lesssim&\, \mr \lambda^{(1-2\gamma)(1+\alpha_0)} \left( | \rho_{\ini} |_{\mathfrak{D}(u,\kappa)} + | \rho_{\ini} |_{\mathfrak{P}(q+1)} \right).    \label{e:6:unconEnergyP}
\end{align}
Moreover, for $\bbeta \in \mathcal{I}_*$ with $[\bbeta] \leq 11Q^3 - 3Qn$,
\begin{align}
      \| \mr D^{\bbeta} \ddot\varrho_n \|_\infty
            \lesssim \mr\lambda^{1-2\gamma} \lambda^{ - \big( 2[\bbeta] + \sum_{i=0}^{n-1} s_i \big) b\gamma }
            \left( | \rho_{\ini} |_{\mathfrak{D}(u,\kappa)} + | \rho_{\ini} |_{\mathfrak{P}(q+1)} \right).    
                  \label{e:8:unconEnergyP}
\end{align}
\end{proposition}

\begin{remark}    \label{r:unconEnergyP}
Compare \eqref{e:2:unconEnergyP} and the estimate in \eqref{e:6:homRslv} following Definition \ref{d:resolve}. We know that, for $[\bbeta] > Q$ and $n > Q$, the estimates needed in \eqref{e:6:homRslv} are covered by Proposition \ref{p:unconEnergyP}. Here, we use the fact $Q\gamma > 3$ from Section \ref{ss:parameters}.
\end{remark}

Next, we prove the stability between the chain estimates for $D_{q+1,t} - \kappa_{q+1} \Delta$ and $O_{q+1,t} - \kappa_{q+1} \Delta$.

\begin{proposition}     \label{p:chainStab}
Suppose that, for any $ \rho_{\ini} \in C^\infty (\T^2) $, we have
\begin{align}
      \resN{ \rho_{\ini}}_{v,\kappa} 
            \leq&\, \frac{1}{2} \lambda^{-\alpha} \bigg( | \rho_{\ini} |_{\mathfrak{D}(v,\kappa)} 
                  + \sum_{j=q_*}^{q+1} \lambda_j^{-\sfrac{\gamma}{2}} | \rho_{\ini} |_{\fS(j)} \bigg),   
            \label{e:3:chainStab}     \\ 
      \resN{ \rho_{\ini}}^{(H,\alpha_0)}_{v,\kappa} 
            \leq&\, \left( 1 - \frac{4}{q} + \lambda^{-\sfrac{\gamma}{2}} \right) \bigg( | \rho_{\ini} |_{\mathfrak{D}(v,\kappa)} 
                  + \sum_{j=q_*}^{q+1} \lambda_j^{-\sfrac{\gamma}{2}} | \rho_{\ini} |_{\fS(j)} \bigg).   \label{e:4:chainStab}
\end{align} 
Then for any $\rho_{\ini} \in C^\infty (\T^2)$, we have the stability estimate for dissipation 
\begin{align}     \label{e:7:chainStab}
      \big| | \rho_{\ini} |_{\mathfrak{D}(v, \kappa)} - | \rho_{\ini} |_{\mathfrak{D}(u, \kappa)} \big|      
            \leq&\, \lambda^{ - \alpha - 2 } 
            \bigg( | \rho_{\ini} |_{\mathfrak{D}(v,\kappa)} 
                  + \sum_{j=q_*}^{q+1} \lambda_j^{-\sfrac{\gamma}{2}} | \rho_{\ini} |_{\fS(j)} \bigg)
\end{align}
and the stability estimates for chains
\begin{align}
      \resN{ \rho_{\ini}}_{u,\kappa} 
            \leq&\, \lambda^{-\alpha} \bigg( | \rho_{\ini} |_{\mathfrak{D}(u,\kappa)} 
                  + \sum_{j=q_*}^{q+1} \lambda_j^{-\sfrac{\gamma}{2}} | \rho_{\ini} |_{\fS(j)} \bigg),     \label{e:5:chainStab} \\ 
      \resN{ \rho_{\ini}}^{(H,\alpha_0)}_{u,\kappa} 
            \leq&\, \left( 1 - \frac{4}{q+1} \right) \bigg( | \rho_{\ini} |_{\mathfrak{D}(u,\kappa)} 
                  + \sum_{j=q_*}^{q+1} \lambda_j^{-\sfrac{\gamma}{2}} | \rho_{\ini} |_{\fS(j)} \bigg).     \label{e:6:chainStab}
\end{align}
\end{proposition}

\begin{proof}
Take any $\{ g_n \}_n \subset \mathcal{O}_{q+1}(2NQ)$, $\{ h_n \}_n \subset \mathcal{P}_{q+1}(2NQ)$ and $\{ \bomega_n \}_{n} \subset \mathcal{I}_*$ with $[\bomega_n] \leq 3Q$ satisfying
\begin{align}     \label{e:8:chainStab}
      \vertiii{ g_n }_{q+1} + \vertiii{ h_n }_{q+1} \lesssim&\, \kappa \mr \lambda^2 \lambda^{-s_nb\gamma}.
\end{align}
Suppose that, for any $n \geq 0$, $g_n$ and $h_n$ are congruent. Define $\{\varrho_n\}_n$ and $\{\ddot \varrho_n\}_n$ in Assumption \ref{a:OEq} and Assumption \ref{a:PEq}. To simplify the notation, we introduce
\begin{align*}
      \Box := \bigg( | \rho_{\ini} |_{\mathfrak{D}(v,\kappa)} 
                  + \sum_{j=q_*}^{q+1} \lambda_j^{-\sfrac{\gamma}{2}} | \rho_{\ini} |_{\fS(j)} \bigg).
\end{align*}

\begin{step}[Estimates on derivatives of high order]
We first prove
\begin{align}
      \| \mr D^{\bbeta} \ddot{\varrho}_n \|_\infty 
            \leq&\, \lambda^{ - \alpha - \big( [\bbeta] + \sum_{i=0}^{n-1} (s_i-1) \big) b \gamma } 
                  \bigg( | \rho_{\ini} |_{\mathfrak{D}(u,\kappa)} 
                  + \sum_{j=q_*}^{q+1} \lambda_j^{-\sfrac{\gamma}{2}} | \rho_{\ini} |_{\fS(j)} \bigg),     \label{e:9:chainStab} \\ 
      \| \ddot{\varrho}_n \|_{L^\infty_t C^{\alpha_0}_x} 
            \leq&\, \frac{1}{2} \lambda^{ - \big( [\bbeta] + \sum_{i=0}^{n-1} (s_i-1) \big) b \gamma } 
                  \bigg( | \rho_{\ini} |_{\mathfrak{D}(u,\kappa)} 
                  + \sum_{j=q_*}^{q+1} \lambda_j^{-\sfrac{\gamma}{2}} | \rho_{\ini} |_{\fS(j)} \bigg)      \label{e:10:chainStab}
\end{align}
for $n \geq Q$ or $[\bbeta] \geq Q$. Both follow from \eqref{e:8:unconEnergyP} in Proposition \ref{p:unconEnergyP} and the sufficient gain in exponent $Q\gamma > 3$.
\end{step}

\begin{step}[Estimates on derivatives of relatively low order]

To prove \eqref{e:9:chainStab} and \eqref{e:10:chainStab} at low orders, we first make the following claim.

\begin{claim}     \label{c:2:chainStab}
For any $0 \leq n \leq Q$ and any $\bbeta \in \mathcal{I}_*$ with $[\bbeta] \leq 11Q^3 - 3Qn$, we have
\begin{equation}        \label{c:e:chainStab}
\begin{split}
      \| \mr D^{\bbeta} ( \varrho_n - \ddot{\varrho}_n ) \|
            +&\, \kappa^{\sfrac12} \mr\lambda \| \mr D^{\bbeta} \mr\nabla ( \varrho_n - \ddot{\varrho}_n ) \|_2   \\ 
            \lesssim&\, \lambda^{ - \alpha - 2 - \big( [\bbeta] + \sum_{i=0}^{n-1} (s_i-1) \big) b \gamma } \Box .
\end{split}
\end{equation}
\end{claim}

Using Claim \ref{c:2:chainStab} and standard embedding estimates, we have for any $0 \leq n \leq Q$ and any $\bbeta \in \mathcal{I}_*$ with $[\bbeta] \leq 11Q^3 - 3Qn - 2$,
\begin{equation}  \label{e:11:chainStab}
\begin{split}
      \| \mr D^{\bbeta} ( \varrho_n - \ddot{\varrho}_n ) \|_{\infty}
            \lesssim&\, \lambda^{-\gamma} \lambda^{ - \alpha - \big( [\bbeta] + \sum_{i=0}^{n-1} (s_i-1) \big) b \gamma } \Box,      \\ 
      \| \varrho_n - \ddot{\varrho}_n \|_{L^\infty_t C^{\alpha_0}_x} 
            \lesssim&\, \lambda^{-\gamma} \lambda^{ - \big( [\bbeta] + \sum_{i=0}^{n-1} (s_i-1) \big) b \gamma } \Box.
\end{split}
\end{equation}
Now we use \eqref{e:11:chainStab} and \eqref{e:3:chainStab}-\eqref{e:4:chainStab} to conclude \eqref{e:9:chainStab}-\eqref{e:10:chainStab} for $0 \leq n \leq Q$ and $[\bbeta] \leq Q$, which yields \eqref{e:5:chainStab} and \eqref{e:6:chainStab}. Here, we use the factor $\lambda^{-\gamma}$ in \eqref{e:11:chainStab} to absorb the constants in $\lesssim$.

\end{step}

\begin{proof}[Proof of \eqref{e:7:chainStab} and Claim \ref{c:2:chainStab}]
From \eqref{e:3:chainStab}, we have for any $0 \leq n \leq 2Q^2$, $\bbeta \in \mathcal{I}_*$ with $[\bbeta] \leq 11Q^3 - 3Q n$
\begin{align}     \label{e:12:chainStab}
      \| \mr D^{\bbeta} \varrho_n \| \lesssim
            \lambda^{ - \alpha - \big( [\bbeta] + \sum_{i=0}^{n-1} (s_i-1) \big) b \gamma } \Box .
\end{align}
Subsequently by Lemma \ref{l:streamEst_rig_nobar}, Lemma \ref{l:strDiffEst2} and Lemma \ref{l:strDiffEst3}, we have that, for any $0 \leq n \leq 2Q^2$ and $\bbeta \in \mathcal{I}_*$ with $[\bbeta] \leq 11Q^3 - 3Qn$,
\begin{align}
      \| \mr D^{\bbeta} ( (u-v) \varrho_n ) \| \lesssim&\, 
            \lambda^{ - \alpha - 3 - \big( [\bbeta] + \sum_{i=0}^{n-1} (s_i-1) \big) b \gamma } \Box ,      \label{e:13:chainStab} \\
      \| \mr D^{\bbeta} ( (u-v) \varrho_n ) \| \lesssim&\, 
            \lambda^{ - \alpha - 3 - \big( [\bbeta] + \sum_{i=0}^{n-1} (s_i-1) \big) b \gamma } \Box .      \label{e:14:chainStab}
\end{align}

Observe that $(\varrho_0 - \ddot{\varrho}_0)(\cdot,0) = 0$ and
\begin{align}     \label{e:16:chainStab}
      D_t ( \varrho_0 - \ddot{\varrho}_0 ) - \kappa \Delta ( \varrho_0 - \ddot{\varrho}_0 ) 
            = ( D_t - O_t ) \varrho_0 = \divr( (u-v) \varrho_0 ).
\end{align}
Applying Lemma \ref{l:L2_F_energy}, we deduce, for $0 \leq [\bbeta] \leq 11Q^3$,
\begin{equation}     \label{e:18:chainStab}
\begin{split}
      \| \mr D^{\bbeta} ( \varrho_0 - \ddot{\varrho}_0 ) \|
            +&\, \kappa^{\sfrac12} \mr\lambda \| \mr D^{\bbeta} \mr\nabla ( \varrho_0 - \ddot{\varrho}_0 ) \|_2      \\ 
            \lesssim&\, \kappa^{-\sfrac12} \lambda^{ - \alpha - 3 - [\bbeta] b \gamma } \Box .
\end{split}
\end{equation}
Then from the definition of dissipation in \eqref{e:dissipDef} and using negative power of $\lambda$ to absorb the constants, we have
\begin{equation}     \label{e:19:chainStab}
\begin{split}
      \big| | \rho_{\ini} |_{\mathfrak{D}(v, \kappa)} - | \rho_{\ini} |_{\mathfrak{D}(u, \kappa)} \big|      
            \leq \lambda^{ - \alpha - 2 } \Box .
\end{split}
\end{equation}

Now we proceed with the induction. Assume following estimates for $\bbeta \in \mathcal{I}_*$ with $[\bbeta] \leq 11Q^3 - 3Qn$,
\begin{equation}        \label{e:20:chainStab}
\begin{split}
      \| \mr D^{\bbeta} ( \varrho_n - \ddot{\varrho}_n ) \|
            +&\, \kappa^{\sfrac12} \mr\lambda \| \mr D^{\bbeta} \mr\nabla ( \varrho_n - \ddot{\varrho}_n ) \|_2   \\ 
            \lesssim&\, \lambda^{ - \alpha - 2 - \big( [\bbeta] + \sum_{i=0}^{n-1} (s_i-1) \big) b \gamma } \Box .
\end{split}
\end{equation}
It is not difficult to see $\varrho_{n+1} - \ddot{\varrho}_{n+1}$ satisfies
\begin{equation}  \label{e:22:chainStab}
\begin{split}
      D_t ( \varrho_{n+1} - \ddot{\varrho}_{n+1} ) -&\, \kappa \Delta ( \varrho_{n+1} - \ddot{\varrho}_{n+1} )    \\ 
            =&\, \mr \divr \big( \mr\lambda (u-v) \varrho_{n+1} \big)
            + \mr \divr \big( (g_n-h_n) \mr O^{\bomega_n} \varrho_n \big)     \\ 
            +&\, \mr \divr \big( h_n ( \mr O^{\bomega_n} - \mr D^{\bomega_n} ) \varrho_n \big)
            + \mr \divr \big( h_n \mr D^{\bomega_n} ( \varrho_n - \ddot{\varrho}_n ) \big).
\end{split}
\end{equation}
Now we have for $\bbeta \in \mathcal{I}_*$ with $[\bbeta] \leq 11Q^3 - 3Q (n+1)$,
\begin{align}
      \kappa^{-\sfrac12} \| \mr D^{\bbeta} \big( (u-v) \varrho_{n+1} \big) \|_2
            \lesssim&\, \lambda^{ - \alpha - 2 - \big( [\bbeta] + \sum_{i=0}^{n} (s_i-1) \big) b \gamma } \Box , 
            \label{e:24:chainStab} \\ 
      \kappa^{-\sfrac12} \mr\lambda^{-1} \| \mr D^{\bbeta} \big( (g_n-h_n) \mr O^{\bomega_n} \varrho_n \big) \|_2
            \lesssim&\, \lambda^{ - \alpha - 2 - \big( [\bbeta] + \sum_{i=0}^{n} (s_i-1) \big) b \gamma } \Box , 
            \label{e:26:chainStab} \\ 
      \kappa^{-\sfrac12} \mr\lambda^{-1} \| \mr D^{\bbeta} \big( h_n ( \mr O^{\bomega_n} - \mr D^{\bomega_n} ) \varrho_n \big) \|_2
            \lesssim&\, \lambda^{ - \alpha - 2 - \big( [\bbeta] + \sum_{i=0}^{n} (s_i-1) \big) b \gamma } \Box , 
            \label{e:28:chainStab} \\ 
      \kappa^{-\sfrac12} \mr\lambda^{-1} \| \mr D^{\bbeta} \big( h_n \mr D^{\bomega_n} ( \varrho_n - \ddot{\varrho}_n ) \big) \|_2
            \lesssim&\, \lambda^{ - \alpha - 2 - \big( [\bbeta] + \sum_{i=0}^{n} (s_i-1) \big) b \gamma } \Box .
            \label{e:30:chainStab}
\end{align}
Here, \eqref{e:24:chainStab} follows from \eqref{e:14:chainStab}. \eqref{e:26:chainStab} follows from \eqref{e:12:chainStab} and Corollary \ref{c:CongruEst}. \eqref{e:28:chainStab} follows from \eqref{e:12:chainStab} and Lemma \ref{l:strDiffEst3}. \eqref{e:30:chainStab} follows from \eqref{e:20:chainStab} and Corollary \ref{c:CalPEst}.

Applying energy estimates in Lemma \ref{l:L2_F_energy} with (\ref{e:22:chainStab}-\ref{e:30:chainStab}), we can deduce that, for $\bbeta \in \mathcal{I}_*$ with $[\bbeta] \leq 11Q^3 - 3Q (n+1)$,
\begin{align*}
      \| \mr D^{\bbeta} ( \varrho_{n+1} - \ddot{\varrho}_{n+1} ) \|
            +&\, \kappa^{\sfrac12} \mr\lambda \| \mr D^{\bbeta} \mr\nabla ( \varrho_{n+1} - \ddot{\varrho}_{n+1} ) \|_2       \\ 
            \lesssim&\, \lambda^{ - \alpha - 2 - \big( [\bbeta] + \sum_{i=0}^{n} (s_i-1) \big) b \gamma } \Box .
\end{align*}
This finishes the proof of Claim \ref{c:2:chainStab}.
\end{proof}

\end{proof}

\subsection{Statement of expansion lemmas}      \label{ss:stateExpInertial}

In this section, we state two core expansion lemmas for homogenization. These lemmas aim to study a single equation in (\ref{e:5:OEq_Idct}-\ref{e:8:OEq_Idct}).

\begin{lemma}[The expansion lemma I]        \label{l:EXP1}
For $N,Q$ chosen in Section \ref{ss:parameters}, consider the equation \eqref{e:hom:meq}-\eqref{e:hom:mat}, or equivalently \eqref{e:5:OEq_Idct}-\eqref{e:6:OEq_Idct}, under Definition \ref{d:homParaI} and Definition \ref{d:homFlowmap}. Then the solution $\varrho$ admits the expansion
\begin{equation}  \label{e:0:EXP1}
\begin{split}
      \varrho(x,t) =& \sum_{n = 0}^{Q} \frac{1}{\lambda^n} 
            \Big( \hat \rho_{n} (x,t, \lambda \Phi(x,t), \mu t ) 
            + \check \rho_{n} (x, t, \mu t) + \bar \rho_{n} (x, t) \Big) 
            + \tilde \rho(x,t), 
\end{split}
\end{equation}
with $\hat \rho_n, \check \rho_n, \bar \rho_n, \tilde \rho$ given by the following.

\begin{enumerate}[leftmargin=*,label=\textsc{(C.\arabic*)},align=left]
\item \label{c:l:exp1:spt} \textbf{\textup{(Spatial correctors)}} $\hat \rho_n : \T^2 \times [0,1] \times \T^2 \times \T \rightarrow \R$ has the form
\begin{align}
      \hat \rho_0 =&\, 0,       \label{e:EXP1:hrho0} \\ 
      \hat \rho_{n}(x,t,\xi,\tau) =&\, \sum_{a=0}^{n-1} \sum_{ [\balpha]=1,\balpha \in \mathcal{I}_* }^{2(n-a)-1} \sum_{l=1}^{\hat I_{n-a}^{(\balpha)}}
            \hat h^{(\balpha, l)}_{n-a} (x,t) \hat \eta^{(\balpha, l)}_{n-a} (\tau) \hat \chi^{(\balpha, l)}_{n-a} (\xi) \mr {\bar D}^{\balpha} \bar \rho_a(x,t)       \label{e:EXP1:hrhon} \\ 
      \hat h^{(\balpha, l)}_{n} \in&\, \mathcal{\bar P} ( N(n-1)+2-[\balpha] )       \label{e:EXP1:hh}
\end{align}
with $\hat h^{(\cdot,\cdot)}_{\cdot} : \T^2 \times [0,1] \rightarrow \R$, $\hat \chi_{\cdot}^{(\cdot,\cdot)} : \T^2 \rightarrow \R$ and $\hat \eta_{\cdot}^{(\cdot,\cdot)} : \T \rightarrow \R $ satisfying the estimates for $k \geq 0$
\begin{align}
      \hat I^{(\balpha)}_{k+1} &\, \leq 2(3N)^{k},   \label{e:EXP1:hI_est} \\ 
      \vertiii{ \hat h^{(\balpha, l)}_{k+1} } &\, \lesssim \frac{\delta^{\sfrac12}\mr{\bar\lambda}}{\kappa\lambda} 
            \bigg( \frac{ \delta^{\sfrac12}\mr{\bar\lambda} } {\kappa\lambda} \cdot \lambda^{\gamma_S} \bigg)^{k}         \label{e:EXP1:hh_Pest}   \\ 
      \sum_{[\balpha]=1}^{2k+1} \sum_{l=1}^{\hat I^{(\balpha)}_{k+1}} 
            \vertiii{ \hat h^{(\balpha, l)}_{k+1} } &\, \lesssim \frac{\delta^{\sfrac12}\mr{\bar\lambda}}{\kappa\lambda} 
            \bigg( \frac{ \delta^{\sfrac12}\mr{\bar\lambda} } {\kappa\lambda} \cdot \lambda^{\gamma_S} \bigg)^{k}         \label{e:EXP1:hhS_Pest}   \\ 
      \| \nabla_\xi^p \hat\chi^{(\balpha,l)}_{k+1} \|_{\infty} &\, \lesssim 1,
            \quad p \leq 11Q^3, 
            \quad \langle \hat\chi^{(\balpha,l)}_{k+1} \rangle_\xi = 0,     \label{e:EXP1:hchi_est} \\
      \| \partial_\tau^p \hat\eta^{(\balpha,l)}_{k+1} \|_{\infty} &\, \lesssim 1,
            \quad p \leq 11Q^3-k-1.       \label{e:EXP1:heta_est}
\end{align}
Here the sum of differentiation index $\balpha$ in \eqref{e:EXP1:hrhon} is among all indices in $\mathcal{I}_*$ such that $1 \leq [\balpha] \leq 2(n-a)-1$. Moreover, $(\hat \chi^{(\balpha, l)}_{n}, \hat \eta^{(\balpha, l)}_{n})$ satisfy the shear condition in Definition \ref{d:shear_structure} for any $n$, $\balpha$ and $l$.

\item \label{c:l:exp1:tpr} \textbf{\textup{(Temporal correctors)}} $\check \rho_n : \T^2 \times [0,1] \times \T \rightarrow \R$ has the form
\begin{align}
      \check \rho_0 =&\, 0,     \label{e:EXP1:crho0} \\ 
      \check \rho_{n}(x,t,\tau) =&\, \sum_{a=0}^{n-1} \sum_{ [\balpha]=1, \balpha \in \mathcal{I}_* }^{2(n-a)} \sum_{l=1}^{\check I_{n-a}^{(\balpha)}}
            \check h^{(\balpha, l)}_{n-a} (x,t) \check \eta^{(\balpha, l)}_{n-a} (\tau) \mr {\bar D}^{\balpha} \bar \rho_a(x,t)    \label{e:EXP1:crhon} \\ 
      \check h^{(\balpha, l)}_{n} \in&\, \mathcal{\bar P} ( N(n-1)+3-[\balpha] \big)    \label{e:EXP1:ch}
\end{align}
with $\check h^{(\cdot,\cdot)}_{\cdot} : \T^2 \times [0,1] \rightarrow \R$, $\check\eta^{(\cdot,\cdot)}_{\cdot} : \T \rightarrow \R $ satisfy the estimates for $k \geq 0$
\begin{align}
      \check I^{(\balpha)}_{k+1} &\, \leq 8(3N)^{k},   \label{e:EXP1:cI_est} \\ 
            \vertiii{ \check h^{(\balpha,l)}_{k+1} } 
            &\, \lesssim \frac{\delta^{\sfrac12}\mr{\bar\lambda}}{\mu} \frac{ \delta^{\sfrac12}\mr{\bar\lambda} } {\kappa\lambda} 
            \bigg( \frac{ \delta^{\sfrac12}\mr{\bar\lambda} } {\kappa\lambda} \cdot \lambda^{\gamma_S} \bigg)^{k},    \label{e:EXP1:ch_Pest} \\ 
      \sum_{[\balpha]=1}^{2k+2} \sum_{l=1}^{\check I^{(\balpha)}_{k+1}}
            \vertiii{ \check h^{(\balpha,l)}_{k+1} } 
            &\, \lesssim \frac{\delta^{\sfrac12}\mr{\bar\lambda}}{\mu} \frac{ \delta^{\sfrac12}\mr{\bar\lambda} } {\kappa\lambda} 
            \bigg( \frac{ \delta^{\sfrac12}\mr{\bar\lambda} } {\kappa\lambda} \cdot \lambda^{\gamma_S} \bigg)^{k},    \label{e:EXP1:chS_Pest} \\
      \| \partial_\tau^p \check\eta^{(\balpha,l)}_{k+1} \|_{\infty} &\, \lesssim 1,   
            \quad p \leq 11Q^3-k,
            \quad \langle \check\eta^{(\balpha,l)}_{k+1} \rangle_\tau = 0.       \label{e:EXP1:ceta_est}
\end{align}

\item \label{c:l:exp1:rsd} \textbf{\textup{(Residual correctors)}} $\bar \rho_0, \{\bar \rho_n\}_{1 \leq n \leq N}: \T^2 \times [0,1] \rightarrow \R $ solve
\begin{align}
      \bar L_0 \bar \rho_0 =&\, 0,   \label{e:EXP1:brho0} \\ 
      \bar L_0 \bar \rho_n =&\, \sum_{a=0}^{n-1} \sum_{ [\balpha]=1, \balpha \in \mathcal{I}_* }^{2(n-a)+1} 
      \mr{\bar \divr} \Big( \bar h^{(\balpha)}_{n-a} \mr {\bar D}^{\balpha} \bar \rho_a \Big)   \label{e:EXP1:brhon} \\ 
      \bar h^{(\balpha)}_{n} \in&\, \mathcal{\bar P} ( Nn+3-[\balpha] )       \label{e:EXP1:bh} 
\end{align}
with $\bar L_0$ given by
\begin{align}   \label{e:EXP1:homedOp}  
      \bar L_0 = \bar D_t - \bar \kappa \Delta
\end{align}
and initial datum given by
\begin{align}
      \bar \rho_0(\cdot, 0) =&\, \bS \rho_{\ini},     \label{e:EXP1:brho0_ini} \\ 
      \bar \rho_n(\cdot, 0) =&\, - \bS [ \check \rho_n(\cdot, 0, 0) ].      \label{e:EXP1:brhon_ini}
\end{align}
Here, \underline{the vector-valued functions} $\bar h^{(\cdot)}_{\cdot} : \T^2 \times [0,1] \rightarrow \R^2$ satisfy the estimates for $k \geq 0$
\begin{align}
      \vertiii{ \bar h^{(\balpha)}_{k+1} } \lesssim&\, 
      \kappa \bigg( \frac{ \delta^{\sfrac12}\mr{\bar\lambda} } {\kappa\lambda} \bigg)^2
      \bigg( \frac{ \delta^{\sfrac12}\mr{\bar\lambda} } {\kappa\lambda} \cdot \lambda^{\gamma_S} \bigg)^{k+1}.       \label{e:EXP1:bh_Pest} \\ 
      \sum_{[\balpha]=1}^{2k+3}
      \vertiii{ \bar h^{(\balpha)}_{k+1} } \lesssim&\, 
      \kappa \bigg( \frac{ \delta^{\sfrac12}\mr{\bar\lambda} } {\kappa\lambda} \bigg)^2
      \bigg( \frac{ \delta^{\sfrac12}\mr{\bar\lambda} } {\kappa\lambda} \cdot \lambda^{\gamma_S} \bigg)^{k+1}.       \label{e:EXP1:bhS_Pest}
\end{align}

\item \label{c:l:exp1:rmd} \textbf{\textup{(remainder)}} $\tilde \rho$ solves the following equation
\begin{align}
      \bar D_t \tilde\rho - \divr( A \nabla \tilde\rho ) 
      =&\, \tilde f    \label{e:exp1:tildeRho}
\end{align}
with initial datum
\begin{align}
      \tilde \rho(x,0) =&\, \rho_{\ini} (x) - \mathbb{S} \rho_{\ini} (x)
            + \sum_{a=1}^Q \frac{1}{\lambda^a} \big( \bS[ \check\rho_a (\cdot,0,0) ] - \check\rho_a (\cdot,0,0) \big).      \label{e:0:tildeRho}
\end{align}
Here, $\tilde f: \T^2 \times [0,1] \rightarrow \R$ is given by
\begin{align}
      \tilde f := &\, -\frac{1}{\lambda^{Q}} 
            \Big( L_1 ( \hat\rho_{Q} + \check\rho_{Q} + \bar\rho_{Q} ) 
                  + L_0 ( \hat\rho_{Q-1} + \check\rho_{Q-1} + \bar\rho_{Q-1} )      \label{e:2:tildeRho} \\
                  &\, \quad + L_{-1} \hat \rho_{Q-2} \Big) (x,t, \lambda \Phi(x,t), \mu t )     \label{e:4:tildeRho} \\ 
            &\,- \frac{1}{\lambda^{Q+1}} \Big( ( L_0 ( \hat\rho_{Q} + \check\rho_{Q} + \bar\rho_{Q} ) + L_{-1} \hat \rho_{Q-1} \Big) (x,t, \lambda \Phi(x,t), \mu t )         \label{e:6:tildeRho} \\ 
            &\,- \frac{1}{\lambda^{Q+2}} \Big( L_{-1} \hat \rho_{Q} \Big) (x,t, \lambda \Phi(x,t), \mu t )        \label{e:8:tildeRho}
\end{align}
and differential operators $L_1,L_0,L_{-1}$ acting on general functions $\rho$ taking arguments $(x,t,\xi,\tau)$ are given by
\begin{align}
      L_1 \rho =&\, - \frac{\delta^{\sfrac12}}{\lambda} \partial_{\xi_i} H_{ij} \partial_{x_j} \rho + \frac{\mu}{\lambda} \partial_\tau \rho
            + \varepsilon \kappa \lambda \sum_m S_{m,ij} \vartheta_m \partial_{\xi_i\xi_j} \rho       \label{e:10:tildeRho} \\ 
      L_0 \rho =&\, \bar D_t \rho - \kappa \Delta_x \rho - \varepsilon \frac{\delta^{\sfrac12}}{\lambda} \sum_m \sigma_m \partial_{x_i} \Omega_{m,ij} H_{12} \partial_{x_j} \rho
            - \varepsilon \delta^{\sfrac12} \sum_m \varphi_m \partial_{\xi_i} ( H E_m^T )_{ij} \partial_{x_j} \rho      \label{e:12:tildeRho} \\ 
            +&\, \varepsilon \delta^{\sfrac12} \sum_m \varphi_m \partial_{x_j} ( HE_m^T )_{ij} \partial_{\xi_i} \rho 
            - \varepsilon \kappa \lambda \sum_m \phi_m \partial_{x_j} B_{m,ij} \partial_{\xi_i} \rho
            - 2 \kappa \lambda \partial_{\xi_i x_i} \rho    \label{e:14:tildeRho} \\ 
            -&\, 2 \varepsilon \kappa \lambda \sum_m \phi_m B_{m,ij} \partial_{\xi_i x_j} \rho        \label{e:16:tildeRho} \\ 
      L_{-1} \rho =&\, \frac{1}{\lambda} \phi_* \omega_* \cdot \nabla_\xi \rho.     \label{e:18:tildeRho}
\end{align}
\end{enumerate}
\end{lemma}

\begin{remark}
For spatial and temporal correctors, we may have none, one, or multiple correctors for fixed $n,a$ and differentiation index $\balpha$. For spatial correctors in \eqref{e:EXP1:hrhon}, the cardinality of such correctors is denoted by $\hat I_{n-a}^{(\balpha)} \in \N$. When $\hat I_{n-a}^{(\balpha)} = 0$, this means we do not have correctors for such $n,a$ and $\balpha$. The same also applies to temporal correctors.
\end{remark}

\begin{remark}[Stationarity of the correctors]  \label{r:stationarity}
Note that the structural information of these correctors, such as $\hat I_{n-a}^{(\balpha)}$, $\hat h_{n-a}^{(\balpha,l)}$, $\hat \eta_{n-a}^{(\balpha,l)}$ and $\hat \chi_{n-a}^{(\balpha,l)}$ depend on the integer $n-a$, instead of $n$ or $a$. For example, this means the structural information in first order corrector $\hat \rho_1$ is shared in higher order correctors $\hat \rho_n$ with $n \geq 2$. We refer this property as the \textit{stationarity} of the correctors. We state the expansion lemma in a way that stationarity is claimed implicitly. However, this property does not come for free, and we shall prove it later.
\end{remark}

We also have the following lemmas to estimate the initial datum of macroscopic states $\bar\rho_n$ in Lemma \ref{l:EXP1}.

\begin{corollary}     \label{c:t=0Ref:EXP1}
Given the assumptions in Lemma \ref{l:EXP1}, then we have, for any $1 \leq n \leq Q$,
\begin{align}
      \frac{1}{\lambda^n} \check \rho_{n} (x,0,0) =&\, \sum_{a=0}^n \sum_{ [\balpha] = 1, \balpha \in \mathcal{I}_x }^{4n+1}
                  \mathfrak{C}_{n,a}^{(\balpha)} \bigg( \frac{1}{\lambda^a} \mr {\bar \partial}^{\balpha} \bar \rho_{a}(x,0) \bigg) 
                  \label{e:2:t=0Ref:EXP1} \\
      \frac{1}{\lambda^n} \bar \rho_{n} (x,0) =&\, \sum_{a=0}^n \sum_{ [\balpha] = 1, \balpha \in \mathcal{I}_x }^{4n+1}
                  \mathfrak{C}_{n,a}^{(\balpha)} \, \bS \bigg( \frac{1}{\lambda^a} \mr {\bar \partial}^{\balpha} \bar \rho_{a}(x,0) \bigg).     
                  \label{e:3:t=0Ref:EXP1}
\end{align}
For $p \geq 1$, we have
\begin{align}
      \frac{1}{\lambda^n} \mr {\bar D}^p_t \bar \rho_{n} (x,0) =&\, \sum_{a=0}^n \sum_{ [\balpha] = 1, \balpha \in \mathcal{I}_x }^{4n+2p+1}
                  \mathfrak{C}_{n,a}^{(p,\balpha)} \bigg( \frac{1}{\lambda^a} \mr {\bar \partial}^{\balpha} \bar \rho_{a}(x,0) \bigg),      
                  \label{e:4:t=0Ref:EXP1}
\end{align}
with the following estimates for the constant coefficients
\begin{align}
      | \mathfrak{C}_{n,a}^{(p,\balpha)} | \lesssim \lambda^{-(n-a)b\gamma}.  \label{e:6:t=0Ref:EXP1}
\end{align}
\end{corollary}

\begin{proof}
This corollary is a direct consequence of \eqref{e:EXP1:brhon}. When $t \in [0, 8\bar\mu^{-1}]$, the velocity field in $\bar D_t$ is zero from \eqref{e:u_zero_initialT}, hence $\nabla$ and $\bar D_t$ commute at time $t=0$. For $n=0$, we can deduce that, from \eqref{e:EXP1:brho0},
\begin{align*}
      \mr{\bar D}_t^p \bar \rho_0 (x, 0) = 
            \bigg( \frac{\bar\kappa \mr{\bar\lambda}^2}{\mr{\bar\mu}} \bigg)^p \mr{\bar \Delta}^p \bar\rho_0 ( x, 0 ).
\end{align*}
Using \eqref{e:1:mScaleRela}, we can deduce \eqref{e:4:t=0Ref:EXP1}-\eqref{e:6:t=0Ref:EXP1} for $n=0$.

For $n \geq 1$, we prove \eqref{e:4:t=0Ref:EXP1}-\eqref{e:6:t=0Ref:EXP1} by induction. Assume the desired estimates hold for $a \leq n$. We use \eqref{e:EXP1:crhon} to derive the desired equality in \eqref{e:4:t=0Ref:EXP1} for $p=0$. We write
\begin{align}     \label{e:8:t=0Ref:EXP1}
      \frac{1}{\lambda^{n+1}} \check \rho_{n+1} (x,0,0) = \sum_{a=0}^{n} \sum_{[\balpha]=1, \balpha \in \mathcal{I}_* }^{2(n+1-a)} \sum_{l=1}^{\check I_{n+1-a}^{(\balpha)}}
            \bigg( \frac{1}{\lambda^{n+1-a}} \check h^{(\balpha, l)}_{n+1-a} \bigg) (x,0) \check \eta^{(\balpha, l)}_{n+1-a} (0) 
            \bigg( \frac{1}{\lambda^a} \mr {\bar D}^{\balpha} \bar \rho_a \bigg) (x,0).
\end{align}
Using the induction hypothesis for $a \leq n$, \eqref{e:EXP1:ch_Pest} and \eqref{e:EXP1:ceta_est}, this yields a representation of $\lambda^{-(n+1)} \check \rho_{n+1} (x,0,0)$ analogous to \eqref{e:4:t=0Ref:EXP1}, with constants satisfying \eqref{e:6:t=0Ref:EXP1}. Here, we use Remark \ref{r:PolyDeriP3} to deduce that $\check h^{(\balpha, l)}_{n+1-a}$ in \eqref{e:8:t=0Ref:EXP1} is constant in $x$ when $t=0$. Then from \eqref{e:EXP1:brhon_ini} and Remark \ref{r:2:SProj}, we deduce \eqref{e:4:t=0Ref:EXP1} and \eqref{e:6:t=0Ref:EXP1} for $n+1$ and $p=0$.

When $p \geq 1$, we can deduce from \eqref{e:EXP1:brhon}
\begin{align}     \label{e:10:t=0Ref:EXP1}
      \mr{\bar D}_t^p \bar \rho_{n+1} (x,0) - \frac{\bar\kappa \mr{\bar\lambda}^2}{\mr{\bar\mu}} \mr{\bar D}_t^{p-1} \mr{\bar \Delta} \bar\rho_{n+1} (x,0) 
            = \frac{1}{\mr{\bar\mu}} \sum_{a=0}^{n} \sum_{ [\balpha]=1, \balpha \in \mathcal{I}_* }^{2(n-a)+3} 
            \mr{\bar \divr} \Big( \bar h^{(\balpha)}_{n-a+1} (x,0) \mr{\bar D}_t^{p-1} \mr {\bar D}^{\balpha} \bar \rho_a (x,0) \Big).
\end{align}
Here, note that $\bar h^{(\balpha)}_{n+1-a} (x,t)$ is constant when $t$ is suitably small, from Remark \ref{r:PolyDeriP3}. Now we can deduce the representation in \eqref{e:4:t=0Ref:EXP1} for $\mr{\bar D}_t^p \bar \rho_{n+1} (x,0)$ by another induction on $p$. The estimates on constants \eqref{e:6:t=0Ref:EXP1} follow from induction assumptions \eqref{e:homGammaRela} and \eqref{e:EXP1:bh_Pest}.

From the parameters relations, including \eqref{e:1:mScaleRela}, \eqref{e:EXP1:bhS_Pest} and induction assumptions, we can deduce \eqref{e:6:t=0Ref:EXP1} for $n+1$ from the induction hypothesis.

\end{proof}

\begin{lemma}     \label{l:t=0:EXP1}
Given the assumptions in Lemma \ref{l:EXP1}, we have that, for any $j \leq q$,
\begin{align}
      \frac{1}{\lambda^{n}} \| \bar \rho_{n} (\cdot,0) \|_2
            + \frac{1}{\lambda^{n}} \| \check \rho_{n} (\cdot,0,0) \|_2
            &\,\lesssim \lambda^{-nb\gamma} \mr{\lambda}_q^{-2b\gamma} | \rho_{\ini} |_{\mathfrak{S}(q)},     \label{e:EXP1:brhoL2Ini}  \\ 
      \frac{1}{\lambda^{n}} | \bar \rho_{n} (\cdot,0) |_{\mathfrak{S}(j)}
            + \frac{1}{\lambda^{n}} | \check \rho_{n} (\cdot,0,0) |_{\mathfrak{S}'(j)}
            &\,\lesssim \lambda^{-nb\gamma} \Bigg( \frac{\mr\lambda_{j}} {\mr{\lambda}_q} \Bigg) 
                  | \rho_{\ini} |_{\mathfrak{S}(j)}.     \label{e:EXP1:brhoSqStrIni}  \\ 
      \frac{1}{\lambda^{n}} \big| \bS [ \check\rho_n (\cdot,0,0) ] - \check\rho_n (\cdot,0,0) \big|_{\mathfrak{P}(q+1)}
            &\,\lesssim \lambda^{-(Q+n)b\gamma} 
                  | \rho_{\ini} |_{\mathfrak{S}(q+1)}.     \label{e:EXP1:brhoSq2Ini}      
\end{align}
\end{lemma}

\begin{proof}
The $\check \rho_n$ part of \eqref{e:EXP1:brhoSqStrIni} follow from \eqref{e:2:t=0Ref:EXP1}, Remark \ref{r:datumProj}. Then the $\bar \rho_n$ part of \eqref{e:EXP1:brhoSqStrIni} follows from Definition \ref{d:SProj} and \eqref{e:EXP1:brhon_ini}. \eqref{e:EXP1:brhoSq2Ini} follows from \eqref{e:d6IniDatNorm} and Remark \ref{r:datumProj}.

\end{proof}

To study equations like \eqref{e:7:OEq_Idct1}-\eqref{e:8:OEq_Idct}, we also need to consider $O_t - \kappa \Delta$ with forces. Here, we specify a typical force containing small scales $\lambda \Phi$ and $\mu t$.

\begin{assumption}[Small scale force]   \label{a:homSmallscl}
Define a force $F: \T^2 \times [0,1] \rightarrow \R^2$ via
\begin{align*}
      F(x,t) = \sum_{k=1}^K \sum_{ [\balpha] = 1, \balpha \in \mathcal{I}_* }^Q
            \mr {\bar D}^{\balpha} \bar \rho_{F,k}(x,t) \sum_{l = 1}^{J} h_{F,k}^{(\balpha,l)}(x,t) \eta_{F,k}^{(\balpha,l)}(\mu t) \chi_{F,k}^{(\balpha,l)} (\lambda \Phi(x,t))
\end{align*}
The collection of functions $\{\bar \rho_{F,k}\}_k$ satisfies transmission condition at step $q$. For some $ \Upsilon > 0 $, $P \in \N$, the \underline{vector-valued functions} $h_F^{(\cdot,\cdot)}: \T^2 \times [0,1] \rightarrow \R^2 $ satisfying
\begin{align} 
      h_{F,k}^{(\balpha,l)} \in [\mathcal{\bar P} ( Q )]^2,       \label{e:2:homSmallscl} \\ 
      \sum_{k=1}^K \sum_{[\balpha] = 1}^Q \sum_{l=1}^{J} \vertiii{ h_{F,k}^{(\balpha,l)} } \leq \kappa \lambda \mr\lambda \Upsilon
            \label{e:4:homSmallscl}
\end{align}
and other scalar-valued functions satisfying
\begin{align}
      \| \partial^p_\tau \eta_{F,k}^{(\balpha,l)} \|_{\infty} \leq&\, 1
            \quad p \leq P,    \label{e:8:homSmallscl} \\ 
      \| \nabla^p_\xi \chi_{F,k}^{(\balpha,l)} \|_{\infty} \leq&\, 1
            \quad p \leq P.    \label{e:10:homSmallscl}
\end{align}
Moreover, we assume $\big( \chi_{F,k}^{(\balpha,l)}, \eta_{F,k}^{(\balpha,l)} \big)$ satisfies generalized shear condition for any $\balpha$ and $l$.
\end{assumption}

\begin{lemma}[The expansion lemma II]        \label{l:EXP2}
For $N,Q$ chosen in Section \ref{ss:parameters} and any force $F$ given by Assumption \ref{a:homSmallscl}, consider the equation 
\begin{equation} \label{e:hom:eqFast}
\begin{split}
    \bar D_t \varrho - \divr ( A \nabla \varrho ) =&\, \mr\lambda^{-1} \divr F, \\ 
    \varrho( \cdot, 0 ) =& 0 
\end{split}
\end{equation}
with \eqref{e:hom:mat} in $\T^2 \times [0,1]$ under Definition \ref{d:homParaI} and Definition \ref{d:homFlowmap}. Then the solution $\varrho$ admits the following expansion
\begin{equation}  \label{e:0:EXP2}
\begin{split}
      \varrho (x,t) =& \sum_{n = 0}^{Q} \frac{1}{\lambda^n} 
            \Big( \hat \rho_{n} (x,t, \lambda \Phi(x,t), \mu t ) 
            + \check \rho_{n} (x, t, \mu t) + \bar \rho_{n} (x, t) \Big) 
            + \tilde \rho(x,t), 
\end{split}
\end{equation}
with $\hat \rho_n, \check \rho_n, \bar \rho_n, \tilde \rho$ given by the following information.

\begin{enumerate}[leftmargin=*,label=\textsc{(D.\arabic*)},align=left]

\item \label{c:l:exp2:spt} \textbf{\textup{(Spatial correctors)}} $\hat \rho_n : \T^2 \times [0,1] \times \T^2 \times \T \rightarrow \R$ has the form
\begin{align}
      \hat \rho_0 =&\, 0,       \label{e:EXP2:hrho0} \\ 
      \hat \rho_{n}(x,t,\xi,\tau) =&\, \sum_{a=0}^{n-1} \sum_{ [\balpha]=1, \balpha \in \mathcal{I}_* }^{2(n-a)-1} \sum_{l=1}^{\hat I_{n-a}^{(\balpha)}}
            \hat h^{(\balpha, l)}_{n-a} (x,t) \hat \eta^{(\balpha, l)}_{n-a} (\tau) \hat \chi^{(\balpha, l)}_{n-a} (\xi) \mr {\bar D}^{\balpha} \bar \rho_a(x,t)       \label{e:EXP2:hrhon} \\ 
            +&\, \sum_{k=1}^K \sum_{ [\balpha]=1, \balpha \in \mathcal{I}_* }^{2n-2+Q} 
            \sum_{l=1}^{\hat I^{(\balpha)}_{F,n}}
            \hat h^{(\balpha,l)}_{F,k,n} (x,t) \hat \eta^{(\balpha,l)}_{F,k,n} (\tau) \hat \chi^{(\balpha,l)}_{F,k,n} (\xi) \mr {\bar D}^{\balpha} \bar \rho_{F,k} (x,t)     \label{e:EXP2:hrhon:F} \\ 
      \hat h^{(\balpha, l)}_{F,n} \in&\, \mathcal{\bar P} ( Q + N(n-1) + 1 )       \label{e:EXP2:hhF}
\end{align}
with $\hat h^{(\cdot,\cdot)}_{\cdot} : \T^2 \times [0,1] \rightarrow \R$, $\hat \chi_{\cdot}^{(\cdot,\cdot)} : \T^2 \rightarrow \R$ and $\hat \eta_{\cdot}^{(\cdot,\cdot)} : \T \rightarrow \R $ satisfy the estimates for $n \geq 1$
\begin{align}
      \hat I^{(\balpha)}_{F,n} &\, \leq J (3N)^{n-1}       \label{e:EXP2:hIF_est} \\ 
      \vertiii{ \hat h^{(\balpha,l)}_{F,k,n} } &\,\lesssim \lambda^2 \bigg( \frac{ \delta^{\sfrac12}\mr{\bar\lambda} } {\kappa\lambda} \cdot \lambda^{\gamma_S} \bigg)^{n-2} \Upsilon,    \label{e:EXP2:hhF_Pest} \\
      \sum_{k=1}^K \sum_{[\balpha]=1}^{2n-2+Q} \sum_{l=1}^{\hat I^{(\balpha)}_{F,n}}
            \vertiii{ \hat h^{(\balpha,l)}_{F,k,n} } &\,\lesssim \lambda^2 \bigg( \frac{ \delta^{\sfrac12}\mr{\bar\lambda} } {\kappa\lambda} \cdot \lambda^{\gamma_S} \bigg)^{n-2} \Upsilon,    \label{e:EXP2:hhSF_Pest} \\
      \| \nabla_\xi^p \hat\chi^{(\balpha,l)}_{F,k,n} \|_{\infty} &\, \lesssim 1,
            \quad p \leq P,
            \quad \langle \hat\chi^{(\balpha,l)}_{F,k,n} \rangle_\xi = 0,   
             \label{e:EXP2:hchi_est} \\
      \| \partial_\tau^p \hat\eta^{(\balpha,l)}_{F,k,n} \|_{\infty} &\, \lesssim
      1,    \quad p \leq P + 1 - n.       \label{e:EXP2:heta_est}
\end{align}
Here the sum of differentiation index $\balpha$ in \eqref{e:EXP2:hrhon} is among all indices in $\mathcal{I}_*$ such that $1 \leq [\balpha] \leq 2(n-a)-1$. Moreover, these correctors satisfy the shear condition in Definition \ref{d:shear_structure}.

\item \label{c:l:exp2:tpr} \textbf{\textup{(Temporal correctors)}} $\check \rho_n : \T^2 \times [0,1] \times \T \rightarrow \R$ has the form
\begin{align}
      \check \rho_0 =&\, 0,     \label{e:EXP2:crho0} \\ 
      \check \rho_{n}(x,t,\tau) =&\, \sum_{a=0}^{n-1} \sum_{ [\balpha]=1, \balpha \in \mathcal{I}_* }^{2(n-a)} \sum_{l=1}^{\check I_{n-a}^{(\balpha)}}
            \check h^{(\balpha, l)}_{n-a} (x,t) \check \eta^{(\balpha, l)}_{n-a} (\tau) \mr {\bar D}^{\balpha} \bar \rho_a(x,t)    \label{e:EXP2:crhon} \\ 
      +&\, \sum_{k=1}^K \sum_{ [\balpha]=1, \balpha \in \mathcal{I}_* }^{2n-1+Q} 
            \sum_{l=1}^{\check I^{(\balpha)}_{F,n}}  \check h^{(\balpha,l)}_{F,k,n} (x,t) \check \eta^{(\balpha,l)}_{F,k,n} (\tau) \mr {\bar D}^{\balpha} \bar \rho_{F,k}(x,t)    \label{e:EXP2:crhon:F} \\ 
      \check h^{(\balpha, l)}_{F,n} \in&\, \mathcal{\bar P} ( Q + N(n-1) + 2 )    \label{e:EXP2:ch}
\end{align}
with $\check h^{(\cdot,\cdot)}_{\cdot} : \T^2 \times [0,1] \rightarrow \R$, $\check\eta^{(\cdot,\cdot)}_{\cdot} : \T \rightarrow \R $ satisfy the estimates for $n \geq 1$
\begin{align}
      \check I^{(\balpha)}_{F,n} &\, \leq 4 J (3N)^{n-1},     \label{e:EXP2:cIF_est}  \\ 
      \vertiii{ \check h^{(\balpha,l)}_{F,k,n} } &\, \lesssim
            \frac{ \kappa\lambda^2 \mr{\bar\lambda} } { \mu } 
            \bigg( \frac{ \delta^{\sfrac12}\mr{\bar\lambda} } {\kappa\lambda} \cdot \lambda^{\gamma_S} \bigg)^{n-1} \Upsilon,        \label{e:EXP2:chF_Pest}  \\ 
      \sum_{k=1}^K \sum_{[\balpha]=1}^{2n-1+Q} \sum_{l=1}^{\check I^{(\balpha)}_{F,n}}
            \vertiii{ \check h^{(\balpha,l)}_{F,k,n} } &\, \lesssim
            \frac{ \kappa\lambda^2 \mr{\bar\lambda} } { \mu } 
            \bigg( \frac{ \delta^{\sfrac12}\mr{\bar\lambda} } {\kappa\lambda} \cdot \lambda^{\gamma_S} \bigg)^{n-1} \Upsilon,        \label{e:EXP2:chSF_Pest}  \\ 
      \| \partial_\tau^p \check\eta^{(\balpha,l)}_{F,k,n} \|_{\infty} &\, \lesssim 1,
            \quad p \leq P + 2 - n.       \label{e:EXP2:ceta_est}
\end{align}

\item \label{c:l:exp2:rsd} \textbf{\textup{(Residual correctors)}} $\bar \rho_0, \{\bar \rho_n\}_{1 \leq n \leq N}: \T^2 \times [0,1] \rightarrow \R $ solve
\begin{align}
      \bar L_0 \bar \rho_0 =&\, \frac{\mr{\bar\lambda}}{\mr\lambda} \sum_{k=1}^K \sum_{[\balpha] = 1, \balpha \in \mathcal{I}_*}^Q \sum_{l = 1}^{J} 
            \mr{\bar \divr} \Big( h_{F,k}^{(\balpha,l)} \mr {\bar D}^{\balpha} \bar \rho_{F,k} \Big) 
                  \langle \eta_{F,k}^{(\balpha,l)} \rangle_{\tau} \langle \chi_{F,k}^{(\balpha,l)} \rangle_{\xi}    \label{e:EXP2:brho0} \\ 
      \bar L_0 \bar \rho_n =&\, \sum_{a=0}^{n-1} \sum_{ [\balpha]=1, \balpha \in \mathcal{I}_* }^{2(n-a)+1} 
            \mr{\bar \divr} \Big( \bar h^{(\balpha)}_{n-a}\mr {\bar D}^{\balpha} \bar \rho_a \Big)
      + \sum_{k=1}^K \sum_{ [\balpha]=1, \balpha \in \mathcal{I}_* }^{2n+Q} 
            \mr{\bar \divr} \Big( \bar h_{F,k,n}^{(\balpha)} \mr {\bar D}^{\balpha} \bar \rho_{F,k} \Big)    \label{e:EXP2:brhon} \\ 
      \bar h^{(\balpha)}_{F,n} \in&\, \mathcal{\bar P} ( Q + Nn + 2 ),       \label{e:EXP2:bhF} 
\end{align}
with $\bar L_0$ given by \eqref{e:EXP1:homedOp} and initial datum given by
\begin{align}
      \bar \rho_0(\cdot, 0) =&\, 0,      \label{e:EXP2:brho0_ini} \\ 
      \bar \rho_n(\cdot, 0) =&\, - \bS [ \check \rho_n(\cdot, 0, 0) ].      \label{e:EXP2:brhon_ini}
\end{align}
Here, \underline{the vector-valued functions} $\bar h^{(\cdot)}_{\cdot} : \T^2 \times [0,1] \rightarrow \R^2$ satisfy the estimates
\begin{align} 
      \vertiii{ \bar h^{(\balpha)}_{F,k,n} } 
      &\,\lesssim \delta^{\sfrac12} \lambda \mr{\bar\lambda} \bigg( \frac{ \delta^{\sfrac12}\mr{\bar\lambda} } {\kappa\lambda} \cdot \lambda^{\gamma_S} \bigg)^{n-1} \Upsilon,          \label{e:EXP2:bhF_Pest} \\ 
      \sum_{k=1}^K \sum_{[\balpha]=1}^{2n+Q}
      \vertiii{ \bar h^{(\balpha)}_{F,k,n} } 
      &\,\lesssim \delta^{\sfrac12} \lambda \mr{\bar\lambda} \bigg( \frac{ \delta^{\sfrac12}\mr{\bar\lambda} } {\kappa\lambda} \cdot \lambda^{\gamma_S} \bigg)^{n-1} \Upsilon.          \label{e:EXP2:bhSF_Pest}
\end{align}

\item \label{c:l:exp2:rmd} \textbf{\textup{(remainder)}} $\tilde \rho$ solves the equation \eqref{e:exp1:tildeRho}-\eqref{e:0:tildeRho}, with initial datum given by
\begin{align}
      \tilde \rho(x,0) =&\, \sum_{a=1}^Q \frac{1}{\lambda^a} \big( \bS[ \check\rho_a (\cdot,0,0) ] - \check\rho_a (\cdot,0,0) \big).      \label{e:EXP2:trho_ini}
\end{align}
The force $\tilde f$ is given by \eqref{e:2:tildeRho}-\eqref{e:8:tildeRho}, and $L_1, L_0, L_{-1}$ are given by \eqref{e:10:tildeRho}-\eqref{e:18:tildeRho}. However, $\hat \rho_{Q-2}, \hat \rho_{Q-1}, \check \rho_{Q-1}, \bar \rho_{Q-1}, \hat \rho_{Q}, \check \rho_{Q}, \bar \rho_{Q}$ are given in \ref{c:l:exp2:spt}, \ref{c:l:exp2:tpr} and \ref{c:l:exp2:rsd}.
\end{enumerate}

\end{lemma}

\begin{remark}    \label{r:shareCorrector}
The $\hat I_{k}^{(\balpha)}$, $\hat h_{k}^{(\balpha,l)}$, $\hat \eta_{k}^{(\balpha,l)}$, $\hat \chi_{k}^{(\balpha,l)}$, $\check I_{k}^{(\balpha)}$, $\check h_{k}^{(\balpha,l)}$, $\check \eta_{k}^{(\balpha,l)}$, $\bar h_{k}^{(\balpha)}$ in Lemma \ref{l:EXP1} and Lemma \ref{l:EXP2} are the same, i.e. they do not depend on $\rho_{\ini}$ in Lemma \ref{l:EXP1} and the forcing terms in Lemma \ref{l:EXP2}. We shall prove this in Section \ref{ss:proof_expansion}. However, the $\bar \rho_n$ in Lemma \ref{l:EXP1} and Lemma \ref{l:EXP2} are different, since they solve different equations, \eqref{e:EXP1:brhon} and \eqref{e:EXP2:brhon} respectively.
\end{remark}

\begin{corollary}       \label{c:t=0Ref:EXP2}
Taking assumptions in Lemma \ref{l:EXP2}, suppose
\begin{align*}
      \Upsilon = \frac{\delta^{\sfrac12}\mr{\bar\lambda}} {\kappa\lambda^2} \cdot \lambda^{ - 2b \gamma }.
\end{align*}
Then for any $0 \leq n \leq Q$, we have
\begin{align}
      \frac{1}{\lambda^n} \check \rho_{n} (x,0,0) =&\, \sum_{a=0}^n \sum_{ [\balpha] = 1, \balpha \in \mathcal{I}_x }^{4n+1}
                  \mathfrak{C}_{n,a}^{(\balpha)} \bigg( \frac{1}{\lambda^a} \mr {\bar \partial}^{\balpha} \bar \rho_{a}(x,0) \bigg)  \nonumber \\
            +&\, \sum_{k=1}^K \sum_{ [\balpha] = 1, \balpha \in \mathcal{I}_x }^{4Q+1}
                  \mathfrak{C}_{F,n,k}^{(\balpha)} \mr {\bar \partial}^{\balpha} \bar \rho_{F,k}(x,0)      \label{e:3:t=0Ref:EXP2} \\
      \frac{1}{\lambda^n} \bar \rho_{n} (x,0) =&\, \sum_{a=0}^n \sum_{ [\balpha] = 1, \balpha \in \mathcal{I}_x }^{4n+1}
                  \mathfrak{C}_{n,a}^{(\balpha)} \, \bS \bigg( \frac{1}{\lambda^a} \mr {\bar \partial}^{\balpha} \bar \rho_{a}(x,0) \bigg) \nonumber \\
            +&\, \sum_{k=1}^K \sum_{ [\balpha] = 1, \balpha \in \mathcal{I}_x }^{4Q+1}
                  \mathfrak{C}_{F,n,k}^{(\balpha)} \, \bS \Big( \mr {\bar \partial}^{\balpha} \bar \rho_{F,k}(x,0)  \Big).     \label{e:5:t=0Ref:EXP2}
\end{align}
For $p \geq 1$, we have
\begin{align}
      \frac{1}{\lambda^n} \mr {\bar D}^p_t \bar \rho_{n} (x,0) =&\, \sum_{a=0}^n \sum_{ [\balpha] = 1, \balpha \in \mathcal{I}_x }^{4n+2p+1}
                  \mathfrak{C}_{n,a}^{(p,\balpha)} \bigg( \frac{1}{\lambda^a} \mr {\bar \partial}^{\balpha} \bar \rho_{a}(x,0) \bigg)     \nonumber \\ 
            +&\, \sum_{k=1}^K \sum_{ [\balpha] = 1, \balpha \in \mathcal{I}_x }^{4Q+2p+1}
                  \mathfrak{C}_{F,n,k}^{(p,\balpha)} \mr {\bar \partial}^{\balpha} \bar \rho_{F,k}(x,0),      \label{e:7:t=0Ref:EXP2}
\end{align}
with the following estimates for the constant coefficients
\begin{align}
      | \mathfrak{C}_{n,a}^{(\balpha)} | + | \mathfrak{C}_{n,a}^{(p,\balpha)} | 
            \lesssim&\, \lambda^{-(n-a)b\gamma},      \label{e:9:t=0Ref:EXP2} \\ 
      | \mathfrak{C}_{F,n,k}^{(\balpha)} | + | \mathfrak{C}_{F,n,k}^{(p,\balpha)} | 
            \lesssim&\, \lambda^{-(n+2)b\gamma}.    \label{e:11:t=0Ref:EXP2}
\end{align}
\end{corollary}

\begin{proof}
The proof is inductive and analogous to the proof of Corollary \ref{c:t=0Ref:EXP1}. Similar to Corollary \ref{c:t=0Ref:EXP1}, $\nabla$ and $\bar D_t$ commute at $t=0$. For $n=0$, (\ref{e:3:t=0Ref:EXP2}-\ref{e:5:t=0Ref:EXP2}) and (\ref{e:9:t=0Ref:EXP2}-\ref{e:11:t=0Ref:EXP2}) follow from \eqref{e:EXP2:crho0} and \eqref{e:EXP2:brho0_ini}. From \eqref{e:EXP2:brho0} and \eqref{e:EXP2:brho0_ini}, we have
\begin{align}     \label{e:21:t=0Ref:EXP2}
      \mr{\bar D}_t^p \bar \rho_0 (x, 0)
            - \frac{ \bar\kappa \mr{\bar\lambda}^2 }{\mr{\bar\mu}} \mr{\bar D}_t^{p-1} \mr{\bar \Delta} \bar\rho_0 ( x, 0 )
            = \frac{ \mr{\bar\lambda} } { \mr{\bar\mu}\mr\lambda } \sum_{k, \balpha, l} 
            \mr{\bar \divr} \Big( h_{F,k}^{(\balpha,l)} (x, 0) \mr{\bar D}_t^{p-1} \mr {\bar D}^{\balpha} \bar \rho_{F,k} (x, 0) \Big) 
                  \langle \eta_{F,k}^{(\balpha,l)} \chi_{F,k}^{(\balpha,l)} \rangle_{\xi,\tau}.
\end{align}
From \eqref{e:1:mScaleRela}, \eqref{e:EXP2:bhSF_Pest} and Assumption \ref{a:homSmallscl}, we can derive 
\begin{align}
      \frac{ \mr{\bar\lambda} } { \mr{\bar\mu}\mr\lambda } \big| h_{F,k}^{(\balpha,l)}(x, 0) \big|
            \lesssim&\, \frac{\delta^{\sfrac12}}{\bar\kappa\lambda} \cdot \lambda^{-2b\gamma}.
\end{align}
Note that $\balpha$ in \eqref{e:21:t=0Ref:EXP2} may contain both space derivatives and transport derivatives. To derive (\ref{e:7:t=0Ref:EXP2}-\ref{e:11:t=0Ref:EXP2}) for $n=0$, we do an induction on $p$ and use transmission condition of $\{\bar\rho_{F,k}\}_k$ in Assumption \ref{a:homSmallscl} to express transport derivative in space derivatives.

Next, assume the desired equalities and estimates for $a \leq n$. We use \eqref{e:EXP2:crhon}-\eqref{e:EXP2:crhon:F} to derive the desired equality in \eqref{e:3:t=0Ref:EXP2}. We write
\begin{equation}     \label{e:25:t=0Ref:EXP2}
      \begin{split}
      \frac{1}{\lambda^{n+1}} \check \rho_{n+1} (x,0,0) =&\, \sum_{a=0}^{n} \sum_{[\balpha]=1, \balpha \in \mathcal{I}_* }^{2(n+1-a)} \sum_{l} 
            \bigg( \frac{1}{\lambda^{n+1-a}} \check h^{(\balpha, l)}_{n+1-a} \bigg) (x,0) \check \eta^{(\balpha, l)}_{n+1-a} (0) 
            \bigg( \frac{1}{\lambda^a} \mr {\bar D}^{\balpha} \bar \rho_a \bigg) (x,0)    \\ 
            +&\, \sum_{k=1}^K \sum_{ [\balpha]=1, \balpha \in \mathcal{I}_* }^{2n+1+Q} \sum_{l}  
            \bigg( \frac{1}{\lambda^{n+1}}\check h^{(\balpha,l)}_{F,k,n+1} \bigg) (x,0) \check \eta^{(\balpha,l)}_{F,k,n+1} (0) \mr {\bar D}^{\balpha} \bar \rho_{F,k}(x,0).
      \end{split}
\end{equation}
Using the induction hypothesis, \eqref{e:EXP1:ch_Pest} and \eqref{e:EXP2:chF_Pest}, we deduce \eqref{e:3:t=0Ref:EXP2}. Here, we use Remark \ref{r:PolyDeriP3}. To express the transport derivatives $\balpha$ in \eqref{e:25:t=0Ref:EXP2} into space derivatives, we use transmission condition of $\{\bar\rho_{F,k}\}_k$ in Assumption \ref{a:homSmallscl} and \eqref{e:7:t=0Ref:EXP2}. \eqref{e:5:t=0Ref:EXP2} for $n+1$ follows from \eqref{e:EXP2:brhon_ini} and \eqref{e:3:t=0Ref:EXP2}.

From \eqref{e:EXP2:brhon} and \eqref{e:EXP2:brhon_ini}, we have
\begin{equation}     \label{e:31:t=0Ref:EXP2}
\begin{split}
      \mr{\bar D}_t^p \bar \rho_{n+1} (x,0) - \frac{\bar\kappa \mr{\bar\lambda}^2}{\mr{\bar\mu}} \mr{\bar D}_t^{p-1} \mr{\bar \Delta} \bar\rho_{n+1} (x,0) 
            =&\, \frac{1}{\mr{\bar\mu}} \sum_{a=0}^{n} \sum_{ [\balpha]=1, \balpha \in \mathcal{I}_* }^{2(n-a)+3} 
            \mr{\bar \divr} \Big( \bar h^{(\balpha)}_{n-a+1} (x,0) \mr{\bar D}_t^{p-1} \mr {\bar D}^{\balpha} \bar \rho_a (x,0) \Big)     \\ 
            +&\, \frac{1}{\mr{\bar\mu}} \sum_{k=1}^K \sum_{ [\balpha]=1, \balpha \in \mathcal{I}_* }^{2n+Q} 
            \mr{\bar \divr} \Big( \bar h_{F,k,n}^{(\balpha)} (x,0) \mr{\bar D}_t^{p-1} \mr {\bar D}^{\balpha} \bar \rho_{F,k} (x,0) \Big).
\end{split}
\end{equation}
From the parameter relations in \eqref{e:1:mScaleRela}, \eqref{e:EXP1:bhS_Pest} and \eqref{e:EXP2:bhSF_Pest}, we can derive
\begin{align}
      \frac{1}{\mr{\bar\mu}} \, \big| \bar h^{(\balpha)}_{n-a+1} (x,0) \big|
            \lesssim&\, \lambda^{-(n - a) b \gamma},        \label{e:35:t=0Ref:EXP2} \\ 
      \frac{1}{\mr{\bar\mu}} \, \big| \bar h_{F,k,n}^{(\balpha)} (x,0) \big|
            \lesssim&\, \lambda^{-(n+2) b \gamma}.          \label{e:37:t=0Ref:EXP2}
\end{align}
Then \eqref{e:7:t=0Ref:EXP2} for $n+1$ follows from an induction on $p$ based on (\ref{e:7:t=0Ref:EXP2}-\ref{e:11:t=0Ref:EXP2}) and \eqref{e:31:t=0Ref:EXP2}. Again, we use that $\nabla$ and $\bar D_t$ commute at $t=0$ and $\bar h^{(\cdot)}_{\cdot} (\cdot,0)$ are constant functions. We also use the transmission condition of $\{\bar\rho_{F,k}\}_k$ in Assumption \ref{a:homSmallscl}.
\end{proof}

\begin{lemma}     \label{l:t=0:EXP2}
Given the assumptions in Lemma \ref{l:EXP2}, for some $\Upsilon_{\ini}>0$, suppose we have
\begin{align*}
      | \bar \rho_{F,k} (\cdot, 0) |_{ \fS(q_*) } + | \bar \rho_{F,k} (\cdot, 0) |_{ \fS(q+1) }
            \lesssim \lambda^{-(s-2)b\gamma} \Upsilon_{\ini},   \quad 
      \Upsilon =&\, \frac{\delta^{\sfrac12}\mr{\bar\lambda}} {\kappa\lambda^2} \cdot \lambda^{ - 2b \gamma }.     \\ 
      \| \bar \rho_{F,k} (\cdot, 0) \|_2 + \sum_{j=q_*}^{q} \lambda_j^{-\gamma} | \bar \rho_{F,k} (\cdot, 0) |_{\fS(j)}
            \lesssim \lambda^{-(s-2)b\gamma} \Upsilon_{\ini},   \quad 
      \Upsilon =&\, \frac{\delta^{\sfrac12}\mr{\bar\lambda}} {\kappa\lambda^2} \cdot \lambda^{ - 2b \gamma }.
\end{align*}
Then we have for $n \geq 0$,
\begin{align}
      \Big\| \frac{1}{\lambda^{n}} \bar \rho_n (\cdot, 0) \Big\|_2 
            + \sum_{j=q_*}^{q} \lambda_j^{-\sfrac{\gamma}{2}} \Big| \frac{1}{\lambda^{n}} \bar \rho_n (\cdot, 0) \Big|_{\fS(j)}
            &\,\lesssim \lambda^{- (n+s) b \gamma } \Upsilon_{\ini},      \label{e:5:t=0:EXP2} \\
      \sum_{j=q_*}^{q} \lambda_j^{-\sfrac{\gamma}{2}} \Big| \frac{1}{\lambda^{n}} \check\rho_n (\cdot,0,0) \Big|_{\fS'(j)}
            &\,\lesssim \lambda^{- (n+s) b \gamma } \Upsilon_{\ini},      \label{e:7:t=0:EXP2} \\
      \frac{1}{\lambda^{n}} \big| \bS[ \check\rho_n (\cdot,0,0) ] - \check\rho_n (\cdot,0,0) \big|_{\mathfrak{P}(q+1)}
            &\,\lesssim \lambda^{- (Q+n+s) b \gamma } \Upsilon_{\ini}.       \label{e:8:t=0:EXP2}
\end{align}
Here, $\{\bar \rho_a\}_a$ is given in \eqref{e:EXP2:brho0} and \eqref{e:EXP2:brhon} with initial conditions in \eqref{e:EXP2:brho0_ini}-\eqref{e:EXP2:brhon_ini}.
\end{lemma}

\begin{proof}
The proof is completely analogous to the proof of Lemma \ref{l:t=0:EXP1} and follows from Lemma \ref{l:EXP2} and Corollary \ref{c:t=0Ref:EXP2}.
\end{proof}

\subsection{Proof of the expansion lemmas}      \label{ss:proof_expansion}

In this section, we prove the expansion Lemma \ref{l:EXP1} and Lemma \ref{l:EXP2}, also Remark \ref{r:stationarity} and Remark \ref{r:shareCorrector}, in three steps.

The first step is given in Section \ref{sss:setUpAsymp}, we show that it suffices to solve two asymptotic systems of differential equations, respectively for Lemma \ref{l:EXP1} and Lemma \ref{l:EXP2}.

The asymptotic systems, obtained in Section \ref{sss:setUpAsymp} are solved by induction. The second step is to solve the initial steps of these asymptotic systems for Lemma \ref{l:EXP1} and Lemma \ref{l:EXP2}. This is done separately in Section \ref{sss:solInitExp1} and Section \ref{sss:solInitExp2}.

The third step is the induction for solving the asymptotic systems for both expansion lemmas. Specifically, we solve for residual correctors, temporal correctors and spatial correctors, respectively in Section \ref{sss:solAsptRSD}, Section \ref{sss:solAsptTPR} and Section \ref{sss:solAsptSPT}. The proofs of Lemma \ref{l:EXP1} and Lemma \ref{l:EXP2} are completely analogous, so we only present it for Lemma \ref{l:EXP1} and comment on their minor difference.

\subsubsection{Set up the asymptotic system for Lemma \ref{l:EXP1} and Lemma \ref{l:EXP2}}      \label{sss:setUpAsymp}

We write \eqref{e:hom:meq} in the following form
\begin{align}     \label{e:aspt:meq}
    L \varrho = \bar D_t \varrho - \divr ( A \nabla \varrho ) = 0,
\end{align}
We make an \textit{ansatz} that $\varrho$ admits an expansion of the following form
\begin{align}     \label{e:aspt:ansatzPre}
      \varrho (x,t) &= \sum_{n = 0}^{Q} \frac{1}{\lambda^n} 
            \rho_{n} (x,t, \lambda \Phi(x,t), \mu t ) + \tilde \rho(x,t).
\end{align}
for some function $\rho_{n} : \T^2 \times [0,1] \times \T^2 \times \T \rightarrow \R $. 

Define $ \hat \rho_{n} : \T^2 \times [0,1] \times \T^2 \times \T \rightarrow \R $ with arguments $(x,t,\xi,\tau)$ and $ \check \rho_{n} : \T^2 \times [0,1] \times \T \rightarrow \R $ with arguments $(x,t,\tau)$ via
\begin{align}
      \bar \rho_n(x,t) :=&\, \langle \rho_{n} (x,t, \cdot, \cdot ) \rangle_{\xi,\tau},    \label{e:1:asptDecomp} \\ 
      \check \rho_n(x,t,\tau) :=&\, \langle \rho_{n} (x,t, \cdot, \tau ) \rangle_{\xi} - \langle \rho_{n} (x,t, \cdot, \cdot ) \rangle_{\xi,\tau},    \label{e:2:asptDecomp} \\ 
      \hat \rho_n(x,t,\xi,\tau) :=&\, \rho_{n} (x,t, \xi, \tau ) - \langle \rho_{n} (x,t, \cdot, \tau ) \rangle_{\xi}.        \label{e:3:asptDecomp}
\end{align}
Then we have the natural decomposition for $\rho_n$ with $n \geq 0$
\begin{align}     \label{e:4:asptDecomp}
      \rho_n(x,t,\xi,\tau) = \hat \rho_n(x,t,\xi,\tau) + \check \rho_n(x,t,\tau) + \bar \rho_n(x,t)
\end{align}
and subsequently
\begin{align}     \label{e:aspt:ansatz}
      \varrho (x,t) &= \sum_{n = 0}^{Q} \frac{1}{\lambda^n} 
            \big( \hat \rho_{n} (x,t, \lambda \Phi(x,t), \mu t ) 
            + \check \rho_{n} (x, t, \mu t) + \bar \rho_{n} (x, t) \big) 
            + \tilde \rho(x,t). 
\end{align}
We have, for any $x,t$,
\begin{align}     \label{e:16:SetAspt}
      \langle \hat \rho_{n}( x,t, \cdot, \tau ) \rangle_\xi = 0 \text{ for any } \tau, \quad
      \langle \check \rho_{n} ( x,t, \cdot ) \rangle_\tau = 0.
\end{align}

\begin{lemma}     \label{l:operator}
Under Definition \ref{d:homParaI} and Definition \ref{d:homFlowmap}, for a general function $\varrho: \T^2 \times [0,1] \rightarrow \R$ given by 
\begin{align}     \label{e:0:operator}
      \varrho(x,t) = \rho( x,t, \lambda \Phi(x,t), \mu t ),
\end{align}
not necessarily solving \eqref{e:hom:meq}, we have
\begin{align}     \label{e:1:operator}
      L\varrho = \lambda^2 L_2 \rho + \lambda L_1 \rho + L_0 \rho
            + \lambda^{-1} L_{-1} \rho
\end{align}
with
\begin{align*}
      L_2 \rho =&\, - \kappa \Delta_\xi \rho - \frac{\delta^{\sfrac12}}{\lambda} (\det \nabla\Phi)^2 \divr_\xi H \cdot \nabla_\xi \rho 
            + \frac{1}{\lambda} Z \cdot \nabla_\xi \rho , \\
      L_1 \rho =&\, - \frac{\delta^{\sfrac12}}{\lambda} \partial_{\xi_i} H_{ij} \partial_{x_j} \rho + \frac{\mu}{\lambda} \partial_\tau \rho
            + \varepsilon \kappa \lambda \sum_m S_{m,ij} \vartheta_m \partial_{\xi_i\xi_j} \rho 
            := L_{11} \rho + L_{12} \rho + L_{13} \rho \\ 
      L_0 \rho =&\, \bar D_t \rho - \kappa \Delta_x \rho - \varepsilon \frac{\delta^{\sfrac12}}{\lambda} \sum_m \sigma_m \partial_{x_i} \Omega_{m,ij} H_{12} \partial_{x_j} \rho
            - \varepsilon \delta^{\sfrac12} \sum_m \varphi_m \partial_{\xi_i} ( H E_m^T )_{ij} \partial_{x_j} \rho \\ 
            +&\, \varepsilon \delta^{\sfrac12} \sum_m \varphi_m \partial_{x_j} ( HE_m^T )_{ij} \partial_{\xi_i} \rho 
            - \varepsilon \kappa \lambda \sum_m \phi_m \partial_{x_j} B_{m,ij} \partial_{\xi_i} \rho
            - 2 \kappa \lambda \partial_{\xi_i x_i} \rho \\ 
            -&\, 2 \varepsilon \kappa \lambda \sum_m \phi_m B_{m,ij} \partial_{\xi_i x_j} \rho  
            := \sum_{1 \leq k \leq 7} L_{0k} \rho \\ 
      L_{-1} \rho =&\, \frac{1}{\lambda} \phi_* \omega_* \cdot \nabla_\xi \rho.
\end{align*}
Note that $L_{01} \rho = D_t \rho - \kappa \Delta_x \rho$, and $L_{0k}, k \geq 2$ is defined by the order of relevant terms.
\end{lemma}

\begin{proof}
We omit the argument $(x,t, \lambda \Phi(x,t), \mu t)$ for the function $\rho$ in this proof. From \eqref{e:0:operator} and \eqref{e:hom:mat}, we compute the elliptic part 
\begin{equation}        \label{e:5:operator}
\begin{split}
      \divr( A \nabla \varrho ) =&\, \lambda^2 \divr_\xi ( \nabla \Phi A \nabla \Phi^T \nabla_\xi \rho ) + \lambda \divr_\xi ( \nabla \Phi A \nabla_x \rho ) \\ 
            &\,+ \lambda \divr_x ( A \nabla \Phi^T \nabla_\xi \rho )
            + \divr_x ( A \nabla_x \rho ).
\end{split}
\end{equation}
Using \eqref{e:hom:mat}, we compute
\begin{align}
      \divr_\xi ( \nabla \Phi A \nabla \Phi^T \nabla_\xi \rho ) 
            =&\, \kappa ( \nabla\Phi\nabla\Phi^T )_{ij} \partial_{\xi_i\xi_j} \rho + \frac{\delta^{\sfrac12}}{\lambda} (\det \nabla\Phi)^2 \partial_{\xi_i} H_{ij} \partial_{\xi_j} \rho       \label{e:6:operator} \\ 
            =&\, \kappa \Delta_\xi \rho + \frac{\delta^{\sfrac12}}{\lambda} (\det \nabla\Phi)^2 \partial_{\xi_i} H_{ij} \partial_{\xi_j} \rho 
            + \kappa ( \nabla\Phi\nabla\Phi^T - \Id )_{ij} \partial_{\xi_i\xi_j} \rho,        \label{e:9:operator}
\end{align}
\begin{align}
      \divr_\xi ( \nabla \Phi A &\, \nabla_x \rho )
            + \divr_x ( A \nabla \Phi^T \nabla_\xi \rho )    \nonumber \\ 
            =&\, ( \nabla\Phi A + \nabla\Phi A^T )_{ij} \partial_{\xi_i x_j} \rho
            + \partial_{\xi_i} (\nabla\Phi A)_{ij} \partial_{x_j} \rho 
            + \partial_{x_j} ( A \nabla\Phi^T )_{ji} \partial_{\xi_i} \rho    \label{e:10:operator} \\ 
            =&\, 2\kappa \nabla\Phi_{ij} \partial_{\xi_i x_j} \rho 
            + \kappa \partial_{x_j} \nabla \Phi_{ij}  \partial_{\xi_i} \rho
            + \frac{\delta^{\sfrac12}}{\lambda} \partial_{\xi_i} ( \det \nabla\Phi H \adj \nabla\Phi^T )_{ij} \partial_{x_j} \rho       \label{e:11:operator} \\ 
            &\,+ \frac{\delta^{\sfrac12}}{\lambda} \partial_{x_j} ( \det \nabla\Phi \adj \nabla\Phi H )_{ji} \partial_{\xi_i} \rho,      \label{e:12:operator}
\end{align}
and
\begin{align}
      \divr_x ( A \nabla_x \rho ) = \kappa \Delta_x \rho + \frac{\delta^{\sfrac12}}{\lambda}
      \partial_{x_i} ( \adj\nabla\Phi H \adj\nabla\Phi^T )_{ij} \partial_{x_j} \rho.     \label{e:14:operator}
\end{align}

For the transport term, we deduce from \eqref{e:0:operator} and \eqref{e:12:homFlowmap}, that
\begin{align}
      \bar D_t \varrho =&\, \mu \partial_\tau \rho + \partial_t \rho + \lambda \partial_t \Phi_i \partial_{\xi_i} \rho + \bar u \cdot \nabla_x \rho + \lambda \bar u \cdot \nabla \Phi_i \partial_{\xi_i} \rho    \nonumber \\ 
      =&\, \mu \partial_\tau \rho + \bar D_t \rho + \lambda \bar D_t \Phi \cdot \nabla_\xi \rho       \nonumber \\ 
      =&\, \mu \partial_\tau \rho + \bar D_t \rho + \lambda Z \cdot \nabla_\xi \rho + \lambda^{-2} \Omega \cdot \nabla_\xi \rho.     \label{e:18:operator} 
\end{align}

From \eqref{e:aspt:meq}, \eqref{e:0:operator}, (\ref{e:5:operator}-\ref{e:18:operator}), (\ref{e:6:homFlowmap}-\ref{e:10:homFlowmap}) and (\ref{e:12:homFlowmap}-\ref{e:18:homFlowmap}), we have \eqref{e:1:operator} with
\begin{align}
      L_2 \rho =&\, - \kappa \Delta_\xi \rho - \frac{\delta^{\sfrac12}}{\lambda} (\det \nabla\Phi)^2 \divr_\xi H \cdot \nabla_\xi \rho 
            + \frac{1}{\lambda} Z \cdot \nabla_\xi \rho ,   \label{e:20:operator} \\
      L_1 \rho =&\, - \frac{\delta^{\sfrac12}}{\lambda} \partial_{\xi_i} H_{ij} \partial_{x_j} \rho + \frac{\mu}{\lambda} \partial_\tau \rho
            + \varepsilon \kappa \lambda \divr_\xi ( \Xi \nabla_\xi \rho ), 
                   \label{e:22:operator} \\ 
      L_0 \rho =&\, \bar D_t \rho - \kappa \Delta_x \rho - \frac{\delta^{\sfrac12}}{\lambda} \partial_{x_i} ( \adj \nabla \Phi H \adj \nabla \Phi^{T} )_{ij} \partial_{x_j} \rho 
            - \varepsilon \delta^{\sfrac12} \sum_m \varphi_m \partial_{\xi_i} ( H E_m^T )_{ij} \partial_{x_j} \rho      \label{e:24:operator} \\ 
            &+ \varepsilon \delta^{\sfrac12} \sum_m \varphi_m \partial_{x_j} ( HE_m^T )_{ij} \partial_{\xi_i} \rho 
            - \kappa \lambda \partial_{x_j} \nabla \Phi_{ij} \partial_{\xi_i} \rho - 2 \kappa \lambda \nabla \Phi_{ij} \partial_{\xi_i x_j} \rho,     \label{e:26:operator} \\ 
      L_{-1} \rho =&\, \frac{1}{\lambda} \Omega \cdot \nabla_\xi \rho.        \label{e:28:operator}
\end{align}

In particular, we remark that, in above computation, we decompose the last term in line \eqref{e:11:operator} into two terms, with one being the first term in line \eqref{e:22:operator}, the other one being the last term in line \eqref{e:24:operator}. Then we use \eqref{e:6:homFlowmap}, \eqref{e:16:homFlowmap} and \eqref{e:18:homFlowmap} and above to conclude the proof.
\end{proof}

\underline{For the proof of Lemma \ref{l:EXP1}}, combining the ansatz \eqref{e:aspt:ansatzPre} and Lemma \ref{l:operator}, we have
\begin{equation}        \label{e:asptExp}
\begin{split}
      L \varrho =&\, \lambda^2 L_2 \rho_0 + \lambda ( L_2 \rho_1 + L_1 \rho_0 ) 
            + ( L_2 \rho_2 + L_1 \rho_1 + L_0 \rho_0 ) \\ 
            &\,+ \frac{1}{\lambda} ( L_2 \rho_3 + L_1 \rho_2 + L_0 \rho_1 + L_{-1} \rho_0 ) + \ldots \\ 
            &\,+ \frac{1}{\lambda^{n+1}} ( L_2 \rho_{n+3} + L_1 \rho_{n+2} + L_0 \rho_{n+1} + L_{-1} \rho_n ) + \ldots \\ 
            &\,+ \frac{1}{\lambda^{Q-1}} ( L_2 \rho_{Q} + L_1 \rho_{Q-1} + L_0 \rho_{Q-2} + L_{-1} \rho_{Q-3} ) \\ 
            &\,+ \frac{1}{\lambda^{Q}} ( L_1 \rho_{Q} + L_0 \rho_{Q-1} + L_{-1} \rho_{Q-2} ) \\ 
            &\,+ \frac{1}{\lambda^{Q+1}} ( L_0 \rho_{Q} + L_{-1} \rho_{Q-1} ) \\ 
            &\,+ \frac{1}{\lambda^{Q+2}}  L_{-1} \rho_{Q} + L\tilde\rho
\end{split}
\end{equation}

Now we seek for $\{\rho_n\}_{0 \leq n \leq Q}$ such that the following equations are satisfied
\begin{align}
      L_2 \rho_0 =& 0, \label{e:aspt:hrc1} \\ 
      L_2 \rho_1 + L_1 \rho_0 =& 0, \label{e:aspt:hrc2} \\ 
      L_2 \rho_2 + L_1 \rho_1 + L_0 \rho_0 =& 0, \label{e:aspt:hrc3} \\ 
      L_2 \rho_{n+3} + L_1 \rho_{n+2} + L_0 \rho_{n+1} + L_{-1} \rho_n =& 0. \label{e:aspt:hrc4}
\end{align}

We use $\tilde \rho: \T^2 \times [0,1] \rightarrow \R$ to denote the solution of the following equation
\begin{equation}  \label{e:eTildeRho}
\begin{split}
      D_t \tilde\rho - \divr( A \nabla \tilde\rho ) 
      = \tilde f := &\, -\frac{1}{\lambda^{Q}} ( L_1 \rho_{Q} + L_0 \rho_{Q-1} + L_{-1} \rho_{Q-2} ) \\ 
            &\,- \frac{1}{\lambda^{Q+1}} ( L_0 \rho_{Q} + L_{-1} \rho_{Q-1} ) \\ 
            &\,- \frac{1}{\lambda^{Q+2}}  L_{-1} \rho_{Q} , \\ 
      \tilde \rho(x,0) =&\, 0.
\end{split}
\end{equation}

Note from Remark \ref{r:2:homFlowmap} that $\Phi(x,0) = x$. We impose the following initial conditions
\begin{align}
      \rho_0(x, 0, \lambda x, 0) =&\, \mathbb{S} \rho_{\ini} (x),      \label{e:0:initialC} \\   
      \bar\rho_n(x, 0) + \bS[ \check\rho_n (\cdot,0,0) ] (x) =&\, 0, \quad \text{for any } n \geq 1,      \label{e:2:initialC} \\ 
      \tilde \rho(x,0) + \sum_{a=0}^Q \frac{1}{\lambda^a} 
            \rho_a (x,0,\lambda x, 0) =&\, \rho_{\ini} (x),      \label{e:4:initialC} \\
      \tilde \rho(x,0) + \mathbb{S} \rho_{\ini} (x) - \rho_{\ini} (x)
            + \sum_{a=0}^Q \frac{1}{\lambda^a} \big( \check\rho_a (\cdot,0,0) -&\, \bS[ \check\rho_a (\cdot,0,0) ] \big) = 0.      \label{e:6:initialC}
\end{align}

\textit{If we have $\{\rho_n\}_{0 \leq n \leq Q}$ and $\tilde \rho$ such that \eqref{e:aspt:hrc1}-\eqref{e:2:initialC} hold, the function $\varrho$ defined by the ansatz \eqref{e:aspt:ansatzPre} solves \eqref{e:hom:meq}. Thus, to prove Lemma \ref{l:EXP1}, it remains to solve \eqref{e:aspt:hrc1}-\eqref{e:2:initialC}.}

\underline{For the proof of Lemma \ref{l:EXP2}}, there are only two difference with Lemma \ref{l:EXP1}. For the first difference, the system for $\hat \rho_n, \check \rho_n$ and $\bar \rho_n$ is now given by
\begin{align}
      L_2 \rho_0 =&\, 0, \label{e:asptFF:hrc1} \\ 
      L_2 \rho_1 + L_1 \rho_0 =&\, 0, \label{e:asptFF:hrc2} \\ 
      L_2 \rho_2 + L_1 \rho_1 + L_0 \rho_0 =&\, G, \label{e:asptFF:hrc3} \\ 
      L_2 \rho_{n+3} + L_1 \rho_{n+2} + L_0 \rho_{n+1} + L_{-1} \rho_n =&\, 0. \label{e:asptFF:hrc4}
\end{align}
For the second difference, the set of initial conditions are given by
\begin{align}
      \rho_0(x, 0, \lambda x, 0) =&\, 0,      \label{e:0:initialCFF} \\   
      \rho_n(x, 0, \lambda x, 0) =&\, 0, \quad \text{for any } n \geq 1.      \label{e:2:initialCFF} \\ 
      \tilde \rho(x,0) =&\, 0.      \label{e:6:initialCFF}
\end{align}
Here, the function $G: \T^2 \times [0,1] \times \T^2 \times \T \rightarrow \R$ is given by computing the forcing term in \eqref{e:hom:eqFast} with $F$ given by Assumption \ref{a:homSmallscl}. Then we have 
\begin{align}
      \mr \lambda^{-1} \divr F =&\, G (x,t, \lambda \Phi(x,t), \mu t ),   \label{e:4:compFHom} \\ 
      G(x,t,\xi,\tau) =&\, \frac{ \mr{\bar\lambda} } { \mr\lambda } 
      \sum_{k=1}^K \sum_{ [\balpha] = 1, \balpha \in \mathcal{I}_* }^Q \sum_{l = 1}^{J} 
      \mr {\bar\divr}_x \Big( h_{F,k}^{(\balpha,l)} \mr {\bar D}^{\balpha} \bar \rho_{F,k} \Big) 
            \eta_{F,k}^{(\balpha,l)}(\tau) \chi_{F,k}^{(\balpha,l)} (\xi)        \label{e:6:compFHom} \\ 
      +&\, \frac{\lambda} { \mr\lambda }
      \sum_{k=1}^K \sum_{ [\balpha] = 1, \balpha \in \mathcal{I}_* }^Q \sum_{l = 1}^{J} 
      \Big[ \Id + \varepsilon \sum_m B_m \phi_m (\tau) \Big]_{ij} h_{F,k,j}^{(\balpha,l)} \mr {\bar D}^{\balpha} \bar \rho_{F,k} 
            \eta_{F,k}^{(\balpha,l)}(\tau) \partial_{\xi_i} \chi_{F,k}^{(\balpha,l)} (\xi)    \label{e:8:compFHom}
\end{align}
Here, $h_{F,k,j}^{(\cdot,\cdot)}$ refers to the $j$-th component of the vector-valued function $h_{F,k}^{(\cdot,\cdot)}$, and we omit the argument $(x,t)$ for $h_{F,k}^{(\cdot,\cdot)}$, $\bar\rho_{F,k}$ and $B_m$.

\textit{Similar to above, to prove Lemma \ref{l:EXP1}, it remains to solve \eqref{e:asptFF:hrc1}-\eqref{e:6:initialCFF}.}

\begin{remark}    \label{r:hatcheckbar}
For $\rho(x,t,\xi,\tau) = \hat\rho(x,t,\xi,\tau) + \check\rho(x,t,\tau) + \bar\rho(x,t)$ with
\begin{align}     \label{e:0:hatcheckbar}
      \langle \hat \rho( x,t, \cdot, \tau ) \rangle_\xi = 0 \text{ for any } x,t,\tau, \quad
      \langle \check \rho ( x,t, \cdot ) \rangle_\tau = 0 \text{ for any } x,t,
\end{align}
we have
\begin{align}
      \langle L_2 \rho \rangle_\xi =&\, 0,      \label{e:2:hatcheckbar} \\ 
      L_1 \rho =&\, L_{11} ( \hat\rho + \check\rho + \bar\rho ) + L_{12} ( \hat\rho + \check\rho ) + L_{13} \hat \rho,        \label{e:4:hatcheckbar} \\ 
      \langle L_1 \rho \rangle_\xi =&\, \langle L_{11} \hat\rho \rangle_\xi + L_{12} \check \rho,     \label{e:6:hatcheckbar} \\ 
      \langle L_1 \rho \rangle_{\xi,\tau} =&\, \langle L_{11} \hat\rho \rangle_{\xi,\tau}.      \label{e:8:hatcheckbar}
\end{align}
Note that $\langle H(\xi,\tau) \rangle_\xi = 0$, then we also have
\begin{align}
      L_0 \rho =&\, ( L_{01} + L_{02} + L_{03} ) ( \hat\rho + \check\rho + \bar\rho ) + L_{04} \hat\rho + L_{05} \hat\rho + L_{06} \hat\rho + L_{07} \hat\rho    \label{e:12:hatcheckbar} \\ 
      \langle L_0 \rho \rangle_\xi =&\, L_{01} \check\rho + L_{01} \bar\rho
            + \langle L_{02} \hat\rho \rangle_{\xi} + \langle L_{03} \hat\rho \rangle_{\xi} + \langle L_{04} \hat\rho \rangle_{\xi},      \label{e:14:hatcheckbar} \\ 
      \langle L_0 \rho \rangle_{\xi,\tau} =&\, L_{01} \bar\rho
            + \langle L_{02} \hat\rho \rangle_{\xi,\tau} + \langle L_{03} \hat\rho \rangle_{\xi,\tau} + \langle L_{04} \hat\rho \rangle_{\xi,\tau}.    \label{e:16:hatcheckbar}
\end{align}
Moreover, we have
\begin{align}    \label{e:18:hatcheckbar} 
      \langle L_{-1} \rho \rangle_\xi = 0.
\end{align}
\end{remark}

\begin{remark}[Disjoint support of $\eta_1$ and $\eta_2$]   \label{r:alter_shear}
When solving an differential equation in $\xi$ of form
\begin{align}     \label{e:alter_shear:2}
      L_2 \rho = f_1(x,t) g_1(\xi_1) \eta_1(\tau) + f_2(x,t) g_2(\xi_2) \eta_2(\tau),
\end{align}
with general smooth functions $f_1, f_2 : \T^2 \times [0,1]$, $g_1, g_2: \T \rightarrow \R$ and $\eta_1, \eta_2$ given by Lemma \ref{l:cutoff}. We use the fact that $\eta_1$ and $\eta_2$ has disjoint support. We claim that a solution to \eqref{e:alter_shear:2} equal to the solution of 
\begin{align}     \label{e:alter_shear:3}
      - \kappa \Delta_\xi \rho_* = f_1(x,t) g_1(\xi_1) \eta_1(\tau) + f_2(x,t) g_2(\xi_2) \eta_2(\tau).
\end{align}
Indeed, the solution of \eqref{e:alter_shear:3} is given by
\begin{align}     \label{e:alter_shear:4}
      \rho_* = - f_1 \eta_1 \partial_{\xi_1}^{-2} g_1 - f_2 \eta_2 \partial_{\xi_2}^{-2} g_1.
\end{align}
From Definition \ref{d:homFlowmap}, \eqref{e:hom:MatEnt}, \eqref{e:hom:sin}, we have
\begin{align}
      \divr H \cdot \nabla_\xi =&\, \eta_2 \partial_{\xi_2} \Pi_2 \partial_{\xi_1}
      - \eta_1 \partial_{\xi_1} \Pi_1 \partial_{\xi_2}, \\
      Z(x,t,\tau) \eta_1 =&\, Z(x,t,\tau) \eta_2 = 0
\end{align}
therefore we have $\divr H \cdot \nabla_\xi \rho_* = Z \cdot \nabla_\xi \rho_* = 0$. Thus the solution to \eqref{e:alter_shear:2} is given by $\rho = \rho_*$.
\end{remark}

\subsubsection{Solve the asymptotic system: initial steps for Lemma \ref{l:EXP1}}   \label{sss:solInitExp1}

In this step, we solve \eqref{e:aspt:hrc1}, \eqref{e:aspt:hrc2} and \eqref{e:aspt:hrc3} with initial conditions \eqref{e:0:initialC} and \eqref{e:2:initialC}, $n=1$. The goal is to deduce
\begin{align}     \label{e:SolAsptInit:14}
      \hat \rho_0, \hat \rho_1, \hat \rho_2,
            \check \rho_0, \check \rho_1, \bar \rho_0. 
\end{align}

We view \eqref{e:aspt:hrc1} as an elliptic equation in $\xi$. Recall the decomposition (\ref{e:1:asptDecomp}-\ref{e:4:asptDecomp}), the solution to \eqref{e:aspt:hrc1} is given by $\hat \rho_{0} = 0$. For \eqref{e:aspt:hrc2}, we first take $\langle \rangle_\xi$ of \eqref{e:aspt:hrc2}. From \eqref{e:0:hatcheckbar}, \eqref{e:2:hatcheckbar} and \eqref{e:6:hatcheckbar}, we have $L_{12} \check \rho_0 = 0$ and the solution is given by
\begin{align*}
      \check \rho_0 =&\, 0.
\end{align*}
Then \eqref{e:aspt:hrc2} can be written as 
\begin{align}     \label{e:SolAsptInit:16}
      L_2 \hat \rho_1 + L_{11} \bar\rho_0 = 0.
\end{align}
Recalling Remark \ref{r:alter_shear}, a solution to \eqref{e:SolAsptInit:16} is given by
\begin{align}
      \hat \rho_1(x,t,\xi,\tau) =&\, \frac{\delta^{\sfrac12}}{\kappa\lambda} \eta_1(\tau) \hat \chi^{(2,1)}_1(\xi) \partial_{2} \bar \rho_0(x,t) - \frac{\delta^{\sfrac12}}{\kappa\lambda} \eta_2(\tau) \hat \chi^{(1,1)}_1(\xi) \partial_{1} \bar \rho_0(x,t)    \nonumber \\ 
      =&\, \frac{ \delta^{\sfrac12} \mr{\bar\lambda} }{\kappa\lambda} \eta_1(\tau) \hat \chi^{(2,1)}_1(\xi) \mr {\bar D}_{2} \bar \rho_0(x,t) - \frac{ \delta^{\sfrac12} \mr{\bar\lambda} }{\kappa\lambda} \eta_2(\tau) \hat \chi^{(1,1)}_1(\xi) \mr {\bar D}_{1} \bar \rho_0(x,t)        \label{e:SolAsptInit:17}
\end{align}
with
\begin{align*}
      \hat \chi^{(2,1)}_{1} = -\partial_{\xi_1}^{-1} \Pi_1, \quad 
      \hat \chi^{(1,1)}_{1} = -\partial_{\xi_2}^{-1} \Pi_2.
\end{align*}
Then we set $\hat I^{(2)}_{1} = \hat I^{(1)}_{1} = 1$. The coefficients in \eqref{e:EXP1:hrhon} for $n=1$ are given by
\begin{equation}        \label{e:SolAsptInit:18}
\begin{split}
      \hat \eta^{(2,1)}_{1}(\tau) = \eta_1(\tau) ,&\, \quad 
      \hat \eta^{(1,1)}_{1}(\tau) = \eta_2(\tau) , \\ 
      \hat h^{(2,1)}_{1}(x,t) = \frac{ \delta^{\sfrac12} \mr{\bar\lambda} } {\kappa\lambda},&\, \quad 
      \hat h^{(1,1)}_{1}(x,t) = - \frac{ \delta^{\sfrac12} \mr{\bar\lambda} }{\kappa\lambda}.
\end{split}
\end{equation}
These lead to \eqref{e:EXP1:hrho0}, \eqref{e:EXP1:crho0}, (\ref{e:EXP1:hrhon}-\ref{e:EXP1:hh_Pest}) for $n=1$ and $k=0$. The estimates \eqref{e:EXP1:hchi_est} and \eqref{e:EXP1:heta_est} follow from \eqref{e:hom:sin} and Lemma \ref{l:cutoff}. Also, it is obvious that $(\hat \chi^{(2,1)}_{1}, \hat \eta^{(2,1)}_{1})$ and $(\hat \chi^{(1,1)}_{1}, \hat \eta^{(1,1)}_{1})$ satisfy the shear condition.

Next we consider the equation \eqref{e:aspt:hrc3} and derive the structure of $\check \rho_1$. With Remark \ref{r:hatcheckbar}, taking $\langle \rangle_\xi$ and $\langle \rangle_{\xi,\tau}$ of \eqref{e:aspt:hrc3} gives
\begin{align}
      \langle L_{11} \hat \rho_1 \rangle_\xi 
            + L_{12} \check\rho_1 + L_{01} \bar \rho_0 =& 0,      \label{e:SolAsptInit:19} \\ 
      \langle L_{11} \hat \rho_1 \rangle_{\xi,\tau} 
            + L_{01} \bar \rho_0 =& 0,          \label{e:SolAsptInit:20} 
\end{align}
From \eqref{e:SolAsptInit:19} and \eqref{e:SolAsptInit:20}, we have
\begin{align*}
      L_{12} \check\rho_1 =&\, \langle L_{11} \hat \rho_1 \rangle_{\xi,\tau}
            - \langle L_{11} \hat \rho_1 \rangle_\xi \\ 
            =&\, \kappa \frac{\delta}{\kappa^2\lambda^2} ( \eta_1^2 - \langle \eta_1^2 \rangle_\tau ) \langle \Pi_1 \partial_{\xi_1} \hat \chi_1^{(2,1)} \rangle_\xi \partial_{22} \bar \rho_0 \\ 
            +&\, \kappa \frac{\delta}{\kappa^2\lambda^2} ( \eta_2^2 - \langle \eta_2^2 \rangle_\tau ) \langle \Pi_2 \partial_{\xi_2} \hat \chi_1^{(1,1)} \rangle_\xi \partial_{11} \bar \rho_0,
\end{align*}
therefore we have
\begin{align}     
      \check \rho_1(x,t,\tau) 
            =&\, - \frac{\delta^{\sfrac12} \mr{\bar\lambda}}{\mu} \frac{\delta^{\sfrac12} \mr{\bar\lambda}}{\kappa\lambda} \partial_\tau^{-1} ( \eta_1^2 - \langle \eta_1^2 \rangle_\tau ) (\tau) \langle \Pi_1 \partial_{\xi_1} \hat \chi_1^{(2,1)} \rangle_\xi \mr {\bar D}_{22} \bar \rho_0(x,t)    \nonumber \\ 
            &\, - \frac{\delta^{\sfrac12} \mr{\bar\lambda}}{\mu} \frac{\delta^{\sfrac12} \mr{\bar\lambda}}{\kappa\lambda} \partial_\tau^{-1} ( \eta_2^2 - \langle \eta_2^2 \rangle_\tau ) (\tau) \langle \Pi_2 \partial_{\xi_2} \hat \chi_1^{(1,1)} \rangle_\xi \mr {\bar D}_{11} \bar \rho_0(x,t)    \nonumber \\ 
            =&\, \check h_{1}^{(22,1)}(x,t) \check\eta_{1}^{(22,1)} (\tau) \mr {\bar D}_{22} \bar \rho_0 (x,t) + \check h_{1}^{(11,1)} \check\eta_{1}^{(11,1)} \mr {\bar D}_{11} \bar \rho_0 (x,t)     \label{e:SolAsptInit:22} 
\end{align}
with 
\begin{align}
      \check I^{(22)}_{1} = \check I^{(11)}_{1} = 1,       \quad
      \check I^{(\balpha)}_{1} = 0, \text{ for } \balpha \neq 11 \text{ and }
            \balpha \neq 22
\end{align}
and coefficients given by
\begin{align*}
      \check h_{1}^{(11,1)} = \frac{\delta^{\sfrac12} \mr{\bar\lambda}}{\mu} \frac{\delta^{\sfrac12} \mr{\bar\lambda}}{\kappa\lambda} \langle \Pi_2 \partial_{\xi_2} \hat \chi_1^{(1,1)} \rangle_\xi, &\,   \quad  
      \check h_{1}^{(22,1)} = \frac{\delta^{\sfrac12} \mr{\bar\lambda}}{\mu} \frac{\delta^{\sfrac12} \mr{\bar\lambda}}{\kappa\lambda} \langle \Pi_1 \partial_{\xi_1} \hat \chi_1^{(2,1)} \rangle_\xi, \\ 
      \check \eta_{1}^{(22,1)} = -\partial_\tau^{-1} ( \eta_1^2 - \langle \eta_1^2 \rangle_\tau ), &\, \quad 
      \check \eta_{1}^{(11,1)} = -\partial_\tau^{-1} ( \eta_2^2 - \langle \eta_2^2 \rangle_\tau ). 
\end{align*}
These information gives (\ref{e:EXP1:crhon}-\ref{e:EXP1:ch_Pest}) for $n=1$ and $k=0$, and \eqref{e:EXP1:ceta_est} follows from Lemma \ref{l:cutoff}.

\underline{For $\bar \rho_0$}, \eqref{e:EXP1:brho0} follows from \eqref{e:SolAsptInit:17}, \eqref{e:SolAsptInit:20} and the definition of $L_{11}$ in Lemma \ref{l:operator}. Here we also use \eqref{e:2:cutoff} and the definition of $\bar\kappa$ in \eqref{d:eddyDiff}.

\underline{Next we derive the structure of $\hat \rho_2$.} We make an ansatz, that $\hat \rho_2$ can be decomposed as 
\begin{align}     \label{e:SolAsptInit:29}
      \hat \rho_2 (x,t,\xi,\tau) =& \sum_{a=0}^{1} \sum_{[\balpha]=1, \balpha \in \mathcal{I}_*}^{3-2a} \sum_{l=1}^{\hat I_{2,a}^{(\balpha)}}
      \hat h^{(\balpha, l)}_{2,a} (x,t) \hat \eta^{(\balpha, l)}_{2,a} (\tau) \hat \chi^{(\balpha, l)}_{2,a} (\xi) \mr {\bar D}^{\balpha} \bar \rho_a (x,t). 
\end{align}

\underline{First, we verify the stationarity} mentioned in Remark \ref{r:stationarity}, i.e. to prove that for $[\balpha]=1, \balpha \in \mathcal{I}_*$ and all admissible $l$, the coefficients in \eqref{e:SolAsptInit:29} are given by
\begin{equation}  \label{e:SolAsptInit:30}
\begin{split}
      \hat I_{2,1}^{(\balpha)} =&\, \hat I_{1}^{(\balpha)}, \quad 
      \hat h^{(\balpha, l)}_{2,1} = \hat h_{1}^{(\balpha,l)}, \\ 
      \hat \eta^{(\balpha, l)}_{2,1} =&\, \hat \eta_{1}^{(\balpha,l)}, \quad 
      \hat \chi^{(\balpha, l)}_{2,1} = \hat \chi_{1}^{(\balpha,l)}.
\end{split}
\end{equation}

Indeed, the difference of \eqref{e:aspt:hrc3} and \eqref{e:SolAsptInit:19} gives
\begin{align}
      L_2 \hat \rho_2 = \langle L_1 \rho_1 \rangle_\xi - L_1 \rho_1 + \langle L_0 \rho_0 \rangle_\xi -& L_0 \rho_0.   \label{e:SolAsptInit:31}
\end{align}
From the equation above and Remark \ref{r:hatcheckbar}, we have
\begin{align}
      L_2 \hat \rho_2 =&\, - L_{11} \bar\rho_1  \label{e:SolAsptInit:32} \\
            &\,- L_{11} \check\rho_1 + \langle L_{11} \hat\rho_1 \rangle_{\xi} - L_{11} \hat\rho_1
            - ( L_{12} + L_{13} ) \hat\rho_1
            - ( L_{02} + L_{03} ) \bar\rho_0.   \label{e:SolAsptInit:33}
\end{align}
Given the information on $\hat\rho_1$ and $\check\rho_1$, i.e. \eqref{e:SolAsptInit:17} and \eqref{e:SolAsptInit:22}, the terms in the line \eqref{e:SolAsptInit:33} does not involve $\bar \rho_1$. Compare the contribution from the line \eqref{e:SolAsptInit:32} with the equation \eqref{e:SolAsptInit:16} for $\hat \rho_1$, i.e.
\begin{align*}
      L_2 \hat \rho_1 =&\, - L_{11} \bar\rho_0.
\end{align*}
We deduce that the coefficients for $a=1$ satisfies the stationarity property mentioned in Remark \ref{r:stationarity}, i.e. \eqref{e:SolAsptInit:30} holds.

\underline{Now we compute the coefficients for $a=0$}, i.e. the coefficient fields $\hat I^{(\cdot)}_{2,0}, \hat h_{2,0}^{(\cdot,\cdot)}, \hat \eta_{2,0}^{(\cdot,\cdot)}, \hat \chi_{2,0}^{(\cdot,\cdot)}$ in \eqref{e:SolAsptInit:29} of $\hat \rho_2$. This comes from the contribution of \eqref{e:SolAsptInit:33}. First, we denote
\begin{align}     \label{e:SolAsptInit:36}
      \hat I^{(\balpha)}_{2,0} := \hat I^{(\balpha)}_{2},   \quad 
      \hat h_{2,0}^{(\balpha,l)} := \hat h_{2}^{(\balpha,l)},           \quad
      \hat \eta_{2,0}^{(\balpha,l)} := \hat \eta_{2}^{(\balpha,l)}      \quad
      \hat \chi_{2,0}^{(\balpha,l)} := \hat \chi_{2}^{(\balpha,l)}.      
\end{align}

Recalling Remark \ref{r:alter_shear}, we can write down the solution 
\begin{align}
      L_2^{-1} ( -L_{11} \check\rho_1 ) = 
      &\, - \frac{\delta^{\sfrac12}}{\mu} \frac{\delta}{\kappa^2\lambda^2} \eta_1 \partial_\tau^{-1} ( \eta_1^2 - \langle \eta_1^2 \rangle_\tau ) \langle \partial_{\xi_1} \Pi_1 \Delta_\xi^{-1} \partial_{\xi_1} \Pi_1 \rangle_\xi \Delta_\xi^{-1} \partial_{\xi_1} \Pi_1 \partial_{222} \bar \rho_0       \nonumber  \\ 
      &\, + \frac{\delta^{\sfrac12}}{\mu} \frac{\delta}{\kappa^2\lambda^2} \eta_2 \partial_\tau^{-1} ( \eta_2^2 - \langle \eta_2^2 \rangle_\tau ) \langle \partial_{\xi_2} \Pi_2 \Delta_\xi^{-1} \partial_{\xi_2} \Pi_2 \rangle_\xi \Delta_\xi^{-1} \partial_{\xi_2} \Pi_2 \partial_{111} \bar \rho_0       \nonumber  \\ 
      =&\, \frac{\delta^{\sfrac12}\mr{\bar\lambda}}{\mu} \bigg( \frac{\delta^{\sfrac12}\mr{\bar\lambda}}{\kappa\lambda} \bigg)^2 \eta_1 \partial_\tau^{-1} ( \eta_1^2 - \langle \eta_1^2 \rangle_\tau ) \langle |\Pi_1|^2 \rangle_\xi \Delta_\xi^{-1} \partial_{\xi_1} \Pi_1 \mr {\bar D}_{222} \bar \rho_0       \label{te:CompHRho2:2} \\ 
      -&\, \frac{\delta^{\sfrac12}\mr{\bar\lambda}}{\mu} \bigg( \frac{\delta^{\sfrac12}\mr{\bar\lambda}}{\kappa\lambda} \bigg)^2 \eta_2 \partial_\tau^{-1} ( \eta_2^2 - \langle \eta_2^2 \rangle_\tau ) \langle |\Pi_2|^2 \rangle_\xi \Delta_\xi^{-1} \partial_{\xi_2} \Pi_2 \mr {\bar D}_{111} \bar \rho_0      \label{te:CompHRho2:4}
\end{align}
\begin{align}
      L_2^{-1} ( \langle L_{11} \hat\rho_1 \rangle_{\xi} - L_{11} \hat\rho_1 )
      =&\, \frac{\delta}{\kappa^2\lambda^2} \eta_1^2 \Delta_\xi^{-1} ( \partial_{\xi_1} \Pi_1 \Delta_\xi^{-1} \partial_{\xi_1} \Pi_1 - \langle \rangle_\xi ) \partial_{22} \bar\rho_0     \nonumber \\ 
      &\, + \frac{\delta}{\kappa^2\lambda^2} \eta_2^2 \Delta_\xi^{-1} ( \partial_{\xi_2} \Pi_2 \Delta_\xi^{-1} \partial_{\xi_2} \Pi_2 - \langle \rangle_\xi ) \partial_{11} \bar\rho_0      \nonumber \\ 
      =&\, \bigg( \frac{\delta^{\sfrac12}\mr{\bar\lambda}}{\kappa\lambda} \bigg)^2 \eta_1^2 \Delta_\xi^{-1} ( \partial_{\xi_1} \Pi_1 \Delta_\xi^{-1} \partial_{\xi_1} \Pi_1 - \langle \rangle_\xi ) \mr {\bar D}_{22} \bar\rho_0    \label{te:CompHRho2:6} \\ 
      +&\, \bigg( \frac{\delta^{\sfrac12}\mr{\bar\lambda}}{\kappa\lambda} \bigg)^2 \eta_2^2 \Delta_\xi^{-1} ( \partial_{\xi_2} \Pi_2 \Delta_\xi^{-1} \partial_{\xi_2} \Pi_2 - \langle \rangle_\xi ) \mr {\bar D}_{11} \bar\rho_0    \label{te:CompHRho2:8}    
\end{align}
In above computation, we omit the formula in $\langle \rangle_\xi$ making the expression zero $\xi$-average. We shall also use this convention below.
\begin{align}
      L_2^{-1} ( -L_{12} \hat\rho_1 ) 
      =&\,- \frac{\mu} {\delta^{\sfrac12}} \frac{\delta}{\kappa^2\lambda^2} \partial_\tau \eta_1 \Delta_\xi^{-1} \partial_{\xi_1} \Pi_1 \partial_{x_2} \bar\rho_0       \nonumber  \\ 
      &\,+ \frac{\mu} {\delta^{\sfrac12}} \frac{\delta}{\kappa^2\lambda^2} \partial_\tau \eta_2 \Delta_\xi^{-1} \partial_{\xi_2} \Pi_2 \partial_{x_1} \bar\rho_0       \nonumber  \\ 
      =&\,- \frac{\mu} {\delta^{\sfrac12}\mr{\bar\lambda}} \bigg( \frac{\delta^{\sfrac12}\mr{\bar\lambda}}{\kappa\lambda} \bigg)^2 \partial_\tau \eta_1 \Delta_\xi^{-1} \partial_{\xi_1} \Pi_1 \mr {\bar D}_{2} \bar\rho_0     \label{te:CompHRho2:10}  \\ 
      &\,+ \frac{\mu} {\delta^{\sfrac12}\mr{\bar\lambda}} \bigg( \frac{\delta^{\sfrac12}\mr{\bar\lambda}}{\kappa\lambda} \bigg)^2 \partial_\tau \eta_2 \Delta_\xi^{-1} \partial_{\xi_2} \Pi_2 \mr {\bar D}_{1} \bar\rho_0     \label{te:CompHRho2:12}
\end{align}
\begin{align}
      L_2^{-1} ( -L_{13} \hat\rho_1 ) = 
      -&\, \varepsilon\lambda \frac{\delta^{\sfrac12}}{\kappa\lambda} \sum_m S_{m,11} \vartheta_m \eta_1 \Delta_\xi^{-1} \partial_{\xi_1} \Pi_1 \partial_{x_2} \bar\rho_0     \nonumber \\ 
      +&\, \varepsilon\lambda \frac{\delta^{\sfrac12}}{\kappa\lambda} \sum_m S_{m,22} \vartheta_m \eta_2 \Delta_\xi^{-1} \partial_{\xi_2} \Pi_2 \partial_{x_1} \bar\rho_0     \nonumber \\ 
      = -&\, \frac{\varepsilon\kappa\lambda^2}{\delta^{\sfrac12} \mr{\bar\lambda}} \bigg( \frac{\delta^{\sfrac12}\mr{\bar\lambda}}{\kappa\lambda} \bigg)^2 \sum_m S_{m,11} \vartheta_m \eta_1 \Delta_\xi^{-1} \partial_{\xi_1} \Pi_1 \mr {\bar D}_{2} \bar\rho_0        \label{te:CompHRho2:14} \\ 
      +&\, \frac{\varepsilon\kappa\lambda^2}{\delta^{\sfrac12} \mr{\bar\lambda}} \bigg( \frac{\delta^{\sfrac12}\mr{\bar\lambda}}{\kappa\lambda} \bigg)^2 \sum_m S_{m,22} \vartheta_m \eta_2 \Delta_\xi^{-1} \partial_{\xi_2} \Pi_2 \mr {\bar D}_{1} \bar\rho_0    \label{te:CompHRho2:16}
\end{align}
\begin{align}
      L_2^{-1} ( - L_{02} \bar\rho_0 ) = 
      &\,- \varepsilon \frac{\delta^{\sfrac12}}{\kappa\lambda} \sum_m \sigma_m \partial_{x_i} \Omega_{m,ij} \big( \eta_1 \Delta_\xi^{-1} \Pi_1 + \eta_2 \Delta_\xi^{-1} \Pi_2 \big) \partial_{x_j} \bar\rho_0       \nonumber \\ 
      =&\, - \frac{\varepsilon\kappa\lambda}{\delta^{\sfrac12}} \bigg( \frac{\delta^{\sfrac12}\mr{\bar\lambda}}{\kappa\lambda} \bigg)^2 
            \sum_m \mr {\bar D}_i \Omega_{m,i1} \sigma_m \eta_1 \Delta_\xi^{-1} \Pi_1 \mr {\bar D}_{1} \bar\rho_0       \label{te:CompHRho2:18}  \\ 
      &\, - \frac{\varepsilon\kappa\lambda}{\delta^{\sfrac12}} \bigg( \frac{\delta^{\sfrac12}\mr{\bar\lambda}}{\kappa\lambda} \bigg)^2
            \sum_m \mr {\bar D}_i \Omega_{m,i2} \sigma_m \eta_1 \Delta_\xi^{-1} \Pi_1 \mr {\bar D}_{2} \bar\rho_0       \label{te:CompHRho2:18_1}  \\ 
      &\, - \frac{\varepsilon\kappa\lambda}{\delta^{\sfrac12}} \bigg( \frac{\delta^{\sfrac12}\mr{\bar\lambda}}{\kappa\lambda} \bigg)^2
            \sum_m \mr {\bar D}_i \Omega_{m,i1} \sigma_m \eta_2 \Delta_\xi^{-1} \Pi_2 \mr {\bar D}_{1} \bar\rho_0       \label{te:CompHRho2:19} \\ 
      &\, - \frac{\varepsilon\kappa\lambda}{\delta^{\sfrac12}} \bigg( \frac{\delta^{\sfrac12}\mr{\bar\lambda}}{\kappa\lambda} \bigg)^2
            \sum_m \mr {\bar D}_i \Omega_{m,i2} \sigma_m \eta_2 \Delta_\xi^{-1} \Pi_2 \mr {\bar D}_{2} \bar\rho_0       \label{te:CompHRho2:19_1}
\end{align}
\begin{align}
      L_2^{-1} ( - L_{03} \bar\rho_0 ) = 
      -&\, \frac{\varepsilon\kappa\lambda^2}{\delta^{\sfrac12} \mr{\bar\lambda}} \bigg( \frac{\delta^{\sfrac12}\mr{\bar\lambda}}{\kappa\lambda} \bigg)^2
            \sum_m \varphi_m \eta_1 \partial_{\xi_1} \Pi_1 E_{m,12} \mr {\bar D}_{1} \bar\rho_0       \label{te:CompHRho2:20} \\ 
      -&\, \frac{\varepsilon\kappa\lambda^2}{\delta^{\sfrac12} \mr{\bar\lambda}} \bigg( \frac{\delta^{\sfrac12}\mr{\bar\lambda}}{\kappa\lambda} \bigg)^2
            \sum_m \varphi_m \eta_1 \partial_{\xi_1} \Pi_1 E_{m,22} \mr {\bar D}_{2} \bar\rho_0       \label{te:CompHRho2:21} \\ 
      +&\, \frac{\varepsilon\kappa\lambda^2}{\delta^{\sfrac12} \mr{\bar\lambda}} \bigg( \frac{\delta^{\sfrac12}\mr{\bar\lambda}}{\kappa\lambda} \bigg)^2
            \sum_m \varphi_m \eta_2 \partial_{\xi_2} \Pi_2 E_{m,11} \mr {\bar D}_{1} \bar\rho_0       \label{te:CompHRho2:22} \\
      +&\, \frac{\varepsilon\kappa\lambda^2}{\delta^{\sfrac12} \mr{\bar\lambda}} \bigg( \frac{\delta^{\sfrac12}\mr{\bar\lambda}}{\kappa\lambda} \bigg)^2
            \sum_m \varphi_m \eta_2 \partial_{\xi_2} \Pi_2 E_{m,21} \mr {\bar D}_{2} \bar\rho_0.       \label{te:CompHRho2:23}
\end{align}

The terms in (\ref{te:CompHRho2:2}-\ref{te:CompHRho2:23}) are from the line \eqref{e:SolAsptInit:33}. Now we need to renumber all these terms and write them in the form of \eqref{e:SolAsptInit:29}. 

Recall $a=0$. Consider the case $\balpha = 2$. Then the contributions from \eqref{te:CompHRho2:10}, \eqref{te:CompHRho2:14}, \eqref{te:CompHRho2:18_1}, \eqref{te:CompHRho2:19_1}, \eqref{te:CompHRho2:21} and \eqref{te:CompHRho2:23} give
\begin{align}     \label{e:SolAsptInit:40}
      \hat I_2^{(2)} = 5N+1
\end{align}
correctors. Taking the contributions in \eqref{te:CompHRho2:10}, \eqref{te:CompHRho2:14} and \eqref{te:CompHRho2:18_1} as an example, we renumber these $2N+1$ correctors as follows. For \eqref{te:CompHRho2:10}, we set
\begin{align*}
      \hat h_2^{(2,1)} = - \frac{\mu} {\delta^{\sfrac12}\mr{\bar\lambda}} \bigg( \frac{\delta^{\sfrac12}\mr{\bar\lambda}}{\kappa\lambda} \bigg)^2,    \quad
      \hat \eta_2^{(2,1)} = \partial_\tau \eta_1,   \quad 
      \hat \chi_2^{(2,1)} = \Delta_\xi^{-1} \partial_{\xi_1} \Pi_1.
\end{align*}
For \eqref{te:CompHRho2:14}, we set, for $2 \leq l \leq N+1$,
\begin{align*}
      \hat h_2^{(2,l)} = - \frac{\varepsilon\kappa\lambda^2}{\delta^{\sfrac12} \mr{\bar\lambda}} \bigg( \frac{\delta^{\sfrac12}\mr{\bar\lambda}}{\kappa\lambda} \bigg)^2 S_{l-1,11},      \quad
      \hat \eta_2^{(2,l)} = \vartheta_{l-1} \eta_1,      \quad 
      \hat \chi_2^{(2,l)} = \Delta_\xi^{-1} \partial_{\xi_1} \Pi_1.
\end{align*}
For \eqref{te:CompHRho2:18_1}, we set, for $N+2 \leq l \leq 2N+1$,
\begin{align*}
      \hat h_2^{(2,l)} = - \varepsilon \frac{\delta^{\sfrac12} \mr{\bar\lambda}^2} {\kappa\lambda} \mr {\bar D}_i \Omega_{l-N-1,i2}, \quad
      \hat \eta_2^{(2,l)} = \sigma_{l-N-1} \eta_1, \quad 
      \hat \chi_2^{(2,l)} = \Delta_\xi^{-1} \Pi_1.
\end{align*}
Other terms are renumbered analogously. Here, it is also clear that $\big( \hat \chi_2^{(\cdot,\cdot)}, \hat \eta_2^{(\cdot,\cdot)} \big)$ satisfies the shear condition in Definition \ref{d:shear_structure}.

For the correctors with $\balpha = 2$, we also estimate
\begin{align}     \label{e:SolAsptInit:46}
      \max_l \vertiii{ \hat h_2^{(2,l)} } \leq \bigg( \frac{\delta^{\sfrac12}\mr{\bar\lambda}}{\kappa\lambda} \bigg)^2 \max \Bigg\{ 
            \frac{\mu} {\delta^{\sfrac12}\mr{\bar\lambda}}, 
            \frac{\varepsilon\kappa\lambda^2} {\delta^{\sfrac12} \mr{\bar\lambda}}
            \Bigg\}.
\end{align}
Here we use Remark \ref{r:opPolynomial} and \eqref{e:mtrEstHom}. Then \eqref{e:EXP1:hh_Pest} follows from \eqref{d:lambda_gamma} and \eqref{e:SolAsptInit:46}. The cardinality estimate \eqref{e:EXP1:hI_est} with $k=1, \balpha = 2$ follows from \eqref{e:SolAsptInit:40}. The property \eqref{e:EXP1:hh} with $n=2$ follows from Remark \ref{r:opPolynomial} and \eqref{e:mtrStrHom}. The estimates \eqref{e:EXP1:hchi_est} and \eqref{e:EXP1:heta_est} follow from Lemma \ref{l:cutoff} and \eqref{e:prdEstHom}. Here, since we have $\partial_\tau \eta_1$ in \eqref{te:CompHRho2:10} and \eqref{te:CompHRho2:12}, the bound on derivative order $p$ in \eqref{e:EXP1:heta_est} may decrease.

The structures of the correctors with $\balpha = 1$ are completely symmetric. The structures of the correctors with $[\balpha] = 2$ or $[\balpha] = 3$ are similar, and the number of correctors is $1$ or $0$, i.e.
\begin{align*}
      \hat I_2^{(\balpha)} =&\, 1, \text{ for }
            \balpha = 11, 22, 111, 222, \\ 
      \hat I_2^{(\balpha)} =&\, 0, \text{ for other }
            \balpha \text{ with } [\balpha] = 2,3.
\end{align*}
The estimates \eqref{e:EXP1:hh}, \eqref{e:EXP1:hchi_est} and \eqref{e:EXP1:heta_est} follow similarly as above. After we have \eqref{e:EXP1:hh_Pest} for all $\balpha$ and $l$, \eqref{e:EXP1:hhS_Pest} follows.

This concludes all the information for the terms in \eqref{e:SolAsptInit:14}.

\subsubsection{Solve the asymptotic system: initial steps for Lemma \ref{l:EXP2}}   \label{sss:solInitExp2} In this section, we solve \eqref{e:asptFF:hrc1}, \eqref{e:asptFF:hrc2}, \eqref{e:asptFF:hrc3} with initial conditions in \eqref{e:0:initialCFF} and \eqref{e:2:initialCFF}, $n=1$. The analysis is largely similar to Section \ref{sss:solInitExp2}. We shall only present the different part.

Similar to the analogue of Lemma \ref{l:EXP1}, we first deduce $\hat\rho_0 = \check\rho_0 = 0$ from \eqref{e:asptFF:hrc1} and \eqref{e:asptFF:hrc2}, verifying \eqref{e:EXP2:hrho0} and \eqref{e:EXP2:crho0}. Again from \eqref{e:asptFF:hrc2}, we have
\begin{align}     \label{e:6:solInitExp2}
      \hat \rho_1 =&\, \frac{ \delta^{\sfrac12} \mr{\bar\lambda} }{\kappa\lambda} \eta_1 \hat \chi^{(2,1)}_1 \mr {\bar D}_{2} \bar \rho_0 - \frac{ \delta^{\sfrac12} \mr{\bar\lambda} }{\kappa\lambda} \eta_2 \hat \chi^{(1,1)}_1 \mr {\bar D}_{1} \bar \rho_0
\end{align}
with the structural information given in \eqref{e:SolAsptInit:17} and \eqref{e:SolAsptInit:18} of Section \ref{sss:solInitExp1}. Setting all forcing related correctors to zero for $n=1$, this verifies (\ref{e:EXP2:hrhon}-\ref{e:EXP2:heta_est}) for $n=1$.

Now we consider \eqref{e:asptFF:hrc3}. Taking $\langle \rangle_{\xi,\tau}$ of \eqref{e:asptFF:hrc3} gives
\begin{align}     \label{e:12:solInitExp2}
      \langle L_{11} \hat \rho_1 \rangle_{\xi,\tau} +&\, L_{01} \bar \rho_0 
      = \frac{ \mr{\bar\lambda} } { \mr\lambda }
            \sum_{k=1}^K \sum_{ [\balpha] = 1, \balpha \in \mathcal{I}_* }^Q \sum_{l = 1}^{J} 
            \mr {\bar\divr}_x \Big( h_{F,k}^{(\balpha,l)} \mr {\bar D}^{\balpha} \bar \rho_{F,k} \Big) 
            \langle \eta_{F,k}^{(\balpha,l)} \chi_{F,k}^{(\balpha,l)} \rangle_{\xi,\tau}
\end{align}
Here, we use the formula for $G$ in \eqref{e:6:compFHom} and \eqref{e:8:compFHom}. Now, combining \eqref{e:6:solInitExp2}, $\hat\rho_0 = \check\rho_0 = 0$ and \eqref{e:0:initialCFF}, we can verify \eqref{e:EXP2:brho0} and \eqref{e:EXP2:brho0_ini}.

\underline{For the structure of $\check \rho_1$}, we take the difference of $\langle \eqref{e:asptFF:hrc3} \rangle_{\xi,\tau}$ and $\langle \eqref{e:asptFF:hrc3} \rangle_{\xi}$. We deduce 
\begin{align*}
      L_{12} \check \rho_1 =&\, \frac{ \mr{\bar\lambda} } { \mr\lambda } 
            \sum_{k=1}^K \sum_{ [\balpha] = 1, \balpha \in \mathcal{I}_* }^Q \sum_{l = 1}^{J} 
            \mr {\bar\divr}_x \Big( h_{F,k}^{(\balpha,l)} \mr {\bar D}^{\balpha} \bar \rho_{F,k} \Big) \big( \eta_{F,k}^{(\balpha,l)} - \langle \eta_{F,k}^{(\balpha,l)} \rangle_\tau \big) \langle \chi_{F,k}^{(\balpha,l)} \rangle_{\xi} \\
            &\, +\langle L_{11} \hat \rho_1 \rangle_{\xi,\tau}
            - \langle L_{11} \hat \rho_1 \rangle_\xi
\end{align*}
and thus
\begin{align}
      \check \rho_1 =&\, - \frac{\delta^{\sfrac12} \mr{\bar\lambda}}{\mu} \frac{\delta^{\sfrac12} \mr{\bar\lambda}}{\kappa\lambda} \partial_\tau^{-1} \big( \eta_1^2 - \langle \eta_1^2 \rangle_\tau \big) \langle \Pi_1 \partial_{\xi_1} \hat \chi_1^{(2,1)} \rangle_\xi \mr {\bar D}_{22} \bar \rho_0       \label{e:13:solInitExp2} \\ 
      &\, - \frac{\delta^{\sfrac12} \mr{\bar\lambda}}{\mu} \frac{\delta^{\sfrac12} \mr{\bar\lambda}}{\kappa\lambda} \partial_\tau^{-1} \big( \eta_2^2 - \langle \eta_2^2 \rangle_\tau \big) \langle \Pi_2 \partial_{\xi_2} \hat \chi_1^{(1,1)} \rangle_\xi \mr {\bar D}_{11} \bar \rho_0    \label{e:14:solInitExp2} \\ 
      &\, + \frac{\lambda}{\mu} \frac{ \mr{\bar\lambda} } { \mr\lambda } 
            \sum_{k=1}^K \sum_{[\balpha] = 1, \balpha \in \mathcal{I}_*}^Q \sum_{l = 1}^{J} 
      h_{F,k,j}^{(\balpha,l)} \partial_{\tau}^{-1} \big( \eta_{F,k}^{(\balpha,l)} - \langle \eta_{F,k}^{(\balpha,l)} \rangle_\tau \big) \langle \chi_{F,k}^{(\balpha,l)} \rangle_{\xi} \mr {\bar\partial}_j \mr {\bar D}^{\balpha} \bar \rho_{F,k}      \label{e:15:solInitExp2} \\ 
      &\, + \frac{\lambda}{\mu} \frac{ \mr{\bar\lambda} } { \mr\lambda } 
            \sum_{k=1}^K \sum_{[\balpha] = 1, \balpha \in \mathcal{I}_*}^Q \sum_{l = 1}^{J} 
      \mr {\bar\divr} h_{F,k}^{(\balpha,l)} \partial_{\tau}^{-1} \big( \eta_{F,k}^{(\balpha,l)} - \langle \eta_{F,k}^{(\balpha,l)} \rangle_\tau \big) \langle \chi_{F,k}^{(\balpha,l)} \rangle_{\xi} \mr {\bar D}^{\balpha} \bar \rho_{F,k}.       \label{e:16:solInitExp2}
\end{align}

Now we write above expressions for $\check \rho_1$ in the form of \eqref{e:EXP2:crhon}. The contributions in (\ref{e:13:solInitExp2}-\ref{e:14:solInitExp2}) go to \eqref{e:EXP2:crhon}. The functions $\check h_{1}^{(\balpha,l)}$ and $\check \eta_{1}^{(\balpha,l)}$ are the same as Lemma \ref{l:EXP1} and are given in Section \ref{sss:solInitExp1}. For the $F$-terms in (\ref{e:15:solInitExp2}-\ref{e:16:solInitExp2}), we have for $\balpha$ with $[\balpha] = 1$ and $\balpha \in \mathcal{I}_*$,
\begin{align}     \label{e:18:solInitExp2}
      \check I_{F,1}^{(\balpha)} = \check I_{F,1}^{(1\balpha)} = \check I_{F,1}^{(2\balpha)} = J,
      \quad       \check I_{F,1}^{(t\balpha)} = 0.
\end{align}
And $\check h_{F,k,1}^{(\bbeta,l)}$, $\check \eta_{F,k,1}^{(\bbeta,l)}$ can be written as 
\begin{align} 
      \check h_{F,k,1}^{(\bbeta,l)} = \frac{\lambda}{\mu} \frac{ \mr{\bar\lambda} } { \mr\lambda }
            \mr {\bar\divr}_x h_{F,k}^{(\balpha,l)},     \quad 
      \check \eta_{F,k,1}^{(\bbeta,l)} = \partial_{\tau}^{-1} \big( \eta_{F,k}^{(\balpha,l)} - \langle \eta_{F,k}^{(\balpha,l)} \rangle_\tau \big) \langle \chi_{F,k}^{(\balpha,l)} \rangle_{\xi} 
            \quad \text{for } \bbeta = \balpha,       \label{e:19:solInitExp2} \\ 
      \check h_{F,k,1}^{(\bbeta,l)} = \frac{\lambda}{\mu} \frac{ \mr{\bar\lambda} } { \mr\lambda }
            h_{F,k,1}^{(\balpha,l)},     \quad 
      \check \eta_{F,k,1}^{(\bbeta,l)} = \partial_{\tau}^{-1} \big( \eta_{F,k}^{(\balpha,l)} - \langle \eta_{F,k}^{(\balpha,l)} \rangle_\tau \big) \langle \chi_{F,k}^{(\balpha,l)} \rangle_{\xi} 
            \quad \text{for } \bbeta = 1\balpha,      \label{e:20:solInitExp2} \\ 
      \check h_{F,k,1}^{(\bbeta,l)} = \frac{\lambda}{\mu} \frac{ \mr{\bar\lambda} } { \mr\lambda }
            h_{F,k,2}^{(\balpha,l)},     \quad 
      \check \eta_{F,k,1}^{(\bbeta,l)} = \partial_{\tau}^{-1} \big( \eta_{F,k}^{(\balpha,l)} - \langle \eta_{F,k}^{(\balpha,l)} \rangle_\tau \big) \langle \chi_{F,k}^{(\balpha,l)} \rangle_{\xi} 
            \quad \text{for } \bbeta = 2\balpha.      \label{e:21:solInitExp2}
\end{align}
Note that $h_{F,k}^{(\balpha,l)}$ is a vector-valued function and $h_{F,k,j}^{(\balpha,l)}$ is its $j$-th component. The information \eqref{e:EXP2:ch} for $n=1$ and (\ref{e:EXP2:cIF_est}-\ref{e:EXP2:ceta_est}) for $n=1$ follow from Assumption \ref{a:homSmallscl}.

\underline{For the structure of $\hat \rho_2$}, we take the differene of $\langle \eqref{e:asptFF:hrc3} \rangle_{\xi}$ and $\eqref{e:asptFF:hrc3}$     
\begin{align}
      L_2 \hat \rho_2 =&\, \langle L_1 \rho_1 \rangle_\xi - L_1 \rho_1 + \langle L_0 \rho_0 \rangle_\xi - L_0 \rho_0    \label{e:22:2:solInitExp2} \\ 
      +&\, \frac{ \mr{\bar\lambda} } { \mr\lambda }
            \sum_{k=1}^K \sum_{[\balpha] = 1, \balpha \in \mathcal{I}_*}^Q \sum_{l = 1}^{J} 
            h_{F,j}^{(\balpha,l)} \eta_{F,k}^{(\balpha,l)} \big( \chi_{F,k}^{(\balpha,l)} - \langle \chi_{F,k}^{(\balpha,l)} \rangle_{\xi} \big) \mr {\bar\partial}_j \mr {\bar D}^{\balpha} \bar \rho_{F,k}    \label{e:22:4:solInitExp2} \\ 
      +&\, \frac{ \mr{\bar\lambda} } { \mr\lambda }
            \sum_{k=1}^K \sum_{[\balpha] = 1, \balpha \in \mathcal{I}_*}^Q \sum_{l = 1}^{J} 
            \mr {\bar\divr}_x h_{F,k}^{(\balpha,l)} \eta_{F,k}^{(\balpha,l)} \big( \chi_{F,k}^{(\balpha,l)} - \langle \chi_{F,k}^{(\balpha,l)} \rangle_{\xi} \big) \mr {\bar D}^{\balpha} \bar \rho_{F,k}    \label{e:22:6:solInitExp2} \\ 
      +&\, \frac{\lambda} {\mr\lambda} 
            \sum_{k=1}^K \sum_{ [\balpha] = 1, \balpha \in \mathcal{I}_* }^Q \sum_{l = 1}^{J} 
            h_{F,k,i}^{(\balpha,l)} \eta_{F,k}^{(\balpha,l)} \partial_{\xi_i} \chi_{F,k}^{(\balpha,l)} \mr {\bar D}^{\balpha} \bar \rho_{F,k}      \label{e:22:8:solInitExp2} \\ 
      +&\, \frac{\varepsilon \lambda} {\mr\lambda} 
            \sum_{k=1}^K \sum_{ [\balpha] = 1, \balpha \in \mathcal{I}_* }^Q \sum_{l = 1}^{J} \sum_m 
            B_{m,ij} h_{F,k,j}^{(\balpha,l)} \phi_m \eta_{F,k}^{(\balpha,l)} \partial_{\xi_i} \chi_{F,k}^{(\balpha,l)} \mr {\bar D}^{\balpha} \bar \rho_{F,k},       \label{e:22:10:solInitExp2}
\end{align}
and thus
\begin{align}
      \hat \rho_2 =&\, 
            \frac{ \delta^{\sfrac12} \mr{\bar\lambda} }{\kappa\lambda} \eta_1 \hat\chi^{(2,1)}_1 \mr {\bar D}_{2} \bar \rho_1 
            - \frac{ \delta^{\sfrac12} \mr{\bar\lambda} }{\kappa\lambda} \eta_2 \hat\chi^{(1,1)}_1 \mr {\bar D}_{1} \bar \rho_1      \label{e:24:solInitExp2} \\ 
      -&\, \frac{\delta^{\sfrac12}\mr{\bar\lambda}}{\kappa\lambda} \frac{\lambda}{\mu} \frac{ \mr{\bar\lambda} } { \mr\lambda }
            \sum_{k=1}^K \sum_{[\balpha] = 1, \balpha \in \mathcal{I}_*}^{Q+1} \sum_{l = 1}^{\check I_{F,1}^{(\balpha)}}
            \check h_{F,1}^{(\balpha,l)} \eta_1 \check \eta_{F,1}^{(\balpha,l)} \Delta_\xi^{-1} \partial_{\xi_1} \Pi_1 \mr {\bar D}_{2} \mr {\bar D}^{\balpha} \bar\rho_F     \label{e:25:solInitExp2} \\ 
      +&\, \frac{\delta^{\sfrac12}\mr{\bar\lambda}}{\kappa\lambda} \frac{\lambda}{\mu} \frac{ \mr{\bar\lambda} } { \mr\lambda }
            \sum_{k=1}^K \sum_{[\balpha] = 1, \balpha \in \mathcal{I}_*}^{Q+1} \sum_{l = 1}^{\check I_{F,1}^{(\balpha)}}
            \check h_{F,1}^{(\balpha,l)} \eta_2 \check \eta_{F,1}^{(\balpha,l)} \Delta_\xi^{-1} \partial_{\xi_2} \Pi_2 \mr {\bar D}_{1} \mr {\bar D}^{\balpha} \bar\rho_F     \label{e:26:solInitExp2} \\ 
      -&\, \frac{\delta^{\sfrac12}\mr{\bar\lambda}}{\kappa\lambda} \frac{\lambda}{\mu} \frac{ \mr{\bar\lambda} } { \mr\lambda } 
            \sum_{k=1}^K \sum_{[\balpha] = 1, \balpha \in \mathcal{I}_*}^{Q+1} \sum_{l = 1}^{\check I_{F,1}^{(\balpha)}}
            \mr {\bar D}_{2} \check h_{F,1}^{(\balpha,l)} \eta_1 \check \eta_{F,1}^{(\balpha,l)} \Delta_\xi^{-1} \partial_{\xi_1} \Pi_1 \mr {\bar D}^{\balpha} \bar\rho_F     \label{e:27:solInitExp2} \\ 
      +&\, \frac{\delta^{\sfrac12}\mr{\bar\lambda}}{\kappa\lambda} \frac{\lambda}{\mu} \frac{ \mr{\bar\lambda} } { \mr\lambda } 
            \sum_{k=1}^K \sum_{[\balpha] = 1, \balpha \in \mathcal{I}_*}^{Q+1} \sum_{l = 1}^{\check I_{F,1}^{(\balpha)}}
            \mr {\bar D}_{1} \check h_{F,1}^{(\balpha,l)} \eta_2 \check \eta_{F,1}^{(\balpha,l)} \Delta_\xi^{-1} \partial_{\xi_2} \Pi_2 \mr {\bar D}^{\balpha} \bar\rho_F     \label{e:28:solInitExp2} \\ 
      -&\, \frac{1}{\kappa} \frac{ \mr{\bar\lambda} } { \mr\lambda }
            \sum_{k=1}^K \sum_{[\balpha] = 1, \balpha \in \mathcal{I}_* }^Q \sum_{l = 1}^{J}
            h_{F,j}^{(\balpha,l)} \eta_{F,k}^{(\balpha,l)} \Delta_\xi^{-1} ( \chi_{F,k}^{(\balpha,l)} - \langle \chi_{F,k}^{(\balpha,l)} \rangle_{\xi} ) \mr {\bar\partial}_j \mr {\bar D}^{\balpha} \bar\rho_F   \label{e:29:solInitExp2} \\ 
      -&\, \frac{1}{\kappa} \frac{ \mr{\bar\lambda} } { \mr\lambda }
            \sum_{k=1}^K \sum_{[\balpha] = 1, \balpha \in \mathcal{I}_* }^Q \sum_{l = 1}^{J} 
            \mr {\bar\divr}_x h_{F}^{(\balpha,l)} \eta_{F,k}^{(\balpha,l)} \Delta_\xi^{-1} \big( \chi_{F,k}^{(\balpha,l)} - \langle \chi_{F,k}^{(\balpha,l)} \rangle_{\xi} \big) \mr {\bar D}^{\balpha} \bar\rho_F   \label{e:30:solInitExp2} \\ 
      -&\, \frac{\lambda} {\kappa \mr\lambda} 
            \sum_{k=1}^K \sum_{ [\balpha] = 1, \balpha \in \mathcal{I}_* }^Q \sum_{l = 1}^{J} 
            h_{F,i}^{(\balpha,l)} \eta_{F,k}^{(\balpha,l)} \Delta_\xi^{-1} \partial_{\xi_i} \chi_{F,k}^{(\balpha,l)} \mr {\bar D}^{\balpha} \bar \rho_{F,k}      \label{e:31:solInitExp2} \\ 
      -&\, \frac{\varepsilon \lambda} {\kappa \mr\lambda} 
            \sum_{k=1}^K \sum_{ [\balpha] = 1, \balpha \in \mathcal{I}_* }^Q \sum_{l = 1}^{J} \sum_m 
            B_{m,ij} h_{F,j}^{(\balpha,l)} \phi_m \eta_{F,k}^{(\balpha,l)} \Delta_\xi^{-1} \partial_{\xi_i} \chi_{F,k}^{(\balpha,l)} \mr {\bar D}^{\balpha} \bar \rho_{F,k}       \label{e:32:solInitExp2} \\ 
      +&\, \text{ terms presented in } (\ref{te:CompHRho2:2}-\ref{te:CompHRho2:23}).       \label{e:33:solInitExp2}
\end{align}

The terms in \eqref{e:33:solInitExp2} has been analyzed in Section \eqref{sss:solInitExp1}. This also verifies Remark \ref{r:shareCorrector} for such correctors. Then we renumber the terms from (\ref{e:25:solInitExp2}-\ref{e:32:solInitExp2}), and we have the following estimates on the cardinality 
\begin{align}
      \hat I^{(\balpha)}_{F,2} \leq&\, 2N J,   \quad  [\balpha] \leq Q, \\ 
      \hat I^{(\balpha)}_{F,2} \leq&\, 5 J,     \quad  [\balpha] \geq Q+1.
\end{align}
We recall that $\big( \chi_{F,k}^{(\balpha,l)}, \eta_{F,k}^{(\balpha,l)} \big)$ satisfies generalized shear condition for any $\balpha$ and any $l$, then $\big( \hat\chi_{F,k,2}^{(\balpha,l)}, \hat\eta_{F,k,2}^{(\balpha,l)} \big)$ satisfies shear condition for any $\balpha$ and any $l$. From Assumption \ref{a:homSmallscl} and \eqref{d:bigMu}, we can deduce \eqref{e:EXP2:hhF} for $n=2$, (\ref{e:EXP2:hhF_Pest}-\ref{e:EXP2:heta_est}) for $n=2$. This is similar to Section \ref{sss:solInitExp1}.

This concludes deriving the structural information on
\begin{align}     \label{e:40:solInitExp2}
      \hat \rho_0, \hat \rho_1, \hat \rho_2,
            \check \rho_0, \check \rho_1, \bar \rho_0. 
\end{align}

\subsubsection{Induction assumptions and residual correctors}      \label{sss:solAsptRSD}

For the induction step, we solve \eqref{e:aspt:hrc4} with initial conditions in \eqref{e:2:initialC}. Together with the analysis in Section \ref{sss:setUpAsymp} and Section \ref{sss:solInitExp1}, this would concludes the proof of Lemma \ref{l:EXP2}. The proof of Lemma \ref{l:EXP1} is the same, and we comment on the difference.

For Lemma \ref{l:EXP2}, we impose the following inductive assumptions
\begin{equation}        \label{e:1:SolAspt}
\begin{split}
      \hat \rho_{n}(x,t,\xi,\tau) =& \sum_{a=0}^{n-1} \sum_{[\balpha]=1, \balpha \in \mathcal{I}_*}^{2(n-a)-1} \sum_{l=1}^{\hat I_{n-a}^{(\balpha)}}
            \hat h^{(\balpha, l)}_{n-a} (x,t) \hat \eta^{(\balpha, l)}_{n-a} (\tau) \hat \chi^{(\balpha, l)}_{n-a} (\xi) \mr {\bar D}^{\balpha} \bar \rho_a(x,t) \\ 
      +&\, \sum_{k=1}^K \sum_{[\balpha]=1, \balpha \in \mathcal{I}_*}^{2n-2+Q} \sum_{l=1}^{\hat I_{F,n}^{(\balpha)}}
            \hat h^{(\balpha,l)}_{F,k,n} (x,t) \hat \eta^{(\balpha,l)}_{F,k,n} (\tau) \hat \chi^{(\balpha,l)}_{F,k,n} (\xi) \mr {\bar D}^{\balpha} \bar \rho_{F,k}(x,t) \\ 
      \hat \rho_{n+1}(x,t,\xi,\tau) =& \sum_{a=0}^n \sum_{[\balpha] = 1, \balpha \in \mathcal{I}_*}^{2(n-a)+1} \sum_{l=1}^{\hat I_{n-a+1}^{(\balpha)}}
            \hat h^{(\balpha, l)}_{n-a+1} (x,t) \hat \eta^{(\balpha, l)}_{n-a+1} (\tau) \hat \chi^{(\balpha, l)}_{n-a+1} (\xi) \mr {\bar D}^{\balpha} \bar \rho_a(x,t)  \\ 
      +&\, \sum_{k=1}^K \sum_{[\balpha]=1, \balpha \in \mathcal{I}_*}^{2n+Q} \sum_{l=1}^{\hat I_{F,n+1}^{(\balpha)}}
            \hat h^{(\balpha,l)}_{F,k,n+1} (x,t) \hat \eta^{(\balpha,l)}_{F,k,n+1} (\tau) \hat \chi^{(\balpha,l)}_{F,k,n+1} (\xi) \mr {\bar D}^{\balpha} \bar \rho_{F,k}(x,t) \\ 
      \hat \rho_{n+2}(x,t,\xi,\tau) =& \sum_{[\balpha] = 1, \balpha \in \mathcal{I}_*}^{2n+3} \sum_{l=1}^{ \hat I_{n+2}^{(\balpha)} }
            \hat h^{(\balpha, l)}_{n+2} (x,t) \hat \eta^{(\balpha, l)}_{n+2} (\tau) \hat \chi^{(\balpha, l)}_{n+2} (\xi) \mr {\bar D}^{\balpha} \bar \rho_0(x,t) \\  
            +&\, \sum_{a = 1}^{n} \sum_{[\balpha] = 1, \balpha \in \mathcal{I}_*}^{2(n-a)+3} \sum_{l=1}^{ \hat I_{n-a+2}^{(\balpha)} }
            \hat h^{(\balpha, l)}_{n-a+2} (x,t) \hat \eta^{(\balpha, l)}_{n-a+2} (\tau) \hat \chi^{(\balpha, l)}_{n-a+2} (\xi) \mr {\bar D}^{\balpha} \bar \rho_a(x,t) \\ 
            +&\, \sum_{[\balpha] = 1, \balpha \in \mathcal{I}_*} \sum_{l=1}^{ \hat I_{1}^{(\balpha)} }
            \hat h^{(\balpha, l)}_{1} (x,t) \hat \eta^{(\balpha, l)}_{1} (\tau) \hat \chi^{(\balpha, l)}_{1} (\xi) \mr {\bar D}^{\balpha} \bar \rho_{n+1}(x,t) \\ 
            +&\, \sum_{k=1}^K \sum_{[\balpha]=1, \balpha \in \mathcal{I}_*}^{2n+2+Q} \sum_{l=1}^{\hat I^{(\balpha)}_{F,n+2}}
            \hat h^{(\balpha,l)}_{F,k,n+2} (x,t) \hat \eta^{(\balpha,l)}_{F,k,n+2} (\tau) \hat \chi^{(\balpha,l)}_{F,k,n+2} (\xi) \mr {\bar D}^{\balpha} \bar \rho_{F,k}(x,t) \\ 
            :=& \hat \rho_{n+2,0} + \hat \rho_{n+2,M} + \hat \rho_{n+2,n+1} + \hat \rho_{n+2,F} \\ 
            :=& \hat \rho_{n+2,0} + \hat \rho_{n+2,E} + \hat \rho_{n+2,F}
                  := \hat \rho_{n+2,R} + \hat \rho_{n+2,F}.
\end{split}
\end{equation}
and
\begin{equation}        \label{e:2:SolAspt}
\begin{split}
      \check \rho_{n}(x,t,\tau) =& \sum_{a=0}^{n-1} \sum_{[\balpha]=1, \balpha \in \mathcal{I}_*}^{2(n-a)} \sum_{l=1}^{\check I_{n-a}^{(\balpha)}}
            \check h^{(\balpha, l)}_{n-a} (x,t) \check \eta^{(\balpha, l)}_{n-a} (\tau) \mr {\bar D}^{\balpha} \bar \rho_a(x,t) \\ 
      +&\, \sum_{k=1}^K \sum_{[\balpha]=1, \balpha \in \mathcal{I}_*}^{2n-1+Q} \sum_{l=1}^{\check I_{F,n}^{(\balpha)}}  \check h^{(\balpha,l)}_{F,k,n} (x,t) \check \eta^{(\balpha,l)}_{F,k,n} (\tau) \mr {\bar D}^{\balpha} \bar \rho_{F,k} (x,t) 
      := \check\rho_{n,R} + \check\rho_{n,F}, \\ 
      \check \rho_{n+1}(x,t,\tau) =& \sum_{a=0}^n \sum_{[\balpha] = 1, \balpha \in \mathcal{I}_*}^{2(n-a)+2} \sum_{l=1}^{\check I_{n-a+1}^{(\balpha)}}
            \check h^{(\balpha, l)}_{n-a+1} (x,t) \check \eta^{(\balpha, l)}_{n-a+1} (\tau) \mr {\bar D}^{\balpha} \bar \rho_a(x,t) \\ 
      +&\, \sum_{k=1}^K \sum_{[\balpha]=1, \balpha \in \mathcal{I}_*}^{2n+1+Q} \sum_{l=1}^{\check I^{(\balpha)}_{F,n+1}} \check h^{(l)}_{F,k,n+1} (x,t) \check \eta^{(l)}_{F,k,n+1} (\tau) \mr {\bar D}^{\balpha} \bar \rho_{F,k} (x,t) 
      :=  \check\rho_{n+1,R} + \check\rho_{n+1,F}.
\end{split}
\end{equation}

We also assume all relevant estimates in Lemma \ref{l:EXP1} and Lemma \ref{l:EXP2} on 
\begin{align}     \label{e:3:SolAspt}
      \hat \rho_n, \hat \rho_{n+1}, \hat \rho_{n+2}, \check \rho_n, \check \rho_{n+1}.
\end{align}
In above formula, we also assume the stationarity property in Remark \ref{r:stationarity} up to the terms in \eqref{e:3:SolAspt} implicitly. The goal in this induction step is to derive the information in $\bar \rho_{n+1}, \check \rho_{n+2}, \hat \rho_{n+3} $, which are done in Section \ref{sss:solAsptRSD}, Section \ref{sss:solAsptTPR} and Section \ref{sss:solAsptSPT}, respectively.

In \eqref{e:1:SolAspt}, we introduce the decomposition notation for $\hat\rho_{n+2}$, where $\hat \rho_{n+2,0}$, $\hat \rho_{n+2,M}$, $\hat \rho_{n+2,n+1}$ and $\hat \rho_{n+2,F}$ contain, respectively, the terms involving $\bar \rho_{0}$, the terms involving $\bar \rho_{a}$ with $1 \leq a \leq n$, the terms involving $\bar \rho_{n+1}$ and the terms involving $\bar \rho_F$. As defined in \eqref{e:1:SolAspt}, we also denote
\begin{align}     \label{e:4:2:SolAspt}
      \hat\rho_{n+2,R} := \hat \rho_{n+2,0} + \hat \rho_{n+2,E}
      := \hat \rho_{n+2,0} + \hat \rho_{n+2,M} + \hat \rho_{n+2,n+1}.
\end{align}
We apply analogous decomposition notation to $\hat\rho_{n+1}$ and $\hat\rho_{n}$,
\begin{equation}  \label{e:4:6:SolAspt}
\begin{split}
      \hat \rho_{n+1} :=&\, \hat \rho_{n+1,0} + \hat \rho_{n+1,M} + \hat \rho_{n+1,n} + \hat \rho_{n+1,F} \\ 
            :=&\, \hat \rho_{n+1,0} + \hat \rho_{n+1,E} + \hat \rho_{n+1,F}
            := \hat \rho_{n+1,R} + \hat \rho_{n+1,F}, \\ 
      \hat \rho_{n} :=&\, \hat \rho_{n,0} + \hat \rho_{n,M} + \hat \rho_{n,n-1} + \hat \rho_{n,F} \\      
            :=&\, \hat \rho_{n,0} + \hat \rho_{n,E} + \hat \rho_{n,F} 
            := \hat \rho_{n,R} + \hat \rho_{n,F}
\end{split}
\end{equation}
Similar notations have been introduced to $\check \rho_n$ and $\check \rho_{n+1}$ in \eqref{e:2:SolAspt}.

\underline{The goal in this section is to derive the structures of $\bar \rho_{n+1}$}, i.e. computing the coefficients in
\begin{align}
      \bar L_0 \bar \rho_{n+1} =&\, \sum_{a=0}^{n} \sum_{[\bbeta]=1, \bbeta \in \mathcal{I}_*}^{2(n-a)+3} 
            \mr{\bar \divr} \Big( \bar h^{(\bbeta)}_{n+1,a} \mr {\bar D}^{\bbeta} \bar \rho_a \Big) 
      + \sum_{k=1}^K \sum_{[\bbeta]=1, \bbeta \in \mathcal{I}_*}^{2n+2+Q} 
            \mr{\bar \divr} \Big( \bar h_{F,k,n+1}^{(\bbeta)} \mr {\bar D}^{\bbeta} \bar \rho_{F,k} \Big),     \label{e:5:SolAspt} \\ 
      \bar \rho_{n+1} (\cdot, 0) =&\, -\check \rho_{n+1} (\cdot, 0,0),     \label{e:5:0:SolAspt}
\end{align}
The initial conditions for the residual correctors $\bar \rho_{n+1}$ are given by (\ref{e:0:initialC}-\ref{e:2:initialC}) for Lemma \ref{l:EXP1} and (\ref{e:0:initialCFF}-\ref{e:6:initialCFF}) for Lemma \ref{l:EXP2}. Here, we use the fact that $(\hat\chi_{n+1}^{(\cdot,\cdot)}, \hat\eta_{n+1}^{(\cdot,\cdot)})$ and $(\hat\chi_{F,k,n+1}^{(\cdot,\cdot)}, \hat\eta_{F,k,n+1}^{(\cdot,\cdot)})$ satisfy the shear condition, which leads to
\begin{align*}
      \hat \rho_{n+1} (x, 0, \xi, 0) = 0, \quad 
      \forall \, x, \xi.
\end{align*}

We decompose
\begin{align}     \label{e:5:1:SolAspt}
      \bar \rho_{n+1} := \bar \rho_{n+1,R} + \bar \rho_{n+1,F}    
\end{align}
with
\begin{align}
      \bar L_0 \bar \rho_{n+1,R} =&\, \sum_{a=0}^{n} \sum_{[\bbeta]=1, \bbeta \in \mathcal{I}_*}^{2(n-a)+3} \mr{\bar \divr} \Big( \bar h^{(\bbeta)}_{n+1,a} \mr {\bar D}^{\bbeta} \bar \rho_a \Big),        \label{e:5:2:SolAspt} \\ 
      \bar L_0 \bar \rho_{n+1,F} =&\, \sum_{k=1}^K \sum_{[\bbeta]=1, \bbeta \in \mathcal{I}_*}^{2n+2+Q} \mr{\bar \divr} \Big( \bar h_{F,k,n+1}^{(\bbeta)} \mr {\bar D}^{\bbeta} \bar \rho_{F,k} \Big).       \label{e:5:3:SolAspt} \\ 
      \bar \rho_{n+1,R} (\cdot, 0) =&\, -\check \rho_{n+1,R} (\cdot, 0,0),    \label{e:5:4:SolAspt}  \\ 
      \bar \rho_{n+1,F} (\cdot, 0) =&\, -\check \rho_{n+1,F} (\cdot, 0,0).    \label{e:5:5:SolAspt}
\end{align}

\underline{We first prove the stationarity property mentioned in Remark \ref{r:stationarity}}, i.e. 

Claim: We have, for $a \geq 1$ and $\bbeta$ with $1 \leq [\bbeta] \leq 2(n-a)+4$, that
\begin{equation}     \label{e:6:SolAspt}
\begin{split}
      \bar h^{(\bbeta)}_{n+1,a} =&\, \bar h^{(\bbeta)}_{n-a+1}.
\end{split}
\end{equation}

\begin{proof}[Proof of stationarity for $\bar \rho_{n+1}$, i.e. the claim above]
When $a \geq 1$, we have $n \geq 1$. From Remark \ref{r:hatcheckbar} and taking $\langle \rangle_{\xi,\tau}$ of \eqref{e:aspt:hrc4}, we deduce
\begin{equation}        \label{e:7:SolAspt}
\begin{split}
      \langle L_{11} (\hat\rho_{n+2,0}+\hat\rho_{n+2,M}) \rangle_{\xi,\tau} + \langle L_{11} \hat\rho_{n+2,n+1} \rangle_{\xi,\tau} &\,+ \langle L_{11} \hat\rho_{n+2,F} \rangle_{\xi,\tau} + L_{01} \bar\rho_{n+1} \\ 
      + \langle (L_{02} + L_{03} + L_{04}) \hat\rho_{n+1,R} \rangle_{\xi,\tau}
            &\,+ \langle (L_{02} + L_{03} + L_{04}) \hat\rho_{n+1,F} \rangle_{\xi,\tau} = 0.       
\end{split}
\end{equation}
Using the coefficients $\hat h_1^{(\cdot,\cdot)}$, $\hat \eta_1^{(\cdot,\cdot)}$ and $\hat \chi_1^{(\cdot,\cdot)}$ computed in Section \ref{sss:solInitExp1}, we have
\begin{align}     \label{e:7:1:SolAspt}
      \langle L_{11} \hat\rho_{n+2,n+1} \rangle_{\xi,\tau} + L_{01} \bar\rho_{n+1}
            = \bar L_0 \bar\rho_{n+1}.
\end{align}
This is the same as deriving the equation for $\bar \rho_0$ in Section \ref{sss:solInitExp1}.

With the decomposition (\ref{e:5:1:SolAspt}-\ref{e:5:5:SolAspt}) and \eqref{e:7:1:SolAspt}, to solve \eqref{e:7:SolAspt}, it suffices to solve
\begin{align}
      \bar L_0 \bar \rho_{n+1,R} + \langle L_{11} (\hat\rho_{n+2,0}+\hat\rho_{n+2,M}) \rangle_{\xi,\tau} + \langle (L_{02} + L_{03} + L_{04}) \hat\rho_{n+1,R} \rangle_{\xi,\tau} =&\, 0,       \label{e:7:2:SolAspt}   \\ 
      \bar L_0 \bar \rho_{n+1,F} + \langle L_{11} \hat\rho_{n+2,F} \rangle_{\xi,\tau} + \langle (L_{02} + L_{03} + L_{04}) \hat\rho_{n+1,F} \rangle_{\xi,\tau} =&\, 0.        \label{e:7:3:SolAspt}
\end{align}

On the other hand, consider the equations at the last induction step, i.e.
\begin{align}
      L_2 \rho_{n+2} + L_1 \rho_{n+1} + L_0 \rho_{n} =&\, 0, \quad \text{if } n=0,        \label{e:8:SolAspt} \\ 
      L_2 \rho_{n+2} + L_1 \rho_{n+1} + L_0 \rho_{n} + L_{-1} \rho_{n-1} =&\, 0, \quad \text{if } n \geq 1.       \label{e:9:SolAspt}
\end{align}
Recalling \eqref{e:18:hatcheckbar}, taking $\langle \rangle_{\xi,\tau}$ of \eqref{e:8:SolAspt} or \eqref{e:9:SolAspt} gives
\begin{align}
      \langle L_{11} \hat\rho_{n+1} \rangle_{\xi,\tau} + L_{01} \bar\rho_{n} 
      + \langle (L_{02} + L_{03} + L_{04}) \hat\rho_{n} \rangle_{\xi,\tau} =&\, 0.  \label{e:10:SolAspt}
\end{align}

Analogous to obtaining \eqref{e:7:2:SolAspt} and \eqref{e:7:3:SolAspt} from \eqref{e:7:SolAspt}, we obtain, from \eqref{e:10:SolAspt},
\begin{align}
      \bar L_0 \bar \rho_{n,R} + \langle L_{11} (\hat\rho_{n+1,0}+\hat\rho_{n+1,M}) \rangle_{\xi,\tau} + \langle (L_{02} + L_{03} + L_{04}) \hat\rho_{n,R} \rangle_{\xi,\tau} =&\, 0,     \label{e:10:2:SolAspt} \\ 
      \bar L_0 \bar \rho_{n,F} + \langle L_{11} \hat\rho_{n+1,F} \rangle_{\xi,\tau} + \langle (L_{02} + L_{03} + L_{04}) \hat\rho_{n,F} \rangle_{\xi,\tau} =&\, 0.      \label{e:10:6:SolAspt}
\end{align}

Then we have from \eqref{e:10:2:SolAspt} and \eqref{e:7:2:SolAspt},
\begin{align}
      - \bar L_{0} \bar\rho_{n,R} =&\, \langle L_{11} (\hat\rho_{n+1,0}+\hat\rho_{n+1,M}) \rangle_{\xi,\tau} + \langle (L_{02} + L_{03} + L_{04}) \hat\rho_{n,R} \rangle_{\xi,\tau} ,  \label{e:12:SolAspt} \\ 
      - \bar L_{0} \bar\rho_{n+1,R} =&\, \langle L_{11} \hat\rho_{n+2,M} \rangle_{\xi,\tau} + \langle (L_{02} + L_{03} + L_{04}) ( \hat\rho_{n+1,M}+\hat\rho_{n+1,n} ) \rangle_{\xi,\tau},      \label{e:14:SolAspt} \\ 
      +&\, \langle L_{11} \hat\rho_{n+2,0} \rangle_{\xi,\tau} + \langle (L_{02} + L_{03} + L_{04}) \hat\rho_{n+1,0} \rangle_{\xi,\tau}.     \label{e:16:SolAspt}
\end{align}
The coefficients $\bar h^{(\bbeta)}_{n+1,a}$ with $1 \leq a \leq n$ are determined by the contributions from the line \eqref{e:14:SolAspt}. Note that the coefficients of $\hat\rho_{n}$, $\hat\rho_{n+1}$ and $\hat\rho_{n+2}$ satisfy the stationarity property in Remark \ref{r:stationarity}, hence the claim \eqref{e:6:SolAspt} follows by carefully comparing the contributions of lines \eqref{e:12:SolAspt} and \eqref{e:14:SolAspt} above for $\bar \rho_n$ and $\bar \rho_{n+1}$ respectively.

\end{proof}

The terms in \eqref{e:16:SolAspt} contribute to those in line \eqref{e:5:2:SolAspt} for $a = 0$. The terms in \eqref{e:7:3:SolAspt} contribute to those in line \eqref{e:5:3:SolAspt}. Now we compute both. And we define
\begin{align}     \label{e:18:SolAspt}
      \bar h^{(\bbeta)}_{n+1} := \bar h^{(\bbeta)}_{n+1,0}.
\end{align}
Next we compute the terms in \eqref{e:16:SolAspt}. Then we have
\begin{align}
      \langle L_{11} \hat\rho_{n+2,0} \rangle_{\xi,\tau}
      = -&\, \frac{\delta^{\sfrac12} \mr{\bar \lambda}}{\lambda} \sum_{[\balpha] = 1, \balpha \in \mathcal{I}_*}^{2n+3} \sum_{l=1}^{ \hat I_{n+2}^{(\balpha)} } 
            \langle \eta_1 \hat \eta^{(\balpha, l)}_{n+2} \rangle_{\tau} \langle \partial_{\xi_1} \Pi_1 \hat \chi^{(\balpha, l)}_{n+2} \rangle_\xi \mr {\bar\partial}_{2} \big( \hat h^{(\balpha, l)}_{n+2} \mr {\bar D}^{\balpha} \bar \rho_0 \big)   \label{e:24:SolAspt} \\ 
      +&\, \frac{\delta^{\sfrac12} \mr{\bar \lambda}}{\lambda} \sum_{[\balpha] = 1, \balpha \in \mathcal{I}_*}^{2n+3} \sum_{l=1}^{ \hat I_{n+2}^{(\balpha)} } 
            \langle \eta_2 \hat \eta^{(\balpha, l)}_{n+2} \rangle_{\tau} \langle \partial_{\xi_2} \Pi_2 \hat \chi^{(\balpha, l)}_{n+2} \rangle_\xi \mr {\bar\partial}_{1} \big( \hat h^{(\balpha, l)}_{n+2} \mr {\bar D}^{\balpha} \bar \rho_0 \big).  \label{e:25:SolAspt}
\end{align}
and 
\begin{align}
      &\, \langle L_{11} \hat\rho_{n+2,F} \rangle_{\xi,\tau} 
            = \bigg( \sum_{k=1}^K \sum_{[\balpha] = 1, \balpha \in \mathcal{I}_*}^{2n+2+Q} \sum_{l=1}^{ \hat I_{F,n+2}^{(\balpha)} }  \text{ below} \bigg)  \nonumber \\
      = -&\, \frac{\delta^{\sfrac12} \mr{\bar \lambda}}{\lambda} 
            \sum \sum \sum
            \langle \eta_1 \hat \eta^{(\balpha, l)}_{F,k,n+2} \rangle_{\tau} \langle \partial_{\xi_1} \Pi_1 \hat \chi^{(\balpha, l)}_{F,k,n+2} \rangle_\xi \mr {\bar\partial}_{2} \big( \hat h^{(\balpha, l)}_{F,k,n+2} \mr {\bar D}^{\balpha} \bar \rho_{F,k} \big)   \label{e:26:SolAspt} \\ 
      +&\, \frac{\delta^{\sfrac12} \mr{\bar \lambda}}{\lambda} 
            \sum \sum \sum
            \langle \eta_2 \hat \eta^{(\balpha, l)}_{F,k,n+2} \rangle_{\tau} \langle \partial_{\xi_2} \Pi_2 \hat \chi^{(\balpha, l)}_{F,k,n+2} \rangle_\xi \mr {\bar\partial}_{1} \big( \hat h^{(\balpha, l)}_{F,k,n+2} \mr {\bar D}^{\balpha} \bar \rho_{F,k} \big).  \label{e:27:SolAspt}
\end{align}

For next term, we have
\begin{align}
      \langle L_{02} \hat\rho_{n+1,0} \rangle_{\xi,\tau} =&\,
      - \varepsilon \frac{\delta^{\sfrac12}}{\lambda} \sum_{[\balpha] = 1, \balpha \in \mathcal{I}_*}^{2n+1} \sum_{l=1}^{ \hat I_{n+1}^{(\balpha)} } \sum_{m=1}^N \partial_{i} \Omega_{m,ij} H_{12} \sigma_m \hat\eta^{(\balpha, l)}_{n+1} \hat\chi^{(\balpha, l)}_{n+1} \partial_{j} \big( \hat h^{(\balpha, l)}_{n+1} \mr {\bar D}^{\balpha} \bar \rho_0 \big)        \nonumber  \\ 
      =- &\, \varepsilon \frac{ \delta^{\sfrac12} \mr{\bar \lambda}^2 } {\lambda} \sum\sum\sum 
      \langle \sigma_m \eta_1 \hat\eta^{(\balpha, l)}_{n+1} \rangle_\tau \langle \partial_{\xi_1} \Pi_1 \hat \chi^{(\balpha, l)}_{n+1} \rangle_{\xi} 
      \mr {\bar\partial}_{i} \Big( \Omega_{m,ij} \mr {\bar\partial}_{j} \big( \hat h^{(\balpha, l)}_{n+1} \mr {\bar D}^{\balpha} \bar \rho_0 \big) \Big)    \label{e:30:SolAspt} \\ 
      -&\, \varepsilon \frac{ \delta^{\sfrac12} \mr{\bar \lambda}^2 } {\lambda} \sum\sum\sum 
      \langle \sigma_m \eta_2 \hat\eta^{(\balpha, l)}_{n+1} \rangle_\tau \langle \partial_{\xi_2} \Pi_2 \hat\chi^{(\balpha, l)}_{n+1} \rangle_{\xi}
      \mr {\bar\partial}_{i} \Big( \Omega_{m,ij} \mr {\bar\partial}_{j} \big( \hat h^{(\balpha, l)}_{n+1} \mr {\bar D}^{\balpha} \bar \rho_0 \big) \Big)    \label{e:32:SolAspt} 
\end{align}
and
\begin{align}
      &\, \langle L_{02} \hat\rho_{n+1,F} \rangle_{\xi,\tau} \nonumber \\
      =&\, - \varepsilon \frac{\delta^{\sfrac12}}{\lambda} \sum_{k=1}^K \sum_{[\balpha] = 1, \balpha \in \mathcal{I}_*}^{2n+Q} \sum_{l=1}^{ \hat I_{F,n+1}^{(\balpha)} } \sum_{m=1}^N 
      \partial_{i} \Omega_{m,ij} H_{12} \sigma_m \hat\eta^{(\balpha, l)}_{F,k,n+1} \hat\chi^{(\balpha, l)}_{F,k,n+1} \partial_{j} \big( \hat h^{(\balpha, l)}_{F,k,n+1} \mr {\bar D}^{\balpha} \bar \rho_{F,k} \big)        \nonumber \\ 
      =&\,- \varepsilon \frac{ \delta^{\sfrac12} \mr{\bar \lambda}^2 } {\lambda} \sum\sum\sum\sum 
      \langle \sigma_m \eta_1 \hat\eta^{(\balpha, l)}_{F,k,n+1} \rangle_\tau \langle \partial_{\xi_1} \Pi_1 \hat \chi^{(\balpha, l)}_{F,k,n+1} \rangle_{\xi} 
      \mr {\bar\partial}_{i} \Big( \Omega_{m,ij} \mr {\bar\partial}_{j} \big( \hat h^{(\balpha, l)}_{F,k,n+1} \mr {\bar D}^{\balpha} \bar \rho_{F,k} \big) \Big)    \label{e:34:SolAspt} \\ 
      &\,- \varepsilon \frac{ \delta^{\sfrac12} \mr{\bar \lambda}^2 } {\lambda} \sum\sum\sum\sum 
      \langle \sigma_m \eta_2 \hat\eta^{(\balpha, l)}_{F,k,n+1} \rangle_\tau \langle \partial_{\xi_2} \Pi_2 \hat\chi^{(\balpha, l)}_{F,k,n+1} \rangle_{\xi}
      \mr {\bar\partial}_{i} \Big( \Omega_{m,ij} \mr {\bar\partial}_{j} \big( \hat h^{(\balpha, l)}_{F,k,n+1} \mr {\bar D}^{\balpha} \bar \rho_{F,k} \big) \Big)    \label{e:36:SolAspt} 
\end{align}

We also have
\begin{align}
      \langle L_{03} \hat\rho_{n+1,0} \rangle_{\xi,\tau} =&\,
      \bigg( \sum_{[\balpha] = 1, \balpha \in \mathcal{I}_*}^{2n+1} \sum_{l=1}^{ \hat I_{n+1}^{(\balpha)} } \sum_{m=1}^N \text{ below} \bigg)         \nonumber \\ 
      =&\, - \varepsilon \delta^{\sfrac12} \mr{\bar \lambda} \sum\sum\sum 
      E_{m,j2} \langle \varphi_m \eta_1 \hat\eta^{(\balpha, l)}_{n+1} \rangle_\tau \langle \partial_{\xi_1} \Pi_1 \hat\chi^{(\balpha, l)}_{n+1} \rangle_{\xi} \mr {\bar\partial}_{j} \big( \hat h^{(\balpha, l)}_{n+1} \mr {\bar D}^{\balpha} \bar \rho_0 \big)      \label{e:39:SolAspt} \\ 
      &\, + \varepsilon \delta^{\sfrac12} \mr{\bar \lambda} \sum\sum\sum 
      E_{m,j1} \langle \varphi_m \eta_2 \hat\eta^{(\balpha, l)}_{n+1} \rangle_\tau \langle \partial_{\xi_2} \Pi_2 \hat\chi^{(\balpha, l)}_{n+1} \rangle_{\xi} \mr {\bar\partial}_{j} \big( \hat h^{(\balpha, l)}_{n+1} \mr {\bar D}^{\balpha} \bar \rho_0 \big)       \label{e:41:SolAspt} 
\end{align}
and
\begin{align}
      \langle L_{04} \hat\rho_{n+1,0} \rangle_{\xi,\tau} =&\, \bigg( \sum_{[\balpha] = 1, \balpha \in \mathcal{I}_*}^{2n+1} \sum_{l=1}^{ \hat I_{n+1}^{(\balpha)} } \sum_{m=1}^N \text{ below} \bigg)       \nonumber \\ 
      =&\,+ \varepsilon \delta^{\sfrac12} \mr{\bar \lambda} \sum\sum\sum \mr {\bar\partial}_{j} E_{m,j2} \hat h^{(\balpha, l)}_{n+1} \langle \varphi_m \eta_1 \hat\eta^{(\balpha, l)}_{n+1} \rangle_\tau \langle \Pi_1 \partial_{\xi_1} \hat\chi^{(\balpha, l)}_{n+1} \rangle_{\xi} \mr {\bar D}^{\balpha} \bar \rho_0       \label{e:46:SolAspt} \\ 
      &\,- \varepsilon \delta^{\sfrac12} \mr{\bar \lambda} \sum\sum\sum \mr {\bar\partial}_{j} E_{m,j1} \hat h^{(\balpha, l)}_{n+1} \langle \varphi_m \eta_2 \hat\eta^{(\balpha, l)}_{n+1} \rangle_\tau \langle \Pi_2 \partial_{\xi_2} \hat\chi^{(\balpha, l)}_{n+1} \rangle_{\xi} \mr {\bar D}^{\balpha} \bar \rho_0.        \label{e:47:SolAspt}
\end{align}

Similarly, we have
\begin{align}
      &\, \langle (L_{03} + L_{04}) \hat\rho_{n+1,F} \rangle_{\xi,\tau} =
            \bigg( \sum_{k=1}^K \sum_{[\balpha] = 1, \balpha \in \mathcal{I}_*}^{2n+Q} \sum_{l=1}^{ \hat I_{F,n+1}^{(\balpha)} } \sum_{m=1}^N \text{ below} \bigg) \nonumber \\ 
      =&\, - \varepsilon \delta^{\sfrac12} \mr{\bar \lambda}
      \sum \sum \sum \sum
            \langle \varphi_m \eta_1 \hat\eta^{(\balpha, l)}_{F,k,n+1} \rangle_\tau \langle \partial_{\xi_1} \Pi_1 \hat\chi^{(\balpha, l)}_{F,k,n+1} \rangle_{\xi} 
            \mr {\bar\partial}_{j} \big( E_{m,j2} \hat h^{(\balpha, l)}_{F,k,n+1} \mr {\bar D}^{\balpha} \bar \rho_{F,k} \big)      \label{e:52:SolAspt} \\ 
      &\, + \varepsilon \delta^{\sfrac12} \mr{\bar \lambda}
      \sum \sum \sum \sum
            \langle \varphi_m \eta_2 \hat\eta^{(\balpha, l)}_{F,k,n+1} \rangle_\tau \langle \partial_{\xi_2} \Pi_2 \hat\chi^{(\balpha, l)}_{F,k,n+1} \rangle_{\xi} \mr {\bar\partial}_{j} \big( E_{m,j1} \hat h^{(\balpha, l)}_{F,k,n+1} \mr {\bar D}^{\balpha} \bar \rho_{F,k} \big)      \label{e:53:SolAspt} 
\end{align}

Now we go back to \eqref{e:16:SolAspt} to compute the coefficients, from the contributions in (\ref{e:24:SolAspt}-\ref{e:47:SolAspt}). We need to write these contributions in the form of the terms in \eqref{e:5:SolAspt} and \eqref{e:18:SolAspt}, which leads to the following cases.

\begin{case}      \label{c:compBar:0}
If $1 \leq [\bbeta] \leq 2n+3$, $\bbeta = \balpha$, this contribution from lines \eqref{e:24:SolAspt} and \eqref{e:25:SolAspt} concerning $L_{11} \hat \rho_{n+2}$ gives
\begin{align*}
      \bar h^{(\bbeta)}_{n+1} = \frac{\delta^{\sfrac12}\mr{\bar \lambda}}{\lambda} \sum_{l=1}^{\hat I^{(\balpha)}_{n+2}} \hat h^{(\balpha, l)}_{n+2} 
      \begin{bmatrix} \langle \eta_2 \hat \eta^{(\balpha, l)}_{n+2} \rangle_{\tau} \langle \partial_{\xi_2} \Pi_2 \hat \chi^{(\balpha, l)}_{n+2} \rangle_\xi \\ 
      - \langle \eta_1 \hat \eta^{(\balpha, l)}_{n+2} \rangle_{\tau} \langle \partial_{\xi_1} \Pi_1 \hat \chi^{(\balpha, l)}_{n+2} \rangle_\xi  \end{bmatrix}.
\end{align*}
The estimate \eqref{e:EXP1:bh_Pest} with $k=n$ follows from (\ref{e:EXP1:hh_Pest}-\ref{e:EXP1:heta_est}) with $k=n+1$.   
\end{case}

\begin{case}      \label{c:compBar:1}
If $1 \leq [\bbeta] \leq 2n+1$, $\bbeta = \balpha$, one contribution from line \eqref{e:30:SolAspt} (with another contribution given in the next case) concerning $L_{02} \hat \rho_{n+1}$ gives
\begin{align*}
      \bar h^{(\bbeta)}_{n+1} = - \varepsilon \frac{ \delta^{\sfrac12} \mr{\bar \lambda}^2 } {\lambda} \sum_{l=1}^{\hat I^{(\balpha)}_{n+1}} \sum_m 
      \Omega_{m} \mr {\bar \nabla} \hat h^{(\balpha, l)}_{n+1} \langle \sigma_m \eta_1 \hat\eta^{(\balpha, l)}_{n+1} \rangle_\tau \langle \partial_{\xi_1} \Pi_1 \hat\chi^{(\balpha, l)}_{n+1} \rangle_{\xi} . 
\end{align*}
The terms involving $\eta_2$ and $\Pi_2$ in \eqref{e:32:SolAspt} are symmetric.
The estimate \eqref{e:EXP1:bh_Pest} with $k=n$ follows from \eqref{e:mtrEstHom} and (\ref{e:EXP1:hh_Pest}-\ref{e:EXP1:heta_est}) with $k=n$.
\end{case}

\begin{case}      \label{c:compBar:2}
If $2 \leq [\bbeta] \leq 2n+2$, $\bbeta = 1 \balpha$, another contribution from line concerning \eqref{e:30:SolAspt} concerning $L_{02} \hat \rho_{n+1}$ gives 
\begin{align*}
      \bar h^{(\bbeta)}_{n+1} = - \varepsilon \frac{ \delta^{\sfrac12} \mr{\bar \lambda}^2 } {\lambda} \sum_{l=1}^{\hat I^{(\balpha)}_{n+1}} \sum_m 
            \hat h^{(\balpha, l)}_{n+1} \langle \sigma_m \eta_1 \hat\eta^{(\balpha, l)}_{n+1} \rangle_\tau \langle \partial_{\xi_1} \Pi_1 \hat\chi^{(\balpha, l)}_{n+1} \rangle_{\xi} 
            \begin{bmatrix} \Omega_{m,11} \\ \Omega_{m,21} \end{bmatrix} . 
\end{align*}
\eqref{e:30:SolAspt} also contribute to symmetric terms for $\bbeta = 2 \balpha$. The terms involving $\eta_2$ and $\Pi_2$ in \eqref{e:32:SolAspt} are symmetric. The estimate \eqref{e:EXP1:bh_Pest} with $k=n$ follows from \eqref{e:mtrEstHom} and (\ref{e:EXP1:hh_Pest}-\ref{e:EXP1:heta_est}) with $k=n$.
\end{case}

\begin{case}      \label{c:compBar:3}
If $1 \leq [\bbeta] \leq 2n+1$, $\bbeta = \balpha$, using product rule, the contribution from line \eqref{e:39:SolAspt} concerning $L_{03} \hat \rho_{n+1}$ and line \eqref{e:46:SolAspt} concerning $L_{04} \hat \rho_{n+1}$ give
\begin{align*}
      \bar h^{(\bbeta)}_{n+1} = - \varepsilon \delta^{\sfrac12} \mr{\bar \lambda} 
            \sum_{l=1}^{\hat I^{(\balpha)}_{n+1}} \sum_m 
            \langle \varphi_m \eta_1 \hat\eta^{(\balpha, l)}_{n+1} \rangle_\tau \langle \partial_{\xi_1} \Pi_1 \hat\chi^{(\balpha, l)}_{n+1} \rangle_{\xi} 
            \begin{bmatrix} 
            E_{m,12} \hat h^{(\balpha, l)}_{n+1} \\ 
            E_{m,22} \hat h^{(\balpha, l)}_{n+1}
            \end{bmatrix}. 
\end{align*}
The lines \eqref{e:41:SolAspt} and \eqref{e:47:SolAspt} give symmetric terms involving $\eta_2$ and $\Pi_2$. The estimate \eqref{e:EXP1:bh_Pest} with $k=n$ follows from \eqref{e:mtrEstHom} and (\ref{e:EXP1:hh_Pest}-\ref{e:EXP1:heta_est}) with $k=n$. Indeed, we just need
\begin{align*}
      \varepsilon \kappa \lambda \leq 
      \kappa \cdot \frac{\delta^{\sfrac12} \mr{\bar \lambda}}{\kappa\lambda} \cdot \lambda^{\gamma_S}.
\end{align*}
\end{case}

\eqref{e:EXP1:bhS_Pest} follows from \eqref{e:EXP1:bh_Pest} for all admissible $\balpha$. \eqref{e:EXP1:bh} follows from \eqref{e:EXP1:hh} for $n+1$ and $n+2$.

Now we compare \eqref{e:7:3:SolAspt} and \eqref{e:16:SolAspt}. We also compare the $F$-related terms in \eqref{e:26:SolAspt}, \eqref{e:27:SolAspt}, \eqref{e:30:SolAspt}, \eqref{e:32:SolAspt}, \eqref{e:52:SolAspt}, \eqref{e:53:SolAspt} and other terms in (\ref{e:24:SolAspt}-\ref{e:53:SolAspt}). It is clear that analyzing the coefficients in \eqref{e:5:3:SolAspt} is the same as above. We omit such details for Lemma \ref{l:EXP2}.

\subsubsection{Inductive step for temporal correctors}      \label{sss:solAsptTPR}

\underline{Next, we derive the structures of $\check \rho_{n+2}$.} We take the difference of $\langle \eqref{e:aspt:hrc4} \rangle_\xi$ and $\langle \eqref{e:aspt:hrc4} \rangle_{\xi,\tau}$. Then we apply Remark \ref{r:hatcheckbar} to deduce
\begin{equation}  \label{e:78:SolAspt}
\begin{split}
      L_{12} \check\rho_{n+2} =&\, - L_{01} \check\rho_{n+1} 
      + \langle L_{11} \hat\rho_{n+2} \rangle_{\xi,\tau} - \langle L_{11} \hat\rho_{n+2} \rangle_{\xi}
      + \langle L_{02} \hat\rho_{n+1} \rangle_{\xi,\tau} - \langle L_{02} \hat\rho_{n+1} \rangle_{\xi} \\
      &\,+ \langle L_{03} \hat\rho_{n+1} \rangle_{\xi,\tau} - \langle L_{03} \hat\rho_{n+1} \rangle_{\xi}
      + \langle L_{04} \hat\rho_{n+1} \rangle_{\xi,\tau} - \langle L_{04} \hat\rho_{n+1} \rangle_{\xi},
\end{split}
\end{equation}
from which we shall compute $\check \rho_{n+2}$.

We claim the stationary property holds for the coefficients in the following formula
\begin{equation}     \label{e:79:SolAspt}
\begin{split}
      \check \rho_{n+2}(x,t,\tau) =& \sum_{a=0}^{n+1} \sum_{[\bbeta] = 1, \bbeta \in \mathcal{I}_*}^{2(n-a)+4} \sum_{l=1}^{\check I_{n+2,a}^{(\bbeta)}}
            \check h^{(\bbeta, l)}_{n+2,a} (x,t) \check \eta^{(\bbeta, l)}_{n+2,a} (\tau) \mr {\bar D}^{\bbeta} \bar \rho_a(x,t)  \\ 
      +&\, \sum_{[\bbeta]=n_0, \bbeta \in \mathcal{I}_*}^{2n+3+n_0} \sum_{l=1}^{\check I_{F,n+2}^{(\bbeta)}} 
            \check h^{(\bbeta,l)}_{F,n+2} (x,t) \check \eta^{(\bbeta,l)}_{F,n+2} (\tau) \mr {\bar D}^{\bbeta} \bar \rho_F (x,t) 
            :=  \check\rho_{n+2,R} + \check\rho_{n+2,F}.
\end{split}
\end{equation}

\begin{proof}[Proof of stationarity for $\check \rho_{n+2}$, i.e. the claim above]
Indeed, we take the difference of $\langle \eqref{e:8:SolAspt} \rangle_\xi$ and $\langle \eqref{e:8:SolAspt} \rangle_{\xi,\tau}$ for $n=0$, or the difference of $\langle \eqref{e:9:SolAspt} \rangle_\xi$ and $\langle \eqref{e:9:SolAspt} \rangle_{\xi,\tau}$ for $n \geq 1$. Either of them gives
\begin{align}
      L_{12} \check\rho_{n+1} =&\, - L_{01} \check\rho_{n} 
      + \langle L_{11} \hat\rho_{n+1} \rangle_{\xi,\tau} - \langle L_{11} \hat\rho_{n+1} \rangle_{\xi},    \label{e:80:0:SolAspt} \\
      &\,+ \langle ( L_{02} + L_{03} + L_{04} ) \hat\rho_{n} \rangle_{\xi,\tau} - \langle ( L_{02} + L_{03} + L_{04} ) \hat\rho_{n} \rangle_{\xi}.    \label{e:80:2:SolAspt}
\end{align}
Subsequently, from the decomposition notation introduced in \eqref{e:2:SolAspt}, we have 
\begin{align}
      L_{12} \check\rho_{n+1,R} =&\, - L_{01} \check\rho_{n,R} 
      + \langle L_{11} \hat\rho_{n+1,R} \rangle_{\xi,\tau} - \langle L_{11} \hat\rho_{n+1,R} \rangle_{\xi}    \label{e:80:4:SolAspt} \\
      &\,+ \langle ( L_{02} + L_{03} + L_{04} ) \hat\rho_{n,R} \rangle_{\xi,\tau} - \langle ( L_{02} + L_{03} + L_{04} ) \hat\rho_{n,R} \rangle_{\xi},   \label{e:81:0:SolAspt} \\ 
      L_{12} \check\rho_{n+1,F} =&\, - L_{01} \check\rho_{n,F} 
      + \langle L_{11} \hat\rho_{n+1,F} \rangle_{\xi,\tau} - \langle L_{11} \hat\rho_{n+1,F} \rangle_{\xi}    \label{e:81:2:SolAspt} \\
      &\,+ \langle ( L_{02} + L_{03} + L_{04} ) \hat\rho_{n,F} \rangle_{\xi,\tau} - \langle ( L_{02} + L_{03} + L_{04} ) \hat\rho_{n,F} \rangle_{\xi}.   \label{e:81:4:SolAspt} 
\end{align}

On the other hand, following the similar decomposition notation in \eqref{e:79:SolAspt} and (\ref{e:1:SolAspt},\ref{e:4:2:SolAspt},\ref{e:4:6:SolAspt}), we deduce, from \eqref{e:78:SolAspt}, that
\begin{align}
      L_{12} \check\rho_{n+2,F}
      =&\, - L_{01} \check\rho_{n+1,F} 
      + \langle L_{11} \hat\rho_{n+2,F} \rangle_{\xi,\tau} - \langle L_{11} \hat\rho_{n+2,F} \rangle_{\xi}  \label{e:81:6:SolAspt} \\
      &\,+ \langle ( L_{02} + L_{03} + L_{04} ) \hat\rho_{n+1,F} \rangle_{\xi,\tau} - \langle ( L_{02} + L_{03} + L_{04} ) \hat\rho_{n+1,F} \rangle_{\xi},  \label{e:81:8:SolAspt} \\ 
      L_{12} \check\rho_{n+2,R}
      =&\, - L_{01} \check\rho_{n+1,E} 
      + \langle L_{11} \hat\rho_{n+2,E} \rangle_{\xi,\tau} - \langle L_{11} \hat\rho_{n+2,E} \rangle_{\xi}      \label{e:82:SolAspt} \\
      &\,+ \langle ( L_{02} + L_{03} + L_{04} ) \hat\rho_{n+1,E} \rangle_{\xi,\tau} - \langle ( L_{02} + L_{03} + L_{04} ) \hat\rho_{n+1,E} \rangle_{\xi}     \label{e:83:SolAspt}  \\
      &\, - L_{01} \check\rho_{n+1,0} 
      + \langle L_{11} \hat\rho_{n+2,0} \rangle_{\xi,\tau} - \langle L_{11} \hat\rho_{n+2,0} \rangle_{\xi}        \label{e:84:SolAspt} \\
      &\,+ \langle ( L_{02} + L_{03} + L_{04} ) \hat\rho_{n+1,0} \rangle_{\xi,\tau} - \langle ( L_{02} + L_{03} + L_{04} ) \hat\rho_{n+1,0} \rangle_{\xi}.       \label{e:85:SolAspt}
\end{align}
\eqref{e:82:SolAspt} and \eqref{e:83:SolAspt} contribute to the terms in \eqref{e:79:SolAspt} for $a \geq 1$, and \eqref{e:84:SolAspt} and \eqref{e:85:SolAspt} contribute to the terms in \eqref{e:79:SolAspt} for $a = 0$. Now we compare the terms in (\ref{e:80:4:SolAspt}-\ref{e:81:0:SolAspt}) and the terms in (\ref{e:82:SolAspt}-\ref{e:83:SolAspt}). Note that the coefficients of $\hat \rho_n$, $\hat \rho_{n+1}$, $\hat \rho_{n+2}$, $\check \rho_n$, $\check \rho_{n+1}$ satisfy the stationary property in Remark \ref{r:stationarity}, as stated in the induction assumptions. Therefore, it is clear that the stationary property still holds for the contributions of (\ref{e:82:SolAspt}-\ref{e:83:SolAspt}) in $\check \rho_{n+2}$.

\end{proof}

Now it remains to compute the terms in \eqref{e:84:SolAspt} and \eqref{e:85:SolAspt}, i.e. the coefficients in \eqref{e:79:SolAspt} for $a=0$. Denoting
\begin{align*}
      \check I_{n+2}^{(\bbeta)} := \check I_{n+2,0}^{(\bbeta)}, \quad
      \check h^{(\bbeta, l)}_{n+2} := \check h^{(\bbeta, l)}_{n+2,0}, \quad
      \check \eta^{(\bbeta, l)}_{n+2} := \check \eta^{(\bbeta, l)}_{n+2,0}, 
\end{align*}
we can write \eqref{e:79:SolAspt} as 
\begin{equation}     \label{e:86:SolAspt}
\begin{split}
      \check \rho_{n+2}(x,t,\tau) =& \sum_{a=0}^{n+1} \sum_{[\bbeta] = 1, \bbeta \in \mathcal{I}_*}^{2(n-a)+4} \sum_{l=1}^{\check I_{n-a+2}^{(\bbeta)}}
            \check h^{(\bbeta, l)}_{n-a+2} (x,t) \check \eta^{(\bbeta, l)}_{n-a+2} (\tau) \mr {\bar D}^{\bbeta} \bar \rho_a(x,t)   \\ 
      +&\, \sum_{[\bbeta]=n_0, \bbeta \in \mathcal{I}_*}^{2n+3+n_0} \sum_{l=1}^{\check I_{F,n+2}^{(\bbeta)}} \check h^{(\bbeta,l)}_{F,n+2} (x,t) \check \eta^{(\bbeta,l)}_{F,n+2} (\tau) \mr {\bar D}^{\bbeta} \bar \rho_F (x,t). 
\end{split}
\end{equation}

\underline{For $a=0$, we now compute the terms in \eqref{e:84:SolAspt} and \eqref{e:85:SolAspt}.} Using the induction assumptions, we have
\begin{align}
      \langle L_{11} \hat\rho_{n+2,0} \rangle_{\xi,\tau} -&\, \langle L_{11} \hat\rho_{n+2,0} \rangle_{\xi} = \Bigg( \sum_{[\balpha] = 1, \balpha \in \mathcal{I}_*}^{2n+3} \sum_{l=1}^{ \hat I_{n+2}^{(\balpha)} } \, \text{below} \Bigg)    \nonumber \\ 
      =&\, \frac{\delta^{\sfrac12} \mr{\bar \lambda}} {\lambda} \sum\sum
            \hat h^{(\balpha, l)}_{n+2} \big( \eta_1 \hat \eta^{(\balpha, l)}_{n+2} - \langle \eta_1 \hat \eta^{(\balpha, l)}_{n+2} \rangle_{\tau} \big) \langle \partial_{\xi_1} \Pi_1 \hat \chi^{(\balpha, l)}_{n+2} \rangle_\xi \mr {\bar \partial}_{2} \mr {\bar D}^{\balpha} \bar \rho_0,    \label{e:90:SolAspt} \\ 
      -&\, \frac{\delta^{\sfrac12} \mr{\bar \lambda}} {\lambda} \sum \sum 
            \hat h^{(\balpha, l)}_{n+2} \big( \eta_2 \hat \eta^{(\balpha, l)}_{n+2} - \langle \eta_2 \hat \eta^{(\balpha, l)}_{n+2} \rangle_{\tau} \big) \langle \partial_{\xi_2} \Pi_2 \hat \chi^{(\balpha, l)}_{n+2} \rangle_\xi \mr {\bar \partial}_{1} \mr {\bar D}^{\balpha} \bar \rho_0,    \label{e:91:SolAspt}
\end{align}

and
\begin{align}
      \langle L_{02} &\, \hat\rho_{n+1,0} \rangle_{\xi,\tau} - \langle L_{02} \hat\rho_{n+1,0} \rangle_{\xi} = \Bigg( \sum_{[\balpha] = 1, \balpha \in \mathcal{I}_*}^{2n+1} \sum_{l=1}^{ \hat I_{n+1}^{(\balpha)} } \sum_m \text{ below} \Bigg)   \nonumber \\ 
      =&\, \varepsilon \frac{\delta^{\sfrac12} \mr{\bar \lambda}^2} {\lambda} \sum\sum\sum 
            \mr {\bar \partial}_{i} \Omega_{m,ij} \mr {\bar \partial}_{j} \hat h^{(\balpha, l)}_{n+1} \big( \sigma_m \eta_1 \hat\eta^{(\balpha, l)}_{n+1} - \langle \sigma_m \eta_1 \hat\eta^{(\balpha, l)}_{n+1} \rangle_\tau \big) \langle \partial_{\xi_1} \Pi_1 \hat\chi^{(\balpha, l)}_{n+1} \rangle_{\xi} \mr {\bar D}^{\balpha} \bar \rho_0     \label{e:92:SolAspt} \\ 
      +&\, \varepsilon \frac{\delta^{\sfrac12} \mr{\bar \lambda}^2} {\lambda} \sum\sum\sum 
            \mr {\bar \partial}_{i} \Omega_{m,ij} \hat h^{(\balpha, l)}_{n+1} \big( \sigma_m \eta_1 \hat\eta^{(\balpha, l)}_{n+1} - \langle \sigma_m \eta_1 \hat\eta^{(\balpha, l)}_{n+1} \rangle_\tau \big) \langle \partial_{\xi_1} \Pi_1 \hat\chi^{(\balpha, l)}_{n+1} \rangle_{\xi} \mr {\bar \partial}_{j} \mr {\bar D}^{\balpha} \bar \rho_0     \label{e:93:SolAspt} \\ 
      +&\, \varepsilon \frac{\delta^{\sfrac12} \mr{\bar \lambda}^2} {\lambda} \sum\sum\sum 
            \mr {\bar \partial}_{i} \Omega_{m,ij} \mr {\bar \partial}_{j} \hat h^{(\balpha, l)}_{n+1} \big( \sigma_m \eta_2 \hat\eta^{(\balpha, l)}_{n+1} - \langle \sigma_m \eta_2 \hat\eta^{(\balpha, l)}_{n+1} \rangle_\tau \big) \langle \partial_{\xi_2} \Pi_2 \hat\chi^{(\balpha, l)}_{n+1} \rangle_{\xi} \mr {\bar D}^{\balpha} \bar \rho_0     \label{e:94:SolAspt} \\ 
      +&\, \varepsilon \frac{\delta^{\sfrac12} \mr{\bar \lambda}^2} {\lambda} \sum\sum\sum 
            \mr {\bar \partial}_{i} \Omega_{m,ij} \hat h^{(\balpha, l)}_{n+1} \big( \sigma_m \eta_2 \hat\eta^{(\balpha, l)}_{n+1} - \langle \sigma_m \eta_2 \hat\eta^{(\balpha, l)}_{n+1} \rangle_\tau \big) \langle \partial_{\xi_2} \Pi_2 \hat\chi^{(\balpha, l)}_{n+1} \rangle_{\xi} \mr {\bar \partial}_{j} \mr {\bar D}^{\balpha} \bar \rho_0,     \label{e:95:SolAspt}
\end{align}

and
\begin{align}
      \langle L_{03}&\, \hat\rho_{n+1,0} \rangle_{\xi,\tau} - \langle L_{03} \hat\rho_{n+1,0} \rangle_{\xi} 
      = \Bigg( \sum_{[\balpha] = 1, \balpha \in \mathcal{I}_*}^{2n+1} \sum_{l=1}^{ \hat I_{n+1}^{(\balpha)} } \sum_m \text{ below} \Bigg)       \nonumber \\ 
      =&\, \varepsilon \delta^{\sfrac12} \mr{\bar \lambda} \sum\sum\sum
            E_{m,j2} \mr {\bar \partial}_{j} \hat h^{(\balpha, l)}_{n+1} \big( \varphi_m \eta_1 \hat\eta^{(\balpha, l)}_{n+1} - \langle \varphi_m \eta_1 \hat\eta^{(\balpha, l)}_{n+1} \rangle_\tau \big) \langle \partial_{\xi_1} \Pi_1 \hat\chi^{(\balpha, l)}_{n+1} \rangle_{\xi} \mr {\bar D}^{\balpha} \bar \rho_0       \label{e:97:SolAspt} \\ 
      +&\, \varepsilon \delta^{\sfrac12} \mr{\bar \lambda} \sum\sum\sum
            E_{m,j2} \hat h^{(\balpha, l)}_{n+1} \big( \varphi_m \eta_1 \hat\eta^{(\balpha, l)}_{n+1} - \langle \varphi_m \eta_1 \hat\eta^{(\balpha, l)}_{n+1} \rangle_\tau \big) \langle \partial_{\xi_1} \Pi_1 \hat\chi^{(\balpha, l)}_{n+1} \rangle_{\xi} \mr {\bar \partial}_{j} \mr {\bar D}^{\balpha} \bar \rho_0       \label{e:98:SolAspt} \\ 
      -&\, \varepsilon \delta^{\sfrac12} \mr{\bar \lambda} \sum\sum\sum
            E_{m,j1} \mr {\bar \partial}_{j} \hat h^{(\balpha, l)}_{n+1} \big( \varphi_m \eta_2 \hat\eta^{(\balpha, l)}_{n+1} - \langle \varphi_m \eta_2 \hat\eta^{(\balpha, l)}_{n+1} \rangle_\tau \big) \langle \partial_{\xi_2} \Pi_2 \hat\chi^{(\balpha, l)}_{n+1} \rangle_{\xi} \mr {\bar D}^{\balpha} \bar \rho_0       \label{e:99:SolAspt} \\ 
      -&\, \varepsilon \delta^{\sfrac12} \mr{\bar \lambda} \sum\sum\sum
            E_{m,j1} \hat h^{(\balpha, l)}_{n+1} \big( \varphi_m \eta_2 \hat\eta^{(\balpha, l)}_{n+1} - \langle \varphi_m \eta_2 \hat\eta^{(\balpha, l)}_{n+1} \rangle_\tau \big) \langle \partial_{\xi_2} \Pi_2 \hat\chi^{(\balpha, l)}_{n+1} \rangle_{\xi} \mr {\bar \partial}_{j} \mr {\bar D}^{\balpha} \bar \rho_0,      \label{e:100:SolAspt} 
\end{align}

and
\begin{align}
      \langle L_{04} &\, \hat\rho_{n+1,0} \rangle_{\xi,\tau} - \langle L_{04} \hat\rho_{n+1,0} \rangle_{\xi} 
      = \Bigg( \sum_{[\balpha] = 1, \balpha \in \mathcal{I}_*}^{2n+1} \sum_{l=1}^{ \hat I_{n+1}^{(\balpha)} } \sum_m \text{ below} \Bigg)       \nonumber \\ 
      =&\, - \varepsilon \delta^{\sfrac12} \mr{\bar \lambda} \sum\sum\sum
            \mr {\bar \partial}_{j} E_{m,j2} \hat h^{(\balpha, l)}_{n+1} \big( \varphi_m \eta_1 \hat\eta^{(\balpha, l)}_{n+1} - \langle \varphi_m \eta_1 \hat\eta^{(\balpha, l)}_{n+1} \rangle_\tau \big) \langle \Pi_1 \partial_{\xi_1} \hat\chi^{(\balpha, l)}_{n+1} \rangle_{\xi} \mr {\bar D}^{\balpha} \bar \rho_0        \label{e:101:SolAspt} \\ 
      &\,+ \varepsilon \delta^{\sfrac12} \mr{\bar \lambda} \sum\sum\sum
            \mr {\bar \partial}_{j} E_{m,j1} \hat h^{(\balpha, l)}_{n+1} \big( \varphi_m \eta_2 \hat\eta^{(\balpha, l)}_{n+1} - \langle \varphi_m \eta_2 \hat\eta^{(\balpha, l)}_{n+1} \rangle_\tau \big) \langle \Pi_2 \partial_{\xi_2} \hat\chi^{(\balpha, l)}_{n+1} \rangle_{\xi} \mr {\bar D}^{\balpha} \bar \rho_0,      \label{e:102:SolAspt}
\end{align}

and
\begin{align}
      - L_{01} \check \rho_{n+1,0} =&\, - \mr{\bar\mu} \sum_{[\balpha] = 1, \balpha \in \mathcal{I}_*}^{2(n+1)} \sum_{l=1}^{\check I_{n+1}^{(\balpha)}}
            \mr {\bar D}_t \check h^{(\balpha, l)}_{n+1} \check \eta^{(\balpha, l)}_{n+1} \mr {\bar D}^{\balpha} \bar \rho_0        \label{e:103:SolAspt} \\ 
      &\, + \kappa \mr{\bar \lambda}^2 \sum_{[\balpha] = 1, \balpha \in \mathcal{I}_*}^{2(n+1)} \sum_{l=1}^{\check I_{n+1}^{(\balpha)}}
            \mr {\bar \partial}_{ii} \check h^{(\balpha, l)}_{n+1} \check \eta^{(\balpha, l)}_{n+1} \mr {\bar D}^{\balpha} \bar \rho_0        \label{e:104:SolAspt} \\ 
      &\,- \mr {\bar\mu} \sum_{[\balpha] = 1, \balpha \in \mathcal{I}_*}^{2(n+1)} \sum_{l=1}^{\check I_{n+1}^{(\balpha)}}
            \check h^{(\balpha, l)}_{n+1} \check \eta^{(\balpha, l)}_{n+1} \mr {\bar D}_t \mr {\bar D}^{\balpha} \bar \rho_0            \label{e:105:SolAspt} \\ 
      &\,+ \kappa \mr{\bar \lambda}^2 \sum_{[\balpha] = 1, \balpha \in \mathcal{I}_*}^{2(n+1)} \sum_{l=1}^{\check I_{n+1}^{(\balpha)}}
            \check h^{(\balpha, l)}_{n+1} \check \eta^{(\balpha, l)}_{n+1} \mr {\bar \partial}_{ii} \mr {\bar D}^{\balpha} \bar \rho_0            \label{e:106:SolAspt} \\ 
      &\,+ 2\kappa \mr{\bar \lambda}^2 \sum_{[\balpha] = 1, \balpha \in \mathcal{I}_*}^{2(n+1)} \sum_{l=1}^{\check I_{n+1}^{(\balpha)}}
            \mr {\bar \partial}_{i} \check h^{(\balpha, l)}_{n+1} \check \eta^{(\balpha, l)}_{n+1} \mr {\bar \partial}_{i} \mr {\bar D}^{\balpha} \bar \rho_0.       \label{e:107:SolAspt} 
\end{align}

We solve \eqref{e:78:SolAspt}, and we write the solution in the form of \eqref{e:86:SolAspt}. For $a=0$ and fixed $\bbeta$, we renumber the contributions from (\ref{e:90:SolAspt}-\ref{e:107:SolAspt}) with index $1 \leq l \leq \check I^{(\bbeta)}_{n+2}$. There are differentiation indices $\bbeta$ for which there are no correctors. In this case, we write $\check I_{n+2}^{(\bbeta)} = 0$ and set all corrector components to zero.  In the end, we also estimate the cardinality of correctors, denoted by $\check I^{(\bbeta)}_{n+2}$.

\begin{case}      \label{c:compCheck:0}
If $2 \leq [\bbeta] \leq 2n+4$, $\bbeta = 1 \balpha$, this contribution from $L_{11} \hat \rho_{n+2}$, i.e. \eqref{e:91:SolAspt}, gives at most $\hat I^{(\balpha)}_{n+2}$ correctors with 
\begin{align*}
      \check h^{(\bbeta, \cdot)}_{n+2} =&\, - \frac{\delta^{\sfrac12} \mr{\bar \lambda}} {\mu} \hat h^{(\balpha, l)}_{n+2} \cdot \langle \partial_{\xi_2} \Pi_2 \hat \chi^{(\balpha, l)}_{n+2} \rangle_\xi,      \\  
      \check \eta^{(\bbeta, \cdot)}_{n+2} =&\, \partial_\tau^{-1} \big( \eta_2 \hat \eta^{(\balpha, l)}_{n+2} - \langle \eta_2 \hat \eta^{(\balpha, l)}_{n+2} \rangle_{\tau} \big).  
\end{align*}
\eqref{e:90:SolAspt} gives the same number of correctors involving $\eta_1$. The estimate \eqref{e:EXP1:ch_Pest} with $k=n+1$ follows from \eqref{e:EXP1:hh_Pest} and \eqref{e:EXP1:hchi_est} with $k=n+1$. The estimate \eqref{e:EXP1:ceta_est} with $k=n+1$ follows from \eqref{e:EXP1:heta_est} with $k=n+1$. 
\end{case}

\begin{case}      \label{c:compCheck:1}
If $1 \leq [\bbeta] \leq 2n+1$, $\bbeta = \balpha$, the contribution \eqref{e:92:SolAspt} from $L_{02} \hat \rho_{n+1}$ gives at most $N \hat I^{(\balpha)}_{n+1}$ correctors with 
\begin{align*}
      \check h^{(\bbeta, \cdot)}_{n+2} =&\, \varepsilon \frac{\delta^{\sfrac12} \mr{\bar \lambda}^2} {\mu} \mr {\bar \partial}_{i} \Omega_{m,ij} \mr {\bar \partial}_{j} \hat h^{(\balpha, l)}_{n+1} \cdot \langle \partial_{\xi_1} \Pi_1 \hat \chi^{(\balpha, l)}_{n+1} \rangle_\xi,     \\ 
      \check \eta^{(\bbeta, \cdot)}_{n+2} =&\, \partial_\tau^{-1} \big( \sigma_m \eta_1 \hat\eta^{(\balpha, l)}_{n+1} - \langle \sigma_m \eta_1 \hat\eta^{(\balpha, l)}_{n+1} \rangle_\tau \big).
\end{align*}
\eqref{e:94:SolAspt} gives the same number of correctors involving $\eta_2$. The estimate \eqref{e:EXP1:ch_Pest} with $k=n+1$ follows from \eqref{d:eddyDiff}, \eqref{e:6:homParaI}, \eqref{e:mtrEstHom}, \eqref{e:EXP1:hh_Pest} and \eqref{e:EXP1:hchi_est} with $k=n$. In particular, we just need
\begin{align*}
      \varepsilon < 1, \quad 
      \frac{\delta^{\sfrac12}}{\kappa\lambda} > 1, \quad 
      \gamma_S > 0.
\end{align*}
The estimate \eqref{e:EXP1:ceta_est} with $k=n+1$ follows from \eqref{e:EXP1:heta_est} with $k=n$ and \eqref{e:prdEstHom}.
\end{case}

\begin{case}      \label{c:compCheck:2}
If $2 \leq [\bbeta] \leq 2n+2$, $\bbeta = 1 \balpha$, the contribution \eqref{e:93:SolAspt} from $L_{02} \hat \rho_{n+1}$ gives at most $N \hat I^{(\balpha)}_{n+1}$ correctors with
\begin{align*}
      \check h^{(\bbeta, \cdot)}_{n+2} =&\, \varepsilon \frac{\delta^{\sfrac12} \mr{\bar \lambda}^2} {\mu} \mr {\bar \partial}_{i} \Omega_{m,i1} \hat h^{(\balpha, l)}_{n+1} 
            \cdot \langle \partial_{\xi_1} \Pi_1 \hat \chi^{(\balpha, l)}_{n+1} \rangle_\xi,     \\ 
      \check \eta^{(\bbeta, \cdot)}_{n+2} =&\, \partial_\tau^{-1} \big( \sigma_m \eta_1 \hat\eta^{(\balpha, l)}_{n+1} - \langle \sigma_m \eta_1 \hat\eta^{(\balpha, l)}_{n+1} \rangle_\tau \big).
\end{align*}
\eqref{e:95:SolAspt} gives the same number of correctors involving $\eta_2$. The estimate \eqref{e:EXP1:ch_Pest} with $k=n+1$ follows in a similar way to Case \ref{c:compCheck:1}. The proof of \eqref{e:EXP1:ceta_est} with $k=n+1$ is also similar to Case \ref{c:compCheck:1}.
\end{case}

\begin{case}      \label{c:compCheck:3}
If $1 \leq [\bbeta] \leq 2n+1$, $\bbeta = \balpha$, the contribution \eqref{e:97:SolAspt} from $L_{03} \hat \rho_{n+1}$ and \eqref{e:101:SolAspt} from $L_{04} \hat \rho_{n+1}$ give at most $N \hat I^{(\balpha)}_{n+1}$ correctors with 
\begin{align*}
      \check h^{(\bbeta, \cdot)}_{n+2} =&\, \varepsilon \frac{ \delta^{\sfrac12} \lambda \mr{\bar \lambda} }{\mu} \mr {\bar \partial}_{j} \big( E_{m,j2} \hat h^{(\balpha, l)}_{n+1} \big)
            \cdot \langle \partial_{\xi_1} \Pi_1 \hat \chi^{(\balpha, l)}_{n+1} \rangle_\xi,     \\ 
      \check \eta^{(\bbeta, \cdot)}_{n+2} =&\, \partial_\tau^{-1} \big( \varphi_m \eta_1 \hat\eta^{(\balpha, l)}_{n+1} - \langle \varphi_m \eta_1 \hat\eta^{(\balpha, l)}_{n+1} \rangle_\tau \big).
\end{align*}
\eqref{e:99:SolAspt} and \eqref{e:102:SolAspt} give the same number of correctors involving $\eta_2$. The estimate \eqref{e:EXP1:ch_Pest} with $k=n+1$ follows from \eqref{d:lambda_gamma}, \eqref{e:mtrEstHom}, \eqref{e:EXP1:hh_Pest} and \eqref{e:EXP1:hchi_est} with $k=n$. In particular, here we need
\begin{align}
      \varepsilon \lambda \leq \frac{\delta^{\sfrac12} \mr{\bar \lambda}}{\kappa\lambda} \cdot \lambda^{\gamma_S}.
\end{align}
The proof of \eqref{e:EXP1:ceta_est} with $k=n+1$ is also similar to Case \ref{c:compCheck:1}.
\end{case}

\begin{case}      \label{c:compCheck:4}
If $2 \leq [\bbeta] \leq 2n+2$, $\bbeta = 1 \balpha$, the contribution \eqref{e:98:SolAspt} from $L_{03} \hat \rho_{n+1}$ gives at most $N \hat I^{(\balpha)}_{n+1}$ correctors with 
\begin{align*}
      \check h^{(\bbeta, \cdot)}_{n+2} =&\, \varepsilon \frac{ \delta^{\sfrac12} \lambda \mr{\bar \lambda} }{\mu}
            E_{m,12} \hat h^{(\balpha, l)}_{n+1} \cdot \langle \partial_{\xi_1} \Pi_1 \hat \chi^{(\balpha, l)}_{n+1} \rangle_\xi,    \\ 
      \check \eta^{(\bbeta, \cdot)}_{n+2} =&\, \partial_\tau^{-1} \big( \varphi_m \eta_1 \hat\eta^{(\balpha, l)}_{n+1} - \langle \varphi_m \eta_1 \hat\eta^{(\balpha, l)}_{n+1} \rangle_\tau \big).
\end{align*}
\eqref{e:100:SolAspt} gives the same number of correctors involving $\eta_2$. The estimates \eqref{e:EXP1:ch_Pest} and \eqref{e:EXP1:ceta_est} with $k=n+1$ follow similarly to Case \ref{c:compCheck:3}.
\end{case}

\begin{case}      \label{c:compCheck:5} 
If $1 \leq [\bbeta] \leq 2n+2$, $\bbeta = \balpha$, the contribution (\ref{e:103:SolAspt}-\ref{e:104:SolAspt}) from $L_{01} \check \rho_{n+1}$ give at most $\check I^{(\balpha)}_{n+1}$ correctors with 
\begin{align*}
      \check h^{(\bbeta, \cdot)}_{n+2} =&\, \frac{\lambda}{\mu} \bigg( - \mr{\bar\mu} \mr {\bar D}_t \check h^{(\balpha, l)}_{n+1} + \kappa \mr{\bar \lambda}^2 \mr {\bar \partial}_{ii} \check h^{(\balpha, l)}_{n+1} \bigg),      \quad 
      \check \eta^{(\bbeta, \cdot)}_{n+2} = \partial_\tau^{-1} \check\eta^{(\balpha, l)}_{n+1}. 
\end{align*}
The estimate \eqref{e:EXP1:ch_Pest} with $k=n+1$ follows from \eqref{d:lambda_gamma} and \eqref{e:EXP1:ch_Pest} with $k=n$. Here, we need 
\begin{align}
      \frac{ \mr{\bar\mu}\lambda } {\mu} + \frac{ \kappa \mr{\bar \lambda}^2\lambda } {\mu} 
      \leq \frac{\delta^{\sfrac12} \mr{\bar \lambda}}{\kappa\lambda} \cdot \lambda^{\gamma_S}.
\end{align}
\eqref{e:EXP1:ceta_est} with $k=n+1$ follows from \eqref{e:EXP1:ceta_est} with $k=n$.
\end{case}

\begin{case}      \label{c:compCheck:6} 
If $2 \leq [\bbeta] \leq 2n+3$, $\bbeta = t \balpha$, the contribution \eqref{e:105:SolAspt} from $L_{01} \check \rho_{n+1}$ gives at most $\check I^{(\balpha)}_{n+1}$ correctors with 
\begin{align*}
      \check h^{(\bbeta, \cdot)}_{n+2} =&\, - \frac{ \mr{\bar\mu}\lambda } {\mu} \check h^{(\balpha, l)}_{n+1},  \quad 
      \check \eta^{(\bbeta, \cdot)}_{n+2} = \partial_\tau^{-1} \check\eta^{(\balpha, l)}_{n+1}. 
\end{align*}
The proofs of estimates \eqref{e:EXP1:ch_Pest} and \eqref{e:EXP1:ceta_est} with $k=n+1$ are similar to Case \ref{c:compCheck:5}.
\end{case}

\begin{case}      \label{c:compCheck:7} 
If $2 \leq [\bbeta] \leq 2n+3$, $\bbeta = 1 \balpha$, the contribution \eqref{e:107:SolAspt} from $L_{01} \check \rho_{n+1}$ gives at most $\check I^{(\balpha)}_{n+1}$ correctors with 
\begin{align*}
      \check h^{(\bbeta, \cdot)}_{n+2} =&\, 2 \frac{ \kappa\lambda\mr{\bar \lambda}^2 } {\mu} \mr {\bar \partial}_{1} \check h^{(\balpha, l)}_{n+1},  \quad 
      \check \eta^{(\bbeta, \cdot)}_{n+2} = \partial_\tau^{-1} \check\eta^{(\balpha, l)}_{n+1}. 
\end{align*}
\eqref{e:107:SolAspt} also gives the same number of correctors involving $\mr {\bar \partial}_{2}$. The proofs of estimates \eqref{e:EXP1:ch_Pest} and \eqref{e:EXP1:ceta_est} with $k=n+1$ are similar to Case \ref{c:compCheck:5}.
\end{case}

\begin{case}      \label{c:compCheck:8}
If $3 \leq [\bbeta] \leq 2n+4$, $\bbeta = 11 \balpha$, the contribution \eqref{e:106:SolAspt} from $L_{01} \check \rho_{n+1}$ gives at most $\check I^{(\balpha)}_{n+1}$ correctors with 
\begin{align*}
      \check h^{(\bbeta, \cdot)}_{n+2} =&\, \frac{ \kappa\lambda\mr{\bar \lambda}^2 }{\mu} \check h^{(\balpha, l)}_{n+1},      \quad 
      \check \eta^{(\bbeta, \cdot)}_{n+2} = \partial_\tau^{-1} \check\eta^{(\balpha, l)}_{n+1}. 
\end{align*}
\eqref{e:106:SolAspt} also gives the same number of correctors involving $\mr {\bar \partial}_{22}$. The proofs of estimates \eqref{e:EXP1:ch_Pest} and \eqref{e:EXP1:ceta_est} with $k=n+1$ are similar to Case \ref{c:compCheck:5}.
\end{case}

In all cases above, the information \eqref{e:EXP1:ch} for $n+2$ follows from \eqref{e:mtrStrHom}, Remark \ref{r:opPolynomial}, \eqref{e:EXP1:hh} for $n+2$ and \eqref{e:EXP1:ch} for $n+1$. \eqref{e:EXP1:chS_Pest} follows from \eqref{e:EXP1:ch_Pest} for all admissible $\balpha$ and $l$.

Now, we estimate the number of correctors $\check I_{n+2}^{(\bbeta)}$. We consider two cases. 
\begin{case}      \label{c:numCheck:0}
For index with form $\bbeta = 1 \balpha = 11\bgamma$, we have
\begin{equation}     \label{e:115:SolAspt}
\begin{split}
      \check I_{n+2}^{(\bbeta)} \leq&\, 
      \hat I^{(\balpha)}_{n+2} + 2 N \hat I^{(\bbeta)}_{n+1} + N \hat I^{(\balpha)}_{n+1} + 2 N \hat I^{(\bbeta)}_{n+1}       \\ 
      &\,+ N \hat I^{(\balpha)}_{n+1} + \check I^{(\bbeta)}_{n+1} + \check I^{(\balpha)}_{n+1} + \check I^{(\bgamma)}_{n+1}
\end{split}
\end{equation}
The contributions in the right hand side of \eqref{e:115:SolAspt} are respectively from Case \ref{c:compCheck:0}, Case \ref{c:compCheck:1}, Case \ref{c:compCheck:2}, Case \ref{c:compCheck:3}, Case \ref{c:compCheck:4}, Case \ref{c:compCheck:5}, Case \ref{c:compCheck:7} and Case \ref{c:compCheck:8}. The cases with forms 
\begin{align*}
      \bbeta = 1 \balpha = 12\bgamma,     \quad 
      \bbeta = 2 \balpha = 21\bgamma,     \quad 
      \bbeta = 2 \balpha = 22\bgamma
\end{align*}
are symmetric, hence the same estimate holds.
\end{case}

\begin{case}      \label{c:numCheck:1}
For index with form $\bbeta = t \balpha$, we have
\begin{align}     \label{e:116:SolAspt}
      \check I_{n+2}^{(\bbeta)} \leq 
      2 N \hat I^{(\bbeta)}_{n+1} + 2 N \hat I^{(\bbeta)}_{n+1} + \check I^{(\bbeta)}_{n+1} + \check I^{(\balpha)}_{n+1}. 
\end{align}
The contributions in the right hand side of \eqref{e:116:SolAspt} are respectively from Case \ref{c:compCheck:1}, Case \ref{c:compCheck:3}, Case \ref{c:compCheck:5} and Case \ref{c:compCheck:6}.
\end{case}

In a similar way, all cases other than Case \ref{c:numCheck:0} and Case \ref{c:numCheck:1} obeys the estimate \eqref{e:115:SolAspt}. Therefore, for $\bbeta, \balpha, \bgamma$ with $[\bbeta] = [\balpha]+1 = [\bgamma]+2$, we have
\begin{align}     \label{e:117:SolAspt}
      \check I_{n+2}^{(\bbeta)} \leq 
      \hat I^{(\balpha)}_{n+2} + 4 N \hat I^{(\bbeta)}_{n+1} + 2 N \hat I^{(\balpha)}_{n+1} + \check I^{(\bbeta)}_{n+1} + \check I^{(\balpha)}_{n+1} + \check I^{(\bgamma)}_{n+1}.
\end{align}
Therefore, we can verify \eqref{e:EXP1:cI_est} for $k=n+1$ from \eqref{e:EXP1:hI_est} with $k \leq n+1$ and \eqref{e:EXP1:cI_est} with $k \leq n$.

\underline{Now we turn to the terms in \eqref{e:81:6:SolAspt} and \eqref{e:81:8:SolAspt}}, i.e. the contribution of $F$-related terms in \eqref{e:86:SolAspt}. Comparing the structures of the lines (\ref{e:81:6:SolAspt}-\ref{e:81:8:SolAspt}) and the lines (\ref{e:84:SolAspt}-\ref{e:85:SolAspt}), it is clear that their computations are analogous.

\subsubsection{Inductive step for spatial correctors}      \label{sss:solAsptSPT}

\underline{Finally, we derive the structures of $\hat \rho_{n+3}$.} We take the difference of $\langle \eqref{e:aspt:hrc4} \rangle_\xi$ and $\eqref{e:aspt:hrc4}$. Then we apply Remark \ref{r:hatcheckbar} to deduce
\begin{equation}  \label{e:118:SolAspt}
\begin{split}
      L_2 \hat\rho_{n+3} =&\, \langle L_{11} \hat\rho_{n+2} \rangle_\xi - L_{11} \hat\rho_{n+2}
      - L_{11} ( \check\rho_{n+2} + \bar\rho_{n+2} ) - ( L_{12} + L_{13} ) \hat\rho_{n+2} \\ 
      &\, - L_{01} \hat\rho_{n+1} + \langle L_{02} \hat\rho_{n+1} \rangle_\xi - L_{02} \hat\rho_{n+1}
      + \langle L_{03} \hat\rho_{n+1} \rangle_\xi - L_{03} \hat\rho_{n+1} \\ 
      &\, + \langle L_{04} \hat\rho_{n+1} \rangle_\xi - L_{04} \hat\rho_{n+1} 
      - ( L_{02} + L_{03} ) ( \check \rho_{n+1} + \bar \rho_{n+1} ) \\ 
      &\, - L_{05} \hat\rho_{n+1} - L_{06} \hat\rho_{n+1} - L_{07} \hat\rho_{n+1} - L_{-1} \hat \rho_n.
\end{split}
\end{equation}
The goal is to compute the coefficients in 
\begin{equation}     \label{e:119:SolAspt}
\begin{split}
      \hat \rho_{n+3} (x,t,\xi,\tau) =&\, \sum_{a = 0}^{n+2} \sum_{[\bbeta] = 1, \bbeta \in \mathcal{I}_*}^{2(n-a)+5} \sum_{l=1}^{ \hat I_{n+3,a}^{(\bbeta)} }
            \hat h^{(\bbeta, l)}_{n+3,a} (x,t) \hat \eta^{(\bbeta, l)}_{n+3,a} (\tau) \hat \chi^{(\bbeta, l)}_{n+3,a} (\chi) \mr {\bar D}^{\bbeta} \bar \rho_a (x,t) \\ 
      +&\, \sum_{[\bbeta]=n_0, \bbeta \in \mathcal{I}_*}^{2n+4+n_0} \sum_{l=1}^{\hat I^{(\bbeta)}_{F,n+3}}
            \hat h^{(\bbeta,l)}_{F,n+3} (x,t) \hat \eta^{(\bbeta,l)}_{F,n+3} (\tau) \hat \chi^{(\bbeta,l)}_{F,n+3} (\xi) \mr {\bar D}^{\bbeta} \bar \rho_F (x,t). 
\end{split}
\end{equation}

Similar to $\check \rho_{n+2}$, we can prove the stationary property in Remark \ref{r:stationarity}. 
\begin{proof}[Proof of stationarity for $\hat \rho_{n+3}$]
This proof has a similar nature to the analogue for $\check \rho_{n+2}$. We consider two cases $n=0$ and $n \geq 1$. In the following computations, we use the decomposition similar to \eqref{e:1:SolAspt}, \eqref{e:4:2:SolAspt} and \eqref{e:4:6:SolAspt}, that $\hat \rho_{\cdot,E}$ contains the correctors involving $\bar \rho_a$ with $a \geq 1$, and $\hat \rho_{\cdot,0}$ contains the correctors involving $\bar \rho_0$.
\begin{case}      \label{c:0:hatn+3_sta}
If $n=0$, similar to \eqref{e:118:SolAspt}, we can derive from \eqref{e:aspt:hrc3} and \eqref{e:aspt:hrc4}, that
\begin{align}
      L_2 \hat \rho_{2,R} =&\, \langle L_{11} \hat\rho_{1,R} \rangle_{\xi} - L_{11} \hat\rho_{1,R} - L_{11} ( \check\rho_{1,R} + \bar\rho_1 ) - ( L_{12} + L_{13} ) \hat\rho_{1,R}  - ( L_{02} + L_{03} ) \bar\rho_0,   \label{e:120:SolAspt} \\ 
      L_2 \hat \rho_{2,F} =&\, \langle L_{11} \hat\rho_{1,F} \rangle_{\xi} - L_{11} \hat\rho_{1,F} - L_{11} \check\rho_{1,F} - ( L_{12} + L_{13} ) \hat\rho_{1,F} ,   \label{e:120:2:SolAspt} \\ 
      L_2 \hat \rho_{3,R} =&\, \langle L_{11} \hat\rho_{2,E} \rangle_\xi - L_{11} \hat\rho_{2,E}
      - L_{11} ( \check\rho_{2,E} + \bar\rho_{2} ) - ( L_{12} + L_{13} ) \hat\rho_{2,E} - ( L_{02} + L_{03} ) \bar \rho_{1}   \label{e:121:SolAspt} \\ 
      &\,+ \langle L_{11} \hat\rho_{2,0} \rangle_\xi - L_{11} \hat\rho_{2,0} - L_{11} \check\rho_{2,0} - ( L_{12} + L_{13} ) \hat\rho_{2,0}    \label{e:122:SolAspt} \\ 
      &\,- L_{01} \hat\rho_{1} + \langle ( L_{02} + L_{03} + L_{04} ) \hat\rho_{1} \rangle_\xi - ( L_{02} + L_{03} + L_{04} ) \hat\rho_{1}    \label{e:122+:SolAspt} \\
      &\,- ( L_{02} + L_{03} ) \check \rho_{1} - (L_{05}+L_{06}+L_{07}) \hat\rho_{1}.      \label{e:123:SolAspt} \\ 
      L_2 \hat \rho_{3,F} =&\, \langle L_{11} \hat\rho_{2,F} \rangle_{\xi} - L_{11} \hat\rho_{2,F} - L_{11} \check\rho_{2,F} - ( L_{12} + L_{13} ) \hat\rho_{2,F} - ( L_{02} + L_{03} ) \check \rho_{1,F}.   \label{e:123:2:SolAspt} 
\end{align}
Here, we use the decomposition $\hat\rho_{\cdot} = \hat\rho_{\cdot,R} + \hat\rho_{\cdot,F}$, in which $\hat\rho_{\cdot,R}$ contains all terms related to $\bar\rho_a, a \geq 0$ and $\hat\rho_{\cdot,F}$ contains all $F$-related terms. (\ref{e:122:SolAspt}-\ref{e:123:SolAspt}) contain the terms in \eqref{e:119:SolAspt} with $a=0$. Compare \eqref{e:120:SolAspt} and \eqref{e:121:SolAspt} and use the stationarity for the coefficients in $\hat\rho_1, \hat\rho_2$ and $\check \rho_1$. We can conclude the stationarity property for correctors in $\hat\rho_3$.
\end{case}

\begin{case}      \label{c:1:hatn+3_sta}
If $n \geq 1$, we consider \eqref{e:aspt:hrc4} for $n-1$ and $n$. Similar to above, we obtain
\begin{align}
      L_2 \hat\rho_{n+2,R} =&\, \langle L_{11} \hat\rho_{n+1,R} \rangle_\xi - L_{11} \hat\rho_{n+1,R}
      - L_{11} ( \check\rho_{n+1,R} + \bar\rho_{n+1} ) - ( L_{12} + L_{13} ) \hat\rho_{n+1,R}       \label{e:124:SolAspt} \\ 
      -&\, L_{01} \hat\rho_{n,R} + \langle ( L_{02} + L_{03} + L_{04} ) \hat\rho_{n,R} \rangle_\xi - ( L_{02} + L_{03} + L_{04} ) \hat\rho_{n,R}   \label{e:125:SolAspt} \\ 
      -&\, ( L_{02} + L_{03} ) ( \check \rho_{n,R} + \bar \rho_{n} ) - (L_{05}+L_{06}+L_{07}) \hat\rho_{n,R} - L_{-1} \hat \rho_{n-1,R},  \label{e:126:SolAspt} \\
      L_2 \hat\rho_{n+3,R} =&\, \langle L_{11} \hat\rho_{n+2,E} \rangle_\xi - L_{11} \hat\rho_{n+2,E} - L_{11} ( \check\rho_{n+2,E} + \bar\rho_{n+2} ) - ( L_{12} + L_{13} ) \hat\rho_{n+2,E}    \label{e:127:SolAspt} \\ 
      -&\, L_{01} \hat\rho_{n+1,E} + \langle ( L_{02} + L_{03} + L_{04} ) \hat\rho_{n+1,E} \rangle_\xi - ( L_{02} + L_{03} + L_{04} ) \hat\rho_{n+1,E}  \label{e:128:SolAspt}  \\ 
      -&\, ( L_{02} + L_{03} ) ( \check \rho_{n+1,E} + \bar \rho_{n+1} ) - (L_{05}+L_{06}+L_{07}) \hat\rho_{n+1,E} - L_{-1} \hat \rho_{n,E}  \label{e:129:SolAspt} \\ 
      +&\, \langle L_{11} \hat\rho_{n+2,0} \rangle_\xi - L_{11} \hat\rho_{n+2,0}
      - L_{11} \check\rho_{n+2,0} - ( L_{12} + L_{13} ) \hat\rho_{n+2,0}      \label{e:130:SolAspt} \\ 
      -&\, L_{01} \hat\rho_{n+1,0} + \langle ( L_{02} + L_{03} + L_{04} ) \hat\rho_{n+1,0} \rangle_\xi - ( L_{02} + L_{03} + L_{04} ) \hat\rho_{n+1,0}  \label{e:131:SolAspt} \\ 
      -&\, ( L_{02} + L_{03} ) \check \rho_{n+1,0} - (L_{05}+L_{06}+L_{07}) \hat\rho_{n+1,0} - L_{-1} \hat \rho_{n,0}.      \label{e:132:SolAspt}
\end{align}
We then compare (\ref{e:124:SolAspt}-\ref{e:126:SolAspt}) and (\ref{e:127:SolAspt}-\ref{e:129:SolAspt}). As usual, we conclude the stationarity property for $\hat\rho_{n+3}$ from the stationarity for $\hat\rho_{n+2}, \hat\rho_{n+1}, \hat\rho_{n}, \hat\rho_{n-1}, \check\rho_{n+1}, \check\rho_{n}$ in the induction assumptions and the stationarity of $\check\rho_{n+2}$ proved before.
\end{case}

\end{proof}

With the stationarity verified above, we compute the coefficients in \eqref{e:119:SolAspt} for terms involving $\bar \rho_0$. Define
\begin{align*}
      \hat I_{n-a+3}^{(\bbeta)} := \hat I_{n+3,a}^{(\bbeta)}, \quad 
      \hat h^{(\bbeta, l)}_{n-a+3} := \hat h^{(\bbeta, l)}_{n+3,a}, \quad 
      \hat \eta^{(\bbeta, l)}_{n-a+3} := \hat \eta^{(\bbeta, l)}_{n+3,a}, \quad 
      \hat \chi^{(\bbeta, l)}_{n-a+3} := \hat \chi^{(\bbeta, l)}_{n+3,a}. 
\end{align*}
Therefore, similar to $\check \rho_{n+2}$, we can write \eqref{e:119:SolAspt} as
\begin{equation}        \label{e:135:SolAspt} 
\begin{split}
      \hat \rho_{n+3}(x,t,\xi,\tau) =&\, \sum_{a=0}^{n+2} \sum_{[\bbeta] = 1, \bbeta \in \mathcal{I}_*}^{2(n-a)+5} \sum_{l=1}^{ \hat I_{n-a+3}^{(\bbeta)} }
            \hat h^{(\bbeta, l)}_{n-a+3} (x,t) \hat \eta^{(\bbeta, l)}_{n-a+3} (\tau) \hat \chi^{(\bbeta, l)}_{n-a+3} (\xi) \mr {\bar D}^{\bbeta} \bar \rho_a (x,t)      \\ 
      +&\, \sum_{[\bbeta]=n_0, \bbeta \in \mathcal{I}_*}^{2n+4+n_0} \sum_{l=1}^{\hat I^{(\bbeta)}_{F,n+3}}
            \hat h^{(\bbeta,l)}_{F,n+3} (x,t) \hat \eta^{(\bbeta,l)}_{F,n+3} (\tau) \hat \chi^{(\bbeta,l)}_{F,n+3} (\xi) \mr {\bar D}^{\bbeta} \bar \rho_F (x,t).
\end{split}
\end{equation}

Now we compute the terms in (\ref{e:122:SolAspt}-\ref{e:123:SolAspt}) or (\ref{e:130:SolAspt}-\ref{e:132:SolAspt}), the coefficient fields for $\bar\rho_0$ in \eqref{e:135:SolAspt}. There are differentiation indices $\bbeta$ for which there is no corrector. In this case, we set $\hat I_{n+3}^{(\bbeta)} = 0$ and all corresponding coefficient fields to zero.

For $\langle L_{11} \hat\rho_{n+2,0} \rangle_{\xi} - L_{11} \hat\rho_{n+2,0}$, we have
\begin{align}
      \langle L_{11} \hat\rho_{n+2,0} \rangle_{\xi} -&\, L_{11} \hat\rho_{n+2,0} 
            = \Bigg( \sum_{[\balpha] = 1, \balpha \in \mathcal{I}_*}^{2n+3} \sum_{l=1}^{ \hat I_{n+2}^{(\balpha)} } \, \text{below} \Bigg)   \nonumber \\ 
      =&\, \frac{\delta^{\sfrac12}\mr{\bar \lambda}} {\lambda} \sum\sum 
            \hat h^{(\balpha, l)}_{n+2} \eta_1 \hat \eta^{(\balpha, l)}_{n+2} \big( \partial_{\xi_1} \Pi_1 \hat \chi^{(\balpha, l)}_{n+2} - \langle \partial_{\xi_1} \Pi_1 \hat \chi^{(\balpha, l)}_{n+2} \rangle_\xi \big) \mr {\bar \partial}_{2} \mr {\bar D}^{\balpha} \bar \rho_0      \label{e:136:SolAspt} \\ 
      +&\, \frac{\delta^{\sfrac12}\mr{\bar \lambda}} {\lambda} \sum\sum 
            \mr {\bar \partial}_{2} \hat h^{(\balpha, l)}_{n+2} \eta_1 \hat \eta^{(\balpha, l)}_{n+2} \big( \partial_{\xi_1} \Pi_1 \hat \chi^{(\balpha, l)}_{n+2} - \langle \partial_{\xi_1} \Pi_1 \hat \chi^{(\balpha, l)}_{n+2} \rangle_\xi \big) \mr {\bar D}^{\balpha} \bar \rho_0      \label{e:137:SolAspt} \\ 
      -&\, \frac{\delta^{\sfrac12}\mr{\bar \lambda}} {\lambda} \sum\sum 
            \hat h^{(\balpha, l)}_{n+2} \eta_2 \hat \eta^{(\balpha, l)}_{n+2} \big( \partial_{\xi_2} \Pi_2 \hat \chi^{(\balpha, l)}_{n+2} - \langle \partial_{\xi_2} \Pi_2 \hat \chi^{(\balpha, l)}_{n+2} \rangle_\xi \big) \mr {\bar \partial}_{1} \mr {\bar D}^{\balpha} \bar \rho_0      \label{e:138:SolAspt} \\ 
      -&\, \frac{\delta^{\sfrac12}\mr{\bar \lambda}} {\lambda} \sum\sum 
            \mr {\bar \partial}_{1} \hat h^{(\balpha, l)}_{n+2} \eta_2 \hat \eta^{(\balpha, l)}_{n+2} \big( \partial_{\xi_2} \Pi_2 \hat \chi^{(\balpha, l)}_{n+2} - \langle \partial_{\xi_2} \Pi_2 \hat \chi^{(\balpha, l)}_{n+2} \rangle_\xi \big) \mr {\bar D}^{\balpha} \bar \rho_0.     \label{e:139:SolAspt}
\end{align}

For the term $- L_{11} \check\rho_{n+2,0}$, we have
\begin{align}
      - L_{11} \check\rho_{n+2,0} &\,= \Bigg( \sum_{[\balpha] = 1, \balpha \in \mathcal{I}_*}^{2n+4} \sum_{l=1}^{ \check I_{n+2}^{(\balpha)} } \, \text{below} \Bigg)    \nonumber \\ 
      &\,= \frac{\delta^{\sfrac12}\mr{\bar \lambda}} {\lambda} \sum\sum 
            \check h^{(\balpha, l)}_{n+2} \eta_1 \check \eta^{(\balpha, l)}_{n+2} \partial_{\xi_1} \Pi_1 \mr {\bar \partial}_{2} \mr {\bar D}^{\balpha} \bar \rho_0      \label{e:140:SolAspt} \\ 
      &\,+ \frac{\delta^{\sfrac12}\mr{\bar \lambda}} {\lambda} \sum\sum 
            \mr {\bar \partial}_{2} \check h^{(\balpha, l)}_{n+2} \eta_1 \check \eta^{(\balpha, l)}_{n+2} \partial_{\xi_1} \Pi_1 \mr {\bar D}^{\balpha} \bar \rho_0     \label{e:141:SolAspt} \\ 
      &\,- \frac{\delta^{\sfrac12}\mr{\bar \lambda}} {\lambda} \sum\sum 
            \check h^{(\balpha, l)}_{n+2} \eta_2 \check \eta^{(\balpha, l)}_{n+2} \partial_{\xi_2} \Pi_2 \mr {\bar \partial}_{1} \mr {\bar D}^{\balpha} \bar \rho_0  \label{e:142:SolAspt} \\ 
      &\,- \frac{\delta^{\sfrac12}\mr{\bar \lambda}} {\lambda} \sum\sum 
            \mr {\bar \partial}_{1} \check h^{(\balpha, l)}_{n+2} \eta_2 \check \eta^{(\balpha, l)}_{n+2} \partial_{\xi_2} \Pi_2 \mr {\bar D}^{\balpha} \bar \rho_0.       \label{e:143:SolAspt}
\end{align}

For the term $-L_{12} \hat \rho_{n+2,0}$ and $-L_{13} \hat \rho_{n+2,0}$, we have
\begin{align}
      -L_{12} \hat \rho_{n+2,0} =&\, - \frac{\mu}{\lambda} 
            \sum_{[\balpha]=1, \balpha \in \mathcal{I}_*}^{2n+3} \sum_{l=1}^{\hat I_{n+2}^{(\balpha)}} 
            \hat h^{(\balpha, l)}_{n+2} \partial_\tau \hat \eta^{(\balpha, l)}_{n+2} \hat \chi^{(\balpha, l)}_{n+2} \mr {\bar D}^{\balpha} \bar \rho_0,       \label{e:144:SolAspt} \\ 
      -L_{13} \hat \rho_{n+2,0} =&\, - \varepsilon \kappa \lambda
            \sum_{[\balpha]=1, \balpha \in \mathcal{I}_*}^{2n+3} \sum_{l=1}^{\hat I_{n+2}^{(\balpha)}} \sum_m 
            S_{m,11} \hat h^{(\balpha, l)}_{n+2} \vartheta_m \hat \eta^{(\balpha, l)}_{n+2} \partial_{\xi_1\xi_1} \hat \chi^{(\balpha, l)}_{n+2} \mr {\bar D}^{\balpha} \bar \rho_0    \label{e:145:SolAspt} \\ 
      &\, - \varepsilon \kappa \lambda 
            \sum_{[\balpha]=1, \balpha \in \mathcal{I}_*}^{2n+3} \sum_{l=1}^{\hat I_{n+2}^{(\balpha)}} \sum_m 
            S_{m,22} \hat h^{(\balpha, l)}_{n+2} \vartheta_m \hat \eta^{(\balpha, l)}_{n+2} \partial_{\xi_2\xi_2} \hat \chi^{(\balpha, l)}_{n+2} \mr {\bar D}^{\balpha} \bar \rho_0.       \label{e:146:SolAspt}
\end{align}

For the term $-L_{01} \hat\rho_{n+1,0}$, we have
\begin{align}
      -L_{01} \hat\rho_{n+1,0} = 
      -&\, \mr{\bar\mu} \sum_{[\balpha]=1, \balpha \in \mathcal{I}_*}^{2n+1} \sum_{l=1}^{\hat I_{n+1}^{(\balpha)}}
            \mr{\bar D}_t \hat h^{(\balpha, l)}_{n+1} \hat \eta^{(\balpha, l)}_{n+1} \hat \chi^{(\balpha, l)}_{n+1} D^{\balpha} \bar \rho_0        \label{e:147:SolAspt} \\ 
      -&\, \mr{\bar\mu} \sum_{[\balpha]=1, \balpha \in \mathcal{I}_*}^{2n+1} \sum_{l=1}^{\hat I_{n+1}^{(\balpha)}}
            \hat h^{(\balpha, l)}_{n+1} \hat \eta^{(\balpha, l)}_{n+1} \hat \chi^{(\balpha, l)}_{n+1} \mr {\bar D}_t \mr {\bar D}^{\balpha} \bar \rho_0      \label{e:148:SolAspt} \\ 
      +&\, \kappa \mr{\bar\lambda}^2 \sum_{[\balpha]=1, \balpha \in \mathcal{I}_*}^{2n+1} \sum_{l=1}^{\hat I_{n+1}^{(\balpha)}}
            \mr {\bar \partial}_{ii} \hat h^{(\balpha, l)}_{n+1} \hat \eta^{(\balpha, l)}_{n+1} \hat \chi^{(\balpha, l)}_{n+1} \mr {\bar D}^{\balpha} \bar \rho_0      \label{e:149:SolAspt} \\ 
      +&\, \kappa \mr{\bar\lambda}^2 \sum_{[\balpha]=1, \balpha \in \mathcal{I}_*}^{2n+1} \sum_{l=1}^{\hat I_{n+1}^{(\balpha)}}
            \hat h^{(\balpha, l)}_{n+1} \hat \eta^{(\balpha, l)}_{n+1} \hat \chi^{(\balpha, l)}_{n+1} \mr {\bar \partial}_{ii} \mr {\bar D}^{\balpha} \bar \rho_0      \label{e:149+:SolAspt} \\ 
      +&\, 2 \kappa \mr{\bar\lambda}^2 \sum_{[\balpha]=1, \balpha \in \mathcal{I}_*}^{2n+1} \sum_{l=1}^{\hat I_{n+1}^{(\balpha)}}
            \mr {\bar \partial}_{i} \hat h^{(\balpha, l)}_{n+1} \hat \eta^{(\balpha, l)}_{n+1} \hat \chi^{(\balpha, l)}_{n+1} \mr {\bar \partial}_{i} \mr {\bar D}^{\balpha} \bar \rho_0.  \label{e:150:SolAspt}
\end{align}

We also have, for the term $L_{02} \hat\rho_{n+1,0}$,
\begin{align}
      \langle L_{02} \hat\rho_{n+1,0} \rangle_{\xi} -&\, L_{02} \hat\rho_{n+1,0} 
      = \Bigg( \sum_{[\balpha]=1, \balpha \in \mathcal{I}_*}^{2n+1} \sum_{l=1}^{\hat I_{n+1}^{(\balpha)}} \sum_m \text{ below} \Bigg)        \nonumber \\ 
      =&\, \varepsilon \frac{\delta^{\sfrac12}\mr{\bar \lambda}^2} {\lambda} \sum\sum\sum 
            \mr {\bar \partial}_{i} \Omega_{m,ij} \mr {\bar \partial}_{j} \hat h^{(\balpha, l)}_{n+1} \sigma_m \eta_1 \hat \eta^{(\balpha, l)}_{n+1} \big(\partial_{\xi_1} \Pi_1 \hat \chi^{(\balpha, l)}_{n+1} - \langle \cdot \rangle_\xi \big) \mr {\bar D}^{\balpha}\bar \rho_0      \label{e:151:SolAspt} \\ 
      +&\, \varepsilon \frac{\delta^{\sfrac12}\mr{\bar \lambda}^2} {\lambda} \sum\sum\sum 
            \mr {\bar \partial}_{i} \Omega_{m,ij} \hat h^{(\balpha, l)}_{n+1} \sigma_m \eta_1 \hat \eta^{(\balpha, l)}_{n+1} \big(\partial_{\xi_1} \Pi_1 \hat \chi^{(\balpha, l)}_{n+1} - \langle \cdot \rangle_\xi \big) \mr {\bar \partial}_{j} \mr {\bar D}^{\balpha}\bar \rho_0       \label{e:152:SolAspt} \\ 
      +&\,  \varepsilon \frac{\delta^{\sfrac12}\mr{\bar \lambda}^2} {\lambda} \sum\sum\sum 
            \mr {\bar \partial}_{i} \Omega_{m,ij} \mr {\bar \partial}_{j} \hat h^{(\balpha, l)}_{n+1} \sigma_m \eta_2 \hat \eta^{(\balpha, l)}_{n+1} \big(\partial_{\xi_2} \Pi_2 \hat \chi^{(\balpha, l)}_{n+1} - \langle \cdot \rangle_\xi \big) \mr {\bar D}^{\balpha}\bar \rho_0       \label{e:153:SolAspt} \\ 
      +&\, \varepsilon \frac{\delta^{\sfrac12}\mr{\bar \lambda}^2} {\lambda} \sum\sum\sum 
            \mr {\bar \partial}_{i} \Omega_{m,ij} \hat h^{(\balpha, l)}_{n+1} \sigma_m \eta_2 \hat \eta^{(\balpha, l)}_{n+1} \big(\partial_{\xi_2} \Pi_2 \hat \chi^{(\balpha, l)}_{n+1} - \langle \cdot \rangle_\xi \big) \mr {\bar \partial}_{j} \mr {\bar D}^{\balpha}\bar \rho_0.       \label{e:154:SolAspt}
\end{align}
For the term $L_{03} \hat\rho_{n+1,0}$, we have
\begin{align}
      \langle L_{03} \hat\rho_{n+1,0} \rangle_{\xi} -&\, L_{03} \hat\rho_{n+1,0} 
      = \Bigg( \sum_{[\balpha]=1, \balpha \in \mathcal{I}_*}^{2n+1} \sum_{l=1}^{\hat I_{n+1}^{(\balpha)}} \sum_m \text{ below} \Bigg)       \nonumber \\ 
      =&\, \varepsilon \delta^{\sfrac12} \mr{\bar \lambda} \sum\sum\sum 
            E_{m,j2} \mr {\bar \partial}_{j} \hat h^{(\balpha, l)}_{n+1} \varphi_m \eta_1 \hat \eta^{(\balpha, l)}_{n+1} \big( \partial_{\xi_1} \Pi_1 \hat \chi^{(\balpha, l)}_{n+1} - \langle \cdot \rangle_\xi \big) \mr {\bar D}^{\balpha}\bar \rho_0      \label{e:155:SolAspt} \\ 
      +&\, \varepsilon \delta^{\sfrac12} \mr{\bar \lambda} \sum\sum\sum 
            E_{m,j2} \hat h^{(\balpha, l)}_{n+1} \varphi_m \eta_1 \hat \eta^{(\balpha, l)}_{n+1} \big( \partial_{\xi_1} \Pi_1 \hat \chi^{(\balpha, l)}_{n+1} - \langle \cdot \rangle_\xi \big) \mr {\bar \partial}_{j} \mr {\bar D}^{\balpha}\bar \rho_0      \label{e:156:SolAspt} \\ 
      -&\, \varepsilon \delta^{\sfrac12} \mr{\bar \lambda} \sum\sum\sum 
            E_{m,j1} \mr {\bar \partial}_{j} \hat h^{(\balpha, l)}_{n+1} \varphi_m \eta_2 \hat \eta^{(\balpha, l)}_{n+1} \big( \partial_{\xi_2} \Pi_2 \hat \chi^{(\balpha, l)}_{n+1} - \langle \cdot \rangle_\xi \big) \mr {\bar D}^{\balpha}\bar \rho_0      \label{e:157:SolAspt} \\ 
      -&\, \varepsilon \delta^{\sfrac12} \mr{\bar \lambda} \sum\sum\sum 
            E_{m,j1} \hat h^{(\balpha, l)}_{n+1} \varphi_m \eta_2 \hat \eta^{(\balpha, l)}_{n+1} \big( \partial_{\xi_2} \Pi_2 \hat \chi^{(\balpha, l)}_{n+1} - \langle \cdot \rangle_\xi \big) \mr {\bar \partial}_{j} \mr {\bar D}^{\balpha}\bar \rho_0.      \label{e:158:SolAspt}
\end{align}
For the term $L_{04} \hat\rho_{n+1,0}$, we have 
\begin{align}
      \langle L_{04} \hat\rho_{n+1,0} \rangle_{\xi} -&\, L_{04} \hat\rho_{n+1,0} 
      = \Bigg( \sum_{[\balpha]=1, \balpha \in \mathcal{I}_*}^{2n+1} \sum_{l=1}^{\hat I_{n+1}^{(\balpha)}} \sum_m \text{ below} \Bigg)     \nonumber \\ 
      =-&\, \varepsilon \delta^{\sfrac12} \mr{\bar \lambda} \sum\sum\sum 
            \mr {\bar \partial}_{j} E_{m,j2} \hat h^{(\balpha, l)}_{n+1} \varphi_m \eta_1 \hat \eta^{(\balpha, l)}_{n+1} \big( \Pi_1 \partial_{\xi_1} \hat \chi^{(\balpha, l)}_{n+1} - \langle \cdot \rangle_\xi \big) \mr {\bar D}^{\balpha}\bar \rho_0      \label{e:159:SolAspt} \\ 
      +&\, \varepsilon \delta^{\sfrac12} \mr{\bar \lambda} \sum\sum\sum 
            \mr {\bar \partial}_{j} E_{m,j1} \hat h^{(\balpha, l)}_{n+1} \varphi_m \eta_2 \hat \eta^{(\balpha, l)}_{n+1} \big( \Pi_2 \partial_{\xi_2} \hat \chi^{(\balpha, l)}_{n+1} - \langle \cdot \rangle_\xi \big) \mr {\bar D}^{\balpha}\bar \rho_0.     \label{e:160:SolAspt}
\end{align}

Next, for $-L_{02} \check \rho_{n+1,0}$, we have
\begin{align}
      -L_{02} \check \rho_{n+1,0} =&\, \Bigg( \sum_{[\balpha]=1, \balpha \in \mathcal{I}_*}^{2n+2} \sum_{l=1}^{\check I_{n+1}^{(\balpha)}} \sum_m \text{ below} \Bigg)     \nonumber \\ 
      =&\, \varepsilon \frac{\delta^{\sfrac12}\mr{\bar \lambda}^2}{\lambda} \sum\sum\sum 
            \mr {\bar \partial}_{i} \Omega_{m,ij} \mr {\bar \partial}_{j} \check h^{(\balpha, l)}_{n+1} \sigma_m \check \eta^{(\balpha, l)}_{n+1} ( \eta_1 \Pi_1 + \eta_2 \Pi_2 ) \mr {\bar D}^{\balpha}\bar \rho_0     \label{e:161:SolAspt} \\ 
      +&\, \varepsilon \frac{\delta^{\sfrac12}\mr{\bar \lambda}^2}{\lambda} \sum\sum\sum 
            \mr {\bar \partial}_{i} \Omega_{m,ij} \check h^{(\balpha, l)}_{n+1} \sigma_m \check \eta^{(\balpha, l)}_{n+1} ( \eta_1 \Pi_1 + \eta_2 \Pi_2 ) \mr {\bar \partial}_{j} \mr {\bar D}^{\balpha}\bar \rho_0.     \label{e:162:SolAspt} 
\end{align}

For the term $-L_{03} \check \rho_{n+1,0}$, we have
\begin{align}
      -L_{03} \check \rho_{n+1,0} =&\, \Bigg( \sum_{[\balpha]=1, \balpha \in \mathcal{I}_*}^{2n+2} \sum_{l=1}^{\check I_{n+1}^{(\balpha)}} \sum_m \text{ below} \Bigg)        \nonumber \\ 
      =&\, \varepsilon \delta^{\sfrac12} \mr{\bar \lambda} \sum\sum\sum 
            E_{m,j2} \mr {\bar \partial}_{j} \check h^{(\balpha, l)}_{n+1} \varphi_m \check \eta^{(\balpha, l)}_{n+1} \eta_1 \partial_{\xi_1} \Pi_1 \mr {\bar D}^{\balpha}\bar \rho_0     \label{e:163:SolAspt} \\ 
      +&\, \varepsilon \delta^{\sfrac12} \mr{\bar \lambda} \sum\sum\sum 
            E_{m,j2} \check h^{(\balpha, l)}_{n+1} \varphi_m \check \eta^{(\balpha, l)}_{n+1} \eta_1 \partial_{\xi_1} \Pi_1 \mr {\bar \partial}_{j} \mr {\bar D}^{\balpha}\bar \rho_0     \label{e:164:SolAspt} \\ 
      -&\, \varepsilon \delta^{\sfrac12} \mr{\bar \lambda} \sum\sum\sum 
            E_{m,j1} \mr {\bar \partial}_{j} \check h^{(\balpha, l)}_{n+1} \varphi_m \check \eta^{(\balpha, l)}_{n+1} \eta_2 \partial_{\xi_2} \Pi_2 \mr {\bar D}^{\balpha}\bar \rho_0     \label{e:165:SolAspt}\\ 
      -&\, \varepsilon \delta^{\sfrac12} \mr{\bar \lambda} \sum\sum\sum 
            E_{m,j1} \check h^{(\balpha, l)}_{n+1} \varphi_m \check \eta^{(\balpha, l)}_{n+1} \eta_2 \partial_{\xi_2} \Pi_2 \mr {\bar \partial}_{j} \mr {\bar D}^{\balpha}\bar \rho_0.     \label{e:166:SolAspt}
\end{align}

For $- (L_{05}+L_{06}+L_{07}) \hat \rho_{n+1,0}$ and $-L_{-1} \hat\rho_{n,0}$, we have
\begin{align}
      -L_{05} \hat \rho_{n+1,0} =&\, \varepsilon \kappa \lambda \mr{\bar \lambda} \sum_{[\balpha]=1, \balpha \in \mathcal{I}_*}^{2n+1} \sum_{l=1}^{\hat I_{n+1}^{(\balpha)}} \sum_m 
      \mr {\bar \partial}_{j} B_{m,ij} \hat h^{(\balpha, l)}_{n+1} \phi_m \hat \eta^{(\balpha, l)}_{n+1} \partial_{\xi_i} \hat \chi^{(\balpha, l)}_{n+1} \mr {\bar D}^{\balpha}\bar \rho_0,     \label{e:170:SolAspt} \\ 
      -L_{06} \hat \rho_{n+1,0} =&\, 2\kappa \lambda \mr{\bar \lambda} \sum_{[\balpha]=1, \balpha \in \mathcal{I}_*}^{2n+1} \sum_{l=1}^{\hat I_{n+1}^{(\balpha)}} \sum_m  
            \mr {\bar \partial}_i \hat h^{(\balpha, l)}_{n+1} \hat \eta^{(\balpha, l)}_{n+1} \partial_{\xi_i} \hat \chi^{(\balpha, l)}_{n+1} \mr {\bar D}^{\balpha}\bar \rho_0      \label{e:171:SolAspt} \\ 
      +&\, 2\kappa \lambda \mr{\bar \lambda} \sum_{[\balpha]=1, \balpha \in \mathcal{I}_*}^{2n+1} \sum_{l=1}^{\hat I_{n+1}^{(\balpha)}} \sum_m  
            \hat h^{(\balpha, l)}_{n+1} \hat \eta^{(\balpha, l)}_{n+1} \partial_{\xi_i} \hat \chi^{(\balpha, l)}_{n+1} \mr {\bar \partial}_i \mr {\bar D}^{\balpha}\bar \rho_0      \label{e:172:SolAspt} \\ 
      -L_{07} \hat \rho_{n+1,0} =&\, 2\varepsilon \kappa \lambda \mr{\bar \lambda} \sum_{[\balpha]=1, \balpha \in \mathcal{I}_*}^{2n+1} \sum_{l=1}^{\hat I_{n+1}^{(\balpha)}} \sum_m 
            B_{m,ij} \mr {\bar \partial}_{j} \hat h^{(\balpha, l)}_{n+1} \phi_m \hat \eta^{(\balpha, l)}_{n+1} \partial_{\xi_i} \hat \chi^{(\balpha, l)}_{n+1} \mr {\bar D}^{\balpha}\bar \rho_0,     \label{e:173:SolAspt} \\ 
      +&\, 2\varepsilon \kappa \lambda \mr{\bar \lambda} \sum_{[\balpha]=1, \balpha \in \mathcal{I}_*}^{2n+1} \sum_{l=1}^{\hat I_{n+1}^{(\balpha)}} \sum_m 
            B_{m,ij} \hat h^{(\balpha, l)}_{n+1} \phi_m \hat \eta^{(\balpha, l)}_{n+1} \partial_{\xi_i} \hat \chi^{(\balpha, l)}_{n+1} \mr {\bar \partial}_{j} \mr {\bar D}^{\balpha}\bar \rho_0,     \label{e:174:SolAspt} \\ 
      -L_{-1} \hat\rho_{n,0} =&\, - \frac{1}{\lambda} \sum_{[\balpha]=1, \balpha \in \mathcal{I}_*}^{2n-1} \sum_{l=1}^{\hat I_{n}^{(\balpha)}} \omega_{*,i} \hat h^{(\balpha, l)}_{n} \phi_* \hat \eta^{(\balpha, l)}_{n} \partial_{\xi_i} \hat \chi^{(\balpha, l)}_{n} \mr {\bar D}^{\balpha}\bar \rho_0.        \label{e:175:SolAspt}
\end{align}

For $a=0$ and fixed $\bbeta$, we need to consider all contributions from (\ref{e:136:SolAspt}-\ref{e:175:SolAspt}) and renumber them with index $1 \leq l \leq \hat I^{(\bbeta)}_{n+3}$. In the end, we also estimate the cardinality of such correctors, denoted by $\hat I^{(\bbeta)}_{n+3}$. These are presented in the following cases.

\begin{case}      \label{c:compHat:0}  
If $2 \leq [\bbeta] \leq 2n+4$ and $\bbeta = 1 \balpha$, the contribution from \eqref{e:138:SolAspt} concerning $L_{11} \hat \rho_{n+2}$ gives at most $\hat I^{(\balpha)}_{n+2}$ correctors with
\begin{align*}
      \hat h^{(\bbeta, \cdot)}_{n+3} =&\, - \frac{\delta^{\sfrac12}\mr{\bar \lambda}} {\kappa\lambda} \hat h^{(\balpha, l)}_{n+2},      \quad 
      \hat \eta^{(\bbeta, \cdot)}_{n+3} = \eta_2 \hat \eta^{(\balpha, l)}_{n+2}, \\
      \hat \chi^{(\bbeta, \cdot)}_{n+3} =&\, (-\Delta_\xi)^{-1} \big( \partial_{\xi_2} \Pi_2 \hat \chi^{(\balpha, l)}_{n+2} - \langle \partial_{\xi_2} \Pi_2 \hat \chi^{(\balpha, l)}_{n+2} \rangle_\xi \big). 
\end{align*}
\eqref{e:136:SolAspt} gives the same number of correctors involving $\eta_1$ for the differentiation index $\bbeta = 2 \balpha$. The estimate \eqref{e:EXP1:hh_Pest} with $k=n+2$ follows from \eqref{e:EXP1:hh_Pest} with $k=n+1$. The estimates \eqref{e:EXP1:hchi_est} and \eqref{e:EXP1:heta_est} for $k=n+2$ follow from (\ref{e:EXP1:hchi_est}-\ref{e:EXP1:heta_est}) for $k=n+1$.
\end{case}

\begin{case}      \label{c:compHat:1}  
If $1 \leq [\bbeta] \leq 2n+3$ and $\bbeta = \balpha$, the contribution from \eqref{e:139:SolAspt} concerning $L_{11} \hat \rho_{n+2}$ gives at most $\hat I^{(\balpha)}_{n+2}$ correctors with
\begin{align*}
      \hat h^{(\bbeta, \cdot)}_{n+3} =&\, - \frac{\delta^{\sfrac12}\mr{\bar \lambda}} {\kappa\lambda} \mr {\bar \partial}_{1} \hat h^{(\balpha, l)}_{n+2},      \quad 
      \hat \eta^{(\bbeta, \cdot)}_{n+3} = \eta_2 \hat \eta^{(\balpha, l)}_{n+2}, \\ 
      \hat \chi^{(\bbeta, \cdot)}_{n+3} =&\, (-\Delta_\xi)^{-1} \big( \partial_{\xi_2} \Pi_2 \hat \chi^{(\balpha, l)}_{n+2} - \langle \partial_{\xi_2} \Pi_2 \hat \chi^{(\balpha, l)}_{n+2} \rangle_\xi \big).      
\end{align*}
\eqref{e:137:SolAspt} gives the same number of correctors involving $\eta_1$. The estimate \eqref{e:EXP1:hh_Pest} with $k=n+2$ follows from \eqref{e:EXP1:hh_Pest} with $k=n+1$. The estimates \eqref{e:EXP1:hchi_est} and \eqref{e:EXP1:heta_est} are similar to Case \ref{c:compHat:0}.
\end{case}

\begin{case}      \label{c:compHat:2}  
If $ 2 \leq [\bbeta] \leq 2n+5$ and $\bbeta = 1 \balpha$, the contribution from \eqref{e:142:SolAspt} concerning $L_{11} \check \rho_{n+2}$ gives at most $\check I^{(\balpha)}_{n+2}$ correctors with
\begin{align*}
      \hat h^{(\bbeta, \cdot)}_{n+3} = - \frac{\delta^{\sfrac12}\mr{\bar \lambda}} {\kappa\lambda} \check h^{(\balpha, l)}_{n+2},       \quad 
      \hat \eta^{(\bbeta, \cdot)}_{n+3} = \eta_2 \check \eta^{(\balpha, l)}_{n+2},      \quad 
      \hat \chi^{(\bbeta, \cdot)}_{n+3} = (-\Delta_\xi)^{-1} \partial_{\xi_2} \Pi_2.
\end{align*}
\eqref{e:140:SolAspt} gives the same number of correctors involving $\eta_1$. The estimate \eqref{e:EXP1:hh_Pest} with $k=n+2$ follows from \eqref{d:lambda_gamma} and \eqref{e:EXP1:ch_Pest} with $k=n+1$. The estimate \eqref{e:EXP1:heta_est} is similar to Case \ref{c:compHat:0}. The estimate \eqref{e:EXP1:hchi_est} is obvious from the definition of $\Pi_2$.
\end{case}

\begin{case}      \label{c:compHat:3}  
If $1 \leq [\bbeta] \leq 2n+4$ and $\bbeta = \balpha$, the contribution from \eqref{e:143:SolAspt} concerning $L_{11} \check \rho_{n+2}$ gives at most $\check I^{(\balpha)}_{n+2}$ correctors with
\begin{align*}
      \hat h^{(\bbeta, \cdot)}_{n+3} = - \frac{\delta^{\sfrac12}\mr{\bar \lambda}} {\kappa\lambda} \mr {\bar \partial}_{1} \check h^{(\balpha, l)}_{n+2},    \quad 
      \hat \eta^{(\bbeta, \cdot)}_{n+3} = \eta_2 \check \eta^{(\balpha, l)}_{n+2},      \quad 
      \hat \chi^{(\bbeta, \cdot)}_{n+3} = (-\Delta_\xi)^{-1} \partial_{\xi_2} \Pi_2.
\end{align*}
\eqref{e:141:SolAspt} gives the same number of correctors involving $\eta_1$. The estimates \eqref{e:EXP1:hh_Pest}, \eqref{e:EXP1:hchi_est} and \eqref{e:EXP1:heta_est} with $k=n+2$ are similar to Case \ref{c:compHat:2}.
\end{case}

\begin{case}      \label{c:compHat:4}  
If $1 \leq [\bbeta] \leq 2n+3$ and $\bbeta = \balpha$, the contribution from \eqref{e:144:SolAspt} concerning $L_{12} \hat \rho_{n+2}$ gives at most $\hat I^{(\balpha)}_{n+2}$ correctors with
\begin{align*}
      \hat h^{(\bbeta, \cdot)}_{n+3} = - \frac{\mu}{\kappa \lambda} \hat h^{(\balpha, l)}_{n+2},     \quad 
      \hat \eta^{(\bbeta, \cdot)}_{n+3} = \partial_\tau \hat \eta^{(\balpha, l)}_{n+2},       \quad 
      \hat \chi^{(\bbeta, \cdot)}_{n+3} = (-\Delta_\xi)^{-1} \hat \chi^{(\balpha, l)}_{n+2}.
\end{align*}
The estimate \eqref{e:EXP1:hh_Pest} with $k=n+2$ follows from \eqref{d:lambda_gamma} and \eqref{e:EXP1:hh_Pest} with $k=n+1$. In particular, we need
\begin{align*}
      \frac{\mu}{\kappa \lambda} \leq \frac{\delta^{\sfrac12} \mr{\bar \lambda}}{\kappa\lambda} \cdot \lambda^{\gamma_S}.
\end{align*}
The estimates \eqref{e:EXP1:hchi_est} and \eqref{e:EXP1:heta_est} with $k=n+2$ are obvious.
\end{case}

\begin{case}      \label{c:compHat:5}  
If $1 \leq [\bbeta] \leq 2n+3$ and $\bbeta = \balpha$, the contribution from \eqref{e:145:SolAspt} concerning $L_{12} \hat \rho_{n+2}$ gives at most $N \hat I^{(\balpha)}_{n+2}$ correctors with
\begin{align*}
      \hat h^{(\bbeta, \cdot)}_{n+3} = -\varepsilon\lambda S_{m,11} \hat h^{(\balpha, l)}_{n+2},    \quad 
      \hat \eta^{(\bbeta, \cdot)}_{n+3} = \vartheta_m \hat \eta^{(\balpha, l)}_{n+2},   \quad
      \hat \chi^{(\bbeta, \cdot)}_{n+3} = (-\Delta_\xi)^{-1} \partial_{\xi_1\xi_1} \hat \chi^{(\balpha, l)}_{n+2}. 
\end{align*}
\eqref{e:146:SolAspt} gives the same number of correctors involving $\xi_2$. The estimate \eqref{e:EXP1:hh_Pest} with $k=n+2$ follows from \eqref{d:lambda_gamma} and \eqref{e:EXP1:hh_Pest} with $k=n+1$. The estimate \eqref{e:EXP1:heta_est} is similar to Case \eqref{c:compHat:0}. The estimate \eqref{e:EXP1:hchi_est} is obvious from the shear condition satisfied by $\hat \chi^{(\balpha, l)}_{n+2}$.
\end{case}

\begin{case}      \label{c:compHat:6}  
If $ 1 \leq [\bbeta] \leq 2n+1$ and $\bbeta = \balpha$, the contributions from \eqref{e:147:SolAspt} and \eqref{e:149:SolAspt} concerning $L_{01} \hat \rho_{n+1}$ give at most $\hat I^{(\balpha)}_{n+1}$ correctors with
\begin{align*}
      \hat h^{(\bbeta, \cdot)}_{n+3} = \Big( - \frac{\mr{\bar\mu}}{\kappa} \mr {\bar D}_t + \mr{\bar \lambda}^2 \mr {\bar \partial}_{ii} \Big) \hat h^{(\balpha, l)}_{n+1},       \quad  
      \hat \eta^{(\bbeta, \cdot)}_{n+3} = \hat \eta^{(\balpha, l)}_{n+1},   \quad 
      \hat \chi^{(\bbeta, \cdot)}_{n+3} = (-\Delta_\xi)^{-1} \hat \chi^{(\balpha, l)}_{n+1}.  
\end{align*}
The estimate \eqref{e:EXP1:hh_Pest} with $k=n+2$ follows from \eqref{d:lambda_gamma} and \eqref{e:EXP1:hh_Pest} with $k=n$. In particular, we need
\begin{align}
      \frac{\mr{\bar\mu}}{\kappa} \leq \bigg( \frac{\delta^{\sfrac12} \mr{\bar \lambda}}{\kappa\lambda} \bigg)^2 \cdot \lambda^{2\gamma_S}.
\end{align}
The estimates \eqref{e:EXP1:hchi_est} and \eqref{e:EXP1:heta_est} with $k=n+2$ are obvious.
\end{case}

\begin{case}      \label{c:compHat:7}  
If $2 \leq [\bbeta] \leq 2n+2$ and $\bbeta = t\balpha$, the contribution from \eqref{e:148:SolAspt} concerning $L_{01} \hat \rho_{n+1}$ gives at most $\hat I^{(\balpha)}_{n+1}$ correctors with
\begin{align*}
      \hat h^{(\bbeta, \cdot)}_{n+3} = - \frac{\mr{\bar\mu}}{\kappa} \hat h^{(\balpha, l)}_{n+1},       \quad 
      \hat \eta^{(\bbeta, \cdot)}_{n+3} = \hat \eta^{(\balpha, l)}_{n+1} ,        \quad
      \hat \chi^{(\bbeta, \cdot)}_{n+3} = (-\Delta_\xi)^{-1} \hat \chi^{(\balpha, l)}_{n+1}.
\end{align*}
The estimate \eqref{e:EXP1:hh_Pest} with $k=n+2$ follows from \eqref{d:lambda_gamma} and \eqref{e:EXP1:hh_Pest} with $k=n$. The estimates \eqref{e:EXP1:hchi_est} and \eqref{e:EXP1:heta_est} with $k=n+2$ are obvious.
\end{case}

\begin{case}      \label{c:compHat:8}  
If $3 \leq [\bbeta] \leq 2n+3$ and $\bbeta = 11 \balpha$, the contribution from \eqref{e:149:SolAspt} concerning $L_{01} \hat \rho_{n+1}$ give at most $\hat I^{(\balpha)}_{n+1}$ correctors with
\begin{align*}
      \hat h^{(\bbeta, \cdot)}_{n+3} = \mr{\bar \lambda}^2 \hat h^{(\balpha, l)}_{n+1},         \quad 
      \hat \eta^{(\bbeta, \cdot)}_{n+3} = \hat \eta^{(\balpha, l)}_{n+1},   \quad 
      \hat \chi^{(\bbeta, \cdot)}_{n+3} = (-\Delta_\xi)^{-1} \hat \chi^{(\balpha, l)}_{n+1},
\end{align*}
and the same number of correctors for the differentiation index $\bbeta = 22 \balpha$. The estimate \eqref{e:EXP1:hh_Pest} with $k=n+2$ follows from \eqref{d:lambda_gamma} and \eqref{e:EXP1:hh_Pest} with $k=n$. (\ref{e:EXP1:hchi_est}-\ref{e:EXP1:heta_est}) with $k=n+2$ are obvious.
\end{case}

\begin{case}      \label{c:compHat:9}  
If $2 \leq [\bbeta] \leq 2n+2$ and $\bbeta = 1 \balpha$, the contribution from \eqref{e:150:SolAspt} concerning $L_{01} \hat \rho_{n+1}$ gives at most $\hat I^{(\balpha)}_{n+1}$ correctors with
\begin{align*}
      \hat h^{(\bbeta, \cdot)}_{n+3} = \mr{\bar \lambda}^2 \mr {\bar \partial}_{1} \hat h^{(\balpha, l)}_{n+1},       \quad 
      \hat \eta^{(\bbeta, \cdot)}_{n+3} = \hat \eta^{(\balpha, l)}_{n+1},         \quad 
      \hat \chi^{(\bbeta, \cdot)}_{n+3} = (-\Delta_\xi)^{-1} \hat \chi^{(\balpha, l)}_{n+1},
\end{align*}
and the same number of correctors for the differentiation index $\bbeta = 2 \balpha$. The estimate \eqref{e:EXP1:hh_Pest} with $k=n+2$ follows from \eqref{d:lambda_gamma} and \eqref{e:EXP1:hh_Pest} with $k=n$. (\ref{e:EXP1:hchi_est}-\ref{e:EXP1:heta_est}) with $k=n+2$ are obvious.
\end{case}

\begin{case}      \label{c:compHat:10}  
If $1 \leq [\bbeta] \leq 2n+1$ and $\bbeta = \balpha$, the contribution from \eqref{e:151:SolAspt} concerning $L_{02} \hat \rho_{n+1}$ gives at most $N\hat I^{(\balpha)}_{n+1}$ correctors with
\begin{align*}
      \hat h^{(\bbeta, \cdot)}_{n+3} =&\, \varepsilon \frac{\delta^{\sfrac12} \mr{\bar \lambda}^2}{\kappa \lambda} \mr {\bar \partial}_{i} \Omega_{m,ij} \mr {\bar \partial}_{j} \hat h^{(\balpha, l)}_{n+1},      \quad 
      \hat \eta^{(\bbeta, \cdot)}_{n+3} = \sigma_m \eta_1 \hat \eta^{(\balpha, l)}_{n+1},     \\
      \hat \chi^{(\bbeta, \cdot)}_{n+3} =&\, (-\Delta_\xi)^{-1} \big( \partial_{\xi_1} \Pi_1 \hat \chi^{(\balpha, l)}_{n+1} - \langle \cdot \rangle_{\xi} \big).      
\end{align*}
\eqref{e:153:SolAspt} gives the same number of correctors involving $\eta_2$. The estimate \eqref{e:EXP1:hh_Pest} with $k=n+2$ follows from \eqref{d:lambda_gamma}, \eqref{e:mtrEstHom} and \eqref{e:EXP1:hh_Pest} with $k=n$. The estimate \eqref{e:EXP1:heta_est} follows from Lemma \ref{l:cutoff}, \eqref{e:prdEstHom} and \eqref{e:EXP1:heta_est} with $k=n$. The proof of \eqref{e:EXP1:hchi_est} is similar to Case \ref{c:compHat:0}.
\end{case}

\begin{case}      \label{c:compHat:11}  
If $2 \leq [\bbeta] \leq 2n+2$ and $\bbeta = 1 \balpha$, the contribution from \eqref{e:152:SolAspt} concerning $L_{02} \hat \rho_{n+1}$ give at most $N\hat I^{(\balpha, l)}_{n+1}$ correctors with
\begin{align*}
      \hat h^{(\bbeta, \cdot)}_{n+3} =&\, \varepsilon \frac{\delta^{\sfrac12}\mr{\bar \lambda}^2}{\kappa \lambda} \mr {\bar \partial}_{i} \Omega_{m,i1} \hat h^{(\balpha, l)}_{n+1},       \quad 
      \hat \eta^{(\bbeta, \cdot)}_{n+3} = \sigma_m \eta_1 \hat \eta^{(\balpha, l)}_{n+1},     \\ 
      \hat \chi^{(\bbeta, \cdot)}_{n+3} =&\, (-\Delta_\xi)^{-1} \big( \partial_{\xi_1} \Pi_1 \hat \chi^{(\balpha, l)}_{n+1} - \langle \cdot \rangle_{\xi} \big).  
\end{align*}
\eqref{e:154:SolAspt} gives the same number of correctors involving $\eta_2$. The estimates (\ref{e:EXP1:hh_Pest},\ref{e:EXP1:hchi_est},\ref{e:EXP1:heta_est}) with $k=n+2$ are similar to Case \ref{c:compHat:10}.
\end{case}

\begin{case}      \label{c:compHat:12}  
If $1 \leq [\bbeta] \leq 2n+1$ and $\bbeta = \balpha$, the contribution from \eqref{e:155:SolAspt} concerning $L_{03} \hat \rho_{n+1}$ gives at most $N\hat I^{(\balpha)}_{n+1}$ correctors with
\begin{align*}
      \hat h^{(\bbeta, \cdot)}_{n+3} =&\, \varepsilon \frac{\delta^{\sfrac12}\mr{\bar \lambda}}{\kappa} E_{m,j2} \mr {\bar \partial}_{j} \hat h^{(\balpha, l)}_{n+1},       \quad 
      \hat \eta^{(\bbeta, \cdot)}_{n+3} = \varphi_m \eta_1 \hat \eta^{(\balpha, l)}_{n+1}     \\
      \hat \chi^{(\bbeta, \cdot)}_{n+3} =&\, (-\Delta_\xi)^{-1} \big( \partial_{\xi_1} \Pi_1 \hat \chi^{(\balpha, l)}_{n+1} - \langle \cdot \rangle_{\xi} \big).
\end{align*}
\eqref{e:157:SolAspt} gives the same number of correctors involving $\eta_2$. The estimate \eqref{e:EXP1:hh_Pest} with $k=n+2$ follows from \eqref{d:lambda_gamma}, \eqref{e:mtrEstHom} and \eqref{e:EXP1:hh_Pest} with $k=n$. In particular, we need
\begin{align*}
      \varepsilon \frac{\delta^{\sfrac12}\mr{\bar \lambda}}{\kappa} 
      \leq \bigg( \frac{\delta^{\sfrac12} \mr{\bar \lambda}}{\kappa\lambda} \bigg)^2 \cdot \lambda^{2\gamma_S}.
\end{align*}
The estimates (\ref{e:EXP1:hchi_est},\ref{e:EXP1:heta_est}) with $k=n+2$ are similar to Case \ref{c:compHat:10}.
\end{case}

\begin{case}      \label{c:compHat:13}  
If $2 \leq [\bbeta] \leq 2n+2$ and $\bbeta = 1 \balpha$, the contribution from \eqref{e:156:SolAspt} concerning $L_{03} \hat \rho_{n+1}$ gives at most $N\hat I^{(\balpha)}_{n+1}$ correctors with
\begin{align*}
      \hat h^{(\bbeta, \cdot)}_{n+3} =&\, \varepsilon \frac{\delta^{\sfrac12}\mr{\bar \lambda}} {\kappa} E_{m,12} \hat h^{(\balpha, l)}_{n+1},         \quad 
      \hat \eta^{(\bbeta, \cdot)}_{n+3} = \varphi_m \eta_1 \hat \eta^{(\balpha, l)}_{n+1},    \\ 
      \hat \chi^{(\bbeta, \cdot)}_{n+3} =&\, (-\Delta_\xi)^{-1} \big( \partial_{\xi_1} \Pi_1 \hat \chi^{(\balpha, l)}_{n+1} - \langle \cdot \rangle_{\xi} \big). 
\end{align*}
\eqref{e:158:SolAspt} gives the same number of correctors involving $\eta_2$. The estimate \eqref{e:EXP1:hh_Pest} with $k=n+2$ follows from \eqref{d:lambda_gamma}, \eqref{e:mtrEstHom} and \eqref{e:EXP1:hh_Pest} with $k=n$. The estimates (\ref{e:EXP1:hchi_est},\ref{e:EXP1:heta_est}) with $k=n+2$ are similar to Case \ref{c:compHat:10}.
\end{case}

\begin{case}      \label{c:compHat:14}  
If $1 \leq [\bbeta] \leq 2n+1$ and $\bbeta = \balpha$, the contribution from \eqref{e:159:SolAspt} concerning $L_{04} \hat \rho_{n+1}$ gives at most $N\hat I^{(\balpha)}_{n+1}$ correctors with
\begin{align*}
      \hat h^{(\bbeta, \cdot)}_{n+3} =&\, \varepsilon \frac{\delta^{\sfrac12}\mr{\bar \lambda}} {\kappa} \mr {\bar \partial}_{j} E_{m,j2} \hat h^{(\balpha, l)}_{n+1},         \quad 
      \hat \eta^{(\bbeta, \cdot)}_{n+3} = \varphi_m \eta_1 \hat \eta^{(\balpha, l)}_{n+1},    \\ 
      \hat \chi^{(\bbeta, \cdot)}_{n+3} =&\, (-\Delta_\xi)^{-1} \big( \Pi_1 \partial_{\xi_1} \hat \chi^{(\balpha, l)}_{n+1} - \langle \cdot \rangle_{\xi} \big).
\end{align*}
\eqref{e:160:SolAspt} gives the same number of correctors involving $\eta_2$. The estimate \eqref{e:EXP1:hh_Pest} with $k=n+2$ follows from \eqref{d:lambda_gamma}, \eqref{e:mtrEstHom} and \eqref{e:EXP1:hh_Pest} with $k=n$. The estimates (\ref{e:EXP1:hchi_est},\ref{e:EXP1:heta_est}) with $k=n+2$ are similar to Case \ref{c:compHat:10}.
\end{case}

\begin{case}      \label{c:compHat:15}  
If $1 \leq [\bbeta] \leq 2n+2$ and $\bbeta = \balpha$, the contribution from \eqref{e:161:SolAspt} concerning $L_{02} \check \rho_{n+1}$ give at most $N\check I^{(\balpha)}_{n+1}$ correctors with
\begin{align*}
      \hat h^{(\bbeta, \cdot)}_{n+3} = \varepsilon \frac{\delta^{\sfrac12}\mr{\bar \lambda}^2} {\kappa\lambda} \mr {\bar \partial}_{i} \Omega_{m,ij} \mr {\bar \partial}_{j} \check h^{(\balpha, l)}_{n+1},     \quad 
      \hat \eta^{(\bbeta, \cdot)}_{n+3} = \sigma_m \eta_1 \check \eta^{(\balpha, l)}_{n+1},   \quad
      \hat \chi^{(\bbeta, \cdot)}_{n+3} = (-\Delta_\xi)^{-1} \Pi_1, 
\end{align*}
and the same number of correctors involving $\eta_2$. The estimate \eqref{e:EXP1:hh_Pest} with $k=n+2$ follows from \eqref{d:lambda_gamma}, \eqref{e:mtrEstHom} and \eqref{e:EXP1:ch_Pest} with $k=n$. The estimate \eqref{e:EXP1:hchi_est} with $k=n+2$ is similar to Case \ref{c:compHat:10}. The estimate \eqref{e:EXP1:heta_est} is obvious.
\end{case}

\begin{case}      \label{c:compHat:16}  
If $2 \leq [\bbeta] \leq 2n+3$ and $\bbeta = 1 \balpha$, the contribution from \eqref{e:162:SolAspt} concerning $L_{02} \check \rho_{n+1}$ gives at most $N\check I^{(\balpha)}_{n+1}$ correctors with
\begin{align*}
      \hat h^{(\bbeta, \cdot)}_{n+3} = \varepsilon \frac{\delta^{\sfrac12}\mr{\bar \lambda}^2} {\kappa\lambda} \mr {\bar \partial}_{i} \Omega_{m,i1} \check h^{(\balpha, l)}_{n+1},     \quad 
      \hat \eta^{(\bbeta, \cdot)}_{n+3} = \sigma_m \eta_1 \check \eta^{(\balpha, l)}_{n+1},   \quad 
      \hat \chi^{(\bbeta, \cdot)}_{n+3} = (-\Delta_\xi)^{-1} \Pi_1,  
\end{align*}
and the same number of correctors involving $\eta_2$. The estimate \eqref{e:EXP1:hh_Pest} with $k=n+2$ follows from \eqref{d:lambda_gamma}, \eqref{e:mtrEstHom} and \eqref{e:EXP1:ch_Pest} with $k=n$. The estimates (\ref{e:EXP1:hchi_est},\ref{e:EXP1:heta_est}) with $k=n+2$ are similar to Case \ref{c:compHat:15}.
\end{case}

\begin{case}      \label{c:compHat:17}  
If $1 \leq [\bbeta] \leq 2n+2$ and $\bbeta = \balpha$, the contribution from line \eqref{e:163:SolAspt} $L_{03} \check \rho_{n+1}$ give at most $N\check I^{(\balpha)}_{n+1}$ correctors with
\begin{align*}
      \hat h^{(\bbeta, \cdot)}_{n+3} = \varepsilon \frac{\delta^{\sfrac12}\mr{\bar \lambda}} {\kappa} E_{m,j2} \mr {\bar \partial}_{j} \check h^{(\balpha, l)}_{n+1},       \quad 
      \hat \eta^{(\bbeta, \cdot)}_{n+3} = \varphi_m \eta_1 \check \eta^{(\balpha, l)}_{n+1},  \quad 
      \hat \chi^{(\bbeta, \cdot)}_{n+3} = (-\Delta_\xi)^{-1} \partial_{\xi_1} \Pi_1.
\end{align*}
\eqref{e:165:SolAspt} gives the same number of correctors involving $\eta_2$. The estimate \eqref{e:EXP1:hh_Pest} with $k=n+2$ follows from \eqref{d:lambda_gamma}, \eqref{e:mtrEstHom} and \eqref{e:EXP1:ch_Pest} with $k=n$. The estimates (\ref{e:EXP1:hchi_est},\ref{e:EXP1:heta_est}) with $k=n+2$ are similar to Case \ref{c:compHat:15}.
\end{case}

\begin{case}      \label{c:compHat:18}  
If $2 \leq [\bbeta] \leq 2n+3$ and $\bbeta = 1 \balpha$, the contribution from line \eqref{e:164:SolAspt} concerning $L_{03} \check \rho_{n+1}$ give at most $N \check I^{(\balpha)}_{n+1}$ correctors with
\begin{align*}
      \hat h^{(\bbeta, \cdot)}_{n+3} = \varepsilon \frac{\delta^{\sfrac12}\mr{\bar \lambda}} {\kappa} E_{m,12} \check h^{(\balpha, l)}_{n+1},       \quad 
      \hat \eta^{(\bbeta, \cdot)}_{n+3} = \varphi_m \eta_1 \check \eta^{(\balpha, l)}_{n+1},        \quad 
      \hat \chi^{(\bbeta, \cdot)}_{n+3} = (-\Delta_\xi)^{-1} \partial_{\xi_1} \Pi_1.
\end{align*}
\eqref{e:166:SolAspt} gives the same number of correctors involving $\eta_2$. The estimates (\ref{e:EXP1:hchi_est},\ref{e:EXP1:heta_est}) with $k=n+2$ are similar to Case \ref{c:compHat:17}.
\end{case}

\begin{case}      \label{c:compHat:19}  
If $1 \leq [\bbeta] \leq 2n+1$ and $\bbeta = \balpha$, the contribution from \eqref{e:170:SolAspt} and \eqref{e:173:SolAspt} concerning $L_{05} \hat \rho_{n+1}$ and $L_{07} \hat \rho_{n+1}$ gives at most $N\hat I^{(\balpha)}_{n+1}$ correctors with
\begin{align*}
      \hat h^{(\bbeta, \cdot)}_{n+3} =&\, \varepsilon \lambda \mr{\bar \lambda} 
            \big( \mr {\bar \partial}_{j} B_{m,1j} \hat h^{(\balpha, l)}_{n+1} 
            + 2B_{m,1j} \mr {\bar \partial}_{j} \hat h^{(\balpha, l)}_{n+1} \big),       \\
      \hat \eta^{(\bbeta, \cdot)}_{n+3} =&\, \phi_m \hat \eta^{(\balpha, l)}_{n+1},        \quad 
      \hat \chi^{(\bbeta, \cdot)}_{n+3} = (-\Delta_\xi)^{-1} \partial_{\xi_1} \hat \chi^{(\balpha, l)}_{n+1}.
\end{align*}
\eqref{e:170:SolAspt} and \eqref{e:173:SolAspt} also give the same number of correctors involving $\eta_2$.
The estimate \eqref{e:EXP1:hh_Pest} with $k=n+2$ follows from \eqref{d:lambda_gamma}, \eqref{e:mtrEstHom} and \eqref{e:EXP1:hh_Pest} with $k=n$. The estimates (\ref{e:EXP1:hchi_est},\ref{e:EXP1:heta_est}) with $k=n+2$ are similar to Case \ref{c:compHat:17}.
\end{case}

\begin{case}      \label{c:compHat:21}
If $2 \leq [\bbeta] \leq 2n+2$ and $\bbeta = 1\balpha$, the contribution from \eqref{e:174:SolAspt} concerning $L_{07} \hat \rho_{n+1}$ gives at most $N\hat I^{(\balpha)}_{n+1}$ correctors with
\begin{align*}
      \hat h^{(\bbeta, \cdot)}_{n+3} = 2\varepsilon \lambda \mr{\bar \lambda} 
            B_{m,1j} \hat h^{(\balpha, l)}_{n+1},       \quad 
      \hat \eta^{(\bbeta, \cdot)}_{n+3} = \phi_m \hat \eta^{(\balpha, l)}_{n+1},        \quad 
      \hat \chi^{(\bbeta, \cdot)}_{n+3} = (-\Delta_\xi)^{-1} \partial_{\xi_1} \hat \chi^{(\balpha, l)}_{n+1}.
\end{align*}
\eqref{e:174:SolAspt} also gives the same number correctors involving $\xi_2$. The estimates (\ref{e:EXP1:hh_Pest},\ref{e:EXP1:hchi_est},\ref{e:EXP1:heta_est}) with $k=n+2$ are similar to Case \ref{c:compHat:19}.
\end{case}

\begin{case}      \label{c:compHat:22}
If $1 \leq [\bbeta] \leq 2n+1$ and $\bbeta = \balpha$, the contribution from \eqref{e:171:SolAspt} concerning $L_{06} \hat \rho_{n+1}$ gives at most $\hat I^{(\balpha)}_{n+1}$ correctors with
\begin{align*}
      \hat h^{(\bbeta, \cdot)}_{n+3} = 2\lambda \mr{\bar \lambda} 
             \mr {\bar \partial}_{1} \hat h^{(\balpha, l)}_{n+1},       \quad 
      \hat \eta^{(\bbeta, \cdot)}_{n+3} = \hat \eta^{(\balpha, l)}_{n+1},        \quad 
      \hat \chi^{(\bbeta, \cdot)}_{n+3} = (-\Delta_\xi)^{-1} \partial_{\xi_1} \hat \chi^{(\balpha, l)}_{n+1}.
\end{align*}
\eqref{e:171:SolAspt} also gives the same number correctors involving $\xi_2$. The estimate \eqref{e:EXP1:hh_Pest} with $k=n+2$ follows from \eqref{d:lambda_gamma} and \eqref{e:EXP1:hh_Pest} with $k=n$. In particular, we need 
\begin{align}
      2\lambda \mr{\bar \lambda} \leq \bigg( \frac{\delta^{\sfrac12} \mr{\bar \lambda}}{\kappa\lambda} \bigg)^2 \cdot \lambda^{2\gamma_S}.
\end{align}
The estimates (\ref{e:EXP1:hchi_est}-\ref{e:EXP1:heta_est}) with $k=n+2$ are obvious.
\end{case}

\begin{case}      \label{c:compHat:23}
If $2 \leq [\bbeta] \leq 2n+2$ and $\bbeta = 1\balpha$, the contribution from \eqref{e:172:SolAspt} concerning $L_{06} \hat \rho_{n+1}$ gives at most $\hat I^{(\balpha)}_{n+1}$ correctors with
\begin{align*}
      \hat h^{(\bbeta, \cdot)}_{n+3} = 2\lambda \mr{\bar \lambda} 
            \hat h^{(\balpha, l)}_{n+1},       \quad 
      \hat \eta^{(\bbeta, \cdot)}_{n+3} = \hat \eta^{(\balpha, l)}_{n+1},        \quad 
      \hat \chi^{(\bbeta, \cdot)}_{n+3} = (-\Delta_\xi)^{-1} \partial_{\xi_1} \hat \chi^{(\balpha, l)}_{n+1}.
\end{align*}
\eqref{e:172:SolAspt} also gives the same number correctors involving $\xi_2$. The estimates (\ref{e:EXP1:hh_Pest},\ref{e:EXP1:hchi_est},\ref{e:EXP1:heta_est}) with $k=n+2$ are similar to Case \ref{c:compHat:22}.
\end{case}

\begin{case}      \label{c:compHat:20}  
If $1 \leq [\bbeta] \leq 2n-1$ and $\bbeta = \balpha$, the contribution from line \eqref{e:175:SolAspt} concerning $L_{-1} \hat \rho_{n}$ gives at most $\hat I^{(\balpha)}_{n}$ correctors with
\begin{align*}
      \hat h^{(\bbeta, \cdot)}_{n+3} = -\frac{1}{\kappa\lambda} \omega_{*,1} \hat h^{(\balpha, l)}_{n},         \quad 
      \hat \eta^{(\bbeta, \cdot)}_{n+3} = \phi_* \hat \eta^{(\balpha, l)}_{n},          \quad 
      \hat \chi^{(\bbeta, \cdot)}_{n+3} = (-\Delta_\xi)^{-1} \partial_{\xi_1} \hat \chi^{(\balpha, l)}_{n}, 
\end{align*}
and the same number of correctors involving $\xi_2$. The estimate \eqref{e:EXP1:hh_Pest} with $k=n+2$ follows from  \eqref{d:lambda_gamma}, \eqref{e:mtrEstHom} and \eqref{e:EXP1:hh_Pest} with $k=n-1$. The estimates (\ref{e:EXP1:hchi_est},\ref{e:EXP1:heta_est}) with $k=n+2$ are similar to Case \ref{c:compHat:17}.
\end{case}

In all cases above, the information \eqref{e:EXP1:hh} for $n+3$ follows from \eqref{e:mtrStrHom}, Remark \ref{r:opPolynomial}, \eqref{e:EXP1:hh} for indices smaller or equal to $n+2$ and \eqref{e:EXP1:ch} for indices smaller or equal to $n+2$. \eqref{e:EXP1:hhS_Pest} follows from \eqref{e:EXP1:hh_Pest} for all admissible $\balpha$ and $l$.

Now we estimate the number of correctors $\hat I^{(\bbeta)}_{n+3}$. Similar to $\check \rho_{n+2}$, we consider two cases.
\begin{case}      \label{c:numHat:0}
For index with form $\bbeta = 1 \balpha = 11\bgamma$, we have
\begin{align}
      \hat I^{(\bbeta)}_{n+3} \leq&\, \hat I^{(\balpha)}_{n+2} + 2\hat I^{(\bbeta)}_{n+2} + \check I^{(\balpha)}_{n+2} + 2 \check I^{(\bbeta)}_{n+2} + \hat I^{(\bbeta)}_{n+2} + 2 N \hat I^{(\bbeta)}_{n+2} + \hat I^{(\bbeta)}_{n+1}      \label{e:2:numHat} \\ 
      +&\, \hat I^{(\bgamma)}_{n+1} + \hat I^{(\balpha)}_{n+1} + 2N \hat I^{(\bbeta)}_{n+1} + N \hat I^{(\balpha)}_{n+1} + 2N \hat I^{(\bbeta)}_{n+1} + N \hat I^{(\balpha)}_{n+1}        \label{e:4:numHat} \\ 
      +&\, 2N \hat I^{(\bbeta)}_{n+1} + 2N \check I^{(\bbeta)}_{n+1} + N \check I^{(\balpha)}_{n+1} + 2N \check I^{(\bbeta)}_{n+1} + N \check I^{(\balpha)}_{n+1} + 2 N \hat I^{(\bbeta)}_{n+1}         \label{e:6:numHat} \\ 
      +&\, N \hat I^{(\balpha)}_{n+1} + 2 \hat I^{(\bbeta)}_{n+1} + \hat I^{(\balpha)}_{n+1} + 2 \hat I^{(\bbeta)}_{n}.     \label{e:8:numHat}
\end{align}
The terms in \eqref{e:2:numHat} come from respectively Case \ref{c:compHat:0}, Case \ref{c:compHat:1}, Case \ref{c:compHat:2}, Case \ref{c:compHat:3}, Case \ref{c:compHat:4}, Case \ref{c:compHat:5}, Case \ref{c:compHat:6}. The terms in \eqref{e:4:numHat} come from Case \ref{c:compHat:8}, Case \ref{c:compHat:9}, Case \ref{c:compHat:10}, Case \ref{c:compHat:11}, Case \ref{c:compHat:12}, Case \ref{c:compHat:13}. The terms in \eqref{e:6:numHat} come from Case \ref{c:compHat:14}, Case \ref{c:compHat:15}, Case \ref{c:compHat:16}, Case \ref{c:compHat:17}, Case \ref{c:compHat:18}, Case \ref{c:compHat:19}. The terms in \eqref{e:8:numHat} come from Case \ref{c:compHat:21}, Case \ref{c:compHat:22}, Case \ref{c:compHat:23}, Case \ref{c:compHat:20}.
\end{case}

\begin{case}      \label{c:numHat:1}
For index with form $\bbeta = t \balpha$, we have
\begin{align}
      \hat I^{(\bbeta)}_{n+3} \leq&\, 2\hat I^{(\bbeta)}_{n+2} + 2 \check I^{(\bbeta)}_{n+2} + \hat I^{(\bbeta)}_{n+2} + 2 N \hat I^{(\bbeta)}_{n+2} + \hat I^{(\bbeta)}_{n+1} + \hat I^{(\balpha)}_{n+1}         \label{e:10:numHat} \\ 
      +&\, 2N \hat I^{(\bbeta)}_{n+1} + 2N \hat I^{(\bbeta)}_{n+1} + 2N \hat I^{(\bbeta)}_{n+1} + 2N \check I^{(\bbeta)}_{n+1} + 2N \check I^{(\bbeta)}_{n+1}     \label{e:12:numHat} \\ 
      +&\, 2N \hat I^{(\bbeta)}_{n+1} + 2I^{(\bbeta)}_{n+1} + \hat I^{(\bbeta)}_{n}.     \label{e:14:numHat}
\end{align}
The terms in \eqref{e:10:numHat} come from Case \ref{c:compHat:1}, Case \ref{c:compHat:3}, Case \ref{c:compHat:4}, Case \ref{c:compHat:5}, Case \ref{c:compHat:6} and Case \ref{c:compHat:7}. The terms in \eqref{e:12:numHat} come from Case \ref{c:compHat:10}, Case \ref{c:compHat:12}, Case \ref{c:compHat:14}, Case \ref{c:compHat:15} and Case \ref{c:compHat:17}. The terms in \eqref{e:14:numHat} come from Case \ref{c:compHat:19}, Case \ref{c:compHat:22} and Case \ref{c:compHat:20}.
\end{case}

Now it is not difficult to see Case \ref{c:numHat:0}, Case \ref{c:numHat:1} and all other cases not listed here, obey the same estimate
\begin{equation}  \label{e:180:SolAspt}
\begin{split}
      \hat I^{(\bbeta)}_{n+3} \leq&\, \hat I^{(\balpha)}_{n+2} + (2N+3) \hat I^{(\bbeta)}_{n+2} + \check I^{(\balpha)}_{n+2} + 2 \check I^{(\bbeta)}_{n+2} \\ 
            +&\, (3N+2) \hat I^{(\balpha)}_{n+1} + (8N+3) \hat I^{(\bbeta)}_{n+1} + \hat I^{(\bgamma)}_{n+1} + 2N \check I^{(\balpha)}_{n+1} + 4N \check I^{(\bbeta)}_{n+1}
            + 2\hat I^{(\bbeta)}_{n},
\end{split}
\end{equation}
from which we can easily verify \eqref{e:EXP1:hI_est} for $k=n+2$, using \eqref{e:EXP1:hI_est} and \eqref{e:EXP1:cI_est} with $k \leq n+1$.

As for the temporal correctors, the $F$-related terms are analogous. We omit the details here.

\subsection{Technical preparation for remainders}    \label{ss:tech_inerH}

One task we shall do later is to estimate the remainders in expansion lemmas in Section \ref{ss:stateExpInertial}. Specifically, we need to estimate $\tilde \rho$ in Lemma \ref{l:EXP1} and Lemma \ref{l:EXP2}. In this section, we prove some technical lemmas to prepare for this task.

\begin{lemma}     \label{l:compDecompedF}
In the notions of Definition \ref{d:homFlowmap} and Definition \ref{d:homParaI}, consider $\rho: \T^2 \times [0,1] \rightarrow \R$, $\chi:\T^2 \rightarrow \R$, $\eta: \T \rightarrow \R $ and $h \in \mathcal{\bar{P}}( 10N^2Q )$. Suppose that, for $\bbeta \in \mathcal{I}_*$ and $p \in \N$ with $[\bbeta], p \leq 11Q^3$, we have the following estimates
\begin{align*}
      \| \mr {\bar D}^{\bbeta} \rho \|_\infty \lesssim &\, \bar \lambda^{-[\bbeta] b \gamma},  \\ 
      \vertiii{ h } + \| \nabla_\xi^p \chi \|_{\infty} + &\, \| \partial_\tau^p \eta \|_{\infty} 
            \lesssim 1.
\end{align*}
Suppose $(\chi, \eta)$ satisfies the generalized shear condition. For some $\balpha \in \mathcal{I}_*$, define $\tilde F: \T^2 \times [0,1] \rightarrow \R$ via
\begin{align*}
      \tilde F(x,t) = h(x,t) \eta(\mu t) \chi(\lambda \Phi(x,t)) \mr {\bar D}^{\balpha} \rho(x,t),
\end{align*}
then we have the following estimates for any $ \bbeta \in \mathcal{I}_* $ with $[\bbeta] + [\balpha] \leq 11Q^3$,
\begin{align}
      \| \mr O^{\bbeta} \tilde F \|
      \lesssim&\, \bar \lambda^{- [\balpha] b\gamma} \bigg( \frac{\lambda}{\mr\lambda} \bigg)^{|\bbeta|_x} \bigg( \frac{ \delta^{\sfrac12}\lambda } {\mr\mu} \bigg)^{|\bbeta|_t}.
            \label{e:10:compDecompedF} 
\end{align}
\end{lemma}

\begin{proof}[Proof of Lemma \ref{l:compDecompedF}]
We deduce from \eqref{e:12:homFlowmap},
\begin{align}
      \bar D_t \big( \eta (\mu t) \chi (\lambda \Phi) \big)
            =&\, \mu \partial_\tau \eta (\mu t) \chi (\lambda \Phi) + \lambda \eta (\mu t) \bar D_t \Phi \cdot \nabla_\xi \chi (\lambda \Phi)      \nonumber \\ 
            =&\, \mu \partial_\tau \eta (\mu t) \chi (\lambda \Phi) + \lambda^{-2} \eta (\mu t) \Omega(x,t,\mu t) \cdot \nabla_\xi \chi (\lambda \Phi)
                  \label{e:12:compDecompedF} \\       
                  +&\, \lambda \eta (\mu t) Z(x,t,\mu t) \cdot \nabla_\xi \chi (\lambda \Phi).
                  \label{e:14:compDecompedF}
\end{align}
Since $(\chi, \eta)$ satisfies the generalized shear condition, we can deduce from \eqref{e:20:homFlowmap} that \eqref{e:14:compDecompedF} is zero. Using \eqref{e:6:homFlowmap}, \eqref{e:8:homFlowmap}, \eqref{e:16:homFlowmap}, \eqref{e:mtrEstHom} and \eqref{e:prdEstHom}, we have
\begin{align}
      \big| \bar D_t \big( \eta (\mu t) \chi (\lambda \Phi) \big) \big| 
            \lesssim&\, \mu,      \label{e:18:compDecompedF} \\ 
      \big| \nabla \big( \eta (\mu t) \chi (\lambda \Phi) \big) \big| 
            \lesssim&\, \lambda.      \label{e:20:compDecompedF}
\end{align}
Given the knowledge in \eqref{e:18:compDecompedF} and \eqref{e:20:compDecompedF}, the following estimates follow from a direct induction argument similar to the proof of Lemma \ref{l:streamEst},
\begin{align*}
      \big| \bar D^{\bbeta} \big( \eta (\mu t) \chi (\lambda \Phi) \big) \big| 
      \lesssim \lambda^{|\bbeta|_x} \mu^{|\bbeta|_t},
            \quad [\bbeta] \leq 11Q^3.
\end{align*}
Note that we have $O_t = \bar D_t + \nabla^\perp \psi_{q+1} \cdot \nabla$, then we have
\begin{align}      \label{e:22:compDecompedF}
      \big| O^{\bbeta} \big( \eta (\mu t) \chi (\lambda \Phi) \big) \big| 
      \lesssim \lambda^{|\bbeta|_x} (\delta^{\sfrac12} \lambda)^{|\bbeta|_t},
            \quad [\bbeta] \leq 11Q^3,
\end{align}
and subsequently 
\begin{align}     \label{e:24:compDecompedF}
      \big\| \mr O^{\bbeta} \big( \eta (\mu t) \chi (\lambda \Phi) \big) \big\|_\infty 
      \lesssim \bigg( \frac{\lambda}{\mr\lambda} \bigg)^{|\bbeta|_x} \bigg( \frac{ \delta^{\sfrac12}\lambda } {\mr\mu} \bigg)^{|\bbeta|_t},
            \quad [\bbeta] \leq 11Q^3.
\end{align}
Indeed, to prove \eqref{e:22:compDecompedF}-\eqref{e:24:compDecompedF}, note that we get additional terms contributed by $\nabla^\perp \psi_{q+1} \cdot \nabla$. Since $\nabla^\perp \psi_{q+1} \cdot \nabla$ gives a factor $\delta^{\sfrac12}\lambda$, we can inductively deduce \eqref{e:22:compDecompedF}-\eqref{e:24:compDecompedF}, similar to the proof of Lemma \ref{l:streamEst}. Here, we also use $\delta^{\sfrac12} \lambda > \mu$ in the parameter relation of Remark \ref{r:auxiliaryScalesRela}.

Similarly, one can show that, for any $\bbeta \in \mathcal{I}_* $,
\begin{align*}
      \| O^{\bbeta} \mr {\bar D}^{\balpha} \rho \|_\infty
      \lesssim \bar \lambda^{-[\balpha] b \gamma} \lambda^{|\bbeta|_x} 
            ( \delta^{\sfrac12} \lambda )^{|\bbeta|_t},
            \quad [\bbeta] + [\balpha] \leq 11Q^3. 
\end{align*}
This leads to 
\begin{align}     \label{e:26:compDecompedF}
      \| \mr O^{\bbeta} \mr {\bar D}^{\balpha} \rho \|_\infty 
      \lesssim \bar \lambda^{-[\balpha] b \gamma} \cdot \bigg( \frac{\lambda}{\mr\lambda} \bigg)^{|\bbeta|_x} \bigg( \frac{ \delta^{\sfrac12}\lambda } {\mr\mu} \bigg)^{|\bbeta|_t},
            \quad [\bbeta] + [\balpha] \leq 11Q^3. 
\end{align}
Applying Corollary \ref{c:CalPEst}, we have
\begin{align}     \label{e:28:compDecompedF}
      \| \mr O^{\bbeta} \tilde F \|_{\infty} 
      &\, \lesssim \bar \lambda^{-[\balpha] b \gamma} \cdot \bigg( \frac{\lambda}{\mr\lambda} \bigg)^{|\bbeta|_x} \bigg( \frac{ \delta^{\sfrac12}\lambda } {\mr\mu} \bigg)^{|\bbeta|_t},
      \quad [\bbeta] + [\balpha] \leq 11Q^3.
\end{align}
\end{proof}

In the next two lemmas, we estimate the force $\tilde f$ for the remainder equations in Lemma \ref{l:EXP1} and Lemma \ref{l:EXP2}.

\begin{lemma}     \label{l:remainderEst_InitS}
Given the assumptions in Lemma \ref{l:EXP1}, suppose that, for any $\bbeta \in \mathcal{I}_*$ with $[\bbeta] \leq 5Q^3$,
\begin{align}     \label{e:2:remainderEst_InitS}
      {\lambda^{-a}} \| \mr{\bar D}^{\bbeta} \bar \rho_{a} \|_\infty
      \lesssim \bar \lambda^{ -\alpha - ( a + [\bbeta] ) b \gamma } \bigg( | \rho_{\ini} |_{\mathfrak{D}(\bar u, \bar \kappa)} 
                  + \sum_{j=q_*}^{q} \lambda_j^{-\sfrac{\gamma}{2}} | \rho_{\ini} |_{\fS(j)} \bigg).
\end{align}
Then for any $\bbeta$ with $[\bbeta] \leq 4Q^2$,
\begin{align}     \label{e:6:remainderEst_InitS}
      \|\mr O^{\bbeta} \tilde f\|_\infty
      \lesssim \lambda^{- (Q\gamma-1)} \bigg( | \rho_{\ini} |_{\mathfrak{D}(\bar u, \bar \kappa)} 
                  + \sum_{j=q_*}^{q} \lambda_j^{-\sfrac{\gamma}{2}} | \rho_{\ini} |_{\fS(j)} \bigg).
\end{align}
\end{lemma}

\begin{proof}
Given all the information in \ref{c:l:exp1:spt}, \ref{c:l:exp1:tpr} and \ref{c:l:exp1:rsd} in Lemma \ref{l:EXP1}, we estimate $\tilde f$ in \eqref{e:eTildeRho} as follows. We expand $\tilde f$ 
\begin{align}     
      \tilde f =&\, -\frac{1}{\lambda^{Q}} \Big( L_{11} ( \hat\rho_{Q} + \check\rho_{Q} + \bar\rho_{Q} ) + L_{12} ( \hat\rho_{Q} + \check\rho_{Q} ) + L_{13} \hat\rho_{Q} \Big)   \label{e:16:remainderEst_InitS} \\ 
            &\, -\frac{1}{\lambda^{Q}} \Big( ( L_{01} + L_{02} + L_{03} ) ( \hat\rho_{Q-1} + \check\rho_{Q-1} + \bar\rho_{Q-1} ) + (L_{04} + L_{05} + L_{06} + L_{07}) \hat\rho_{Q-1} \Big)    \label{e:17:remainderEst_InitS} \\ 
            &\, - \frac{1}{\lambda^{Q+1}} \Big( ( L_{01} + L_{02} + L_{03} ) ( \hat\rho_{Q} + \check\rho_{Q} + \bar\rho_{Q} ) + (L_{04} + L_{05} + L_{06} + L_{07}) \hat\rho_{Q} \Big)     \label{e:18:remainderEst_InitS} \\ 
            &\,- \frac{1}{\lambda^{Q+2}}  \Big( \lambda^2 L_{-1} \hat \rho_{Q-2} + \lambda L_{-1} \hat \rho_{Q-1} + L_{-1} \hat \rho_{Q} \Big)   \label{e:19:remainderEst_InitS}
\end{align}
Note that the argument on the right hand side of (\ref{e:16:remainderEst_InitS}-\ref{e:19:remainderEst_InitS}) is $(x,t,\lambda \Phi(x,t), \mu t)$.

We estimate $\tilde f$ presented in (\ref{e:16:remainderEst_InitS}-\ref{e:19:remainderEst_InitS}). All terms are estimated in the same way. Here we take the term $L_{13} \hat\rho_{Q}$ as a typical example. First, we write 
\begin{align*}
      \frac{1}{\lambda^{Q}} L_{13} \hat\rho_{Q} = \varepsilon \kappa \lambda \cdot \frac{1}{\lambda^{Q}}
      \sum_{m=1}^N S_{m,11} \vartheta_m 
      \sum_{a=0}^{Q-1} \sum_{[\balpha]=1, \balpha \in \mathcal{I}_*}^{2(Q-a)-1} \sum_{l=1}^{\hat I_{Q-a}^{(\balpha)}} 
             \hat h_{Q-a}^{(\balpha,l)} \hat \eta_{Q-a}^{(\balpha,l)} \hat \chi_{Q-a}^{(\balpha,l)} \mr {\bar D}^{\balpha} \bar \rho_a.
\end{align*}
From the estimate \eqref{e:EXP1:hhS_Pest}, we have
\begin{align}
      S_{m,11} \hat h_{Q-a}^{(\balpha,l)}  
            &\,\in \mathcal{\bar P} ( (Q-a)N + m ),    \nonumber \\ 
      \sum_{m=1}^N \sum_{[\balpha]=1, \balpha \in \mathcal{I}_*}^{2(Q-a)-1} \sum_{l=1}^{\hat I_{Q-a}^{(\balpha)}}
      \frac{\varepsilon \kappa \lambda} {\lambda^{Q-a}} 
            &\,\vertiii{ S_{m,11} \hat h_{Q-a}^{(\balpha,l)} }
            \lesssim \varepsilon \kappa \lambda \lambda^{-2(Q-a)\gamma}.      \label{e:12:remainderEst_InitS}
\end{align}
From \eqref{e:prdEstHom}, \eqref{e:EXP1:hchi_est} and \eqref{e:EXP1:heta_est}, we have, for $p \leq 8Q^3$,
\begin{align*}
      \| \partial_\tau^p ( \vartheta_m \hat \eta_{Q-a}^{(\balpha,l)} ) \|_\infty \lesssim&\, 1,    \\ 
      \| \nabla^p_\xi \hat \chi_{Q-a}^{(\balpha,l)} \|_\infty \lesssim&\, 1.
\end{align*}
We apply Lemma \ref{l:compDecompedF} to deduce that, for any $[\bbeta] \leq Q + 3Q^2$,
\begin{align}     \label{e:20:remainderEst_InitS}
      \| \mr O^{\bbeta} \tilde f \|_\infty  
      \lesssim \sum_{a=0}^{Q-1} \bar \lambda^{ - \alpha - ( a + [\bbeta] ) b \gamma } 
            \cdot \varepsilon \kappa \lambda \lambda^{- 2 (Q-a)\gamma} 
            \bigg( | \rho_{\ini} |_{\mathfrak{D}(\bar u, \bar \kappa)} 
                  + \sum_{j=q_*}^{q} \lambda_j^{-\sfrac{\gamma}{2}} | \rho_{\ini} |_{\fS(j)} \bigg).
\end{align}
The major gain comes from the factor $\lambda^{-Q\gamma}$ in \eqref{e:12:remainderEst_InitS}. From direct computation and parameter relation in Section \ref{ss:parameters}, we can deduce \eqref{e:16:remainderEst_InitS}.
\end{proof}

\begin{lemma}     \label{l:remainderEst_IndtS}
Suppose $P \geq Q^3$. Given the assumptions in Lemma \ref{l:EXP2}, suppose we have, for some $s \in \N^+$ and any $\bbeta \in \mathcal{I}_*$ with $[\bbeta] \leq 5Q^3$,
\begin{align*}
      \| \mr{\bar D}^{\bbeta} \bar \rho_{F,k} \| 
            \lesssim&\, \bar \lambda^{ -\alpha-[\bbeta]b\gamma } \lambda^{-(s-1)b\gamma},     \\
      \| \mr{\bar D}^{\bbeta} \bar \rho_0 \| 
            \lesssim&\, \bar \lambda^{ -\alpha-[\bbeta]b\gamma } \lambda^{-(s-1)b\gamma},     \\ 
      {\lambda^{-a}} \| \mr{\bar D}^{\bbeta} \bar \rho_{a} \| 
            \lesssim&\, \bar \lambda^{ -\alpha-[\bbeta]b\gamma } \lambda^{-(s-2+a)b\gamma},
            \quad a \geq 1.
\end{align*}
Then for any $\bbeta$ with $[\bbeta] \leq 4Q^2$,
\begin{align}
      \| \mr O^{\bbeta} \tilde f \|
      \lesssim \lambda^{- ([\bbeta]+s-1) b\gamma } \lambda^{-(Q\gamma-1)}.
\end{align}
\end{lemma}

\begin{proof}
As mentioned in \ref{c:l:exp2:rmd} of Lemma \ref{l:EXP2}, the $\tilde f$ is represented by (\ref{e:16:remainderEst_InitS}-\ref{e:19:remainderEst_InitS}), but the correctors $\hat \rho_Q, \check \rho_Q, \bar \rho_Q, \hat \rho_{Q-1}, \check \rho_{Q-1}, \bar \rho_{Q-1}, \hat \rho_{Q-2}$ are given in \ref{c:l:exp2:spt}, \ref{c:l:exp2:tpr}, \ref{c:l:exp2:rsd}. This proof is analogous to the proof of Lemma \ref{l:remainderEst_InitS}. Again, we take the term $L_{13} \hat\rho_{Q}$ as a typical example. Consider
\begin{align}
      \frac{1}{\lambda^{Q}} L_{13} \hat\rho_{Q} =&\, \varepsilon \kappa \lambda \cdot \frac{1}{\lambda^{Q}}
      \sum_{m=1}^N S_{m,11} \vartheta_m 
      \sum_{a=0}^{Q-1} \sum_{[\balpha]=1, \balpha \in \mathcal{I}_*}^{2(Q-a)-1} \sum_{l=1}^{\hat I_{Q-a}^{(\balpha)}} 
             \hat h_{Q-a}^{(\balpha,l)} \hat \eta_{Q-a}^{(\balpha,l)} \hat \chi_{Q-a}^{(\balpha,l)} \mr {\bar D}^{\balpha} \bar \rho_a   \label{e:10:remainderEst_IndtS} \\ 
      +&\, \varepsilon \kappa \lambda \cdot \frac{1}{\lambda^{Q}}
      \sum_{m=1}^N S_{m,11} \vartheta_m 
      \sum_{k=1}^K \sum_{ [\balpha]=1, \balpha \in \mathcal{I}_* }^{3Q-2} 
            \sum_{l=1}^{\hat I^{(\balpha)}_{F,n}}
            \hat h^{(\balpha,l)}_{F,k,Q} \hat \eta^{(\balpha,l)}_{F,k,Q} \hat \chi^{(\balpha,l)}_{F,k,Q} \mr {\bar D}^{\balpha} \bar \rho_{F,k}.      \label{e:12:remainderEst_IndtS}
\end{align}
The analysis of \eqref{e:10:remainderEst_IndtS} is the same as Lemma \ref{l:remainderEst_InitS}. For \eqref{e:12:remainderEst_IndtS}, we have
\begin{align}
      S_{m,11} \hat h_{F,k,Q}^{(\balpha,l)}  
            &\,\in \mathcal{\bar P} \big( N (Q+1) + 1 \big),    \\ 
      \sum_{m=1}^N \sum_{[\balpha]=1, \balpha \in \mathcal{I}_*}^{2(Q-a)-1} \sum_{l=1}^{\hat I_{Q-a}^{(\balpha)}}
      \frac{\varepsilon \kappa \lambda} {\lambda^Q} 
            &\,\vertiii{ S_{m,11} \hat h_{F,k,Q}^{(\balpha,l)} }
            \lesssim \varepsilon \kappa \lambda \lambda^{-2(Q-2)\gamma}. 
\end{align}
From \eqref{e:prdEstHom}, \eqref{e:EXP1:heta_est} and \eqref{e:EXP1:hchi_est}, we have
\begin{align}
      \| \partial_\tau^p ( \vartheta_m \hat \eta_{F,k,Q}^{(\balpha,l)} ) \|_\infty 
            \lesssim&\, 1,       \quad p \leq P-Q,    \\ 
      \| \nabla^p_\xi \hat \chi_{F,k,Q}^{(\balpha,l)} \|_\infty
            \lesssim&\, 1,       \quad p \leq P.
\end{align}
Then we apply Corollary \ref{c:CalPEst} and Lemma \ref{l:compDecompedF} to deduce the following estimate for any $\bbeta$ with $[\bbeta] \leq 4Q^2$,
\begin{align*}
      \|\mr O^{\bbeta} ( \ref{e:12:remainderEst_IndtS} ) \|_\infty  
      \lesssim&\, \bar \lambda^{ - \alpha } 
            \lambda^{ - (s-1+[\bbeta]) b \gamma }
            \cdot \varepsilon \kappa \lambda \lambda^{- 2 (Q-2) \gamma}.
\end{align*}
Note that $Q$ is a large integer from Section \ref{ss:parameters}, we can deduce the desired estimate.
\end{proof}

\subsection{The first expansion at step $q+1$}        \label{ss:reEXP1}

In this section, we rearrange the expansion formula in Lemma \ref{l:EXP1} to retrieve a simplified expansion. We also prove some related estimates for the simplified expansion.

\begin{lemma}     \label{l:reEXP1}
Suppose that for any $\rho_{\ini} \in C^\infty (\T^2)$, we have
\begin{align}
      \resN{ \rho_{\ini} }_{\bar u,\bar \kappa} 
            \leq&\, \bar\lambda^{-\alpha} \bigg( | \rho_{\ini} |_{\mathfrak{D}(\bar u, \bar \kappa)} 
                  + \sum_{j=q_*}^{q} \lambda_j^{-\sfrac{\gamma}{2}} | \rho_{\ini} |_{\fS(j)} \bigg),     \label{e:1:reEXP1}\\
      \resN{ \rho_{\ini} }_{\bar u,\bar \kappa} ^{(H,\alpha_0)}
            \leq&\, \left( 1 - \frac{4}{q} \right) \bigg( | \rho_{\ini} |_{\mathfrak{D}(\bar u, \bar \kappa)} 
                  + \sum_{j=q_*}^{q} \lambda_j^{-\sfrac{\gamma}{2}} | \rho_{\ini} |_{\fS(j)} \bigg).    \label{e:2:reEXP1}
\end{align}
Now we fix some $ \rho_{\ini} \in C^\infty (\T^2) $ in Assumption \ref{a:OEq}, consider the solution $\varrho_0$ in \eqref{e:5:OEq_Idct}-\eqref{e:6:OEq_Idct}. Then $\varrho_0$ admits the expansion
\begin{align}     \label{e:exp:reEXP1}
      \varrho_{0}(x,t) = \sum_{a=0}^Q \sum_{ [\balpha]=1, \balpha \in \mathcal{I}_* }^{Q} \mr {\bar D}^{\balpha} \bar \rho_{0,a} 
            \sum_{l=1}^{I_0} h_{0,a}^{(\balpha,l)} \eta_{0,a}^{(\balpha,l)} (\mu t) \chi_{0,a}^{(\balpha,l)} (\lambda\Phi) 
            + \sum_{a=0}^Q \bar\rho_{0,a} + \tilde \rho_{0}
\end{align}
with
\begin{align}
      h^{(\balpha,l)}_{0,a} \in&\, \bar{\mathcal{P}} ( Q ),       \qquad 
      \sum_{a,\balpha,l}
      \vertiii{ h^{(\balpha,l)}_{0,a} }
            \lesssim \frac{\delta^{\sfrac12}\mr{\bar\lambda}}{\kappa\lambda^2}.    \label{e:h:reEXP1}
\end{align}
The pair $\big(\chi_{0,a}^{(\balpha,l)}, \eta_{0,a}^{(\balpha,l)} \big)$ satisfies the generalized shear condition for all $a,\balpha,l$. The correctors $\chi_{0,a}^{(\cdot,\cdot)}: \T^2 \rightarrow \R$, $\eta_{0,a}^{(\cdot,\cdot)}: \T \rightarrow \R$ satisfy $I_0 \lesssim 1$ and
\begin{align}
      \| \partial_\tau^p \eta_{0,a}^{(\balpha,l)} \|_\infty &\,\lesssim 1,        
            \quad p \leq 8Q^3-Q,     \label{e:etaEst:reEXP1} \\ 
      \| \nabla_\xi^p \chi_{0,a}^{(\balpha,l)} \|_\infty &\,\lesssim 1,
            \quad p \leq 8Q^3.     \label{e:chiEst:reEXP1} 
\end{align}
The collection $\{\bar\rho_{0,a}\}_a$ satisfies the transmission condition, and for any $0 \leq a \leq Q$,
\begin{align}
      \| \bar\rho_{0,a} (\cdot,0) \|_2 + \sum_{j=q_*}^{q} \lambda_j^{-\sfrac{\gamma}{2}} | \bar\rho_{0,a} (\cdot,0) |_{\fS(j)}
            \lesssim&\, \lambda^{-ab\gamma} \bigg( | \rho_{\ini} |_{\mathfrak{D}(\bar u, \bar \kappa)} 
                  + \sum_{j=q_*}^{q} \lambda_j^{-\sfrac{\gamma}{2}} | \rho_{\ini} |_{\fS(j)} \bigg).        \label{e:6:reEXP1}
\end{align}
The following estimates hold,
\begin{align}
      \resF{ \bar\rho_{0,a} }_{\bar u, \bar \kappa, a} 
            \lesssim&\, \bar \lambda^{-\alpha} \lambda^{-ab\gamma} \bigg( | \rho_{\ini} |_{\mathfrak{D}(\bar u, \bar \kappa)}
                  + \sum_{j=q_*}^{q} \lambda_j^{-\sfrac{\gamma}{2}} | \rho_{\ini} |_{\fS(j)} \bigg),      \label{e:8:reEXP1} \\ 
      \resF{ \bar\rho_{0,a} }_{\bar u, \bar \kappa, a} ^{(H,\alpha_0)}
            \lesssim&\, \lambda^{-ab\gamma} \bigg( | \rho_{\ini} |_{\mathfrak{D}(\bar u, \bar \kappa)}
                  + \sum_{j=q_*}^{q} \lambda_j^{-\sfrac{\gamma}{2}} | \rho_{\ini} |_{\fS(j)} \bigg),      \label{e:holder:reEXP1} \\     
      \resF{ \bar\rho_{0,0} }_{\bar u, \bar \kappa, 0} ^{(H,\alpha_0)}
            \leq&\, \left( 1 - \frac{4}{q} \right) \bigg( | \rho_{\ini} |_{\mathfrak{D}(\bar u, \bar \kappa)} 
                  + \sum_{j=q_*}^{q} \lambda_j^{-\sfrac{\gamma}{2}} | \rho_{\ini} |_{\fS(j)} \bigg).       \label{e:const:reEXP1}
\end{align}
Moreover, for any $\bbeta$ with $[\bbeta] \leq 2Q + 3Q^2$, we have
\begin{align}     \label{e:9:reEXP1}
      \| \mr O^{\bbeta} \tilde \rho_0 \|
      + \kappa^{\sfrac12} \mr\lambda \| \mr O^{\bbeta} \mr \nabla \tilde \rho_0 \|_2  
            \lesssim \lambda^{-3-[\bbeta]b\gamma} \bigg( | \rho_{\ini} |_{\mathfrak{D}(\bar u, \bar \kappa)}
                  + \sum_{j=q_*}^{q+1} \lambda_j^{-\sfrac{\gamma}{2}} | \rho_{\ini} |_{\fS(j)} \bigg).
\end{align}
\end{lemma}

\begin{lemma}     \label{l:enReEXP1}
In the setting of Lemma \ref{l:reEXP1}, we have
\begin{align}
      h_{0,0}^{(\balpha,l)} = - \frac{\delta^{\sfrac12}\mr{\bar\lambda}}{\kappa\lambda^2}, \quad
      \eta_{0,0}^{(\balpha,l)} = \eta_2, \quad
      \chi_{0,0}^{(\balpha,l)} = - \partial_{\xi_2}^{-1} \Pi_2,  \quad \text{for } \balpha = 1, l=1.        \label{e:2:enReEXP1} \\ 
      h_{0,0}^{(\balpha,l)} = \frac{\delta^{\sfrac12}\mr{\bar\lambda}}{\kappa\lambda^2}, \quad
      \eta_{0,0}^{(\balpha,l)} = \eta_1, \quad
      \chi_{0,0}^{(\balpha,l)} = - \partial_{\xi_1}^{-1} \Pi_1,  \quad \text{for } \balpha = 2, l=1.        \label{e:4:enReEXP1}
\end{align}
The following estimate holds
\begin{align}     \label{e:6:enReEXP1}
      \| h_{0,0}^{(\balpha,l)} \|_{\infty} \leq 
            \frac{\delta^{\sfrac12}\mr{\bar\lambda}}{\kappa\lambda^2} \cdot \lambda^{-b\gamma}, \quad \text{for } a = 0, ( \balpha, l ) \notin \{ (1,1), (2,1) \}.
\end{align}
Moreover, $\bar\rho_{0,0}$ satisfies
\begin{align}     \label{e:8:enReEXP1}
      \bar D_t \bar\rho_{0,0} - \bar\kappa \Delta \bar\rho_{0,0} = 0,      \quad \bar\rho_{0,0}(\cdot,0) = \rho_{\ini}(\cdot).
\end{align}
\end{lemma}

\begin{remark}    \label{r:holderConst2}
Note that in \eqref{e:2:reEXP1} and \eqref{e:const:reEXP1}, we have $\leq$ instead of $\lesssim$. We need these two precise estimates for the global iterative scheme to converge. In contrast, the constants in the other estimates are not important, since we can absorb them later.
\end{remark}

In the rest of this section, we prove two lemmas above. Specifically, \eqref{e:exp:reEXP1}-\eqref{e:6:reEXP1} and Lemma \ref{l:enReEXP1} are proved in Section \ref{sss:p1:reEXP1}. \eqref{e:8:reEXP1} and \eqref{e:9:reEXP1} are proved in Section \ref{sss:p2:reEXP1} and Section \ref{sss:p3:reEXP1} respectively.

\subsubsection{Proof of Lemma \ref{l:reEXP1} and Lemma \ref{l:enReEXP1}: Rearranging the expansion}      \label{sss:p1:reEXP1}

In this section, we prove \eqref{e:exp:reEXP1}-\eqref{e:6:reEXP1}. From the notation in Definition \ref{d:homParaI} and Definition \ref{d:homFlowmap}, \eqref{e:hom:meq} and (\ref{e:5:OEq_Idct}-\ref{e:6:OEq_Idct}) are equivalent.

From Lemma \ref{l:EXP1}, we have an expansion of $\varrho_0 := \varrho$ given by \eqref{e:0:EXP1}, \eqref{e:EXP1:hrho0}, \eqref{e:EXP1:hrhon}, \eqref{e:EXP1:crho0} and \eqref{e:EXP1:crhon}. In this proof, we omit the arguments $\mu t$ for $\eta$ and $\lambda \Phi$ for $\chi$. For fixed $0 \leq n \leq Q$ and $0 \leq a \leq n-1$, consider the terms in this expansion involving $\bar \rho_a$, i.e.
\begin{align}
      \sum_{n=a+1}^{Q} &\,\frac{1}{\lambda^n} \Bigg( \sum_{[\balpha]=1, \balpha \in \mathcal{I}_*}^{2(n-a)-1} \sum_{l=1}^{\hat I_{n-a}^{(\balpha)}}
            \hat h^{(\balpha, l)}_{n-a} \hat \eta^{(\balpha, l)}_{n-a} \hat \chi^{(\balpha, l)}_{n-a} \mr {\bar D}^{\balpha} \bar \rho_a 
            + \sum_{[\balpha]=1, \balpha \in \mathcal{I}_*}^{2(n-a)} \sum_{l=1}^{\check I_{n-a}^{(\balpha)}}
            \check h^{(\balpha, l)}_{n-a} \check \eta^{(\balpha, l)}_{n-a} \mr {\bar D}^{\balpha} \bar \rho_a \Bigg)    \label{e:11:reEXP1} \\ 
            + &\, \frac{1}{\lambda^a} \bar \rho_a     \label{e:12:reEXP1}
\end{align}
Letting $i = n-a$ and setting $\hat I_i^{(\balpha)} = 0$ for $[\balpha] = 2i$, we can arrange the terms in \eqref{e:11:reEXP1} as
\begin{align}
      \eqref{e:11:reEXP1}
      &\,=\sum_{i=1}^{Q-a} \frac{1}{\lambda^{i+a}} \Bigg( \sum_{ [\balpha]=1, \balpha \in \mathcal{I}_* }^{2i} \sum_{l=1}^{\hat I_{i}^{(\balpha)}}
            \hat h^{(\balpha, l)}_i \hat \eta^{(\balpha, l)}_i \hat \chi^{(\balpha, l)}_i \mr {\bar D}^{\balpha} \bar \rho_a 
      + \sum_{ [\balpha]=1, \balpha \in \mathcal{I}_* }^{2i} \sum_{l=1}^{\check I_i^{(\balpha)}}
            \check h^{(\balpha, l)}_i \check \eta^{(\balpha, l)}_i \mr {\bar D}^{\balpha} \bar \rho_a \Bigg)      \label{e:14:reEXP1} \\ 
      &\,= \frac{1}{\lambda^a} \sum_{ [\balpha]=1, \balpha \in \mathcal{I}_* }^{2(Q-a)} \mr {\bar D}^{\balpha} \bar \rho_a 
            \sum_{i \geq \frac{[\balpha]}{2}}^{Q-a}
            \Bigg( \sum_{l=1}^{\hat I_{i}^{(\balpha)}} \frac{1}{\lambda^i} \hat h^{(\balpha, l)}_i \hat \eta^{(\balpha, l)}_i \hat \chi^{(\balpha, l)}_i 
            + \sum_{l=1}^{\check I_{i}^{(\balpha)}} \frac{1}{\lambda^i} \check h^{(\balpha, l)}_i \check \eta^{(\balpha, l)}_i \Bigg)      \label{e:16:reEXP1} \\ 
      &\,= \frac{1}{\lambda^a} \sum_{ [\balpha]=1, \balpha \in \mathcal{I}_* }^{2(Q-a)} \mr {\bar D}^{\balpha} \bar \rho_a 
            \sum_{i \geq \frac{[\balpha]}{2}}^{Q-a} 
            \Bigg( \sum_{l=1}^{\hat I_{i}^{(\balpha)}} \frac{1}{\lambda^i} \hat h^{(\balpha, l)}_i \hat \eta^{(\balpha, l)}_i \hat \chi^{(\balpha, l)}_i 
            + \sum_{l=1}^{\check I_{i}^{(\balpha)}} \frac{1}{\lambda^i} \check h^{(\balpha, l)}_i \check \eta^{(\balpha, l)}_i \check \chi^{(\balpha, l)}_i \Bigg).   \label{e:18:reEXP1}
\end{align}
Here, from \eqref{e:16:reEXP1} to \eqref{e:18:reEXP1}, we define periodic functions $\check \chi^{(\balpha, l)}_i : \T^2 \rightarrow \R$ by
\begin{align}     \label{e:22:reEXP1}
      \check \chi^{(\balpha, l)}_i = 1.
\end{align}
Note that $\langle \check \chi^{(\balpha, l)}_i \rangle_\xi = 1$.

Next, we represent the terms in \eqref{e:18:reEXP1}. Define 
\begin{align}     \label{e:24:reEXP1}
      \bar\rho_{0,a} = \frac{1}{\lambda^a} \bar \rho_a. 
\end{align}
For fixed $a \geq 0$ and $\balpha \in \mathcal{I}_*$, we renumber the following collection of function triples
\begin{align*} 
      \bigg \{ 
            \bigg\{\frac{1}{\lambda^i} \hat h^{(\balpha, l)}_i, 
                  \hat \eta^{(\balpha, l)}_i, \hat \chi^{(\balpha, l)}_i 
                  \bigg\}_{1 \leq l \leq \hat I_{i}^{(\balpha)}}, 
            \bigg\{ \frac{1}{\lambda^i} \check h^{(\balpha, l)}_i, 
                  \check \eta^{(\balpha, l)}_i, \check \chi^{(\balpha, l)}_i 
                  \bigg\}_{1 \leq l \leq \check I_{i}^{(\balpha)}} 
                        \bigg\}_{ \frac{[\balpha]}{2} \leq i \leq Q-a }
\end{align*}
to be (we remind that below $h_{a}^{(\balpha,l)}$ is not a misprint)
\begin{align}     \label{e:32:reEXP1}
      \big\{ h_{a}^{(\balpha,l)}, \eta_{0,a}^{(\balpha,l)}, \chi_{0,a}^{(\balpha,l)} \big\}_{1 \leq l \leq I_0}
\end{align}
for some integer $I_0 \in \N$. It is clear from Lemma \ref{l:EXP1} that we can arrange this renumbering such that Lemma \ref{l:enReEXP1} holds. Here, we use the exact formula for the first two correctors in \eqref{e:2:enReEXP1} and \eqref{e:4:enReEXP1}. This information is not stated in Lemma \ref{l:EXP1}, but it is given in Section \ref{sss:solInitExp1}. The cardinality $I_0$ must be sufficiently large. Recalling \eqref{e:EXP1:hI_est} and \eqref{e:EXP1:cI_est}, we have 
\begin{align}     \label{e:34:reEXP1}
      \sum_{\frac{[\balpha]}{2} \leq i \leq Q-a} \hat I_{i}^{(\balpha)} + \check I_{i}^{(\balpha)}
      \leq 11(3N)^{Q-a-1}.
\end{align}
Hence, we set $I_0 = 11(3N)^{Q-1}$. In case that $I_0$ is larger than the left-hand side of \eqref{e:34:reEXP1}, we just set those correctors in additional slots of \eqref{e:32:reEXP1} to be zero.

With above representation of \eqref{e:11:reEXP1}, we can write the expansion in Lemma \ref{l:EXP1} to be
\begin{align}
      \varrho_{0}(x,t) =&\, \sum_{a=0}^Q \sum_{ [\balpha]=1, \balpha \in \mathcal{I}_* }^{2(Q-a)} \mr {\bar D}^{\balpha} \bar \rho_{0,a} 
            \sum_{l=1}^{I_0} h_{a}^{(\balpha,l)} \eta_{0,a}^{(\balpha,l)} (\mu t) \chi_{0,a}^{(\balpha,l)} (\lambda\Phi) 
            + \sum_{a=0}^Q \bar\rho_{0,a} + \tilde \rho     \label{e:38:reEXP1} \\ 
      =&\, \sum_{a=0}^Q \sum_{ [\balpha]=1, \balpha \in \mathcal{I}_* }^{Q} \mr {\bar D}^{\balpha} \bar \rho_{0,a} 
            \sum_{l=1}^{I_0} \bP_Q h_{a}^{(\balpha,l)} \eta_{0,a}^{(\balpha,l)} (\mu t) \chi_{0,a}^{(\balpha,l)} (\lambda\Phi) 
            + \sum_{a=0}^Q \bar\rho_{0,a} + \tilde \rho_{0}       \label{e:40:reEXP1}
\end{align}
with
\begin{align}
      \tilde \rho_{0} =&\, \sum_{a=0}^Q \sum_{ [\balpha]=Q+1, \balpha \in \mathcal{I}_* }^{2(Q-a)} \mr {\bar D}^{\balpha} \bar \rho_{0,a} 
            \sum_{l=1}^{I_0} h_{a}^{(\balpha,l)} \eta_{0,a}^{(\balpha,l)} (\mu t) \chi_{0,a}^{(\balpha,l)} (\lambda\Phi)    \label{e:42:reEXP1} \\ 
      +&\, \sum_{a=0}^Q \sum_{ [\balpha]=1, \balpha \in \mathcal{I}_* }^{Q} \mr {\bar D}^{\balpha} \bar \rho_{0,a} 
            \sum_{l=1}^{I_0} \Big( h_{a}^{(\balpha,l)} - \bP_Q h_{a}^{(\balpha,l)} \Big) \eta_{0,a}^{(\balpha,l)} (\mu t) \chi_{0,a}^{(\balpha,l)} (\lambda\Phi)    \label{e:43:reEXP1}\\ 
      +&\, \tilde \rho.        \label{e:44:reEXP1}
\end{align}

Now, defining $h_{0,a}^{(\balpha,l)} := \bP_Q h_{a}^{(\balpha,l)}$, \eqref{e:exp:reEXP1} follows from \eqref{e:40:reEXP1}. \eqref{e:h:reEXP1} follows from \eqref{e:homGammaRela}, \eqref{d:bigMu}, \eqref{e:EXP1:hh}, \eqref{e:EXP1:hh_Pest}, \eqref{e:EXP1:ch} and \eqref{e:EXP1:ch_Pest}. \eqref{e:etaEst:reEXP1} follows from \eqref{e:EXP1:heta_est} and \eqref{e:EXP1:ceta_est}. \eqref{e:chiEst:reEXP1} follows from \eqref{e:EXP1:hchi_est} and \eqref{e:22:reEXP1}. \eqref{e:6:reEXP1} follows from \eqref{e:24:reEXP1} and Lemma \ref{l:t=0:EXP1}.

\subsubsection{Proof of Lemma \ref{l:reEXP1}: Estimating macroscopic states}     \label{sss:p2:reEXP1}
From \eqref{d:eddyDiff}, \eqref{e:homGammaRela}, \eqref{e:EXP1:bh_Pest} and \eqref{e:EXP1:bhS_Pest}, we have $\bar{h}^{(\balpha)}_{n} \in \mathcal{\bar P} (Nn+3)$ and
\begin{align}     \label{e:62:reEXP1}
      \sum_{[\balpha]=1}^{2n+1} \frac{1}{\lambda^n} \vertiii{ \bar{h}^{(\balpha)}_{n} } \lesssim
            \bar\kappa \mr{\bar\lambda}^2 \lambda^{-2nb\gamma}.
\end{align}
Recalling \eqref{e:24:reEXP1}, \eqref{e:EXP1:brhon} for $n \geq 1$ becomes
\begin{align}     \label{e:64:reEXP1}
      \bar L_0 \bar \rho_{0,n} =&\, \sum_{a=0}^{n-1} \sum_{ [\balpha]=1, \balpha \in \mathcal{I}_* }^{2(n-a)+1} 
      \mr{\bar \divr} \bigg( \frac{1}{\lambda^n} \bar{h}^{(\balpha)}_{n-a} \mr {\bar D}^{\balpha} \bar \rho_{0,a} \bigg).
\end{align}

From Lemma \ref{l:P1resolveToF}, \eqref{e:EXP1:brho0}, \eqref{e:1:reEXP1} and \eqref{e:24:reEXP1}, we deduce
\begin{align}     \label{e:65:reEXP1}
      \resF{ \bar\rho_{0,0} }_{\bar u, \bar \kappa, 0} 
            \leq&\, \bar \lambda^{-\alpha} \bigg( | \rho_{\ini} |_{\mathfrak{D}(\bar u, \bar \kappa)} 
                  + \sum_{j=q_*}^{q} \lambda_j^{-\sfrac{\gamma}{2}} | \rho_{\ini} |_{\fS(j)} \bigg).
\end{align}
Similarly, from the alternative of Lemma \ref{l:P2resolveToF} for $\resF{\cdot}_{\cdot,\cdot,\cdot}^{(H,\alpha_0)}$, \eqref{e:EXP1:brho0}, \eqref{e:2:reEXP1} and \eqref{e:24:reEXP1}, we deduce \eqref{e:const:reEXP1}.

Next we proceed with an induction argument to prove \eqref{e:8:reEXP1}. Suppose \eqref{e:8:reEXP1} holds for any $0 \leq a \leq n$, i.e.
\begin{align}
      \resF{ \bar\rho_{0,a} }_{\bar u, \bar \kappa, a}
            \lesssim&\, \bar\lambda^{-\alpha} \lambda^{-ab\gamma} \bigg( | \rho_{\ini} |_{\mathfrak{D}(\bar u, \bar \kappa)}
                  + \sum_{j=q_*}^{q} \lambda_j^{-\sfrac{\gamma}{2}} | \rho_{\ini} |_{\fS(j)} \bigg).        \label{e:66:reEXP1}
\end{align}
From \eqref{e:1:reEXP1} and \eqref{e:6:reEXP1}, we have
\begin{align}     \label{e:67:reEXP1}
      \resN{ \bar\rho_{0,n+1} ( \cdot, 0 ) }_{\bar u, \bar \kappa}
      \lesssim \bar\lambda^{-\alpha} \lambda^{-(n+1)b\gamma} 
            \bigg( | \rho_{\ini} |_{\mathfrak{D}(\bar u, \bar \kappa)} 
                  + \sum_{j=q_*}^{q} \lambda_j^{-\sfrac{\gamma}{2}} | \rho_{\ini} |_{\fS(j)} \bigg).
\end{align}
Then from \eqref{e:62:reEXP1}, \eqref{e:64:reEXP1}, \eqref{e:66:reEXP1}, \eqref{e:67:reEXP1} and Lemma \ref{l:P2resolveToF}, we can deduce
\begin{align}     \label{e:70:reEXP1}
      \resF{ \bar\rho_{0,n+1} }_{\bar u, \bar \kappa, n+1}
            \lesssim \bar\lambda^{-\alpha} \lambda^{-(n+1)b\gamma} \bigg( | \rho_{\ini} |_{\mathfrak{D}(\bar u, \bar \kappa)}
                  + \sum_{j=q_*}^{q} \lambda_j^{-\sfrac{\gamma}{2}} | \rho_{\ini} |_{\fS(j)} \bigg).
\end{align}
This concludes the proof of \eqref{e:8:reEXP1}. The proof of \eqref{e:holder:reEXP1} is completely analogous.

\subsubsection{Proof of Lemma \ref{l:reEXP1}: Estimating the remainder}     \label{sss:p3:reEXP1}

To estimate the remainder $\tilde \rho_0$, we first estimate terms in \eqref{e:42:reEXP1} and \eqref{e:43:reEXP1}. From \eqref{e:8:reEXP1}, we have that, for any $\bbeta \in \mathcal{I}_* $ with $[\bbeta] \leq 11Q^3 - 3Qn$,
\begin{align}     \label{e:76:reEXP1}
      \| \mr {\bar D}^{\bbeta} \bar \rho_{0,n} \| \lesssim&\, 
            \bar \lambda^{ - \alpha - [\bbeta] b \gamma } \lambda^{ - n b \gamma } 
            \bigg( | \rho_{\ini} |_{\mathfrak{D}(\bar u, \bar \kappa)} 
                  + \sum_{j=q_*}^{q} \lambda_j^{-\sfrac{\gamma}{2}} | \rho_{\ini} |_{\fS(j)} \bigg).
\end{align}
Using \eqref{e:h:reEXP1}, \eqref{e:etaEst:reEXP1}, \eqref{e:chiEst:reEXP1}, \eqref{e:76:reEXP1} and applying Lemma \ref{l:compDecompedF}, we have the following estimate for any $\bbeta$ with $[\bbeta] \leq 5Q^3-3Q$,
\begin{align*}
      \kappa^{\sfrac12} \mr\lambda \| \mr O^{\bbeta} (\ref{e:42:reEXP1}) \| 
            \lesssim \bar \lambda^{ - \alpha } \lambda^{ - ( [\bbeta] + Q ) b\gamma } \cdot \kappa^{\sfrac12} \mr\lambda \cdot \frac{\delta^{\sfrac12}\mr{\bar\lambda}}{\kappa\lambda^2} 
            \bigg( | \rho_{\ini} |_{\mathfrak{D}(\bar u, \bar \kappa)} 
                  + \sum_{j=q_*}^{q} \lambda_j^{-\sfrac{\gamma}{2}} | \rho_{\ini} |_{\fS(j)} \bigg)
\end{align*}
and subsequently for any $\bbeta$ with $[\bbeta] \leq 4Q^3$,
\begin{align}     \label{e:80:reEXP1}
      \| \mr O^{\bbeta} (\ref{e:42:reEXP1}) \|
      + \kappa^{\sfrac12} \mr\lambda \| \mr O^{\bbeta} \mr\nabla (\ref{e:42:reEXP1}) \|_2 
            \lesssim \lambda^{ - 3 - [\bbeta] b\gamma } \bigg( | \rho_{\ini} |_{\mathfrak{D}(\bar u, \bar \kappa)}
                  + \sum_{j=q_*}^{q} \lambda_j^{-\sfrac{\gamma}{2}} | \rho_{\ini} |_{\fS(j)} \bigg).
\end{align}

Estimating \eqref{e:43:reEXP1} is similar. Note that the gain in negative power of $\lambda$ is from the lower bound on the differentiation order of $h_{a}^{(\balpha,l)} - \bP_Q h_{a}^{(\balpha,l)}$. Applying Corollary \ref{c:CalPEst} and Lemma \ref{l:compDecompedF}, we have that, for any $\bbeta$ with $[\bbeta] \leq 4Q^3$,
\begin{align}     \label{e:82:reEXP1}
      \| \mr O^{\bbeta} (\ref{e:43:reEXP1}) \|
      + \kappa^{\sfrac12} \mr\lambda \| \mr O^{\bbeta} \mr \nabla (\ref{e:43:reEXP1}) \|_2 
            \lesssim \lambda^{ - 3 - [\bbeta] b\gamma } \bigg( | \rho_{\ini} |_{\mathfrak{D}(\bar u, \bar \kappa)} 
                  + \sum_{j=q_*}^{q} \lambda_j^{-\sfrac{\gamma}{2}} | \rho_{\ini} |_{\fS(j)} \bigg).
\end{align}

For $\tilde \rho$ in \eqref{e:44:reEXP1}, we deduce, from Remark \ref{r:datumProj} and Lemma \ref{l:t=0:EXP1}, the following estimates on initial datum $\tilde \rho(\cdot,0)$,
\begin{align*}
      \| \mr{\bar\nabla}^p \tilde\rho (\cdot, 0) \| 
            \lesssim&\, \lambda^{-Q\gamma} | \rho_{\ini} |_{\fS(q+1)},     \quad p \leq 10Q^3,
\end{align*}
and from Lemma \ref{l:remainderEst_InitS}, the following estimates on $\tilde f$
\begin{align*}
      \| \mr O^{\bbeta} \tilde f \|_\infty 
            \lesssim \lambda^{-(Q\gamma-1)} \bigg( &\, | \rho_{\ini} |_{\mathfrak{D}(\bar u, \bar \kappa)} 
                  + \sum_{j=q_*}^{q} \lambda_j^{-\sfrac{\gamma}{2}} | \rho_{\ini} |_{\fS(j)} \bigg),   \quad [\bbeta] \leq 4Q^3.
\end{align*}
Then we apply the energy estimates Lemma \ref{l:L2_init_energy} and Lemma \ref{l:L2_F_energy}, in the style of Proposition \ref{p:unconEnergyHom}, to deduce the following estimates for any $\bbeta$ with $[\bbeta] \leq 2Q + 3Q^2$,
\begin{align}     \label{e:84:reEXP1}
      \| \mr O^{\bbeta} \tilde \rho \|
      + \kappa^{\sfrac12} \mr\lambda \| \mr O^{\bbeta} \mr \nabla \tilde \rho \|_2  
            \lesssim \kappa^{-\sfrac12} \lambda^{-(Q\gamma-1)} \bigg( | \rho_{\ini} |_{\mathfrak{D}(\bar u, \bar \kappa)} 
                  + \sum_{j=q_*}^{q+1} \lambda_j^{-\sfrac{\gamma}{2}} | \rho_{\ini} |_{\fS(j)} \bigg).
\end{align}

Finally, (\ref{e:42:reEXP1}-\ref{e:44:reEXP1}) with $Q\gamma > 3$ and (\ref{e:80:reEXP1}-\ref{e:84:reEXP1}) together give \eqref{e:9:reEXP1}.

\subsection{Dissipation estimates: Proof of Theorem \ref{t:homDissip}}    \label{ss:proofHomDissip}

In this section, we prove the estimate on dissipation.

\begin{proof}[Proof of Theorem \ref{t:homDissip}]
In view of \eqref{e:7:chainStab} in Proposition \ref{p:chainStab}, it suffices to prove
\begin{align}
      \big| | \rho_{\ini} |_{\mathfrak{D}(v, \kappa)} - | \rho_{\ini} |_{\mathfrak{D}(\bar u, \bar\kappa)} \big|      
            \leq&\, \lambda^{ - \alpha - 2 } 
            \bigg( | \rho_{\ini} |_{\mathfrak{D}(v,\kappa)} 
                  + \sum_{j=q_*}^{q+1} \lambda_j^{-\sfrac{\gamma}{2}} | \rho_{\ini} |_{\fS(j)} \bigg)
\end{align}
We work directly with the expansion formula \eqref{e:exp:reEXP1} in Lemma \ref{l:reEXP1}. We compute
\begin{align}     
      \nabla \varrho_{0} =&\, \lambda \sum_{a=0}^Q \sum_{ [\balpha]=1, \balpha \in \mathcal{I}_* }^{Q} \mr {\bar D}^{\balpha} \bar \rho_{0,a} 
            \sum_{l=1}^{I_0} h_{0,a}^{(\balpha,l)} \eta_{0,a}^{(\balpha,l)} (\mu t) \partial_{\xi_j} \chi_{0,a}^{(\balpha,l)} (\lambda\Phi) \nabla \Phi_j   \label{e:10:homDissip} \\ 
            +&\, \sum_{a=0}^Q \sum_{ [\balpha]=1, \balpha \in \mathcal{I}_* }^{Q} \nabla \mr {\bar D}^{\balpha} \bar \rho_{0,a} 
            \sum_{l=1}^{I_0} h_{0,a}^{(\balpha,l)} \eta_{0,a}^{(\balpha,l)} (\mu t) \chi_{0,a}^{(\balpha,l)} (\lambda\Phi)   \label{e:12:homDissip} \\
            +&\, \sum_{a=0}^Q \sum_{ [\balpha]=1, \balpha \in \mathcal{I}_* }^{Q} \mr {\bar D}^{\balpha} \bar \rho_{0,a} 
            \sum_{l=1}^{I_0} \nabla h_{0,a}^{(\balpha,l)} \eta_{0,a}^{(\balpha,l)} (\mu t) \chi_{0,a}^{(\balpha,l)} (\lambda\Phi)   \label{e:14:homDissip} \\
            +&\, \sum_{a=0}^Q \nabla \bar\rho_{0,a} + \nabla \tilde \rho_{0} .   \label{e:16:homDissip}
\end{align}
The goal is to show the leading term of $\kappa \iint |\nabla \varrho_0|^2 \, dx dt$ is contributed by
\begin{align}
      &\, \frac{\delta^{\sfrac12}} {\kappa\lambda} \partial_1 \bar \rho_{0,0} 
            \eta_{0,0}^{(1,1)} (\mu t) \partial_{\xi_j} \chi_{0,0}^{(1,1)} (\lambda\Phi) \nabla \Phi_j
      + \frac{\delta^{\sfrac12}} {\kappa\lambda} \partial_2 \bar \rho_{0,0} 
            \eta_{0,0}^{(2,1)} (\mu t) \partial_{\xi_j} \chi_{0,0}^{(2,1)} (\lambda\Phi) \nabla \Phi_j      \label{e:18:homDissip} \\ 
      = &\, \frac{\delta^{\sfrac12}} {\kappa\lambda} \partial_1 \bar \rho_{0,0} 
            \eta_{0,0}^{(1,1)} (\mu t) \partial_{\xi_2} \chi_{0,0}^{(1,1)} (\lambda\Phi) \nabla \Phi_2
      + \frac{\delta^{\sfrac12}} {\kappa\lambda} \partial_2 \bar \rho_{0,0} 
            \eta_{0,0}^{(2,1)} (\mu t) \partial_{\xi_1} \chi_{0,0}^{(2,1)} (\lambda\Phi) \nabla \Phi_1      \label{e:20:homDissip}
\end{align}
and the rest are small error terms. We shall show all error terms in $\kappa \iint |\nabla \varrho_0|^2 \, dx dt$ can be bounded by
\begin{align}     \label{e:21:homDissip}
      \lambda^{-b\gamma} \bigg( | \rho_{\ini} |_{\mathfrak{D}(v,\kappa)} 
                  + \sum_{j=q_*}^{q+1} \lambda_j^{-\sfrac{\gamma}{2}} | \rho_{\ini} |_{\fS(j)} \bigg)^2.
\end{align}

Here, we also remark that, using \eqref{e:alphaDef}, we have
\begin{align}     \label{e:17:homDissip}
      \bar\kappa \mr{\bar\lambda}^2 \bar\lambda^{-2\alpha} \lesssim \bar\lambda^{2b^2\gamma}.
\end{align}

\begin{step}[Estimate the leading term \eqref{e:18:homDissip}]

Notice that $\eta_{0,0}^{(1,1)}$ and $\eta_{0,0}^{(2,1)}$ have disjoint support, and $\nabla \Phi = \Id + \varepsilon O(1)$ given by \eqref{e:6:uEstimate}. We have
\begin{align}     
      \kappa \iint \eqref{e:18:homDissip}^2 \, dx dt 
            =&\, \frac{\delta} {\kappa\lambda^2} \iint |\partial_1 \bar \rho_{0,0}|^2 
            |\eta_{0,0}^{(1,1)} (\mu t)|^2 |\partial_{\xi_2} \chi_{0,0}^{(1,1)} (\lambda\Phi)|^2     \label{e:22:homDissip}         \\
            +&\, \frac{\delta} {\kappa\lambda^2} \iint |\partial_2 \bar \rho_{0,0}|^2 
            |\eta_{0,0}^{(2,1)} (\mu t)|^2 |\partial_{\xi_1} \chi_{0,0}^{(2,1)} (\lambda\Phi)|^2    \label{e:24:homDissip}         \\ 
            +&\, \frac{\delta} {\kappa\lambda^2} \iint |\nabla \bar \rho_{0,0}|^2 \cdot \varepsilon^2 O(1).     \label{e:26:homDissip}                      
\end{align}
Note that \eqref{e:26:homDissip} can be bounded by \eqref{e:21:homDissip} using \eqref{d:eddyDiff}, $\varepsilon^2 < \lambda^{-b\gamma}$ and integration by parts.

Then we estimate
\begin{align}     
      \eqref{e:22:homDissip} - \frac{\delta} {\kappa\lambda^2} \iint |\partial_1 \bar \rho_{0,0}|^2    
      = &\, \frac{\delta} {\kappa\lambda^2} \iint |\partial_1 \bar \rho_{0,0}|^2 
            |\eta_{0,0}^{(1,1)} (\mu t)|^2 \Big( |\partial_{\xi_2} \chi_{0,0}^{(1,1)} (\lambda\Phi)|^2 
            - \langle |\partial_{\xi_2} \chi_{0,0}^{(1,1)}|^2 \rangle_\xi \Big)     \label{e:28:homDissip}         \\ 
      + &\, \frac{\delta} {\kappa\lambda^2} \langle |\partial_{\xi_2} \chi_{0,0}^{(1,1)}|^2 \rangle_\xi \iint |\partial_1 \bar \rho_{0,0}|^2 
            \Big( |\eta_{0,0}^{(1,1)} (\mu t)|^2 - \langle |\eta_{0,0}^{(1,1)}|^2 \rangle_\tau \Big).     \label{e:30:homDissip} 
\end{align}
Notice that there exists a potential $f:\T^2 \rightarrow \R$ with $\partial_{\xi_1} f = 0$ and $\|f\|_\infty \lesssim 1$ such that
\begin{align*}
      \partial_{\xi_2} f = | \partial_{\xi_2} \chi_{0,0}^{(1,1)} ( \xi )|^2 
            - \langle |\partial_{\xi_2} \chi_{0,0}^{(1,1)}|^2 \rangle_\xi,
\end{align*}
hence
\begin{align*}
      \partial_2 \Big( \frac{1}{\lambda} f(\lambda \Phi) \Big) = \partial_{\xi_2} f (\lambda \Phi) \cdot \partial_2 \Phi_2.
\end{align*}
To estimate \eqref{e:28:homDissip}, we use again $\nabla \Phi = \Id + \varepsilon O(1)$ given by \eqref{e:6:uEstimate},
\begin{align}
      \eqref{e:28:homDissip} =&\, \frac{\delta} {\kappa\lambda^2} \iint |\partial_1 \bar \rho_{0,0}|^2 
            |\eta_{0,0}^{(1,1)} (\mu t)|^2 \partial_2 \Big( \frac{1}{\lambda} f(\lambda \Phi) \Big)
      + \frac{\delta} {\kappa\lambda^2} \iint |\partial_1 \bar \rho_{0,0}|^2 \cdot \varepsilon O(1)
            \label{e:36:homDissip}    \\ 
      =&\, - \frac{\delta} {\kappa\lambda^2} \iint 2 \partial_{21} \bar \rho_{0,0} \partial_1 \bar \rho_{0,0} 
            |\eta_{0,0}^{(1,1)} (\mu t)|^2 \cdot \frac{1}{\lambda} f(\lambda \Phi)
      + \frac{\delta} {\kappa\lambda^2} \iint |\partial_1 \bar \rho_{0,0}|^2 \cdot \varepsilon O(1). 
            \label{e:38:homDissip}                                                                                      
\end{align}
Using \eqref{d:eddyDiff} and \eqref{e:8:reEXP1}, the first term in \eqref{e:38:homDissip} can be estimated by
\begin{align}
      &\, \big| \text{the first term in } \eqref{e:38:homDissip} \big| 
      \lesssim \frac{\mr{\bar\lambda}}{\lambda} \cdot \bar\kappa \mr{\bar\lambda}^2 \bar\lambda^{-2\alpha}
            \bigg( | \rho_{\ini} |_{\mathfrak{D}(v,\kappa)} 
                  + \sum_{j=q_*}^{q+1} \lambda_j^{-\sfrac{\gamma}{2}} | \rho_{\ini} |_{\fS(j)} \bigg)^2.     \label{e:42:homDissip}         
\end{align}
Now we can use the parameter relations \eqref{e:smallGammaR}, \eqref{e:12:mScaleRela} to deduce that \eqref{e:42:homDissip} is bounded by \eqref{e:21:homDissip}.

To estimate \eqref{e:30:homDissip}, notice that there exists a potential $g:\T \rightarrow \R$ with $\|g\|_\infty \lesssim 1$ such that
\begin{align*}
      \partial_\tau g = |\eta_{0,0}^{(1,1)} ( \tau )|^2 - \langle |\eta_{0,0}^{(1,1)}|^2 \rangle_\tau,
\end{align*}
then
\begin{align*}
      \bar D_t \Big( \frac{1}{\mu} g(\mu t) \Big) = \partial_{\tau} g (\mu t).
\end{align*}
Now we estimate \eqref{e:30:homDissip} as
\begin{align}
      \eqref{e:30:homDissip} =&\, \frac{\delta} {\kappa\lambda^2} \langle |\partial_{\xi_2} \chi_{0,0}^{(1,1)}|^2 \rangle_\xi \iint |\partial_1 \bar \rho_{0,0}|^2 
            \bar D_t \Big( \frac{1}{\mu} g(\mu t) \Big)     \label{e:44:homDissip}         \\ 
      =&\, - \frac{\delta} {\kappa\lambda^2} \langle |\partial_{\xi_2} \chi_{0,0}^{(1,1)}|^2 \rangle_\xi \iint 2 \bar D_t \partial_1 \bar \rho_{0,0} \partial_1 \bar \rho_{0,0} 
            \cdot \frac{1}{\mu} g(\mu t)     \label{e:46:homDissip}     \\ 
      &\,+ \frac{\delta} {\kappa\lambda^2} \langle |\partial_{\xi_2} \chi_{0,0}^{(1,1)}|^2 \rangle_\xi \int |\partial_1 \bar \rho_{0,0}|^2 \cdot \frac{1}{\mu} g(\mu t) \bigg|_{t=0}^{t=1}.     \label{e:48:homDissip}
\end{align}
Using \eqref{d:eddyDiff} and \eqref{e:8:reEXP1}, we estimate \eqref{e:46:homDissip} and \eqref{e:48:homDissip} as
\begin{align}
      \big| \eqref{e:46:homDissip} \big| 
      \lesssim \frac{\mr{\bar\mu}}{\mu} \cdot \bar\kappa \mr{\bar\lambda}^2 \bar\lambda^{-2\alpha}
            \bigg( | \rho_{\ini} |_{\mathfrak{D}(v,\kappa)} 
                  + \sum_{j=q_*}^{q+1} \lambda_j^{-\sfrac{\gamma}{2}} | \rho_{\ini} |_{\fS(j)} \bigg)^2,     \label{e:50:homDissip} \\ 
      \big| \eqref{e:48:homDissip} \big| 
      \lesssim \frac{1}{\mu} \cdot \bar\kappa \mr{\bar\lambda}^2 \bar\lambda^{-2\alpha}
            \bigg( | \rho_{\ini} |_{\mathfrak{D}(v,\kappa)} 
                  + \sum_{j=q_*}^{q+1} \lambda_j^{-\sfrac{\gamma}{2}} | \rho_{\ini} |_{\fS(j)} \bigg)^2.     \label{e:52:homDissip}
\end{align}
From the parameter relations \eqref{e:smallGammaR}, \eqref{e:mrParameter} and \eqref{e:14:mScaleRela}, we deduce that \eqref{e:50:homDissip} and \eqref{e:52:homDissip} can be bounded by \eqref{e:21:homDissip}.

Combining above, we deduce that \eqref{e:26:homDissip}-\eqref{e:30:homDissip} can be bounded by \eqref{e:21:homDissip}.
\end{step}

\begin{step}[Estimate other error terms]
Now we consider the contributions from \eqref{e:10:homDissip}. Other terms in \eqref{e:12:homDissip}-\eqref{e:16:homDissip} and cross terms can be dealt similarly.
Apart from \eqref{e:18:homDissip}, other contributions in \eqref{e:10:homDissip} come from two type of terms:
\begin{align}
      \lambda \mr {\bar D}^{\balpha} \bar \rho_{0,a} 
            h_{0,a}^{(\balpha,l)} \eta_{0,a}^{(\balpha,l)} (\mu t) \partial_{\xi_j} \chi_{0,a}^{(\balpha,l)} (\lambda\Phi) \nabla \Phi_j,
            &\, \quad \text{for } a \geq 1,      \label{e:60:homDissip} \\
      \lambda \mr {\bar D}^{\balpha} \bar \rho_{0,0} 
            h_{0,0}^{(\balpha,l)} \eta_{0,0}^{(\balpha,l)} (\mu t) \partial_{\xi_j} \chi_{0,0}^{(\balpha,l)} (\lambda\Phi) \nabla \Phi_j,
            &\, \quad \text{for } ( \balpha, l ) \notin \{ (1,1), (2,1) \}.      \label{e:62:homDissip}
\end{align}
For terms like \eqref{e:60:homDissip}, we use \eqref{e:8:reEXP1} and \eqref{d:eddyDiff} to estimate
\begin{align}
      \kappa \iint \eqref{e:60:homDissip}^2 \, dx dt
      \lesssim&\, \lambda^{-4b\gamma} \cdot \bar\kappa \mr{\bar\lambda}^2 \bar\lambda^{-2\alpha}
            \bigg( | \rho_{\ini} |_{\mathfrak{D}(v,\kappa)} 
                  + \sum_{j=q_*}^{q+1} \lambda_j^{-\sfrac{\gamma}{2}} | \rho_{\ini} |_{\fS(j)} \bigg)^2,     \label{e:64:homDissip}   \\ 
      \kappa \iint \eqref{e:62:homDissip}^2 \, dx dt
      \lesssim&\, \lambda^{-4b\gamma} \cdot \bar\kappa \mr{\bar\lambda}^2 \bar\lambda^{-2\alpha}
            \bigg( | \rho_{\ini} |_{\mathfrak{D}(v,\kappa)} 
                  + \sum_{j=q_*}^{q+1} \lambda_j^{-\sfrac{\gamma}{2}} | \rho_{\ini} |_{\fS(j)} \bigg)^2.     \label{e:66:homDissip}
\end{align}
In the factor $\lambda^{-4b\gamma}$ of \eqref{e:64:homDissip}, half of the power comes from $a \geq 1$ and \eqref{e:8:reEXP1}. The other half comes from $[\balpha]_x \geq 1$. For the factor $\lambda^{-4b\gamma}$ of \eqref{e:66:homDissip}, half of the power comes from $( \balpha, l ) \notin \{ (1,1), (2,1) \}$, Lemma \ref{l:enReEXP1} and \eqref{e:8:reEXP1}. The other half comes from $[\balpha]_x \geq 1$.

\end{step}

Combining above and \eqref{e:17:homDissip}, we deduce that all error terms in $\kappa \iint |\nabla \varrho_0|^2 \, dx dt$ can be bounded by \eqref{e:21:homDissip}.

\end{proof}

\subsection{Reformulate the expansion lemma II}    \label{ss:generalExpF}

In this section, we rearrange the expansion Lemma \ref{l:EXP2} and prove some related estimates.

\begin{lemma}       \label{l:reEXP2}
Given \eqref{e:1:reEXP1}-\eqref{e:2:reEXP1}, taking Assumption \ref{a:homSmallscl} with some $P \geq Q^3$ and
\begin{align}
      \Upsilon = \frac{\delta^{\sfrac12}\mr{\bar\lambda}} {\kappa\lambda^2} &\, \cdot \lambda^{ - 2 b \gamma },      \label{e:Up:reEXP2}
\end{align}
suppose that, for some $\Upsilon_{\ini} > 0$, we have
\begin{align}
      \resF{ \bar\rho_{F,k} }_{\bar u, \bar \kappa, n_0} 
            \lesssim&\, \bar \lambda^{ - \alpha } \lambda^{ - (s-2) b\gamma } \Upsilon_{\ini},  \label{e:macro:reEXP2} \\ 
      \resF{ \bar\rho_{F,k} }_{\bar u, \bar \kappa, n_0}^{(H,\alpha_0)} 
            \lesssim&\, \lambda^{ - (s-2) b\gamma } \Upsilon_{\ini},  \label{e:macHld:reEXP2} \\ 
      \| \bar\rho_{F,k} (\cdot,0) \|_2 + \sum_{j=q_*}^{q} \lambda_j^{-\sfrac{\gamma}{2}} &\, | \bar\rho_{F,k} (\cdot,0) |_{\fS(j)}
            \lesssim \lambda^{ - (s-2) b\gamma } \Upsilon_{\ini}.     \label{e:macroIn:reEXP2} 
\end{align}
Then the solution $\varrho$ to the equation \eqref{e:hom:eqFast} admits the expansion
\begin{align}
      \varrho(x,t) =&\, \sum_{a=0}^{K+Q} \sum_{ [\balpha]=1, \balpha \in \mathcal{I}_* }^{Q} 
            \mr {\bar D}^{\balpha} \ddot{\rho}_{a}  
            \sum_{l=1}^{\mathfrak{I}} \ddot{h}_{a}^{(\balpha,l)} \ddot{\eta}_{a}^{(\balpha,l)} (\mu t) \ddot{\chi}_{a}^{(\balpha,l)} (\lambda\Phi) 
            + \sum_{a=0}^Q \ddot{\rho}_{a} + \tilde \rho_F  \label{e:E1:reEXP2} 
\end{align}
with $\mathfrak{I} \lesssim 1$,
\begin{align} 
      \ddot{h}_{k}^{(\balpha,l)} \in&\, \bar{\mathcal{P}} ( 2NQ ),      \qquad 
      \sum_{k,\balpha,l} \vertiii{ \ddot{h}_{k}^{(\balpha,l)} } 
            \lesssim \frac{\delta^{\sfrac12}\mr{\bar\lambda}}{\kappa\lambda^2}.     
            \label{e:hEst:reEXP2} 
\end{align}
The pair $(\ddot{\chi}_{a}^{(\balpha,l)}, \ddot{\eta}_{a}^{(\balpha,l)})$ satisfies the generalized shear condition for all $a,\balpha,l$, and
\begin{align}
      \| \nabla_\xi^p \ddot{\chi}_{F,k}^{(\balpha,l)} \|_\infty 
            \lesssim&\, 1,  \quad p \leq P,      \label{e:chiEst:reEXP2}  \\ 
      \| \partial_\tau^p \ddot{\eta}_{F,k}^{(\balpha,l)} \|_\infty 
            \lesssim&\, 1,  \quad p \leq P-Q.      \label{e:etaEst:reEXP2} 
\end{align}
Moreover, the collection $\{ \ddot{\rho}_{a} \}_a$ satisfies transmission condition, and for any $0 \leq a \leq K+Q$, we have
\begin{align}
      \| \ddot{\rho}_{a} (\cdot,0) \|_2 + \sum_{j=q_*}^{q} \lambda_j^{-\sfrac{\gamma}{2}} | \ddot{\rho}_{a} (\cdot,0) |_{\fS(j)}
            \lesssim&\, \lambda^{ - (a+s) b \gamma } \Upsilon_{\ini},   \label{e:4:reEXP2}
\end{align}
and
\begin{align}
      \resF{ \ddot{\rho}_{a} }_{\bar u, \bar \kappa, n_0+a+1} 
            \lesssim&\, \bar \lambda^{ - \alpha } \lambda^{ - (s-1) b \gamma }
                  \lambda^{ - (b-1)\gamma } \Upsilon_{\ini},       \label{e:6:reEXP2} \\ 
      \resF{ \ddot{\rho}_{a} }_{\bar u, \bar \kappa, n_0+a+1}^{(H,\alpha_0)} 
            \lesssim&\, \lambda^{ - (s-1) b \gamma }
                  \lambda^{ - (b-1)\gamma } \Upsilon_{\ini}.       \label{e:7:reEXP2}
\end{align}
For $\bbeta$ with $[\bbeta] \leq 2Q + 3Q^2$, we have
\begin{align}     \label{e:8:reEXP2}
      \| \mr O^{\bbeta} \tilde \rho_F \|
      + \kappa^{\sfrac12} \mr\lambda \| \mr O^{\bbeta} \mr \nabla \tilde \rho_F \|_2
            \lesssim&\, \lambda^{ - 3 - ( [\bbeta] + s - 1 ) b \gamma } \Upsilon_{\ini}.
\end{align}
\end{lemma}

The proof of Lemma \ref{l:reEXP2} is presented in the rest of this section and is similar to that of Lemma \ref{l:reEXP1} in Section \ref{ss:reEXP1}. In Section \ref{sss:reEXP2:reform}, we give an immediate reformulation of the expansion in Lemma \ref{l:EXP2}. In Section \ref{sss:p1:reEXP2}, we prove \eqref{e:E1:reEXP2}-\eqref{e:4:reEXP2}. In Section \ref{sss:p2:reEXP2}, we estimate the macroscopic states $\ddot \rho_a$, i.e. proving \eqref{e:6:reEXP2}. In Section \ref{sss:p3:reEXP2}, we estimate the remainder term $\tilde \rho_F$. i.e. proving \eqref{e:8:reEXP2}.

\subsubsection{Proof of Lemma \ref{l:reEXP2}: Reformulation}      \label{sss:reEXP2:reform}

Recalling Remark \ref{r:shareCorrector}, there are shared coefficients and functions in Lemma \ref{l:EXP1} and Lemma \ref{l:EXP2}, i.e. $\hat h_{n-a}^{(\balpha,l)}$, $\hat \eta_{n-a}^{(\balpha,l)}$, $\hat \chi_{n-a}^{(\balpha,l)}$, $\check h_{n-a}^{(\balpha,l)}$, $\check \eta_{n-a}^{(\balpha,l)}$ and $\bar h_{n-a}^{(\balpha,l)}$. However, the macroscopic states $\bar \rho_a, 0 \leq a \leq Q$ in Lemma \ref{l:EXP1} and Lemma \ref{l:EXP2} are different.

Collecting the expansion in Lemma \ref{l:EXP2}, we have 
\begin{align}
      \varrho = \sum_{n=1}^Q \frac{1}{\lambda^n} \Bigg( &\, \sum_{a=0}^{n-1} \sum_{ [\balpha]=1, \balpha \in \mathcal{I}_* }^{2(n-a)-1} \sum_{l=1}^{\hat I_{n-a}^{(\balpha)}}
            \hat h^{(\balpha, l)}_{n-a} \hat \eta^{(\balpha, l)}_{n-a} (\mu t) \hat \chi^{(\balpha, l)}_{n-a} (\lambda \Phi) \mr {\bar D}^{\balpha} \bar \rho_{a}   \label{e:24:reEXP2} \\ 
      +&\, \sum_{k=1}^K \sum_{ [\balpha]=1, \balpha \in \mathcal{I}_* }^{2n-2+Q} \sum_{l=1}^{\hat I^{(\balpha)}_{F,n}}
            \hat h^{(\balpha,l)}_{F,k,n} \hat \eta^{(\balpha,l)}_{F,k,n} (\mu t) \hat \chi^{(\balpha,l)}_{F,k,n} (\lambda \Phi) \mr {\bar D}^{\balpha} \bar \rho_{F,k}  \label{e:26:reEXP2} \\ 
      +&\, \sum_{a=0}^{n-1} \sum_{ [\balpha]=1, \balpha \in \mathcal{I}_* }^{2(n-a)} \sum_{l=1}^{\check I_{n-a}^{(\balpha)}}
            \check h^{(\balpha, l)}_{n-a} \check \eta^{(\balpha, l)}_{n-a} (\mu t) \mr {\bar D}^{\balpha} \bar \rho_{a}    \label{e:28:reEXP2} \\ 
      +&\, \sum_{k=1}^K \sum_{ [\balpha]=1, \balpha \in \mathcal{I}_* }^{2n-1+Q} \sum_{l=1}^{\check I^{(\balpha)}_{F,n}}  
            \check h^{(\balpha,l)}_{F,k,n} \check \eta^{(\balpha,l)}_{F,k,n} (\mu t) \mr {\bar D}^{\balpha} \bar \rho_{F,k}  \Bigg)      \label{e:30:reEXP2} \\ 
      +&\, \sum_{n=0}^Q \frac{1}{\lambda^n} \bar \rho_{n} + \tilde \rho.   \label{e:32:reEXP2}
\end{align}
From \eqref{e:EXP2:hhF}, \eqref{e:EXP2:hIF_est}, \eqref{e:EXP2:ch} and \eqref{e:EXP2:cIF_est}, we have $\hat I^{(\balpha)}_{F,n}, \check I^{(\balpha)}_{F,n} \lesssim 1 $ and
\begin{align}
      \hat h^{(\balpha,l)}_{F,k,n}, \check h^{(\balpha,l)}_{F,k,n} 
            \in&\, \bar {\mathcal{P}} ( Q + NQ ).    \label{e:36:reEXP2} 
\end{align}
From \eqref{e:macro:reEXP2}, \eqref{e:EXP2:hhSF_Pest} and \eqref{e:EXP2:chSF_Pest}, we have 
\begin{align}
      \sum_{l,\balpha,k} 
            \frac{1} {\lambda^n} \vertiii{ \hat h^{(\balpha,l)}_{F,k,j,n} } 
            &\, \lesssim \frac{\delta^{\sfrac12}\mr{\bar\lambda}}{\kappa\lambda^2} 
                  \cdot \lambda^{-2(n-1)\gamma},    
                  \quad n \geq 2,         \label{e:38:reEXP2} \\ 
      \sum_{l,\balpha,k} 
            \frac{1} {\lambda^n} \vertiii{ \check h^{(\balpha,l)}_{F,k,j,n} } 
            &\, \lesssim \frac{\delta^{\sfrac12}\mr{\bar\lambda}} {\mu}
            \frac{\mr{\bar\lambda}} {\lambda} 
                  \cdot \lambda^{-2n\gamma},     
                  \qquad n \geq 1.         \label{e:40:reEXP2}
\end{align}
From \eqref{e:EXP2:hchi_est}, \eqref{e:EXP2:heta_est} and \eqref{e:EXP2:ceta_est}, we have
\begin{align}
      \| \nabla_\xi^p \hat \chi^{(\balpha,l)}_{F,k,j,n} \|_\infty 
            &\, \lesssim 1,   \quad p \leq P,    \label{e:42:reEXP2} \\ 
      \| \partial_\tau^p \hat \eta^{(\balpha,l)}_{F,k,j,n} \|_\infty +  \| \partial_\tau^p \check \eta^{(\balpha,l)}_{F,k,j,n} \|_\infty 
            &\, \lesssim 1,   \quad p \leq P-Q.    \label{e:44:reEXP2} 
\end{align}

\subsubsection{Proof of Lemma \ref{l:reEXP2}: Rearrangement}    \label{sss:p1:reEXP2}

For the terms in \eqref{e:24:reEXP2} and \eqref{e:28:reEXP2} involving $\bar \rho_{a}$, we define
\begin{align}
      \ddot{\rho}_{a} := \frac{1}{\lambda^a} \bar \rho_a.   \label{e:96:reEXP2}
\end{align}
The rearrangement of \eqref{e:24:reEXP2} and \eqref{e:28:reEXP2} is completely analogous to the rearrangement for Lemma \ref{l:EXP1} in Section \ref{sss:p1:reEXP1}, so we omit the details. This leads to
\begin{align}
      \eqref{e:24:reEXP2} + \eqref{e:28:reEXP2}
      =&\, \sum_{a=0}^Q \sum_{ [\balpha]=1, \balpha \in \mathcal{I}_* }^{Q} \mr {\bar D}^{\balpha} \ddot \rho_{a} 
      \sum_{l=1}^{I_0} \bP h_{0,a}^{(\balpha,l)} \eta_{0,a}^{(\balpha,l)} (\mu t) \chi_{0,a}^{(\balpha,l)} (\lambda\Phi)    \label{e:98:0:reEXP2} \\ 
      +&\, \sum_{a=0}^Q \sum_{ [\balpha]=1, \balpha \in \mathcal{I}_* }^{Q} \mr {\bar D}^{\balpha} \ddot \rho_{a} 
      \sum_{l=1}^{I_0} \Big( h_{0,a}^{(\balpha,l)} - \bP h_{0,a}^{(\balpha,l)} \Big) \eta_{0,a}^{(\balpha,l)} (\mu t) \chi_{0,a}^{(\balpha,l)} (\lambda\Phi)    \label{e:98:2:reEXP2} \\ 
      +&\,\sum_{a=0}^Q \sum_{ [\balpha]=Q+1, \balpha \in \mathcal{I}_* }^{2(Q-a)} \mr {\bar D}^{\balpha} \ddot \rho_{a} 
      \sum_{l=1}^{I_0} h_{0,a}^{(\balpha,l)} \eta_{0,a}^{(\balpha,l)} (\mu t) \chi_{0,a}^{(\balpha,l)} (\lambda\Phi).    \label{e:98:4:reEXP2}
\end{align}
Here, $I_0$, $h_{0,a}^{(\balpha,l)}$, $\eta_{0,a}^{(\balpha,l)}$ and $\chi_{0,a}^{(\balpha,l)}$ are given in Lemma \ref{l:reEXP1}.

Next, we rearrange the $F$-related terms as follows, i.e. the terms in \eqref{e:26:reEXP2} and \eqref{e:30:reEXP2} involving $\bar\rho_{F,k}$.
\begin{align}
      &\, \eqref{e:26:reEXP2} + \eqref{e:30:reEXP2} \\
      =&\, \sum_{n=1}^{Q} \frac{1}{\lambda^n} \Bigg( \sum_{k=1}^K \sum_{ [\balpha]=1, \balpha \in \mathcal{I}_* }^{2n-2+Q} \sum_{l=1}^{\hat I^{(\balpha)}_{F,n}}
            \hat h^{(\balpha,l)}_{F,k,n} \hat \eta^{(\balpha,l)}_{F,k,n} (\mu t) \hat \chi^{(\balpha,l)}_{F,k,n} (\lambda \Phi) \mr {\bar D}^{\balpha} \bar\rho_{F,k}        \label{e:109:reEXP2} \\ 
      + &\, \sum_{k=1}^K \sum_{ [\balpha]=1, \balpha \in \mathcal{I}_* }^{2n-1+Q} \sum_{l=1}^{\check I^{(\balpha)}_{F,n}}  
            \check h^{(\balpha,l)}_{F,k,n} \check \eta^{(\balpha,l)}_{F,k,n} (\mu t) \mr {\bar D}^{\balpha} \bar\rho_{F,k} \Bigg)        \label{e:110:reEXP2} \\ 
      =&\, \sum_{k=1}^K \sum_{ [\balpha]=1, \balpha \in \mathcal{I}_* }^{2n-1+Q}
                  \mr {\bar D}^{\balpha} \bar\rho_{F,k} 
            \sum_{ n \geq \max\{ 1, \frac{[\balpha]+1-Q}{2}\} }^{Q} \frac{1}{\lambda^n} 
                  \Bigg( \sum_{l=1}^{\check I^{(\balpha)}_{F,n}} \check h^{(\balpha,l)}_{F,k,n} \check \eta^{(\balpha,l)}_{F,k,n} (\mu t)     \label{e:112:reEXP2} \\ 
      + &\, \sum_{l=1}^{\hat I^{(\balpha)}_{F,k,n}}
                  \hat h^{(\balpha,l)}_{F,n} \hat \eta^{(\balpha,l)}_{F,k,n} (\mu t) \hat \chi^{(\balpha,l)}_{F,k,n} (\lambda \Phi) \Bigg)    \label{e:114:reEXP2}
\end{align}
From (\ref{e:109:reEXP2}-\ref{e:110:reEXP2}) to (\ref{e:112:reEXP2}-\ref{e:114:reEXP2}), we set 
\begin{align}     \label{e:118:reEXP2}
      \hat h^{(\balpha,l)}_{F,k,n}, \hat \eta^{(\balpha,l)}_{F,k,n}, \hat \chi^{(\balpha,l)}_{F,k,n} = 0, \quad 
            \text{if } [\balpha] = 2n-1+Q.
\end{align} 
Then, we renumber the following correctors indexed by $n$ and $l$
\begin{align*}
      \bigg\{ 
      \bigg\{ \frac{1}{\lambda^n} \hat h^{(\balpha,l)}_{F,k,n}, \hat \eta^{(\balpha,l)}_{F,k,n}, \hat \chi^{(\balpha,l)}_{F,k,n} \bigg\}_{ 1 \leq l \leq \hat I^{(\balpha)}_{F,k,n} }, 
      \bigg\{ \frac{1}{\lambda^n} \check h^{(\balpha,l)}_{F,k,n}, \check \eta^{(\balpha,l)}_{F,k,n}, \check \chi^{(\balpha,l)}_{F,k,n} \bigg\}_{ 1 \leq l \leq \check I^{(\balpha)}_{F,k,n} }
      \bigg\}_{n}
\end{align*}
with $\check \chi^{(\balpha,l)}_{F,k,n} = 1$ in a single index, denoted again by $l$, and write them in the following form
\begin{align*}
      \big\{ \ddot{h}^{(\balpha,l)}_{F,k}, \ddot{\eta}^{(\balpha,l)}_{F,k}, \ddot{\chi}^{(\balpha,l)}_{F,k} \big\}_{1 \leq l \leq \mathfrak{I}}.
\end{align*}
Hence, the cardinality has the estimate
\begin{align}     \label{e:122:reEXP2}
      \mathfrak{I} \leq \max_{k,\balpha} \sum_{ n \geq \max\{ 1, \frac{[\balpha]+1-Q}{2}\} }^{Q}
            \Big( \hat I^{(\balpha)}_{F,k,n} + \check I^{(\balpha)}_{F,k,n} \Big)
\end{align}
In case where $\leq$ in \eqref{e:122:reEXP2} turns into a strict inequality, we set all correctors in additional slots to be zero. This leads to $\mathfrak{I} \lesssim 1$.

With the renumbering above, we deduce
\begin{align}     
      &\, \eqref{e:26:reEXP2} + \eqref{e:30:reEXP2}       \nonumber \\ 
      =&\, \sum_{k=0}^K \sum_{ [\balpha]=1, \balpha \in \mathcal{I}_* }^{2n-1+Q}
            \mr {\bar D}^{\balpha} \ddot{\rho}_{F,k}
            \sum_{l=1}^{\mathfrak{I}} 
            \ddot{h}^{(\balpha,l)}_{F,k} \ddot{\eta}^{(\balpha,l)}_{F,k} (\mu t) \ddot{\chi}^{(\balpha,l)}_{F,k} (\lambda \Phi)     \label{e:124:0:reEXP2} \\ 
      =&\, \sum_{k=0}^K \sum_{ [\balpha]=1, \balpha \in \mathcal{I}_* }^{Q}
            \mr {\bar D}^{\balpha} \ddot{\rho}_{F,k}
            \sum_{l=1}^{\mathfrak{I}} 
            \bP_Q \ddot{h}^{(\balpha,l)}_{F,k} \ddot{\eta}^{(\balpha,l)}_{F,k} (\mu t) \ddot{\chi}^{(\balpha,l)}_{F,k} (\lambda \Phi)     \label{e:124:2:reEXP2} \\ 
      +&\, \sum_{k=0}^K \sum_{ [\balpha]=1, \balpha \in \mathcal{I}_* }^{Q}
            \mr {\bar D}^{\balpha} \ddot{\rho}_{F,k}
            \sum_{l=1}^{\mathfrak{I}} 
            \Big( \ddot{h}^{(\balpha,l)}_{F,k} - \bP_Q \ddot{h}^{(\balpha,l)}_{F,k} \Big) \ddot{\eta}^{(\balpha,l)}_{F,k} (\mu t) \ddot{\chi}^{(\balpha,l)}_{F,k} (\lambda \Phi)    \label{e:124:4:reEXP2} \\ 
      +&\, \sum_{k=0}^K \sum_{ [\balpha]=Q+1, \balpha \in \mathcal{I}_* }^{2n-1+Q}
            \mr {\bar D}^{\balpha} \ddot{\rho}_{F,k}
            \sum_{l=1}^{\mathfrak{I}} 
            \ddot{h}^{(\balpha,l)}_{F,k} \ddot{\eta}^{(\balpha,l)}_{F,k} (\mu t) \ddot{\chi}^{(\balpha,l)}_{F,k} (\lambda \Phi)     \label{e:124:6:reEXP2}
\end{align}
From \eqref{e:38:reEXP2} and \eqref{e:40:reEXP2}, we have the following estimates
\begin{align}     \label{e:125:reEXP2}
      \sum_{k=0}^K \sum_{ [\balpha]=1, \balpha \in \mathcal{I}_* }^{2n-1+Q} \sum_{l=1}^{\mathfrak{I}} \vertiii{ \ddot{h}_{F,k}^{(\balpha,l)} } 
            \lesssim&\, \frac{\delta^{\sfrac12}\mr{\bar\lambda}}{\kappa\lambda^2} \cdot \lambda^{ - \gamma }. 
\end{align}
From \eqref{e:EXP2:hchi_est}, \eqref{e:EXP2:heta_est} and \eqref{e:EXP2:ceta_est}, we have
\begin{align}
      \| \nabla_\xi^p \ddot{\chi}^{(\balpha,l)}_{F,n} \|_\infty 
            \lesssim&\, 1,    \quad p \leq P,     \label{e:126:reEXP2} \\ 
      \| \partial_\tau^p \ddot{\eta}^{(\balpha,l)}_{F,n} \|_\infty 
            \lesssim&\, 1,    \quad p \leq P-Q.     \label{e:128:reEXP2}
\end{align}

Now we return to (\ref{e:24:reEXP2}-\ref{e:32:reEXP2}). The terms in \eqref{e:24:reEXP2} and \eqref{e:28:reEXP2} are represented in (\ref{e:98:0:reEXP2}-\ref{e:98:4:reEXP2}). The terms in \eqref{e:26:reEXP2} and \eqref{e:30:reEXP2} are represented in (\ref{e:124:0:reEXP2}-\ref{e:124:6:reEXP2}). We rename the terms in \eqref{e:98:0:reEXP2} and \eqref{e:124:2:reEXP2} such that
\begin{align}     \label{e:130:reEXP2}
      (\ref{e:98:0:reEXP2}) + (\ref{e:124:2:reEXP2})
      = \sum_{a=0}^{K+Q} \sum_{ [\balpha]=1, \balpha \in \mathcal{I}_* }^{Q} 
            \mr {\bar D}^{\balpha} \ddot{\rho}_{a}  
            \sum_{l=1}^{\mathfrak{I}} \ddot{h}_{a}^{(\balpha,l)} \ddot{\eta}_{a}^{(\balpha,l)} (\mu t) \ddot{\chi}_{a}^{(\balpha,l)} (\lambda\Phi),
\end{align}
which gives \eqref{e:E1:reEXP2} with the remainder term given by
\begin{align}     \label{e:132:reEXP2}
      \tilde \rho_F := (\ref{e:98:2:reEXP2}) + (\ref{e:98:4:reEXP2}) + (\ref{e:124:4:reEXP2}) + (\ref{e:124:6:reEXP2}) + \tilde \rho.
\end{align}
This proves \eqref{e:E1:reEXP2}. The estimates \eqref{e:hEst:reEXP2}, \eqref{e:chiEst:reEXP2} and \eqref{e:etaEst:reEXP2} follow from \eqref{e:125:reEXP2}, \eqref{e:126:reEXP2} and \eqref{e:128:reEXP2}. The estimate \eqref{e:4:reEXP2} and transmission condition follow from \eqref{e:96:reEXP2}, Lemma \ref{c:t=0Ref:EXP2} and Lemma \ref{l:t=0:EXP2}.

\subsubsection{Proof of Lemma \ref{l:reEXP2}: Estimating macroscopic states}    \label{sss:p2:reEXP2}

Define
\begin{align}
      \ddot{h}_{F,k,0}^{(\balpha)} :=&\, \frac{\mr{\bar\lambda}} {\mr\lambda}
            \sum_{l = 1}^{J} h_{F,k}^{(\balpha,l)} \langle \eta_{F,k}^{(\balpha,l)} \rangle_{\tau} \langle \chi_{F,k}^{(\balpha,l)} \rangle_{\xi},    \label{e:56:reEXP2} \\ 
      \ddot{h}_{F,k,n}^{(\balpha)} :=&\, 
      \frac{1}{\lambda^n} \bar h_{F,k,n}^{(\balpha)},       \qquad n \geq 1.   \label{e:58:reEXP2}
\end{align}
The macroscopic states $\{\ddot{\rho}_n\}_n$ satisfy the following equations
\begin{align}
      \bar L_0 \ddot{\rho}_{0}
      =&\, \sum_{k=0}^K \sum_{ [\balpha] = 1, \balpha \in \mathcal{I}_* }^Q
            \mr{\bar \divr} \Big( \ddot{h}_{F,k,0}^{(\balpha)} \mr {\bar D}^{\balpha} \bar \rho_{F,k} \Big),      \label{e:60:reEXP2} \\ 
      \bar L_0 \ddot{\rho}_{n}
      =&\, \sum_{a=0}^{n-1} \sum_{ [\balpha]=1, \balpha \in \mathcal{I}_* }^{2(n-a)+1} 
            \mr{\bar \divr} \Big( \ddot h^{(\balpha)}_{n-a} \mr {\bar D}^{\balpha} \ddot{\rho}_{a} \Big) 
      + \sum_{k=1}^K \sum_{ [\balpha]=1, \balpha \in \mathcal{I}_* }^{2n+Q} 
            \mr{\bar \divr} \Big( \ddot h_{F,k,n}^{(\balpha)} \mr {\bar D}^{\balpha} \bar \rho_{F,k} \Big)     \label{e:62:reEXP2}
\end{align}
with $\ddot h_{F,k,n}^{(\balpha)} \in \bar {\mathcal{P}} ( 2Q + NQ )$ and
\begin{align}
       \sum_{k,\balpha} \vertiii{ \ddot{h}_{F,k,0}^{(\balpha)} }
            \lesssim&\, \bar\kappa \mr{\bar\lambda}^2 \lambda^{-2b\gamma},      
                  \label{e:66:reEXP2} \\
      \sum_{k,\balpha} \vertiii{ \ddot h_{F,k,n}^{(\balpha)} } 
            \lesssim&\, \bar\kappa \mr{\bar\lambda}^2 \lambda^{-2nb\gamma},
                  \quad n \geq 1.     \label{e:68:reEXP2}
\end{align}
Here, \eqref{e:66:reEXP2} follows from \eqref{e:4:homSmallscl}, \eqref{e:56:reEXP2} and $ \bar \kappa > \frac{\delta^{\sfrac12}} {\lambda} $. \eqref{e:68:reEXP2} follows from \eqref{e:EXP2:bhF_Pest}, \eqref{e:Up:reEXP2} and \eqref{e:58:reEXP2}.

Note that from \eqref{e:EXP2:brho0_ini} and \eqref{e:96:reEXP2}, $\ddot\rho_0(x,0) = 0$. We apply Lemma \ref{l:P2resolveToF} with \eqref{e:66:reEXP2} to deduce
\begin{align}
      \resF{ \ddot{\rho}_0 }_{\bar u, \bar \kappa, n_0+1} 
            \lesssim&\, \bar \lambda^{ - \alpha } \lambda^{ - (s-1)b\gamma }
                  \lambda^{ - (b-1)\gamma } \Upsilon_{\ini},      \label{e:70:reEXP2} \\ 
      \resF{ \ddot{\rho}_0 }_{\bar u, \bar \kappa, n_0+1}^{(H,\alpha_0)}
            \lesssim&\, \lambda^{ - (s-1)b\gamma }
                  \lambda^{ - (b-1)\gamma } \Upsilon_{\ini}.      \label{e:70_2:reEXP2}
\end{align}
Here, note that \eqref{e:66:reEXP2} involves $\mr{\bar\lambda}$ and $\lambda$, while the estimate in Lemma \ref{l:P2resolveToF} involves $\mr{\bar\lambda}$ and $\bar\lambda$. This gives the rise of the factor $\lambda^{ - (s-1)b\gamma }$. It will help us to control certain constants.
From \eqref{e:1:reEXP1}, \eqref{e:2:reEXP1} and \eqref{e:4:reEXP2}, we can deduce
\begin{align}
      \resN{ \ddot{\rho}_1(\cdot,0) }_{\bar u, \bar \kappa} 
            \lesssim&\, \bar \lambda^{ - \alpha } \lambda^{ - sb\gamma } \Upsilon_{\ini},       \label{e:71:reEXP2} \\ 
      \resN{ \ddot{\rho}_1(\cdot,0) }_{\bar u, \bar \kappa}^{(H,\alpha_0)} 
            \lesssim&\, \lambda^{ - sb\gamma } \Upsilon_{\ini}.       \label{e:71_2:reEXP2}
\end{align}
Applying Lemma \ref{l:P2resolveToF} again with \eqref{e:68:reEXP2} and \eqref{e:71:reEXP2}, we have
\begin{align}
      \resF{ \ddot{\rho}_1 }_{\bar u, \bar \kappa, n_0+2} 
            \lesssim&\, \bar \lambda^{ - \alpha } \lambda^{ - (s-1)b\gamma }
                  \lambda^{ - (b-1)\gamma } \Upsilon_{\ini},     \label{e:72:reEXP2} \\
      \resF{ \ddot{\rho}_1 }_{\bar u, \bar \kappa, n_0+2}^{(H,\alpha_0)}
            \lesssim&\, \lambda^{ - (s-1)b\gamma }
                  \lambda^{ - (b-1)\gamma } \Upsilon_{\ini}.     \label{e:72_2:reEXP2}
\end{align}
Here, note that $s-1 \geq 0$ and $\ddot{h}^{(\balpha)}_{n-a}$ is given in Lemma \ref{l:reEXP1}, according to Remark \ref{r:shareCorrector}.

Similar to the proof of Lemma \ref{l:reEXP1} in Section \ref{sss:p2:reEXP1}, inductively we can show
\begin{align}
      \resF{ \ddot{\rho}_n }_{\bar u, \bar \kappa, n_0+n+1} 
            \lesssim&\, \bar \lambda^{ - \alpha } \lambda^{ - (s+n-2)b\gamma } \lambda^{ - (b-1)\gamma } \Upsilon_{\ini},
            \quad \forall \, n \geq 2.    \label{e:74:reEXP2} \\ 
      \resF{ \ddot{\rho}_n }_{\bar u, \bar \kappa, n_0+n+1}^{(H,\alpha_0)} 
            \lesssim&\, \lambda^{ - (s+n-2)b\gamma } \lambda^{ - (b-1)\gamma } \Upsilon_{\ini},
            \quad \forall \, n \geq 2.    \label{e:76:reEXP2}
\end{align}
Then \eqref{e:70:reEXP2}, \eqref{e:72:reEXP2} and \eqref{e:74:reEXP2} give \eqref{e:6:reEXP2}.

\subsubsection{Proof of Lemma \ref{l:reEXP2}: Estimating the remainder}    \label{sss:p3:reEXP2}

This is analogous to Section \ref{sss:p3:reEXP1}. To estimate $\tilde \rho$, we appeal to the equation satisfied by $\tilde \rho$ with initial datum given in \eqref{e:EXP2:trho_ini}. Hence, the initial datum is estimated by Remark \ref{r:datumProj} and Lemma \ref{l:t=0:EXP2}, and the force $\tilde f$ is estimated by Lemma \ref{l:remainderEst_IndtS}. With such estimates, we use the energy estimates in Lemma \ref{l:L2_init_energy} and Lemma \ref{l:L2_F_energy} to estimate $\tilde \rho$. We remark here we use \eqref{e:74:reEXP2} to verify the assumptions in Lemma \ref{l:remainderEst_IndtS}. The estimates on \eqref{e:98:2:reEXP2}, \eqref{e:98:4:reEXP2}, \eqref{e:124:4:reEXP2} and \eqref{e:124:6:reEXP2} are similar to Section \ref{sss:p3:reEXP1}, so we omit the details.

Above together with \eqref{e:132:reEXP2} gives \eqref{e:8:reEXP2}.

\subsection{Inductive expansion at step $q+1$: Proof of Theorem \ref{t:homRslv}}    \label{ss2:estQ+1}

In this section we prove Theorem \ref{t:homRslv}. The key element is the following proposition, which shows how to propagate the expansion structure from step $n$ to step $n+1$. Then, Theorem \ref{t:homRslv} follow from the expansions and direct computations.

\begin{proposition}     \label{p:ind_resolvability}
Fix some $\rho_{\ini} \in C^\infty (\T^2)$ and $\Upsilon_{\ini} > 0$. Given \eqref{e:1:reEXP1}-\eqref{e:1:reEXP1} and Assumption \ref{a:OEq}, consider the equations \eqref{e:7:OEq_Idct1} and \eqref{e:8:OEq_Idct}. Fix any $n \leq Q-1$. Suppose $\varrho_n$ admits the following expansion
\begin{align}     \label{e:16:estQ+1}
      \varrho_{n} =&\, \sum_{a=0}^{(Q+1)(n+1)} \sum_{[\balpha]=1, \balpha \in \mathcal{I}_*}^{Q} \mr {\bar D}^{\balpha} \bar \rho_{n,a} 
      \sum_{l=1}^{I_{n}} h_{n,a}^{(\balpha,l)} \eta_{n,a}^{(\balpha,l)} (\mu t) \chi_{n,a}^{(\balpha,l)} (\lambda\Phi) + \sum_{a=0}^{Q} \bar\rho_{n,a} + \tilde \rho_{n},
\end{align}
with $I_{n} \lesssim 1$,
\begin{align}
      h_{n,a}^{(\balpha,l)} \in \bar{\mathcal{P}} ( Q ),    \quad 
      \sum_{a,\balpha,l} \vertiii{ h_{n,a}^{(\balpha,l)} }
            \lesssim&\, \frac{\delta^{\sfrac12}\mr{\bar\lambda}} {\kappa\lambda^2},      
            \label{e:18:estQ+1}  
\end{align}
and
\begin{align}
      \| \nabla_\xi^p \chi_{n,a}^{(\balpha,l)} \|_{\infty} \lesssim&\, 1
            \quad p \leq 8Q^3 - 6Q (n+1),     \label{e:24:estQ+1} \\ 
      \| \partial_\tau^p \eta_{n,a}^{(\balpha,l)} \|_{\infty} \lesssim&\, 1     
            \quad p \leq 8Q^3 - 6Q (n+1).     \label{e:26:estQ+1} 
\end{align}
For any $a,\balpha,l$, $\big( \chi_{n,a}^{(\balpha,l)}, \eta_{n,a}^{(\balpha,l)} \big)$ satisfies the generalized shear condition. The collection $\{\bar\rho_{n,a}\}_a$ satisfy the transmission condition. Suppose that for any $0 \leq a \leq (Q+1)(n+1)$,
\begin{align}
      \| \bar\rho_{n,a} (\cdot,0) \|_2 + \sum_{j=q_*}^{q} \lambda_j^{-\sfrac{\gamma}{2}} | \bar\rho_{n,a} (\cdot,0) |_{\fS(j)} 
            \lesssim&\, \lambda^{ - \big( \sum_{i=0}^{n-1} (s_i-1) \big) b \gamma } \Upsilon_{\ini},     \label{e:27:estQ+1} \\ 
      \resF{ \bar\rho_{n,a} }_{ \bar u, \bar \kappa, (Q+1)(n+1) } 
            \lesssim&\, \bar \lambda^{ - \alpha } \lambda^{ - \big( \sum_{i=0}^{n-1} (s_i-1) \big) b \gamma } 
                  \Upsilon_{\ini},     \label{e:28:estQ+1} \\ 
      \resF{ \bar\rho_{n,a} }_{ \bar u, \bar \kappa, (Q+1)(n+1) } ^{ (H,\alpha_0) }
            \lesssim&\, \lambda^{ - \big( \sum_{i=0}^{n-1} (s_i-1) \big) b \gamma } 
                  \Upsilon_{\ini}      \label{e:29:estQ+1}
\end{align}
and
\begin{equation}     \label{e:30:estQ+1} 
\begin{split}
      \| \mr O^{\bbeta} \tilde \rho_{n} \| + \kappa \mr\lambda \| \mr O^{\bbeta} \mr\nabla \tilde \rho_{n} \|_2 
            \lesssim&\, \lambda^{ - 3 - \big( [\bbeta] + \sum_{i=0}^{n-1} (s_i-1) \big) b \gamma } \Upsilon_{\ini}  \quad 
            \text{for } [\bbeta] \leq 2Q + 3Q(Q-n).
\end{split}
\end{equation}
Then $\varrho_{n+1}$ admits the following expansion
\begin{align}     \label{e:32:estQ+1}
      \varrho_{n+1} = \sum_{a=0}^{(Q+1)(n+2)} \sum_{[\balpha]=1, \balpha \in \mathcal{I}_*}^{Q} \mr {\bar D}^{\balpha} \bar \rho_{n+1,a} 
      \sum_{l=1}^{I_{n+1}} h_{n+1,a}^{(\balpha,l)} \eta_{n+1,a}^{(\balpha,l)} (\mu t) \chi_{n+1,a}^{(\balpha,l)} (\lambda\Phi) + \sum_{a=0}^{Q} \bar\rho_{n+1,a} + \tilde \rho_{n+1}
\end{align}
with $I_{n+1} \lesssim 1$,
\begin{align}
      h_{n+1,a}^{(\balpha,l)} \in \bar{\mathcal{P}} ( Q ),  \quad 
      \sum_{a,\balpha,l} \vertiii{ h_{n+1,a}^{(\balpha,l)} } 
            \lesssim&\, \frac{\delta^{\sfrac12}\mr{\bar\lambda}} {\kappa\lambda^2}
                  \label{e:36:estQ+1}  
\end{align}
and
\begin{align}
      \| \nabla_\xi^p \chi_{n+1,a}^{(\balpha,l)} \|_{\infty}
            \lesssim&\, 1,  \quad p \leq 8Q^3 - 6Q (n+2),     \label{e:42:estQ+1} \\ 
      \| \partial_\tau^p \eta_{n+1,a}^{(\balpha,l)} \|_{\infty}
            \lesssim&\, 1,  \quad p \leq 8Q^3 - 6Q (n+2).     \label{e:44:estQ+1} 
\end{align}
For any $a,\balpha,l$, $\big( \chi_{n+1,a}^{(\balpha,l)}, \eta_{n+1,a}^{(\balpha,l)} \big)$ satisfies the generalized shear condition. Moreover, for any $0 \leq a \leq (Q+1)(n+2)$,
\begin{align}
      \| \bar\rho_{n+1,a} (\cdot,0) \|_2 + \sum_{j=q_*}^{q} \lambda_j^{-\sfrac{\gamma}{2}} | \bar\rho_{n+1,a} (\cdot,0) |_{\fS(j)}
            \lesssim&\, \lambda^{ - \big( \sum_{i=0}^{n} (s_i-1) \big) b \gamma } \Upsilon_{\ini},     \label{e:47:estQ+1} \\ 
      \resF{ \bar\rho_{n+1,a} }_{\bar u, \bar \kappa, (Q+1)(n+2)} 
            \lesssim&\, \bar \lambda^{ - \alpha } \lambda^{ - \big( \sum_{i=0}^{n} (s_i-1) \big) b \gamma } 
                  \lambda^{ - (b-1)\gamma } \Upsilon_{\ini}    \label{e:48:estQ+1} \\ 
      \resF{ \bar\rho_{n+1,a} }_{\bar u, \bar \kappa, (Q+1)(n+2)} ^{ (H,\alpha_0) }
            \lesssim&\, \lambda^{ - \big( \sum_{i=0}^{n} (s_i-1) \big) b \gamma } 
                  \lambda^{ - (b-1)\gamma } \Upsilon_{\ini}      \label{e:49:estQ+1}
\end{align}
and
\begin{align}     \label{e:50:estQ+1} 
      \| \mr O^{\bbeta} \tilde \rho_{n+1} \| + \kappa \mr\lambda \| \mr O^{\bbeta} \mr\nabla \tilde \rho_{n+1} \|_2 
            \lesssim&\, \lambda^{ - 3 - \big( [\bbeta] + \sum_{i=0}^{n} (s_i-1) \big) b \gamma } \Upsilon_{\ini}  \quad 
            \text{for } [\bbeta] \leq 2Q + 3Q(Q-n-1).
\end{align}
\end{proposition}

\begin{remark}
Comparing \eqref{e:28:estQ+1}-\eqref{e:29:estQ+1} with \eqref{e:48:estQ+1}-\eqref{e:49:estQ+1}, we see that the factor $\lambda^{ - (b-1)\gamma }$ appears in the estimates for $\bar\rho_{n+1,a}$, which is not present in the estimates for $\bar\rho_{n,a}$. This is an actual gain in each step of the expansion, which helps to control the constants in the estimates.
\end{remark}

We also need the following corollaries.

\begin{corollary}       \label{c:reEXP1Est}
For $\varrho_0$ in Lemma \ref{l:reEXP1}, we have that, for any $\bbeta \in \mathcal{I}_*$ with $[\bbeta] \leq Q$,
\begin{align}
      \| \mr O^{\bbeta} \varrho_0 \|_\infty \leq&\, \frac{1}{2} \lambda^{ - \alpha - [\bbeta] b \gamma } 
            \bigg( | \rho_{\ini} |_{\mathfrak{D}(v, \kappa)}
                  + \sum_{j=q_*}^{q+1} \lambda_j^{-\sfrac{\gamma}{2}} | \rho_{\ini} |_{\fS(j)} \bigg),      \label{e:2:reEXP1Est} \\ 
      \| \varrho_0 \|_{L^\infty_t C^{\alpha_0}_x} \leq&\, \bigg( 1 - \frac{4}{q} + \lambda^{-\sfrac{\gamma}{2}} \bigg) 
            \bigg( | \rho_{\ini} |_{\mathfrak{D}(v, \kappa)}
                  + \sum_{j=q_*}^{q+1} \lambda_j^{-\sfrac{\gamma}{2}} | \rho_{\ini} |_{\fS(j)} \bigg).   \label{e:4:reEXP1Est}
\end{align}
\end{corollary}

\begin{corollary}       \label{c:ind_resolve}
For $n \leq Q-1$ and the expansion \eqref{e:32:estQ+1} with \eqref{e:36:estQ+1}-\eqref{e:50:estQ+1} in Proposition \ref{p:ind_resolvability}, we have that, for any $\bbeta \in \mathcal{I}_*$ with $[\bbeta] \leq Q$,
\begin{align}
      \| \mr O^{\bbeta} \varrho_{n+1} \|_\infty \lesssim&\, \lambda^{-b\gamma} \lambda^{ - \alpha - \big( [\bbeta] + \sum_{i=0}^{n} (s_i-1) \big) b \gamma } \Upsilon_{\ini},      \label{e:3:ind_resolve} \\ 
      \| \varrho_{n+1} \|_{L^\infty_t C^{\alpha_0}_x} \lesssim&\, \lambda^{-(b-1)\gamma} \lambda^{ - \big( \sum_{i=0}^{n} (s_i-1) \big) b \gamma } \Upsilon_{\ini}.     \label{e:4:ind_resolve}
\end{align}
\end{corollary}

From above results, we can prove Theorem \ref{t:homRslv}. We postpone the proof of Proposition \ref{p:ind_resolvability}, Corollary \ref{c:reEXP1Est} and Corollary \ref{c:ind_resolve} to the end of this section.

\begin{proof}[Proof of Theorem \ref{t:homRslv}]
Given Proposition \ref{p:chainStab}, it suffices to prove that, for any $ \rho_{\ini} \in C^\infty (\T^2)$, we have
\begin{align*}
      \resN{ \rho_{\ini} }_{v,\kappa} 
            \leq&\, \frac{1}{2} \lambda^{-\alpha} \bigg( | \rho_{\ini} |_{\mathfrak{D}(v, \kappa)}
                  + \sum_{j=q_*}^{q+1} \lambda_j^{-\sfrac{\gamma}{2}} | \rho_{\ini} |_{\fS(j)} \bigg),    \\ 
      \resN{ \rho_{\ini} }^{(H,\alpha_0)}_{v,\kappa} 
            \leq&\, \bigg( 1 - \frac{4}{q} + \lambda^{-\frac{(b-1)\gamma}{2}} \bigg) \bigg( | \rho_{\ini} |_{\mathfrak{D}(v, \kappa)}
                  + \sum_{j=q_*}^{q+1} \lambda_j^{-\sfrac{\gamma}{2}} | \rho_{\ini} |_{\fS(j)} \bigg).
\end{align*}
In this proof, we only prove the estimate for $\resN{ \cdot }_{v,\kappa}$ in the setting of Assumption \ref{a:OEq}. The H\"older estimate for $\resN{ \cdot }^{(H,\alpha_0)}_{v,\kappa}$ is similar.

It suffices to prove that for any $0 \leq n \leq 2Q^2$ and any $[\bbeta] \leq 11Q^3 - 3Qn$, we have 
\begin{align}
      \| \mr O^{\bbeta} \varrho_{n} \|_\infty \leq&\, \frac{1}{2} \lambda^{ - \alpha - \big( [\bbeta] + \sum_{i=0}^{n-1} (s_i-1) \big) b \gamma } \bigg( | \rho_{\ini} |_{\mathfrak{D}(v, \kappa)}
                  + \sum_{j=q_*}^{q+1} \lambda_j^{-\sfrac{\gamma}{2}} | \rho_{\ini} |_{\fS(j)} \bigg).   \label{e:18:homRslv}
\end{align}

When $n=0$ and $[\bbeta] \leq Q$, (\ref{e:18:homRslv}) follows from Lemma \ref{l:reEXP1} and Corollary \ref{c:reEXP1Est}. 

When $n \leq Q$ and $[\bbeta] \leq Q$, notice that Lemma \ref{l:reEXP1} and Proposition \ref{p:ind_resolvability} form an inductive proof, which verifies that the solution $\varrho_n$ to (\ref{e:5:OEq_Idct}-\ref{e:8:OEq_Idct}) admits the expansion \eqref{e:16:estQ+1}. Here, we set 
\begin{align*}
      \Upsilon_{\ini} = | \rho_{\ini} |_{\mathfrak{D}(v, \kappa)}
                  + \sum_{j=q_*}^{q+1} \lambda_j^{-\sfrac{\gamma}{2}} | \rho_{\ini} |_{\fS(j)}.
\end{align*} Then (\ref{e:18:homRslv}) with $[\bbeta] \leq Q$ follows from Corollary \ref{c:ind_resolve}.

When $n \leq Q$ or $[\bbeta] \geq Q$, (\ref{e:18:homRslv}) follows from Proposition \ref{p:unconEnergyP} and Remark \ref{r:unconEnergyP}. Here, the constants in $\lesssim$ of \eqref{e:4:unconEnergyP} are absorbed by the additional negative power of $\lambda$. This finishes the proof of \eqref{e:18:homRslv} for $0 \leq n \leq 2Q^2$.
\end{proof}

The proofs of Corollary \ref{c:ind_resolve} and Corollary \ref{c:reEXP1Est} given below are analogous.

\begin{proof}[Proof of Corollary \ref{c:ind_resolve}]
To prove \eqref{e:3:ind_resolve} and \eqref{e:4:ind_resolve}, we use the expansion \eqref{e:32:estQ+1}. We need to estimate the three terms in \eqref{e:32:estQ+1}. Moreover, here we need to control the constants in $\lesssim$ carefully.

\begin{step}[The term $\tilde\rho_{n+1}$]       \label{c:s:2:ind_resolve}
The $L^\infty_t L^2_x$-estimate on $\nabla^p \tilde\rho_{n+1}$ is obvious from \eqref{e:50:estQ+1}, and we have gained a factor $\lambda^{-3}$. The $L^\infty_{x,t}$-estimate and $L^\infty_t C^{\alpha_0}_x$-estimate on $\tilde\rho_{n+1}$ follow from $L^\infty_t L^2_x$-estimate on $\nabla^p \tilde\rho_{n+1}$ and interpolation inequality in the spirit of Proposition \ref{p:unconEnergyHom}. The substantial gain $\lambda^{-3}$ compensates the derivative loss in interpolation inequality and also absorbs constants in $\lesssim$.
\end{step}

\begin{step}[The terms $\bar\rho_{n+1,a}$]      \label{c:s:4:ind_resolve}

For terms like $\bar\rho_{n+1,a}$, we use \eqref{e:48:estQ+1} and Lemma \ref{l:compDecompedF} to deduce that for any $\bbeta \in \mathcal{I}_*$,
\begin{align}     \label{e:12:ind_resolve}
      \| \mr O^{\bbeta} \bar\rho_{n+1,a} \|_\infty \lesssim \bar \lambda^{ - \alpha } \lambda^{ - \big( \sum_{i=0}^{n} (s_i-1) \big) b \gamma } \cdot
            \bigg( \frac{\mr{\bar\lambda}} {\mr\lambda} \bigg)^{|\bbeta|_x} 
            \bigg( \frac{\delta^{\sfrac12}{\lambda}} {\mr\mu} \bigg)^{|\bbeta|_t} \Upsilon_{\ini}.
\end{align}
From the parameter setting, in particular \eqref{e:alphaDef} in Section \ref{ss:parameters}, we have
\begin{align}     \label{e:14:ind_resolve}
      \bar \lambda^{ - \alpha } \cdot \frac{\mr{\bar\lambda}} {\mr\lambda} 
            \leq \lambda^{ - \alpha - b\gamma }.
\end{align}
\eqref{e:12:ind_resolve} and \eqref{e:14:ind_resolve} together yield the desired $L^\infty_{x,t}$-estimate for $\bar\rho_{n+1,a}$.

The $L^\infty_t C^{\alpha_0}_x$-estimate of $\bar\rho_{n+1,a}$ follows from \eqref{e:49:estQ+1}. The factor $\lambda^{-(b-1)\gamma}$ in \eqref{e:49:estQ+1} absorbs the constant in $\lesssim$.

\end{step}

\begin{step}[The main oscillatory terms]    \label{c:s:6:ind_resolve}

For terms like 
\begin{align}     \label{e:6:ind_resolve}
      \mr {\bar D}^{\balpha} \bar \rho_{n+1,a} h_{n+1,a}^{(\balpha,l)} \eta_{n+1,a}^{(\balpha,l)} (\mu t) \chi_{n+1,a}^{(\balpha,l)} (\lambda\Phi),
\end{align}
we apply Lemma \ref{l:compDecompedF} to deduce, for any $\bbeta \in \mathcal{I}_*$,
\begin{align}     \label{e:8:ind_resolve}
      \| \mr O^{\bbeta} (\ref{e:6:ind_resolve}) \|_\infty
            \lesssim \bar \lambda^{ - \alpha - [\balpha] \gamma } \lambda^{ - \big( \sum_{i=0}^{n} (s_i-1) \big) b \gamma } \cdot \frac{\delta^{\sfrac12}\mr{\bar\lambda}} {\kappa\lambda^2} 
            \bigg( \frac{\lambda}{\mr\lambda} \bigg)^{|\bbeta|_x} 
            \bigg( \frac{\delta^{\sfrac12}{\lambda}} {\mr\mu} \bigg)^{|\bbeta|_t} \Upsilon_{\ini}.
\end{align}
From the parameter setting including \eqref{e:alphaDef} in Section \ref{ss:parameters}, we have
\begin{align}     \label{e:10:ind_resolve}
      \bar \lambda^{ - \alpha - b\gamma } \cdot \frac{\delta^{\sfrac12}\mr{\bar\lambda}} {\kappa\lambda^2} \cdot \frac{\lambda}{\mr\lambda} 
            \leq \lambda^{ - \alpha - b\gamma } \cdot \bar\lambda^{-(b-1)^2b\gamma}.
\end{align}
Now we use $\bar\lambda^{-(b-1)^2b\gamma}$ to absorb the constant in $\lesssim$ of \eqref{e:8:ind_resolve}. This finishes the $L^\infty_{x,t}$-estimate for the terms of form \eqref{e:6:ind_resolve}.

To derive the $L^\infty_t C^{\alpha_0}_x$-estimate of $(\ref{e:6:ind_resolve})$, we have that, similar to \eqref{e:8:ind_resolve} and \eqref{e:10:ind_resolve}
\begin{equation}  \label{e:16:ind_resolve}
\begin{split}
      \| (\ref{e:6:ind_resolve}) \|_{\infty} \lesssim&\, \lambda^{ - \alpha - b\gamma } \lambda^{ - \big( \sum_{i=0}^{n} (s_i-1) \big) b \gamma } \Upsilon_{\ini},   \\ 
      \| \nabla (\ref{e:6:ind_resolve}) \|_{\infty} \lesssim&\, \mr \lambda \lambda^{ - \alpha - b\gamma } \lambda^{ - \big( \sum_{i=0}^{n} (s_i-1) \big) b \gamma } \Upsilon_{\ini}.
\end{split}
\end{equation}
Then we interpolate between $\| (\ref{e:6:ind_resolve}) \|_{\infty}$ and $\| \nabla (\ref{e:6:ind_resolve}) \|_{\infty}$ to get the $L^\infty_t C^{\alpha_0}_x$-estimate, i.e.
\begin{equation}
\begin{split}
      \| (\ref{e:6:ind_resolve}) \|_{L^\infty_t C^{\alpha_0}_x} 
            \leq&\, \| \nabla (\ref{e:6:ind_resolve}) \|_{\infty}^{\alpha_0} \| (\ref{e:6:ind_resolve}) \|_{\infty}^{1-\alpha_0}    \\ 
            \lesssim&\, \lambda^{ - \alpha - b\gamma } \lambda^{ - \big( \sum_{i=0}^{n} (s_i-1) \big) b \gamma } \mr \lambda^{\alpha_0} \Upsilon_{\ini},
\end{split}
\end{equation}
which indeed follows from \eqref{e:16:ind_resolve} and the parameter relations in Section \ref{ss:parameters}, in particular \eqref{e:2:bConsParameter}, \eqref{e:0:mScaleRela} and \eqref{e:smallGammaR}.
\end{step}

\end{proof}

\begin{proof}[Proof of Corollary \ref{c:reEXP1Est}]

The proof of \eqref{e:2:reEXP1Est} is analogous to \eqref{e:3:ind_resolve} in Corollary \ref{c:ind_resolve}. For the proof of \eqref{e:4:reEXP1Est}, the only difference is the control of constants in $\lesssim$. In the expansion formula \eqref{e:exp:reEXP1} of $\varrho_0$ in Lemma \ref{l:reEXP1}, we deduce from \eqref{e:const:reEXP1} and \eqref{e:holder:reEXP1} that
\begin{align}
      \| \bar \rho_{0,0} \|_{L^\infty_t C^{\alpha_0}_x} 
            \leq&\, \bigg( 1 - \frac{4}{q} \bigg) \bigg( | \rho_{\ini} |_{\mathfrak{D}(v, \kappa)}
                  + \sum_{j=q_*}^{q+1} \lambda_j^{-\sfrac{\gamma}{2}} | \rho_{\ini} |_{\fS(j)} \bigg),      \label{e:12:reEXP1Est} \\
      \| \bar \rho_{0,a} \|_{L^\infty_t C^{\alpha_0}_x}
            \lesssim&\, \lambda^{ - ab\gamma } \bigg( | \rho_{\ini} |_{\mathfrak{D}(v, \kappa)}
                  + \sum_{j=q_*}^{q+1} \lambda_j^{-\sfrac{\gamma}{2}} | \rho_{\ini} |_{\fS(j)} \bigg).      \label{e:14:reEXP1Est}
\end{align}
The $\| \cdot \|_{L^\infty_t C^{\alpha_0}_x}$ of other terms in \eqref{e:exp:reEXP1} are estimated similarly as in Step \ref{c:s:6:ind_resolve} and Step \ref{c:s:2:ind_resolve}.
\end{proof}

\begin{proof}[Proof of Proposition \ref{p:ind_resolvability}]

We apply Lemma \ref{l:decompRho} to $\mr O^{\bomega_n} \varrho_n$. Note that $[\bomega_n] \leq 3Q$. Then we get
\begin{align}
      \mr O^{\bomega_n} \varrho_{n} =&\, \sum_{a=0}^{(Q+1)(n+1)} \sum_{|\balpha|=1, \balpha \in \mathcal{I}_*}^{Q} \sum_{l=1}^{I_n}
      \sum_{[\bbeta] \leq [\bomega_n]}
            \mr {\bar D}^{\bbeta \balpha} \bar \rho_{n,a} 
      \sum_{k=1}^{ (3N^2)^{[\bomega_n]} }
            h_{n,a,k}^{(\balpha,\bbeta,l)} \eta_{n,a,k}^{(\balpha,\bbeta,l)} (\mu t) \chi_{n,a,k}^{(\balpha,\bbeta,l)} (\lambda\Phi)      \label{e:60:estQ+1} \\
      +&\, \mr D^{\bomega_n} \tilde \rho_{n}    \label{e:62:estQ+1}
\end{align}
with $|\bbeta\balpha| \leq 4Q$, $h_{n,a,k}^{(\balpha,\bbeta,l)} \in \bar{\mathcal{P}}(7NQ)$ and
\begin{align}
      \sum_{\bbeta,k} \vertiii{ h_{n,a,k}^{(\balpha,\bbeta,l)} }
            \lesssim&\, \vertiii{ h_{n,a}^{(\balpha,l)} } \bigg( \frac{\lambda}{\mr\lambda} \bigg)^{|\bomega_n|_{x}} + \frac{\mr{\bar\lambda}} {\mr\lambda}.      \label{e:66:estQ+1}
\end{align}
The corrector components satisfy
\begin{align}
      \| \nabla_\xi^p \chi_{n,a,k}^{(\balpha,\bbeta,l)} \|_{\infty}
            \lesssim&\, 1,
            \quad p \leq 8Q^3 - 6Q (n+1) - 3Q,       \label{e:68:estQ+1} \\ 
      \| \partial_\tau^p \eta_{n,a,k}^{(\balpha,\bbeta,l)} \|_{\infty}
            \lesssim&\, 1,
            \quad p \leq 8Q^3 - 6Q (n+1) - 3Q.       \label{e:70:estQ+1}
\end{align}
Here, for any $a,k,\balpha,\bbeta,l$, $\big( \chi_{n,a,k}^{(\balpha,\bbeta,l)}, \eta_{n,a,k}^{(\balpha,\bbeta,l)} \big)$ satisfies the generalized shear condition.

Then we apply Corollary \ref{c:decomp} to $g_n \in \mathcal{\bar P}(2NQ)$, to deduce
\begin{align}     \label{e:74:estQ+1} 
      g_n = g_{n,0} + \sum_{j = 1}^{(10N)^{8NQ}} \ddot{\eta}_{n,j} (\mu t) \ddot{\chi}_{n,j} (\lambda \Phi) g_{n,j}
\end{align}
with $g_{n,j} \in \bar{\mathcal{P}} \big( 8N^2Q \big)$ and
\begin{align}
      \sum_j \vertiii{ g_{n,j} }
            &\, \lesssim \kappa \mr\lambda^2 \lambda^{-s_nb\gamma},    \label{e:78:estQ+1} \\ 
      \| \nabla_\xi^p \ddot{\chi}_{n,j} \|_{\infty} 
            &\, \lesssim 1,       \quad p \leq 8Q^3,        \label{e:80:estQ+1} \\ 
      \| \partial_\tau^p \ddot{\eta}_{n,j} \|_{\infty} 
            &\, \lesssim 1,       \quad p \leq 8Q^3.        \label{e:82:estQ+1} 
\end{align}
Here, for any $j$, $\big( \ddot{\chi}_{n,j}, \ddot{\eta}_{n,j} \big)$ satisfies the generalized shear condition.

Next we expand the product of $g_n$ and $\mr D^{\bomega_n} \varrho_n$ and renumber all terms in the following form
\begin{align}
      g_n \mr D^{\bomega_n} \varrho_n = F_n + \tilde f_n + g_n \mr D^{\bomega_n} \tilde \varrho_n     \label{e:88:estQ+1} \\ 
\end{align}
with
\begin{align}
      F_n =&\, \sum_{a=0}^{(Q+1)(n+1)} \sum_{|\balpha|=1, \balpha \in \mathcal{I}_*}^{Q} 
            \mr{\bar D}^{\balpha} \bar \rho_{n,a} \sum_{l=1}^{I'_n} \mathcal{\bP}_Q h_{F,n,a}^{(\balpha,l)} \eta_{F,n,a}^{(\balpha,l)} (\mu t) \chi_{F,n,a}^{(\balpha,l)} (\lambda\Phi),   \label{e:90:estQ+1} \\ 
      \tilde f_n = &\, \sum_{a=0}^{(Q+1)(n+1)} \sum_{|\balpha|=1, \balpha \in \mathcal{I}_*}^{Q} 
            \mr{\bar D}^{\balpha} \bar \rho_{n,a} \sum_{l=1}^{I'_n} 
            \Big( h_{F,n,a}^{(\balpha,l)} - \mathcal{\bP}_Q h_{F,n,a}^{(\balpha,l)} \Big) \eta_{F,n,a}^{(\balpha,l)} (\mu t) \chi_{F,n,a}^{(\balpha,l)} (\lambda\Phi)         \label{e:92:estQ+1} \\ 
      &\,+ \sum_{a=0}^{(Q+1)(n+1)} \sum_{|\balpha|=Q, \balpha \in \mathcal{I}_*}^{Q+[\bomega_n]} 
            \mr{\bar D}^{\balpha} \bar \rho_{n,a} \sum_{l=1}^{I'_n} h_{F,n,a}^{(\balpha,l)} \eta_{F,n,a}^{(\balpha,l)} (\mu t) \chi_{F,n,a}^{(\balpha,l)} (\lambda\Phi),        \label{e:94:estQ+1} 
\end{align}
$I'_n \lesssim 1$, $h_{F,n,a}^{(\balpha,l)} \in \bar{\mathcal{P}}(9NQ)$ and
\begin{align}
      \vertiii{ h_{F,n,a}^{(\balpha,l)} } 
            &\,\lesssim \kappa \lambda \mr\lambda \cdot \frac{\delta^{\sfrac12}\mr{\bar\lambda}} {\kappa\lambda^2} \cdot \lambda^{-s_nb\gamma},
            \label{e:100:estQ+1} \\ 
      \| \nabla_\xi^p \chi_{F,n,a}^{(\balpha,l)} \|_{\infty} 
            &\,\lesssim 1,       \quad p \leq 8Q^3 - 6Q (n+1) - 3Q,     \label{e:101:estQ+1} \\ 
      \| \partial_\tau^p \eta_{F,n,a}^{(\balpha,l)} \|_{\infty} 
            &\,\lesssim 1,       \quad p \leq 8Q^3 - 6Q (n+1) - 3Q.     \label{e:102:estQ+1}
\end{align}
Naturally, for any $a,\balpha,l$, $\big( \chi_{F,n,a}^{(\balpha,l)}, \eta_{F,n,a}^{(\balpha,l)} \big)$ satisfies the generalized shear condition.

Now we define $\varrho_{n+1,0}$, $\varrho_{n+1,1}$ and $\varrho_{n+1,2}$ to be the solutions of 
\begin{align}
      O_t \varrho_{n+1,0} - \kappa \Delta \varrho_{n+1,0} =&\, \mr{\divr} F_n,      \label{e:104:estQ+1} \\ 
      \varrho_{n+1,0} (\cdot, 0) =&\, 0,    \label{e:106:estQ+1} \\
      O_t \varrho_{n+1,1} - \kappa \Delta \varrho_{n+1,1} =&\, \mr{\divr} \tilde f_n ,      \label{e:108:estQ+1} \\ 
      \varrho_{n+1,1} (\cdot, 0) =&\, 0,    \label{e:110:estQ+1} \\ 
      O_t \varrho_{n+1,2} - \kappa \Delta \varrho_{n+1,2} =&\, \mr{\divr} \big( g_n \mr D^{\bomega_n} \tilde \varrho_n \big), \label{e:112:estQ+1} \\ 
      \varrho_{n+1,2} (\cdot, 0) =&\, 0.    \label{e:114:estQ+1}
\end{align}
From Corollary \ref{c:CalPEst}, we have
\begin{align}
      \| \mr {\bar D}^{\bbeta} h_{F,n,a}^{(\balpha,l)} \|_\infty \leq&\, \kappa \lambda \mr\lambda \lambda^{-s_nb\gamma}
            \quad \text{for } [\bbeta] \leq Q^3       \label{e:116:estQ+1}
\end{align}
and subsequently from \eqref{e:28:estQ+1}, Lemma \ref{l:compDecompedF},
\begin{align}     \label{e:118:estQ+1}
      \| \mr O^{\bbeta} \tilde f_{n} \|_\infty 
            \leq&\, \lambda^{ - \frac{Q\gamma}{2} - \alpha - \big( [\bbeta] + \sum_{i=0}^{n} (s_i-1) \big) b \gamma }  \Upsilon_{\ini}
            \quad \text{for } [\bbeta] \leq Q + 3Q(Q-n-1).
\end{align}
Here, we gain the negative power $-\frac{Q\gamma}{2}$ for the terms in \eqref{e:92:estQ+1} and \eqref{e:94:estQ+1}. For \eqref{e:92:estQ+1}, we gain this power because of Corollary \ref{c:CalPEst} and the lower bound $Q$ on the differentiation order of $h_{F,n,a}^{(\balpha,l)} - \mathcal{\bP}_Q h_{F,n,a}^{(\balpha,l)}$. For \eqref{e:94:estQ+1}, we gain this power because of the lower bound $Q$ on $|\balpha|$ and Lemma \ref{l:compDecompedF}.

With \eqref{e:118:estQ+1}, we apply energy estimates in Lemma \ref{l:L2_F_energy} to deduce
\begin{align}
      \| \mr O^{\bbeta} \varrho_{n+1,1} \|
            \leq&\, \lambda^{ - 3 - \big( [\bbeta] + \sum_{i=0}^{n} (s_i-1) \big) b \gamma } \Upsilon_{\ini}  
            \quad \text{for } [\bbeta] \leq 2Q + 3Q(Q-n-1). 
\end{align}
From \eqref{e:30:estQ+1} and energy estimate in Corollary \ref{c:L2_FHom_energy}, we also have
\begin{align}     \label{e:124:estQ+1}
      \| \mr O^{\bbeta} \varrho_{n+1,2} \|
            \leq&\, \lambda^{ - 3 - \big( [\bbeta] + \sum_{i=0}^{n} (s_i-1) \big) b \gamma }  \Upsilon_{\ini} 
            \quad \text{for } [\bbeta] \leq 2Q + 3Q(Q-n-1). 
\end{align}
For $\varrho_{n+1,0}$, we apply the expansion in Lemma \ref{l:reEXP2} to (\ref{e:104:estQ+1}-\ref{e:106:estQ+1}). Here we briefly verify the assumptions. The conditions in Assumption \ref{a:homSmallscl} follow from \eqref{e:100:estQ+1}-\eqref{e:102:estQ+1} and the definition of the projection $\mathbb{P}$. Suitably redistributing constants for $h_{F,\cdot}^{(\cdot,\cdot)}$ and $\bar\rho_{F,\cdot}$, we can verify \eqref{e:Up:reEXP2}-\eqref{e:macroIn:reEXP2} from \eqref{e:18:estQ+1} and \eqref{e:27:estQ+1}-\eqref{e:29:estQ+1}. Then we get the following expansion
\begin{align*}
      \varrho_{n+1,0} 
      =&\, \sum_{a=0}^{(Q+1)(n+2)} \sum_{|\balpha|=1, \balpha \in \mathcal{I}_*}^{Q} 
            \mr {\bar D}^{\balpha} \bar \rho_{n+1,a} 
            \sum_{l=1}^{I_{n+1}} h_{n+1,a}^{(\balpha,l)} \eta_{n+1,a}^{(\balpha,l)} (\mu t) \chi_{n+1,a}^{(\balpha,l)} (\lambda\Phi) 
            + \sum_{a=0}^Q \bar\rho_{n+1,a} + \tilde \rho_{F} 
\end{align*}
with $I_{n+1} \lesssim 1$, $h_{n+1,a}^{(\balpha,l)} \in \bar{\mathcal{P}} \big( Q \big)$ and
\begin{align*}
      \sum \vertiii{h_{n+1,a}^{(\balpha,l)}} 
            \lesssim \frac{\delta^{\sfrac12}\mr{\bar\lambda}} {\kappa\lambda^2}.
\end{align*}
For the correctors, we have
\begin{align*}
      \| \nabla_\xi^p \chi_{n+1,a}^{(\balpha,l)} \|_{\infty} 
            &\, \lesssim 1,       \quad p \leq 8Q^3 - 6Q (n+1) - 4Q,    \\ 
      \| \partial_\tau^p \eta_{n+1,a}^{(\balpha,l)} \|_{\infty} 
            &\, \lesssim 1,       \quad p \leq 8Q^3 - 6Q (n+1) - 4Q.
\end{align*}
Here, for any $a, \balpha, l$, $\big( \chi_{n+1,a}^{(\balpha,l)}, \eta_{n+1,a}^{(\balpha,l)} \big)$ satisfies the generalized shear condition. $\bar\rho_{n+1,a}$ satisfies \eqref{e:48:estQ+1} for any $0 \leq a \leq (Q+1)(n+2)$. We also have
\begin{align*}
      \| \mr O^{\bbeta} \tilde \rho_F \|_\infty \leq&\, 
            \lambda^{ - 3 - \big( [\bbeta] + \sum_{i=0}^{n} (s_i-1) \big) b \gamma } \Upsilon_{\ini} 
            \quad \text{for } [\bbeta] \leq Q + 4Q^2(Q-n-1).  \\ 
\end{align*}
Finally, setting
\begin{align*}
      \tilde \rho_{n+1} = \tilde \rho_F + \varrho_{n+1,1} + \varrho_{n+1,2}
\end{align*}
concludes the proof.

\end{proof}

\newpage

\section{Homogenization in dissipative range}         \label{s:homDissip}

In this section, we prove the inductive homogenization estimates from step $q$ to step $q+1$ for dissipative range. In dissipative range, the difference between eddy diffusivity and molecular diffusivity is very small and can be treated as a perturbation. Throughout this section, we work with diffusivity $\kappa \geq 0$ satisfying
\begin{align}     \label{e:D:ConDiffusivity}
      \kappa \geq \kappa_q.
\end{align}

\subsection{Main results of this section}       \label{ss:mainR_dissH}

Using the notations introduced in Section \ref{ss:mainR_inerH}, we state the main results as follows

\begin{theorem}   \label{t:D:homDissip}
Fix some $q_* \leq q$ with \eqref{e:D:ConDiffusivity}. Suppose that, for any $\rho_{\ini} \in C^\infty (\T^2)$, we have
\begin{align}     
      \resN{ \rho_{\ini} }_{u_q,\kappa} 
            \leq&\, \lambda_q^{-\alpha} \bigg( | \rho_{\ini} |_{\mathfrak{D}(u_q,\kappa)} 
                  + \sum_{j=q_*}^{q} \lambda_j^{-\sfrac{\gamma}{2}} | \rho_{\ini} |_{\fS(j)} \bigg),   \label{e:1:D:homDissip} 
\end{align}
then
\begin{align}     \label{e:2:D:homDissip}
      \left| | \rho_{\ini} |_{\mathfrak{D}(u_q,\kappa)} 
            - | \rho_{\ini} |_{\mathfrak{D}(u_{q+1},\kappa)} \right|
            \leq&\, \lambda_{q+1}^{-\gamma} \bigg( | \rho_{\ini} |_{\mathfrak{D}(u_q,\kappa)} 
                  + \sum_{j=q_*}^{q+1} \lambda_j^{-\sfrac{\gamma}{2}} | \rho_{\ini} |_{\fS(j)} \bigg).
\end{align}
\end{theorem}

\begin{theorem}   \label{t:D:homRslv}
Fix some $q_* \leq q$ with \eqref{e:D:ConDiffusivity}. Suppose that, for any $\rho_{\ini} \in C^\infty (\T^2)$, we have
\begin{align}     
      \resN{ \rho_{\ini} }_{u_q,\kappa} 
            \leq&\, \lambda_q^{-\alpha} \bigg( | \rho_{\ini} |_{\mathfrak{D}(u_q,\kappa)} 
                  + \sum_{j=q_*}^{q} \lambda_j^{-\sfrac{\gamma}{2}} | \rho_{\ini} |_{\fS(j)} \bigg),   \label{e:2:D:homRslv} \\ 
      \resN{ \rho_{\ini} }^{(H,\alpha_0)}_{u_q,\kappa} 
            \leq&\, \left( 1 - \frac{4}{q} \right) \bigg( | \rho_{\ini} |_{\mathfrak{D}(u_q,\kappa)} 
                  + \sum_{j=q_*}^{q} \lambda_j^{-\sfrac{\gamma}{2}} | \rho_{\ini} |_{\fS(j)} \bigg), \label{e:4:D:homRslv}
\end{align}
then for any $\rho_{\ini} \in C^\infty (\T^2)$, we have
\begin{align}
      \resN{ \rho_{\ini} }_{u_{q+1},\kappa} 
            \leq&\, \lambda_{q+1}^{-\alpha} \bigg( | \rho_{\ini} |_{\mathfrak{D}(u_{q+1},\kappa)} 
                  + \sum_{j=q_*}^{q+1} \lambda_j^{-\sfrac{\gamma}{2}} | \rho_{\ini} |_{\fS(j)} \bigg),      \label{e:6:D:homRslv} \\
      \resN{ \rho_{\ini} }^{(H,\alpha_0)}_{u_{q+1},\kappa} 
            \leq&\, \left( 1 - \frac{4}{q+1} \right) \bigg( | \rho_{\ini} |_{\mathfrak{D}(u_{q+1},\kappa)}
                  + \sum_{j=q_*}^{q+1} \lambda_j^{-\sfrac{\gamma}{2}} | \rho_{\ini} |_{\fS(j)} \bigg).     \label{e:8:D:homRslv}
\end{align}
\end{theorem}

The proof strategy of Theorem \ref{t:D:homRslv} is the same as Theorem \ref{t:homRslv}. Since we have a relatively large diffusivity $\kappa$, the expansion lemma and its proof can be simplified. At heuristic level, since the difference between eddy diffusivity and molecular diffusivity is very small, we treat this difference as a perturbation. Hence, we do not need to solve cell problem in time in homogenization.

\subsection{Notations of this section} Consider the following equation for $\varrho: \T^2 \times [0,1] \rightarrow \R$,
\begin{equation}  \label{e:D:q+1:meq}
\begin{split}
      O_{q+1,t} \varrho - \kappa \Delta \varrho =&\, 0, \\ 
      \varrho (\cdot, 0) =&\, \rho_{\ini}.
\end{split}
\end{equation}

As in Section \ref{ss:nota_inerH}, we use the notation introduced in Definition \ref{d:homFlowmap}. We also introduce the following simplified notations, slightly different from those in Definition \ref{d:homParaI}.

\begin{definition}[Parameter]   \label{d:D:homParaI}
Define
\begin{align}     \label{d:D:mainPara}
      \lambda := \lambda_{q+1},            \quad
      \delta := \delta_{q+1},            \quad
      \mu := \mu_{q+1},            \quad 
      \bar\lambda := \lambda_{q},      \quad
      \bar\delta := \delta_{q},            \quad
      \bar\mu := \mu_{q},  
\end{align}
\begin{align}     \label{d:D:circleLambda}
      \mr {\lambda} := \mr {\lambda}_{q+1},      \quad
      \mr {\mu} := \mr {\mu}_{q+1},       \quad
      \mr {\bar\lambda} := \mr {\lambda}_{q},      \quad
      \mr {\bar\mu} := \mr {\mu}_{q},  
\end{align} 
\begin{align}     \label{d:D:varepsilon} 
      \varepsilon := \frac{\bar\delta^{\sfrac12} \bar\lambda^{1+2b\gamma_R}} {\mu}
\end{align}
From \eqref{e:4:mScaleRela}, we have
\begin{align}     \label{d:D:lambda_gamma}
      \frac{\mu} {\delta^{\sfrac12}\mr{\bar\lambda}} \leq \lambda^{\gamma_S}.
\end{align}
From \eqref{e:smallGammaR} and \eqref{e:2:mScaleRela}, we have
\begin{align}     \label{e:D:homGammaRela}
      \frac{\mu}{\kappa\lambda^2}
            + \frac{ \delta^{\sfrac12}\mr{\bar\lambda} } {\kappa\lambda^2} 
            \cdot \lambda^{2b\gamma+\gamma_S} \leq 1.
\end{align}
\end{definition}

\subsection{Expansion lemma}

\begin{lemma}[The expansion lemma]        \label{l:D:EXP1}
For $N,Q$ chosen in Section \ref{ss:parameters}, consider the equation \eqref{e:D:q+1:meq}, under Definition \ref{d:D:homParaI} and Definition \ref{d:homFlowmap}. Then the solution $\varrho$ admits the expansion
\begin{equation}  \label{e:D:0:EXP1}
\begin{split}
      \varrho(x,t) =& \sum_{n = 0}^{Q} \frac{1}{\lambda^n} 
            \Big( \hat \rho_{n} (x,t, \lambda \Phi(x,t) ) 
            + \bar \rho_{n} (x, t) \Big) 
            + \tilde \rho(x,t), 
\end{split}
\end{equation}
with $\hat \rho_n, \bar \rho_n, \tilde \rho$ given by the following.

\begin{enumerate}[leftmargin=*,label=\textsc{(C.\arabic*)},align=left]
\item \label{c:D:l:exp1:spt} \textbf{\textup{(Spatial correctors)}} $\hat \rho_n : \T^2 \times [0,1] \times \T^2 \rightarrow \R$ has the form
\begin{align}
      \hat \rho_0 =&\, 0,       \label{e:D:EXP1:hrho0} \\ 
      \hat \rho_{n}(x,t,\xi) =&\, \sum_{a=0}^{n-1} \sum_{ [\balpha]=1,\balpha \in \mathcal{I}_* }^{2(n-a)-1} \sum_{l=1}^{\hat I_{n-a}^{(\balpha)}}
            \hat h^{(\balpha, l)}_{n-a} (x,t) \hat \eta^{(\balpha, l)}_{n-a} ( \mu t ) \hat \chi^{(\balpha, l)}_{n-a} (\xi) \mr {\bar D}^{\balpha} \bar \rho_a(x,t)       \label{e:D:EXP1:hrhon} \\ 
      \hat h^{(\balpha, l)}_{n} \in&\, \mathcal{\bar P} ( N(n-1)+2-[\balpha], 0 )       \label{e:D:EXP1:hh}
\end{align}
with $\hat h^{(\cdot,\cdot)}_{\cdot} : \T^2 \times [0,1] \rightarrow \R$, $\hat \chi_{\cdot}^{(\cdot,\cdot)} : \T^2 \rightarrow \R$ and $\hat \eta_{\cdot}^{(\cdot,\cdot)} : \T \rightarrow \R $ satisfying the estimates for $k \geq 0$
\begin{align}
      \hat I^{(\balpha)}_{k+1} &\, \leq 2(3N)^{k},   \label{e:D:EXP1:hI_est} \\ 
      \vertiii{ \hat h^{(\balpha, l)}_{k+1} } &\, \lesssim \frac{\delta^{\sfrac12}\mr{\bar\lambda}}{\kappa\lambda} ( \varepsilon\lambda )^{k},         \label{e:D:EXP1:hh_Pest}   \\ 
      \sum_{[\balpha]=1}^{2k+1} \sum_{l=1}^{\hat I^{(\balpha)}_{k+1}} 
            \vertiii{ \hat h^{(\balpha, l)}_{k+1} } &\, \lesssim \frac{\delta^{\sfrac12}\mr{\bar\lambda}}{\kappa\lambda} ( \varepsilon\lambda )^{k},         \label{e:D:EXP1:hhS_Pest}   \\ 
      \| \nabla_\xi^p \hat\chi^{(\balpha,l)}_{k+1} \|_{\infty} &\, \lesssim 1,
            \quad p \leq 11Q^3, 
            \quad \langle \hat\chi^{(\balpha,l)}_{k+1} \rangle_\xi = 0,     \label{e:D:EXP1:hchi_est} \\
      \| \partial_\tau^p \hat\eta^{(\balpha,l)}_{k+1} \|_{\infty} &\, \lesssim 1,
            \quad p \leq 11Q^3-k-1.       \label{e:D:EXP1:heta_est}
\end{align}
Here the sum of differentiation index $\balpha$ in \eqref{e:D:EXP1:hrhon} is among all indices in $\mathcal{I}_*$ such that $1 \leq [\balpha] \leq 2(n-a)-1$. Moreover, $(\hat \chi^{(\balpha, l)}_{n}, \hat \eta^{(\balpha, l)}_{n})$ satisfy the shear condition in Definition \ref{d:shear_structure} for any $n$, $\balpha$ and $l$.

\item \label{c:D:l:exp1:rsd} \textbf{\textup{(Residual correctors)}} $\bar \rho_0, \{\bar \rho_n\}_{1 \leq n \leq N}: \T^2 \times [0,1] \rightarrow \R $ solve
\begin{align}
      \bar L_0 \bar \rho_0 =&\, 0,   \label{e:D:EXP1:brho0} \\ 
      \bar L_0 \bar \rho_n =&\, \sum_{a=0}^{n-1} \sum_{ [\balpha]=1, \balpha \in \mathcal{I}_* }^{2(n-a)+1} \sum_{l=1}^{\bar I_{n}^{(\balpha)}}
      \mr{\bar \divr} \Big( \bar h^{(\balpha,l)}_{n-a} \odot \bar \eta^{(\balpha,l)}_{n-a}(\mu t) \mr {\bar D}^{\balpha} \bar \rho_a \Big)   \label{e:D:EXP1:brhon} \\ 
      \bar h^{(\balpha)}_{n} \in&\, \mathcal{\bar P} ( Nn+3-[\balpha], 0 )       \label{e:D:EXP1:bh} 
\end{align}
with $\bar L_0$ given by
\begin{align}   \label{e:D:EXP1:homedOp}  
      \bar L_0 = \bar D_t - \kappa \Delta
\end{align}
and initial datum given by
\begin{align}
      \bar \rho_0(\cdot, 0) =&\, \rho_{\ini},     \label{e:D:EXP1:brho0_ini} \\ 
      \bar \rho_n(\cdot, 0) =&\, 0,       \quad n \geq 1.      \label{e:D:EXP1:brhon_ini}
\end{align}
Here, \underline{the vector-valued functions} $\bar h^{(\cdot,\cdot)}_{\cdot} : \T^2 \times [0,1] \rightarrow \R^2$ and $\bar \eta^{(\cdot,\cdot)}_{\cdot} : \T \rightarrow \R^2$ satisfy the estimates for $k \geq 0$
\begin{align}
      \vertiii{ \bar h^{(\balpha,l)}_{k+1} } \lesssim&\, 
      \kappa \bigg( \frac{ \delta^{\sfrac12}\mr{\bar\lambda} } {\kappa\lambda} \bigg)^2
      ( \varepsilon\lambda )^{k}.       \label{e:D:EXP1:bh_Pest} \\ 
      \sum_{[\balpha]=1}^{2k+3}
      \vertiii{ \bar h^{(\balpha,l)}_{k+1} } \lesssim&\, 
      \kappa \bigg( \frac{ \delta^{\sfrac12}\mr{\bar\lambda} } {\kappa\lambda} \bigg)^2
      ( \varepsilon\lambda )^{k}.       \label{e:D:EXP1:bhS_Pest} \\ 
      \| \partial_\tau^p \bar \eta^{(\balpha,l)}_{k+1} \|_{\infty} &\, \lesssim 1,
            \quad p \leq 11Q^3-k-1.
\end{align}

\item \label{c:D:l:exp1:rmd} \textbf{\textup{(remainder)}} $\tilde \rho$ solves the following equation
\begin{align}
      \bar D_t \tilde\rho - \divr( A \nabla \tilde\rho ) 
      =&\, \tilde f,    \label{e:D:exp1:tildeRho} \\ 
      \tilde \rho(x,0) =&\, 0.      \label{e:D:0:tildeRho}
\end{align}
Here, $\tilde f: \T^2 \times [0,1] \rightarrow \R$ is given by
\begin{align}
      \tilde f := &\, -\frac{1}{\lambda^{Q}} 
            L_0 ( \hat\rho_{Q} + \bar\rho_{Q} )       \label{e:D:2:tildeRho} 
\end{align}
and differential operators $L_1,L_0,L_{-1}$ acting on general functions $\rho$ taking arguments $(x,t,\xi,\tau)$ are given by
\begin{align}
      L_0 \rho =&\, - \delta^{\sfrac12} \partial_{\xi_i} H_{ij} \partial_{x_j} \rho + \bar D_t \rho - \kappa \Delta_x \rho - \varepsilon \frac{\delta^{\sfrac12}}{\lambda} \sum_m \sigma_m (\mu \cdot) \partial_{x_i} \Omega_{m,ij} H_{12} \partial_{x_j} \rho    \\ 
            -&\, \varepsilon \delta^{\sfrac12} \sum_m \varphi_m (\mu \cdot) \partial_{\xi_i} ( H E_m^T )_{ij} \partial_{x_j} \rho + \varepsilon \delta^{\sfrac12} \sum_m \varphi_m (\mu \cdot) \partial_{x_j} ( HE_m^T )_{ij} \partial_{\xi_i} \rho      \\ 
            -&\, \varepsilon \kappa \lambda \sum_m \phi_m (\mu \cdot) \partial_{x_j} B_{m,ij} \partial_{\xi_i} \rho - 2 \kappa \lambda \partial_{\xi_i x_i} \rho - 2 \varepsilon \kappa \lambda \sum_m \phi_m (\mu \cdot) B_{m,ij} \partial_{\xi_i x_j} \rho      \\ 
            +&\, \varepsilon \kappa \lambda^2 \sum_m S_{m,ij} \vartheta_m (\mu \cdot) \partial_{\xi_i\xi_j} \rho + \frac{1}{\lambda^2} \phi_* (\mu \cdot) \omega_* \cdot \nabla_\xi \rho 
\end{align}
\end{enumerate}
\end{lemma}

\begin{proof}[Sketch of proof to Theorem \ref{t:D:homRslv}]
Compared with Theorem \ref{t:homRslv}, the difference in Theorem \ref{t:D:homRslv} is only due to the larger diffusivity $\kappa$.

To prove Theorem \ref{t:D:homRslv}, it suffices to verify two main ingredients. The first ingredient is unconditional energy estimates, in the form of Proposition \ref{p:unconEnergyHom} and Proposition \ref{p:chainStab}. From Lemma \ref{l:L2_init_energy}, Lemma \ref{l:L2_F_energy} and Lemma \ref{c:L2_FHom_energy}, these obvious hold, since the larger diffusivity yields stronger energy estimates.

The second ingredient is the homogenization estimates via asymptotic expansion. The whole homogenization analysis is essentially consequences of Lemma \ref{l:EXP1} and Lemma \ref{l:EXP2}. And Lemma \ref{l:EXP2} is a variant of Lemma \ref{l:EXP1}. Therefore, it boils down to prove Lemma \ref{l:D:EXP1}, which is the variant of Lemma \ref{l:EXP1} in the dissipative range.
\end{proof}

\subsection{Proof of Lemma \ref{l:D:EXP1}}      \label{ss:D:proofExp}

In this section, we prove Lemma \ref{l:D:EXP1}. The proof is largely similar to the proof of Lemma \ref{l:EXP1}, except for two different points. The first point is that we choose to deal with easier cell problems, which are set up in Section \ref{sss:D:setUpAsymp} below. The second point is the different parameter relations, hence we need to verify the expansion converges under these parameter relations. This is done in Section \ref{sss:D:solInitAspt}, Section \ref{sss:D:expIndRsd} and Section \ref{sss:D:expIndSpt}. We shall omit many other similar details.

\subsubsection{Set up the asymptotic system for Lemma \ref{l:D:EXP1}}      \label{sss:D:setUpAsymp}

As in inertial range homogenization, we have
\begin{align}     \label{e:D:aspt:meq}
    L \varrho = \bar D_t \varrho - \divr ( A \nabla \varrho ) = 0.
\end{align}
We make the following \textit{ansatz}
\begin{align}     \label{e:D:aspt:ansatzPre}
      \varrho (x,t) &= \sum_{n = 0}^{Q} \frac{1}{\lambda^n} 
            \rho_{n} (x,t, \lambda \Phi(x,t) ) + \tilde \rho(x,t).
\end{align}
for some function $ \rho_{n} : \T^2 \times [0,1] \times \T^2 \rightarrow \R $. 

Define $ \hat \rho_{n} : \T^2 \times [0,1] \times \T^2 \rightarrow \R $ with arguments $(x,t,\xi)$ via
\begin{align}
      \bar \rho_n(x,t) :=&\, \langle \rho_{n} ( x,t, \cdot ) \rangle_{\xi},    \label{e:D:1:asptDecomp} \\
      \hat \rho_n(x,t,\xi) :=&\, \rho_{n} (x,t, \xi) - \langle \rho_{n} (x,t, \cdot) \rangle_{\xi}.        \label{e:D:3:asptDecomp}
\end{align}
Then we have the natural decomposition for $\rho_n$ with $n \geq 0$
\begin{align}     \label{e:D:4:asptDecomp}
      \rho_n(x,t,\xi) = \hat \rho_n(x,t,\xi) + \bar \rho_n(x,t)
\end{align}
and subsequently
\begin{align}     \label{e:D:aspt:ansatz}
      \varrho (x,t) &= \sum_{n = 0}^{Q} \frac{1}{\lambda^n} 
            \big( \hat \rho_{n} (x,t, \lambda \Phi(x,t) ) + \bar \rho_{n} (x, t) \big) 
            + \tilde \rho(x,t). 
\end{align}
We have, for any $x,t$,
\begin{align}     \label{e:D:16:SetAspt}
      \langle \hat \rho_{n}( x,t, \cdot ) \rangle_\xi = 0.
\end{align}

An analogue of Lemma \ref{l:operator} is given as follows. There are two crucial differences. Since we do not homogenize in time variable $\mu t$, the derivatives in variable $\mu t$ are treated as slow/macroscopic derivatives. The second difference is that the operators in cell problems are different, because of the simplification given by stronger diffusion.

\begin{lemma}     \label{l:D:operator}
Under Definition \ref{d:D:homParaI} and Definition \ref{d:homFlowmap}, for a general function $\varrho: \T^2 \times [0,1] \rightarrow \R$ given by 
\begin{align}     \label{e:D:0:operator}
      \varrho(x,t) = \rho( x,t, \lambda \Phi(x,t) ),
\end{align}
not necessarily solving \eqref{e:D:aspt:meq}, we have
\begin{align}     \label{e:D:1:operator}
      L \varrho = \lambda L_1 \rho + L_0 \rho
\end{align}
with
\begin{align*}
      L_1 \rho =&\, - \kappa \lambda \Delta_\xi \rho - \delta^{\sfrac12} (\det \nabla\Phi)^2 \divr_\xi H \cdot \nabla_\xi \rho 
            + \sum_m z_m (\mu \cdot) \omega_m \cdot \nabla_\xi \rho ,   \\ 
      L_0 \rho =&\, - \delta^{\sfrac12} \partial_{\xi_i} H_{ij} \partial_{x_j} \rho + \bar D_t \rho - \kappa \Delta_x \rho - \varepsilon \frac{\delta^{\sfrac12}}{\lambda} \sum_m \sigma_m (\mu \cdot) \partial_{x_i} \Omega_{m,ij} H_{12} \partial_{x_j} \rho \\ 
            -&\, \varepsilon \delta^{\sfrac12} \sum_m \varphi_m (\mu \cdot) \partial_{\xi_i} ( H E_m^T )_{ij} \partial_{x_j} \rho + \varepsilon \delta^{\sfrac12} \sum_m \varphi_m (\mu \cdot) \partial_{x_j} ( HE_m^T )_{ij} \partial_{\xi_i} \rho 
            - \varepsilon \kappa \lambda \sum_m \phi_m (\mu \cdot) \partial_{x_j} B_{m,ij} \partial_{\xi_i} \rho        \\ 
            -&\,2 \kappa \lambda \partial_{\xi_i x_i} \rho - 2 \varepsilon \kappa \lambda \sum_m \phi_m (\mu \cdot) B_{m,ij} \partial_{\xi_i x_j} \rho + \varepsilon \kappa \lambda^2 \sum_m S_{m,ij} \vartheta_m (\mu \cdot) \partial_{\xi_i\xi_j} \rho + \frac{1}{\lambda^2} \phi_* (\mu \cdot) \omega_* \cdot \nabla_\xi \rho      \\ 
            :=&\, \sum_{0 \leq k \leq 9} L_{0k} \rho 
\end{align*}
Here, we omit the argument $(\xi,\mu t)$ of $H$. We define $L_{00} := - \delta^{\sfrac12} \partial_{\xi_i} H_{ij} \partial_{x_j} \rho$, $L_{01} \rho := D_t \rho - \kappa \Delta_x \rho$, and $L_{0k}, k \geq 2$ is defined by the order of relevant terms.
\end{lemma}

\begin{proof}
The proof is similar to the proof of Lemma \ref{l:operator}. Computing the terms $\divr( A \nabla \varrho )$, namely $\divr_\xi ( \nabla \Phi A \nabla \Phi^T \nabla_\xi \rho )$, $\divr_\xi ( \nabla \Phi A \nabla_x \rho )$, $\divr_x ( A \nabla \Phi^T \nabla_\xi \rho )$ and $\divr_x ( A \nabla_x \rho )$ is done in (\ref{e:5:operator}-\ref{e:14:operator}).

For the transport term, we deduce from \eqref{e:0:operator} and \eqref{e:12:homFlowmap}, that
\begin{align}
      \bar D_t \varrho =&\, \partial_t \rho + \lambda \partial_t \Phi_i \partial_{\xi_i} \rho + \bar u \cdot \nabla_x \rho + \lambda \bar u \cdot \nabla \Phi_i \partial_{\xi_i} \rho    \nonumber \\ 
      =&\, \bar D_t \rho + \lambda \bar D_t \Phi \cdot \nabla_\xi \rho       \nonumber \\ 
      =&\, \bar D_t \rho + \lambda Z \cdot \nabla_\xi \rho + \lambda^{-2} \Omega \cdot \nabla_\xi \rho.     \label{e:D:18:operator} 
\end{align}

From \eqref{e:D:aspt:meq}, \eqref{e:D:0:operator}, (\ref{e:5:operator}-\ref{e:14:operator}) and \eqref{e:18:homFlowmap}, we have \eqref{e:1:operator} with
\begin{align}
      L_1 \rho =&\, - \kappa \lambda \Delta_\xi \rho - \delta^{\sfrac12} (\det \nabla\Phi)^2 \divr_\xi H \cdot \nabla_\xi \rho 
            + Z \cdot \nabla_\xi \rho ,   \label{e:D:20:operator} \\ 
      L_0 \rho =&\, - \delta^{\sfrac12} \partial_{\xi_i} H_{ij} \partial_{x_j} \rho + \bar D_t \rho - \kappa \Delta_x \rho - \frac{\delta^{\sfrac12}}{\lambda} \partial_{x_i} ( \adj \nabla \Phi H \adj \nabla \Phi^{T} )_{ij} \partial_{x_j} \rho      \label{e:D:22:operator} \\ 
            -&\, \varepsilon \delta^{\sfrac12} \sum_m \varphi_m \partial_{\xi_i} ( H E_m^T )_{ij} \partial_{x_j} \rho
            + \varepsilon \delta^{\sfrac12} \sum_m \varphi_m \partial_{x_j} ( HE_m^T )_{ij} \partial_{\xi_i} \rho 
            - \kappa \lambda \partial_{x_j} \nabla \Phi_{ij} \partial_{\xi_i} \rho        \label{e:D:24:operator} \\ 
            -&\, 2 \kappa \lambda \nabla \Phi_{ij} \partial_{\xi_i x_j} \rho + \varepsilon \kappa \lambda^2 \divr_\xi ( \Xi \nabla_\xi \rho ) + \frac{1}{\lambda^2} \Omega \cdot \nabla_\xi \rho,     \label{e:D:26:operator} 
\end{align}

Now we use \eqref{e:6:homFlowmap}, \eqref{e:16:homFlowmap} and \eqref{e:18:homFlowmap} and above to conclude the proof.

\end{proof}

From Lemma \ref{l:D:operator}, we have
\begin{equation}        \label{e:D:asptExp} 
      \begin{split}
      L \varrho =&\, \lambda L_1 \rho_0 + L_1 \rho_1 + L_0 \rho_0 
            + \frac{1}{\lambda} ( L_1 \rho_2 + L_0 \rho_1 ) + \ldots       \\ 
            +&\, \frac{1}{\lambda^{Q-1}} ( L_1 \rho_{Q} + L_0 \rho_{Q-1} )
            + \frac{1}{\lambda^Q} L_0 \rho_Q + L \tilde \rho
      \end{split}
\end{equation}

Then we seek for $\{\rho_n\}_{0 \leq n \leq Q}$ such that 
\begin{align}
      L_2 \rho_0 =&\, 0,     \label{e:D:aspt:hrc1} \\ 
      L_1 \rho_1 + L_0 \rho_0 =&\, 0,     \label{e:D:aspt:hrc2} \\ 
      L_1 \rho_{n+1} + L_0 \rho_n =&\, 0    \label{e:D:aspt:hrc3}
\end{align}
with initial conditions
\begin{align}
      \rho_0( x,0,\xi ) =&\, \rho_{\ini} (x),     \label{e:D:0:initialC} \\ 
      \rho_n( x,0,\xi ) =&\, 0, \quad \text{for any } n \geq 1.   \label{e:D:2:initialC}
\end{align}
Then $\tilde \rho$ satisfies
\begin{equation}        \label{e:D:eTildeRho}
      \begin{split}
      \bar D_t \tilde \rho - \divr ( A \nabla \tilde \rho ) =&\, \tilde f 
            := - \frac{1}{\lambda^Q} (L_0 \rho_Q) (x,t,\lambda\Phi(x,t)), \\ 
      \tilde \rho (x,0) =&\, 0.
      \end{split}
\end{equation}

\textit{If we have $\{\rho_n\}_{0 \leq n \leq Q}$ and $\tilde \rho$ such that \eqref{e:D:aspt:hrc1}-\eqref{e:D:eTildeRho} hold, the function $\varrho$ given by \eqref{e:D:aspt:ansatz} solves \eqref{e:D:aspt:meq}. Thus, to prove Lemma \ref{l:D:EXP1}, it remains to solve \eqref{e:D:aspt:hrc1}-\eqref{e:D:eTildeRho}. }

\begin{remark}    \label{r:D:hatcheckbar}
For $\rho(x,t,\tau) = \hat\rho(x,t,\tau) + \bar\rho(x,t)$ with
\begin{align*}
      \langle \hat \rho (x,t,\cdot) \rangle_\xi = 0 \quad \text{for any }x,t,
\end{align*}
we have
\begin{align}
      \langle L_1 \rho \rangle_\xi =&\, 0,     \label{e:D:6:hatcheckbar} \\ 
      \langle L_0 \rho \rangle_\xi =&\, \langle L_{00} \hat\rho \rangle_\xi + L_{01} \bar \rho + \langle L_{02} \hat\rho \rangle_\xi + \langle L_{03} \hat\rho \rangle_\xi + \langle L_{04} \hat\rho \rangle_\xi,       \label{e:D:8:hatcheckbar} \\ 
      L_0 \rho =&\, \sum_{k=0}^3 L_{0k} ( \hat\rho + \bar\rho ) + \sum_{k=4}^9 L_{0k} \hat\rho.        \label{e:D:10:hatcheckbar}
\end{align}
\end{remark}

\subsubsection{Solve the asymptotic system: initial steps}     \label{sss:D:solInitAspt}

In this step, we solve for
\begin{align}
      \hat \rho_0, \hat \rho_1, \bar \rho_0.
\end{align}

\eqref{e:D:aspt:hrc1} gives $\hat \rho_0 = 0$. Taking $\xi$-average of \eqref{e:D:aspt:hrc2}, we have 
\begin{align}
      \bar D_t \bar\rho_0 - \kappa \Delta_x \bar\rho_0 = 0. 
\end{align}

Then we have
\begin{align}
      L_1 \hat\rho_1 =&\, \langle L_0\rho_0 \rangle_\xi - L_0\rho_0 
            = - L_{00} \bar \rho_0 - L_{02} \bar \rho_0 - L_{03} \bar \rho_0
\end{align}
and
\begin{align}
      - L_{00} \bar \rho_0 =&\, \delta^{\sfrac12} \partial_{\xi_i} H_{ij} \partial_{x_j} \bar \rho_0,       \\ 
      - L_{02} \bar \rho_0 =&\, \varepsilon \frac{\delta^{\sfrac12}}{\lambda} \sum_m \sigma_m (\mu \cdot) \partial_{x_i} \Omega_{m,ij} H_{12} \partial_{x_j} \bar \rho_0,     \\ 
      - L_{03} \bar \rho_0 =&\, \varepsilon \delta^{\sfrac12} \sum_m \varphi_m (\mu \cdot) \partial_{\xi_i} ( H E_m^T )_{ij} \partial_{x_j} \bar \rho_0 
\end{align}
and subsequently
\begin{align}
      \hat \rho_1 =&\, J_1 + J_2 + J_3,      \\ 
      J_1 (x,t,\xi) =&\, - \frac{ \delta^{\sfrac12} \mr{\bar\lambda} }{\kappa\lambda} \eta_1 (\mu t) \partial_{\xi_1}^{-1} \Pi_1 (\xi) \mr {\bar \partial}_{2} \bar \rho_0(x,t)      \\ 
            &\,+ \frac{ \delta^{\sfrac12} \mr{\bar\lambda} }{\kappa\lambda} \eta_2 (\mu t) \partial_{\xi_2}^{-1} \Pi_2 (\xi) \mr {\bar \partial}_{1} \bar \rho_0(x,t), \\ 
      J_2 (x,t,\xi) =&\, - \frac{ \varepsilon \delta^{\sfrac12} \mr{\bar\lambda}^2 }{\kappa\lambda^2} \sum_m [\eta_1\sigma_m] (\mu t) \mr {\bar\partial}_i \Omega_{m,ij} (x,t) \partial_{\xi_1}^{-2} \Pi_1 (\xi) \mr {\bar\partial}_j \bar \rho_0(x,t)      \\ 
            &\,- \frac{ \varepsilon \delta^{\sfrac12} \mr{\bar\lambda}^2 }{\kappa\lambda^2} \sum_m [\eta_2\sigma_m] (\mu t) \mr {\bar\partial}_i \Omega_{m,ij} (x,t) \partial_{\xi_2}^{-2} \Pi_2 (\xi) \mr {\bar\partial}_j \bar \rho_0(x,t)    \\ 
      J_3 (x,t,\xi) =&\, - \frac{ \varepsilon \delta^{\sfrac12} \mr{\bar\lambda} }{\kappa\lambda} \sum_m [\eta_1\varphi_m] (\mu t) E_{m,j2} (x,t) \partial_{\xi_1}^{-1} \Pi_1 (\xi) \mr {\bar\partial}_j \bar \rho_0(x,t)     \\ 
      &\, + \frac{ \varepsilon \delta^{\sfrac12} \mr{\bar\lambda} }{\kappa\lambda} \sum_m [\eta_2\varphi_m] (\mu t) E_{m,j1} (x,t) \partial_{\xi_1}^{-1} \Pi_1 (\xi) \mr {\bar\partial}_j \bar \rho_0(x,t)
\end{align}

\subsubsection{Induction step for residual correctors}      \label{sss:D:expIndRsd}

In this induction step, we assume we have $\bar \rho_0$, $\bar \rho_1, \ldots,$ $\bar \rho_{n-1}$ and
\begin{align}
      \hat \rho_{n} (x,t,\xi) =&\, \sum_{a=0}^{n-1} \sum_{ [\balpha] = 1, \balpha \in \mathcal{I}_* }^{2(n-a)-1} \sum_{l=1}^{\hat I_{n-a}^{(\balpha)}} \hat h_{n-a}^{(\balpha,l)}(x,t) \hat \eta_{n-a}^{(\balpha,l)} (\mu t) \hat \chi_{n-a}^{(\balpha,l)} (\xi) \mr{ \bar D }^{\balpha} \bar\rho_a (x,t)      \label{e:D:6:expInd} \\ 
      =&\, \hat \rho_{n,0} + \hat \rho_{n,E}.         \label{e:D:7:expInd}
\end{align}
The goal is to deduce the structures of $\bar \rho_n$, i.e. the coefficient fields in the following equation
\begin{align}     \label{e:D:8:expInd}
      \bar L_0 \bar \rho_n =&\, \sum_{a=0}^{n-1} \sum_{ [\bbeta]=1, \bbeta \in \mathcal{I}_* }^{2(n-a)+1} \sum_{l=1}^{ \bar I_n^{(\bbeta)} } 
      \mr{\bar \divr} \Big( \bar h^{(\bbeta)}_{n-a} \odot \bar \eta^{(\bbeta)}_{n-a} (\mu \cdot) \mr {\bar D}^{\bbeta} \bar \rho_a \Big).
\end{align}
Similar to the proof of Lemma \ref{l:EXP1} and Lemma \ref{l:EXP2}, we have stationarity for coefficient fields and we omit the proof here.

Because of the stationarity property, it suffices to derive the coefficient fields in \eqref{e:D:8:expInd} for $a=0$. We take $\langle \cdot \rangle_\xi$ of \eqref{e:D:aspt:hrc3}, giving
\begin{align}     \label{e:D:10:expInd}
      \langle L_{00} \hat\rho_n \rangle_\xi + L_{01} \bar \rho_n + \langle L_{02} \hat\rho_n \rangle_\xi + \langle L_{03} \hat\rho_n \rangle_\xi + \langle L_{04} \hat\rho_n \rangle_\xi = 0.
\end{align}

Now we have
\begin{align} 
      \langle L_{00} \hat\rho_{n,0} \rangle_\xi 
            =&\, - \delta^{\sfrac12} \mr {\bar\lambda} \sum_{ [\balpha] = 1, \balpha \in \mathcal{I}_* }^{2n-1} \sum_{l=1}^{\hat I_{n}^{(\balpha)}} [ \eta_1 \hat \eta_{n}^{(\balpha,l)} ] (\mu \cdot) \langle \partial_{\xi_1} \Pi_1 \hat \chi_{n}^{(\balpha,l)} \rangle_\xi \mr{ \bar\partial }_2 \big( \hat h_{n}^{(\balpha,l)} \mr{ \bar D }^{\balpha} \bar\rho_0 \big)       \label{e:D:16:expInd} \\ 
            &\, + \delta^{\sfrac12} \mr {\bar\lambda} \sum_{ [\balpha] = 1, \balpha \in \mathcal{I}_* }^{2n-1} \sum_{l=1}^{\hat I_{n}^{(\balpha)}} [ \eta_2 \hat \eta_{n}^{(\balpha,l)} ] (\mu \cdot) \langle \partial_{\xi_2} \Pi_2 \hat \chi_{n}^{(\balpha,l)} \rangle_\xi \mr{ \bar\partial }_1 \big( \hat h_{n}^{(\balpha,l)} \mr{ \bar D }^{\balpha} \bar\rho_0 \big)       \label{e:D:18:expInd}
\end{align}

\begin{align}
      \langle L_{02} \hat\rho_{n,0} \rangle_\xi 
            =&\, - \varepsilon \frac{\delta^{\sfrac12}\mr{\bar\lambda}^2}{\lambda} \sum_{ [\balpha] = 1, \balpha \in \mathcal{I}_* }^{2n-1} \sum_{l=1}^{\hat I_{n}^{(\balpha)}} \sum_m [ \eta_1 \sigma_m \hat \eta_{n}^{(\balpha,l)} ] (\mu \cdot) \langle \Pi_1 \hat \chi_{n}^{(\balpha,l)} \rangle_\xi \mr{\bar \partial}_{i} \Omega_{m,ij} \mr{\bar \partial}_{j} \big( \hat h_{n}^{(\balpha,l)} \mr{ \bar D }^{\balpha} \bar\rho_0 \big)       \label{e:D:20:expInd} \\ 
            &\, + \varepsilon \frac{\delta^{\sfrac12}\mr{\bar\lambda}^2} {\lambda} \sum_{ [\balpha] = 1, \balpha \in \mathcal{I}_* }^{2n-1} \sum_{l=1}^{\hat I_{n}^{(\balpha)}} \sum_m [ \eta_2 \sigma_m \hat \eta_{n}^{(\balpha,l)} ] (\mu \cdot) \langle \Pi_2 \hat \chi_{n}^{(\balpha,l)} \rangle_\xi \mr{\bar \partial}_{i} \Omega_{m,ij} \mr{\bar \partial}_{j} \big( \hat h_{n}^{(\balpha,l)} \mr{ \bar D }^{\balpha} \bar\rho_0 \big)       \label{e:D:22:expInd}
\end{align}

\begin{align}
      \langle ( L_{03} + L_{04} ) \hat\rho_{n,0} \rangle_\xi
            =&\, - \varepsilon \delta^{\sfrac12} \mr {\bar\lambda} \sum_{ [\balpha] = 1, \balpha \in \mathcal{I}_* }^{2n-1} \sum_{l=1}^{\hat I_{n}^{(\balpha)}} \sum_m [ \eta_1 \varphi_m \hat \eta_{n}^{(\balpha,l)} ] (\mu \cdot) \langle \partial_{\xi_1} \Pi_1 \hat \chi_{n}^{(\balpha,l)} \rangle_\xi E_{m,j2} \mr{\bar \partial}_{j} \big( \hat h_{n}^{(\balpha,l)} \mr{ \bar D }^{\balpha} \bar\rho_0 \big)      \nonumber \\ 
            &\, + \varepsilon \delta^{\sfrac12} \mr {\bar\lambda} \sum_{ [\balpha] = 1, \balpha \in \mathcal{I}_* }^{2n-1} \sum_{l=1}^{\hat I_{n}^{(\balpha)}} \sum_m [ \eta_2 \varphi_m \hat \eta_{n}^{(\balpha,l)} ] (\mu \cdot) \langle \partial_{\xi_2} \Pi_2 \hat \chi_{n}^{(\balpha,l)} \rangle_\xi E_{m,j1} \mr{\bar \partial}_{j} \big( \hat h_{n}^{(\balpha,l)} \mr{ \bar D }^{\balpha} \bar\rho_0 \big)    \nonumber \\ 
            &\, + \varepsilon \delta^{\sfrac12} \mr {\bar\lambda} \sum_{ [\balpha] = 1, \balpha \in \mathcal{I}_* }^{2n-1} \sum_{l=1}^{\hat I_{n}^{(\balpha)}} \sum_m [ \eta_1 \varphi_m \hat \eta_{n}^{(\balpha,l)} ] (\mu \cdot) \langle \Pi_1 \partial_{\xi_1} \hat \chi_{n}^{(\balpha,l)} \rangle_\xi \mr{\bar \partial}_{j} E_{m,j2} \hat h_{n}^{(\balpha,l)} \mr{ \bar D }^{\balpha} \bar\rho_0    \nonumber \\ 
            &\, - \varepsilon \delta^{\sfrac12} \mr {\bar\lambda} \sum_{ [\balpha] = 1, \balpha \in \mathcal{I}_* }^{2n-1} \sum_{l=1}^{\hat I_{n}^{(\balpha)}} \sum_m [ \eta_2 \varphi_m \hat \eta_{n}^{(\balpha,l)} ] (\mu \cdot) \langle \partial_{\xi_2} \Pi_2 \hat \chi_{n}^{(\balpha,l)} \rangle_\xi \mr{\bar \partial}_{j} E_{m,j1} \hat h_{n}^{(\balpha,l)} \mr{ \bar D }^{\balpha} \bar\rho_0    \nonumber \\ 
            =&\, - \varepsilon \delta^{\sfrac12} \mr {\bar\lambda} \sum_{ [\balpha] = 1, \balpha \in \mathcal{I}_* }^{2n-1} \sum_{l=1}^{\hat I_{n}^{(\balpha)}} \sum_m [ \eta_1 \varphi_m \hat \eta_{n}^{(\balpha,l)} ] (\mu \cdot) \langle \partial_{\xi_1} \Pi_1 \hat \chi_{n}^{(\balpha,l)} \rangle_\xi \mr{\bar \partial}_{j} \big( E_{m,j2} \hat h_{n}^{(\balpha,l)} \mr{ \bar D }^{\balpha} \bar\rho_0 \big)      \label{e:D:24:expInd} \\ 
            &\, + \varepsilon \delta^{\sfrac12} \mr {\bar\lambda} \sum_{ [\balpha] = 1, \balpha \in \mathcal{I}_* }^{2n-1} \sum_{l=1}^{\hat I_{n}^{(\balpha)}} \sum_m [ \eta_2 \varphi_m \hat \eta_{n}^{(\balpha,l)} ] (\mu \cdot) \langle \partial_{\xi_2} \Pi_2 \hat \chi_{n}^{(\balpha,l)} \rangle_\xi \mr{\bar \partial}_{j} \big( E_{m,j1} \hat h_{n}^{(\balpha,l)} \mr{ \bar D }^{\balpha} \bar\rho_0 \big)    \label{e:D:26:expInd}
\end{align}

\begin{case}      \label{c:D:0:expInd}
For $1 \leq [\bbeta] \leq 2n-1$ and $\bbeta = \balpha$, the contribution from \eqref{e:D:16:expInd} and \eqref{e:D:18:expInd} gives at most $\bar I_n^{(\balpha)}$ terms with
\begin{align*}
      \bar h_{n+1}^{(\bbeta,\cdot)} =&\, \delta^{\sfrac12} \mr {\bar\lambda} 
      \begin{bmatrix}
            \langle \partial_{\xi_2} \Pi_2 \hat \chi_{n}^{(\balpha,l)} \rangle_\xi \hat h_{n}^{(\balpha,l)}       \\ 
            - \langle \partial_{\xi_1} \Pi_1 \hat \chi_{n}^{(\balpha,l)} \rangle_\xi \hat h_{n}^{(\balpha,l)}
      \end{bmatrix},    \quad
      \bar \eta_{n+1}^{(\bbeta,\cdot)} =
      \begin{bmatrix}
            \eta_2 \hat \eta_{n}^{(\balpha,l)}       \\ 
            \eta_1 \hat \eta_{n}^{(\balpha,l)}
      \end{bmatrix}.
\end{align*}
\end{case}

The contributions from \eqref{e:D:20:expInd}, \eqref{e:D:22:expInd}, \eqref{e:D:24:expInd} and \eqref{e:D:26:expInd} are analyzed by the same way as Case \ref{c:D:0:expInd}. Here we note that the constants in (\ref{e:D:20:expInd}-\ref{e:D:26:expInd}) are smaller than those in (\ref{e:D:16:expInd}-\ref{e:D:18:expInd}). This closes our derivation on $\bar \rho_{n}$.

\subsubsection{Induction step for spatial correctors}      \label{sss:D:expIndSpt}

We compute $\eqref{e:D:aspt:hrc3} - \langle \eqref{e:D:aspt:hrc3} \rangle_\xi$, giving
\begin{align}     \label{e:D:40:expInd}
      L_1 \hat \rho_{n+1} = \langle L_0 \hat \rho_n \rangle_\xi - L_0 \hat \rho_n. 
\end{align}

In this section, the goal is to derive the coefficients in 
\begin{align}     \label{e:D:42:expInd}
      \hat\rho_{n+1} (x,t,\xi) = \sum_{a=0}^n \sum_{[\bbeta] = 1, \bbeta \in \mathcal{I}_*}^{2(n-a)+1} \sum_{l=1}^{\hat I_{n+1,a}^{(\bbeta)}}
            \hat h_{n+1,a}^{(\bbeta,l)} (x,t) \hat \eta_{n+1,a}^{(\bbeta,l)} (\mu t) \hat \chi_{n+1,a}^{(\bbeta,l)} (\xi) \mr {\bar D}^{\bbeta} \bar \rho_a (x,t).
\end{align}

The stationarity follows in the same way as in Lemma \ref{l:EXP1} and Lemma \ref{l:EXP2}. We omit the proof of stationarity. Define
\begin{align}     \label{e:D:44:expInd}
      \hat I_{n-a+1}^{(\bbeta)} := \hat I_{n+1,a}^{(\bbeta)},           \quad 
      \hat h_{n-a+1}^{(\bbeta,l)} := \hat h_{n+1,a}^{(\bbeta,l)},       \quad 
      \hat \eta_{n-a+1}^{(\bbeta,l)} := \hat \eta_{n+1,a}^{(\bbeta,l)}, \quad 
      \hat \chi_{n-a+1}^{(\bbeta,l)} := \hat \chi_{n+1,a}^{(\bbeta,l)}.
\end{align}
Then we can write
\begin{align}     \label{e:D:46:expInd}
      \hat\rho_{n+1} (x,t,\xi) = \sum_{a=0}^n \sum_{[\bbeta] = 1, \bbeta \in \mathcal{I}_*}^{2(n-a)+1} \sum_{l=1}^{\hat I_{n-a+1}^{(\bbeta)}}
            \hat h_{n-a+1}^{(\bbeta,l)} (x,t) \hat \eta_{n-a+1}^{(\bbeta,l)} (\mu t) \hat \chi_{n-a+1}^{(\bbeta,l)} (\xi) \mr {\bar D}^{\bbeta} \bar \rho_a (x,t).
\end{align}

Note that from Remark \ref{r:D:hatcheckbar},
\begin{equation}  \label{e:D:48:expInd}
\begin{split}
      \langle L_0 \hat \rho_n \rangle_\xi - L_0 \hat \rho_n 
      = &\, \langle L_{00} \hat\rho_n \rangle_\xi - L_{00} \hat\rho_n - L_{00} \bar\rho_n - L_{01} \hat \rho_n + \sum_{k=2}^3 \langle L_{0k} \hat\rho_n \rangle_\xi - L_{0k} \hat\rho_n - L_{0k} \bar\rho_n     \\ 
      +&\, \langle L_{04} \hat\rho_n \rangle_\xi - L_{04} \hat\rho_n - \sum_{k=5}^9 L_{0k} \hat\rho_n.
\end{split}
\end{equation}

Now we compute the contributions from \eqref{e:D:48:expInd}. We have
\begin{align}
      \langle L_{00} \hat\rho_{n,0} \rangle_\xi - L_{00} \hat\rho_{n,0} 
            =&\, \delta^{\sfrac12} \mr {\bar\lambda} \sum_{ [\balpha] = 1, \balpha \in \mathcal{I}_* }^{B} \sum_{l=1}^{\hat I_{n}^{(\balpha)}} [ \eta_1 \hat \eta_{n}^{(\balpha,l)} ] (\mu \cdot) \big( \partial_{\xi_1} \Pi_1 \hat \chi_{n}^{(\balpha,l)} - \langle \partial_{\xi_1} \Pi_1 \hat \chi_{n}^{(\balpha,l)} \rangle_\xi \big) \mr{ \bar\partial }_2 \hat h_{n}^{(\balpha,l)} \mr{ \bar D }^{\balpha} \bar\rho_0   \label{e:D:50:expInd} \\ 
            +&\, \delta^{\sfrac12} \mr {\bar\lambda} \sum_{ [\balpha] = 1, \balpha \in \mathcal{I}_* }^{B} \sum_{l=1}^{\hat I_{n}^{(\balpha)}} [ \eta_1 \hat \eta_{n}^{(\balpha,l)} ] (\mu \cdot) \big( \partial_{\xi_1} \Pi_1 \hat \chi_{n}^{(\balpha,l)} - \langle \partial_{\xi_1} \Pi_1 \hat \chi_{n}^{(\balpha,l)} \rangle_\xi \big) \hat h_{n}^{(\balpha,l)} \mr{ \bar\partial }_2 \mr{ \bar D }^{\balpha} \bar\rho_0   \label{e:D:52:expInd} \\ 
      +&\, \text{symmetric terms}   \label{e:D:54:expInd}
\end{align}

\begin{align}
      -L_{01} \hat\rho_{n,0} =
            &\, - \mr{\bar\mu} \sum_{ [\balpha] = 1, \balpha \in \mathcal{I}_* }^{B} \sum_{l=1}^{\hat I_n^{(\balpha)}} \mr{\bar D}_t \hat h_{n}^{(\balpha,l)} \hat \eta_{n}^{(\balpha,l)} (\mu \cdot) \hat \chi_{n}^{(\balpha,l)} \mr{ \bar D }^{\balpha} \bar\rho_0    \label{e:D:56:expInd} \\ 
            &\, + \kappa \mr {\bar\lambda}^2 \sum_{ [\balpha] = 1, \balpha \in \mathcal{I}_* }^{B} \sum_{l=1}^{\hat I_n^{(\balpha)}} \mr {\bar \partial}_{ii} \hat h_{n}^{(\balpha,l)} \hat \eta_{n}^{(\balpha,l)} (\mu \cdot) \hat \chi_{n}^{(\balpha,l)} \mr{ \bar D }^{\balpha} \bar\rho_0   \label{e:D:58:expInd} \\   
            &\, - \mr{\bar\mu} \sum_{ [\balpha] = 1, \balpha \in \mathcal{I}_* }^{B} \sum_{l=1}^{\hat I_n^{(\balpha)}} \hat h_{n}^{(\balpha,l)} \hat \eta_{n}^{(\balpha,l)} (\mu \cdot) \hat \chi_{n}^{(\balpha,l)} \mr{\bar D}_t \mr{ \bar D }^{\balpha} \bar\rho_0    \label{e:D:60:expInd} \\ 
            &\, + \kappa \mr {\bar\lambda}^2 \sum_{ [\balpha] = 1, \balpha \in \mathcal{I}_* }^{B} \sum_{l=1}^{\hat I_n^{(\balpha)}} \hat h_{n}^{(\balpha,l)} \hat \eta_{n}^{(\balpha,l)} (\mu \cdot) \hat \chi_{n}^{(\balpha,l)} \mr {\bar \partial}_{ii} \mr{ \bar D }^{\balpha} \bar\rho_0   \label{e:D:62:expInd} \\ 
            &\, + \kappa \mr {\bar\lambda}^2 \sum_{ [\balpha] = 1, \balpha \in \mathcal{I}_* }^{B} \sum_{l=1}^{\hat I_n^{(\balpha)}} \mr {\bar \partial}_i \hat h_{n}^{(\balpha,l)} \hat \eta_{n}^{(\balpha,l)} (\mu \cdot) \hat \chi_{n}^{(\balpha,l)} \mr {\bar \partial}_i \mr{ \bar D }^{\balpha} \bar\rho_0        \label{e:D:64:expInd} \\ 
            &\, - \mu \sum_{ [\balpha] = 1, \balpha \in \mathcal{I}_* }^{B} \sum_{l=1}^{\hat I_n^{(\balpha)}} \hat h_{n}^{(\balpha,l)} \partial_\tau \hat \eta_{n}^{(\balpha,l)} (\mu \cdot) \hat \chi_{n}^{(\balpha,l)} \mr{ \bar D }^{\balpha} \bar\rho_0       \label{e:D:65:expInd}
\end{align}

\begin{align}
      &\, \langle L_{02} \hat\rho_{n,0} \rangle_\xi - L_{02} \hat\rho_{n,0}       \nonumber \\
      =&\, \varepsilon \frac{\delta^{\sfrac12}\mr{\bar\lambda}^2} {\lambda} \sum_{ [\balpha] = 1, \balpha \in \mathcal{I}_* }^{B} \sum_{l=1}^{\hat I_{n}^{(\balpha)}} \sum_m [ \eta_1 \sigma_m \hat \eta_{n}^{(\balpha,l)} ] (\mu \cdot) \big( \Pi_1 \hat \chi_{n}^{(\balpha,l)} - \langle \Pi_1 \hat \chi_{n}^{(\balpha,l)} \rangle_\xi \big) \mr{\bar \partial}_{i} \Omega_{m,ij} \mr{\bar \partial}_{j} \hat h_{n}^{(\balpha,l)} \mr{ \bar D }^{\balpha} \bar\rho_0       \label{e:D:66:expInd} \\ 
      +&\, \varepsilon \frac{\delta^{\sfrac12}\mr{\bar\lambda}^2} {\lambda} \sum_{ [\balpha] = 1, \balpha \in \mathcal{I}_* }^{B} \sum_{l=1}^{\hat I_{n}^{(\balpha)}} \sum_m [ \eta_1 \sigma_m \hat \eta_{n}^{(\balpha,l)} ] (\mu \cdot) \big( \Pi_1 \hat \chi_{n}^{(\balpha,l)} - \langle \Pi_1 \hat \chi_{n}^{(\balpha,l)} \rangle_\xi \big) \mr{\bar \partial}_{i} \Omega_{m,ij} \hat h_{n}^{(\balpha,l)} \mr{\bar \partial}_{j} \mr{ \bar D }^{\balpha} \bar\rho_0      \label{e:D:68:expInd} \\ 
      +&\, \text{symmetric terms}         \label{e:D:69:expInd}
\end{align}

\begin{align}
      &\, \langle L_{03} \hat\rho_{n,0} \rangle_\xi - L_{03} \hat\rho_{n,0}       \nonumber \\ 
      =&\, \varepsilon \delta^{\sfrac12} \mr {\bar\lambda} \sum_{ [\balpha] = 1, \balpha \in \mathcal{I}_* }^{B} \sum_{l=1}^{\hat I_{n}^{(\balpha)}} \sum_m [ \eta_1 \varphi_m \hat \eta_{n}^{(\balpha,l)} ] (\mu \cdot) \big( \partial_{\xi_1} \Pi_1 \hat \chi_{n}^{(\balpha,l)} - \langle \partial_{\xi_1} \Pi_1 \hat \chi_{n}^{(\balpha,l)} \rangle_\xi \big) E_{m,j2} \mr{\bar \partial}_{j} \hat h_{n}^{(\balpha,l)} \mr{ \bar D }^{\balpha} \bar\rho_0     \label{e:D:70:expInd} \\ 
      +&\, \varepsilon \delta^{\sfrac12} \mr {\bar\lambda} \sum_{ [\balpha] = 1, \balpha \in \mathcal{I}_* }^{B} \sum_{l=1}^{\hat I_{n}^{(\balpha)}} \sum_m [ \eta_1 \varphi_m \hat \eta_{n}^{(\balpha,l)} ] (\mu \cdot) \big( \partial_{\xi_1} \Pi_1 \hat \chi_{n}^{(\balpha,l)} - \langle \partial_{\xi_1} \Pi_1 \hat \chi_{n}^{(\balpha,l)} \rangle_\xi \big) E_{m,j2} \hat h_{n}^{(\balpha,l)} \mr{\bar \partial}_{j} \mr{ \bar D }^{\balpha} \bar\rho_0     \label{e:D:72:expInd}\\ 
      +&\, \text{symmetric terms}   \label{e:D:74:expInd}
\end{align}

\begin{align}
      - L_{06} \hat\rho_{n,0} 
            =&\, 2\kappa\lambda\mr {\bar\lambda} \sum_{ [\balpha] = 1, \balpha \in \mathcal{I}_* }^{B} \sum_{l=1}^{\hat I_{n-a}^{(\balpha)}} \mr{\bar \partial}_{i} \hat h_{n-a}^{(\balpha,l)} \hat \eta_{n-a}^{(\balpha,l)} (\mu \cdot) \partial_{\xi_i} \hat \chi_{n-a}^{(\balpha,l)} \mr{ \bar D }^{\balpha} \bar\rho_0     \label{e:D:76:expInd} \\ 
            +&\, 2\kappa\lambda\mr {\bar\lambda} \sum_{ [\balpha] = 1, \balpha \in \mathcal{I}_* }^{B} \sum_{l=1}^{\hat I_{n-a}^{(\balpha)}} \hat h_{n-a}^{(\balpha,l)} \hat \eta_{n-a}^{(\balpha,l)} (\mu \cdot) \partial_{\xi_i} \hat \chi_{n-a}^{(\balpha,l)} \mr{\bar \partial}_{i} \mr{ \bar D }^{\balpha} \bar\rho_0     \label{e:D:78:expInd}
\end{align}

\begin{align}
      - L_{08} \hat\rho_{n,0} 
            =&\, \varepsilon\kappa\lambda^2 \sum_{ [\balpha] = 1, \balpha \in \mathcal{I}_* }^{B} \sum_{l=1}^{\hat I_{n-a}^{(\balpha)}} \mr{\bar \partial}_{i} \hat h_{n-a}^{(\balpha,l)} \hat \eta_{n-a}^{(\balpha,l)} (\mu \cdot) \partial_{\xi_i} \hat \chi_{n-a}^{(\balpha,l)} \mr{ \bar D }^{\balpha} \bar\rho_0     \label{e:D:80:expInd} \\ 
            +&\, \varepsilon\kappa\lambda^2 \sum_{ [\balpha] = 1, \balpha \in \mathcal{I}_* }^{B} \sum_{l=1}^{\hat I_{n-a}^{(\balpha)}} \hat h_{n-a}^{(\balpha,l)} \hat \eta_{n-a}^{(\balpha,l)} (\mu \cdot) \partial_{\xi_i} \hat \chi_{n-a}^{(\balpha,l)} \mr{\bar \partial}_{i} \mr{ \bar D }^{\balpha} \bar\rho_0.     \label{e:D:82:expInd}
\end{align}

For $a=0$ and a fixed $\bbeta$, we consider all contributions above and renumber them with index $1 \leq l \leq \hat I_{n+1}^{(\bbeta)}$, in the form of \eqref{e:D:46:expInd}.

\begin{case}      \label{c:D:12:expInd}
If $2 \leq [\bbeta] \leq 2n$ and $\bbeta = 2\balpha$, the contribution from \eqref{e:D:52:expInd} gives at most $\hat I_n^{(\balpha)}$ correctors with
\begin{align*}
      \hat h_{n+1}^{(\bbeta,\cdot)} =&\, - \frac{\delta^{\sfrac12} \mr {\bar\lambda}} {\kappa\lambda} \hat h_n^{(\balpha,l)},   \quad 
      \hat \eta_{n+1}^{(\bbeta,\cdot)} = \eta_1 \hat \eta_n^{(\balpha,l)},    \\ 
      \hat \chi_{n+1}^{(\bbeta,\cdot)} =&\, \Delta_\xi^{-1} \big( \partial_{\xi_1} \Pi_1 \hat \chi_{n}^{(\balpha,l)} - \langle \partial_{\xi_1} \Pi_1 \hat \chi_{n}^{(\balpha,l)} \rangle_\xi \big).
\end{align*}
The contributions from \eqref{e:D:50:expInd} and \eqref{e:D:54:expInd} are similar to above.
\end{case}

\begin{case}      \label{c:D:14:expInd}
If $1 \leq [\bbeta] \leq 2n-1$ and $\bbeta = \balpha$, the contribution from \eqref{e:D:56:expInd} gives at most $\hat I_n^{(\balpha)}$ correctors with
\begin{align*}
      \hat h_{n+1}^{(\bbeta,\cdot)} =&\, \frac{\mr{\bar \mu}}{\kappa\lambda} \mr {\bar D}_t \hat h_n^{(\balpha,l)},   \quad 
      \hat \eta_{n+1}^{(\bbeta,\cdot)} = \hat \eta_n^{(\balpha,l)},    \\ 
      \hat \chi_{n+1}^{(\bbeta,\cdot)} =&\, \Delta_\xi^{-1} \hat \chi_{n}^{(\balpha,l)} .
\end{align*}
\end{case}

\begin{case}      \label{c:D:16:expInd}
If $1 \leq [\bbeta] \leq 2n-1$ and $\bbeta = \balpha$, the contribution from \eqref{e:D:58:expInd} gives at most $\hat I_n^{(\balpha)}$ correctors with
\begin{align*}
      \hat h_{n+1}^{(\bbeta,\cdot)} =&\, - \frac{\mr {\bar\lambda}^2}{\lambda} \mr{\bar\partial}_{ii} \hat h_n^{(\balpha,l)},   \quad 
      \hat \eta_{n+1}^{(\bbeta,\cdot)} = \hat \eta_n^{(\balpha,l)},    \\ 
      \hat \chi_{n+1}^{(\bbeta,\cdot)} =&\, \Delta_\xi^{-1} \hat \chi_{n}^{(\balpha,l)} .
\end{align*}
\end{case}

\begin{case}      \label{c:D:18:expInd}
If $1 \leq [\bbeta] \leq 2n-1$ and $\bbeta = \balpha$, the contribution from \eqref{e:D:65:expInd} gives at most $\hat I_n^{(\balpha)}$ correctors with
\begin{align*}
      \hat h_{n+1}^{(\bbeta,\cdot)} =&\, \frac{\mu}{\kappa\lambda} \hat h_n^{(\balpha,l)},   \quad 
      \hat \eta_{n+1}^{(\bbeta,\cdot)} = \partial_\tau \hat \eta_n^{(\balpha,l)},    \\ 
      \hat \chi_{n+1}^{(\bbeta,\cdot)} =&\, \Delta_\xi^{-1} \hat \chi_{n}^{(\balpha,l)} .
\end{align*}
\end{case}

\begin{case}      \label{c:D:20:expInd}
If $1 \leq [\bbeta] \leq 2n-1$ and $\bbeta = \balpha$, the contribution from \eqref{e:D:70:expInd} gives at most $\hat I_n^{(\balpha)}$ correctors with
\begin{align*}
      \hat h_{n+1}^{(\bbeta,\cdot)} =&\, \varepsilon \frac{\delta^{\sfrac12} \mr{\bar\lambda}}{\kappa\lambda} E_{m,j2} \mr{\bar \partial}_{j} \hat h_{n}^{(\balpha,l)},   \quad 
      \hat \eta_{n+1}^{(\bbeta,\cdot)} = \eta_1 \varphi_m \hat \eta_{n}^{(\balpha,l)},    \\ 
      \hat \chi_{n+1}^{(\bbeta,\cdot)} =&\, \Delta_\xi^{-1} \big( \partial_{\xi_1} \Pi_1 \hat \chi_{n}^{(\balpha,l)} - \langle \partial_{\xi_1} \Pi_1 \hat \chi_{n}^{(\balpha,l)} \rangle_\xi \big) .
\end{align*}
\end{case}

\begin{case}
If $1 \leq [\bbeta] \leq 2n-1$ and $\bbeta = \balpha$, the contribution from \eqref{e:D:76:expInd} gives at most $\hat I_n^{(\balpha)}$ correctors with
\begin{align*}
      \hat h_{n+1}^{(\bbeta,\cdot)} =&\, 2 \mr {\bar\lambda} \mr{\bar \partial}_{i} \hat h_{n}^{(\balpha,l)},   \quad 
      \hat \eta_{n+1}^{(\bbeta,\cdot)} = \hat \eta_{n}^{(\balpha,l)},    \\ 
      \hat \chi_{n+1}^{(\bbeta,\cdot)} =&\, \Delta_\xi^{-1} \partial_{\xi_i} \hat \chi_{n}^{(\balpha,l)} .
\end{align*}
\end{case}

\begin{case}
If $1 \leq [\bbeta] \leq 2n-1$ and $\bbeta = \balpha$, the contribution from \eqref{e:D:80:expInd} gives at most $\hat I_n^{(\balpha)}$ correctors with
\begin{align*}
      \hat h_{n+1}^{(\bbeta,\cdot)} =&\, \varepsilon \lambda \mr{\bar \partial}_i \hat h_{n}^{(\balpha,l)},   \quad 
      \hat \eta_{n+1}^{(\bbeta,\cdot)} = \hat \eta_{n}^{(\balpha,l)},    \\ 
      \hat \chi_{n+1}^{(\bbeta,\cdot)} =&\, \Delta_\xi^{-1} \big( \partial_{\xi_1} \Pi_1 \hat \chi_{n}^{(\balpha,l)} - \langle \partial_{\xi_1} \Pi_1 \hat \chi_{n}^{(\balpha,l)} \rangle_\xi \big) .
\end{align*}
\end{case}

This closes our derivation on $\hat \rho_{n+1}$.

\newpage

\section{Proof of the main result}      \label{s:mainRes}

\subsection{Choice of parameters}
Let us start with $\alpha_0,\beta_0$ as specified in the Theorem, i.e.~$0<\beta_0<\frac{1}{3}$ and $\alpha_0>0$ such that $2\alpha_0 + \beta_0 < 1$, and $s >\frac{1+\beta_0}{2}$, $\alpha' \in (\alpha_0, 1)$. Choose $\beta>\beta_0$ so that $\beta<\frac{1}{3}$, $2\alpha_0 + \beta < 1$ and $\frac{1+\beta}{2}<s$. We then choose the parameters $b$, $\gamma_I$, $\gamma_R$, $\gamma_S$, $\gamma$, $\lambda_0$, $N_*$, $N$ and $Q$ according to Section \ref{ss:parameters}. Our aim is then to construct a divergence-free velocity field $u \in C^{\beta_0}(\T^2\times [0,1])$ such that \eqref{e:2:t:mainHld} and \eqref{e:4:t:mainHld} with $\alpha=\alpha_0$ holds.

\subsection{The velocity field}\label{ss:main_pf_u}

 We construct the vector field $u$ as in Section \ref{ss:vectorfield}. In particular, estimate \eqref{e:streamEs_rig2_nobar_psi} (with $q+1$ replaced with $q$) yields 
\[
  \begin{aligned}
      \|  \nabla^\perp \psi_{q} \|_{\infty}  
      &\lesssim \delta^{\sfrac12}_{q} = \lambda^{-\beta}_{q} \\
      \|  \nabla \nabla^\perp   \psi_{q} \|_{\infty}  
      &\lesssim \delta^{\sfrac12}_{q} \lambda_{q} = \lambda^{1-\beta}_{q} \\
      \| D_{q} \nabla^\perp \psi_{q} \|_{\infty} 
      &\lesssim \delta_{q} \lambda_{q} =  \lambda^{1-2\beta}_{q} 
 \end{aligned}
\]
Since \[
\partial_t \nabla^\perp \psi_q = D_q \nabla^\perp \psi_q - u_q \cdot \nabla \nabla^\perp \psi_q,
\] the above estimates give also
\[
 \| \partial_t \nabla^\perp \psi_q \|_{\infty} \lesssim \lambda^{1-2\beta}_{q} +  \| u_q\|_{\infty}  \lambda^{1-\beta}_{q} \lesssim \lambda^{1-\beta}_{q},
\]
with the last estimate valid in view of \eqref{e:streamEs_rig2_nobar}.
Using half of the gap  $g:=\beta- \beta_0$ and the base of our superexponential sequence to consume the constants of $\lesssim$, we obtain
\[
  \begin{aligned}
      \|  \nabla^\perp \psi_{q} \|_{\infty}  
      &\le \lambda^{-\sfrac{g}{2}}_{q} \lambda^{-\beta_0}_{q} \\
      \|  \nabla \nabla^\perp   \psi_{q} \|_{\infty}  +  \| \partial_t \nabla^\perp \psi_q \|_{\infty}
      &\le \lambda^{-\sfrac{g}{2}}_{q}  \lambda^{1-\beta_0}_{q}.
 \end{aligned}
\]
Consequently by interpolation
\[ \|  \nabla^\perp \psi_{q} \|_{C_{xt}^{\beta_0}}  \le \lambda^{-\sfrac{g}{2}}_{q} 
\]
Since $u_m = \sum_{q=0}^m\nabla^\perp \psi_{q}$ and our $\{\lambda_{q}\}_q$ is at least geometric, this gives for any $n \ge m \ge q_0$ 
\[\|  u_m - u_n \|_{C_{xt}^{\beta_0}}  
      \le \lambda^{-\sfrac{g}{4}}_{m},
\]
so $\{u_{q}\}_{q \ge q_0}$ is a Cauchy sequence in ${C_{xt}^{\beta_0}}$.

\subsection{Proof of \eqref{e:4:t:mainHld} in Theorem \ref{t:mainHld}}
For the sequence $\{\kappa_q\}_q$ set in Section \ref{ss:parameters}, our goal is to estimate the solution to the following equation
\begin{equation}\label{e:5:p:mainHld}
\begin{aligned}
      \partial_t \rho_n + u \cdot \nabla \rho_n - \kappa_n \Delta \rho_n =&\, 0,       \\ 
      \rho_n (\cdot, 0) =&\, \rho_{\ini}
\end{aligned}
\end{equation}
for any initial datum $\rho_{\ini} \in H^{1+\alpha'} (\T^2)$ and for an arbitrary $n \geq 1$. In the rest of the proof, we fix an arbitrary $n \geq 1$. For this fixed $n$, we will decompose the problem with respect to parameter $k \in \N$.
In this proof, we will use the subscript $q$ as an additional parameter needed for iterations. Thus the final result concerns \eqref{e:5:p:mainHld} with $n$, which is replaced in formulation of Theorem \ref{t:mainHld} by $q$.

\begin{step}[Decomposition of the initial datum]      \label{st:2:mainHld}
Consider the Littlewood-Paley decomposition of the initial datum
\begin{align}
      \rho_{\ini} = \sum_{j \geq -1} \Delta_j \rho_{\ini}.      \label{e:10:p:mainHld}
\end{align}
We define
\begin{align}
      \rho_{\ini,k} :=&\, \sum_{j \in \Lambda_k} \Delta_j \rho_{\ini}.      \label{e:11:p:mainHld} 
\end{align}
Here, $\Lambda_k = \{ j: \lambda_{k-1} < 2^j \leq \lambda_k \}$ for $k \geq 1$ and $\Lambda_0 = \{ j: 2^j \leq \lambda_0 \}$. From Lemma \ref{l:a:littlewood_paley}, we have the following standard estimate for any integer,
\begin{align}
      \| \mr\nabla_k^p \rho_{\ini,k} \|_{2} 
            \lesssim_p&\, \bigg( \frac{\lambda_k}{\mr\lambda_k} \bigg)^p \sum_{j \in \Lambda_k} \| \Delta_j \rho_{\ini} \|_{2},   \label{e:12:p:mainHld} 
\end{align}
then we have, using the quantities defined in Section \ref{ss:char_datum}
\begin{align}
      \| \rho_{\ini,k} \|_{2} + | \rho_{\ini,k} |_{\mathfrak{P}(k)} + \sum_{j=k}^\infty \lambda_j^{-\sfrac{\gamma}{2}} | \rho_{\ini,k} |_{\fS(j)}
            \lesssim&\, \sum_{j \in \Lambda_k} \| \Delta_j \rho_{\ini} \|_{2}.      \label{e:14:p:mainHld}
\end{align}

\end{step}

\begin{step}[Iteration for large-scale datum ($k \leq n$)]      \label{st:3:mainHld}

Fix $k \leq n$. From direct energy estimates in Proposition \ref{p:unconEnergyP}, we have that, for any $\tilde \rho_{\ini} \in C^\infty(\T^2)$
\begin{align}
      \resN{ \tilde \rho_{\ini} }_{u_k, \kappa_k}
            \lesssim&\, \mr\lambda_k^{1-2\gamma} \bigg( | \tilde\rho_{\ini} |_{\mathfrak{D}(u_k,\kappa_k)} 
                  + | \tilde\rho_{\ini} |_{\mathfrak{P}(k)} \bigg),       \label{e:17:p:mainHld} \\ 
      \resN{ \tilde \rho_{\ini} }_{u_k, \kappa_k}^{(H,\alpha_0)}
            \lesssim&\,   \mr\lambda_k^{(1-2\gamma)(1+\alpha_0)} 
                \bigg( | \tilde\rho_{\ini} |_{\mathfrak{D}(u_k,\kappa_k)} 
                  + | \tilde\rho_{\ini} |_{\mathfrak{P}(k)} \bigg).       \label{e:18:p:mainHld}
\end{align}
Then we use \eqref{e:2:unw_IniDatNorm} and absorb the constant in $\lesssim$ by a small negative power of $\lambda_k$ to deduce the following: for any $\tilde \rho_{\ini} \in C^\infty(\T^2)$,
\begin{align}
      \resN{ \tilde \rho_{\ini} }_{u_k, \kappa_k}
            \leq&\, \mr\lambda_k^{1-\gamma} \bigg( | \tilde\rho_{\ini} |_{\mathfrak{D}(u_k,\kappa_k)} 
                  + \lambda_k^{-\sfrac{\gamma}{2}} | \tilde\rho_{\ini} |_{\fS(k)} \bigg),       \label{e:20:p:mainHld} \\ 
      \resN{ \tilde \rho_{\ini} }_{u_k, \kappa_k}^{(H,\alpha_0)}
            \leq&\,   \mr\lambda_k^{(1-\gamma)(1+\alpha_0)} 
                  \cdot \bigg( 1 - \frac{4}{k} \bigg)\bigg( | \tilde\rho_{\ini} |_{\mathfrak{D}(u_k,\kappa_k)} 
                  + \lambda_k^{-\sfrac{\gamma}{2}} | \tilde\rho_{\ini} |_{\fS(k)} \bigg).       \label{e:22:p:mainHld}
\end{align}
Now we apply homogenization estimates of Theorem \ref{t:homRslv} inductively over $q$ from $q=k$ (so current $k=q_*$ of Theorem \ref{t:homRslv}) up to $q=n$ to deduce the following: for any $q$ with $k \leq q \leq n$, we have
\begin{align}
      \resN{ \tilde \rho_{\ini} }_{u_q, \kappa_q}
            \leq&\, \mr\lambda_k^{1-\gamma} \bigg( | \tilde\rho_{\ini} |_{\mathfrak{D}(u_q,\kappa_q)} 
                  + \sum_{j=k}^{q} \lambda_j^{-\sfrac{\gamma}{2}} | \tilde\rho_{\ini} |_{\fS(j)} \bigg),    \label{e:24:p:mainHld} \\ 
      \resN{ \tilde \rho_{\ini} }_{u_q, \kappa_q}^{(H,\alpha_0)}
            \leq&\,   \mr\lambda_k^{(1-\gamma)(1+\alpha_0)} 
                  \bigg( 1 - \frac{4}{q} \bigg)\bigg( | \tilde\rho_{\ini} |_{\mathfrak{D}(u_q,\kappa_q)} 
                  + \sum_{j=k}^{q} \lambda_j^{-\sfrac{\gamma}{2}} | \tilde\rho_{\ini} |_{\fS(j)} \bigg).    \label{e:26:p:mainHld}
\end{align}
Here, we use the linearity in Remark \ref{r:linearityHomRslv}.

Using \eqref{e:24:p:mainHld} and \eqref{e:26:p:mainHld} with $q=n$, we apply homogenization estimates of Theorem \ref{t:D:homRslv} inductively to deduce the following: for any $q \geq n$,
\begin{align}
      \resN{ \tilde \rho_{\ini} }_{u_q, \kappa_n}
            \leq&\, \mr\lambda_k^{1-\gamma} \bigg( | \tilde\rho_{\ini} |_{\mathfrak{D}(u_q,\kappa_n)} 
                  + \sum_{j=k}^{q} \lambda_j^{-\sfrac{\gamma}{2}} | \tilde\rho_{\ini} |_{\fS(j)} \bigg),    \label{e:28:p:mainHld} \\ 
      \resN{ \tilde \rho_{\ini} }_{u_q, \kappa_n}^{(H,\alpha_0)}
            \leq&\, \mr\lambda_k^{(1-\gamma)(1+\alpha_0)} 
                   \bigg( 1 - \frac{4}{q} \bigg)\bigg( | \tilde\rho_{\ini} |_{\mathfrak{D}(u_q,\kappa_n)} 
                  + \sum_{j=k}^{q} \lambda_j^{-\sfrac{\gamma}{2}} | \tilde\rho_{\ini} |_{\fS(j)} \bigg).    \label{e:30:p:mainHld}
\end{align}
Next, we apply the above estimates for $\tilde \rho_{\ini} = \rho_{\ini,k}$ to the solution to
\begin{align}
      \partial_t \bar\rho_{k,q} + u_q \cdot \nabla \bar\rho_{k,q} - \kappa_n \Delta \bar\rho_{k,q} =&\, 0,      \label{e:34:p:mainHld} \\ 
      \bar\rho_{k,q} (\cdot, 0) =&\, \rho_{\ini,k}.   \label{e:36:p:mainHld}
\end{align}
This yields, using definition of  $\resN{ \cdot }_{u_q, \kappa_n}$ (preceding Theorem \ref{t:homDissip}) and \eqref{e:trivialDissipEst}, the following estimate uniformly in $q$ for any $q \geq n$,
\begin{align}
      \| \bar\rho_{k,q} \|_{L_t^\infty C_x^{\alpha_0}} 
            \leq&\, \mr\lambda_k^{(1-\gamma)(1+\alpha_0)} \sum_{j \in \Lambda_k} \| \Delta_j \rho_{\ini} \|_{2}.       \label{e:38:p:mainHld}
\end{align}
Similarly as in Section \ref{ss:main_pf_u}, Lemma \ref{l:streamEst_rig_nobar} gives $u_q \rightarrow u$ in $C^{\sfrac{\gamma}{2}}_t C^\gamma_x$. In particular, we have  that  $\{u_q\}_q$ is uniformly in $C^{\sfrac{\gamma}{2}}_t C^\gamma_x$ (i.e.\ with q-independent bound). Therefore, since $\kappa_n$ of \eqref{e:38:p:mainHld} is a fixed positive number, from standard Schauder theory we know that $\{\bar\rho_{k,q}\}_q$ is $q$-uniformly bounded in $C_t^{\sfrac{\gamma}{2}} C_x^{1,\gamma}$. Thus, sending $q \to \infty$ we see that $\{\bar\rho_{k,q}\}_q$, up to a subsequence, converges strongly in $C_t^0 C_x^{\alpha_0}$ to $\rho_{k,n}$ uniquely solving the following equation
\begin{align}
      \partial_t \rho_{k,n} + u \cdot \nabla \rho_{k,n} - \kappa_n \Delta \rho_{k,n} =&\, 0,    \label{e:58:p:mainHld} \\ 
      \rho_{k,n} (\cdot, 0) =&\, \rho_{\ini,k}.       \label{e:60:p:mainHld}
\end{align}
Knowing the strong convergence, we have from \eqref{e:38:p:mainHld} the following estimate for the limiting $\rho_{k,n}$
\begin{align}     \label{e:62:p:mainHld}
      \| \rho_{k,n} \|_{L_t^\infty C_x^{\alpha_0}} 
            \leq&\, \mr\lambda_k^{(1-\gamma)(1+\alpha_0)} \sum_{j \in \Lambda_k} \| \Delta_j \rho_{\ini} \|_{2},
\end{align}
which is uniform in $k$ for any $k \leq n$.

Let us remark that  we used Schauder theory only to obtain the qualitative convergence, whilst the estimate \eqref{e:62:p:mainHld} follows from the uniform estimate \eqref{e:38:p:mainHld} and the qualitative convergence. Obtaining \eqref{e:62:p:mainHld} directly from Schauder estimates is unlikely, since they have typically bad constants dependence on $\kappa_n$ and $u$.

\end{step}

\begin{step}[Iteration for small-scale datum ($k \geq n$)]      \label{st:4:mainHld}

Fix $k \geq n$. We proceed as in last step. Instead of applying Theorem \ref{t:homRslv} and Theorem \ref{t:D:homRslv}, in this step we only apply Theorem \ref{t:D:homRslv}. We deduce the estimate \eqref{e:62:p:mainHld} uniform in $k$ for any $k \geq n$.

\end{step}

\begin{step}[Summing up all frequencies]  \label{st:5:mainHld}
From Step \ref{st:3:mainHld} and Step \ref{st:4:mainHld}, for fixed $n$, we have the estimate \eqref{e:62:p:mainHld} for $\rho_{k,n}$ uniform in $k$ for any $k \geq 0$. From the definition of $\rho_n$, we have
\[
\begin{aligned}
      \| \rho_n \|_{L_t^\infty C_x^{\alpha_0}}
            \leq&\, \sum_{k \geq 0} \| \rho_{k,n} \|_{L_t^\infty C_x^{\alpha_0}}       
            \\ 
            \leq&\, \sum_{k \geq 0} \mr\lambda_k^{(1-\gamma)(1+\alpha_0)} 
                  \sum_{j \in \Lambda_k} \| \Delta_j \rho_{\ini} \|_{2}       
                  \\ 
            \leq&\, \sum_{j \geq -1} 2^{jb(1-\gamma)(1+\alpha_0)} \| \Delta_j \rho_{\ini} \|_{2}       
            \\ 
            \lesssim&\, \| \rho_{\ini} \|_{H^{(1-\gamma)(1+\alpha_0)b}}       
\end{aligned}
\]
\end{step}
Inspecting the preceding computations, we see that all of the functions are pointwisely defined in time, so we can replace $L_t^\infty$ norm appearing in the preceding computations with $C_t^0$ norm.
This concludes the proof of \eqref{e:4:t:mainHld} in Theorem \ref{t:mainHld}. 

\subsection{Proof of \eqref{e:2:t:mainHld} in Theorem \ref{t:mainHld}}
Our goal now is the justify the remaining \eqref{e:2:t:mainHld}.
Recall that the velocity field $u$ was constructed in Section \ref{ss:vectorfield}. For sufficiently large $n$, our goal now is to study the dissipation
\begin{align}     \label{e:4:p:mainEne}
      \kappa_n \iint | \nabla \rho_n |^2 dxdt
\end{align}
for the solution to
\begin{equation}  \label{e:6:p:mainEne}
\begin{split}
      \partial_t \rho_n + u \cdot \nabla \rho_n - \kappa_n \Delta \rho_n =&\, 0,       \\ 
      \rho_n (\cdot, 0) =&\, \rho_{\ini}.
\end{split}
\end{equation}

\begin{step}[Mollification of the initial datum]
Define the ratio given by the initial datum
\begin{equation*}
    \ell := \frac{\|\rho_{\ini}\|_2} {\|\nabla^s \rho_{\ini}\|_2},
\end{equation*}
and mollify the initial datum at length scale $r\ell$:
\begin{align}     \label{e:12:p:mainEne}
      \tilde \rho_{\ini} = \rho_{\ini} * m_{r\ell},
\end{align}
where $m$ is a standard symmetric mollifier (see e.g.\ Appendix) and a suitably small $0<r<1$ is to be fixed. Then, for all $0 \leq n \leq 12Q^3$, we have
\begin{align}
      \| \tilde\rho_{\ini} - \rho_{\ini} \|_{2} 
            &\leq C (r\ell)^{s} \| \nabla^s \rho_{\ini} \|_{2} = C r^s \ell^{s-1} \|\rho_{\ini}\|_{2} \leq \frac{1}{2} \|\rho_{\ini}\|_{2},      \label{e:14:p:mainEne} \\ 
      \| \nabla^{n+1} \tilde\rho_{\ini} \|_{2} 
            &\leq C (r\ell)^{-(n+1-s)} \| \nabla^s \rho_{\ini} \|_{2} 
            \leq C (r\ell)^{-(n+1)} r^s \ell^{s-1} \| \rho_{\ini} \|_{2}.        \label{e:16:p:mainEne}
\end{align}
In (\ref{e:14:p:mainEne}), we impose the first condition on $r$:
\begin{align}     \label{e:18:p:mainEne}
      2 C r^s \ell^{s-1} \leq 1.
\end{align}

\end{step}

\begin{step}[The initial step of homogenization: choosing $q_0$, $\ell$ and $r$]
For $q_0$ to be fixed, we consider the solution to
\begin{align*}
      \partial_t \tilde\rho_q + u_q \cdot \nabla \tilde\rho_q - \kappa_q \Delta \tilde\rho_q =&\, 0,       \\ 
      \tilde\rho_q (\cdot, 0) =&\, \tilde\rho_{\ini}
\end{align*}
for any $q \geq q_0$. From Poincar\'e inequality, there exists absolute constant $c_0>0$ such that
\begin{align}  \label{e:24:p:mainEne}
	\frac{d}{dt} \| \tilde\rho_{q}(\cdot,t) \|^2_{2} dt 
		\leq -2 \kappa_{q} \| \nabla \tilde\rho_{q}(\cdot,t) \|^2_{2}
		\leq - c_0 \kappa_q\| \tilde\rho_{q}(\cdot,t) \|^2_{2}.
\end{align}
Integrating \eqref{e:24:p:mainEne} from $0$ to $T$, we have that, for any $q \geq q_0$,
\begin{equation}  \label{e:26:p:mainEne}
\begin{split}
	\kappa_{q} \iint | \nabla \tilde\rho_{q} |^2 dxdt 
		\geq& \frac{1}{2} \left( 1-\exp( -c_0\kappa_{q}T ) \right) \| \tilde\rho_{\ini} \|_{2}^2 
		\geq 16c_1 \kappa_{q} T \|\tilde\rho_{\ini}\|_{2}^2 \geq 4c_1 \kappa_qT\|\rho_{\ini}\|^2_{2}.
\end{split}
\end{equation}
Here $c_0,c_1>0$ are absolute constants.

For a small constant $c>0$ to be fixed, we require
\begin{align}
      (r\ell)^{-1} \leq&\, c \mr\lambda_{q_0}^{1-b\gamma},     \label{e:28:p:mainEne} \\
      r^s \ell^{s-1} \leq&\, c  (\kappa_{q_0} T)^{\sfrac{1}{2}}. \label{e:30:p:mainEne}
\end{align}
Note that \eqref{e:18:p:mainEne} follows from \eqref{e:30:p:mainEne} for suitably small $c$. And an equivalent set of conditions for \eqref{e:28:p:mainEne}-\eqref{e:30:p:mainEne} is
\begin{align}
      c^{-s} \mr\lambda_{q_0}^{(b\gamma-1)s} \ell^{-1}
            \leq r^s \ell^{s-1}
            \leq c  (\kappa_{q_0} T)^{\sfrac{1}{2}}.    \label{e:32:p:mainEne}
\end{align}
We use \eqref{e:32:p:mainEne} to determine $r$ and $q_0$. First, notice that that our choices of section \ref{ss:parameters} yield 
\begin{align}
      c^{-1}  (\kappa_{q} T)^{-\sfrac{1}{2}} 
            < c^{s} \mr\lambda_{q}^{(1-b\gamma)s} \ell
            \quad \text{for sufficiently large } q,    \label{e:33:p:mainEne}
\end{align}
provided $s>\sfrac{(1 + \beta)}{2}$, which is our assumption.
Then we choose $q_0$ sufficiently large so that \eqref{e:33:p:mainEne} holds for $q=q_0$. Increasing $q_0$ if needed we ensure that $\ell^{s-1} > c  (\kappa_{q_0} T)^{\sfrac{1}{2}}$ holds; consequently we can choose $r$ such that \eqref{e:32:p:mainEne} holds.

As a result of \eqref{e:28:p:mainEne} and \eqref{e:30:p:mainEne}, we have the following relations for any $q \geq q_0$ and $p \geq 1$,
\begin{align}
      (r\ell)^{-1} \leq\, c \mr\lambda_{q}^{1-b\gamma},  \quad
      (r\ell)^{-p} r^s \ell^{s-1} \leq\, c ( \kappa_{q} T )^{\sfrac{1}{2}} \mr\lambda_{q}^{(1-b\gamma)p}.         \label{e:34:p:mainEne}
\end{align}
\end{step}

\begin{step}[Homogenization in inertial range]
First, we estimate the quantity $\mathfrak{P}(q)$ for any $q \geq q_0$. Using \eqref{e:16:p:mainEne}, \eqref{e:26:p:mainEne} and \eqref{e:32:p:mainEne}, we have for any $p \geq 1$,
\begin{align}     \label{e:38:p:mainEne}
      \| \nabla^p \tilde\rho_{\ini} \|_{2} 
            \lesssim_p \mr\lambda_{q_0}^{(1-b\gamma)p} 
            \bigg( \kappa_{q} \iint | \nabla \tilde\rho_{q} |^2 dxdt \bigg)^{\frac{1}{2}},
\end{align}
from which we deduce, for any $q \geq q_0$,
\begin{align}     \label{e:40:p:mainEne}
      | \tilde\rho_{\ini} |_{\mathfrak{P}(q)} 
            + \sum_{j=q_0}^q \lambda_j^{-\sfrac{\gamma}{2}} | \tilde\rho_{\ini} |_{\fS(j)}
            \leq&\, | \tilde\rho_{\ini} |_{\mathfrak{D}(u_{q},\kappa_{q})}.
\end{align}
Here, we impose the first condition on $c$ in \eqref{e:28:p:mainEne}-\eqref{e:30:p:mainEne} to ensure \eqref{e:40:p:mainEne} holds, i.e. absorbing the constant when deriving \eqref{e:40:p:mainEne}.

Using \eqref{e:40:p:mainEne}, we apply Theorem \ref{t:homDissip} inductively to deduce the following: for any $q \geq q_0$,
\begin{align}     \label{e:44:p:mainEne}
      \kappa_{q_0} \iint | \nabla \tilde\rho_{q_0} |^2 dxdt \prod_{j=q_0+1}^{q} ( 1 - \lambda_j^{-\gamma} )
            \leq\, \kappa_q \iint | \nabla \tilde\rho_q |^2 dxdt.
\end{align}

\end{step}

\begin{step}[Homogenization in dissipative range]
We consider the solution to
\begin{equation}  \label{e:50:p:mainEne}
\begin{split}
      \partial_t \tilde\rho^{(n)} + u \cdot \nabla \tilde\rho^{(n)} - \kappa_n \Delta \tilde\rho^{(n)} =&\, 0,       \\ 
      \tilde\rho^{(n)} (\cdot, 0) =&\, \tilde\rho_{\ini}.
\end{split}
\end{equation}
Below $\psi$ and $\psi_n$ are the stream functions associated to $u$ and $u_n$ respectively, as defined in Section \ref{ss:vectorfield}. Note that from Lemma \ref{l:streamEst_rig_nobar} we have
\begin{align}
      \| \psi - \psi_n \|_{L^\infty} \leq&\, \frac{\delta_{n+1}^{\sfrac12}}{\lambda_{n+1}} 
            \leq \kappa_n \lambda_n^{-2b\gamma}.
            \label{e:52:p:mainEne}
\end{align}
Then we apply basic $L^2$-stability estimate (see e.g.\ Proposition 4.3 of \cite{burczak2023anomalous}) to deduce that
\begin{align}     \label{e:54:p:mainEne}
      \bigg| \kappa_n \iint | \nabla \tilde\rho^{(n)} |^2 dxdt - \kappa_n \iint | \nabla \tilde\rho_n |^2 dxdt \bigg| 
            \lesssim&\, \lambda_n^{-b\gamma} \kappa_n \iint | \nabla \tilde\rho_n |^2 dxdt.
\end{align}

Now we use \eqref{e:26:p:mainEne} for $q_0$, \eqref{e:44:p:mainEne} for $q=n$ and \eqref{e:54:p:mainEne} to deduce 
\begin{align}
      \kappa_n \iint | \nabla \tilde\rho^{(n)} |^2 dxdt
            \geq 2c_1 \kappa_{q_0} T \|\rho_{\ini}\|^2_{2}.     \label{e:56:p:mainEne}
\end{align}
Finally, we compare \eqref{e:6:p:mainEne} and \eqref{e:50:p:mainEne}. The basic energy estimate and \eqref{e:14:p:mainEne} yields
\begin{align}
      \kappa_n \iint | \nabla ( \tilde\rho^{(n)} - \rho_n ) |^2 dxdt
            \leq&\, \frac{1}{2} \| \tilde\rho_{\ini} - \rho_{\ini} \|_{2}^2
            \leq C^2 r^{2s} \ell^{2s-2} \|\rho_{\ini}\|_{2}^2.     \label{e:60:p:mainEne}
\end{align}
From \eqref{e:56:p:mainEne}, \eqref{e:60:p:mainEne} we have
\begin{align}
      \bigg( \kappa_n \iint | \nabla \rho_n |^2 dxdt \bigg)^{\frac{1}{2}}
            \geq&\, \big( (c_1 \kappa_{q_0} T)^{\sfrac12} - C r^{s} \ell^{s-1}  \big) \|\rho_{\ini}\|_{2}
            \geq \frac{1}{4} (c_1 \kappa_{q_0} T)^{\sfrac12} \|\rho_{\ini}\|_{2},     \label{e:62:p:mainEne}
\end{align}
where the latter inequality of \eqref{e:62:p:mainEne} is valid by imposing in \eqref{e:30:p:mainEne} the second condition on $c$.

Notice that \eqref{e:62:p:mainEne} is independent of $n \geq q_0$, hence we finish the proof of \eqref{e:2:t:mainHld}.
\end{step}
Proof of Theorem \ref{t:mainHld} is now complete.

\newpage 

\section{Appendix}

\subsection{Littlewood-Paley decomposition}

Here we collect some standard estimates. We refer to Chapter 2 in \cite{BaChDa11} and Chapter 3 in \cite{GuPer15} for Littlewood-Paley theory and Besov spaces. For a function $f$ defined on $\mathbb{T}^d$, we denote its Littlewood-Paley decomposition by
\begin{align}     \label{e:littlewood_paley}
      f = \sum_{j \geq -1} \Delta_j f.
\end{align}

Below we cite the characterization of H\"older spaces by Littlewood-Paley decomposition, see page 38 in \cite{GuPer15}.

\begin{lemma}[Littlewood-Paley decomposition]     \label{l:a:littlewood_paley}
For the decomposition \eqref{e:littlewood_paley}, we have the following estimates for any $m \in \N$ and $1 \leq p \leq \infty$,
\begin{align*}             
      &\, \| \nabla^m \Delta_j f \|_{L^p} \lesssim 2^{mj} \| \Delta_j f \|_{L^p}.
\end{align*}
Moreover, the following characterization of H\"older spaces holds for any $\alpha \in (0,1)$,
\begin{align*}
      \| f \|_{C^\alpha} \sim \sup_{j \geq -1} 2^{j\alpha} \| \Delta_j f \|_{ L^\infty }.
\end{align*}
\end{lemma}

\subsection{Interpolation inequality}

We also cite the following Gagliardo-Nirenberg interpolation inequality without proof.

\begin{lemma}[Gagliardo-Nirenberg interpolation inequality]     \label{l:a:interpolate}
For any $m \in \N$, we have that for any $f \in C^\infty (\mathbb{T}^n)$,
\begin{align}
      \| \nabla f \|_{L^\infty} \lesssim&\, \| \nabla^m f \|_{L^2}^{a_0} \| f \|_{L^\infty}^{1-a_0},
            \quad a_0 = \frac{2}{2m-n}, m > \frac{n}{2},          \label{e:l:2:interpolate} \\ 
      \| \nabla f \|_{L^\infty} \lesssim&\, \| \nabla^m f \|_{L^2}^{a_1} \| f \|_{L^2}^{1-a_1},
            \quad a_1 = \frac{n+2}{2m}, m \geq \frac{n+2}{2},     \label{e:l:4:interpolate} \\ 
      \| f \|_{L^\infty} \lesssim&\, \| \nabla^m f \|_{L^2}^{a_2} \| f \|_{L^2}^{1-a_2},
            \quad a_2 = \frac{n}{2m}, m > \frac{n}{2}.            \label{e:l:6:interpolate}
\end{align}
Moreover, we have
\begin{align}
      [ f ]_{C^\alpha} \lesssim&\, \| \nabla^m f \|_{L^2}^{a_3} \| f \|_{L^2}^{1-a_3},
            \quad a_3 = \frac{n+2\alpha}{2m}, m > \alpha + \frac{n}{2}.       \label{e:l:8:interpolate} \\ 
      [ f ]_{C^\alpha} \lesssim&\, \| \nabla^m f \|_{L^2}^{a_4} \| f \|_{L^\infty}^{1-a_4},
            \quad a_4 = \frac{2\alpha}{2m-n}, m > \alpha + \frac{n}{2}.       \label{e:l:10:interpolate}
\end{align}
\end{lemma}

\subsection{Mollification}

For $N_m \in \N^+$, we need mollifiers which are deep smoothing operators of depth $N_m$, see \cite[Section 2.3.4]{GromovBook}. The case of infinite depth was introduced by Nash \cite{Nash56}. For this purpose, we fix a mollification kernel $\mathfrak{m} \in C^\infty_c(\R^n)$ which, in addition to $ \int_{\R^n} \mathfrak{m} \,dx = 1 $, also satisfies 
\begin{equation}\label{e:deepmollifier}
      \int_{\R^n} \mathfrak{m}(x) x^{\balpha} \,dx = 0,
            \quad \textrm{ for any multiindex $\balpha$ with $1 \leq |\balpha| \leq N_m$ }.
\end{equation}
If $N_m\geq 2$, then $\mathfrak{m}$ cannot be nonnegative. As a variant of Lemma 7.2 in \cite{burczak2023anomalous}, we have the following mollification estimates.

\begin{lemma}[Mollification]     \label{l:a:mollify}
Let  $\mathfrak{m} \in C^\infty_c(\R^n)$ be a smoothing operator of depth $N_m \geq 1$ and such that $\int_{\R^n} \mathfrak{m} \,dx = 1$. Then for any real $r,s\geq 0$ and any $f \in C^\infty (\mathbb{T}^n)$
\begin{equation}
	\| f * \mathfrak{m}_\ell \|_{H^{r+s}} \lesssim \ell^{-s} \|f\|_{H^r}    \label{e:mollify1}
\end{equation}
and for any $r\geq 0$, $0\leq s\leq N_m$
\begin{equation}
      \|f-f*\mathfrak{m}_\ell\|_{H^r} \lesssim \ell^{s} \|f\|_{H^{r+s}}   \label{e:mollify2} \\
\end{equation}
The absolute constants depend on $r,s$ and $\mathfrak{m}$.
\end{lemma}

The proof is similar to that of Lemma 7.2 in \cite{burczak2023anomalous}.

\bibliographystyle{plain}
\bibliography{scalar.bib}

\end{document}